%% file: Book-nonlocal.tex
\documentclass{gsm-l}

\usepackage{color}
\usepackage{framed}
\usepackage{empheq}
\usepackage{enumerate}
\usepackage{amsmath, amsthm, amssymb}
\usepackage{hyperref}
\usepackage{array}
\usepackage{mathrsfs} 
\usepackage[shortlabels]{enumitem}

\usepackage{background}
\backgroundsetup{
  position=current page.east,
  angle=-90,
  vshift=-4mm,
  hshift=2mm,
  opacity=0.8,
  scale=3,
  contents= {\rm --- $\overline{\underline{\mbox{DRAFT}}}$ --- }
}

\usepackage{chngcntr}
\counterwithin{figure}{section}

\newlist{steps}{enumerate}{1}
\setlist[steps, 1]{label = {\it \underline{\smash{Step~\arabic*:}}}, ref={\it Step~\arabic*}, wide, itemindent = 3pt}

\newtheorem{thm}{Theorem}[section]
\newtheorem*{thm*}{Theorem}
\newtheorem{lem}[thm]{Lemma}
\newtheorem{cor}[thm]{Corollary}
\newtheorem{prop}[thm]{Proposition}

\theoremstyle{definition}
\newtheorem{defi}[thm]{Definition}
\newtheorem{example}[thm]{Example}

\theoremstyle{remark}
\newtheorem{rem}[thm]{Remark}
\newtheorem*{rem*}{Remark}

\newtheorem*{claim*}{Claim}

\numberwithin{section}{chapter}
\numberwithin{equation}{section}

\newcommand{\ssubset}{\subset\joinrel\subset}

\renewcommand{\tocsection}[3]{\indentlabel{{\ignorespaces#1\makebox[\widthof{000.}][l]{#2.}\quad}}#3\dotfill}

\newcommand{\ubar}[1]{\text{\b{$#1$}}}

\newcommand{\liminfs}[1]{\underset{{\smash{#1}}}{{\rm lim\,inf}_*}\,}
\newcommand{\limsups}[1]{\underset{{\smash{#1}}}{{\rm lim\,sup}^*}\,}

\newcommand{\tosu}{\xrightarrow{\vphantom{*}\smash{\raisebox{-1ex}{{\rm *}}}}}
\newcommand{\tosd}{\xrightarrow[\vphantom{*}\smash{\raisebox{-0.1ex}{{\rm *}}}]{}}

\newcommand{\average}{{\mathchoice {\kern1ex\vcenter{\hrule height.4pt
width 6pt depth0pt} \kern-9.7pt} {\kern1ex\vcenter{\hrule
height.4pt width 4.3pt depth0pt} \kern-7pt} {} {} }}
\newcommand{\ave}{\average\int}

\def\fls{{(-\Delta)^s}}
\newcommand{\sqrtl}{{\sqrt{-\Delta}\,}}
\def\be{{\boldsymbol{e}}}

\def\R{\mathbb{R}}
\def\C{\mathbb{C}}
\def\Z{\mathbb{Z}}
\def\N{\mathbb{N}}
\def\E{\mathbb{E}}
\def\L{\mathcal{L}}
\def\MR{\mathcal{R}}
\def\LL{\mathfrak{L}}
\def\LLL{\LL_s(\lambda, \Lambda)}
\def\I{\mathcal{I}}
\def\J{\mathcal{J}}
\def\II{\mathfrak{I}}
\def\III{\II_s(\lambda, \Lambda)}
\def\IIL{\mathfrak{I}_s(\lambda, \Lambda)}
\def\A{\mathcal{A}}
\def\B{\mathcal{B}}
\def\D{\mathcal{D}}
\def\dr{\mathfrak{d}}
\def\M{\mathcal{M}}
\def\MMp{\M^+_{\LLL}}
\def\MMm{\M^-_{\LLL}}
\def\MMpm{\M^\pm_{\LLL}}

\DeclareRobustCommand{\S}{\ifmmode\mathsection\else\textsection\fi}
\renewcommand\mathsection{\operatorname{\mathbb{S}}}

\newcommand{\eps}{\varepsilon}
\newcommand{\de}{\partial}

\DeclareMathOperator{\dist}{dist}

\newcommand{\G}{{\mathfrak{G}}}
\newcommand{\GL}{{\mathfrak{G}_s(\lambda, \Lambda)}}

\newcommand{\GLh}{{\mathfrak{G}^{\hspace{-0.3mm} \rm \scriptscriptstyle h\hspace{-0.1mm}o\hspace{-0.1mm}m}_s \hspace{-0.4mm} (\lambda, \Lambda)}}

\newcommand{\LLh}{{\mathfrak{L}^{\hspace{-0.3mm} \rm \scriptscriptstyle h\hspace{-0.1mm}o\hspace{-0.1mm}m}_s \hspace{-0.4mm} (\lambda, \Lambda)}}

\newcommand{\Mp}{{\M^+}}
\newcommand{\Mm}{{\M^-}}
\newcommand{\Mpm}{{\M^\pm}}

\DeclareMathOperator*{\osc}{osc}

\makeindex

\begin{document}
\frontmatter
\title{\hspace{-3mm}{\Huge \mbox{\hspace{-20mm}Integro-Differential Elliptic~Equations}}\hspace{-5mm}\\ \vspace{1cm}}

\author{Xavier Fern\'andez-Real}
\address{EPFL SB MATH, Institute of Mathematics, Station 8, CH-1015 Lausanne, Switzerland}
\email{xavier.fernandez-real@epfl.ch}

\author{Xavier Ros-Oton}
\address{ICREA, Pg. Lluis Companys 23, 08010 Barcelona, Spain; and}
\address{Universitat de Barcelona, Gran Via de les Corts Catalanes 585, 08007 Barcelona, Spain; and}
\address{Centre de Recerca Matem\`atica, Barcelona, Spain.}
\email{xros@icrea.cat}

\makeatletter
\@namedef{subjclassname@2020}{%
  \textup{2020} Mathematics Subject Classification}
\makeatother

\date{\today}

\begin{abstract}
This book aims to provide a self-contained introduction to the regularity theory for integro-differential elliptic equations, mostly developed in the 21st century.
Such a class of equations often arises in analysis, probability theory, mathematical physics, and in several contexts in applied sciences.
The authors give a detailed presentation of all the necessary techniques, primarily focusing on the main ideas rather than proving all results in their greatest generality.

The book starts from the very basics, studying the square root of the Laplacian and weak solutions to linear equations.
Then, the authors develop the theory of viscosity solutions to nonlinear equations and prove the main known results in this context.
Finally, they study obstacle problems for integro-differential operators and establish the regularity of solutions and free boundaries.

Almost all the covered material appears in book form for the first time, and several proofs are different (and shorter) than those in the original papers. Moreover, several open problems are listed throughout the book. 
\end{abstract}

\maketitle

\setcounter{page}{2}
\setcounter{tocdepth}{2}
\tableofcontents

\include{preface}

\mainmatter
\include{chap1}

\include{chap2}

\include{chap3}

\include{chap4}

\appendix
\include{appendixA}

\include{appendixB}

\backmatter
\include{notation}

\include{biblio}

\printindex

\end{document}

%% file: preface.tex
%

\chapter*{Preface}

The aim of this book is to study nonlocal equations of the form
\begin{equation}\label{eq-preface} \tag{$\ast$}
\int_{\R^n}\bigl(u(x)-u(x+y)\bigr)K(y)dy=0\qquad \textrm{for}\quad x\in\Omega\subset\R^n,
\end{equation}
with kernels $K\geq0$.
These equations are \emph{elliptic}, and they share many properties with elliptic PDE: maximum principle, existence and uniqueness of solutions, regularity results, etc.

The equations we will study are of   type \eqref{eq-preface}, with a kernel $K(y)$ that is \emph{not integrable at the origin}, while nice and integrable at infinity.
The first and simplest example is 
\[K(y)=\frac{1}{|y|^{n+1}}.\]
In this particular case, the operator in \eqref{eq-preface} is (a multiple of) the square root of the Laplacian $\sqrt{-\Delta}$.

Since the kernel $K$ is not integrable at the origin, the operator in \eqref{eq-preface} is in some sense differentiating the function $u$, and this is why it is called an \emph{integro-differential} equation.
Such a type of equations appears in several contexts, and has been studied for many years in:
\begin{itemize}
 \item Probability theory (stochastic processes with jumps).
 \item Fluid mechanics (for example, in the SQG equation, or even the Boltzmann equation).
 \item Mathematical physics (relativistic Schr\"odinger operators).
 \item Applied sciences (anomalous diffusions).
\end{itemize}
From an analytical point of view, the general theory for such  equations has been mostly developed in the 21st century, including both \emph{linear} and \emph{nonlinear} equations.

This area of research has attracted great interest in the PDE community in the last 15 years, especially since the first works of Caffarelli and Silvestre on the topic.
However, there is no book yet providing a systematic study of the (existence and) regularity properties of solutions to general integro-differential elliptic equations.

The goal of this book is to fill this void by developing the regularity theory for (both linear and nonlinear) integro-differential elliptic operators of order $2s$, with $s\in(0,1)$.
As we will see, there is a strong parallelism with the theory of elliptic PDE, which corresponds to the limiting case $s=1$.
Still, quite often the proofs of the results for $s\in(0,1)$ are completely independent of the ``local case'' $s=1$, and new ideas and techniques have been developed in order to treat such nonlocal equations.

Most of the material we present here is accessible in book form for the first time, and we provide several simplified proofs compared to the original papers.
Furthermore, we establish some new results as well, which generalize or unify previously known results.
For example, the interior regularity theory for linear operators that we develop in Chapter~\ref{ch:2} includes for the first time the most general scale-invariant class of operators of order $2s$, the existence and uniqueness of viscosity solutions in Chapter~\ref{ch:fully_nonlinear} is established for a \emph{very} wide class of domains, and the boundary Harnack inequality in Chapter~\ref{ch:obst_pb} is proved here under a very mild assumption on the domain $\Omega$.

We hope that this book will be useful to some of the many active researchers working in this field, while at the same time being self-contained and accessible to graduate students interested in this topic.

We wish to thank Gerd Grubb,  Joaquim Serra, and Marvin Weidner for several comments and suggestions on this book.

Finally, we acknowledge the support from the following funding agencies. 
X.F. was supported by the Swiss National Science Foundation (SNF grants 
200021\_182565 and PZ00P2\_208930), and by the Swiss State Secretariat for Education, Research and lnnovation (SERI) under contract number MB22.00034;
X.R. was supported by the European Research Council under the Grant Agreement No. 801867 ``Regularity and singularities in elliptic PDE (EllipticPDE)'', by AEI project PID2021-125021NA-I00 (Spain),  by the grant RED2022-134784-T funded by AEI/10.13039/501100011033, by AGAUR Grant 2021 SGR 00087 (Catalunya), and by the Spanish State Research Agency through the Mar\'ia de Maeztu Program for Centers and Units of Excellence in R{\&}D (CEX2020-001084-M).


%% file: chap1.tex
%
%
%


\chapter{The square root of the Laplacian}
\label{ch:fract_Lapl}
\addtocontents{toc}{\protect\setcounter{tocdepth}{1}}

In this chapter we focus our attention on the square root of the Laplacian, denoted $\sqrtl$. This operator, when acting on a smooth function $u\in C^\infty_c(\R^n)$, can be defined as follows:\index{Square root of the Laplacian}
\[
\begin{split}
{\sqrtl} u(x) & := c_n \textrm{P.V.}\int_{\R^n} \frac{u(x) - u(x+y)}{|y|^{n+1}}\, dy,
\end{split}
\]
for some positive constant $c_n$.

\section{Motivation}

Diffusive processes form the backbone of the quantitative study of many-particle systems. 
The mathematical description of the Brownian motion has a vast range of applications that basically touch any area in which a macroscopic description of microscopic phenomena is needed.
In the last decades, however, there has been increasing observational evidence that many systems that were previously thought to behave like ordinary diffusion processes are actually better explained by assuming an anomalous behavior of the forming particles. 
Namely, the Brownian motion is formulated under the assumption that the movement of such particles is \emph{continuous} or, alternatively, that in short time intervals it can be described (thanks to the central limit theorem) by means of Gaussian distributions. 
However, such a description fails to account for the situations in which the underlying random variables have infinite variance. 
This is precisely the setting where \emph{anomalous diffusion} appears, the \emph{Cauchy process} being  a first important example. 
In this setting, particles can \emph{jump} towards another point and the increments of the process follow a Cauchy distribution (see Section~\ref{sec:Levy} for more details).
Macroscopically, the evolution of the density of particles is not only characterized by its local profile, but rather it is influenced by the whole distribution (being further objects less relevant).
In terms of the corresponding transition function, the operator governing the evolution of a given density is no longer a local operator (contrary to the Brownian motion, where we find the Laplacian $-\Delta$), but instead we find the \emph{square root of the Laplacian}, denoted by $\sqrtl$, a nonlocal operator. 
More precisely, the evolution of a density of particles is governed by the fractional heat equation: 
\[
\partial_t u + \sqrtl u = 0\quad\text{in}\quad \R^n. 
\]
This type of nonlocal equations have received an increasing amount of attention in the last decades, mainly motivated and driven by numerous applications: starting from the observation of Mandelbrot in the 1960s, on the deviation of stock prices from normality, \cite{Man63}, they also appear in \emph{physics},  \cite{MK, MK2, Las00, ZK94}, \emph{ecology} and \emph{biology}, \cite{RR09, Hum10, Kla90, Vis96}, and \emph{finance}, \cite{Merton, Sch03, AB88}. 

As we will see, the square root of the Laplacian $\sqrtl$ can also be interpreted as the Dirichlet-to-Neumann map for the Laplace equation in the half-space (see Section~\ref{sec:harm_extension}). 
As such, this operator appears in \emph{fluid mechanics}, for example in the surface quasi-geostrophic equation, which is used to model the temperature on the surface of a fluid in oceanography, \cite{Con,CCV,CV}, or the Benjamin-Ono equation \cite{Ben67, Ono75},
\[
\sqrtl u = -u+u^2\quad\text{in}\quad \R,
\]
used to describe one-dimensional internal waves in deep water, \cite{AT91, FL}, and which plays an important role in the understanding of the gravity of water waves' equations in dimensions 2 and 3 \cite{GMS12}. 

The same operator $\sqrtl$ arises in \emph{elasticity}, too: the Peierls-Nabarro equation is a model in crystal dislocation to study microscopical deformations of a material, \cite{BV, Tol97, Lu05, CS-M, DPV15}; and the Signorini or thin obstacle problem can be used to model the equilibrium configuration of an elastic body on a frictionless surface \cite{Sig59, F-survey} (see Section~\ref{sec:motivations}). 
More generally, we find  nonlocal operators such as $\sqrtl$ when trying to model phenomena that takes into account long-range interactions of a given system: in many materials the stress points depend on the strains of surrounding regions, \cite{Kro67, Eri72}, and nonlocal forces have been observed to propagate along fibers or laminae in composite materials \cite{INW81,Dru03, MD04}. 

On the other hand, in \emph{quantum mechanics}, nonlocal operators like $\sqrtl$ appear as the relativistic momentum operator. 
In particular, this operator arises in dispersive equations describing the dynamics and gravitational collapse of boson stars, \cite{LY87, ES07}, or in the study of stability of relativistic matter, \cite{FLS08,Lieb}; see also  \cite{FL,Las00, Las02, Sei09}.

This kind of operators appears as well in \emph{image processing}, where nonlocal denoising algorithms are able to detect patterns and contours in a better way than PDE-based models, and can be used for image reconstruction and enhancement \cite{Yar85, GO08, ELB08, KFEA09, Zha10, DSL12, Yue16}.

Fractional powers of the Laplacian also appear naturally in \emph{conformal geometry}, where the fractional Paneitz operators are conformally covariant operators which encode geometric and topological information about the manifold; see \cite{geometry0, geometry1,Geometry2,Mar}.

Finally, we refer to \cite{mwiki,CR,CRSire,CRSavin,CSVaz,DG,DGLZ,Eri02,PK,Serfaty1,IS,IS2,IS3,Oskendal,SiICM,S-Boltzmann,Vaz,ZD} for other models and motivations to study nonlocal operators like the square root of the Laplacian.

We next start studying the basic properties of $\sqrt{-\Delta}$, which will serve as an introduction to integro-differential operators and their corresponding elliptic equations.

\section{Basic properties}

The square root of the Laplacian may be defined as
\begin{equation}
\label{eq:square_root_Laplacian}\index{Square root of the Laplacian}
\begin{split}
{\sqrtl} u(x) & := c_n \textrm{P.V.}\int_{\R^n} \frac{u(x) - u(x+y)}{|y|^{n+1}}\, dy \\
& =  c_n \textrm{P.V.}\int_{\R^n} \frac{u(x) - u(z)}{|x-z|^{n+1}}\, dz,
\end{split}
\end{equation}
where $c_n$ is a positive constant, given by \eqref{eq:cn} below.

 The notation $\mathrm{P.V.}$ stands for \emph{Principal Value}, which is a way (coming from distribution theory) to assign values to an integral of an a priori not absolutely integrable function. Indeed, observe that, when $u$ is smooth, $u(x) - u(x+y) = \nabla u (0)\cdot y + O(|y|^2)$. If $\nabla u$ is nonzero, then $|u(x) - u(x+y)|\cdot |y|^{-n-1}$ is comparable to $\frac{1}{|y|^n}$, which is not integrable at the origin! That is, the function to be integrated\footnote{Integrability at infinity is not a problem at this point, we can assume  $u$ to be compactly supported.} is not~$L^1$.

The integral is then taken in principal value sense, which takes advantage of cancellations to assign a value to the integral:
\[\index{Principal Value}
\textrm{P.V.}\int_{\R^n} \frac{u(x) - u(x+y)}{|y|^{n+1}}\, dy := \lim_{\eps \downarrow 0} \int_{\R^n\setminus B_\eps} \frac{u(x) - u(x+y)}{|y|^{n+1}}\, dy.
\]
Notice that, with this definition, by symmetry we have 
\[
{\rm P.V.} \int_{B_1}\frac{\nabla u \cdot y }{|y |^{n+1}}\, dy = 0,\]
 and hence the integral in \eqref{eq:square_root_Laplacian} is well defined for any $u\in C^\infty_c(\R^n)$. 

The square root of the Laplacian \eqref{eq:square_root_Laplacian} is a \emph{nonlocal} operator: when acting on $u$ at a point $x\in \R^n$, it uses information about $u$ that is far from~$x$. Also, we have observed that the kernel $\frac{1}{|y|^{n+1}}$ is singular (not integrable) at the origin, and hence it requires a certain regularity from $u$ near $x\in \R^n$ in order to evaluate $\sqrtl u(x)$. The singularity of the kernel makes it such that the operator ${\sqrtl}$ is ``differentiating'', in some sense, the function $u$, and this is why it is called an \emph{integro-differential} operator. 

Sometimes, it is convenient to write the following alternative expression for~\eqref{eq:square_root_Laplacian}:
\begin{equation}
\label{eq:square_root_Laplacian_incr_quotients}
{\sqrtl} u(x) = \frac{c_n}{2}\int_{\R^n} \frac{2 u(x) - u(x+y)-u(x-y)}{|y|^{n+1}}\, dy.
\end{equation}
Now there is no need to write $\textrm{P.V.}$, since the numerator of the integrand is a second order incremental quotient, and thus it has norm comparable to $|y|^2$ (and $|y|^{-n+1}$ is integrable in $\R^n$ around the origin). In particular, this expression implies that in order to evaluate ${\sqrtl}$ it is enough to assume $u\in C^2(\R^n)\cap L^\infty(\R^n)$:
\[
\begin{split}
\big|{\sqrtl} u(x)\big|& \le \frac{c_n}{2}\int_{B_1}\frac{\|D^2u\|_{L^\infty(\R^n)}}{|y|^{n-1}}\, dy + \frac{c_n}{2}\int_{\R^n\setminus B_1} \frac{4\|u\|_{L^\infty(\R^n)}}{|y|^{n+1}}\, dy \\
& \le C \left(\|u\|_{L^\infty(\R^n)}+\|D^2 u\|_{L^\infty(\R^n)}\right),
\end{split}
\]
for some $C$ depending only on $n$. The regularity of $u$ takes care of the integrability of the kernel around the origin.

On the other hand, when evaluating the square root of the Laplacian at a point $x\in\R^n$, it is also essential to have a global integrability assumption on $u$ to account for the contributions far from $x$. In the previous inequality, this global integrability assumption is given by the fact that $u\in L^\infty(\R^n)$. More generally, though, it is enough to consider the following space:
\begin{defi}
\label{defi:L1w}\index{L1w@$L^1_\omega$}
We say that $v\in L^1_\omega(\R^n)$ if 
\[
\|v\|_{L^1_\omega(\R^n)} := \int_{\R^n}\frac{|v(y)|}{1+|y|^{n+1}}\, dy < +\infty.
\]
\end{defi}
\begin{rem}
The previous definition introduces the minimum global integrability requirement so that ${{\sqrtl}} u$ makes sense. Indeed, we can split into two terms, 
\[
\begin{split}
{\sqrtl} u(x)  & = c_n \textrm{P.V.}\int_{B_1} \frac{u(x) - u(x+y)}{|y|^{n+1}}\, dy + c_n \int_{\R^n\setminus B_1} \frac{u(x) - u(x+y)}{|y|^{n+1}}\, dy.
\end{split}
\] 
Then, in order to evaluate the first term we just need a local regularity assumption on $u$ around $x$, whereas the second term is bounded (in absolute value) by $\|u\|_{L^1_\omega(\R^n)}$. 
\end{rem}

On the other hand, if we want to evaluate ${\sqrtl} u$ pointwise at a point~$x$, it is enough for $u$ to be $C^{1,\alpha}$ around $x$, for some $\alpha > 0$:

\begin{lem}
\label{lem:laplu}\index{Strong solutions!Square root of the Laplacian}
Let $u\in C^{1,\alpha}(B_1)\cap L^1_\omega(\R^n)$ for some $\alpha\in(0,1)$. Then ${\sqrtl}u(0)$ is well-defined, and
\begin{equation}
\label{eq:flsat0}
\big|{\sqrtl}u (0)\big| \le C \left(\|u\|_{C^{1,\alpha}(B_1)} + \|u\|_{L^1_\omega(\R^n)}\right)
\end{equation}
for some $C$ depending only on $n$ and $\alpha$. Moreover, ${\sqrtl}u \in C_{\rm loc}^{0,\alpha}(B_1)$ and 
\[
\big\|{\sqrtl} u\big\|_{C^{0,\alpha}(B_{1/2})} \le C \left(\|u\|_{C^{1,\alpha}(B_1)} + \|u\|_{L_\omega^1(\R^n)} \right) 
\]
for some $C$ depending only on $n$ and $\alpha$.
\end{lem}

\begin{proof}  By considering $v = u/{C_\circ}$ with $C_\circ = \|u\|_{C^{1,\alpha}(B_1)} + \|u\|_{L^1_\omega(\R^n)}$ we may assume that 
\begin{equation}
\label{eq:weassume1}
\|u\|_{C^{1,\alpha}(B_1)} + \|u\|_{L^1_\omega(\R^n)} = 1.
\end{equation} 
 
 Moreover, this implies 
\[
\big|2u(x)-u(x+y)-u(x-y)\big|\le 2|y|^{1+\alpha}
\]
for all $x\in B_1$ and $y\in B_{1-|x|}$ (see \eqref{eq:proc_as_in}).

From \eqref{eq:square_root_Laplacian_incr_quotients}, and integrating separately in $B_1$ and $\R^n\setminus B_1$ (recall \eqref{eq:weassume1}),
\[
\big|{\sqrtl} u(0)\big| \le C\int_{B_1}\frac{dy}{|y|^{n-\alpha}} + C\int_{\R^n\setminus B_1} \frac{|u(0)|+|u(y)|+|u(-y)|}{|y|^{n+1}}\, dy.
\]

Observe that now, on the one hand (again using \eqref{eq:weassume1}),
\[
\int_{\R^n\setminus B_1}\frac{|u(0)|}{|y|^{n+1}}\, dy \le C\int_{\R^n\setminus B_1}\frac{dy}{|y|^{n+1}}\le C 
\]
for some $C$ depending only on $n$. And on the other hand, 
\[
\int_{\R^n\setminus B_1}\frac{|u(y)|+|u(-y)|}{|y|^{n+1}}\, dy = 2\int_{\R^n\setminus B_1}\frac{|u(y)|}{1+| y|^{n+1}} \frac{1+| y|^{n+1}}{|y|^{n+1}}\, dy \le 4,
\]
where we are using that $1+|y|^{n+1}\le 2|y|^{n+1}$ since $|y|\ge 1$, and $\|u\|_{L^1_\omega(\R^n)}\le 1$. Putting everything together, we get \eqref{eq:flsat0}. 

More generally, repeating for any point $z\in B_{1/2}$ (and taking integrals in $B_{1/2}$ instead of $B_1$) we obtain
\begin{equation}
\label{eq:sqrtlLinf}
\big\|{\sqrtl} u\big\|_{L^\infty(B_{1/2})} \le C
\end{equation}
as we wanted. 

Let us now take a cut-off function $\eta\in C^\infty_c(\R^n)$ such that $\eta \ge 0$, $\eta \equiv 0$ in $\R^n \setminus B_{3/4}$ and $\eta\equiv 1$ in $B_{2/3}$; and let us define $u_1 := \eta u$ compactly supported in $B_{3/4}$, and $u_2 := u-u_1 = (1-\eta) u$ such that $u_2 \equiv 0$ in $B_{2/3}$. Let us obtain an estimate for each $u_1$ and $u_2$ separately. 

We prove first that $\sqrtl u_1\in C^{0,\alpha}(B_{1/2})$, with $\|\sqrtl u_1\|_{C^{0,\alpha}(B_{1/2})}\le C$ (where $C$ is universal, explicit in terms of the cut-off $\eta$). Observe that the $L^\infty$ bound follows as in \eqref{eq:sqrtlLinf}, since $u_1\in C^{1,\alpha}(\R^n)$ with $\|u_1\|_{C^{1,\alpha}(\R^n)}\le C\|u\|_{C^{1,\alpha}(B_1)}\le C$. If we let now $x\in B_{1/2}$ and $r := |x|$, we split
\[
\begin{split}
\sqrtl u_1(x) & = \frac{c_n}{2} \int_{B_r} \frac{2u_1(x) -u_1(x+y)-u_1(x-y)}{|y|^{n+1}}\, dy \\
& \quad +c_n  \int_{\R^n\setminus B_r} \frac{u_1(x) -u_1(x+y)}{|y|^{n+1}}\, dy,
\end{split}
\]
so that 
\[
\begin{split}
& \left|\sqrtl u_1(x)+\sqrtl u_1(-x) - 2\sqrtl u_1(0)\right| \leq  \\
& \qquad \le  C\int_{B_r} |y|^{\alpha-n}\, dy  + C\int_{\R^n\setminus B_r} \frac{|x|^{1+\alpha}}{|y|^{n+1}}\, dy \le  C r^\alpha,
\end{split}	
\]
where we used 
\[
\begin{split}
|u_1(x+y)+u_1(x-y)-2u_1(x)| &\le  C |y|^{1+\alpha},\\
|u_1(x)+u_1(-x)-2u_1(0)| &\le  C |x|^{1+\alpha},\\
|u_1(x+y)+u_1(-x+y)-2u_1(y)| &\le  C |x|^{1+\alpha}.
\end{split}
\]
 Repeating it around any point in $B_{1/2}$ and together with the $L^\infty$ bound, we get  $\|\sqrtl u_1\|_{C^{0,\alpha}(B_{1/2})} \le C$ by Lemma~\ref{it:H7}.

Secondly, let us consider $u_2$. 
Given any $x_1, x_2 \in B_{1/2}$, since $u_2 \equiv 0$ in~$B_{2/3}$,  
\[
\begin{split}
& \left|\sqrtl  u_2(x_1)-\sqrtl u_2(x_2)\right| =  \\
&\qquad  =  c_n \left|\int_{D} u_2(z) \left(\frac{1}{|z-x_1|^{n+1}}-\frac{1}{|z-x_2|^{n+1}}\right)\, dz\right|,
\end{split}
\]
where $D$ can be taken to be $\R^n\setminus (B_{1/6}(x_1)\cup B_{1/6}(x_2))$. Observe that, in this domain, 
\[
\left|\frac{1}{|z-x_1|^{n+1}}-\frac{1}{|z-x_2|^{n+1}}\right|\le C\frac{|x_1-x_2|}{1+|z|^{n+2}}\quad\text{for all}\quad z \in D,
\]
by regularity of $|z|^{-n-1}$ outside of the origin. In particular, 
\[
\big|\sqrtl  u_2(x_1)-\sqrtl u_2(x_2)\big|   \le C |x_1-x_2|\|u_2\|_{L^1_\omega(\R^n)},
\]
and hence $[\sqrtl u_2]_{C^{0, 1}(B_{1/2})} \le C \|u_2\|_{L^1_\omega(\R^n)}\le C \|u\|_{L^1_\omega(\R^n)} \le C$.

Combining the estimates for $u_1$ and $u_2$ with \eqref{eq:sqrtlLinf}, we get the desired result. 
\end{proof}

We also observe that the operator ${\sqrtl}$ is positive definite:
\begin{lem}\index{Integration by parts!Square root of the Laplacian}
\label{lem:pos_def}
Let $u, v \in C^\infty_c(\R^n)$. Then 
\[
\int_{\R^n} u(x) {\sqrtl} v(x)\, dx = \frac{c_n}{2}\int_{\R^n}\int_{\R^n} \frac{\left(u(x) - u(z)\right)\left(v(x)-v(z)\right)}{|x-z|^{n+1}}\, dx\, dz.
\]
In particular, 
\[
\int_{\R^n} u(x) {\sqrtl} u(x)\, dx = \frac{c_n}{2}\int_{\R^n}\int_{\R^n} \frac{\left(u(x) - u(z)\right)^2}{|x-z|^{n+1}}\, dx\, dz \ge 0,
\]
 with strict inequality unless $u\equiv 0$. That is, ${\sqrtl}$ is positive definite. 
\end{lem}
\begin{proof}
By definition, 
\[
\int_{\R^n} u(x) {\sqrtl} v(x)\, dx =c_n \int_{\R^n} u(x)  \textrm{P.V.}\int_{\R^n} \frac{v(x) - v(z)}{|x-z|^{n+1}}\, dz\, dx.
\]
Changing the roles of the variables we also have 
\[
\int_{\R^n} u(x) {\sqrtl} v(x)\, dx =c_n \int_{\R^n} u(z)  \textrm{P.V.}\int_{\R^n} \frac{v(z) - v(x)}{|x-z|^{n+1}}\, dx\, dz.
\]
Adding the previous two expressions and dividing by two, we obtain the desired result.
\end{proof}

 Let us finish this short introduction to the square root of the Laplacian noticing the invariances satisfied by the operator, and how these invariances actually characterize it: 
 \begin{lem}
 \label{lem:invariances}
 Let $u\in C^2(\R^n)\cap L^\infty(\R^n)$. Then, the following properties hold for all $x\in \R^n$,
 \begin{enumerate}[(i)]
 \item \label{it:prop_trans}\index{Translation invariance!Square root of the Laplacian}Translation invariance: \[
 {\sqrtl} [u(x+x_\circ)] = \big({\sqrtl} u\big)(x+x_\circ)
 \] for any $x_\circ\in \R^n$. 
 \item \index{Rotation invariance!Square root of the Laplacian}\label{it:prop_rot} Rotation/Orthogonal invariance: \[{\sqrtl} [u(Ox)] = \big({\sqrtl} u\big)(Ox)\]
  for any $O\in \mathcal{O}(n)$ orthogonal transformation.
 \item \index{Scale invariance!Square root of the Laplacian}\label{it:prop_scale}Scale invariance/1-homogeneity: 
 \[{\sqrtl} [u(\lambda x)] = |\lambda|\big({\sqrtl} u\big)(\lambda x)\] for any $\lambda \in \R$.
 \end{enumerate}
 \end{lem}
\begin{proof}
It follows using the definition of  ${\sqrtl}$, \eqref{eq:square_root_Laplacian}, and changing variables appropriately for each case. 
\end{proof}

In particular, estimates like the ones in Lemma~\ref{lem:laplu} can be rescaled to any ball $B_r(x_\circ)\subset \R^n$. 

The properties in Lemma~\ref{lem:invariances} are not mere consequences of the definition of the operator ${\sqrtl}$, but in fact, they fully characterize it. In the following lemma, the set $\mathcal{S}(\R^n)$ denotes the \emph{Schwartz space}, that is, the set of $C^\infty(\R^n)$ functions that have all derivatives decaying faster than any power. We also denote by $\mathcal{F}$ the Fourier transform, acting on functions $f\in \mathcal{S}(\R^n)$ (in fact, it is an automorphism on $\mathcal{S}$):
\[
\mathcal{F}(f)(\xi) := \int_{\R^n} f(x) e^{-i\xi\cdot x}\, dx. 
\]

\begin{lem}
\label{lem:sqrtl_unique_invariances}\index{Fourier symbol!Square root of the Laplacian}
Let $\L:\mathcal{S}\to \mathcal{S}$ be any linear operator satisfying properties \ref{it:prop_trans}-\ref{it:prop_rot}-\ref{it:prop_scale} from Lemma~\ref{lem:invariances}. Then, there exists a constant $\kappa\in \R$ such that $\L = \kappa {\sqrtl}$.
\end{lem}
\begin{proof}
Since $\L$ is translation invariant by property \ref{it:prop_trans}, it is given by a Fourier multiplier $\A(\xi)$ (see, for example, \cite[Chapter I, Theorem 3.16]{SW71}). Namely, 
\[
\L u = \mathcal{F}^{-1} (\A(\xi)\mathcal{F}(u)(\xi))
\]
where $\mathcal{F}$ denotes the Fourier transform, and $\mathcal{F}^{-1}$ is the inverse Fourier transform.

Rotation and scale invariances, \ref{it:prop_rot}-\ref{it:prop_scale}, now imply that $\A(\xi) = \kappa |\xi|$. In particular, all operators satisfying \ref{it:prop_trans}-\ref{it:prop_rot}-\ref{it:prop_scale} from Lemma~\ref{lem:invariances} are multiple of each other, and the lemma follows. 
\end{proof}
\begin{rem}
\label{rem:sqrtlsqrtl}
As a consequence of the previous proof, we also showed that the Fourier multiplier of ${\sqrtl}$ is $\kappa |\xi|$ for some $\kappa\in \R$. In Section~\ref{sec:harm_extension} below we will show that ${\sqrtl}\circ {\sqrtl} = -\Delta$, and since the Fourier multiplier of $-\Delta$ is simply $|\xi|^2$, together with the positive definiteness of ${\sqrtl}$ we obtain that the Fourier multiplier of ${\sqrtl}$ is $|\xi|$ (namely, $\kappa = 1$). 
\end{rem}

\section{Harmonic extension}
\label{sec:harm_extension}

For any $u\in C^\infty_c(\R^n)$, let us consider the operator $Lu$ defined as follows: 

Let ${\tilde u}(x, y):\R^n\times\R\to  \R$ be the harmonic extension of $u(x)$ towards $\R^{n+1}_+ := \{(x, y)\in \R^n\times\R : y > 0\}$.  That is, $\tilde u(x, y)$ is the only solution to 
\[
\left\{
\begin{array}{rcll}
\Delta_{x, y} {\tilde u} & = & 0& \quad\text{in}\quad \{y > 0\}\\
{\tilde u}(x, 0) & = & u(x)& \quad\text{for}\quad x\in \R^n
\end{array}
\right.
\]
which decays at infinity, where $\Delta_{x, y} := \Delta_x + \partial_{yy} = \sum_{i = 1}^n \partial_{x_ix_i} + \partial_{yy}$ denotes the Laplacian in the $(x, y)$ coordinates. We then define
\[
L u(x) := -\partial_y\big|_{y = 0} {\tilde u}(x, y).
\]

The operator $L$ is called a \emph{Dirichlet-to-Neumann operator}: given a function $u$, we use it as a Dirichlet datum for the Laplace equation in the upper half-space in $\R^{n+1}$, and then compute its Neumann condition on $\R^n\times\{0\}$. Thus, since  $u$ is a function defined on $\R^n$, we have that $Lu$ is a new function defined on $\R^n$ as well. Notice, also, that $Lu$ is nonlocal, in the sense that its value is influenced by the values of $u$ in all of~$\R^n$; and that roughly speaking, it has one derivative \emph{less} than $u$. 

We can compute ${\tilde u}$ explicitly by means of the Poisson kernel of the upper-half space for the Laplacian in $\R^{n+1}$. That is, since $u\in C^\infty_c(\R^n)$,  $\tilde u$ is given by 
\begin{equation}
\label{eq:poissonkernel0}
{\tilde u}(x, y) = [P(\cdot, y) * u](x) := \int_{\R^n} P(x-z, y) u(z)\, dz
\end{equation}
where $P(x, y)$ is the Poisson kernel for the upper half-space, 
\begin{equation}
\label{eq:poissonkernel}
P(x, y) = c_n \frac{y}{(|x|^2 +y^2)^{\frac{n+1}{2}}},
\end{equation}
with 
\begin{equation}
\label{eq:cn}
c_n := \Gamma\left(\textstyle{\frac{n+1}{2}}\right) \pi^{-\frac{n+1}{2}}.
\end{equation}
Notice that $u\in L^1_\omega(\R^n)$ is the necessary and sufficient condition for the integral in \eqref{eq:poissonkernel0} to exist. In this case ($u\in L^1_\omega(\R^n)$) there exists a unique well-defined extension (decaying at infinity). 

Then, for any $u\in C^\infty_c(\R^n)$, we can compute $Lu$:
\[
\begin{split}
Lu(x) = -\lim_{y \downarrow 0} \frac{{\tilde u}(x, y) - u(x)}{y} & = - \lim_{y\downarrow 0} \frac{1}{y}\left\{\int_{\R^n} P(x-z, y) (u(z)-u(x))\, dz \right\}\\
& = c_n \textrm{P.V.} \int_{\R^n} \frac{u(x) - u(z)}{|x-z|^{n+1}}\, dz \\
& = {\sqrtl} u(x),
\end{split}
\]
where we are using that, for each $y > 0$, $\int_{\R^n} P(z, y)\, dz = 1$, and we need to take the $\textrm{P.V.}$ in order to make sense of the integral. That is, the operator $Lu$ above coincides with the square root of the Laplacian!

Alternatively, we can also compute $L(Lu)$. Indeed, if $\tilde u$ is the harmonic extension of $u$, then $-\partial_y \tilde u$ is the harmonic extension of $Lu$ (since it coincides with $L u$ on $\{y = 0\}$ and it is harmonic). Thus,
\[
L(Lu)(x) = -\partial_y\big|_{y = 0}(-\partial_y\big|_{y = 0} {\tilde u}) = \partial_{yy}\big|_{y = 0} {\tilde u}(x, y).
\]
Now, since $\Delta_{x,y} {\tilde u} = \Delta_x {\tilde u} + \partial_{yy} {\tilde u} = 0$, 
\[
L(Lu) (x) = -\Delta_x {\tilde u}(x, 0) = -\Delta u(x). 
\]
That is, $L^2 u (x) = -\Delta u(x)$, which justifies the notation $L = {\sqrtl}$. 

In all, we have:
\begin{thm}
\label{thm:extension}\index{Harmonic extension!Square root of the Laplacian}
Let $u\in C^{1,\eps}(B_1)\cap L_\omega^1(\R^n)$ for some $\eps>0$, and let ${\tilde u}:\R^{n+1}\cap\{y > 0\} \to \R$ be the unique solution to
\begin{equation}
\label{eq:extension}
\left\{
\begin{array}{rcll}
\Delta_{x, y} {\tilde u} & = & 0& \quad\text{in}\quad \{y > 0\}\\
{\tilde u}(x, 0) & = & u(x)& \quad\text{for}\quad x\in \R^n
\end{array}
\right.
\end{equation}
with sublinear growth, given by the Poisson kernel representation \eqref{eq:poissonkernel0}-\eqref{eq:poissonkernel}. Then, 
\[
{\sqrtl}u(x) = -\partial_y\big|_{y = 0} {\tilde u}(x, y)\qquad\text{for all}\qquad x\in B_1,
\]
where ${\sqrtl}$ is given by \eqref{eq:square_root_Laplacian}.
\end{thm}
\begin{proof}
The proof is the same as the reasoning above, which also works for $u \in C^{1,\eps}(B_1)\cap L^1_\omega(\R^n)$ (thanks to Lemma~\ref{lem:laplu}). 
\end{proof}

As a corollary, we obtain the local $C^\infty$ regularity of solutions to the equation ${\sqrtl} u = 0$:
\begin{cor}
\label{cor:smoothness}
Let $u\in C^{1,\eps}(B_1)\cap L^1_\omega(\R^n)$ for some $\eps>0$ satisfy 
\[
{\sqrtl} u = 0\quad\text{in}\quad B_{1}.
\]
Then $u\in C^\infty(B_1)$ (in fact, it is real analytic).
\end{cor}
\begin{proof}
Let us use Theorem~\ref{thm:extension}, and consider ${\tilde u}$ to be the solution to \eqref{eq:extension}. Consider its even extension ${\tilde u}_e$ to $\mathcal{B}_1 := \{(x, y) \in \R^n\times\R: |(x, y)| < 1\}$, namely 
\[
{\tilde u}_e (x, y) = \left\{
\begin{array}{ll}
{\tilde u}(x, y)& \quad\text{if}\quad y \ge 0\\
{\tilde u}(x, -y)& \quad\text{if}\quad y < 0.
\end{array}
\right.
\]
Observe that ${\tilde u}_e$ is continuous in $\mathcal{B}_1$ and $C^1$ as well, since $\partial_y\big|_{y = 0} {\tilde u}(x, y) = 0$. Hence, ${\tilde u}_e$ is harmonic in the whole $\mathcal{B}_1$, and by interior regularity for harmonic functions we obtain that ${\tilde u}_e\in C^\infty(\mathcal{B}_1)$ and it is real analytic as well.
In particular, $u(x) = {\tilde u}_e(x, 0) \in C^\infty(B_1)$. 
\end{proof}

\section{Heat kernel and fundamental solution}

\label{sec:heatkernel_sqrtl}

\subsection*{Heat kernel}\index{Heat kernel!Square root of the Laplacian}
The heat kernel for the square root of the Laplacian is a function $p(t, x)$ such that, for any $\varphi\in C^\infty_c(\R^n)$,
\[
\Phi(t, x) := [p(t,\cdot)*\varphi](x)
\]
 satisfies the heat equation for the square root of the Laplacian in $\R^n$ with initial value $\varphi(x)$, that is, 
\begin{equation}
\label{eq:fract_heat_equation}
\left\{
\begin{array}{rcll}
\partial_t \Phi + {\sqrtl} \Phi & = & 0& \quad\text{in}\quad (0, +\infty) \times \R^n\\
\Phi(0, x) & = &   \varphi(x)&\quad\text{for}\quad x\in \R^n. 
\end{array}
\right.
\end{equation}
Formally, $p(t, x)$ satisfies 
\[
\left\{
\begin{array}{rcll}
\partial_t p+ {\sqrtl} p & = & 0& \quad\text{in}\quad (0, +\infty) \times \R^n\\
p(0,x)& = &  \delta_{0}&\quad\text{for}\quad x\in \R^n,
\end{array}
\right.
\]
where $\delta_{0}$ denotes a Dirac delta at $0$. 

Interestingly, the heat kernel $p(t, x)$ coincides with the Poisson kernel of the upper half-space for the Laplacian, $P(t, x)$ given by \eqref{eq:poissonkernel}:

\begin{lem} 
\label{lem:heatkernel}
Let $p(t, x)$ be given by 
\[
p(t, x) := c_n \frac{t}{(|x|^2+t^2)^{\frac{n+1}{2}}}.
\]
Then, for all $\varphi\in C^\infty_c(\R^n)$, the function $\Phi(t, x) := [p(t, \cdot)*\varphi](x)$ satisfies \eqref{eq:fract_heat_equation}. 
\end{lem}
\begin{proof}
Observe that $p(t, x) = P(t, x)$ where $P$ is the Poisson kernel in the half-space for the Laplacian, \eqref{eq:poissonkernel}. Then, on the one hand, by definition of the Poisson kernel $[\varphi*P(0, \cdot)] = \varphi$ for any $\varphi\in C^\infty_c(\R^n)$. On the other hand, observe that if $\tilde u(x, y)$ is the harmonic extension of a function $u(x)$, then $\tilde u(x, y+t)$ is the harmonic extension of $\tilde u(x, t)$ for each fixed $t > 0$. In particular, for each $t >0$ fixed, the harmonic extension of $\Phi(t, x) = [P(t, \cdot)*\varphi](x)$ is given by $\tilde\Phi(t, x, y) = [P(t+y, \cdot)*\varphi](x)$. Hence, we can compute (using Theorem~\ref{thm:extension})
\[
{\sqrtl} \Phi(t, x) = -\partial_y\big|_{y = 0}\tilde\Phi(t, x, y) = -\partial_t [P(t, \cdot)*\varphi](x) = - \partial_t \Phi(t, x),
\]
that is, $\Phi$ satisfies the fractional heat equation \eqref{eq:fract_heat_equation}. 
\end{proof}
\subsection*{Fundamental solution}
\index{Fundamental solution!Square root of the Laplacian}
Recall that, for the Laplace operator $-\Delta$, the inverse operator is given by the Riesz potential $I_2$. Indeed, a solution to 
\[
-\Delta v = f\quad\text{in}\quad \R^n
\]
(with $f$ decaying at infinity) is given by convolution against a locally integrable function,
\[
v(x) = (I_2 f)(x) =(K_2 * f)(x) =  \kappa_{2,n} \int_{\R^n} \frac{f(z)}{|x-z|^{n-2}}\, dz \qquad \text{(when $n\neq  2$)},
\]
where 
\begin{equation}
\label{eq:K2}
K_2(x) := \left\{
\begin{array}{lll}
 \kappa_{2,n} |x|^{2-n},\qquad &\quad\textrm{if $n \neq 2 $},&\textrm{with}\quad  \kappa_{2,n} =  \frac{\Gamma\left(\frac{n-2}{2}\right)}{4}\pi^{-\frac{n}{2}},\\
-\frac{1}{2\pi}\log|x|,\qquad &\quad\textrm{if $n =2$}, &
\end{array}
\right.
\end{equation}
is called the \emph{fundamental solution} for the Laplacian in $\R^n$. In other words, $I_2 f = (-\Delta)^{-1} f$. Formally, the kernel $K_2$ satisfies $-\Delta K_2 = \delta_0$, where $\delta_0$ is a Dirac  delta at the origin.

For the square root of the Laplacian, ${\sqrtl}$, we need to introduce the Riesz potential $I_1$. That is, a solution to 
\[
{\sqrtl} u = f\quad\text{in}\quad \R^n
\]
(when $f$ has sufficient decay) is given by 
\[
u(x) = (I_1 f)(x) =(K_1 * f)(x) = \kappa_{1,n} \int_{\R^n} \frac{f(z)}{|x-z|^{n-1}}\, dz \qquad \text{(when $n\ge 2$)},
\]
where 
\begin{equation}
\label{eq:K1}
K_1(x) := \left\{
\begin{array}{lll}
\kappa_{1,n} |x|^{1-n},\qquad &\quad\textrm{if $n \ge 2$},&\textrm{with}\quad  \kappa_{1,n} =  \frac{\Gamma\left(\frac{n-1}{2}\right)}{2}\pi^{-\frac{n+1}{2}},\\
-\frac{1}{\pi}\log|x|,\qquad &\quad\textrm{if $n =1$}. &
\end{array}
\right.
\end{equation}

In the distributional sense, we have
\[
{\sqrtl} K_1 = \delta_0 \quad \textrm{in}\quad \R^n,
\]
where $\delta_0$ is a Dirac delta at the origin. The function $K_1$ is called the \emph{fundamental solution} of the square root of the Laplacian in $\R^n$.

\begin{lem}[Fundamental solution]
\label{lem:fundamental_solution}
Let $f\in C^\infty_c(\R^n)$, and let $K_1$ be given by \eqref{eq:K1}. Then 
\[
{\sqrtl} (K_1*f) = f\quad\textrm{in}\quad\R^n.
\]
\end{lem}
\begin{proof}[First proof]
We use the harmonic extension from Theorem~\ref{thm:extension}. That is, we consider $\tilde u:\R^{n+1}_+\to \R$ to be the extension towards $\{y > 0\}$ of $K_1*f$, which is explicit in terms of the Poisson kernel \eqref{eq:poissonkernel0}-\eqref{eq:poissonkernel}
\[
\tilde u(x, y) = \left(P(\cdot, y) * [K_1 * f]\right)(x) = \left([P(\cdot, y) * K_1] * f\right)(x) 
\]
where we are using the associative property of the convolution. Observe that now $P(\cdot, y) * K_1$ should be the harmonic extension of $K_1$, which is the fundamental solution of the Laplacian in dimension $n+1$ (up to a multiplicative factor 2): indeed, it is the unique harmonic function vanishing at infinity  (or with sublinear growth, for $n=1$) that coincides with $K_1$ on $\{y = 0\}$ and is harmonic in $\{y > 0\}$. Thus, 
\[
\tilde u(x, y) = [\tilde K_2(\cdot, y) * f](x) \quad\text{with}\quad \tilde K_2(x, y) :=  \kappa_{1,n} \left(|x|^2 + y^2\right)^{-\frac{n-1}{2}}
\]
when $n \ge 2$, and $\tilde K_2(x, y) = -\frac{1}{2\pi} \log\left(|x|^2+y^2\right)$ if $n = 1$. 

In particular, we can compute 
\[
{\sqrtl} (K_1*f)(x) = -\partial_y\big|_{y = 0} \tilde u(x, y) = [-\partial_y \tilde K_2(\cdot, y) * f]\big|_{y = 0}(x).
\]
We now notice that 
\[
-\partial_y \tilde K_2(\cdot, y) =  P(x, y),
\]
and so $[P(\cdot, y)*f](x) = f(x)$ when $y \downarrow 0$, as we wanted to see. 
\end{proof}
\begin{proof}[Second proof]
Let us do the case $n\ge 2$. Observe that the function $K_1$ satisfies
\[
K_1(x) = \int_0^\infty p(t, x) \, dt,
\]
where $p(t, x)$ is the heat kernel (see Lemma~\ref{lem:heatkernel}). Then, for all $f\in C^\infty_c(\R^n)$ and using Lemma~\ref{lem:heatkernel}, 
\[
\begin{split}
{\sqrtl}(K_1*f) & = {\sqrtl}\left(\int_0^\infty p(t, \cdot)* f\, dt \right) = \int_0^\infty {\sqrtl}\left(p(t, \cdot) *f\right) \, dt \\
& = -\int_0^\infty \partial_t \left(p(t, \cdot) *f\right) \, dt =  p(0, \cdot) *f = f,
\end{split}
\]
as wanted.
\end{proof}

\section{Maximum principle}
\label{ssec:max_sqrtl}
Super- and subsolutions to equations of the form $\sqrtl u=0 $ satisfy a maximum principle and a comparison principle: 
\begin{lem}[Maximum Principle]
\label{lem:max_principle}\index{Maximum principle!Square root of the Laplacian}
Let $\Omega\subset\R^n$ be a bounded open set, and let $u\in C^{1,\eps}_{\rm loc}(\Omega)\cap C(\overline{\Omega})\cap L^1_\omega(\R^n)$ for some $\eps>0$. Let us assume that 
\[
\left\{
\begin{array}{rcll}
\sqrtl u & \ge & 0 & \quad\text{in}\quad \Omega,\\
u & \ge & 0 & \quad\text{in}\quad \R^n\setminus \Omega.
\end{array}
\right.
\]
Then $u\ge 0$ in $\R^n$. Moreover, either $u > 0$ in $\Omega$ or $u \equiv 0$ in $\R^n$. 
\end{lem}
\begin{proof}
Let us suppose that it is not true, and that $\inf_{\R^n} u = u(x_\circ) \le 0$ for some $x_\circ\in \Omega$. Then
\[
\sqrtl u (x_\circ) = c_n \textrm{P.V.}\int_{\R^n} \frac{u(x_\circ) - u(x_\circ+y)}{|y|^{n+1}}\, dy  < 0 
\]
if $u \not\equiv 0$ in $\R^n$, a contradiction. Thus, either $u \equiv 0$ in $\R^n$ or $u > 0$ in $\Omega$. 
\end{proof}

As a consequence we obtain the comparison principle for strong solutions. 
\begin{cor}[Comparison Principle]
\label{cor:comp_principle}\index{Comparison principle!Square root of the Laplacian}
Let $\Omega\subset\R^n$ be a bounded open set, and let $u_1, u_2\in C^{1,\eps}_{\rm loc}(\Omega)\cap C(\overline{\Omega})\cap L^1_\omega(\R^n)$ for some $\eps>0$. Let us assume that 
\[
\left\{
\begin{array}{rcll}
\sqrtl u_1 & \ge & \sqrtl u_2 & \quad\text{in}\quad \Omega,\\
u_1 & \ge & u_2 & \quad\text{in}\quad \R^n\setminus \Omega.
\end{array}
\right.
\]
Then $u_1\ge u_2$ in $\R^n$. 
\end{cor}
\begin{proof}
This is an immediate consequence of Lemma~\ref{lem:max_principle}, applied to $w := u_1 - u_2$. 
\end{proof}

We can also deduce the following: 

\begin{cor}[Uniqueness of strong solutions]
\label{cor:uniqueness_sqrtl}\index{Uniqueness!Square root of the Laplacian}
Let $\Omega\subset\R^n$ be a bounded open set, and let $u_1, u_2\in C_{\rm loc}^{1,\eps}(\Omega)\cap C(\overline{\Omega})\cap L^1_\omega(\R^n)$ for some $\eps>0$. Let us assume that 
\[
\left\{
\begin{array}{rcll}
\sqrtl u_1 & =& \sqrtl u_2 & \quad\text{in}\quad \Omega,\\
u_1 & = & u_2 & \quad\text{in}\quad \R^n\setminus \Omega.
\end{array}
\right.
\]
Then $u_1= u_2$ in $\R^n$. 
\end{cor}
\begin{proof}
It is a consequence of applying the comparison principle (Corollary~\ref{cor:comp_principle}) twice to $u_1$ and $u_2$ by exchanging their roles. 
\end{proof}

\section{Poisson kernel and the mean value property}

\subsection*{Poisson kernel}
\index{Poisson kernel!Square root of the Laplacian}
As   used in Section~\ref{sec:harm_extension}, the Poisson kernel for the upper half-space $\R^{n}_+$ for the Laplacian is the kernel that allows us to compute the harmonic extension to a boundary datum on $\R^{n} \cap \{x_n = 0\}$.  Namely, given a (sufficiently smooth) function $g\in C^\infty_c(\{x_n = 0\})$, we can explicitly solve the Dirichlet problem in the upper half-space with boundary datum given by $g$ (among bounded solutions), 
\begin{equation}
\label{eq:poisson_kernel_Laplacian}
\begin{array}{c}
\left\{
\begin{array}{rcll}
\Delta u & = & 0& \quad\text{in}\quad \{x_n > 0\}\\
u & = & g & \quad\text{on}\quad \{x_n = 0\}
\end{array}
\right.
\\[0.4cm]
\Downarrow \\[0.2cm]
u(x', x_n) = c_{n-1} \displaystyle{\int_{\{x_n = 0\}}}  \frac{x_n g(x'-z')}{(|z'|^2 + x_n^2)^{\frac{n}{2}}}\, dz',
\end{array}
\end{equation}
where we have denoted $x = (x', x_n) \in \R^{n-1}\times \R$. 

We now want  an analogous result for the square root of the Laplacian. That is, given a domain $\Omega = \{x_n > 0\}$, we want to find the  solution $u$ being $\frac12$-harmonic in $\Omega$ (i.e., ${\sqrtl} u = 0$ in $\Omega$) and with the corresponding boundary condition. Observe that now, however, due to the nonlocality of the operator, in order to compute ${\sqrtl} u$ we need to know the value of $u$ at \emph{all} points in $\R^n$. That is, we need to impose a condition in the whole $\R^n\setminus \Omega$. This is called an \emph{exterior condition} (contrary to a \emph{boundary condition} in the previous problem, where the operator is local). 

The Poisson kernel in the half-space then reads as follows:
\begin{prop}[Poisson kernel in a half-space]
\label{prop:poisson_half_space}\index{Poisson kernel!Square root of the Laplacian!Half-space}
Let $g\in L^\infty(\{x\cdot e \le 0\})$ for some $e\in \mathbb{S}^{n-1}$.
Let $u$ be defined as
\[u(x)=a_n\int_{\{z\cdot e\le 0\}}\frac{\sqrt{x\cdot e}\,g(z)}{\sqrt{|z\cdot e|}|x-z|^n}\,dz\qquad \textrm{if}\quad \{x\cdot e > 0\},\]
 where $a_n=\Gamma\left(\frac{n}{2}\right)\pi^{-\frac{n}{2}-1}$, and $u(x) = g(x)$ if $x\cdot e \le 0$. Then $u\in C^\infty(\{x\cdot e > 0\})\cap L^\infty(\R^n)$ and it solves 
\[
\left\{
\begin{array}{rcll}
{\sqrtl}u &=&0&\quad\textrm{in}\quad \{x\cdot e >0\}\\
u&=&g&\quad\textrm{in}\quad \{x\cdot e\le 0\}.
\end{array}
\right.
\]
\end{prop}

\begin{proof}
We give a short proof of the result  for $n = 1$, and we refer to \cite{BMR09} for a detailed proof of the general case $n \ge 2$. 
 
Observe that, for any $x > 0$ and by a scaling argument, 
\[
a_1\int_{\{z\le 0\}}\frac{\sqrt{x}}{\sqrt{|z|}|x-z|}\,dz = a_1\int_0^\infty\frac{dz}{\sqrt{z}(1+z)}= \frac{2}{\pi}\int_0^\infty\frac{dx}{1+x^2} = 1.
\]
Hence, $\|u\|_{L^\infty(\{x > 0\})}\le \|g\|_{L^\infty(\{x\le 0\})}$; and $u\in C^\infty(\{x  > 0\})$ by the smoothness of the integrand in $x$. It remains to prove that $\sqrtl u = 0$ in $\{x  >0\}$:

Let $\tilde u(x, y):\R^2\to \R$ be the harmonic extension of $u(-x)$ towards $y > 0$, extended evenly to $y \le 0$; and let us define
\[
w(\bar x, \bar y) := \tilde u(\bar x^2-\bar y^2, 2\bar x\bar y). 
\]

A simple computation shows that $w$ is harmonic in the half-space $\R^{2}_+=\{\bar y>0\}$ if and only if $\tilde u = \tilde u( x, y)$ is harmonic  in $\R^2\setminus \{ x \ge 0,  y = 0\}$ (alternatively, $w$ is the composition of a function $\tilde u$ with a conformal map in $\R^2$, so it is harmonic if and only if $\tilde u$ is harmonic). Thus, we need to impose $w$ to be harmonic in $\{\bar y > 0\}$. For this, we can use the Poisson kernel representation \eqref{eq:poisson_kernel_Laplacian} to derive
\[w(\bar x,\bar y)=\frac{1}{\pi}\int_\R \frac{w(z,0)\bar y}{(\bar x-z)^2+\bar y^2}\,dz.\]
In particular, for $x>0$ we recover an expression for $u$,
\[u(x)=w(0,\sqrt{x})=\frac{1}{\pi}\int_\R \frac{u(-z^2)\sqrt{x}}{z^2+x}\,dz=\frac{1}{\pi}\int_{\{z<0\}} \frac{g(z)\sqrt{x}}{\sqrt{|z|}(x-z)}\,dz.\]
Thanks to the extension theorem, Theorem~\ref{thm:extension}, such $u$ satisfies $\sqrtl u = 0$ for $x > 0$.
\end{proof}

Thanks to the previous result, and by means of an inversion in a sphere, we can also obtain the Poisson kernel for a ball. When the operator is the Laplacian, $-\Delta$, a similar argument yields that
\[
\begin{array}{c}
\left\{
\begin{array}{rcll}
\Delta u & = & 0& \quad\text{in } B_1\\
u & = & g & \quad\text{on }\partial B_1
\end{array}
\right.
\\[0.4cm]
\Downarrow
\\[0.2cm]
 u(x) = \displaystyle{\frac{c_{n-1}}{2} \int_{\partial B_1} \frac{(1-|x|^2)g(z)}{|x-z|^n}\, dz'}.
\end{array}
\]

In our case, in order to obtain the Poisson kernel in a ball for ${\sqrtl}$, we need to define the \emph{Kelvin transform of the square root of the Laplacian} for a function $u$. That is, given the inversion in $B_1$ 
\[
x\longmapsto x^* = \frac{x}{|x|^2}
\]
we define $u^*$ the Kelvin transform of $u$ as 
\[
u^*(x) := |x|^{1-n}u(x^*).
\]

The following lemma says that if $u$   satisfies $\sqrtl u = 0$ in $D$, then $u^*$ does so in $D^*$, where $D^*$ is the image of $D$ through $x\mapsto x^*$. We provide two proofs; one by direct computation and one using the extension variable instead. 
\begin{lem}[Kelvin transform for $\sqrtl$]
\label{lem:kelvin_transform}\index{Kelvin transform}
Let $u\in L^1_\omega(\R^n)$. Then $u^*\in L^1_\omega(\R^n)$, with 
\begin{equation}
\label{eq:L1KT}
\|u^*\|_{L^1_\omega(\R^n)} = \|u\|_{L^1_\omega(\R^n)}.
\end{equation}
Moreover, if $x\neq 0$ and $u\in C^{1,\eps}(B_r(x^*))$  for some $\eps>0$ and $r > 0$, then ${\sqrtl} u^*$ is well-defined at $x = (x^*)^*\in \R^n$ and 
\[{\sqrtl}u^*(x)=|x|^{-1-n}{\sqrtl}u(x^*).\]
\end{lem}
\begin{proof}
By changing variables $z\mapsto \zeta = z^*$, so that $\frac{dz}{|z|^{2n}} = d\zeta$, we get
\[
\int_{\R^n} \frac{|u^*(z)|}{1+|z|^{n+1}}\, dz = \int_{\R^n} \frac{|u^*(\zeta^*)|}{1+|\zeta^*|^{n+1}}\, |\zeta|^{-2n}d\zeta = \int_{\R^n} \frac{ |u(\zeta)|}{1+|\zeta|^{n+1}}\, d\zeta 
\]
and hence \eqref{eq:L1KT} holds. 

The fact that $\sqrtl u^*$ is well-defined at $x^*$ is due to Lemma~\ref{lem:laplu}, since $u^*$ is $C^{1,\eps}$ around $x^*$ and $u^*\in L^1_\omega(\R^n)$. 

Changing variables again, $z\mapsto \zeta = z^*$,  we find
\[
\begin{split}
{\sqrtl} u(x^*) & = c_n\textrm{P.V.} \int_{\R^n} \frac{u(x^*)- u(z)}{|x^*-z|^{n+1}}\, dz \\
& = c_n\textrm{P.V.} \int_{\R^n} \frac{u(x^*)- u(\zeta^*)}{|x^*-\zeta^*|^{n+1}}\, |\zeta|^{-2n}\,d\zeta.
\end{split}
\]
Observe now
\begin{equation}
\label{eq:obs_xstar_zstar}
|x^* - z^*|\cdot |x|\cdot |z| = |x-z|,
\end{equation}
so that
\[
\begin{split}
{\sqrtl} u(x^*) & = |x|^{n+1} c_n\textrm{P.V.} \int_{\R^n} \frac{u^*(x)\left(\frac{|x|}{|z|}\right)^{n-1} - u^*(\zeta)}{|x-\zeta|^{n+1}} \, d\zeta\\
&  \hspace{-1.25cm}= |x|^{n+1} c_n\left(\textrm{P.V.} \int_{\R^n} \frac{u^*(x) - u^*(\zeta)}{|x-\zeta|^{n+1}} \, d\zeta+|x|^{n-1} u^*(x){\sqrtl} \left(|x|^{1-n}\right)\right)\\ 
&\hspace{-1.25cm} = |x|^{n+1}{\sqrtl}u^*(x),
\end{split}
\]
where in the last inequality we are using that $|x|^{1-n}$ is  (a multiple of) the fundamental solution (see Lemma~\ref{lem:fundamental_solution}), and so ${\sqrtl} \left(|x|^{1-n}\right) = 0$ since $x\neq 0$. Taking $y = x^*$ we get the desired result. 
\end{proof}
\begin{proof}[Second proof]
We can give a proof of the second part by means of the extension method, Theorem~\ref{thm:extension}. Let us see that if $\tilde u(x, y)$ and ${\widetilde {u^*}}(x, y)$ denote the harmonic extensions towards $\{y > 0\}$ of $u$ and $u^*$ respectively, then $\tilde u^\circ = \widetilde {u^*}$, where 
\[
w^\circ(X) := |X|^{1-n} w(X^*),\qquad\text{for}\quad X = (x, y)\in \R^n\times \R. 
\]
Here, $w^\circ$ is the Kelvin transform for the (classical) Laplacian: observe that the only difference from before is in the power $1-n$, since now we are in dimension $n+1$. It is now a direct computation to check that
\[
\Delta w^\circ(X) = |X|^{-n-2} \Delta w (X^*),
\]
and so $w^\circ$ is harmonic in $D$ if and only if $w$ is harmonic in $D^*$. In particular, since half-spaces and hyperplanes that go through the origin are invariant under the inversion $x\mapsto x^*$, we obtain that $\tilde u^\circ = \widetilde {u^*}$. Now, we can compute (using Theorem~\ref{thm:extension})
\[
{\sqrtl} u^*(x) = \partial_{y}\big|_{y = 0} \widetilde {u^*}(x, y) = \partial_{y}\big|_{y = 0} \tilde {u}^\circ (x, y) = \partial_{y}\big|_{y = 0} |X|^{1-n} \tilde {u}(X^*). 
\]
Observe now that $\partial_{y}\big|_{y = 0} |X|^{1-n} = 0$ (if $X\neq 0$) and 
\[
\partial_{y}\big|_{y = 0} \tilde {u}(X^*) = \nabla_X \tilde u(X^*)\cdot \partial_y (X/|X|^2) \big|_{y = 0} =  (\partial_y \tilde u)(X^*)\big|_{y =0} |x|^{-2}.
\]
Hence, we obtain
\[
{\sqrtl} u^*(x) = [\partial_y\big|_{y = 0}\tilde u](x^*)|x|^{-1-n} = |x|^{-1-n}{\sqrtl} u(x^*),
\]
as we wanted to see.
\end{proof}

The inversion in $B_1$, $x\mapsto x^*$, maps (conformally) the interior of $B_1$ into its exterior, and vice versa. Observe that this transformation maps any ball whose boundary goes through the origin into a half-space, and in particular $(B_{1}(\boldsymbol{e}_n))^* = \{x_{n}> 1/2\}$. Thanks to this observation, we can compute the Poisson kernel in a ball.\footnote{We refer to \cite{Buc16} for a direct proof of Proposition~\ref{prop:poisson_kernel_ball}, not based on Kelvin transforms nor the Poisson kernel of the half-space. 
}
\begin{prop}[Poisson kernel in a ball]
\label{prop:poisson_kernel_ball}\index{Poisson kernel!Square root of the Laplacian!Ball}
Let $g\in L^1_\omega(\R^n)$ such that it is continuous on $\partial B_1$. Let $u$ be defined as
\[u(x)=a_n \int_{\R^n\setminus B_1}\frac{\sqrt{1-|x|^2}\,g(z)}{\sqrt{|z|^2-1}|x-z|^n}\,dz\qquad\textrm{if}\quad x\in B_1,\]
 where $a_n=\Gamma\left(\frac{n}{2}\right)\pi^{-\frac{n}{2}-1}$, and $u(x) = g(x)$ in  $\R^n \setminus B_1$. Then $u\in C^\infty(B_1)\cap C(\overline{B_1})$ and it solves 
\begin{equation}
\label{eq:whatwewanted}
\left\{
\begin{array}{rcll}
{\sqrtl} u &=&0&\quad\textrm{in}\quad B_1\\
u&=&g&\quad\textrm{in}\quad \R^n\setminus B_1.
\end{array}
\right.
\end{equation}
\end{prop}

\begin{proof}
Given $g_0\in C^\infty_c(\R^n\setminus \overline{B_1(\be_n)})$. we look for  a function $u_0$ that satisfies
\begin{equation}
\label{eq:uenough}
\left\{
\begin{array}{rcll}
{\sqrtl} u_0 &=&0&\quad\textrm{in}\quad B_{1}(\be_n)\\
u_0&=&g_0&\quad\textrm{in}\quad \R^n\setminus B_{1}(\be_n).
\end{array}
\right.
\end{equation}
By Lemma~\ref{lem:kelvin_transform}, $u_0^*(x) := |x|^{1-n}u_0(x^*)$ would  satisfy
\[
{\sqrtl} u_0^* (x) = |x|^{-1-n} {\sqrtl} u_0(x^*) = 0
\]
for all $x$ such that $x^*\in B_{1}(\be_n)$. Since $x\mapsto x^*$ is an involution, and we can check that $(B_{1}(\be_n))^* = \{x_n > 1/2\}$, we have that $u_0$ satisfying the above equation \eqref{eq:uenough} is equivalent to $u_0^*$ satisfying
\[
\left\{
\begin{array}{rcll}
{\sqrtl} u_0^* &=&0&\quad\textrm{in}\quad \{x_n > 1/2\}\\
u_0^*&=&g_0^*&\quad\textrm{in}\quad \{x_n \le 1/2\}.
\end{array}
\right.
\]

By Proposition~\ref{prop:poisson_half_space} (and Lemma~\ref{lem:invariances}-\eqref{it:prop_trans},  
we can take:
\[
\begin{split}
u_0^* (x) & = a_n \int_{\{z\cdot \be_n \le 1/2\}} \frac{\sqrt{x\cdot \be_n-1/2}\, g_0^*(z)}{\sqrt{1/2-z\cdot \be_n }|x-z|^n}\, dz \\
& = a_n \int_{\R^n\setminus B_{1}(\be_n)} \frac{\sqrt{x\cdot \be_n-1/2}\, g_0^*(\zeta^*)}{\sqrt{1/2-\zeta^*\cdot \be_n }|x-\zeta^*|^n}\, |\zeta|^{-2n} d\zeta,
\end{split}
\]
where we have used the change of variables $z\mapsto \zeta = z^*$. Hence, also recalling \eqref{eq:obs_xstar_zstar}, we find that 
\[
u_0 (x) = |x|^{1-n} u_0^*(x^*) = a_n \int_{\R^n\setminus B_{1}(\be_n)} \frac{\sqrt{x\cdot \be_n-|x|^2/2}\,g_0(\zeta)}{\sqrt{|\zeta|^2/2-\zeta\cdot \be_n }|x-\zeta|^n}\,   d\zeta
\]
satisfies \eqref{eq:uenough}.

Observe now, that by translation invariance of the square root of the Laplacian (Lemma~\ref{lem:invariances}-\eqref{it:prop_trans}), if we define $u(x) = u_0\left({x+\be_n}\right)$ with $g_0(y) = g\left(y-\be_n\right)$, changing variables in the expression above we deduce that
\[
u (x) = a_n \int_{\R^n\setminus B_{1}} \frac{ \sqrt{1-|x|^2} g(z)}{\sqrt{|z|^2-1}|x-z|^n}\,    dz
\]
satisfies \eqref{eq:whatwewanted} whenever $g\in C^\infty_c(\R^n\setminus \overline{B_1})$. By approximation, it also works for any $g\in L^1_\omega(\R^n)$ that is continuous on $\partial B_1$. Here, the continuity on $\partial B_1$ ensures that the integral defining $u(x)$ is uniformly bounded in $x$. 

Finally, the smoothness of $u$ inside $B_1$ is by definition (since the kernel is smooth in $B_1$), and the continuity of $u$ up to $\partial B_1$ holds because 
\[
a_n\int_{\R^n\setminus B_1} \frac{\sqrt{1-|x|^2}}{\sqrt{|z|^2-1}|x-z|^2}\, dz = 1
\]
independently of $x\in B_1$, so that taking the limit $B_1\ni x \to x_\circ\in \partial B_1$, 
\[
\lim_{B_1\ni x\to x_\circ} |u(x) - g(x_\circ)|  \le  \lim_{B_1\ni x\to x_\circ}\int_{\R^n\setminus B_1} a_n \frac{\sqrt{1-|x|^2}\,\big|g(z) - g(x_\circ)\big|}{\sqrt{|z|^2-1}|x-z|^2}\, dz =0,
\]
by the continuity of $g$ and the fact that the integrand concentrates all the mass around $x_\circ$. 
\end{proof}

As a corollary, we obtain the unique representation of strong solutions in the unit ball:
\begin{cor}
\label{cor:uniqueness_strong_solutions}
Let $u\in C_{\rm loc}^{1,\eps}(B_1) \cap C(\overline{B_1})\cap L^1_\omega(\R^n)$ for some $\eps>0$ such that 
\[
{\sqrtl}  u = 0\quad\text{in}\quad B_1. 
\]
Then, $u$ satisfies
\[u(x)=a_n \int_{\R^n\setminus B_1}\frac{\sqrt{r^2-|x|^2}\,u(z)}{\sqrt{|z|^2-r^2}|x-z|^n}\,dz\]
for all $x\in B_1$.
\end{cor}
\begin{proof}
If we denote
\[u_0(x)=a_n \int_{\R^n\setminus B_r}\frac{\sqrt{r^2-|x|^2}\,u(z)}{\sqrt{|z|^2-r^2}|x-z|^n}\,dz\]
for $x\in B_1$, and $u_0(x) = u(x)$ if $x\in \R^n\setminus B_1$, then, by Proposition~\ref{prop:poisson_kernel_ball}, $u_0\in C^\infty(B_1)\cap C(\overline{B_1})\cap L^1_\omega(\R^n)$ and ${\sqrtl}  u_0 = 0$ in $B_1$. The result now follows by the uniqueness of strong solutions, Corollary~\ref{cor:uniqueness_sqrtl}.
\end{proof}

\begin{rem}[Poisson kernel in $B_r$]
\label{rem:rescaledpoisson}
By rescaling the Poisson kernel in a ball $B_1$ from Proposition~\ref{prop:poisson_kernel_ball} (see Lemma~\ref{lem:invariances}) we obtain the Poisson kernel in any ball $B_r$. That is, the solution of
\[
\left\{
\begin{array}{rcll}
{\sqrtl}u &=&0&\quad\textrm{in}\quad B_r\\
u&=&g&\quad\textrm{in}\quad \R^n\setminus B_r
\end{array}
\right.
\]
is given by 
\[u(x)=a_n \int_{\R^n\setminus B_r}\frac{\sqrt{r^2-|x|^2}\,g(z)}{\sqrt{|z|^2-r^2}|x-z|^n}\,dz,\]
for any $r > 0$. 
\end{rem}

\subsection*{Mean value property}
\index{Mean value property!Square root of the Laplacian}
A direct computation for the Laplace operator yields the mean value property,
\[
\Delta u = 0\quad\textrm{in}\quad \Omega \quad \Longrightarrow \quad u (x) = \ave_{\partial B_r(x)} u := \frac{1}{|\partial B_r|} \int_{\partial B_r(x)} u
\]
for any $x\in \Omega$, and $B_r(x) \ssubset \Omega$. Moreover, integrating in $r$ we also get 
\[
u(x) = \ave_{B_r(x)} u := \frac{1}{|B_r|} \int_{B_r(x)} u\qquad\textrm{whenever}\quad B_r(x) \subset \Omega. 
\]

We obtain the analogy of the mean value property for functions $u$ that satisfy $\sqrtl u = 0$ in $\Omega$. Observe that, unsurprisingly, we cannot restrict ourselves to the values of $u$ on $\partial B_r$, and we need to take exterior domains instead.

\begin{prop}[Mean Value Property]
\label{prop:meanvalue}
Let $\Omega\subset \R^n$ be an open domain, and let $u\in C_{\rm loc}^{1,\eps}(\Omega)\cap L^1_\omega(\R^n)$ for some $\eps>0$. If ${\sqrtl} u = 0$ in $\Omega$, then for every $B_r\ssubset \Omega$,
\[
u(0) = a_n \int_{\R^n\setminus B_r} \frac{ru(z)}{\sqrt{|z|^2-r^2}|z|^n}\, dz.
\]
\end{prop}
\begin{proof}
This is a consequence of Remark~\ref{rem:rescaledpoisson} and Corollary~\ref{cor:uniqueness_strong_solutions}. 
\end{proof}

As a corollary, we obtain:
\begin{cor}[Mean Value Property II]
\label{cor:meanvalue}
There exists a nonincreasing function $\omega(t):[0,+\infty)\to [0,+\infty)$  such that if $u\in C_{\rm loc}^{1,\eps}(B_1)\cap L^1_\omega(\R^n)$ for some $\eps>0$ satisfies 
\[
{\sqrtl} u = 0\quad\textrm{in}\quad B_1 
\]
then 
\[
u(0) = \int_{\R^n} u(z) \omega(|z|)\, dz.
\]
Moreover,
\begin{equation}
\label{eq:omega_estimates}
\frac{C^{-1}}{1+|t|^n} \le \omega(t) \le \frac{C}{1+|t|^{n+1}}
\end{equation}
for some $C$ depending only on $n$. 
\end{cor}
\begin{proof}
By Proposition~\ref{prop:meanvalue}, 
\[
u(0) = n\int_0^1 r^{n-1} u(0) \, dr = a_n \int_0^1\int_{|z|\ge r} \frac{r^n u(z)}{\sqrt{|z|^2-r^2}|z|^n}\, dz\,dr = a_nI_1 + a_nI_2,
\]
where 
\[
I_1 := \int_0^1 \int_{|z|\ge 1} \frac{r^n u(z)}{\sqrt{|z|^2-r^2}|z|^n}\, dz\,dr,\quad I_2 = \int_{|z|\le 1} \int_0^{|z|} \frac{r^n u(z)}{\sqrt{|z|^2-r^2}|z|^n}\, dr\,dz.
\]
Now, we compute separately
\[
I_1 = \int_{|z|\ge 1} \frac{u(z) }{|z|^n}\int_0^1 \frac{r^n}{\sqrt{|z|^2-r^2}}\, dr\, dz = \int_{|z| > 1} u(z) \int_0^{\frac{1}{|z|}} \frac{\rho^n\, d\rho}{\sqrt{1-\rho^2}}\, dz
\]
and 
\[
I_2 = \int_{|z|\le 1}\frac{u(z) }{|z|^n}\int_0^{|z|} \frac{r^n}{\sqrt{|z|^2-r^2}}\, dr\, dz = \int_0^1 \frac{\rho^n\, d\rho}{\sqrt{1-\rho^2}} \int_{|z|\le 1} u(z)\, dz.
\]
Hence, we obtain 
\[
\omega(t) := \left\{
\begin{array}{ll}
\displaystyle a_n \int_0^1 \frac{\rho^n\, d\rho}{\sqrt{1-\rho^2}} \qquad & \textrm{if}\quad 0\le t\le 1\\
\displaystyle a_n \int_0^{\frac{1}{t}} \frac{\rho^n\, d\rho}{\sqrt{1-\rho^2}} \qquad & \textrm{if}\quad t\ge 1,
\end{array}
\right.
\]
from which \eqref{eq:omega_estimates} follows. 
\end{proof}
\begin{rem}\label{rem:meanvalue_super}
The analogue for super- and subsolutions of the mean value property also holds. Namely, as a direct consequence of the comparison principle, Corollary~\ref{cor:comp_principle}, 
under the hypotheses of Proposition~\ref{prop:meanvalue}, if we had instead that
\[
{\sqrtl} u \ge 0\quad\textrm{in}\quad \Omega,
\]
then 
\[
u(0) \ge a_n \int_{\R^n\setminus B_r} \frac{ru(z)}{\sqrt{|z|^2-r^2}|z|^n}\, dz,
\]
for every $B_r\subset \Omega$.

Hence, repeating the proof of Corollary~\ref{cor:meanvalue} in this case we obtain that if
\[
{\sqrtl} u \ge 0\quad\textrm{in}\quad B_1 
\]
then 
\[
u(0) \ge \int_{\R^n} u(z) \omega(|z|)\, dz,
\]
where $\omega$ satisfies \eqref{eq:omega_estimates}. 
\end{rem}

\section{The Harnack inequality}
\index{Harnack's inequality!Square root of the Laplacian}

The Harnack inequality is an essential tool in the study of solutions to elliptic PDE. It has many important consequences and generalizations, especially when studying regularity and convergence properties of solutions. 

For the Laplacian, it says that the supremum and the infimum of any nonnegative harmonic function are comparable. That is, 
\begin{equation}
\label{eq:harnack_laplacian}
\left\{
\begin{array}{rcll}
u &\ge& 0& \quad\textrm{in}\quad{B_1}\\
\Delta u& =& 0&\quad\textrm{in}\quad{B_1}
\end{array}
\right.
\quad
\Longrightarrow
\quad 
\sup_{B_{1/2}} u \le C \inf_{B_{1/2}} u
\end{equation}
for some $C$ depending only on $n$. 

For the square root of the Laplacian we have an analogous Harnack inequality, but in this case, the nonnegativity condition needs to be imposed in the whole space.

\begin{prop}
\label{prop:Harnack_sqrt}
Let $u\in C_{\rm loc}^{1,\eps}(B_1)\cap L^1_\omega(\R^n)$ for some $\eps>0$ be such that 
\[
\left\{
\begin{array}{rcll}
u &\ge& 0& \quad\textrm{in}\quad \R^n\\
{\sqrtl} u& =& 0&\quad\textrm{in}\quad{B_1}.
\end{array}
\right.
\]
Then, 
\begin{equation}
\label{eq:harnack}
\sup_{B_{1/2}} u \le C \inf_{B_{1/2}} u
\end{equation}
for some constant $C$ depending only on $n$. 
\end{prop}
\begin{proof}
Let $\tilde u(x, y)$ be the harmonic extension of $u$ given by Theorem~\ref{thm:extension}, extended evenly to $\R^{n+1}$. Now, since $u \ge 0$, we have $\tilde u \ge 0$ as well (it is explicit in terms of the corresponding Poisson kernel, \eqref{eq:poissonkernel0}-\eqref{eq:poissonkernel}). Moreover, $\partial_y \tilde u(x, 0) = 0$ for all $x\in B_1$, and so $\tilde u$ is harmonic across $\{y = 0\}$ in $B_1$:
\[
\Delta \tilde u = 0\quad\textrm{in}\quad \mathcal{B}_1,
\]
where $\mathcal{B}_1\subset \R^{n+1}$ denotes the unit ball in $\R^{n+1}$.
Thus, we can apply the Harnack inequality for harmonic functions, \eqref{eq:harnack_laplacian}, to deduce 
\[
\sup_{B_{1/2}} u \le \sup_{\mathcal{B}_{1/2}} \tilde u \le C \inf_{\mathcal{B}_{1/2}} \tilde u \le C \inf_{B_{1/2}} u,
\]
as we wanted to see.
\end{proof}

We can alternatively show the Harnack inequality by proving the two half-Harnack inequalities for super- and subsolutions to the equation.

\begin{prop}[Half-Harnack for supersolutions]
\label{prop:half_Harnack_super}\index{Half-Harnack!Square root of the Laplacian!Supersolutions}
Let $u\in C_{\rm loc}^{1,\eps}(B_1)\cap L^1_\omega(\R^n)$  for some $\eps>0$ be such that 
\[
\left\{
\begin{array}{rcll}
u &\ge& 0& \quad\textrm{in}\quad \R^n\\
{\sqrtl} u& \ge & 0&\quad\textrm{in}\quad{B_1}.
\end{array}
\right.
\]
Then, 
\[
\inf_{B_{1/2}} u \ge c \int_{\R^n}\frac{u(z)}{1+|z|^{n+1}}\, dz,
\]
for some constant $c$ depending only on $n$. 
\end{prop}
\begin{proof} It is a direct consequence of Remark~\ref{rem:meanvalue_super} together with \eqref{eq:omega_estimates} and the fact that $u\ge 0$. 
\end{proof}

The second ``half-Harnack'' is an $L^\infty$ bound for subsolutions. 

\begin{prop}[Half-Harnack for subsolutions]
\label{prop:half_Harnack_sub}\index{Half-Harnack!Square root of the Laplacian!Subsolutions}
Let $u\in C_{\rm loc}^{1,\eps}(B_1)\cap L^1_\omega(\R^n)$  for some $\eps>0$ be such that $\sqrtl u \le 0$ in $B_1$. Then, 
\[
\sup_{B_{1/2}} u \le C \int_{\R^n}\frac{|u(z)|}{1+|z|^{n+1}}\, dz,
\]
for some constant $C$ depending only on $n$. 
\end{prop}
\begin{proof} It is a direct consequence of Remark~\ref{rem:meanvalue_super} applied to $-u$ together with \eqref{eq:omega_estimates}. 
\end{proof}

By using the previous two results, we obtain an alternative proof of the Harnack inequality: 

\begin{proof}[Second proof of Proposition~\ref{prop:Harnack_sqrt}]
The result follows by concatenating Proposition~\ref{prop:half_Harnack_super} and Proposition~\ref{prop:half_Harnack_sub}.
\end{proof}
\begin{rem}
From this last proof we actually obtain stronger information. Namely, both quantities $\sup_{B_{1/2}} u $ and $\inf_{B_{1/2}} u$ are actually comparable to $\|u\|_{L^1_\omega(\R^n)}$, that is, 
\[
\frac{1}{C_2} \int_{\R^n}\frac{u(z)}{1+|z|^{n+1}}\, dz \le \sup_{B_{1/2}} u \le C \inf_{B_{1/2}} u\le C_2 \int_{\R^n}\frac{u(z)}{1+|z|^{n+1}}\, dz 
\]
for some $ C_2$ depending only on $n$. 
\end{rem}

\begin{rem} In Proposition~\ref{prop:Harnack_sqrt}, the requirement $u \ge 0$ in $\R^n$ is necessary in order to achieve the conclusion, contrary to the requirement in \eqref{eq:harnack_laplacian}, where only nonnegativity in $B_1$ was needed. Let us show it through a counterexample:

Let $v$ be the (unique) bounded solution to 
\[
\left\{
\begin{array}{rcll}
{\sqrtl} v &=&0&\quad\textrm{in}\quad B_1\\
v&=&\chi_{B_3\setminus B_2}(x)&\quad\textrm{in}\quad \R^n\setminus B_1,
\end{array}
\right.
\]
where $\chi_A(x)$ is the characteristic function of the set $A\subset\R^n$. Then, by Proposition~\ref{prop:poisson_kernel_ball},
\[
v(x)=a_n \sqrt{1-|x|^2} \int_{B_3\setminus B_2} \frac{dz}{\sqrt{|z|^2-1}|x-z|^n}.
\]
In particular, $v$ is continuous in $B_{3/2}$, and $v \ge 0$ everywhere. Let $0<M := \max_{B_1} v(x) = v(x_\circ)$, which is attained for some $x_\circ \in B_1$. Moreover, since $\sqrtl(1-v) = 0$ in $B_1$ and $1-v \ge 0$ in $\R^n\setminus B_1$, by the strict maximum principle (Lemma~\ref{lem:max_principle}) $1-v > 0$ in $B_1$ and by continuity ($1-v = 1$ on $\partial B_1$) we have that $M < 1$. 

 The function $u := 1-\frac{v}{M}$ satisfies that ${\sqrtl} u = 0$ in $B_1$, $u \ge 0$ in $B_1$ (and $u \not\equiv 0$ in $B_1$), but nonetheless $\inf_{B_1} u = u(x_\circ) = 0$; and hence, it does not satisfy \eqref{eq:harnack}. In this case, $u \ge 0$ in $B_3^c\cup B_2$,  but $u = 1-M^{-1} < 0$ in $B_3\setminus B_2$. 

This observation is a actually a consequence of a more general fact: any smooth function $f$ in $B_1$ can be approximated by functions $\tilde f$ satisfying $\sqrtl \tilde f = 0$ in $B_1$.  Namely, for any $\eps > 0$ there exists $f_\eps\in C^{\infty}_c(\R^n)\cap L^1_\omega(\R^n)$ satisfying $\sqrtl f_\eps = 0$ in $B_1$ such that $\|f - f_\eps\|_{L^\infty(B_1)}\le \eps$; see \cite{DSV, BV}. In particular, any function $u\ge 0$ in $B_1$ can be approximated by functions $u_\eps$ with $\sqrtl u_\eps  = 0$ in $B_1$ (that, in general, might be negative somewhere in $\R^n\setminus B_1$), and thus a ``local'' Harnack inequality cannot hold. 
\end{rem}

\section{Interior regularity}
\label{sec:interior_regularity}

Notice that, since $\sqrtl$ is an elliptic operator of order one (see Lemma~\ref{lem:laplu}),   we expect solutions to $\sqrtl u = f$ to ``gain a derivative''\footnote{Alternatively, a function that satisfies $\sqrtl u = f$ in $B_1$, as seen in the extension, corresponds to a harmonic function $\tilde u$ in $\R^{n+1}_+$ such that $\partial_{x_{n+1}} \tilde u = f$ on $\{x_{n+1} = 0\}$ . Hence, we expect $\tilde u$ (and, therefore, $u$) to have one more ``derivative'' of regularity than $f$. 
} with respect to $f$. The goal of this section is to prove the following a priori interior regularity estimate:

\begin{thm}[Interior regularity]
\label{thm:interior_regularity_sqrtl}\index{Interior regularity!Square root of the Laplacian}
Let $u\in C^{1,\eps}(B_1)\cap L^1_\omega(\R^n)$ for some $\eps>0$ satisfy
\[
\sqrtl u = f\quad\text{in}\quad B_1,
\]
for some $f\in C^{k,\alpha}(B_1)$ with $k\in \N_0$ and $\alpha \in (0, 1)$. 
Then, $u\in C_{\rm loc}^{k+1, \alpha}(B_1)$ and 
\[
\|u\|_{C^{k+1,\alpha}(B_{1/2})}\le C \left(\|f\|_{C^{k,\alpha}(B_1)} + \|u\|_{L^1_\omega(\R^n)}\right),
\]
for some $C$ depending only on $n$, $k$, and $\alpha$.
\end{thm}

Let us first prove two intermediate lemmas, which will then yield Theorem~\ref{thm:interior_regularity_sqrtl}. 

The first lemma refers to the interior regularity of functions $u$ satisfying $\sqrtl u = 0$ in $B_1$, and follows by means of the Poisson kernel representation and the mean value formula:

\begin{lem}
\label{lem:interior_reg_zero_sqrtl}
Let $u\in C^{1,\eps}(B_1)\cap L^1_\omega(\R^n)$  for some $\eps>0$  satisfy 
\[
\sqrtl u = 0\quad\text{in}\quad B_1. 
\]
Then, $u\in C^\infty(B_1)$ and 
\[
\|u\|_{C^\nu(B_{1/2})}\le C_\nu \|u\|_{L^1_\omega(\R^n)}\quad\text{for all}\quad \nu \in \N,
\]
for some constant $C_\nu$ depending only on $n$ and $\nu$. 
\end{lem}
\begin{proof}
By the uniqueness of strong solutions and the Poisson kernel representation in $B_{3/4}$ (Corollary~\ref{cor:uniqueness_strong_solutions} and Remark~\ref{rem:rescaledpoisson}) we have that 
\[
u(x) = a_n \int_{\R^n\setminus B_{3/4}}\frac{\sqrt{9/16-|x|^2}\,u(z)}{\sqrt{|z|^2-9/16}|x-z|^n}\,dz\,\quad\text{for all}\quad x\in B_{3/4}. 
\]
In particular, we can differentiate with respect to $x$ under the integral sign and use that 
\[
\left|D^\nu_x  \left(\frac{\sqrt{9/16-|x|^2}}{|x-z|^n}\right)\right|\le \frac{C_\nu}{1+|z|^n}
\]
for any $x\in B_{1/2}$, $z\in \R^n\setminus B_{3/4}$, to get 
\[
|D^\nu u(x)|\le C_\nu \int_{\R^n\setminus B_{5/6}} \frac{|u(z)|}{1+|z|^{n+1}}\, dz + C_\nu \int_{B_{5/6}\setminus B_{3/4}} \frac{|u(z)|}{\sqrt{|z|-3/4}}\, dz.
\]
We finish by observing that by Proposition~\ref{prop:half_Harnack_sub} applied around any point $z\in B_{5/6}$ we get
\[
\|u\|_{L^\infty(B_{5/6})}\le C \|u\|_{L^1_\omega(\R^n)}.
\]
Hence, 
\[
|D^\nu u(x)|\le C_\nu \|u\|_{L^1_\omega(\R^n)}
\]
for all $x\in B_{1/2}$, as we wanted to see. 
\end{proof}

The second intermediate lemma refers to the regularity of global solutions. Observe that, under the hypotheses below, $\sqrtl u = f$ in $\R^n$. 
\begin{lem}
\label{lem:global_interior_reg}
Let $k\in \N_0$, $\alpha \in (0, 1)$, and let $f\in C^{k,\alpha}_c(B_2)$. Let 
\[
u(x) := (K_1 * f)(x) = \int_{\R^n} f(y)K_1(x-y)\, dy,
\]
where $K_1$ is the fundamental solution for $\sqrtl$, \eqref{eq:K1}. Then $u\in C_{\rm loc}^{k+1,\alpha}(B_2)\cap L^1_\omega(\R^n)$ and 
\[
\|u\|_{C^{k+1,\alpha}(B_1)} + \|u\|_{L^1_\omega(\R^n)}\le C \|f\|_{C^{k,\alpha}(B_2)},
\]
for some $C$ depending only on $n$, $k$, and $\alpha$. 
\end{lem}
\begin{proof}
We start observing that, since $f$ is compactly supported, 
\begin{equation}
\label{eq:firstbound}
\|u\|_{L^\infty(B_1)}\le C\|f\|_{L^\infty(B_2)},
\end{equation}
for some $C$ (due to the local integrability of $K_1$). More generally, for $n \ge 2$ and since $K_1$ is radially decreasing, we have
\[
\|u\|_{L^\infty(\R^n)}\le C\|f\|_{L^\infty(B_2)}.
\]
If $n = 1$, instead, we have 
\[
\|u\|_{L^\infty(B_R)}\le C\log (R+1) \|f\|_{L^\infty(B_2)}\quad\text{for}\quad R \ge 1,
\]
for some $C$ depending only on $n$. In both cases we obtain 
\begin{equation}
\label{eq:firstbound0}
\|u\|_{L^1_\omega(\R^n)}\le C \|f\|_{L^\infty(B_2)}
\end{equation}
for some $C$ depending only on $n$. 

Let us assume first that $k = 0$.  We take the gradient under the integral sign and use that 
\[
\nabla K_1(y)  = - C \frac{y}{|y|^{n+1}},
\]
to obtain,  for  $0 \le r \le 1$ fixed,
\[
 \nabla u(x) = \int_{B_r}  \frac{f(x+y)}{|y|^{n+1}}y \, dy + I^r(x)= \int_{B_r} \left(f(x+y) -f(x)\right)\frac{y\, dy}{|y|^{n+1}} + I^r(x),
\]
where we have used that $\text{P.V.} \int_{B_r} \frac{y}{|y|^{n+1}}\, dy = 0$ and we have denoted  $I^r(x) = (I^r_1(x),\dots, I^r_n(x))$ with
\[
I^r_i(x) := \int_{\R^n\setminus B_r}  f(x+y) \frac{y_i}{|y|^{n+1}}\, dy = \int_{\R^n\setminus B_r(x)}  f(z) \frac{z_i-x_i}{|z-x|^{n+1}}\, dz 
\]
for $i \in \{1,\dots,n\}$.
By the mean value theorem and using that $f\in C^{0,\alpha}_c(B_2)$ we get, for any $\bar x\in B_r$,
\begin{equation}
\label{eq:plugbacknablau}
\begin{split}
|\nabla u(\bar x) - \nabla u(0)|& \le C\|f\|_{C^\alpha(B_2)} \int_{B_r} |y|^{-n+\alpha}\, dy + r\|D_x I^r\|_{L^\infty(B_r)}\\
& \le C\|f\|_{C^\alpha(B_2)} r^\alpha +r\|D_x I^r\|_{L^\infty(B_r)}.
\end{split}
\end{equation}
 Moreover,
\[
\partial_{x_i} I^r_j(\bar x) = -\int_{\R^n\setminus B_r} f(\bar  x +z) \partial_{z_i} \left(\frac{z_j}{|z|^{n+1}}\right) \, dz- \int_{\partial B_r} f(\bar  x+z) \frac{z_iz_j}{r^{n+2}}\, dz.
\]
We can now use  that 
\[
\frac{1}{r^{n+2}}\int_{\partial B_r} z_i z_j\, dz = -\int_{\R^n\setminus B_r} \partial_{z_i}\left(\frac{z_j}{|z|^{n+1}}\right)
\]
(which is equal to zero if $i \neq j$) to rewrite it as
\[
\begin{split}
-\partial_{x_i} I^r_j(\bar  x) & =  \int_{\R^n\setminus B_r} \left( f(\bar x +z)-f(\bar  x)\right) \partial_{z_i} \left(\frac{z_j}{|z|^{n+1}}\right) \, dz\\
& \quad + \int_{\partial B_r} \left( f( \bar x+z)-f(\bar  x)\right) \frac{z_iz_j}{r^{n+2}}\, dz,
\end{split}
\]
and bound it, using $f\in C^{0,\alpha}(B_2)$, by 
\[
\begin{split}
\|D_x I^r\|_{L^\infty(B_r)}&  \le C\|f\|_{C^\alpha(B_2)}\left(\int_{\R^n\setminus B_r} |z|^{-n-1+\alpha}\, dz + \frac{1}{r^{n+2}}\int_{\partial B_r} |z|^{2+\alpha}\, dz\right) \\
& \le Cr^{\alpha - 1}\|f\|_{C^\alpha(B_2)}.
\end{split}
\]
Hence, from \eqref{eq:plugbacknablau} we get 
\[
|\nabla u(x) - \nabla u(0)| \le C\|f\|_{C^\alpha(B_2)} r^\alpha,\quad\text{for all}\quad x\in B_r.
\]
Since we can repeat the argument in any ball $B_r(x_\circ)$ with $x_\circ \in B_1$,  we obtain
\[
[\nabla u]_{C^{\alpha}(B_1)}\le C \|f\|_{C^\alpha(B_2)}. 
\]
Finally, when $k\ge 1$, we notice that 
\[
D^{k} u(x) = (K_1 * D^k f)(x)
\]
and $D^{k} f\in C^\alpha_c(B_2)$, so that by the previous estimate we have 
\[
[D^{k+1} u]_{C^\alpha(B_1)}\le C  \|D^k f\|_{C^\alpha(B_2)}.
\]
Together with \eqref{eq:firstbound} and \eqref{eq:firstbound0} we get the desired result. 
\end{proof}

Combining the previous two lemmas, we obtain the proof of the interior regularity in Theorem~\ref{thm:interior_regularity_sqrtl}:

\begin{proof}[Proof of Theorem~\ref{thm:interior_regularity_sqrtl}]
Let $\tilde f$ be a globally defined $C^{k,\alpha}$ extension of $f$; namely, $\tilde f\in C^{k,\alpha}_c(B_2)$ such that $\tilde f = f$ in $B_1$ and 
\[
\|\tilde f\|_{C^{k,\alpha}(\R^n)} \le C \|f\|_{C^{k,\alpha}(B_1)}. 
\]

Let us define $w_1 := K_1 * \tilde f$ so that, thanks to the fundamental solution representation from Lemma~\ref{lem:fundamental_solution} (and an approximation argument), we have that $\sqrtl w_1 = \tilde f$ in $\R^n$. Then, by Lemma~\ref{lem:global_interior_reg}) we have
\begin{equation}
\label{eq:w1}
\|w_1\|_{C^{k+1,\alpha}(B_{1})}\le C \|\tilde f\|_{C^{k,\alpha}(\R^n)} \le C \|f\|_{C^{k,\alpha}(B_1)},
\end{equation}
for some $C$ depending only on $n$, $k$, and $\alpha$. 

On the other hand, we define $w_2 := u - w_1\in C^{1,\alpha_\circ}(B_1)$ where $\alpha_\circ  =\min\{\eps, \alpha\}$ and (by Lemma~\ref{lem:global_interior_reg})
\[
\|w_2\|_{L^1_\omega(\R^n)}\le \|u\|_{L^1_\omega(\R^n)} + \|w_1\|_{L^1_\omega(\R^n)} \le \|u\|_{L^1_\omega(\R^n)} + C\|f\|_{C^{k,\alpha}(B_1)} . 
\]

Notice that $\sqrtl w_2 = 0$ in $B_1$ so that, by Lemma~\ref{lem:interior_reg_zero_sqrtl} $w_2\in C^\infty(B_1)$ and 
\begin{equation}
\label{eq:w2}
\|w_2\|_{C^\nu(B_{1/2})}\le C_\nu \|w_2\|_{L^1_\omega(\R^n)}\le C_\nu (\|u\|_{L^1_\omega(\R^n)} + \|f\|_{C^{k,\alpha}(B_1)})
\end{equation}
for any $\nu\in \N$ and for some constant $C_\nu$ depending only on $\nu$ and $n$. Combining \eqref{eq:w1}-\eqref{eq:w2} and taking $\nu = k+2$ we get
\[
\begin{split}
\|u\|_{C^{k+1,\alpha}(B_{1/2})}& \le \|w_1\|_{C^{k+1,\alpha}(B_{1})}+\|w_2\|_{C^\nu(B_{1/2})}\\
&  \le C(\|u\|_{L^1_\omega(\R^n)} + C\|f\|_{C^{k,\alpha}(B_1)}),
\end{split}
\]
for some $C$ depending only on $n$, $k$, and $\alpha$, as wanted. 
\end{proof}

\section{Some explicit solutions}
\label{sec:explicit_solutions}\index{Explicit solutions!Square root of the Laplacian}
Let us give some explicit solutions. In the following statements, we denote by $t_+ := \max\{t , 0\}$ the positive part of any number $t\in \R$. 

\begin{figure}
\centering
\makebox[\textwidth][c]{\includegraphics[scale = 1]{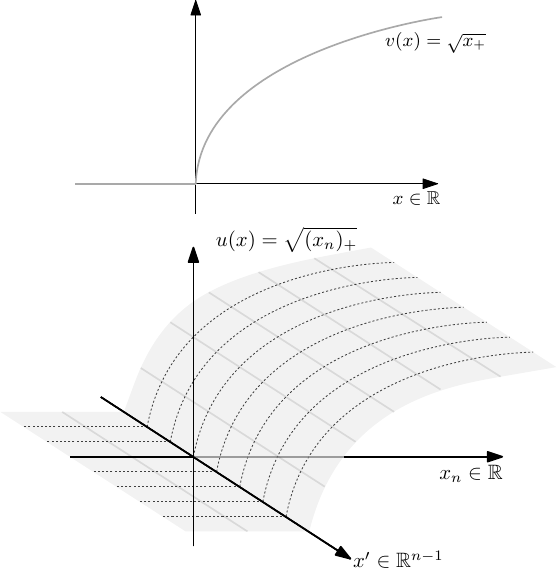}}
\caption{\label{fig:01} Graphical representation of the solutions $v$ and $u$ from Proposition~\ref{prop:explicit_sol_half_space_sqrtl}. They are not differentiable at $x_1 = 0$. }
\end{figure}

We first construct half-space solutions (see Figure~\ref{fig:01}).

\begin{prop}
\label{prop:explicit_sol_half_space_sqrtl}
The function $v:\R\to \R$ defined as 
\[
v(x) = \sqrt{x_+}
\]
satisfies 
\[
{\sqrtl} v = 0\quad\text{in}\quad \{x > 0\}. 
\]
More generally, the function $u:\R^n\to \R$ defined as $u(x) = \sqrt{(x\cdot \be)_+}$ for $\be\in \mathbb{S}^{n-1}$ fixed, satisfies ${\sqrtl} u = 0$ in $\{x\cdot \be > 0\}$. 
\end{prop}

\begin{proof}
We use the extension, Theorem~\ref{thm:extension}. Observe that the extension of $v$ towards $\R^2\cong \C$ is given by $\tilde v(x, y) = \textrm{Re}(z^{\frac12})$, with $z = x+iy$ and where $z^{\frac12}$ is the principal square root of $z$, using the nonpositive real axis as branch cut. Alternatively, in polar coordinates $x = r\cos(\theta)$ and $y = r\sin(\theta)$, we have $\tilde v = r^{\frac12}\cos\left(\frac{\theta}{2}\right)$, with $ r > 0$ and $-\pi < \theta <\pi$. Then, since $\partial_\theta\big|_{\theta = 0} \tilde v = 0$, we get that ${\sqrtl}  v(x) = -\partial_y\big|_{y = 0} \tilde v(x, y) = 0$ when $x > 0$, as we wanted to see.

To see the result in $\R^n$ we just need to realize that the harmonic extension of $u(x)$ in this case is given by $\tilde v(x\cdot \be, y)$. 
\end{proof}

In the second construction (and its corollary) we present solutions in a ball, with constant right-hand side and zero exterior condition. We start with the one-dimensional solution:

\begin{prop}
\label{prop:sqrtlv11D}
The function $v:\R\to \R$ defined as 
\[
v(x) = \sqrt{(1-x^2)_+}
\]
satisfies 
\[
{\sqrtl} v = 1\quad\text{in}\quad |x|< 1. 
\]
\end{prop}
\begin{proof}
We can search for the extension towards $\R^2\cong \C$ of $v$. In this case, it is given by\footnote{This is the unique harmonic extension given by the Poisson kernel, because it is sublinear at infinity. A priori, if we just required harmonicity of the extension we could have subtracted $\lambda y$ for any $\lambda\in \R$, but only $\lambda = 1$ gives the sublinear growth at infinity.}
\[
\tilde v(x, y) = \textrm{Re}\left(\sqrt{1-z^2} + iz\right) = \textrm{Re}\left(\sqrt{1-z^2}\right) - y,
\]
where again we take the principal square root. Observe that the branch cut is now $\{|z|\ge 1 : z = \bar z\}$, that is, real numbers with modulus greater or equal than 1.

Using that $ \textrm{Re}\big(\sqrt{1-z^2}\big)$ is even with respect to the real axis in $B_1$, we get 
\[
{\sqrtl} v(x) = -\partial_y\big|_{y = 0} \tilde v(x, y) = 1\qquad\text{for}\quad -1< x < 1,
\]
as we wanted to see.
\end{proof}

\begin{figure}
\centering
\makebox[\textwidth][c]{\includegraphics[scale = 1]{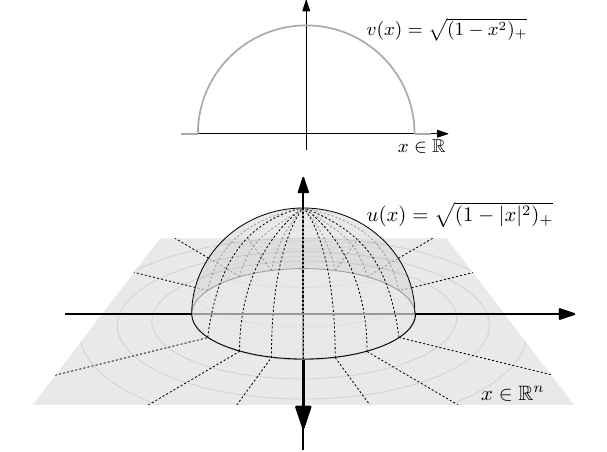}}
\caption{\label{fig:02} Graphical representation of the solutions $v$ and $u$ from Proposition~\ref{prop:sqrtlv11D} and Corollary~\ref{cor:sqrtlv1}.}
\end{figure}

Using polar coordinates, we find the corresponding solution in $\R^n$ (see Figure~\ref{fig:02} for a graphical representation of the solutions): 

\begin{cor}
\label{cor:sqrtlv1}
The function $u:\R^n\to \R$ defined as $u(x) = \sqrt{(1-|x|^2)_+}$ satisfies 
\[
{\sqrtl} u = q_n \quad\text{in}\quad B_1,\quad\text{with}\quad q_n := \sqrt{\pi} \frac{\Gamma\left(\frac{n+1}{2}\right)}{\Gamma\left(\frac{n}{2}\right)}.
\]
\end{cor}

\begin{proof}
We use the result in Proposition~\ref{prop:sqrtlv11D}. We observe that in $\R^n$ we can write
\[
\begin{split}
{\sqrtl} u(x) & = c_n \textrm{P.V.}\int_{\R^n} \frac{u(x) - u(x+y)}{|y|^{n+1}}\, dy \\
& = \frac{c_n}{2} \int_{\mathbb{S}^{n-1}} \textrm{P.V.}\int_{\R} \frac{u(x) - u(x+\rho \theta)}{\rho^{2}}\, d\rho\, d\theta.
\end{split}
\]

If we define $w_{x, \theta} = w_{x, \theta}(\rho) :\R\to \R$ as 
\[
w_{x, \theta}(\rho) := u(x+\rho\theta),
\]
then
\[
c_1 \textrm{P.V.} \int_{\R} \frac{u(x) - u(x+\rho \theta)}{|\rho|^{2}}\, d\rho ={\sqrtl} w_{x, \theta}(0).
\]
Finally, we notice 
\[
w_{x,\theta}(\rho) = A_{x, \theta}\, v\left(\frac{\rho+B_{x,\theta}}{A_{x,\theta}}\right),\quad A_{x, \theta} = \sqrt{1+(x\cdot \theta)^2 - |x|^2},\quad B_{x,\theta} = x\cdot \theta,
\]
so that 
\[
{\sqrtl} w_{x, \theta}(0) = \big({\sqrtl} v\big)\left(\frac{B_{x,\theta}}{A_{x,\theta}}\right) = 1
\]
by the 1-dimensional computation, whenever $|B_{x,\theta}| < A_{x,\theta}$, i.e., $|x|^2 < 1$. Thus,
\[
{\sqrtl}u(x) = \frac{c_n}{2c_1}|\mathbb{S}^{n-1}| = \sqrt{\pi} \frac{\Gamma\left(\frac{n+1}{2}\right)}{\Gamma\left(\frac{n}{2}\right)}\quad\text{in}\quad B_1,
\]
as claimed.
\end{proof}

\section{The fractional Laplacian}
\label{sec:fls}

The square root of the Laplacian is a particular instance of a more general class of nonlocal operators: \emph{the fractional Laplacian}, denoted $\fls$ for $s\in(0,1)$, which is an operator of order $2s$. In this section we  extend the properties introduced for ${\sqrtl}$ in Sections~\ref{sec:harm_extension}-\ref{sec:explicit_solutions} to the context of the fractional Laplacian $\fls$; we refer to \cite{AV19, BV, CLM20, Gar19, MRS} for further properties and results.

We start with the definition. When acting on a smooth function $u\in C^\infty_c(\R^n)$, we define the fractional Laplacian $\fls$ as follows:
\begin{equation}
\label{eq:fractional_Laplacian}\index{Fractional Laplacian}
\begin{split}
{\fls} u(x) &  := c_{n,s} \textrm{P.V.}\int_{\R^n} \frac{u(x) - u(x+y)}{|y|^{n+2s}}\, dy   \\
& \,=  c_{n,s} \textrm{P.V.}\int_{\R^n} \frac{u(x) - u(z)}{|x-z|^{n+2s}}\, dz,
\end{split}
\end{equation}
where
\begin{equation}
\label{eq:cns}
c_{n,s} := 2^{2s} s \frac{\Gamma\left(\frac{n+2s}{2}\right)}{\Gamma(1-s)} \pi^{-\frac{n}{2}}.
\end{equation}
 As for the square root of the Laplacian ($s = \frac12$), due to concerns on integrability at the origin (when $s \ge \frac12$), we use the \emph{Principal Value} of the integral. We also have the alternative definition that takes advantage of the symmetry of the kernel to obtain an expression with second order incremental quotients (that no longer needs $\textrm{P.V.}$):
\[
{\fls} u(x) = \frac{c_{n,s}}{2}\int_{\R^n} \frac{2 u(x) - u(x+y)-u(x-y)}{|y|^{n+2s}}\, dy.
\]
Exactly as for ${\sqrtl}$, if one assumes $u\in C^2(\R^n)\cap L^\infty(\R^n)$ then we can evaluate $\fls u$. In general, though, it is enough to require as global assumption  integrability of $u$ against the appropriate weight:
\begin{defi}
\label{defi:L1omegas}\index{L1ws@$L^1_{\omega_s}$}
We say that $w\in L^1_{\omega_s}(\R^n)$ if 
\[
\|w\|_{L^1_{\omega_s}(\R^n)} := \int_{\R^n}\frac{|w(y)|}{1+|y|^{n+2s}}\, dy<+\infty.
\]
\end{defi}

If we want to evaluate $\fls u$ pointwise at $x\in \R^n$, it is enough for $u$ to be $C^{2s+\alpha}$ around $x$, for some $\alpha > 0$. From now on, we use the convention that if $\beta = k +\alpha$ for some $k\in \N_0$, $\alpha\in (0, 1)$, then $C^\beta := C^{k,\alpha}$ (see Appendix~\ref{app.A}). 

\begin{lem}
\label{lem:laplu_s}\index{Strong solutions!Fractional Laplacian}
Let $s\in (0, 1)$, let $\alpha > 0$ with $\alpha\notin\N$, and let $u\in C^{2s+\alpha}(B_1)\cap L^1_{\omega_s}(\R^n)$. Then ${\fls}u$ is well-defined in $B_1$ and ${\fls}u \in C_{\rm loc}^\alpha(B_1)$, with
\[
\|{\fls} u\|_{C^\alpha(B_{1/2})} \le C \left(\|u\|_{C^{2s+\alpha}(B_1)} + \|u\|_{L_{\omega_s}^1(\R^n)} \right) 
\]
for some $C$ depending only on $n$, $s$, and $\alpha$.
\end{lem}
\begin{proof} It is a small modification of the proof of Lemma~\ref{lem:laplu} (cf. the proof of Lemma~\ref{lem:Lu_2} as well, in a more general context). 
\end{proof}

The operator ${\fls}$ is also positive definite, and the set of invariances for $\fls$ is the same as for ${\sqrtl}$, except that now the $1$-homogeneity becomes $2s$-homogeneity: 
 \begin{lem}
 \label{lem:invariances_s}
 Let $u\in C_{\rm loc}^{2s+\eps}(\R^n)\cap L^1_{\omega_s}(\R^n)$ for some $\eps >0$. Then, the following properties hold for all $x\in \R^n$,
 \begin{enumerate}[(i)]
 \item \label{it:prop_trans_s}\index{Translation invariance!Fractional Laplacian} Translation invariance:
 \[{\fls} [u(x+x_\circ)] = \big({\fls} u\big)(x+x_\circ)\] for any $x_\circ\in \R^n$. 
 \item \label{it:prop_rot_s}\index{Rotation invariance!Fractional Laplacian} Rotation invariance: \[{\fls} [u(Ox)] = \big({\fls} u\big)(Ox)\] for any $O\in \mathcal{O}(n)$ orthogonal transformation.
 \item \label{it:prop_scale_s} \index{Scale invariance!Fractional Laplacian} Scale invariance: \[{\fls} [u(\lambda x)] = |\lambda|^{2s}\big({\fls} u\big)(\lambda x)\] for any $\lambda \in \R$.
 \end{enumerate}
 \end{lem}

Finally, as for the square root of the Laplacian, the previous invariances actually characterize the operator:

\begin{lem}
\label{lem:fls_unique_invariances}
Let $L:\mathcal{S}\to \mathcal{S}$ be any linear operator satisfying properties \ref{it:prop_trans_s}-\ref{it:prop_rot_s}-\ref{it:prop_scale_s} from Lemma~\ref{lem:invariances_s}. Then, there exists a constant $\kappa\in \R$ such that $L = \kappa {\fls}$.
\end{lem}

The proofs of the previous results are analogous to those of Lemmas~\ref{lem:invariances} and~\ref{lem:sqrtl_unique_invariances}. 

\begin{rem}
\label{rem:sqrtlsqrtl_s}\index{Fourier symbol!Fractional Laplacian}
The constant in the definition of $\fls$, \eqref{eq:cns}, is chosen so that the Fourier multiplier of ${\fls}$ is $|\xi|^{2s}$ (see, for example, \cite[Proposition~3.3]{DPV12}, or Example~\ref{example:fract_Lapl} in Section~\ref{sec:Levy}).

Notice also that, thanks to this, the fractional Laplacian can be seen as a pseudo-differential operator (plus a smoothing operator), see \cite{Gru09}, and hence solutions to $\fls u = 0$ are always $C^\infty$, and solutions to $\fls u = f$ gain $2s$-derivatives of regularity (see Theorem~\ref{thm:interior_regularity_fls}).
\end{rem}

\begin{rem}[Semigroup property]
In fact, thanks to the previous remark, we observe that the fractional Laplacian satisfies the semigroup property:
\[
\fls \circ (-\Delta)^t = (-\Delta)^{s+t},
\]
for any $s, t\in (0, 1]$ with $s+t\le 1$. 
\end{rem}

\subsection*{The $a$-harmonic extension}
\index{Harmonic extension!Fractional Laplacian}
For the fractional Laplacian, $\fls$, there is an analogue to the harmonic extension introduced in Section~\ref{sec:harm_extension}. In the PDE community, such an extension is  known as the \emph{Caffarelli--Silvestre extension}, since Caffarelli and Silvestre introduced its use in the context of the fractional Laplacian in \cite{CSext}.  From a probabilistic perspective, it had been observed in \cite{Spi58}  (for $n = 1$ and $s = \frac12$) and in \cite{MO69,Hsu86} (when $s \in (0, 1)$).

In this case, we  deal with the local operator $L_a$  defined for functions $w = w(x, y):\R^n\times \R_+$ as 
\begin{equation}
\label{eq:Ladef}
L_a w := \Delta_x w + \frac{a}{y} \partial_y w + \partial_{yy} w = y^{-a}{\rm div}_{x,y}(y^a \nabla_{x,y} w),\quad a\in (-1,1).
\end{equation}
Observe that $L_0 \equiv \Delta$. Operators of the form $L_a$ with $a\in (-1,1)$ are degenerate elliptic operators, that in its variational form have associated a Muckenhoupt weight $A_2$. Many of the properties that hold for $\Delta$ also hold for this more general class of operators, see \cite{FKS82,FJK83,  HKM93, Kil94}.

The analogue of Theorem~\ref{thm:extension} in the context of the fractional Laplacian $\fls$ is the following. We refer to \cite{CSext, BV, ST10, FL}. 
\begin{thm}
\label{thm:extension_s}
Let $s\in (0,1)$ and let $a := 1-2s \in (-1,1)$. Let $u\in C_{\rm loc}^{2s+\eps}(B_1)\cap L_{\omega_s}^1(\R^n)$ for some $\eps > 0$, and let ${\tilde u}:\R^{n+1}\cap\{y > 0\} \to \R$ be the unique solution to
\begin{equation}
\label{eq:extension_s}
\left\{
\begin{array}{rcll}
L_a \tilde u & = & 0& \quad\text{in}\quad \{y > 0\}\\
\tilde u(x, 0) & = & u(x)& \quad\text{for}\quad x\in \R^n
\end{array}
\right.
\end{equation}
with sublinear growth at infinity, given by the Poisson kernel representation \eqref{eq:poissonkernel0_s}-\eqref{eq:poissonkernel_s}. Then,  
\[
\fls u(x) = -d_s\lim_{y \downarrow 0}y^a \partial_y \tilde u(x, y),\quad \text{for any}~~ x\in B_1,\quad d_s := 2^{2s-1}\frac{\Gamma(s)}{\Gamma(1-s)},
\]
where $\fls$ is given by \eqref{eq:fractional_Laplacian}.
\end{thm}

That is, we   recover the fractional Laplacian $\fls$ by taking the $a$-harmonic extension towards one more dimension (i.e., the extension $\tilde u$ satisfying $L_a\tilde u$ with $L_a$ given by \eqref{eq:Ladef} and $a = 1-2s$), and then taking a ``fractional derivative'' at $y = 0$, up to a constant, $-d_s \lim_{y \downarrow 0}y^a \partial_y \tilde u$. As already observed, the operator $L_a$ is degenerate elliptic  as $y\downarrow 0$; nonetheless, it retains many of the properties of the Laplacian. 

A recurrent heuristic  intuition that it is often useful to guess formulas for the operator $L_a$ is interpreting the expression \eqref{eq:Ladef} as the ``Laplacian'' of a function $u(x, y):\R^n\times \R^{1+a}\to \R$ that is ``radial'' in the last $1+a$ coordinates, namely, $u(x, y) = \tilde u(x, |y|)$ with $\tilde u:\R^n\times\R\to \R$. Whenever $a\in\N_0$, this intuition is exactly true, but for other noninteger values of $a$ only the expression \eqref{eq:Ladef} retains its sense.  Two such examples are as follows:
\begin{itemize}
\item   The fundamental solution for the operator $L_a$ in $\R^{n+1}$ is given by 
\begin{equation}
\label{eq:Ka}
K_{a}(x, y) = C^{(K)}_{n,a} |(x, y)|^{1-a-n},
\end{equation}
for some constant $C^{(K)}_{n,a}$ depending only on $n$ and $a$ (cf. \eqref{eq:K2}).
\item The solution to \eqref{eq:extension_s} is also given in terms of $u$ by means of the Poisson kernel (cf. \eqref{eq:poissonkernel0}-\eqref{eq:poissonkernel}):
\begin{equation}
\label{eq:poissonkernel0_s}
{\tilde u}(x, y) = [P_a(\cdot, y) * u](x) := \int_{\R^n} P_a(x-z, y) u(z)\, dz
\end{equation}
where $P_a(x, y)$ is the Poisson kernel of an upper half-space for the operator $L_a$, \index{Poisson kernel!Fractional Laplacian!Half-space}
\begin{equation}
\label{eq:poissonkernel_s}
P_a(x, y) = C^{(P)}_{n,a}\frac{y^{1-a}}{(|x|^2 +y^2)^{\frac{n+1-a}{2}}},
\end{equation}
for some constant $C^{(P)}_{n,a}$ depending only on $n$ and $a$.  
\end{itemize}

\subsection*{Maximum principle}
\index{Maximum principle!Fractional Laplacian}
Super- and subsolutions to the equation $\fls u=0 $ also satisfy a maximum principle and a comparison principle. The proof is exactly the same as the one of Lemma~\ref{lem:max_principle}.
\begin{lem}[Maximum Principle]
\label{lem:max_principle_s}
Let $\Omega\subset\R^n$ be a bounded open set, and let $u\in C^{2s+\eps}_{\rm loc}(\Omega)\cap C(\overline{\Omega})\cap L^1_{\omega_s}(\R^n)$ for some $\eps>0$. Let us assume that 
\[
\left\{
\begin{array}{rcll}
\fls u & \ge & 0 & \quad\text{in}\quad \Omega,\\
u & \ge & 0 & \quad\text{in}\quad \R^n\setminus \Omega.
\end{array}
\right.
\]
Then $u\ge 0$ in $\R^n$. Moreover, either $u > 0$ in $\Omega$ or $u \equiv 0$ in $\R^n$. 
\end{lem}

As a consequence, we also have a comparison principle (cf. Corollary~\ref{cor:comp_principle}) and the uniqueness of strong solutions (cf. Corollary~\ref{cor:uniqueness_sqrtl}).

\subsection*{Fundamental solution}
\index{Fundamental solution!Fractional Laplacian}

The fundamental solution for the fractional Laplacian, denoted $K_{2s}$, is given by 
\begin{equation}
\label{eq:K2s}
K_{2s}(x) := 
\kappa_{2s,n} |x|^{2s-n},\quad \textrm{with}\quad  \kappa_{2s,n} =  2^{-2s}\frac{\Gamma\left(\frac{n-2s}{2}\right)}{\Gamma(s)} \pi^{-\frac{n}{2}},
\end{equation}
whenever $n\neq 2s$. For $n = 2s = 1$, it is \eqref{eq:K1}. 

\begin{lem}[Fundamental solution]
\label{lem:fundamental_solution_s}
Let $f\in C^\infty_c(\R^n)$, and let $K_{2s}$ be given by \eqref{eq:K2s}. Then 
\[
{\fls} (K_{2s}*f) = f\quad\textrm{in}\quad\R^n.
\]
\end{lem}

The previous lemma follows from \eqref{eq:Ka} by means of the extension theorem, Theorem~\ref{thm:extension_s}. For the value of the constant see, for example, \cite{Buc16}.

\subsection*{Poisson kernel in a ball and mean value property}

The Poisson kernel in a ball in this case is the following. We refer to \cite{Buc16} or \cite{L} for a proof of this result. 

\begin{prop}[Poisson kernel in a ball]
\label{prop:poisson_kernel_ball_s}\index{Poisson kernel!Fractional Laplacian!Ball}
Let $s\in(0,1)$, and let $g\in L^1_{\omega_s}(\R^n)$ be such that $g$ is continuous on $\partial B_1$. Let $u$ be defined as 
\[u(x)=a_{n,s} \int_{\R^n\setminus B_1}\frac{(1-|x|^2)^s\,g(z)}{(|z|^2-1)^s|x-z|^n}\,dz\qquad\textrm{if}\quad x\in B_1,\]
where $a_{n, s} = \Gamma\left(\frac{n}{2}\right)\pi^{-\frac{n}{2}-1}\sin(\pi s)$, and $u(x) = g(x)$ in $\R^n\setminus B_1$. Then $u\in C^\infty(B_1)\cap C(\overline{B_1})$ and it solves
\[
\left\{
\begin{array}{rcll}
\fls u &=&0&\quad\textrm{in}\quad B_1\\
u&=&g&\quad\textrm{in}\quad \R^n\setminus B_1.
\end{array}
\right.
\]
\end{prop}

From the Poisson kernel in a ball, as in the case $s = \frac12$, we can derive the corresponding mean value property for $\fls$: 

\begin{prop}[Mean value property]
\label{prop:meanvalue_s}\index{Mean value property!Fractional Laplacian}
Let $s\in (0, 1)$. Let $u\in C_{\rm loc}^{2s+\eps}(B_1)\cap L^1_{\omega_s}(\R^n)$ for some $\eps> 0$, and assume that $\fls u = 0$ in $B_1$.  Then, for any $r\in (0, 1)$, 
\[
u(0) = a_{n,s} \int_{\R^n\setminus B_r} \frac{r^{2s} u(z)}{(|z|^2-r^2)^s|z|^n}\, dz.
\]
 Moreover, there exists a nonincreasing function $\omega_s(t):[0,+\infty)\to [0,+\infty)$ such that 
\[
u(0) = \int_{\R^n} u(z) \omega_s(|z|)\, dz,
\]
with 
\begin{equation}
\label{eq:omega_estimates_s}
\frac{C^{-1}}{1+|t|^{n+2s}} \le \omega_s(t) \le \frac{C}{1+|t|^{n+2s}}
\end{equation}
for some $C$ depending only on $n$ and $s$. 
\end{prop}

\begin{proof}
The first part is an immediate consequence of the rescaled version of Proposition~\ref{prop:poisson_kernel_ball_s} (cf. Remark~\ref{rem:rescaledpoisson}). For the second part, we can proceed exactly as in Corollary~\ref{cor:meanvalue} to obtain
\[
\omega_s(t) := \left\{
\begin{array}{ll}
a_{n,s} \int_0^1 {\rho^{n+2s-1}}{(1-\rho^2)^{-s}}\, d\rho \qquad & \textrm{if}\quad 0\le t\le 1,\\
a_{n,s} \int_0^{\frac{1}{t}} {\rho^{n+2s-1}}{(1-\rho^2)^{-s}}\, d\rho \qquad & \textrm{if}\quad t\ge 1,
\end{array}
\right.
\]
from which \eqref{eq:omega_estimates_s} follows. 
\end{proof}

\subsection*{The Harnack inequality}

As a consequence of the mean value property, we deduce the Harnack inequality. 

\begin{prop}
\label{prop:Harnack_fls}\index{Harnack's inequality!Fractional Laplacian}
Let $s\in(0,1)$, and let $u\in C^{2s+\eps}(B_1)\cap L^1_{\omega_s}(\R^n)$ for some $\eps >0$ be such that
\[
\left\{
\begin{array}{rcll}
u &\ge& 0& \quad\textrm{in}\quad \R^n\\
\fls \, u& =& 0&\quad\textrm{in}\quad{B_1}.
\end{array}
\right.
\]
Then, 
\[
\sup_{B_{1/2}} u \le C \inf_{B_{1/2}} u
\]
for some constant $C$ depending only on $n$ and $s$.  Moreover, both quantities $\sup_{B_{1/2}} u $ and $\inf_{B_{1/2}} u$ are comparable to $\|u\|_{L^1_{\omega_s}(\R^n)}$. 
\end{prop}
\begin{proof}
The proof is the same as the second proof of Proposition~\ref{prop:Harnack_sqrt} using Proposition~\ref{prop:meanvalue_s}. 
\end{proof}

\subsection*{Interior regularity}

\index{Interior regularity!Fractional Laplacian}
Concerning the interior regularity of solutions, we have the following result. Notice that, since $\fls$ is an elliptic operator of order $2s$ (see Lemmas~\ref{lem:laplu_s} and \ref{lem:invariances_s}),   we expect solutions of $\fls u = f$ to ``gain $2s$ derivatives'' with respect to $f$. 
\begin{thm}[Interior regularity for $\fls$]
\label{thm:interior_regularity_fls}
Let $s\in (0, 1)$, and let $u\in C^{2s+\eps}(B_1)\cap L^1_{\omega_s}(\R^n)$ for some $\eps > 0$ satisfy
\[
\fls u = f\quad\text{in}\quad B_1,
\]
for some $f\in C^\alpha(B_1)$ with $\alpha > 0$ and $\alpha+2s\notin\N$. Then, $u\in C_{\rm loc}^{\alpha+2s}(B_1)$ and 
\[
\|u\|_{C^{\alpha+2s}(B_{1/2})}\le C \left(\|f\|_{C^\alpha(B_1)} + \|u\|_{L^1_{\omega_s}(\R^n)}\right),
\]
for some $C$ depending only on $\alpha$, $s$, and $n$. 
\end{thm}

The proof of this result is a small modification of that of Theorem~\ref{thm:interior_regularity_sqrtl}.

\subsection*{Some explicit solutions} \index{Explicit solutions!Fractional Laplacian} For the fractional Laplacian, we also have some explicit solutions. The first one is an explicit solution which will play an important role in Section~\ref{sec:bdryregularity}. 

\begin{prop}
\label{prop:onedbarrier}
Let $s\in (0, 1)$. The function $u:\R^n\to \R$ defined as 
\[
u(x) = ( x\cdot \be)_+^s
\]
 for $\be\in \mathbb{S}^{n-1}$ fixed satisfies 
\[
\fls u = 0\quad\text{in}\quad \{x\cdot \be > 0\}. 
\]
\end{prop}
\begin{proof}
 As in the proof of Proposition~\ref{prop:explicit_sol_half_space_sqrtl}, we can look for the extension when $n = 1$ for $v(x) = (x_+)^s$, which in this case is (in polar coordinates $(r, \theta)\in [0, \infty)\times (-\pi, \pi)$):
\[
\tilde v (r, \theta) = r^s\left[\cos\left(\theta/2\right)\right]^{2s}.
\]
Then $\tilde v$ satisfies ${\rm div}(y^a \nabla u ) = 0$ in $\{y > 0\}$. Since $\partial_\theta\big|_{\theta = 0}  \tilde v = 0$, by Theorem~\ref{thm:extension_s} we get that $\fls v = 0$ when $x  > 0$.

The result in $\R^n$ now follows as in Proposition~\ref{prop:explicit_sol_half_space_sqrtl}. 
\end{proof}

The second one, is a solution in the unit ball. 

\begin{prop}
\label{prop:count_ex_s}
Let $s\in (0, 1)$. The function $u:\R^n\to \R$ defined as 
\[
u(x) = (1-|x|^2)_+^s
\]
 satisfies 
\[
\fls  u = q_{n,s} \quad\text{in}\quad B_1,\quad\text{with}\quad q_{n,s} := 2^{2s} \Gamma(1+s) \frac{\Gamma\left(\frac{n+2s}{2}\right)}{\Gamma\left(\frac{n}{2}\right)}.
\]
\end{prop}
\begin{proof}
In dimension $n = 1$ this is a (nontrivial) computation, which can be performed explicitly (done, for example, in \cite{G}; see also \cite{BV, Dyd12}). In dimensions $n\ge 2$, we can proceed in polar coordinates as in the proof of Corollary~\ref{cor:sqrtlv1}.
\end{proof}

\subsection*{Classification of 1D solutions}
\index{Classification 1D solutions of $\fls$}
Let us finish this section with the following classification result, that will be useful later on in the study of the higher-order boundary regularity of solutions.

\begin{thm}
\label{thm:class_1d_fls}
Let $s\in (0, 1)$, and let   $u\in C(\R)\cap C^{\infty}((0, \infty))$  satisfy
\[
\left\{
\begin{array}{rcll}
|u(x)| & \le & C\big(1+|x|^{2s-\eps}\big)& \quad\text{in}\quad \R\\
\fls u & = & 0 & \quad \text{in}\quad (0, \infty)\\
u & = & 0 & \quad\text{in}\quad (-\infty, 0]
\end{array}
\right.
\]
for some $C$ and $\eps > 0$. Then,
\[
u(x) = \kappa (x_+)^s
\]
for some $\kappa\in \R$. 
\end{thm}

In order to prove it, we will use the following result:

\begin{lem}
\label{lem:itsolves}
Let $s\in (0, 1)$, and let us consider polar coordinates in~$\R^2$, $(x, y ) = (r\cos\theta, r\sin\theta)$ for $r\ge 0$ and $\theta \in (-\pi, \pi]$. Let us define, for each $m\in \Z$, $m\ge -1$, 
\[
\varphi_m(r, \theta) := r^{s+m} \Theta_m(\theta),\qquad \Theta_m(\theta) = c_m |\sin\theta|^s P_m^s(\cos\theta),
\]
where $P_m^s$ is the associated Legendre function of the first kind, and $c_m$ is chosen so that $\int_{-\pi}^\pi|\Theta_m(\theta)|^2|\sin\theta|^{1-2s}\, d\theta=1$. Then, each $\varphi_m$ satisfies
\[
\left\{
\begin{array}{rcll}
L_a \varphi_m & = & 0& \quad \text{in}\quad \{y \neq 0\},\\
\lim_{y\downarrow 0}\, y^a\partial_y \varphi_m(x, y) & = & 0& \quad \text{for}\quad x > 0,\\
\varphi_m(x, 0) & = & 0& \quad\text{for}\quad x \le 0,
\end{array}
\right.
\]
where $L_a$ is given by \eqref{eq:Ladef} with $a = 1-2s$. Moreover,  the functions $\{\Theta_m(\theta)\}_{m\in \N_0}$ are a complete orthogonal system in the subspace of even functions in the weighted space $L^2\big((-\pi, \pi); |\sin\theta|^{1-2s}\big)$. 
\end{lem}

\begin{proof}
The proof of this lemma is a computation, using that $P_m^s(t)$ satisfies the second order ODE:
\[
(1-t^2) \big(P_m^s(t)\big)''-2t\big(P_m^s(t)\big)' + \left(m+m^2-\frac{s^2}{1-t^2}\right) P_m^s(t) = 0. 
\] 
We refer to \cite[Lemma 6.1]{RS-Duke} for more details on the proof (see also \cite[Proposition A.1]{FS} for another proof). 
\end{proof}

Using the previous lemma, we can proceed with the classification result: 

\begin{proof}[Proof of Theorem~\ref{thm:class_1d_fls}]
Let $\tilde u(x, y)$ be the extension of $u$ towards $\{y > 0\}$ given by the Poisson kernel representation \eqref{eq:poissonkernel0_s} (and extended evenly to $y< 0$),
\[
{\tilde u}(x, y) = \big[P_a(\cdot, y) * u\big](x) \quad
P_a(x, y) = \frac{C^{(P)}_{1,a}y^{1-a}}{(x^2 +y^2)^{\frac{2-a}{2}}}, \quad a = 1-2s. 
\]
Since $u(x) \le C(1+|x|^{2s-\eps})$, we can explicitly bound for $y > 0$
\[
|\tilde u(x,y)|\le C y^{1-a}\int_{-\infty}^\infty \frac{1+|z|^{2s-\eps}}{\left((x-z)^2+y^2\right)^{\frac{2-a}{2}}}\, dz = C\int_{-\infty}^\infty\frac{1+|x-y\xi|^{2s-\eps}}{\left(1+\xi^2\right)^{\frac{2-a}{2}}}\, d\xi.
\]
Using that $|x-y\xi|^{2s-\eps}\le 4(|x|^{2s-\eps}+y^{2s-\eps}|\xi|^{2s-\eps})$, that $2-a-2s+\eps > 1$, and that $|a|^p + |b|^p\le C_p (|a|^2+|b|^2)^{\frac{p}{2}}$ for any $a, b\in \R$, $p \in (0, 2)$, we get
\begin{equation}
\label{eq:growthcontrolparseval}
|\tilde u(x, y)|\le C \left(1+|x|^{2s-\eps}+|y|^{2s-\eps}\right) \le C\left(1+|(x, y)|^{2s-\eps}\right),
\end{equation}
so that $\tilde u$ has the same growth as $u$. 

Moreover, for each $R$ fixed, we can write (thanks to Lemma~\ref{lem:itsolves}) 
\[
\tilde u(R\cos\theta, R\sin\theta) = \sum_{m\ge 0} a_m(R) R^{s+m}\Theta_m(\theta),
\]
and since $L_a \tilde u = 0$ in $\R^2\setminus \{x \le 0, y = 0\}$, by uniqueness of solutions $a_m(R) = a_m$ is independent of $R$. 
Finally, by Parseval's identity (recall $\|\Theta_m\|_{L^2((-\pi, \pi); |\sin\theta|^{1-2s})} = 1$) we have
\begin{equation}
\label{eq:growthcontrolparseval2}
\int_{\partial B_R} |\tilde u(x, y)|^2 |y|^a\, d\sigma = \sum_{m \ge 0} a_m^2 R^{2(s+m)+1+a} .
\end{equation}
From the growth control \eqref{eq:growthcontrolparseval} we also know that for $R \ge 1$,
\[
\int_{\partial B_R} |\tilde u(x, y)|^2 |y|^a\, d\sigma \le C (1+R^{2s-\eps})^2R^{1+a} \le C R^{2(2s-\eps) +1+a}. 
\]
Combined with \eqref{eq:growthcontrolparseval2} this yields that $a_m = 0$ whenever $m+s \ge 2s-\eps$, that is, for any $m\ge 1$. 
Hence, $\tilde u(r\cos\theta, r\sin\theta) = a_0 r^s\Theta_0(\theta)$ and in particular, $u(x) = \kappa (x_+)^s$ for some $\kappa\in \R$. 
\end{proof}

\addtocontents{toc}{\protect\setcounter{tocdepth}{2}}

%% file: chap2.tex
%
%
%

\chapter{Linear integro-differential equations}
\label{ch:2}

 In this chapter, we mostly  study general elliptic operators of order $2s$ that are linear and translation invariant (with no $x$-dependence). 
 In the case of (local) second order elliptic PDE, this would mean that $\L = -\sum_{ij} a_{ij} \partial_{ij} $, an operator with constant coefficients, which after an affine change of variables becomes the Laplacian.
In other words, for second order elliptic PDE, the class of linear and translation invariant operators essentially consists of only one operator: $\L=-\Delta$.

For nonlocal equations the situation is \emph{very} different. The corresponding class of linear operators with ``constant coefficients'' is {extremely} rich, and presents interesting features that do not appear for the fractional Laplacian $\fls$.

We will study linear operators of the form:
  \begin{equation}
 \label{eq:L_def_lin}\index{Linear operator}
 \begin{split}
 \L f(x) & = {\rm P.V.}\int_{\R^n} \bigl(f(x) - f(x+y)\bigr)K (dy) \\
 & =\frac12 \int_{\R^n} \bigl(2f(x) - f(x+y) - f(x-y)\bigr)K (dy),
 \end{split}
 \end{equation}
 for some (nonnegative) measure $K$, the \emph{L\'evy measure} (or \emph{kernel}) of the operator, that will be symmetric (see Definition~\ref{defi:sym_meas}) and satisfy certain ellipticity conditions.

\section{L\'evy processes and classes of kernels}
\label{sec:Levy}
 Operators of the form \eqref{eq:L_def_lin} are integro-differential operators, and they appear as infinitesimal generators of L\'evy processes. 
Let us briefly explain how, and we refer to more specialized books for further insight in this direction,  \cite{Fel71, ST94, Bertoin, Sat99, Eva13, KS16, Sch}.

 As we will see below (and as we have seen in Chapter~\ref{ch:fract_Lapl} for the fractional Laplacian), operators like \eqref{eq:L_def_lin} can be alternatively defined by their \emph{Fourier multiplier} or \emph{Fourier symbol}, given by 
\begin{equation}
\label{eq:fourier_symb0}
\mathcal{F}(\L u) (\xi) = \A(\xi) \mathcal{F}(u)(\xi)
\end{equation}
for some function $\A:\R^n\to \R$. We say that $\mathcal{A}$ is the \emph{Fourier multiplier} or \emph{Fourier symbol} of $\L$, and it is explicit in terms of $K$, 
\begin{equation}
\label{eq:fourier_symb}\index{Fourier symbol!General operators}
\A(\xi) =  \int_{\R^n} \bigl(1-\cos(y\cdot \xi)\bigr)K(dy);
\end{equation}
see \eqref{eq:pseudo}-\eqref{eq:sym_psi} below.

\subsection{L\'evy processes and infinitely divisible distributions}
\label{ssec:Levy}
A stochastic process $(X_t)_{t \in T}$ is a family of random variables in $\R^n$ (on a common probability space) indexed by $t\in T$. In this exposition we will consider $ T = [0, \infty)$ to be regarded as \emph{time}. 

Stochastic processes are widely used as models of systems and phenomena that evolve in a random way, and appear in particular in many areas of physics, biology, or even finance.
 The Brownian motion is the most famous stochastic process, and together with the Poisson process, the Cauchy process, or the Gamma process, they belong to a wider class known as \emph{L\'evy processes}.

Let us start with a definition:
\begin{defi}
\label{defi:Levy}\index{L\'evy process}
A stochastic process $(X_t)_{t\ge 0}$ on $\R^n$ is a \emph{L\'evy process} if it satisfies the following properties:
\begin{enumerate}
\item $X_0 = 0$ almost surely.

\item (Independent increments) For any $0\le t_1 < t_2 < \dots < t_n$, the random variables $X_{t_1},\ X_{t_2}-X_{t_1},\ X_{t_3}-X_{t_2},\ \dots,\ X_{t_n}-X_{t_{n-1}}$ are mutually (or jointly) independent. 

\item (Stationary increments) For any $s < t$, the random variable $X_t-X_s$ is equal in distribution to $X_{t-s}$. 

\item (Stochastic continuity) $X_t$ is \emph{continuous in probability}, that is, for any $\eps > 0$ and $t \ge 0$, $\lim_{s\to t} P(|X_s-X_t| > \eps) = 0$.
\end{enumerate}
\end{defi}

Sometimes, one also adds the condition that $t\mapsto X_t$ is almost surely right-continuous in $t\ge 0$ with left limits in $t > 0$ (alternatively, $t\mapsto X_t$ is \emph{c\`adl\`ag}). In fact, given any \emph{L\'evy process} $X_t$, there always exists a modification of $X_t$ (or a stochastically equivalent process) satisfying this property; see, e.g., \cite{Sat99}. For this reason, we will also assume that our L\'evy processes are c\`adl\`ag.

The study of L\'evy processes and their characterization is done by means of  the concept of \emph{infinitely divisible distributions} (which also play an important role for generalizations of the central limit theorem and in additive processes): 
\begin{defi}[Infinitely divisible distributions]\index{Infinitely divisible distribution}
Let $\mu\in \mathcal{P}(\R^n)$ be a probability measure. We say that $\mu$ is \emph{infinitely divisible} if, for any $m\in \N$, there exists $\mu_m\in \mathcal{P}(\R^n)$ such that
\[
\mu = \underbrace{\mu_m * \dots * \mu_m}_{m}.
\]
Alternatively, we say that a random variable $X$ is infinitely divisible if, for any $m\in \N$, there exist a random variable $Y_m$ such that 
\[
X \overset{d}{=} Y^{(1)}_m +\dots Y^{(m)}_m
\]
where $Y^{(1)}_m,\dots, Y^{(m)}_m$ are independent copies of $Y_m$, and $\overset{d}{=}$ denotes equality in distribution (that is, ${\rm Law}(X)$ is an infinitely divisible distribution).  
\end{defi}

\begin{example}
The most common infinitely divisible distribution is the Gaussian distribution. That is, if $X$ is a random variable with distribution $N(0, 1)$ (normal with zero mean and variance one), then for each $m$ the previous statement holds with $Y  \sim  N(0, 1/m)$. This example will appear recurrently.
\end{example}

Observe that  if $(X_t)_{t \ge 0}$ is a L\'evy process then, for any $t\ge 0$, the random variable $X_t$ (or its distribution) is infinitely divisible, since for every $m \in \N$ we have
\[
X_t \overset{d}{=} \left(X_{\frac{t}{m}}-X_0\right) + \left(X_{t\frac{2}{m}}-X_{\frac{t}{m}}\right) +  \dots+ \left(X_{t}-X_{t\frac{m-1}{m}}\right).
\]

In fact, the relationship between L\'evy processes and infinitely divisible distributions is an equivalence (see \cite[Theorem 7.10]{Sat99}):

\begin{thm}
If $\mu\in \mathcal{P}(\R^n)$ is an infinitely divisible distribution, then there is a L\'evy process $(X_t)_{t \ge 0}$ such that $X_1  \sim  \mu$. Moreover, if $(X_t')_{t\ge 0}$ is another L\'evy process such that $X_1 \overset{d}{=} X_1'$, then $(X_t)_{t\ge 0} \overset{d}{=} (X_t')_{t\ge 0}$. 
\end{thm}

For example, when $\mu = \frac{1}{\sqrt{2\pi}} e^{-\frac{x^2}{2}}\, dx$ (the normal distribution, $N(0, 1)$), the corresponding L\'evy process is the \emph{Brownian motion} (also called \emph{Wiener process}). 
In this case, the process has continuous paths. 

Thanks to the previous result, we have that the characterization of L\'evy processes is equivalent to the characterization of infinitely divisible distributions. In the latter case, such characterization is accomplished through the \emph{L\'evy-Khintchine formula}, which gives a representation of the characteristic function of any infinitely divisible measure. 

We recall that the characteristic function of a probability measure $\mu$, denoted $\hat{\mu}$, is its (inverse) Fourier transform:
\[
\hat{\mu}(\xi)  = \E_{Y\sim \mu}\left[e^{i \xi\cdot Y}\right] =  \int_{\R^n} e^{i \xi\cdot y}\mu(dy). 
\]

\begin{thm}[L\'evy-Khintchine formula]\label{thm-LK}\index{L\'evy-Khintchine formula}
Let $\mu\in \mathcal{P}(\R^n)$ be an infinitely divisible probability measure. Then 
\[
\hat{\mu}(\xi) = e^{\Psi(\xi)},
\]
where $\Psi$ is the \emph{characteristic exponent}, given by 
\begin{equation}
\label{eq:Psi}
\Psi(\xi) = -\frac12 \xi\cdot A\xi + i b\cdot \xi + \int_{\R^n} \left(e^{i \xi\cdot y} - 1- i \xi\cdot y\chi_{B_1}(y)\right) \nu(dy).
\end{equation}
Here, $A\in \R^{n\times n}$ is symmetric and nonnegative-definite, $b\in \R^n$, and $\nu$ is a measure on $\R^n$ satisfying
\begin{equation}
\label{eq:cond_nu}
\nu(\{0\}) =0\qquad \int_{\R^n} \min\{1, |y|^2\}\nu(dy) < \infty.
\end{equation}
Conversely, any $(A, b, \nu)$ with the previous properties generates an infinitely divisible distribution. 
\end{thm}

\begin{rem}
Observe that the conditions on $\nu$ are necessary to make sense of the integral defining $\Psi(\xi)$. This is because, for any $\xi, y\in \R^n$, 
\[
\big|e^{i\xi\cdot y} - 1\big|\le 2\quad\text{in}\quad \R^n\setminus B_1,\qquad\big|e^{i\xi\cdot y} - 1-i\xi\cdot y\big|\le C|\xi \cdot y|^2\quad\text{in}\quad B_1.
\]
This is however not the only way to make sense of the integral. One could take instead $\left(e^{i \xi\cdot y} - 1- i c(y)\xi\cdot y \right)$ as the integrand, under appropriate conditions on $c(y)$, by changing accordingly the value of $b$ (so that, for each $c(y)$, there is a definition of $b$). Other possible $c(y)$ are $c(y) = \chi_{B_\delta}(y)$ for $\delta > 0$ or $c(y) = (1+|y|^2)^{-1}$. 
\end{rem}

\begin{rem}
\label{rem:p_nonnegative} 
As a consequence of Theorem \ref{thm-LK}, if we know that the Fourier transform $\hat{\mu}$ of a distribution is given by $\hat{\mu}(\xi) = e^{\Psi(\xi)}$ with $\Psi$ of the form \eqref{eq:Psi}-\eqref{eq:cond_nu} (observe $\Psi(-\xi)$ is of the same form), then necessarily $\mu$ is a probability measure (in particular, it is nonnegative). 
\end{rem}

\begin{example}
\label{example:gaussian}
When we take as infinitely divisible distribution the Gaussian distribution, $N(0, 1)$, then in $\Psi$ above we have $A = {\rm Id}$, $b = 0$, and $\nu \equiv 0$. 
\end{example}

 The triple $(A, b, \nu)$ is called a \emph{generating triplet} for $\mu$, and characterizes the infinitely divisible distribution $\mu$ (and hence, also characterizes the law of the corresponding L\'evy process). Given a L\'evy process $(X_t)_{t\ge 0}$, the characteristic exponent of the distribution of $X_1$, $\Psi$ in the L\'evy-Khintchine formula, is called the \emph{characteristic exponent} of the process. The characteristic function of $X_t$ for any $t\ge 0$ is then $e^{t\Psi(\xi)}$.

 Justified by Example~\ref{example:gaussian} above, the matrix $A$ is called the \emph{Gaussian covariance matrix}, and the measure $\nu$ is the \emph{L\'evy measure}. When $A = 0$, we say that $(X_t)_{t\ge 0}$ is \emph{purely non-Gaussian} process.  
 Conversely, $\nu$ is zero if and only if $\mu$ is a Gaussian distribution (not necessary centered or with unit variance). 
 
 We will focus our attention on symmetric processes:
 \begin{defi}\index{Symmetric measure}\label{defi:sym_meas}
 We say that a measure $\mu$ on $\R^n$ is \emph{symmetric} if $\mu(B) = \mu(-B)$ for any Borel set $B\subset \R^n$.  We say that a stochastic process $(X_t)_{t\ge 0}$ is \emph{symmetric} if $X_t \overset{d}{=} -X_t$ for all $t\ge 0$. 
 \end{defi}
 
 We will use that a L\'evy process with triple $(A, b, \nu)$ is symmetric if and only if $b = 0$ and $\nu$ is symmetric. In this case, the characteristic exponent \eqref{eq:Psi} can be written as 
\begin{equation}
\label{eq:sym_psi}
\begin{split}
\Psi(\xi) & = -\frac12 \xi\cdot A\xi + {\rm P.V.}\,\int_{\R^n} \big(e^{i \xi\cdot y} - 1\big) \nu(dy)\\
& = -\frac12 \xi\cdot A\xi - \int_{\R^n} \bigl(1-\cos(\xi\cdot y)\bigr) \nu(dy),
\end{split}
\end{equation}
where we have taken advantage of the symmetry of the measure $\nu$ to remove the ${\rm P.V.}$ in the second expression (since $2 \cos(\xi\cdot x) = e^{i\xi \cdot x}+e^{-i\xi\cdot x}$).

\subsection{Infinitesimal generators}
 \label{ssec:inf_gen} \index{Infinitesimal generator}
 
 Let us introduce the notion of infinite\-simal generator of a semigroup, and in particular, of a L\'evy process. We refer to \cite{Bertoin,KS16, Eva13}  for further details. 
 
 Given a L\'evy process $(X_t)_{t \ge 0}$ on $\R^n$,  the measure kernel
 \[
 p_t(x, B) := P(X_t + x\in B) \quad\text{for}\quad t \ge 0, ~~x\in \R^n,~~B\subset \R^n~\text{a Borel set}
 \]
 is what is known as a \emph{Markov transition function} (or \emph{kernel}).  In particular, it defines a linear operator on bounded and measurable functions $f:\R^n\to \R$, 
 \begin{equation}
 \label{eq:Ptdef}
 P_t f(x) = \E[f(x+X_t)] = \int_{\R^n} f(y) p_t(x, dy)\quad\text{for}\quad t \ge 0, ~~x\in \R^n.
 \end{equation}
 
 Observe that, since the characteristic function of $X_t$ is $e^{t\Psi(\xi)}$, we have
 \[
 P_t(e^{i\xi\cdot x}) = e^{i\xi \cdot x} e^{t\Psi(\xi)}, 
 \]
and then $P_t$ is explicit in the Fourier space,
 \begin{equation}
 \label{eq:semigroup_fourier}
 \mathcal{F}(P_t f)(\xi) = \mathcal{F}(f)(\xi) e^{t\Psi(\xi)},
 \end{equation}
 where we recall that $\Psi(\xi)$ is given by \eqref{eq:Psi}. 
 
 The operators $(P_t)_{t\ge 0}$ form a semigroup  (that is, $P_{t+s} = P_t\circ P_s$ and $P_0 = {\rm id}$, which follows from \eqref{eq:semigroup_fourier}) acting on bounded measurable functions. 
Moreover, they are sub-Markovian (i.e., $0 \le f\le 1$ implies $0 \le P_t f \le 1$), contractive ($\|P_t f\|_{L^\infty(\R^n)}\le \|f\|_{L^\infty(\R^n)}$), conservative ($P_t 1 = 1$), and they are strongly continuous Feller operators. That is, if $C_0(\R^n)$ denotes the space of continuous functions vanishing at infinity, then $P_t: C_0(\R^n)\to C_0(\R^n)$ and $\lim_{t\downarrow 0} \|P_t f - f\|_{L^\infty(\R^n)} = 0$ for all $f\in C_0(\R^n)$. 
 
 As a Feller semigroup, one can define for $(P_t)_{t\ge 0}$ the corresponding \emph{infi\-nitesimal generator}, which is the linear operator given by\footnote{We define it with a minus sign in front for future convenience.} 
 \begin{equation}
 \label{eq:def-infinit-generat}
 -\L f := \lim_{t\downarrow 0} \frac{P_t f - f}{t},\quad\text{for all}\quad f\in \D(-\L),
 \end{equation}
 where 
 \[
 \D(-\L):= \left\{ f\in C_0(\R^n) : \lim_{t\downarrow 0}\frac{P_t f - f}{t}~~\text{exists as uniform limit} \right\}.
 \]
 Since $(X_t)_{t\ge 0}$ is a L\'evy process,  we have $C^\infty_c(\R^n)\subset \D(-\L)$.  We have, in particular, that for $t \ge 0$, 
 \begin{equation}
 \label{eq:heat_eq_p}
 \frac{d}{dt} P_t  f = -\L P_t f = -P_t \L f \quad\text{for}\quad f\in C_c^\infty(\R^n)
 \end{equation}
 (this justifies the use of the symbolic notation $P_t = e^{-t\L}$). From \eqref{eq:semigroup_fourier} and the linearity of the Fourier transform, the infinitesimal generator becomes a Fourier multiplier\footnote{General operators of the form \eqref{eq:pseudo} are called \emph{pseudodifferential operators} ($\psi$do), we refer to \cite{Hor07} or \cite{Gru22} for a nice introduction.}
 \begin{equation}
 \label{eq:pseudo}
 \mathcal{F}(\L f) (\xi)  = -\Psi(\xi) \mathcal{F}(f)(\xi) \quad\text{for all}\quad f\in C^\infty_c(\R^n). 
 \end{equation}
 (I.e., $\A = -\Psi$ and $K = \nu$ in \eqref{eq:fourier_symb0}-\eqref{eq:fourier_symb}.)  As a consequence, since we have an explicit expression for $\Psi$, \eqref{eq:Psi}, we deduce how $\L$ acts on $C^\infty_c(\R^n)$ functions:
 \begin{equation}
 \label{eq:L_def}
 \begin{split}
 \L f(x) & = - \frac12 {\rm div}\big(A\nabla f(x)\big) - b\cdot \nabla f(x)\\
 &  \quad - \int_{\R^n} \big(f(x+y) - f(x) - \nabla f(x) \cdot y \chi_{B_1}(y)\big)\nu (dy). 
 \end{split}
 \end{equation}
 That is, $\L$ has three terms: a local diffusion term, a drift term, and a nonlocal diffusion term.  For symmetric L\'evy processes \eqref{eq:sym_psi} we have 
\[
 \L f(x)  =  - \frac12 {\rm div}\big(A\nabla f(x)\big) + {\rm P.V.} \int_{\R^n} \big(f(x) - f(x+y)\big)\nu (dy),
\]
where $\nu$ is symmetric.
 We call any operator of the form \eqref{eq:L_def} (with $A\geq0$ and $\nu$ satisfying \eqref{eq:cond_nu}) a \emph{L\'evy operator}.

We will focus our attention on processes that are purely nonlocal, and as such we will always assume that $A = 0$, so there is no second-order term\footnote{When $A>0$, the corresponding operators are of second order with lower-order nonlocal terms. We refer to \cite{GM02} for a study of such operators.}. In this case, the operator $\L$ is basically a convolution against a measure $\nu(dy)$, which will be often referred as the \emph{kernel} of the operator.

 \begin{example}
 When $(X_t)_{t\ge 0}$ is a Brownian motion (that is, $X_1$ is a  $N(0, 1)$, or alternatively, the characteristic exponent has $ A = {\rm Id}$, $b = 0$, $\nu= 0$), the corresponding infinitesimal generator of the process has $\frac12 |\xi|^2$ as Fourier symbol. That is, $\L = -\frac12 \Delta$ in this case. 
 \end{example}

 \begin{example}
 \label{example:fract_Lapl}\index{Fourier symbol!Fractional Laplacian}
 The most famous non-Gaussian example is in the case $A = 0$, $b = 0$, and $\nu(dy)= \frac{c_\circ}{|y|^{n+2s}}\, dy$ for some $s\in (0, 1)$, and $c_\circ > 0$. Observe that this is an admissible L\'evy measure $\nu$, since conditions \eqref{eq:cond_nu} are satisfied. 
 The corresponding operator in this case is (a multiple of) the \emph{fractional Laplacian}, denoted $\fls$, since the corresponding characteristic exponent is now proportional to $-|\xi|^{2s}$.  (This was already observed in Chapter~\ref{ch:fract_Lapl}, in Remark~\ref{rem:sqrtlsqrtl} or in Remark~\ref{rem:sqrtlsqrtl_s}.) In particular, thanks to \eqref{eq:sym_psi}, we can check that the constant $c_\circ$ such that the characteristic exponent is exactly $-|\xi|^{2s}$ is given by 
 \[
 c_\circ^{-1} = \int_{\R^n}\frac{1-\cos(y_1)}{|y|^{n+2s}}\, dy,
 \]
 which coincides with $c_{n, s}$ in \eqref{eq:cns}.
\end{example}

 \subsection{A characterization via comparison principle}
\index{Comparison principle!L\'evy processes}
Infinitesimal generators of L\'evy processes define a bounded, linear, and translation invariant operator
\[
\L : C^2_c(\R^n) \to C(\R^n). 
\]
Moreover, it follows from the definition \eqref{eq:def-infinit-generat} that such operator satisfies the following \emph{Global Minimum Principle} (GMP): 
\begin{equation}
\label{eq:GMP}
\text{If $f\in C^2_c(\R^n)$ with $f\ge 0$ and $f(x_\circ) = 0$, then $\L f(x_\circ) \le 0$.}
\end{equation}

Notice that we are \emph{not} using the L\'evy-Khintchine formula in order to deduce these properties of $\L$.

Interestingly, it turns out that these properties, combined with the fact that  $\L 1 = 0$, completely  characterize infinitesimal generators of L\'evy processes. 
This is due to the following result:

\begin{thm}[\cite{Cou65}]
Let $\L : C^2_c(\R^n)\to C(\R^n)$ be any bounded, linear, and translation invariant operator that satisfies the Global Minimum Principle \eqref{eq:GMP}.
Then $\L$ is of the form 
\begin{equation}
\label{eq:levytypeoperator}
 \begin{split}
 \L f(x) & = - \frac12 {\rm div}\big(A\nabla f(x)\big) - b\cdot \nabla f(x) - c f(x)\\
 &  \quad - \int_{\R^n} \big(f(x+y) - f(x) - \nabla f(x) \cdot y \chi_{B_1}(y)\big)\nu (dy)
 \end{split}
 \end{equation}
 for a symmetric and nonnegative matrix $A$, $b\in \R^n$, $c\le 0$, and a nonnegative $\nu$ measure satisfying
 \[
\nu(\{0\}) = 0\quad\text{and}\quad \int_{\R^n}\min\{1, |x|^2\}\nu(dx) < \infty.
 \]
\end{thm}

Notice that if we additionally have  $\L 1 = 0$ then it follows that $c=0$.
Thus, this gives an alternative way to characterize the infinitesimal generators of L\'evy processes.

\subsection{Integro-differential equations}
\label{subsect-2.1.4}

Infinitesimal generators are very useful to study analytic properties of the underlying stochastic process. We next list some examples, and refer to \cite{Eva13, Bertoin} for more details.
 
\subsubsection{Expected payoff} Let $(X_t)_{t\ge 0}$ be a L\'evy process with infinitesimal generator $\L$, and let $\Omega\subset \R^n$ be an open domain. Suppose that we have a payoff function $g:\R^n\setminus \Omega\to \R$, and given any $x\in \Omega$, we consider the process $X_t^x := x+ X_t$: that is, we initialize a particle at $x$ that is moving according to $X_t$. 

\begin{figure}
\centering
\makebox[\textwidth][c]{\includegraphics[scale = 1.0]{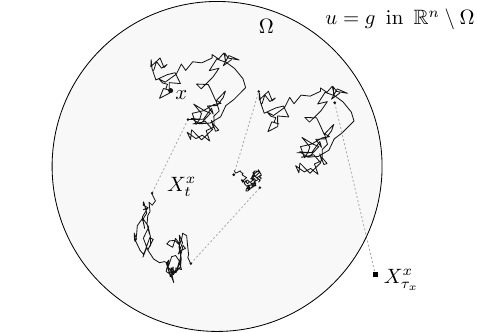}}
\caption{\label{fig:random} Graphical representation of a L\'evy flight inside $\Omega$ and the first time $\tau_x$ where it falls outside. Since the paths are not continuous, we position of $X_{\tau_x}^x$ can be anywhere in $\R^n\setminus \Omega$.}
\end{figure}

If we denote by $\tau_x$ the first time $x+X_t\notin \Omega$ (that is, the first time the particle falls outside of the domain\footnote{Recall from the definition of L\'evy process, Definition~\ref{defi:Levy}, that we do not necessarily have continuous paths, so in the first hitting time it is not necessarily true that $x+X_{\tau_x}\in \partial \Omega$. Still, we are working with c\`adl\`ag modifications of the underlying processes, and hence $\tau_x$ is a stopping time.}), we get paid $g(x+X_{\tau_x})$ (see Figure~\ref{fig:random} for a graphical representation of this setting). 
We want to compute the expected payoff
\[
u(x) := \E[g(x+X_{\tau_x})],
\]
a deterministic function that depends only on the starting point $x$. 
It turns out that such function $u$ solves the Dirichlet problem for the operator $\L$ with exterior datum $g$,
\[
\left\{
\begin{array}{rcll}
\L u & = & 0& \quad\text{in}\quad \Omega\\
u & = & g& \quad\text{in}\quad \R^n\setminus \Omega.
\end{array}
\right.
\]
Indeed, in this case we have that, using \eqref{eq:Ptdef}-\eqref{eq:heat_eq_p},
\[
u(x) = u(x) + \E\left[-\int_0^{\tau_x} \L u(x+X_t)\, dt\right] = \E(g(x+X_{\tau_x})). 
\]
By taking $g(x) = \chi_A(x)$ for some $A\subset \R^n$ this can be used, for example, to compute the probability that the first time the particle falls outside of $\Omega$ it does so in $A$.

When the stochastic process is the Brownian motion (and $\L = -\frac12\Delta$), then $(X_t)_{t\ge 0}$ has continuous paths, $g$ needs to be defined only on $\partial \Omega$, and the previous equation is the Laplace equation with Dirichlet datum $g:\partial \Omega \to \R$.

\subsubsection{Expected hitting time} Under the previous assumptions, one may ask instead what is the first time that the particle starting at $x\in \Omega$ will fall outside of $\Omega$. 
That is, to find the expected time
\[
u(x) = \E[\tau_x]. 
\]
In this case, $u$ solves the Dirichlet problem
\[
\left\{
\begin{array}{rcll}
\L u & = & 1& \quad\text{in}\quad \Omega\\
u & = & 0& \quad\text{in}\quad \R^n\setminus \Omega.
\end{array}
\right.
\]
When putting a general right-hand side $f:\Omega\to \R$ instead, we obtain the average value of $f(x+X_t)$ along a path before the particle falls outside of the domain.

\subsubsection{Surviving particles} Under the previous assumptions, assume moreover that at each time $t\ge 0$ when $x+ X_t\in \Omega$, there is a certain probability such that the particle is being killed (and gives zero payoff). 
That is, there exists some function $c:\Omega\to [0, \infty)$ such that a particle that is at position $\bar x$ for a time $\tau > 0$ is killed with probability proportional to $c(\bar x) \tau$.   
In this case, the average value of $f:\Omega\to \R$ along a path $x+X_t$ before the particle either falls outside of the domain, or is killed (with probability proportional to $c$) is given by the solution to 
\[
\left\{
\begin{array}{rcll}
\L u +cu & = & f& \quad\text{in}\quad \Omega\\
u & = & 0& \quad\text{in}\quad \R^n\setminus \Omega.
\end{array}
\right.
\]

\subsubsection{Nonlinear problems} There are also several  nonlinear problems for L\'evy processes that can be analytically tackled by means of the corresponding infinitesimal generators. 
We refer to Chapter~\ref{ch:fully_nonlinear} for examples on controlled diffusion processes (which yield convex fully nonlinear equations) or on stochastic differential games (yielding general fully nonlinear equations); and to Chapter~\ref{ch:obst_pb} for a description of the optimal stopping problem (which yields the obstacle problem).

\subsection{Stable processes} An important class of L\'evy processes are those given by stable distributions:
\index{Stable process}
\begin{defi}[Stable distribution] \index{Stable distribution}
\label{defi:stable} 
A random variable $X$ on $\R^n$ is \emph{stable} if for any $a, b>0$, there exists $c > 0$ and $d\in \R^n$ such that 
\[
a X_1 + b X_2 \overset{d}{=} c X + d.
\]
Here $X_1$ and $X_2$ are independent copies of $X$. Equivalently, $X$ is stable if for all $m \ge 2$ there exist $c_m \ge 0$ and $d_m \in \R^n$ such that
\begin{equation}
\label{eq:sum}
X_1 + \dots + X_m \overset{d}{=} c_m X + d_m,
\end{equation}
where $X_1,\dots,X_m$ are $m$ independent copies of $X$. We say that $X$ is \emph{strictly stable} if $d_m = 0$. 
\end{defi}

We will say that a L\'evy process $(X_t)_{t\ge 0}$ is \emph{(strictly) stable} if $X_1$ is a (strictly) stable random variable. 
Observe that  any stable distribution is infinitely divisible, but the converse is \emph{not} true (infinitely divisible distributions are not necessarily stable).

Equivalently, a random variable variable $X$ is stable if it has a \emph{domain of attraction}, that is, there exist a sequence of independent identically distributed (i.i.d.) random variables $Z_1, Z_2, \dots$ such that 
\[
\frac{Z_1+Z_2+\dots+Z_n}{C_n} + D_n \ \overset{d}{\longrightarrow} \ X,
\]
for some $C_n > 0$ and where $\overset{d}{\rightarrow}$ indicates convergence in distribution. Hence, stable distributions appear as limits of sums of i.i.d. random variables (after an appropriate rescaling). If such variables have finite variance, the limiting stable distribution  is a Gaussian distribution by the Central Limit Theorem. 
If the variances are not finite, however, we obtain a larger class of stable distributions (and a generalized Central Limit Theorem, for random variables with infinite variance). 

In Definition~\ref{defi:stable}, the value of the constant $c_m$ is explicit. Indeed, by taking $X_1 - X_2$ we can assume the random variables to be symmetric (and hence $d_m = 0$). If we let $m_1, m_2\in \N$, then we can separate $X_1, \dots, X_{m_1 m_2}$ into $m_2$ groups of $m_1$ variables each. Then, by the defining property, we have
\[
c_{m_1m_2} X \overset{d}{=} X_1+\dots X_{m_1 m_2} \overset{d}{=} c_{m_1} c_{m_2} X.
\]
Therefore $c_{m_1m_2} = c_{m_1} c_{m_2}$. From here, it is now easy to deduce $c_m = m^{\frac{1}{\alpha}}$ for some $\alpha > 0$, and from the fact that the normal distribution is stable with $\alpha  =2$, one can deduce $\alpha\in (0, 2]$ (see, for example, \cite[Chapter VI]{Fel71}).

For each stable distribution, hence, there is a unique $\alpha\in (0, 2]$ such that \eqref{eq:sum} holds with $c_m = m^{\frac{1}{\alpha}}$. Such number is called the \emph{characteristic exponent} or the \emph{index of the distribution}, and any random variable $X$ that is stable with index $\alpha$ is called $\alpha$-\emph{stable}. Similarly, from now on, a L\'evy process $(X_t)_{t\ge 0}$ is $\alpha$-\emph{stable} if $X_1$ is $\alpha$-\emph{stable} (respectively, we call it strictly $\alpha$-stable).

More precisely, the characteristic exponent $\Psi$ of a strictly $\alpha$-stable distribution, \eqref{eq:Psi}, is explicit and of the form  (see \cite[Theorems 2.3.1-2.4.1]{ST94})
\[
\Psi(\xi) = \left\{
\begin{array}{ll} \displaystyle \small
-\int_{\S^{n-1}} |\xi \cdot \theta|^\alpha\left(1-i\tan\frac{\pi\alpha}{2}\, {\rm sign}(\xi\cdot \theta)\right) \, \zeta(d\theta)&  ~\text{if}~\alpha \neq 1\\  \displaystyle
-\int_{\S^{n-1}} \left(|\xi \cdot \theta|+i\frac{2}{\pi}\log|\xi\cdot \theta| \, (\xi\cdot \theta)\right) \, \zeta(d\theta)-i\mu_0\cdot \xi& ~\text{if}~\alpha = 1,
\end{array}
\right.
\]
where $0 < \alpha < 2$, ${\rm sign}(t) = \pm 1$ for $\pm t > 0$, and ${\rm sign}(0) = 0$, and for some $\mu_0\in \R^n$, and some finite measure $\zeta$ on $\S^{n-1}$. The measure $\zeta$ is called the \emph{spectral measure} of the $\alpha$-stable distribution or random variable. 

In the case of symmetric distributions, we have a simpler representation,
\[
\Psi(\xi) = -\int_{\S^{n-1}} |\xi\cdot\theta|^\alpha\, \zeta(d\theta)
\]
for $0< \alpha < 2$, 
where $\zeta$ is the (symmetric) spectral measure of the corresponding symmetric $\alpha$-stable distribution.

\begin{example}
As we have seen, 2-stable distributions are Gaussian distributions. On the other hand, an example of a 1-stable distribution that is also rotation invariant is the Cauchy distribution (whose corresponding infinitesimal generator is the square root of the Laplacian, $\sqrt{-\Delta}$).
\end{example} 

At the level of the corresponding strictly $\alpha$-stable L\'evy process, we have that they are self-similar, that is, 
\[
X_t \sim {t^\frac{1}{\alpha}} X_1\quad\text{for all}\quad t > 0. 
\]

In terms of the generating triplet of a symmetric process, $(A, b, \nu)$, we have that if $\alpha = 2$ then $b = 0$ and $\nu = 0$, while if $0 < \alpha < 2$, then $A = 0$, $b = 0$, and the symmetric L\'evy measure $\nu$  is scale invariant
\[
\nu = \kappa^{-\alpha}S_\kappa \nu\quad\text{for all}\quad \kappa > 0
\]
where $S_\kappa \nu(B) = \nu(\kappa^{-1} B)$  for all Borel set $B\subset \R^n$. Alternatively, we can write, for $s = \frac{\alpha}{2}\in (0, 1)$, and where we will   denote $K = \nu$, 
\begin{equation}
\label{eq:Lu_stab0}
\nu(B) = K(B) = c_s \int_{\S^{n-1}} \int_0^\infty \chi_B(r\theta) \frac{dr}{r^{1+2s}}\zeta(d\theta)\quad\text{for all}\quad B\subset \R^n~~\text{Borel},
\end{equation}
for some constant $c_s>0$ depending only on $s$.

Overall, the infinitesimal generator of a symmetric $\alpha$-stable L\'evy process is of the form \index{Stable L\'evy process@$\alpha$-stable L\'evy process}
\begin{equation}
\label{eq:Lu_stab}\index{Stable operator}
\L u(x) = \frac{c_s}{2} \int_{\S^{n-1}}{\rm P.V.}\int_{-\infty}^{\infty} \big(u(x) - u(x+r\theta)\big)\frac{dr}{|r|^{1+2s}} \zeta(d\theta),
\end{equation}
where $\zeta$ is a finite measure on the sphere (the \emph{spectral measure}\index{Spectral measure}). Observe that such operators are scale invariant, that is, if $u_r(x) = u(rx)$, then 
\begin{equation}
\label{eq:scale_invariant_stable}\index{Scale invariance!Stable operators}
(\L u_r)(x) = r^{2s}(\L u)(rx).
\end{equation}

We will also ask for ellipticity conditions, to ensure boundedness and nondegeneracy of the operator. 
On the one hand we consider an upper ellipticity constant $\tilde \Lambda$ such that
\begin{equation}
\label{eq:Stab_L_ellipt}\index{Ellipticity!Stable operators}
\int_{\S^{n-1}}\zeta(d\theta) \le \tilde \Lambda. 
\end{equation}

On the other hand, we also require a nondegeneracy-type condition\footnote{In the local case, when $\alpha = 2$ this is typically required in the form $A \ge \tilde \lambda\, {\rm Id}$ for some $\tilde \lambda > 0$, or alternatively, 
\[
\inf_{\be \in \S^{n-1}} \be \cdot A \be \ge \tilde \lambda > 0.
\]
This ensures that the matrix $A$ does not degenerate in any direction. 
}
 to ensure that the corresponding measure $\zeta$ is full dimensional, that is, it is not concentrated on any hyperplane. 
A way to control it is to ask that, for some $\lambda > 0$, 
\begin{equation}
\label{eq:Stab_l_ellipt}
 \inf_{\be\in \S^{n-1}}\int_{\S^{n-1}}|\be\cdot \theta|^{2s}\zeta(d\theta) \ge \tilde \lambda > 0. 
\end{equation}

Operators of the form \eqref{eq:Lu_stab} satisfying the ellipticity conditions \eqref{eq:Stab_L_ellipt}-\eqref{eq:Stab_l_ellipt} are called \emph{stable operators} of order $2s$,  with ellipticity constants $\tilde \lambda$ and~$\tilde \Lambda$.

\begin{rem}
\label{rem:equiv_ellipt_homog}
Equation \eqref{eq:Stab_l_ellipt} is not the only way to ensure that the spectral measure is not concentrated on a hyperplane. One could instead impose $ \inf_{\be\in \S^{n-1}}\int_{\S^{n-1}}|\be\cdot \theta|^{\beta}\zeta(d\theta) \ge \lambda_\beta > 0$ for any $\beta >0$. Thanks to H\"older's inequality (and \eqref{eq:Stab_L_ellipt}), however, both conditions are equivalent (up to changing the value of $\lambda_\beta$). We have chosen $\beta = 2s$ because then the expression \eqref{eq:Stab_l_ellipt} has a meaning in terms of the Fourier symbol: it says that the Fourier symbol $\A(\xi) = -\Psi(\xi)\ge \tilde \lambda |\xi|^{2s}> 0$.
That is, it is bounded below by the Fourier symbol of the fractional Laplacian. 
Notice also that \eqref{eq:Stab_L_ellipt} is equivalent to $\mathcal A(\xi) \leq \tilde \Lambda |\xi|^{2s}$.
\end{rem}

\subsubsection{The fractional Laplacian} The first example of stable operator is, of course, the case $\alpha = 2$, that is, $\L = -\frac12 \Delta$. In this case, the corresponding stable distribution is the normal distribution. 

The most famous example of stable operator that is nonlocal is the square root of the Laplacian $\sqrt{-\Delta}$, or more generally, the fractional Laplacian. 
It corresponds to the case in which the measure is not only symmetric, but also rotationally invariant\footnote{A measure $\mu$ is rotationally invariant if $\mu(B) = \mu(O B)$ for any orthogonal transformation $O\in \mathcal{O}(n)$ and for any Borel set $B\subset \R^n$.}, so that the infinitesimal generator is of the form
\[
\L u(x) = c_{n,s} {\rm P.V.} \int_{\R^n}\big((u(x) - u(x+y)\big)\frac{dy}{|y|^{n+2s}}
\]
for some $s = \frac{\alpha}{2} \in (0, 1)$. 
The name is justified by the corresponding Fourier multiplier $-\Psi$, that in this case is
\[
-\Psi(\xi)  = |\xi|^{2s}. 
\]

\begin{example}
\label{ex:four_symb_ex}
Not all stable spectral measures satisfying the above mentioned ellipticity conditions need to be as nice as the fractional Laplacian, or even,   be absolutely continuous measures. For example, the L\'evy process $X_t = (X_t^1,\dots, X_t^n)$ with $X_t^i$ independent one-dimensional symmetric stable processes is a symmetric stable L\'evy process, with infinitesimal generator given by 
\[
\L u = (-\partial_{x_1x_1})^s u + \dots + (-\partial_{x_nx_n})^s u,
\]
and spectral measure $\zeta(d\theta) = \delta_{\be_1}(d\theta) +\dots \delta_{\be_n}(d\theta)+\delta_{-\be_1}(d\theta) +\dots \delta_{-\be_n}(d\theta)$, up to a multiplicative constant, where $\be_i$ are the coordinate vectors in $\R^n$. 
Its Fourier symbol is $\mathcal A(\xi)=|\xi_1|^{2s}+...+|\xi_n|^{2s}$, and thus it is still comparable to the one of the fractional Laplacian.
\end{example}

\subsection{Kernels comparable to the fractional Laplacian}\label{ssec:comparable} \index{Comparable to fractional Laplacian}\index{Stable-like operators}

The fractional Laplacian is by far the most studied non-Gaussian stable operator. 
As such, it serves as a motivation to justify the following more general classes of L\'evy processes, that will be of importance throughout this book.

A very typical and natural assumption in the study of integro-differential operators is to consider infinitesimal generators with L\'evy measures (or kernels) \emph{comparable} to the one of the fractional Laplacian, see \eqref{eq:Kcompfls}.
These are  sometimes called \emph{stable-like operators} (see \cite{BL02, BL02b, CH03, Sil06, Bas09, Kas09}), and correspond to the  $\mathcal L_0$ class of Caffarelli and Silvestre \cite{CS,CS2,CS3}. 

Notice that, under this assumption, we drop the homogeneity property that was introduced in stable operators (since they are self-similar), but we now require the kernel to always be positive and uniformly bounded (in particular, it must be absolutely continuous with respect to the Lebesgue measure). 
As we will see, this extra assumption will allow us to prove some results (like the Harnack inequality) which do not hold for general stable operators.

In this case,  the L\'evy measure $K$ in \eqref{eq:L_def_lin} (also denoted $\nu$) can be expressed as an absolutely continuous kernel 
\begin{equation}
\label{eq:abs_cont_K}
\nu(dy) = K(dy) =  K(y) dy,
\end{equation}
with the stronger ellipticity conditions
\begin{equation}
\label{eq:Kcompfls}\index{Ellipticity!Operators comparable to fractional Laplacian}
0< \frac{\lambda}{|y|^{n+2s}}\le K(y) \le \frac{\Lambda}{|y|^{n+2s}} \qquad\text{for all}\quad y \in \R^n,
\end{equation}
with $0<\lambda\le \Lambda$.
The second inequality ensures boundedness of the operator, while the first one ensures its nondegeneracy.

\begin{rem}
Stable operators and operators comparable to the fractional Laplacian  both arise as infinitesimal generators of L\'evy processes, but they encompass different subsets. 
That is, they are not included one into the other: stable-like operators, contrary to stable operators, do not need to be homogeneous; whereas stable operators, contrary to stable-like operators, do not need to be absolutely continuous with respect to the Lebesgue measure. 
\end{rem}

We emphasize that by dropping the homogeneity assumption from stable operators, operators satisfying \eqref{eq:Kcompfls} are not scale invariant anymore. 
Still, the whole class of operators comparable to the fractional Laplacian with given ellipticity constants is scale invariant. 
That is, if $\L$ satisfies  \eqref{eq:Kcompfls} with ellipticity constants $\lambda$ and $\Lambda$, and if we denote  $u_r(x) = u(rx)$, then $(\L_r u_r)(x) = r^{2s}(\L u)(rx)$, where $\L_r$ also satisfies \eqref{eq:Kcompfls} with the same ellipticity constants, $\lambda$ and $\Lambda$.

\subsection{General elliptic operators of order $2s$} 
\label{ssec:general_ellipt}
\index{General operators}

A more general (and very natural) class of integro-differential operators is obtained by, instead of comparing them with the fractional Laplacian in the physical space, comparing them to the fractional Laplacian in the Fourier space. 

Given $s\in (0, 1)$ we consider symmetric operators $\L$ of the form 
\begin{equation}
\label{eq:Lu_nu}
\begin{split}
 \L u(x) &=   {\rm P.V.}\int_{\R^n} \big(u(x) - u(x+y)\big)K (dy) \\
 & =\frac12 \int_{\R^n} \big(2u(x) - u(x+y) - u(x-y)\big)K (dy),
 \end{split}
\end{equation}
where
\begin{equation}
\label{eq:nu_cond}
\text{$K\ge 0\,$ is symmetric, \quad $K(\{0\}) = 0,$}\quad \int_{\R^n} \min\{1, |y|^2\} \,K(dy) <\infty. 
\end{equation}
(cf. \eqref{eq:L_def_lin} and recall Definition~\ref{defi:sym_meas}) with Fourier symbol $\A$ (see \eqref{eq:fourier_symb}) such that, for some ellipticity constants $0 < \tilde \lambda \le \tilde \Lambda$, we have
\begin{equation}
\label{eq:GSLO}
0 < \tilde \lambda |\xi|^{2s} \le \A(\xi) \le \tilde \Lambda |\xi|^{2s}\quad\text{for all}\quad \xi\in \R^n. 
\end{equation}
We will call this class $\G_s$, which contains both stable operators (recall Remark \ref{rem:equiv_ellipt_homog}) and operators comparable to the fractional Laplacian (recall \eqref{eq:fourier_symb}). 
Observe that, as before, the operator $\L$ may not be scale invariant, but the whole class of kernels satisfying \eqref{eq:GSLO} is scale invariant\footnote{This can be seen either through the properties of the Fourier transform, or by means of the equivalent ellipticity conditions \eqref{eq:Kellipt_gen_L}-\eqref{eq:Kellipt_gen_l}.}. 

Typically, we want to work with the operator $\L$ via its L\'evy measure (and not via its Fourier symbol), and so we will need an equivalent formulation in terms only of the measure $K$ (often denoted $\nu$ before). 
We will prove that, given any L\'evy operator~$\L$, the ellipticity condition \eqref{eq:GSLO} is equivalent to the following: there are constants $\lambda, \Lambda>0$ (depending only on $n$, $s$, $\tilde \lambda$, and~$\tilde \Lambda$), such that
\begin{equation}
\label{eq:Kellipt_gen_L} \index{Ellipticity!General operators}
r^{2s} \int_{B_{2r}\setminus B_r} K(dy) \le \Lambda\qquad\text{for all}\quad r  >0,
\end{equation}
and 
\begin{equation}
\label{eq:Kellipt_gen_l}
0 < \lambda \le r^{2s-2}\inf_{e\in \S^{n-1}}\int_{B_{r}}|e\cdot y|^{2} K(dy)\qquad\text{for all}\quad r  >0.
\end{equation}
We refer to Proposition~\ref{prop:equiv_fourier} below for a proof of this fact.

Most of the results we will prove in this chapter hold for the following general class of operators.

\begin{defi}[General elliptic operator]
\label{defi:G}
Let $s\in (0, 1)$, and $\lambda, \Lambda > 0$. 
We define
\[
\GL :=\left\{\L : \begin{array}{l}
\text{$\L$ is an operator of the form \eqref{eq:Lu_nu}-\eqref{eq:nu_cond}}\\
\text{such that \eqref{eq:Kellipt_gen_L}-\eqref{eq:Kellipt_gen_l} hold}
\end{array}
 \right\}.
\]
Equivalently, $\L\in \GL$ for some $\lambda$ and $\Lambda$ if and only if it is of the form \eqref{eq:Lu_nu}-\eqref{eq:nu_cond} and \eqref{eq:GSLO} holds for some $\tilde\lambda$ and $\tilde \Lambda$ (see Proposition~\ref{prop:equiv_fourier}). 
\end{defi}

\begin{rem}[Scale invariance]
\label{rem:Scale_invariance}\index{Scale invariance!General operators}
The whole class $\GL$ is scale invariant. 
That is, if $\L\in \GL$,   and we denote $u_r(x) = u(rx)$ with $r > 0$, then $(\L_r u_r)(x) = r^{2s}(\L u)(rx)$, for some $\L_r\in \GL$.  

More precisely, if $\L$ has kernel  $K(dy)$, then $\L_r$ has kernel $K_r(dy) = r^{2s} K(r\, dy)$. Here, and throughout the work,  we denote by $K(x+\rho \, dy)$ for $\rho\in \R\setminus\{0\}$ and $x\in \R^n$, the measure $\tilde K(dy) = K(x+\rho\, dy)$ such that $\tilde K(B) = K(x+\rho B)$ for any $B\subset \R^n$ Borel (see also \nameref{notation} on page \pageref{not:measures}). 
\end{rem}

\begin{rem}
Observe that the class $\GL$ contains both stable operators and operators comparable with the fractional Laplacian. 
In particular, if $\L \in \GL$ with $K$ homogeneous, then $\L$ is a stable operator and satisfies \eqref{eq:Stab_L_ellipt}-\eqref{eq:Stab_l_ellipt} for some ellipticity constants $\tilde \lambda$ and $\tilde \Lambda$.
\end{rem}

In some results, we will consider those operators in $\GL$ whose L\'evy measure is homogeneous:

\begin{defi}[General stable operator]
\label{defi:Gh}
Let $s\in (0, 1)$, and $\lambda, \Lambda > 0$. 
We define
\[
\GLh :=\big\{\L \in \GL : \text{its L\'evy measure, $K$, is homogeneous} \big\}.
\]
\end{defi}

This is the class of \emph{stable operators}, and they can equivalently be written as \eqref{eq:Lu_stab} (with a homogeneous L\'evy measure, of the form \eqref{eq:Lu_stab0}).  In this case, one can express the ellipticity conditions directly in terms of the spectral measure, \eqref{eq:Stab_L_ellipt}-\eqref{eq:Stab_l_ellipt} (recall Remark~\ref{rem:equiv_ellipt_homog}).

For notational convenience only, we may sometimes write the measure $K$ as if it was absolutely continuous with respect to the Lebesgue measure, 
\begin{equation}
\label{eq:Kabscont}
K(dy) = K(y)\, dy
\end{equation}
where the kernel is such that $K\in L^1_{\rm loc}(\R^n\setminus \{0\})$ and
\begin{equation}
\label{eq:Kint}
K\ge 0\qquad\text{and}\qquad \int_{\R^n}\min\{1, |y|^2\}K(y)\, dy < \infty. 
\end{equation}
By doing so, integrals against $r^n K(ry)\, dy$ or $K(x-y)\, dy$ become less notationally heavy\footnote{Observe that, in this case, given $\rho\in \R\setminus\{0\}$, we have $K(\rho\, dy) = |\rho|^n K(\rho y)\, dy$.}.
With this notation, the symmetry of the measure $K$ (recall Definition~\ref{defi:sym_meas}) is
\begin{equation}
\label{eq:Ksym}
K(y) = K(-y)\qquad\text{for all}\quad y \in \R^n. 
\end{equation}
Observe that, when $\L = \fls$, we have $K(y)dy  = c_{n,s}|y|^{-n-2s} dy$. 
Observe also that, in this notation,  stable operators  have kernels satisfying
\begin{equation}
\label{eq:Khom}
\text{$K(y)$ is homogeneous of degree $-n-2s$ for some $s\in (0, 1)$,}
\end{equation}
or equivalently, 
\begin{equation}
\label{eq:Khom2}
K(y) = \frac{K\left(y/{|y|}\right)}{|y|^{n+2s}}\qquad\text{for some}\quad K\big|_{\mathbb{S}^{n-1}}\in L^1(\mathbb{S}^{n-1}),~~K \ge 0,
\end{equation}
where $K$ is the spectral measure (denoted $\zeta$ in \eqref{eq:Lu_stab}, up to a multiplicative constant, which we are assuming to be absolutely continuous).

\subsection{Regular kernels} \index{General operators!Regular kernels}
Some of the results we will prove  require kernels that are somewhat regular (in an integral sense). 
The exact condition we will need is that, for $\alpha\in (0, 1]$,  the following quantity is finite
\begin{equation}\label{Calpha-assumption}
[K]_\alpha := \sup_{\rho > 0}  \sup_{x, x'\in B_{\rho/2}}\rho^{2s+\alpha}\int_{B_{2\rho}\setminus B_\rho}\frac{\big|K(x-dy)-K(x'-dy)\big|}{|x-x'|^{\alpha}},
\end{equation}
where  we recall that the measure $K(x-dy) = K^x(dy)$ is such that $K^x(B) = K(x-B)$ for all Borel sets $B\subset \R^n$. More generally, we define for any $\mu > 0$ with $\lceil \mu-1\rceil = m$ and for $K$ such that its distributional derivative $D^m K$ is a locally finite signed measure in $\R^n\setminus\{0\}$ (and $K$ is absolutely continuous for $m\ge 1$),
\begin{equation}\label{Cmu-assumption}
[K]_\mu := \sup_{\rho > 0}  \sup_{x, x'\in B_{\rho/2}}\rho^{2s+\mu}\int_{B_{2\rho}\setminus B_\rho}\frac{\big|D^mK(x-dy)-D^mK(x'-dy)\big|}{|x-x'|^{\mu-m}}.
\end{equation}
 If $\L$ is an operator of the form \eqref{eq:Lu_nu} with kernel $K$, we denote
\begin{equation}
\label{eq:LKmu}
[\L]_\mu := [K]_\mu. 
\end{equation}
\begin{defi}[Regular general elliptic operators]
\label{defi:Gmu}
Given $\mu > 0$, we define  the class
\[
\G_s(\lambda, \Lambda; \mu) :=\big\{\L \in \GL : [\L]_\mu < \infty\big\},
\]
that is, those operators in $\GL$ with regular kernels, in the sense of \eqref{Calpha-assumption}-\eqref{Cmu-assumption}-\eqref{eq:LKmu}. 
\end{defi}

Notice that the boundedness of \eqref{Cmu-assumption} is equivalent to (cf. \eqref{eq:prop_kernel} below) 
\begin{equation}
\label{eq:prop_kernel_Calpha}
\int_{B_\rho^c}\big|D^mK(x-dy)-D^m K(x'-dy)\big|  \le C \rho^{-2s-\mu} |x-x'|^{\mu - m}
\end{equation}
for all $x,x'\in B_{\rho/2}$ and all $\rho > 0$, and for some constant $C = C'[K]_\mu$.

It is important to emphasize that this property is preserved by scaling. That is, if we define $u_r(x) = u(rx)$, given $\L$ with kernel $K$ we can consider the operator $\L_r$ with kernel $K_r$ such that $(\L_r  u_r)(x) = r^{2s}(\L u)(rx)$. We then have that if $\L\in \G_s(\lambda, \Lambda; \mu)$ then $\L_r\in \G_s(\lambda, \Lambda; \mu)$ with $[K]_\mu = [K_r]_\mu$.

Furthermore, the boundedness of \eqref{Cmu-assumption} (or \eqref{eq:prop_kernel_Calpha}) is of course implied by the (stronger) assumption
\begin{equation}
\label{eq:easycond_K_smooth}
[K]_{C^\mu(B_\rho^c)} \le C_\circ \rho^{-n-2s-\mu}.
\end{equation}
Notice also that if 
\[[K]_{{\rm BV}(B_{2\rho}\setminus B_\rho)} \leq C\rho^{-2s-1}\]
then \eqref{Calpha-assumption} holds for all $\alpha\in(0,1]$ (where ${\rm BV}$ denotes the \emph{bounded variation} norm).
In particular, kernels satisfying \eqref{Calpha-assumption} (even for $\alpha=1$) may be discontinuous and unbounded.

Finally, let us mention that the quantity in \eqref{Calpha-assumption} (and \eqref{Cmu-assumption}) is essentially a Besov seminorm,
\[
[f]_{B_{p, \infty}^\alpha} = \sup_{h} \frac{\|f(\cdot +h) - f\|_{L^p}}{|h|^\alpha},
\]
that is, there is an equivalence of seminorms for a given kernel $K$:
\[
[K]_\alpha\asymp  \sup_{\rho > 0}\, \rho^{2s+\alpha}[K]_{B_{1,\infty}^\alpha(B_{2\rho\setminus B_\rho})}. 
\]
By classical embeddings for Besov spaces, we have $W^{\alpha, 1}\subset B_{1,\infty}^\alpha$, so it suffices to bound the seminorms in $W^{\alpha, 1}$ (see \cite[Chapter V]{Stein} or \cite{Tri92}).

\begin{rem}
The fractional Laplacian $\fls$ has a kernel of the form $K(y) = c_{n,s}|y|^{-n-2s}$ which, by homogeneity, satisfies \eqref{eq:easycond_K_smooth} for any $\mu > 0$. Hence $\fls \in \G_s(\lambda, \Lambda; \mu)$ for all $\mu > 0$.   More generally, any operator $\L_{K_\zeta}$ with kernel $K_\zeta$ of the form $K_\zeta(y) = \zeta(y/|y|)|y|^{-n-2s}$ for some $\zeta\in C^\infty(\mathbb{S}^{n-1})$   symmetric, nonnegative, and nonzero (namely, smooth homogeneous kernels), satisfies $\L_{K_\zeta}\in \G_s(\lambda, \Lambda; \mu)$ for all $\mu > 0$, and for some $\lambda, \Lambda > 0$ depending only on $\zeta$.  
\end{rem}

\begin{rem}
When dealing with operators with homogeneous kernels $\GLh$, as in \eqref{eq:Khom}-\eqref{eq:Khom2}, it is enough to check the condition for   $\rho = 1$:
\[
[K]_\alpha :=   \sup_{x, x'\in B_{1/2}}\int_{B_{2}\setminus B_1}\frac{\big|K(x-y)-K(x'-y)\big|}{|x-x'|^{\alpha}} dy\asymp  [K]_{B_{1,\infty}^\alpha(B_{2}\setminus B_1)}.
\]
In particular, in this case it suffices to have  $K\in W^{\alpha, 1}(\partial B_1)$ or $K \in C^\alpha(\partial B_1)$.

Here, $W^{\alpha,p}(\Omega)$ or $W^{\alpha,p}(\partial\Omega)$ denote the spaces of functions in $L^p(\Omega)$ or $L^p(\partial\Omega)$ for which the seminorm
\[
[w]_{W^{\alpha,p}(\Omega)} :=\int_{\Omega}\int_{\Omega}\frac{|w(x)-w(y)|^p}{|x-y|^{n+\alpha p}}\,dx\,dy
\]
or
\begin{equation}\label{Walphap}\index{Fractional Sobolev space}
[w]_{W^{\alpha,p}(\partial\Omega)} :=\int_{\partial\Omega}\int_{\partial\Omega}\frac{|w(x)-w(y)|^p}{|x-y|^{n-1+\alpha p}}\,dx\,dy
\end{equation}
is finite.
\end{rem}

\subsection{Brief discussion on the notation}

Throughout this chapter, we consider L\'evy operators of the type \eqref{eq:L_def_lin}, where $K$ denotes the (nonnegative) \emph{L\'evy measure} that in subsection~\ref{ssec:Levy} we have introduced as $\nu$, and satisfies \eqref{eq:cond_nu}. 
The Fourier symbol of these L\'evy operators is given by~\eqref{eq:fourier_symb}.

The notation $\nu$ is typically used in the probability literature, whereas in the PDE community people usually study the case of absolutely continuous kernels, which are commonly denoted $K(y)\, dy$.
In this PDE book, however, we want to allow for general measures as well, so we decided to denote these measures $K(dy)$.
 
Among the reasons to include general measures (and not only absolutely continuous ones), we have:
\begin{itemize}[leftmargin=*]
\item It is the most general class of operators of order $2s$ for which one can develop a regularity theory.
 
\item It includes relevant operators such as $(-\partial^2_{x_1x_1})^{2s}+...+(-\partial^2_{x_nx_n})^{2s}$, whose kernel is not absolutely continuous (see Example~\ref{ex:four_symb_ex}). 

\item The class of operators is then closed under weak limits (see Proposition~\ref{prop:stab_distr}). 
\end{itemize}

In order to use the usual convenient notation for absolutely continuous kernels, we decided to keep it for measures as well (see \nameref{notation} on page \pageref{not:measures}).
This means that, given a measure $\mu(dy)$ (for example,  $K(dy)$ as above), we have defined the associated measures $\mu(x_\circ +r \, dx)$ for $x_\circ\in  \R^n$ and $r\in \R\setminus\{0\}$.

In particular, when $K$ is absolutely continuous, $K(dy) = K(y) \, dy$, we have the equivalences:
\[
K(-dy) = K(-y)\, dy,\qquad K(x_\circ+dy) = K(x_\circ+y)\, dy,
\]
and more generally, 
\[
K(x_\circ + r\, dy) = |r|^n K(x_\circ + r\,y)\, dy \qquad\text{for}\quad x_\circ\in \R^n, \ r\in \R\setminus\{0\}. 
\]

For example, in the case of absolutely continuous kernels\footnote{Notice, though, that there are purely singular measures for which $[K]_\alpha$ is finite.} the seminorm \eqref{Calpha-assumption} is
\[
[K]_\alpha := \sup_{\rho > 0}  \sup_{x, x'\in B_{\rho/2}}\rho^{2s+\alpha}\int_{B_{2\rho}\setminus B_\rho}\frac{\big|K(x-y)-K(x'-y)\big|}{|x-x'|^{\alpha}}\, dy.
\]
On the other hand, notice that if $[K]_\mu <\infty$ for $\mu > 1$ (see \eqref{Cmu-assumption}),  then the fact that $D^m K$ is a locally finite measure implies that $K$ (and $D^{m-1} K$) is absolutely continuous, and in particular, $K(dy) = K(y)\, dy$ in this case.

\section{Basic properties and notions of solution}
\label{sec:basic_prop}

In this section we will prove some basic properties on the class of ope\-rators we consider, some of them analogous to those we have seen for the fractional Laplacian. Moreover, we also introduce the notions of strong, weak, and distributional solution for these linear elliptic integro-differential operators.

\subsection{Preliminaries on the kernels}

Let us start with a preliminary property on the class of kernels $\GL$ itself, given by Definition~\ref{defi:G}. We show here that the conditions \eqref{eq:Kellipt_gen_L}-\eqref{eq:Kellipt_gen_l}   are indeed equivalent to  \eqref{eq:GSLO} (i.e., the Fourier symbol is comparable to the Fourier symbol of the fractional Laplacian). 

In order to do that, we will use that  condition \eqref{eq:Kellipt_gen_L} implies (and up to a constant, is equivalent to)
\begin{equation}
\label{eq:prop_kernel}
\int_{\R^n\setminus B_r} K(dy) \le C \Lambda r^{-2s}\quad\text{for all}\quad r > 0,
\end{equation}
for some $C$ depending only on $s$, since
\[
\int_{\R^n\setminus B_r} K(dy)  = \sum_{k\ge 0} \int_{B_{2^{k+1}r} \setminus B_{2^k r}} K(dy)  \le r^{-2s}\sum_{k\ge 0} 2^{-2ks} = C  r^{-2s}.
\]
On the other hand, we also have
\begin{equation}
\label{eq:prop_kernel2}
\int_{B_r} |y|^{2s+\alpha}K(dy) \le C \Lambda r^\alpha\quad\text{for all}\quad r > 0,
\end{equation}
and 
\begin{equation}
\label{eq:prop_kernel3}
\int_{B_r^c} |y|^{2s-\alpha}K(dy) \le C \Lambda r^{-\alpha}\quad\text{for all}\quad r > 0,
\end{equation}
for any $\alpha > 0$ and some $C$ depending only on $\alpha$. Indeed, 
\[
\begin{split}
\int_{B_r} |y|^{2s+\alpha}K(dy) & \le \sum_{k\ge 0}\int_{B_{r/2^{k}}\setminus B_{r/2^{k+1}}}  (2^{-k}r)^{2s+\alpha}K(dy)\\ & \le \Lambda r^{\alpha} \sum_{k\ge 0} 2^{-k\alpha+2s} \le C \Lambda r^\alpha, 
\end{split}
\]
and 
\[
\begin{split}
\int_{B_r^c} |y|^{2s-\alpha}K(dy) & \le \sum_{k\ge 0}\int_{B_{2^{k+1}r}\setminus B_{2^{k}r}}  (2^{k+1}r)^{2s-\alpha}K(dy)\\ & \le \Lambda r^{-\alpha} \sum_{k\ge 0} 2^{2s-\alpha k -\alpha} \le C \Lambda r^{-\alpha}.
\end{split}
\]
Finally, by a similar reasoning we have that if \eqref{eq:prop_kernel_Calpha} holds, then 
\begin{equation}
\label{eq:prop_kernel4}
\int_{B_\rho^c}|y|^{2s+\mu-\alpha} \big|D^mK(x-dy)-D^m K(x'-dy)\big|  \le C [K]_\mu \rho^{-\alpha} |x-x'|^{\mu - m}
\end{equation}
for any $\alpha > 0$.
 
Properties \eqref{eq:prop_kernel}-\eqref{eq:prop_kernel2}-\eqref{eq:prop_kernel3} will appear recurrently throughout this chapter. Let us use them to show that \eqref{eq:Kellipt_gen_L}-\eqref{eq:Kellipt_gen_l} are equivalent to  \eqref{eq:GSLO}.

\begin{prop}
\label{prop:equiv_fourier}
Let $s\in (0, 1)$, and let $\L$ be an operator of the form \eqref{eq:Lu_nu}-\eqref{eq:Kint}-\eqref{eq:Ksym}. Then, 
\begin{enumerate}[leftmargin=*,label=(\roman*)]
\item \label{it:equiv_fourier_i} If $\L$ satisfies \eqref{eq:Kellipt_gen_L}-\eqref{eq:Kellipt_gen_l} for some $\lambda, \Lambda>0$, then there exist $\tilde\lambda, \tilde \Lambda>0$ depending only on $n$, $s$, $\lambda$, and $\Lambda$, such that \eqref{eq:GSLO} holds. 
\item \label{it:equiv_fourier_ii} If $\L$ satisfies \eqref{eq:GSLO} for some $\tilde\lambda, \tilde\Lambda>0$, then there exist $\lambda, \Lambda>0$ depending only on $n$, $s$, $\tilde\lambda$, and $\tilde\Lambda$, such that \eqref{eq:Kellipt_gen_L}-\eqref{eq:Kellipt_gen_l} holds.
\end{enumerate}
\end{prop}
\begin{rem}
In fact, for any $\alpha > 0$ fixed,  property \eqref{eq:prop_kernel2} (and also \eqref{eq:prop_kernel3}) is equivalent to \eqref{eq:prop_kernel}, and therefore to \eqref{eq:Kellipt_gen_L} (up to a multiplicative constant). Thus, choosing $\alpha = 2-2s>0$, we can equivalently express the ellipticity conditions \eqref{eq:Kellipt_gen_L}-\eqref{eq:Kellipt_gen_l} as 
\[
0 < \lambda \le r^{2s-2}\int_{B_{  r} }|e\cdot y|^{2} K(dy)\le C \Lambda \qquad\text{for all}\ r  >0\ \text{and}\ e\in \S^{n-1},
\]
for some $C$ depending only on $n$ and $s$.
\end{rem}
\begin{proof}
We divide the proof into two steps.
\begin{steps}
\item 
Let us show \ref{it:equiv_fourier_i} first. We recall that (see \eqref{eq:fourier_symb} or \eqref{eq:sym_psi}) 
\[
\A(\xi) = \int_{\R^n} \big(1-\cos(y\cdot \xi)\big)K(dy).
\]
Since $1-\cos(t) \le \frac12 t^2$, 
\[
\A(\xi) \le 2\int_{\R^n}\min\left\{1, \textstyle{\frac14} |y\cdot \xi|^2\right\} K(dy)\le 2\int_{\R^n}\min\left\{1,  |y|^2|\xi|^2\right\} K\left(dy\right).
\]
In particular, 
\[
\A(\xi) \le 2 \int_{\R^n\setminus B_{1/|\xi|}} K(dy) + 2|\xi|^2 \int_{B_{1/|\xi|}}|y|^2 K(dy).
\]
By assumption and thanks to properties \eqref{eq:prop_kernel} and \eqref{eq:prop_kernel2}, we have
\[
\int_{\R^n\setminus B_{1/|\xi|}} K(dy) \le  C \Lambda |\xi|^{2s}\quad\text{and}\quad \int_{B_{1/|\xi|}} |y|^2 K(dy)\le C \Lambda |\xi|^{2s-2}.
\]
Combining the previous inequalities, we get $\A(\xi)\le C \Lambda |\xi|^{2s}$ for all $\xi\in \R^n$, so we can take $\tilde \Lambda = C\Lambda$.

On the other hand,   let us assume that \eqref{eq:Kellipt_gen_l} holds. We use now that, for $t^2 \le 1$, $1-\cos(t) \ge \frac14 t^2$, to obtain
\[
\A(\xi)   \ge \int_{B_{1/|\xi|}} \big(1 - \cos(y\cdot \xi)\big) K(dy)  \ge \frac {|\xi|^2}{4} \int_{B_{1/|\xi|}} \left|\frac{\xi}{|\xi|}\cdot y\right|^2 K(dy)\ge \frac{\lambda}{4}|\xi|^{2s},
\]
so that we can take $\tilde \lambda =  \frac{\lambda}{4}$.

\item Let us now show \ref{it:equiv_fourier_ii}. We know that, for any $r > 0$, 
\[
r^{2s}\int_{\mathbb{S}^{n-1}} \A\left(\frac{\sigma}{r}\right)\, d\sigma\le |\partial B_1| \tilde \Lambda.
\]
As before, we have that, for some $c$ depending only on $n$ and $s$,
\[
\int_{\mathbb{S}^{n-1}} \big(1-\cos(\lambda \cdot \sigma)\big)\, d\sigma \ge c \min\{1, |\lambda|^2\},
\]
so that,  
\[
\begin{split}
r^{2s}\int_{\mathbb{S}^{n-1}} \A\left(\frac{\sigma}{r}\right)\, d\sigma & = r^{2s}\int_{\mathbb{S}^{n-1}} \int_{\R^n} \left(1-\cos(r^{-1}y\cdot \sigma)\right)K(dy)\, d\sigma\\
& \ge c r^{2s}\int_{\R^n} \min\{1, r^{-2}|y|^2\} K(dy).
\end{split}
\]
In all, we have  
\[
  r^{2s}\int_{\R^n\setminus B_1 } K(r \, dy) \le r^{2s}\int_{\R^n} \min\{1, |y|^2\} K(r \, dy)\le c^{-1} \tilde \Lambda,
\]
and in particular, 
\begin{equation}
\label{eq:prev_arg}
\int_{B_{2r}\setminus B_r} K(dy) \le  \int_{\R^n\setminus B_1} K(r \, dy) \le c^{-1} \tilde \Lambda r^{-2s}
\end{equation}
for all $r >0$, so that we can take any $\Lambda \ge c^{-1}\tilde\Lambda$ for some $c$ that depends only on $n$ and $s$. 

In order to determine $\lambda$ now, observe that for any $e\in \mathbb{S}^{n-1}$ and $r > 0$ we have that 
\[
\begin{split}
\tilde \lambda \le r^{2s} \A(e / r) & = r^{ 2s} \int_{\R^n} \big(1-\cos(y\cdot e)\big) K(r \, dy)\\ & \le 2 r^{ 2s} \int_{\R^n} \min\{1, |y\cdot e|^2\}K(r \, dy),
\end{split}
\]
where we are using, as in the first step, that $1-\cos(t) \le 2\min\{1,t^2\}$. Hence, we can split, for any $R> 1$ to be chosen, as
\[
\tilde \lambda \le 2 r^{ 2s} \int_{B_R} |y\cdot e|^2 K(r \, dy)+ 2 r^{2s} \int_{\R^n\setminus B_R} K(r \, dy).
\]
Observe  now that, on the one hand and thanks to \eqref{eq:prev_arg} combined with property \eqref{eq:prop_kernel},
\[
2 r^{ 2s} \int_{\R^n\setminus B_R} K(r \, dy)=2 r^{2s}\int_{\R^n\setminus B_{Rr}} K(dy) \le C R^{-2s} \tilde \Lambda,
\]
for some $C$ depending only on $n$ and $s$. We therefore get
\[
\frac{\tilde \lambda - CR^{-2s}\tilde\Lambda}{2}\le  r^{ 2s} \int_{B_R} |y\cdot e|^2 K(r \, dy) = r^{2s-2} \int_{Rr} |y\cdot e|^2 K(dy).
\]
By denoting $\rho := Rr$ we obtain
\[
\frac{\tilde \lambda - CR^{-2s}\tilde\Lambda}{2R^{2-2s}}\le  \rho^{2s-2}   \int_{B_\rho} |y\cdot e|^2 K(dy).
\]
Optimizing in $R$ we get 
\[
\rho^{2s-2}   \int_{B_\rho} |y\cdot e|^2 K(dy)\ge \lambda
\]
for some $\lambda > 0$, that can be taken to be $\lambda = C {\tilde\lambda}^{\frac{1}{s}}{\tilde \Lambda}^{1-\frac{1}{s}}$, for some $C$ depending only on $n$ and $s$. 
\qedhere
\end{steps}
\end{proof}

\subsection{Strong solutions}\index{Strong solutions!General operators} Before stating some of our estimates, let us introduce the following space, that we denote $L^\infty_{\tau}(\R^n)$:
\begin{defi}
\label{defi:Linfs}\index{Linfty@$L^\infty_{\tau}$}
Given $\tau \ge 0$, we say that $w\in L^\infty_{\tau}(\R^n)$ if 
\[
\|w\|_{L^\infty_{\tau}(\R^n)} := \left\|\frac{w(x)}{1+|x|^\tau}\right\|_{L^\infty(\R^n)}<\infty.
\]
\end{defi}

The following lemma shows that $\L u$ is well-defined in the $L^\infty$ sense around a point $x\in \R^n$ (and hence, we say that $u$ is a \emph{strong solution} around~$x$) whenever $u\in L^\infty_{2s-\eps}(\R^n)$ and $u$ is $C^{2s+\eps}$ around~$x$. 

\begin{lem}
\label{lem:Lu}
Let $s\in (0, 1)$ and $\L\in \GL$. Let   $u\in C^{2s+\eps}(B_1)\cap L^\infty_{2s-\eps}(\R^n)$ for some $\eps > 0$. Then $\L u\in L^\infty_{\rm loc}(B_1)$ with 
\[
\|\L u\|_{L^\infty(B_{1/2})} \le C \Lambda \left(\|u\|_{C^{2s+\eps}(B_1)} + \|u\|_{L^{\infty}_{2s-\eps}(\R^n)} \right) 
\]
for some $C$ depending only on $n$, $s$, and $\eps$.
\end{lem}
\begin{proof}
We proceed as in Lemma~\ref{lem:laplu}, by splitting for any $x\in B_{1/2}$
\[
\begin{split}
\L u(x) & = \frac12 \int_{B_{1/2}} \left(2u(x) - u(x+y)-u(x-y)\right) K(dy)\\
& \quad  +  u(x) \int_{B^c_{1/2}} K(dy)  - \int_{B^c_{1/2}} u(x+y) K(dy).
\end{split}
\]
If we assume without loss of generality $2s+\eps \le 2$, so that (see, for example, Lemma~\ref{it:H7_gen})
\[
\big|u(x+y) + u(x-y) - 2u(x)\big|\le C[u]_{C^{2s+\eps}(B_1)}|y|^{2s+\eps}
\]
for all $x, y\in B_1$ such that $x+y, x-y\in B_1$, then we see that the first term above is well-defined pointwise at every $x\in B_{1/2}$. The same holds for the second term, while for the last term we have   the convolution of a locally bounded function against a finite measure, which by the decay of $K$ is well-defined in an $L_{\rm loc}^\infty$ sense. In this case, if $K$ is not absolutely continuous, then $\L u$ might not be defined at every point in $B_{1/2}$, but nonetheless can be identified with a function in $L^\infty(B_{1/2})$.

Moreover, for any $x\in B_{1/2}$, and using \eqref{eq:prop_kernel} and \eqref{eq:prop_kernel2}, we have the bound
\[
\begin{split}
|\L u (x)|& \le C  \|u\|_{C^{2s+\eps}(B_1)} \int_{B_{1/2}} |y|^{2s+\eps} K(dy)\\
& \quad + C  \left\|\frac{u(y)}{|y|^{2s-\eps}}\right\|_{L^\infty(\R^n\setminus B_{1/2})} \int_{\R^n\setminus B_{1/2}} |y|^{2s-\eps}K(dy)
\\ & \le C\Lambda \|u\|_{C^{2s+\eps}(B_1)} + C_\eps \Lambda  \|u\|_{L^\infty_{2s-\eps}(\R^n)},
\end{split}
\]
which gives the desired result.
\end{proof}

\begin{rem}
In Chapter~\ref{ch:fully_nonlinear} we will see that, in fact, in order to evaluate $u$ pointwise, it is in general enough to impose a one sided regularity condition. 
We refer to Lemma~\ref{lem:Lu_LL} for an analogous pointwise statement in the context of linear operators with kernels comparable to the fractional Laplacian.
\end{rem}

In case of the fractional Laplacian,  the previous lemma holds with norm $\|u\|_{L^1_{\omega_s}(\R^n)}$  instead of the (stronger) $\|u\|_{L^\infty_{2s-\eps}(\R^n)}$ norm. However,  for general operators $\L\in \GL$, the condition $u\in L^1_{\omega_s}(\R^n)$ is not enough to evaluate $\L u$ in an $L^\infty$ sense, and we need the stronger assumption $u\in L^\infty_{2s-\eps}(\R^n)$ (cf. Lemma~\ref{lem:Lu_LL}, when the condition can be relaxed again).

The previous lemma gives the local boundedness of $\L u$. In order to obtain higher regularity for $\L u$ (as in Lemma~\ref{lem:laplu})   we not only need to impose that $u$ is more regular in $B_1$, but we must make sure that the operator $\L$ is sufficiently regular as well (in the sense of Definition~\ref{defi:Gmu}). 
Alternatively, as shown in part \ref{it:lem_Lu2_ii} below, instead of imposing regularity of $\L$, one may assume \emph{global} regularity of $u$ (contrary to regularity only in $B_1$). Later on, in Lemma~\ref{lem:counterexamples_strong} we will actually show that these conditions are not only sufficient but also necessary. 
 
\begin{lem}
\label{lem:Lu_2}
Let $s\in (0, 1)$ and let $\alpha>0$ with $\alpha\not\in\mathbb N$. Then:
\begin{enumerate}[leftmargin=*,label=(\roman*)]
\item \label{it:lem_Lu2_i} If $\L\in \G_s(\lambda, \Lambda; \alpha)$ (see Definition~\ref{defi:Gmu}), then for any $u\in C^{2s+\alpha}(B_1)\cap L^\infty_{2s-\eps}(\R^n)$ with $\eps > 0$ we have  $\L u\in C_{\rm loc}^\alpha(B_{1})$ and
\[
\|\L u\|_{C^\alpha (B_{1/2})} \le C  \left(\|u\|_{C^{2s+\alpha}(B_1)} + \|u\|_{L^{\infty}_{2s-\eps}(\R^n)}\right) ,
\]
with $C$ depending only on $n$, $s$, $\Lambda$, $\eps$, $\alpha$, and $[\L]_\alpha$.
\item \label{it:lem_Lu2_ii} If $\L\in \GL$, then for any $u\in C^{2s+\alpha}(B_1)\cap C^\alpha(\R^n)$ we have $\L u\in C_{\rm loc}^\alpha(B_{1})$ and
\[
\|\L u\|_{C^\alpha (B_{1/2})} \le C \Lambda \left(\|u\|_{C^{2s+\alpha}(B_1)} +\|u\|_{C^{\alpha}(\R^n)}\right),
\]
with $C$ depending only on $n$, $s$, and $\alpha$.
 
\item \label{it:lem_Lu2_iii} If $\gamma\in (0, \alpha)$   and $\L\in \G_s(\lambda, \Lambda; \alpha-\gamma)$, then for any $u\in C^{2s+\alpha}(B_1)\cap C^\gamma(\R^n)$ we have $\L u\in C_{\rm loc}^\alpha(B_{1})$ and
\[
\|\L u\|_{C^\alpha (B_{1/2})} \le C \left(\|u\|_{C^{2s+\alpha}(B_1)} +\|u\|_{C^{\gamma}(\R^n)}\right),
\]
with $C$ depending only on $n$, $s$, $\Lambda$, $\alpha$, $\gamma$, and $[\L]_{\alpha-\gamma}$.
\end{enumerate}
\end{lem}

\begin{proof} In case \ref{it:lem_Lu2_i}, up to dividing by $\|u\|_{C^{2s+\alpha}(B_1)} +\|u\|_{L^{\infty}_{2s-\eps}(\R^n)}$ we may assume $\|u\|_{C^{2s+\alpha}(B_1)} +\|u\|_{L^{\infty}_{2s-\eps}(\R^n)} = 1$. In case \ref{it:lem_Lu2_ii}, up to dividing by $ \|u\|_{C^{2s+\alpha}(B_1)} + \|u\|_{C^{\alpha}(\R^n)}$ we assume $\|u\|_{C^{2s+\alpha}(B_1)} + \|u\|_{C^{\alpha}(\R^n)} = 1$. In case \ref{it:lem_Lu2_iii}, up to dividing by $ \|u\|_{C^{2s+\alpha}(B_1)} + \|u\|_{C^{\gamma}(\R^n)}$ we assume $\|u\|_{C^{2s+\alpha}(B_1)} + \|u\|_{C^{\gamma}(\R^n)} = 1$. 

We proceed analogously to the case of the fractional Laplacian $\fls$ in  Lemma~\ref{lem:laplu}. We fix a cut-off function $\eta\in C^\infty_c(\R^n)$ such that $\eta \ge 0$, $\eta \equiv 0$ in $\R^n\setminus B_{3/4}$ and $\eta \equiv 1$ in $B_{2/3}$, and define 
\[
u_1 := \eta u\qquad\text{and}\qquad u_2 := (1-\eta) u,
\]
so that $u = u_1 + u_2$, with $u_1$ being compactly supported in $B_{3/4}$ and $u_2$ satisfying that $u_2 \equiv 0$ in $B_{2/3}$. We divide the proof into three steps:

\begin{steps}

\item \label{step1:laplu} Let us compute first a bound for $\|\L u_1\|_{C^\alpha(B_{1/2})}$ in all cases \ref{it:lem_Lu2_i}, \ref{it:lem_Lu2_ii}, and \ref{it:lem_Lu2_iii}, by using that $u_1\in C^{2s+\alpha}(\R^n)$ with $\|u_1\|_{C^{2s+\alpha}(\R^n)}\le C$ for some universal $C$. Notice that the $L^\infty$ bound follows from Lemma~\ref{lem:Lu}. 

We start by assuming $ \alpha < 1$. Let $x\in B_{1/2}$ be fixed, and $r := |x|$. We split
\[
\begin{split}
\L u_1(x) & = \frac12 \int_{B_r}\big(2u_1(x)-u_1(x+y)-u_1(x-y)\big)K(dy)\\
& \quad +  \frac12 \int_{\R^n\setminus B_r}\big(2u_1(x)-u_1(x+y)-u_1(x-y)\big)K(dy).
\end{split}
\]
Then, if $2s+\alpha \le 2$,
\[
\begin{split}
\big|\L u_1(x) + \L u_1(-x) - 2 \L u_1(0) \big|& \le  C\int_{B_r} |y|^{2s+\alpha} K(dy) \\
&\quad  + C \int_{\R^n\setminus B_r}|x|^{2s+\alpha} K(dy),
\end{split}
\]
where we are using 
\begin{equation}
\label{eq:toused124}
\begin{split}
|u_1(x+y)+u_1(x-y)-2u_1(x)| &\le  C |y|^{2s+\alpha},\\
|u_1(x)+u_1(-x)-2u_1(0)| &\le  C |x|^{2s+\alpha},\\
|u_1(x\pm y)+u_1(-x\pm y)-2u_1(\pm y)| &\le  C |x|^{2s+\alpha}.
\end{split}
\end{equation}
Since $K$ satisfies \eqref{eq:Kellipt_gen_L}, we have by \eqref{eq:prop_kernel}-\eqref{eq:prop_kernel2}, 
\[
\int_{\R^n\setminus B_r}r^{2s+\alpha} K(dy) \le C\Lambda r^\alpha\quad\text{and}\quad \int_{B_r}|y|^{2s+\alpha} K(dy) \le C\Lambda r^\alpha,
\]
so that 
\begin{equation}
\label{eq:weobtainagain}
\big|\L u_1(x) + \L u_1(-x) - 2 \L u_1(0) \big| \le C \Lambda |x|^\alpha.
\end{equation}
On the other hand, if $2<2s+\alpha<3$,  by Lemma \ref{lem:A_imp_2}-\ref{it:A:3} we have 
\begin{equation}
\label{eq:2salpha1}
\begin{split}
& \big|u_1(x+y)+u_1(x-y)-2u_1(x)-u_1(y)-u_1(-y)+2u_1(0)\big| \leq \\
& \qquad \le |y|^2 |x|^{2s+\alpha-2}.
\end{split}
\end{equation}
and 
\begin{equation}
\label{eq:2salpha2}
\begin{split}
& \big|u_1(x+y)+u_1(-x+y)-2u_1(y)-u_1(x)-u_1(-x)+2u_1(0)\big| \leq \\
& \qquad \le |x|^2 |y|^{2s+\alpha-2}.
\end{split}
\end{equation}

Thus, using \eqref{eq:2salpha1}-\eqref{eq:2salpha2} and thanks to properties \eqref{eq:prop_kernel}-\eqref{eq:prop_kernel2}-\eqref{eq:prop_kernel3},
\[
\begin{split}
 \big|\L u_1(x) + \L u_1(-x) - 2 \L u_1(0) \big|  & \le  C\int_{B_r} |y|^{2}|x|^{2s+\alpha-2}K(dy)  \\
& \quad +C\int_{B_r^c} |x|^{2}|y|^{2s+\alpha-2} K(dy) \le C\Lambda r^\alpha,
\end{split}
\]
that is, we obtain again \eqref{eq:weobtainagain}.  Repeating around any point in $B_{1/2}$ we get $\|\L u_1\|_{C^\alpha(B_{1/2})} \le C\Lambda $ by Lemma~\ref{it:H7}, if $\alpha < 1$. If $\alpha = k +\beta$ with $k \in \N$, $\beta\in (0, 1)$, we take $\ell$ derivatives of $\L u_1$ with $\ell =1,\dots k$ and repeat the arguments above iteratively, to obtain $\|D^\ell \L u_1\|_{C^\beta(B_{1/2})} \le C\Lambda $ for all $\ell\in\{0,\dots,k\}$.  In all, we get
\[
\|\L u_1\|_{C^\alpha(B_{1/2})}\le C\Lambda,
\]
which is the desired result for $u_1$. 

\item \label{it:step2Lu2}  Let us now find a bound for $\|\L u_2\|_{C^\alpha(B_{1/2})}$, where we denote  $\alpha = k+\beta$, with $k \in \N_0$ and $\beta\in (0, 1)$. We first consider case \ref{it:lem_Lu2_ii}, so $u_2\in C^{\alpha}(\R^n)$ with $[D^k u_2]_{C^\beta(\R^n)} \le C$ and $u_2 \equiv 0$ in $B_{2/3}$. In this case, for any $x\in B_{1/2}$ (using property \eqref{eq:prop_kernel} as before), 
\[
\begin{split}
\left|D^k \L u_2(x) - D^k \L u_2(0) \right|& = \left|\L D^k u_2(x) - \L D^k u_2(0) \right|\\
&  \le 
 C \int_{\R^n\setminus B_{1/6}}|x|^{\beta} K(dy) \le C \Lambda |x|^{\beta},
 \end{split}
\]
where we have used
\[
\begin{split}
|D^k u_2(x+y)-D^k u_2(y)| &\le C |x|^{\beta}\quad\text{for all}\quad y \in \R^n,
\end{split}
\]
and the fact that $u_2 \equiv 0$ in $B_{2/3}$, so that $u_2(x+ y) = 0$ for any $x\in B_{1/2}$ and $y\in B_{1/6}$. Combined with the $L^\infty$ bound in Lemma~\ref{lem:Lu} and interpolation inequalities (Proposition~\ref{it:H9}), we obtain the desired result in case \ref{it:lem_Lu2_ii}.  

We now prove the first statement, \ref{it:lem_Lu2_i}, where we have the extra hypothesis on the kernel, \eqref{Calpha-assumption}, and $\|u_2\|_{L^\infty_{2s-\eps}(\R^n)} \le C$ for some $C$ universal.  Let $x_1, x_2\in B_{1/2}$, so that (since $u_2 \equiv 0$ in $B_{2/3}$), 
\[
(D^k \L u_2) (x_i) =  -\int_{B_{1/8}^c(x_i)} u_2(z) D^k K(-x_1+dz)  
\]
and 
\[
\left|D^k \L u_2(x_1) - D^k \L u_2 (x_2) \right| \le \int_D |u_2(z)| \left|D^k K(-x_1+dz) - D^k K(-x_2+dz)\right|,
\]
where $D$ denotes the domain $D = \R^n\setminus (B_{1/8}(x_1)\cup B_{1/8}(x_2))$. We can now assume without loss of generality that $|x_1-x_2|\le \frac{1}{32}$ (otherwise, take a finite chain of inequalities) so that by assumption on the kernel (see \eqref{eq:prop_kernel_Calpha} and \eqref{eq:prop_kernel4}) we have 
\begin{equation}
\label{eq:weneedtoboundlaststep}
\begin{split}
\left|D^k \L u_2(x_1) - D^k \L u_2 (x_2) \right| \leq & \\
& \hspace{-4cm} \le \int_{B_{1/16}^c}|u_2(z)| \left|D^k K(-x_1+dz) - D^k K(-x_2+dz)\right| \\
& \hspace{-4cm}\le C \left\|\frac{u_2(z)}{|z|^{2s-\eps}}\right\|_{L^\infty(B^c_{1/16})}\int_{B_{1/16}^c} |z|^{2s-\eps}\big|D^k K(-x_1+dz) - D^k K(-x_2+dz)\big| \\
& \hspace{-4cm} \le C [K]_\alpha  \|u\|_{L^\infty_{2s-\eps}(\R^n)}  |x_1-x_2|^\beta,
\end{split}
\end{equation}
which combined with Lemma~\ref{lem:Lu} completes the proof of \ref{it:lem_Lu2_i}. For $\alpha > 1$, we proceed as before, by taking derivatives.
\item \label{it:step3Lu2}
  To finish, we show \ref{it:lem_Lu2_iii} by using \ref{it:lem_Lu2_i}. Up to taking $\lfloor\gamma\rfloor$ derivatives and arguing as before, we may assume $\gamma \in (0,1)$. 

Let us consider
\[
D^\gamma_h u_2(x) :=\frac{u_2(x+h) - u_2(x)}{|h|^\gamma},
\]
for some $|h|$ sufficiently small. By assumption, it satisfies $\|D^\gamma_h u_2\|_{L^\infty(\R^n)}\le C$ and $D_h^\gamma u_2 \equiv 0$ in $B_{1/2}$ for $h$ small enough, and we can apply the first result, \ref{it:lem_Lu2_i}, in balls $B_{1/4}$ and $B_{1/2}$ to get:
\[
\|\L D^\gamma_h u_2\|_{C^{\alpha-\gamma} (B_{1/4})} \le C   \left(\|D^\gamma_h u_2\|_{C^{2s+\alpha}(B_{1/2})} + \|D^\gamma_h u_2\|_{L^{\infty}(\R^n)}\right) \le C,
\]
where we have used $u_2 \equiv 0$ in $B_{2/3}$ and $|h|$ small. Using $\L D_h^\gamma u_2 = D_h^\gamma \L u_2$ and Lemma~\ref{it:H8}, we obtain the desired result also in case \ref{it:lem_Lu2_iii}, after a covering argument.
\qedhere
\end{steps}
\end{proof}

 Some remarks are in order:
 
\begin{rem}
\label{rem:fullnorm_seminorm}
From the previous proof, if $\alpha \in (0, 1)$ we could have taken in the right-hand side of Lemma~\ref{lem:Lu_2}-\ref{it:lem_Lu2_ii} above the value $[u]_{C^\alpha(\R^n)}+\|u\|_{L^\infty_{2s-\eps}(\R^n)}$ instead  of the full norm $\|u\|_{C^\alpha(\R^n)}$. When $\alpha > 1$, however, we need some condition to ensure that $\L D^k u$ is well-defined, with $k = \lfloor\alpha\rfloor$. 
\end{rem}

\begin{rem}
\label{rem:Lu_2}
In cases \ref{it:lem_Lu2_i} and \ref{it:lem_Lu2_iii}, the previous proof actually gives a more refined estimate on the bounds of the seminorm $[\L u]_{C^\alpha(B_{1/2})}$, in which we can also make explicit the dependence on $\Lambda$ and $[\L]_\alpha$ or $[\L]_{\alpha-\gamma}$. 

Indeed, by simply inspecting the proof we see that, in fact, we have the  following bounds:
\begin{itemize}[leftmargin=*]
\item In the first case, Lemma~\ref{lem:Lu_2}-\ref{it:lem_Lu2_i}, we showed:
\[
[\L u]_{C^\alpha (B_{1/2})} \le C  \left(\Lambda \|u\|_{C^{2s+\alpha}(B_1)} + [\L]_\alpha \|u\|_{L^{\infty}_{2s-\eps}(\R^n)}\right),
\]
for some $C$ depending only on $n$, $s$, $\eps$, and $\alpha$. 
\item  In the third case, Lemma~\ref{lem:Lu_2}-\ref{it:lem_Lu2_iii}, we showed:
\begin{equation}
\label{eq:rem_Lu22}
[\L u]_{C^\alpha (B_{1/2})} \le C  \left(\Lambda \|u\|_{C^{2s+\alpha}(B_1)} + [\L]_{\alpha-\gamma} \left(\|u\|_{L^\infty_{2s-\eps}(\R^n)} + [u]_{C^\gamma(\R^n)}\right)\right),
\end{equation}
for some $C$ depending only on $n$, $s$, $\eps$,  $\alpha$, and $\gamma$. 
\end{itemize}

We can  then combine these estimates with the one in Lemma~\ref{lem:Lu} to get more refined bounds for the corresponding full norms. 
\end{rem}
\begin{rem}
\label{rem:nonneg_kernels}
Even if Lemmas~\ref{lem:Lu} and \ref{lem:Lu_2} are stated for operators belonging to $\GL$ (with the corresponding regularity assumptions), in their proofs we have never  used the lower ellipticity assumptions, \eqref{eq:Kellipt_gen_l}, which can be seen  from the lack of dependence on $\lambda$. In fact, we have not even used the nonnegativity of the kernel, and we are only using that $K$ is even (i.e., symmetric) and satisfies
\[
r^{2s}\int_{B_{2r}\setminus B_r} |K(dy)|  \le \Lambda < \infty\qquad\text{for all}\quad r > 0,
\]
as well as the corresponding regularity assumption, \eqref{Calpha-assumption} or \eqref{Cmu-assumption}  in Lemma~\ref{lem:Lu_2}-\ref{it:lem_Lu2_i} and \ref{it:lem_Lu2_iii}.
\end{rem}

As shown next, the higher regularity of the kernel (or global regularity for $u$) is necessary if we want to obtain regularity for $\L u$: 

\begin{lem}\label{lem:counterexamples_strong}\index{Counter-example}
There exists $\L \in \GL$ with kernel $K(y)\, dy$ comparable to the fractional Laplacian, \eqref{eq:Kcompfls}, and $u\in L^\infty(\R^n)$  compactly supported and with $u \equiv 0$ in $B_2$, such that $\L u \notin C^\eps(B_1)$ for any $\eps > 0$. 
\end{lem}
\begin{proof}
We prove it in dimension $n=1$. 
Let us first show that there exists a function on $0\le w_0 \le 1$ in $\R$ with $w_0\equiv 0$ in $\R\setminus (0, 1)$, such that 
\begin{equation}
\label{eq:w1cont}
w_1(x) := \int_{-1}^1 w_0(x+y)w_0(y)\,dy \notin C^\eps((-1, 1))\quad\text{for any}\quad \eps > 0. 
\end{equation}
To construct it, let $\varphi$ be a square wave,  $\varphi(x) = 1$ if $x\in (2k, 2k+1)$ for some $k\in \Z$, and $\varphi(x) = 0$ otherwise. That is, 
\[
\varphi = \sum_{k\in \Z} \chi_{(2k, 2k+1)}.
\]
 Let us define, for some $m, t>0$ to be chosen with $m/t\gg 1$, 
\[
\psi_{m, t}(x) = \varphi(x/t)\chi_{(0, m)},
\]
so that 
\begin{equation}
\label{eq:sothatthanks}
\int_0^m \psi_{m, t}(y+t)\psi_{m, t}(y) \, dx= 0\quad\text{and}\quad \int_0^m \psi_{m, t}^2(y)\, dy  \approx \frac{m}{2},
\end{equation}
if $m \gg t$. If we fix $m_k = k^{-2}/4$ and $t_k = e^{-k^2}$, and define 
\[
g_k(x) := \psi_{m_k, t_k}(x+\sigma_k)\quad\text{where}\quad \sigma_k := \frac12 \sum_{i = 1}^k i^{-2}<1,
\]
 we can take 
\[
w_0(x) := \sum_{k \ge 3} g_k(x),
\]
and $w_1$ as in \eqref{eq:w1cont}. In particular, observe that since $g_k$ is supported in $(\sigma_k, \sigma_k+m_k)$ and since $\sigma_k +m_k < \sigma_{k+1}$ for $k \ge 3$, we have that all $g_k$ have disjoint supports, and since $\sigma_k < 1$ for all $k\in \N$, they are all in $(0, 1)$. 

Then $0\le w_0\le 1$, $w_0$ has support in $(0, 1)$, $w_1(0) \ge w_1(x)$ for any $x\in \R$ and 
\[
w_1(0)-w_1(t_k) \ge \frac{m_k}{3}\quad\text{if $k$ is large enough},
\]
thanks to \eqref{eq:sothatthanks}. Since $m_k = 1/|4\log(t_k)|$, we get
\[
w_1(0)-w_1(t_k) \ge 1/|12\log(t_k)|\quad\text{if $k$ is large enough},
\]
and $w_1\notin C^\eps((-1, 1))$ for any $\eps > 0$. 

Let $\L \in \GL$ with $n = 1$ and kernel $K(y)\, dy$ be given by 
\[
K(y) := \frac{1}{|y|^{1+2s}}+w_0(y+4)+w_0(-y-4),
\]
which is comparable to the fractional Laplacian and symmetric (but nonsmooth), and let $u(x) = w_0(x+4)$. Thanks to \eqref{eq:w1cont} we have that $\L u \notin C^\eps(B_1)$ for any $\eps > 0$, and $u$ is compactly supported with $u = 0$ in $B_2$. 
\end{proof}

We also prove the following result, that gives estimates on $\L\varphi$ whenever $\varphi\in C^\infty_c(\R^n)$. In particular, it implies that $\L \varphi\in L^1(\R^n)\cap L^\infty(\R^n)$. 

\begin{lem}
\label{lem:welldefined}
Let $s\in (0, 1)$, $\L \in \GL$,  $\Omega\subset \R^n$ be a bounded domain with $0\in \Omega$, and  $\alpha\in [0, 2s]$. Then, for any $\varphi \in C^\infty_c(\Omega)$, we have
\[
\int_{\R^n\setminus B_r} \big(1+|x|^{2s-\alpha}\big)|\L \varphi|\le C\|\varphi\|_{C^2(\R^n)}(1+r)^{-\alpha}\quad \text{for any}\quad r \ge 0,
\]
and 
\[
|\L \varphi(x)| \le C\|\varphi\|_{C^2(\R^n)}(1+|x|)^{-2s}\quad\text{for all}\quad x\in \R^n,
\]
for some $C$ depending only on $n$, $s$, $\alpha$, ${\rm diam}(\Omega)$, $\lambda$, and $\Lambda$.
\end{lem}
\begin{proof}
Let $\Omega \subset B_R$ with $ R= {\rm diam}(\Omega)$.  

If $r > 2R$, we can compute, using that $\varphi(x+y) = 0$ for all $x\in \R^n\setminus B_r$ and $y \in B_{r-R}$
\[
\begin{split}
\int_{\R^n\setminus B_r} |\L \varphi| & \le  \int_{\R^n \setminus B_r} \int_{ \R^n}|\varphi(x+y)|K(dy)\, dx\\
& \le \|\varphi\|_{L^\infty(\R^n)}  \int_{\R^n \setminus B_{r-R}} K(dy)\\
& \le C\|\varphi\|_{L^\infty(\R^n)}  r^{-2s}. 
\end{split}
\]
where we have also used \eqref{eq:prop_kernel} and $r-R \ge \frac12 r$.
On the other hand, if $r \le 2R$ we immediately have that $\int_{B_{2R}} |\L \varphi|\le C\|\varphi\|_{C^2(\R^n)}$ thanks to Lemma~\ref{lem:Lu}. 

We assume, therefore,  $r > 2R$, in which case
\[
\int_{\R^n\setminus B_r}\big(1+|x|^{2s-\alpha}\big)|\L\varphi|  \le  C \int_{\R^n \setminus B_{r}} |x|^{2s-\alpha}|\L\varphi|
\]
and we can bound it as 
\[
\begin{split}
\int_{\R^n \setminus B_{r}} |x|^{2s-\alpha}|\L\varphi| & \le  \sum_{k \ge 0}2^{(k+1)(2s-\alpha)}r^{2s-\alpha}  \int_{B_{2^{k+1}r}\setminus B_{2^kr}} |\L\varphi| \\
& \le C r^{-\alpha} \sum_{k \ge 0} 2^{-\alpha k}\\
& \le Cr^{-\alpha}. 
\end{split}
\]
Combining the previous inequalities gives the first result. 

For the second part, we assume that $|x|\ge 2R$ (otherwise, we are done by Lemma~\ref{lem:Lu}) and hence
\[
\begin{split}
|\L u(x)| & \le C \int_{\R^n} |\varphi(x+y)|K(dy) \\
& \le \|\varphi\|_{L^\infty(\R^n)}\int_{\R^n \setminus B_{|x|-R}} K(dy)\\
& \le C (|x|-R)^{-2s}\\
& \le C |x|^{-2s},
\end{split}
\]
where we have used that $\varphi(x+y) = 0$ if $x+y\notin B_R$, and $|x|\ge 2R$. This completes the proof. 
\end{proof}

\subsection{Fractional Sobolev spaces}
\index{Fractional Sobolev space}
As we will see, the natural space in which to study \emph{weak solutions} to the equation $\L u = f$ is the \emph{fractional Sobolev space} $H^s(\R^n)$. 

Let us give a very brief introduction to the fractional Sobolev space $H^s(\R^n)$ (we refer to \cite{DPV12, EE22, Leo23} for further details; see also \cite{MRS}). 

\begin{defi}[Fractional Sobolev space] The fractional Sobolev space $H^s(\R^n)$ is defined as 
\[
H^s(\R^n) := \left\{ u\in L^2(\R^n) : \|u\|_{H^s(\R^n)} < \infty  \right\},
\]
with norm
\[
\|u \|_{H^s(\R^n)} = \|u\|_{L^2(\R^n)} + [u]_{H^s(\R^n)},
\]
and seminorm
\[
[u]_{H^s(\R^n)} := \left(\int_{\R^n}\int_{  \R^n} \frac{(u(x) - u(y))^2}{|x-y|^{n+2s}}\, dx\, dy \right)^{\frac12}.
\]
\end{defi}

We could alternatively define fractional Sobolev spaces using the Fourier transform 
\[
H^s(\R^n) := \left\{ u\in L^2(\R^n) : (1+|\xi|^{2})^{\frac{s}{2}}\mathcal{F}(u)(\xi) \in L^2 (\R^n) \right\}
\]
and taking the corresponding equivalent norm in the Fourier space. In fact, we have that
\begin{equation}
\label{eq:ufourier}
[u]^2_{H^s(\R^n)} = \frac{2}{c_{n,s}} \int_{\R^n} (\mathcal{F}(u)(\xi))^2 |\xi|^{2s}\, d\xi. 
\end{equation}
 with $c_{n,s}$ given by \eqref{eq:cns}.

The space $H^s(\R^n)$ is a Hilbert space, with scalar product given by 
\[
\langle u, v \rangle_{H^s(\R^n)} = \int_{\R^n} u v + \iint_{\R^n\times \R^n} \frac{(u(x) - u(y))(v(x) - v(y))}{|x-y|^{n+2s}}\, dx\, dy. 
\]
Moreover, we also have a fractional Sobolev inequality in $H^s(\R^n)$ (we refer to \cite[Proposition 15.5]{Ponce} for a very short proof due to Brezis, and to \cite{SV11} for a different proof): 
\begin{thm}[Fractional Sobolev inequality]\index{Fractional Sobolev inequality}
\label{thm:FSI}
Let $s\in (0, 1)$, $u\in H^s(\R^n)$, and $n> 2s$. Then 
\[
\|u\|_{L^q(\R^n)} \le C [u]_{H^s(\R^n)},
\]
where $q = \frac{2n}{n-2s}$, for some $C$ depending only on $n$ and $s$.
\end{thm}

And a Poincar\'e--Friedrichs inequality: 
 
\begin{prop}[Fractional Poincar\'e inequality]
\label{prop:Poincare}\index{Poincar\'e's inequality}
Let $s\in (0, 1)$, and let $\Omega\subset \R^n$ be a bounded domain. Let $u \in H^s(\R^n)$ with $u \equiv 0$ in $\R^n\setminus \Omega$. Then, 
\[
\|u\|_{L^2(\Omega)} \le C [u]_{H^s(\R^n)},
\]
for some $C$ depending only on $\Omega$, $n$, and $s$.
\end{prop}
\begin{proof}
Since $\Omega$ is bounded, for any $x\in \Omega$ we have 
\[
\int_{\Omega^c} \frac{dy}{|x-y|^{n+2s}} \ge c,
\]
for some $c>0$ depending only on $\Omega$, $n$, and $s$. Hence, 
\[
|u(x)|^2 \le \frac{1}{c} \int_{\Omega^c} \frac{|u(x)|^2}{|x-y|^{n+2s}} \, dy \le \frac{1}{c} \int_{\R^n} \frac{|u(x)-u(y)|^2}{|x-y|^{n+2s}} \, dy,
\]
where we are using that $u(y) = 0$ for $y \in \Omega^c$. Integrating in $x\in \R^n$ we get the desired result. 
\end{proof}

\subsection{Integration by parts and weak solutions}
\index{Integration by parts!General operators}
Any L\'evy operator  (and in particular, any $\L \in \GL$) has an integration by parts formula. In order to write it, we   define the following symmetric and positive semi-definite bilinear form
$
\langle\cdot,\cdot\rangle_K:C_c^\infty(\R^n)\times C^\infty_c(\R^n)\to \R
$
given by 
\begin{equation}
\label{eq:bilin}
\langle u, v \rangle_K := \frac12 \iint_{\R^n\times\R^n} \big(u(x) - u(x+y)\big)\big(v(x)-v(x+y)\big) K(dy)\, dx.
\end{equation}

 Observe that, for notational simplicity only, in this section that  we will sometimes assume 
\[
K(dy) = K(y)\, dy.
\]
In this case, the expression \eqref{eq:bilin} has a symmetric representation as 
\begin{equation}
\label{eq:bilin_sym}
\langle u, v \rangle_K := \frac12 \iint_{\R^n\times\R^n} \big(u(x) - u(z)\big)\big(v(x)-v(z)\big) K(x-z)\, dx\, dz.
\end{equation}
\begin{rem}\label{rem:bilin_sym}
One can also write \eqref{eq:bilin} in a symmetric form like \eqref{eq:bilin_sym} by considering  $\tilde K(d(x, z))$ defined on $\R^{2n}$ to be the measure such that $
\tilde K(A\times B) = K(A-B) 
$
for any Borel sets $A, B\subset\R^n$, where $A-B := \bigcup_{\substack{a\in A\\ b\in B}}\{a-b\}$. Then, \eqref{eq:bilin} becomes 
\[
\langle u, v \rangle_K := \frac12 \iint_{\R^n\times\R^n} \big(u(x) - u(z)\big)\big(v(x)-v(z)\big) \tilde K(d(x, z)).
\]
\end{rem}

\begin{lem}[Integration by parts]
\label{lem:integration_by_parts}
Let $s\in (0, 1)$ and let $\L\in \GL$. Let $u, v\in C^\infty_c(\R^n)$. Then, 
\[
\int_{\R^n}  u \, \L v = \langle u, v \rangle_K  = \int_{\R^n} \L u \, v.
\]
\end{lem}

\begin{proof}
We first observe that, for $u, v\in C^\infty_c(\R^n)$, $\langle u, v\rangle_K$ is well-defined, since the integrand is absolutely integrable. Indeed, 
\[
\begin{split}
&\frac12 \iint_{\R^n\times \R^n}\big|u(x+y) - u(x)\big|\big|v(x+y)-v(x)\big| K(dy)\, dx\\
& \quad \le \frac14 \iint_{\R^n\times \R^n}\big|u(x+y) - u(x)\big|^2 K(dy)\, dx\\
& \qquad +\frac14 \iint_{\R^n\times \R^n}\big|v(x+y) - v(x)\big|^2 K(dy)\, dx.
\end{split}
\]
Since, whenever ${\rm supp}(u) \subset B_R$,  
\[
\big|u(x+y)-u(x)\big|^2\le C\|u\|_{C^2(\R^n)}\min\{1, |y|^2\}\chi_{B_R\cup B_R(-y)}(x),\]
the previous integrals are finite (by \eqref{eq:Kint}) and $\langle u, v\rangle_K$ is well-defined. In particular, we also have
\[
\langle u , v\rangle_K = \frac12 \lim_{\eps\downarrow 0} \iint_{B_\eps^c\times \R^n}\big(u(x+y) - u(x)\big)\big(v(x+y)-v(x)\big) K(dy)\, dx.
\]

By considering the domain $B_\eps^c\times \R^n$, we can split the above integral into two integrable terms, 
\begin{equation}
\label{eq:splitlike}
\begin{split}
  \iint_{B_\eps^c\times \R^n}&\big(u(x+y) - u(x)\big)\big(v(x+y)-v(x)\big) K(dy)\, dx=\\
 &  = \iint_{B_\eps^c\times \R^n}u(x+y) \big(v(x+y)-v(x)\big) K(dy)\, dx\\
 &\quad - \iint_{B_\eps^c\times \R^n}u(x)\big(v(x+y)-v(x)\big) K(dy)\, dx\\
 &   = \iint_{B_\eps^c\times \R^n}u(x)\big(2v(x)-v(x-y)-v(x+y)\big) K(dy)\, dx,
 \end{split}
\end{equation}
where in the second equality we have performed a change of variables   $x\mapsto x-y$ (and Fubini). In particular, we can take the limit $\eps\downarrow 0$, and since all terms are well-defined we obtain 
\[
\langle u, v\rangle_K = \int_{\R^n} u \, \L v. 
\]
By symmetry in $u$ and $v$ we get the desired result. 
\end{proof}

Observe that, in the previous integration by parts formula (and similarly to what occurs in the local case), in order to make sense of the term $\langle u, v \rangle_K$ we require less regularity on $v$ than to compute $\L v$ in the left-hand side (recall Lemma~\ref{lem:Lu}). In particular, we have the following:  

\begin{lem}
\label{lem:u_fourier}
Let $s\in (0, 1)$ and let $\L\in \GL$. Let $u, v\in C^\infty_c(\R^n)$. Then, 
\[
\langle u, v \rangle_K  = \int_{\R^n} \mathcal{F}(u)(\xi) \mathcal{F}(v)(\xi) \A(\xi)\, d\xi, 
\]
where $\A$ is the Fourier symbol of $\L$. In particular, 
\begin{equation}
\label{eq:Plancherel}
\langle u, u\rangle_K  = \int_{\R^n} \big(\mathcal{F}(u)(\xi)\big)^2 \A(\xi)\, d\xi\asymp [u]^2_{H^s(\R^n)}, 
\end{equation}
and the bilinear form $\langle\cdot, \cdot\rangle_K$ is  defined on $H^s(\R^n)\times H^s(\R^n)$. 
\end{lem}
\begin{proof}
The first equality follows by Plancherel's theorem,\footnote{Observe that we can use it, since $\L v \in L^1(\R^n)\cap L^\infty(\R^n)$ by Lemma~\ref{lem:welldefined}, and hence $\L v\in L^2(\R^n)$}
\[
\langle u, v \rangle_K = \int_{\R^n} u \, \L v = \int_{\R^n} \mathcal{F}(u)(\xi) \mathcal{F}(\L v)(\xi)\, d\xi = \int_{\R^n} \mathcal{F}(u)(\xi) \mathcal{F}(v)(\xi)\A(\xi)\, d\xi. 
\]
Then \eqref{eq:Plancherel} is a consequence of the comparability of $\A$ with $|\xi|^{2s}$ (recall \eqref{eq:GSLO} or Proposition~\ref{prop:equiv_fourier}~\ref{it:equiv_fourier_i}), together with \eqref{eq:ufourier}. 

Finally, $\langle\cdot, \cdot\rangle_K$ is a symmetric positive semi-definite bilinear form, and as such it satisfies a Cauchy-Schwarz inequality:
\[
\langle u, v \rangle_K^2 \le \langle u, u\rangle_K \langle v, v\rangle_K \le C[u]^2_{H^s(\R^n)}[v]^2_{H^s(\R^n)}.
\]
By a density argument, $\langle u, v \rangle_K<\infty$ as soon as $u, v\in H^s(\R^n)$. 
\end{proof}

This allows us to introduce the notion of \emph{weak solution} for $u\in H^s(\R^n)$ (in fact, for $u$ satisfying \eqref{eq:u_weak_sol} below); as well as \emph{weak supersolution} and \emph{weak subsolution}.  To do that, we define first a new bilinear form adapted to a given domain $\Omega$:

\begin{defi}\label{defi:bilin_om}
When $K$ is absolutely continuous, we define the following bilinear form,
\begin{equation}
\label{eq:prod_KO}
\langle u, v\rangle_{K;\Omega} := \frac12 \iint_{\R^n\times\R^n \setminus (\Omega^c\times \Omega^c)} \big(u(x) - u(z)\big)\big(v(x) - v(z)\big)K(x-z)\, dx\, dz.
\end{equation}
When $K$ is not absolutely continuous, we can replace $K(x-z)\, dx\, dz$  by $\tilde K(d(x, z))$ from Remark~\ref{rem:bilin_sym}, or alternatively, consider instead
\begin{equation}
\label{eq:prod_KO_y}
\langle u, v\rangle_{K;\Omega} := \frac12 \iint_{\R^n\times\R^n \setminus (x, x+y\in \Omega^c)}\hspace{-3.5mm} \big(u(x) - u(x+y)\big)\big(v(x) - v(x+y)\big)K(dy)\, dx.
\end{equation}
\end{defi}
Observe that $\langle\cdot,\cdot\rangle_{K;\Omega}$ is a modification of $\langle\cdot, \cdot\rangle_K$, in which the integration is performed everywhere except for exterior-exterior interactions, that is, except at $\Omega^c\times\Omega^c$.

\begin{defi}[Weak solution, supersolution, and subsolution]\index{Weak solutions}
\label{defi:weaksol}
Let $s\in (0, 1)$ and let $\L\in \GL$. 
Let $\Omega\subset \R^n$ be any bounded domain and let $f\in L^p(\Omega)$ for $p \ge \frac{2n}{n+2s}$ and $n > 2s$. 
We say that $u$ such that 
\begin{equation}
\label{eq:u_weak_sol}
\langle u, u\rangle_{K;\Omega} 
  <\infty
\end{equation}
is a \emph{weak solution} of 
\begin{equation}
\label{eq:dirpb}
\L u  =  f  \quad\text{in}\quad \Omega,
\end{equation}
if
\begin{equation}
\label{eq:weak_sol}
\frac12 \iint_{\R^n\times\R^n} \big(u(x) - u(z)\big)\big(\varphi(x)-\varphi(z)\big) K(x-z)\, dx\, dz = \int_\Omega f \varphi
\end{equation}
for all $\varphi \in H^s(\R^n)$ with $\varphi \equiv 0$ in $\R^n\setminus \Omega$, when $K$ is absolutely continuous. In particular, $u\in H^s(\R^n)$ is a \emph{weak solution} to \eqref{eq:dirpb} if it satisfies \eqref{eq:weak_sol} for all $\varphi \in H^s(\R^n)$ with $\varphi \equiv 0$ in $\R^n\setminus \Omega$. 

We say that $u$ satisfying \eqref{eq:u_weak_sol} is a \emph{weak supersolution} in $\Omega$ (resp. \emph{weak subsolution}) to the equation $\L u = f$, and we denote it 
$
\L u \ge f
$ in $\Omega$ (resp. $\L u \le f$ in $\Omega$) 
if 
\begin{equation}
\label{eq:weak_subsol}
\begin{split}
\frac12 \iint_{\R^n\times\R^n} \big(u(x) - u(z)\big)\big(\varphi(x)-\varphi(z)\big) K(x-z)\, dx\, dz \underset{\left(\text{resp. $\le$}\right)}{\ge} \int_\Omega f \varphi,
\end{split}
\end{equation}
for all $\varphi \in H^s(\R^n)$ with $\varphi \ge 0$ and $\varphi \equiv 0$ in $\R^n\setminus \Omega$, when $K$ is absolutely continuous. 

For a general measure $K$, we have, instead of \eqref{eq:weak_sol} (analogously for \eqref{eq:weak_subsol}), 
\begin{equation}
\label{eq:weak_sol_meas}
\begin{split}
\frac12 \iint_{\R^n\times\R^n} \big(u(x) - u(x+y)\big)\big(\varphi(x)-\varphi(x+y)\big) K(dy)\, dx  =  \int_\Omega f \varphi.
\end{split}
\end{equation}

\end{defi} 

\begin{rem}
 In the previous definition, by density arguments it suffices to have \eqref{eq:weak_subsol} for all $\varphi\in C^\infty_c(\Omega)$. Indeed, if $\varphi\in H^s(\R^n)$ with $\varphi \equiv 0$ in $\R^n\setminus \Omega$, and $\varphi_\eps\in C^\infty_c(\Omega)$ is such that $\varphi_\eps\to \varphi$ in $H^s(\R^n)$, then 
\[
\big|\langle u, \varphi - \varphi_\eps\rangle_K\big|^2\le C \langle u, u\rangle_{K;\Omega} [\varphi-\varphi_\eps]^2_{H^s(\R^n)}\to 0 
\]
as $\varphi_\eps \to \varphi$ in $H^s(\R^n)$ (where we also used Lemma~\ref{lem:u_fourier}). 
\end{rem}

\begin{rem}\label{rem:fpwhy}
The assumption $p \ge \frac{2n}{n+2s}$ and $n > 2s$ is such that in \eqref{eq:weak_sol} the right-hand side is bounded:
\[
\int_\Omega f \varphi \le \|f\|_{L^p(\Omega)}\|\varphi\|_{L^{p'}(\Omega)} \le \|f\|_{L^p(\Omega)}[\varphi]_{H^{s}(\R^n)} 
\]
where we have used H\"older's inequality and the fractional Sobolev inequality, Theorem~\ref{thm:FSI}. 
\end{rem}
\begin{rem}
\label{rem:allp}
When $n = 2s$, functions in $H^s(\R^n)$ that vanish outside of $\Omega$ are in $L^{p'}(\Omega)$ for all $p' < \infty$ and when $n < 2s$, functions in $H^s(\R^n)$ vanishing outside of $\Omega$ are in $L^\infty(\Omega)$. 
In particular, the previous definition is still well-posed for $n = 2s$ and $p > 1$, and for $n < 2s$ and $p \ge 1$. Similarly, the theorem of existence of weak solutions below, Theorem~\ref{thm:exist_weak_sol}, is also valid in these cases. 
\end{rem}

If $u\in C^\infty_c(\R^n)$ is a weak solution then, by the integration by parts formula above, $\int \L u \, \varphi = \int f \varphi$ for all $\varphi \in C^\infty_c(\Omega)$. Hence it follows that $\L u = f$ in $\Omega$ and thus $u$ is a strong solution. More generally, we have the following: 
\begin{lem}
\label{lem:weak_strong}\index{Equivalent notions!Strong and weak}
Let $s\in (0, 1)$ and let $\L\in \GL$. Let $\Omega\subset\R^n$ be any bounded domain, and let $f\in L^p(\Omega)$ for $p \ge \frac{2n}{n+2s}$ and $n > 2s$. Let $u\in C^{2s+\eps}(\Omega)\cap L^\infty_{2s-\eps}(\R^n)$ for some $\eps > 0$ satisfying \eqref{eq:u_weak_sol}. 

Then, $\L u = f$ in $\Omega$ in the weak sense (see Definition~\ref{defi:weaksol}) if and only if it satisfies it in the strong sense.
\end{lem}
\begin{proof}
Let $\varphi\in C^\infty_c(\Omega)$, and we proceed as in the proof of Lemma~\ref{lem:integration_by_parts}. On the one hand, we have that $\langle u, \varphi\rangle_K$ is well-defined, using $\langle u, u\rangle_{K;\Omega}<\infty$. On the other hand, we can still split like in \eqref{eq:splitlike}, since we claim that each of the two terms appearing are integrable, $\int_{B_\eps^c\times \R^n} \varphi(x)(u(x+y)-u(y)) K(dy)\, dx<\infty$ (the other follows in the same way). 
To see this, we notice that we only integrate for $z\in \Omega$, where $u(z)$ is bounded (since $\Omega$ is bounded). Moreover, $|u(x+z)|\le C(1+|x|^{2s-\eps})$, and so the claim follows thanks to  \eqref{eq:prop_kernel3}. 

In all, we have 
\[
\begin{split}
 & \iint_{B_\eps^c\times \R^n}(u(x+y) - u(x))(\varphi(x+y)-\varphi(x)) K(dy)\, dx = \\
 & \qquad = \iint_{B_\eps^c\times \R^n}\varphi(x)(2u(x)-u(x-y)-u(x+y)) K(dy)\, dx.
 \end{split}
\]
We can take limits $\eps\downarrow 0$ in the first term because $\langle u, \varphi\rangle_K$ is well-defined (it is integrable), and in the second term thanks to Lemma~\ref{lem:Lu}. Hence, since $\L u$ is well-defined (again, by Lemma~\ref{lem:Lu}, in an $L^\infty$ sense) we have obtained
\[
\langle u, \varphi\rangle_K = \int_{\Omega} \L u \, \varphi\quad\text{for any}\quad\varphi\in C^\infty_c(\Omega).
\]
Now, if $u$ is a weak solution to $\L  u = f$ we get
\[
\int_{\Omega} f\varphi = \int_{\Omega}\L u \, \varphi\quad\text{for any}\quad\varphi\in C^\infty_c(\Omega),
\]
and hence, since $\L u$ exists and is bounded at every point, we deduce $f = \L u$ almost everywhere in $\Omega$. On the other hand, if $\L u = f$ a.e., then plugging it in the previous equation we obtain the weak formulation of the equation. 
\end{proof}

Observe that, since in the previous definitions we are asking $\varphi \equiv 0$ in $\R^n \setminus \Omega$, we have that the left-hand side in  \eqref{eq:weak_sol} (and in \eqref{eq:weak_subsol}) in reality can be integrated only in $(\R^n\times \R^n)\setminus (\Omega^c\times \Omega^c)$ (i.e., it is $\frac12 \langle u, \varphi\rangle_{K;\Omega}$). Thus the relaxed  requirement  \eqref{eq:u_weak_sol} instead of $u\in H^s(\R^n)$ in the definition of weak solution. This is relevant, in particular, when the boundary datum $g$ is not regular enough for it to be in the  energy space $H^s(\R^n)$ (or does not decay to zero at infinity).  


The energy functional associated to \eqref{eq:dirpb} is then\footnote{We note that the first term in \eqref{eq:energy_weak} is the analogue of the Dirichlet energy, $\frac12\int|\nabla u|^2$, in the local case $s = 1$ with $\L = -\Delta$. The constant $1/4$ appears here because the integration by parts in Lemma~\ref{lem:integration_by_parts} holds for the scalar product \eqref{eq:bilin}, which has a constant $1/2$ in front.}
\begin{equation}
\label{eq:energy_weak}
\begin{split}
\mathcal{E}(u) & = \frac12 \langle u, u\rangle_{K;\Omega} - \int_\Omega f (u-g)\\
&  = \frac14 \iint_{\R^n\times\R^n \setminus (\Omega^c\times \Omega^c)} \big|u(x) - u(z)\big|^2K(x-z)\, dx\, dz - \int_\Omega f(u-g),\\
\end{split}
\end{equation}
where the second equality holds for $K$ absolutely continuous. The unique minimizer of $\mathcal{E}$ among functions with $u = g$ in $\R^n\setminus \Omega$ will be the unique weak solution to \eqref{eq:dirpb} with $u\equiv g$ in $\R^n\setminus\Omega$. Observe that, when $g \equiv 0$ (or when $g\in H^s(\R^n)$), then we could integrate the first term in the whole $\R^n\times \R^n$. However, when $g\notin H^s(\R^n)$, the term $\iint_{\Omega^c\times \Omega^c} \big|g(x) - g(y)\big|^2 K(x-z)\, dx\, dz$ could be infinite, and this is why in general one has to take the previous expression.  On the other hand, when $\int_\Omega fg$ is finite, we can replace $\int_\Omega f(u-g)$ by simply $\int_\Omega f u$.

\subsection{Existence of weak solutions}
\label{ssec:existence_weak}

Let us now prove the existence theorem of weak solutions. We show it by means of the Riesz representation theorem in the appropriate space (see Remark~\ref{rem:allp} for the case $n\le 2s$).  

\begin{thm}[Existence of weak solutions]\index{Weak solutions!Existence}
\label{thm:exist_weak_sol}
Let $s\in (0, 1)$, let $\L\in \GL$, and let $\Omega\subset \R^n$ be any bounded domain. Let $f\in L^p(\Omega)$ for $p \ge \frac{2n}{n+2s}$ and $n > 2s$, and $g:\R^n \to \R$ be such that (recall Definition~\ref{defi:bilin_om})
\[
 \langle g, g\rangle_{K;\Omega}   < \infty.
\]
Then, there exists a unique weak solution $u$ of 
\[
\left\{
\begin{array}{rcll}
\L u & = & f & \quad\text{in}\quad \Omega,\\
u & = & g& \quad \text{in}\quad \R^n\setminus \Omega
\end{array}
\right.
\]
satisfying \eqref{eq:u_weak_sol}.  Moreover, $u$ is the unique minimizer of the energy $\mathcal{E}$, \eqref{eq:energy_weak}, among functions satisfying $u = g$ in $\R^n\setminus \Omega$.  
\end{thm}
\begin{proof}
Let us define 
\[
H^s_*(\Omega) := \left\{w \in H^s(\R^n) : w \equiv 0 \quad\text{in}\quad  \R^n\setminus \Omega\right\}.
\]
If $u$ is a weak solution, it satisfies \eqref{eq:weak_sol} and equivalently
\begin{equation}
\label{eq:weaksol2} 
\langle u - g, \varphi \rangle_K = \int_\Omega f \varphi - \langle g, \varphi \rangle_K = \int_\Omega f \varphi - \langle g, \varphi \rangle_{K;\Omega} 
\end{equation}
for all $\varphi \in H^s_*(\Omega)$. Since $\langle\cdot, \cdot\rangle_{K;\Omega}$ is a bilinear form, it satisfies the Cauchy-Schwarz inequality and hence (recall Lemma~\ref{lem:u_fourier})
\[
|\langle g, \varphi\rangle_{K;\Omega}|^2 \le |\langle g, g\rangle_{K;\Omega}||\langle \varphi, \varphi\rangle_{K}| \le C|\langle g, g\rangle_{K;\Omega}|[\varphi]^2_{H^s(\R^n)}.
\]
Therefore (see Remark~\ref{rem:fpwhy}), 
\begin{equation}
\label{eq:fromtheprevineqaaaa}
|\langle u - g, \varphi \rangle_K |\le C\left(\|f\|_{L^p(\Omega)} + \sqrt{|\langle g, g\rangle_{K;\Omega}|}\right) \|\varphi\|_{H^s(\R^n)}
\end{equation}
for all $\varphi \in H^s_*(\Omega)$. Observe now that, thanks to the Poincar\'e inequality, Proposition~\ref{prop:Poincare}, $\langle\cdot, \cdot\rangle_K$ defines a scalar product on $H^s_*(\Omega)$  (where we have $\langle\cdot, \cdot\rangle_K=\langle\cdot, \cdot\rangle_{K;\Omega}$) and hence it is a Hilbert space, since $\langle u, u\rangle_K = 0$ implies $u\equiv 0$.

In particular, since the norm induced by $\langle\cdot, \cdot\rangle_K$ on $H_*^s(\Omega)$ is comparable to $H^s(\R^n)$, we get from \eqref{eq:fromtheprevineqaaaa} that $u-g\in H^s(\R^n)$ with $u-g \equiv 0$ in $\R^n\setminus \Omega$. Thus, we can write $ u = g+v$ for some $v\in H^s_*(\Omega)$.

 By \eqref{eq:weaksol2} we are searching for a solution $v\in H^s_*(\Omega)$ to the weak formulation
\begin{equation}
\label{eq:crit_funct}
\langle v, \varphi\rangle_K  = \int_\Omega f\varphi - \langle g, \varphi\rangle_{K;\Omega},
\end{equation}
for all $\varphi\in H^s_*(\Omega)$ (recall that $\varphi\mapsto \langle g, \varphi\rangle_{K;\Omega}$ and $\varphi \mapsto \int_\Omega f\varphi $ are bounded and  linear operators on $H^s_*(\Omega)$). By the Riesz representation theorem, such solution exists and is unique, and hence $u = v+g$ is the unique weak solution to our problem.

Finally, observe that if $v = u - g$ is the unique solution of \eqref{eq:crit_funct}, then 
\[
\mathcal{E}(u+\varphi) - \mathcal{E}(u) = \mathcal{E}(\varphi) + \langle u, \varphi\rangle_{K; \Omega} -\int_{\Omega} f g = \mathcal{E}(\varphi)\ge 0,
\]
for all $\varphi\in H^s_*(\Omega)$, and hence $u$ is the unique minimizer of the energy. 
\end{proof}
 
We refer to \cite{KD, GH22} for a characterization of the functions $g:\R^n\setminus \Omega\to \R$ for which there exists an extension $\tilde g:\R^n\to \R^n$ (with $\tilde g = g$ in $\R^n\setminus \Omega$) such that $\langle \tilde g, \tilde g\rangle_{K;\Omega} < +\infty$, for kernels $K(y)\asymp |y|^{-n-2s}$. See also \cite{FKV} for a very general existence theorem, including non-symmetric and $x$-dependent kernels.

\subsection{Distributional solutions}
\index{Distributional solutions}
From the integration by parts formula, one can also define the following notion of solution to an equation $\L u = f$, requiring even less regularity on the function $u$  (see Remark~\ref{rem:relax_reg_dist}):

\begin{defi}[Distributional solution]
\label{defi:dist}
Let $s\in (0, 1)$ and $\L\in \GL$. Let $\Omega\subset \R^n$ be bounded, and let $f\in L^1_{\rm loc}(\Omega)$. We say that $u\in L^\infty_{2s-\eps}(\R^n)$ for some $\eps > 0$ is a \emph{distributional solution} of
\[
\L u = f \quad\text{in}\quad \Omega
\]
whenever
\[
\int_{\R^n} u \, \L \varphi = \int_{\Omega} f \varphi 
\qquad \textrm{for all} \qquad \varphi\in C^\infty_c(\Omega).
\]
More generally, we say that $u$ is a \emph{distributional supersolution} (resp. \emph{distributional subsolution}) to $\L u = f$ in $\Omega$, and we denote it $\L u \ge f$ in $\Omega$ (resp. $\L u \le f$ in $\Omega$) if
 \[
\int_{\R^n} u \, \L \varphi \ge \int_{\R^n} f \varphi \qquad\left(\text{resp.}\quad \int_{\R^n} u \, \L \varphi \le \int_{\R^n} f \varphi \right)
\]
for all $\varphi\in C^\infty_c(\Omega)$ with $\varphi \ge 0$. 
\end{defi}

\begin{rem}
Observe that   the previous definition is well-posed, since $u\in L^\infty_{2s-\eps}(\R^n)$ and hence $u\,\L \varphi\in L^1(\R^n)$ by Lemma~\ref{lem:welldefined} with $\alpha = \eps$.
\end{rem}

\begin{rem}[Relaxation of regularity]
\label{rem:relax_reg_dist}
To define \emph{distributional solutions} in Definition~\ref{defi:dist} we assume $u\in L^\infty_{2s-\eps}(\R^n)$ in order to force $u\,\L\varphi\in L^1(\R^n)$. This is, however, not necessary: as long as $u\,\L\varphi\in L^1(\R^n)$ for all $\varphi\in C^\infty_c(\R^n)$  we could make sense of the notion of  distributional solution  (and have the same properties). For example, we could ask for $u\in L^1(\R^n)$ or $u\in L^1(\R^n) + L^\infty_{2s-\eps}(\R^n) $ for some $\eps > 0$. Alternatively, we could define 
\[
L^1_\L (\R^n) := \left\{u \in L^1_{\rm loc}(\R^n) :  \sup_{\varphi\in C_c^\infty(B_1)} \int_{\R^n} |u \,\L \varphi|< +\infty\right\},
\]
and when $\L = \fls$ (or more generally, when $\L$ is comparable to the fractional Laplacian, see subsection~\ref{ssec:comparable}) then this space is exactly $L^1_{\omega_s}(\R^n)$ (recall Definition~\ref{defi:L1omegas}). 
\end{rem}

Distributional solutions retain some of the properties of their weak and strong counterparts. In particular, translation and scale invariance (Remark~\ref{rem:Scale_invariance}) follow from the definition itself. Observe, also, that smooth functions that are distributional solutions are strong solutions. This follows from the following lemma, saying that the integration by parts formula from Lemma~\ref{lem:integration_by_parts} works if  $u\in C_{\rm loc}^{2s+\eps}(\Omega)\cap L^\infty_{2s-\eps}(\R^n)$ (so that $\L u$ can be defined in the strong sense):  

\begin{lem}
\label{lem:int_parts_2}
Let $\Omega\subset \R^n$ be any bounded domain. Let  $v\in C^\infty_c(\Omega)$ and $u\in C^{2s+\eps}_{\rm loc}(\Omega)\cap L^\infty_{2s-\eps}(\R^n)$ for some $\eps > 0$. Then,
\[
\int_{\R^n} u \,\L v =   \int_{\Omega} \L u \, v.
\]
\end{lem}
\begin{proof}
We assume here, for notational simplicity, that $K(dy) = K(y)\, dy$. After a rescaling, let us assume $\Omega \subset B_1$. The proof follows  similarly to that of Lemma~\ref{lem:weak_strong}, by observing that the argument from Lemma~\ref{lem:welldefined}  also works to show that, for $v\in C^\infty_c(B_1)$, \footnote{If $K$ is a measure, one can define first the locally finite measure $\mu(d y) = \min\{1, |y|^2\} K(dy)$ (since \eqref{eq:Kint} holds), and then prove instead that 
\[
u(x)  \frac{2v(x) - v(x+y)-v(x-y)}{\min\{1,|y|^2\}}\in L^1(\R^n\times\R^n; dx\otimes \mu(dy)).
\]}
\[
w(x, y) = u(x) (2v(x) - v(x+y)-v(x-y)) K(y) \in L^1(\R^n\times\R^n).
\]
Indeed, 
\[
|w(x, y)|\le C (1+|x|^{2s-\eps})|2v(x) - v(x+y)-v(x-y)|K(y),
\]
and hence, noticing that 
\[
|2v(x) - v(x+y) - v(x-y)|\le C \min\{1, |y|^2\}\quad\text{for}\quad x\in B_2,  ~y\in \R^n,
\]
together with Tonelli's theorem and the symmetry of $K$, we get
\[
\begin{split}
\int_{\R^n}\int_{\R^n}|w(x, y)|\, dx\, dy & \le C \int_{B_2}(1+|x|^{2s-\eps})\, dx \int_{\R^n} \min\{1, |y|^2\} K(y)\, dy \\
& \quad  + C\int_{B_2^c}\int_{\R^n} |x|^{2s-\eps}|v(x+y)|K(y)\, dy\, dx,
\end{split}
\]
where the first term is now bounded (recall \eqref{eq:Kint}). For the second term, we use Tonelli's theorem (and the fact that $v$ is bounded) to bound it by
\[
C\int_{B_2^c}\int_{B_1} |x|^{2s-\eps}K(x-z)\, dz\, dx \le  C\int_{B_1^c}\int_{B_1} |z+\xi|^{2s-\eps}K(\xi)\, dz\, d\xi. 
\]
Using now that in this domain $|z+\xi|^{2s-\eps}\le C|\xi|^{2s-\eps}$ and by yet another application of Tonelli's theorem together with \eqref{eq:prop_kernel3} (with $\alpha = \eps$) we get that $w\in L^1(\R^n\times\R^n)$. In particular, we can apply Fubini's theorem to 
\[
\int_{\R^n} u\,\L v = \frac12 \iint_{\R^n\times\R^n} u(x) (2v(x) - v(x+y) -v(x-y) )K(y) \, dy\, dx
\]
and integrate in $x$ first (which is a well-defined integral, since $u$ is locally bounded and $v$ is compactly supported) to get
\[
\int_{\R^n} \hspace{-1mm} u(x) (2v(x) - v(x+y) -v(x-y) )\, dx \hspace{-0.5mm}=\hspace{-0.5mm} \int_{\R^n}\hspace{-1mm}  v(x) (2u(x) - u(x+y) -u (x-y) )\, dx.
\]
We can now multiply by $K(y)$ and integrate in $y$, and since $\L u$ is well-defined (by Lemma~\ref{lem:Lu}) we have 
\[
\int_{\R^n} u \,\L v = \int_{\R^n} \L u\, v ,
\]
as we wanted to see. 
\end{proof}

As a consequence, we obtain that distributional and strong solutions are the same for smooth functions: 

\begin{cor}
\label{cor:dist_strong_smooth}\index{Equivalent notions!Strong and distributional}
Let $s\in (0, 1)$,   $\L \in \GL$, and   $\Omega\subset \R^n$ be any bounded domain. Let $u\in C_{\rm loc}^{2s+\eps}(\Omega)\cap L^\infty_{2s-\eps}(\R^n)$ for some $\eps > 0$ and $f\in L^1_{\rm loc}(\Omega)$. 

Then, $u$ is a distributional solution to $\L u = f$ in $\Omega$ if and only if it is a strong solution.
\end{cor}

\begin{proof}
The fact that if $u$ is a strong solution then it is a distributional solution follows directly from Lemma~\ref{lem:int_parts_2}.

If $u$ is distributional, also from Lemma~\ref{lem:int_parts_2} we have 
\[
\int_{\R^n} u\L\varphi = \int_{\Omega} \L u \, \varphi = \int_{\Omega} f \varphi\qquad\text{for all}\quad \varphi \in C^\infty_c(\Omega).
\]
In particular, $\L u = f$ almost everywhere in $\Omega$, and $u$ is a strong solution as well. Notice that we also obtain that we must have $f\in L^\infty_{\rm loc}(\Omega)$, by Lemma~\ref{lem:Lu}.
\end{proof}

We observe that distributional solutions are also well behaved under convolution against a smooth mollifier. We denote by $\psi\in C^\infty_c(B_1)$ a smooth function with $\psi \ge 0$ and $\int_{\R^n} \psi = 1$, and we define 
\begin{equation}
\label{eq:mollifier}
\psi_\delta (x) := \delta^{-n}\psi\left(\frac{x}{\delta}\right).
\end{equation}
Then, we have:

\begin{lem}
\label{lem:conv_dist_sol}\index{Convolution}
Let $s\in (0, 1)$, $\L \in \GL$, and let $u\in L^\infty_{2s-\eps}(\R^n)$ for some $\eps > 0$ be a distributional solution to $\L u = f$ in $B_1$. Then the functions 
\[
u_\delta := u * \psi_\delta,\qquad\text{and}\qquad f_\delta := f * \psi_\delta,
\]
are $C^\infty(\R^n)$ and satisfy 
\begin{equation}
\label{eq:Lueps_d}
\L u_\delta = f_\delta\quad\text{in}\quad B_{1-\delta}
\end{equation}
in the strong sense. 
\end{lem}

\begin{proof}
First observe that if $u\in L^\infty_{2s-\eps}(\R^n)$ then $u_\eps\in L^\infty_{2s-\eps}(\R^n)$ as well. Hence, we can evaluate, for any $\varphi\in C^\infty_c(B_1)$,
\[
\begin{split}
\int_{\R^n} u_\eps \L \varphi&  = \int_{\R^n} \int_{\R^n}u(x-y)\psi_\eps(y) \, dy \,\L \varphi(x)\, dx \\
&  = \int_{\R^n} \int_{\R^n}u(x-y)\L \varphi(x)\, dx \,\psi_\eps(y) \, dy  \\
& =  \int_{\R^n} \int_{\R^n}f(x-y)\varphi(x)\, dx\, \psi_\eps(y) \, dy  \\
& =  \int_{\R^n} f_\eps \varphi.
\end{split} 
\]

Finally, since $u_\eps$ is smooth and it satisfies \eqref{eq:Lueps_d} in the distributional sense, by Corollary~\ref{cor:dist_strong_smooth}  it also satisfies it in the strong sense. 
\end{proof}
\begin{rem}\label{rem:conv_dist_sol}
In the whole space, the previous lemma says that if $\L u = f$ distributionally in $\R^n$, then $\L \tilde u = \tilde f$ in the strong sense in $\R^n$, where
\[
\tilde u = u * \zeta,\qquad \tilde f = f * \zeta,
\]
for  $\zeta\in C^\infty(\R^n)$ with $\zeta \ge 0$ and $\int_{\R^n}\zeta = 1$. 
\end{rem}

In particular, we have shown that distributional solutions can be approximated by strong solutions. 

Let us finish this subsection by showing that, when they are in the appropriate space, distributional and weak solutions are equivalent:

\begin{lem}
\label{lem:weak_distr}
Let $s\in (0, 1)$,  $\L\in \GL$,   $\Omega\subset \R^n$ be any bounded domain, and    $f\in L^p(\Omega)$ for $p \ge \frac{2n}{n+2s}$ and $n > 2s$. Let $u\in   L^\infty_{2s-\eps}(\R^n)$ for some $\eps > 0$ be satisfying \eqref{eq:u_weak_sol}. 

Then, $\L u = f$ in $\Omega$ in the weak sense (see Definition~\ref{defi:weaksol}) if and only if it satisfies it in the distributional sense (see Definition~\ref{defi:dist}). \index{Equivalent notions!Weak and distributional}
\end{lem}
\begin{proof}
We assume here, for notational simplicity, that $K(dy) = K(y)\, dy$. It is enough to show
\[
\langle u, \varphi\rangle_K = \int_{\R^n} u\, \L \varphi\qquad\text{for all}\quad\varphi\in C^\infty_c(\Omega).
\]
To do so, let us fix some $\varphi\in C^\infty_c(\Omega)$, and let 
\[
2\eta := \dist({\rm supp}\, \varphi, \R^n\setminus\Omega).
\]
Let us regularize $u$, by considering $u_\delta := u * \psi_\delta$ for some smooth mollifier as in \eqref{eq:mollifier} with $\psi\in C^\infty_c(B_1)$, and with $\delta \le \eta$. 

We denote 
\[
w(x, y, z) = \left(u(x-z)- u(y-z)\right)\psi_\delta(z)\left(\varphi(x)-\varphi(y)\right)K(x-y).
\]
Then $w\in L^1(D)$ where $D := \left[\R^{2n}\setminus(\Omega^c\times \Omega^c)\right]\times \R^n\subset \R^{3n}$. Indeed, if we define $\Omega_\eta := \{x\in \Omega : \dist(x, \R^n\setminus\Omega) \ge \eta\}$, we have
\[
\begin{split}
\|w\|_{L^1(D)}& \le \iint_{\R^{2n}\setminus(\Omega_\eta^c)^2}\int_{B_\delta} |w(x, y, z)|\, dz\, dx\, dy\\
& \le \int_{B_\delta}\psi_\delta(z) \langle u(\,\cdot\, - z), u(\,\cdot\, - z)\rangle_{K,\Omega_\eta}\, dz  + \int_{B_\delta} \psi_\delta(z) \langle \varphi, \varphi\rangle_{K}\, dz,
\end{split}
\]
where we have used $ab\le \frac12 a^2 + \frac12 b^2$ with $a = u(x-z)- u(y-z)$ and $b = \varphi(x)- \varphi(y)$, and the fact that $\varphi$ is supported in $\Omega_\eta$. Now, since $\Omega_\eta + z \subset \Omega$ for $z\in B_\delta$ and $\delta \le \eta$, by changing variables $x-z\mapsto x$ and $y-z\mapsto y$ we know 
\[
\langle u(\,\cdot\, - z), u(\,\cdot\, - z)\rangle_{K,\Omega_\eta}\le \langle u, u\rangle_{K,\Omega} < \infty.
\]
Together with $\langle \varphi, \varphi\rangle_{K} < \infty$ (see Lemma~\ref{lem:u_fourier}) and $\int\psi_\delta = 1$ we obtain that $w\in L^1(D)$. We can therefore apply Fubini to derive the following (also using that $\psi_\delta(-z) = \psi_\delta(z)$ for $z\in B_\delta$)
\begin{equation}
\label{eq:comb1d}
\langle u_\delta, \varphi\rangle_K = \langle u_\delta, \varphi\rangle_{K;\Omega_\eta} = \langle u, \varphi_\delta\rangle_{K;\Omega} = \langle u, \varphi_\delta\rangle_{K}.
\end{equation}

We have used again that $\varphi$ is supported in $\Omega_\eta$ and hence $\varphi_\delta$ is supported in $\Omega$ for $\delta < \eta$. 

Observe that $\langle u_\delta, u_\delta\rangle_{K;\Omega_\eta} \le \langle u, u\rangle_{K;\Omega}< \infty$ by H\"older's inequality (together with the previous arguments) and so $u_\delta$ is a weak solution in $\Omega_\eta$. Since it is also smooth, by Lemma~\ref{lem:weak_strong} it is also a strong solution, and 
\begin{equation}
\label{eq:comb1d2}
\langle u_\delta, \varphi\rangle_K = \int_{\R^n} \L u_\delta \, \varphi  =\int_{\R^n} u_\delta\, \L \varphi, 
\end{equation}
where in the last inequality we are using that $u_\delta$ is  a distributional solution as well (see Corollary~\ref{cor:dist_strong_smooth}). Finally, since
\[
\langle u, \varphi_\delta - \varphi\rangle_{K}\le C \langle u, u\rangle_{K;\Omega}[\varphi_\delta - \varphi]_{H^s(\R^n)}\downarrow 0\quad\text{as}\quad\delta\downarrow 0,
\]
combined with \eqref{eq:comb1d}-\eqref{eq:comb1d2} we obtain
\[
\langle u, \varphi\rangle_K = \lim_{\delta \downarrow 0}\,\langle u, \varphi_\delta \rangle_K = \lim_{\delta \downarrow 0} \int_{\R^n}u_\delta \, \L \varphi = \int_{\R^n}u \, \L \varphi.
\]
Since $\varphi$ is arbitrary, we get the desired result. 
\end{proof}

Finally, as a corollary we have that distributional, weak, and strong solutions are all equivalent notions for functions in the appropriate space:

\begin{cor}
Let $s\in (0, 1)$ and let $\L\in \GL$. Let $\Omega\subset\R^n$ be any bounded domain and let $f\in L^p(\Omega)$ for $p \ge \frac{2n}{n+2s}$ and $n > 2s$. Let $u\in C_{\rm loc}^{2s+\eps}(\Omega)\cap L^\infty_{2s-\eps}(\R^n)$ for some $\eps > 0$ be satisfying \eqref{eq:u_weak_sol}. 

Then, $\L u = f$ in $\Omega$ in the distributional sense (see Definition~\ref{defi:dist}) if and only if it satisfies it in the weak sense (see Definition~\ref{defi:weaksol}), and if and only if it satisfies it in the strong sense.
\end{cor}

\begin{proof}
It is now the result of combining Lemma~\ref{lem:weak_strong}, Corollary~\ref{cor:dist_strong_smooth}, and Lemma~\ref{lem:weak_distr}. 
\end{proof}
\subsection{Stability of distributional solutions}
\index{Stability!Distributional solutions}
The next result establishes the stability under $L^1$ convergence of distributional solutions $\L u = f$. 

In order to prove it, we use the following version of Vitali's convergence theorem, that we state here for the convenience of the reader (see \cite[Corollary 4.5.5]{Bog07}):
\begin{thm}[Vitali's Convergence Theorem]
\label{thm:vit_conv_thm}\index{Vitali's convergence theorem}
Let $g_k, g\in L^1(\R^n)$ with $k\in \N$ be such that $g_k(x) \to g(x)$ for a.e. $x\in \R^n$. Then, the convergence of $g_k$ to $g$ in $L^1(\R^n)$ (as $k\to\infty$) is equivalent to the following: 
\begin{equation}
\label{eq:cond1_vit}
\lim_{|E|\to 0}\sup_{k\in \N} \int_E |g_k|  =0
\end{equation}
and, for every $\eps > 0$, there exists a Borel set $E_\eps$ with $|E_\eps|< \infty$ such that 
\begin{equation}
\label{eq:cond2_vit}
\sup_{k\in \N}\int_{\R^n\setminus E_\eps}|g_k| < \eps. 
\end{equation}
\end{thm}

In the previous statement the first condition is immediately satisfied if $g_k$ are all uniformly bounded in $k$, and the second condition is the tightness of the measures $|g_k(x)|\, dx$.

We will also use Prokhorov's Theorem  (see \cite[Theorem 8.6.2]{Bog07} or \cite[Theorem 2.1.11]{FG21}).

\begin{thm}[Prokhorov's Theorem]\label{thm:Prokhorov}\index{Prokhorov's theorem}
Let $\mathcal{A}\subset\mathcal{P}(\R^n)$ be any family of probability measures.
Then, the following are equivalent:

\begin{itemize}
\item[(i)] For every sequence $(\mu_k)_{k\in \N}$ with $\mu_k\in \mathcal{A}$ there exists a subsequence $(\mu_{k_i})_{i\in \N}$ such that $\mu_{k_i} \rightharpoonup \mu$ as $i \to \infty$ for some $\mu\in \mathcal{P}(\R^n)$.

\item[(ii)] For every $\eps > 0$ there exists a compact set $K_\eps \subset \R^n$ such that $\mu(K_\eps) > 1-\eps$  for all $\mu\in \mathcal{A}$.
\end{itemize}
\end{thm}

Let us now prove the stability of distributional solutions:

\begin{prop}[Stability of distributional solutions]
\label{prop:stab_distr}\index{Compactness of linear operators!Distributional solutions}
Let $s\in (0, 1)$ and let, for each $k \in \N$, $\L_k \in \GL$. Let $\Omega\subset \R^n$ be any bounded domain. 

Let $u_k\in L^\infty_{2s-\eps}(\R^n)$ with $\sup_{k\in \N}\|u_k\|_{L^\infty_{2s-\eps}(\R^n)} < \infty$ for some $\eps > 0$ and $f_k \in L^1_{\rm loc}(\Omega)$ with  
\[
u_k \to u \quad\text{in}\quad L^1_{\rm loc}(\R^n)\qquad \text{and}\qquad f_k \to f \quad \text{in}\quad  L^1_{\rm loc}(\Omega),
\]
for some $u\in L^\infty_{2s-\eps}(\R^n)$ and $f\in L^1_{\rm loc}(\Omega)$, 
be such that 
\[
\L_k u_k = f_k \quad\text{in}\quad \Omega\qquad\left({\rm resp.}\ \ \L_k u_k \le f_k\ \ \text{in}\ \ \Omega\right)
\]
in the distributional sense. Then, there exists some $\L\in \GL$ such that 
\[
\L u = f\quad\text{in}\quad \Omega\qquad\left({\rm resp.}\ \ \L u \le f\ \ \text{in}\ \ \Omega\right)
\]
in the distributional sense. 
\end{prop}
\begin{proof}
Let $\varphi\in C^\infty_c(\Omega)$ arbitrary, so that
\[
\int_{\R^n} u_k \, \L_k \varphi = \int_\Omega f_k \varphi \qquad \left(\text{resp.\ \ $\varphi \ge 0$ \ and \ }  \int_{\R^n} u_k \, \L_k \varphi \le \int_\Omega f_k \varphi\right).
\]
Observe, first, that since $f_k \in L^1_{\rm loc}(\Omega)$ are converging to $f\in L^1_{\rm loc}(\Omega)$ in $L^1_{\rm loc}(\Omega)$, we immediately have 
\[
\int_\Omega f_k \varphi \to \int_\Omega f \varphi.  
\]
We want to show 
\begin{equation}
\label{eq:toshow_stab}
\int_{\R^n} u_k \, \L_k  \varphi \to \int_{\R^n} u \, \L \varphi
\end{equation}
for some $\L\in \GL$. In order to do that, let us consider the sequence of measures 
\[
\mu_k (dy) := c_k \min\{1, |y|^2\} K_k(dy),
\]
where $K_k$ are the L\'evy measures of $\L_k$ and $c_k$ are chosen so that $\mu_k$ is a probability measure, $\mu_k\in \mathcal{P}(\R^n)$. From the ellipticity conditions, we know that $c_k$ remains uniformly bounded and positive (recall that we are assuming \eqref{eq:Kint} and we have \eqref{eq:Kellipt_gen_l} and \eqref{eq:prop_kernel2}). Moreover, from \eqref{eq:prop_kernel}, for any $\delta > 0$ we can find $R \ge 1$ such that 
\[
\int_{B_R} \mu_k(dy) > 1-\delta,
\] 
where $R$ is independent of $k$. 
By  Prokhorov's Theorem (see Theorem~\ref{thm:Prokhorov}),  there exists a subsequence $\mu_{k_m}$ converging weakly $\mu_{k_m}\rightharpoonup \mu$ for some $\mu\in \mathcal{P}(\R^n)$. In particular, up to a further subsequence if necessary, we can assume $c_k \to c_\infty > 0$ and we can define a measure
\[
K(dy) = \frac{\mu(dy)}{c_\infty\min\{1, |y|^2\}},
\]
that will satisfy both properties \eqref{eq:Kellipt_gen_L}-\eqref{eq:Kellipt_gen_l} in the limit by weak convergence outside of the origin (in fact, it will satisfy \eqref{eq:Kellipt_gen_L}-\eqref{eq:Kellipt_gen_l} for some measure in the limit). 

Hence, for any $x\in \R^n$ (and using that $\varphi$ is smooth), 
\[
\begin{split}
\L_k \varphi (x) & = \frac12 \int_{\R^n}\left(2\varphi(x) - \varphi(x+y) - \varphi(x-y) \right) K_k(dy) \\
& = \frac12 \int_{\R^n}\frac{2\varphi(x) - \varphi(x+y) - \varphi(x-y)}{c_k\min\{1, |y|^2\}} c_k\min\{1, |y|^2\}K_k(dy) \\
& \to \frac12 \int_{\R^n}\left(2\varphi(x) - \varphi(x+y) - \varphi(x-y)\right) K(dy)
\end{split} 
\]
and $\L_k \varphi(x) \to \L \varphi(x)$ pointwise for all $x\in \R^n$. We want to apply Vitali's Convergence Theorem (Theorem~\ref{thm:vit_conv_thm}) to $g_k = u_k \L_k \varphi$ which satisfies $g_k(x) \to g(x) = u(x) \L\varphi(x)$ for a.e. $x\in \R^n$, and $g_k, g\in L^1(\R^n)$ by Lemma~\ref{lem:welldefined} (with bounds independent of $k\in \N$).

Observe that, on the one hand, we have by Lemma~\ref{lem:welldefined}
\[
|g_k(x)| = \big|u_k(x)\L_k \varphi(x)\big|\le C\big(1+|x|^{2s-\eps}\big)(1+|x|)^{-2s} \le C\quad\text{for all}\quad x\in \R^n,
\]
for some constant $C$ independent of $k$. That is, condition \eqref{eq:cond1_vit} holds, since $u_k\L_k\varphi$ is uniformly bounded. 

On the other hand, and again thanks to Lemma~\ref{lem:welldefined}, we also have that 
\[
\int_{\R^n\setminus B_r} |g_k|\le C\int_{\R^n\setminus B_r} \big(1+|x|^{2s-\eps}\big)|\L_k \varphi| \le C (1+r)^{-\eps}\quad\text{for any}\quad r > 0,
\]
and some $C$ independent of $k$. That is, condition \eqref{eq:cond2_vit} also holds. 

Hence, we can apply Theorem~\ref{thm:vit_conv_thm} to the sequence $g_k$ to deduce that $g_k\to g$ in $L^1(\R^n)$ and \eqref{eq:toshow_stab} is satisfied, as we wanted to see.
\end{proof}

\section{Maximum principle}
\label{sec:max_principle}
 
Super- and subsolutions to equations of the form $\L u=0 $ satisfy a maximum principle and a comparison principle. 

As for the classical case of the Laplacian $-\Delta$, the maximum principle essentially relies on the fact that $\L u(x_\circ) \ge 0$ whenever $u$ has a maximum at $x_\circ$. For local equations this is true for any local maximum, while in the nonlocal setting this is only true when $u$ has a global maximum at $x_\circ$. In this case, the inequality $\L u(x_\circ) \ge 0$ follows simply from the integral expression of $\L$ and the fact that $u(x_\circ)\ge u(x_\circ\pm  y)$ for any\footnote{When $K > 0$ in $\R^n$, such argument already yields a strong maximum principle, since one has $\L u(x_\circ) > 0$ unless $u$ is constant. However, for $\L \in \GL$, the strong maximum principle turns out to be more delicate, see subsection~\ref{ssec:strongmaximumprinciple}.} $y \in \R^n$.

The previous considerations work when $u$ is regular enough, so that $\L u$ can be evaluated (see Lemma~\ref{lem:Lu}). In the following, we first show how the proof for strong solutions works, and then we prove the maximum principle for weak solutions, and for continuous distributional solutions as well. 

\subsection{Maximum principle for strong solutions}

We start by proving the maximum principle and its consequences for strong solutions with general integro-differential operators. The proofs here follow similarly to those for the fractional Laplacian (see Section~\ref{ssec:max_sqrtl}).

\begin{lem}[Maximum Principle for strong solutions]
\label{lem:max_principle_G}\index{Maximum principle!General operators!Strong solutions}
Let $s\in (0, 1)$ and   $\L\in \GL$. Let $\Omega\subset \R^n$ be any bounded open set, and   $u\in C_{\rm loc}^{2s+\eps}(\Omega)\cap C(\overline{\Omega})\cap L^\infty_{2s-\eps}(\R^n)$ for some $\eps>0$. Let us assume that 
\[
\left\{
\begin{array}{rcll}
\L u & \ge & 0 & \quad\text{in}\quad \Omega,\\
u & \ge & 0 & \quad\text{in}\quad \R^n\setminus \Omega.
\end{array}
\right.
\]
Then, $u\ge 0$ in $\R^n$. 
\end{lem}
\begin{proof}
Let us suppose that it is not true, and that $\inf_{\R^n} u = u(x_\circ) < 0$ for some $x_\circ\in \Omega$ (which is achieved inside $\Omega$ by continuity of $u$ and the fact that $u\ge 0$ in $\R^n\setminus \Omega$). 

Since $K\not\equiv 0$, there exists some $\be \in \mathbb{S}^{n-1}$ and $\eps > 0$ such that 
\begin{equation}
\label{eq:nondegcond}
\int_{B_{\eps/2}(\eps \be)} K(dy) > 0
\end{equation}
 (in other words, take $\eps \be \in {\rm supp}(K)$).

Since $x_\circ\in \Omega$ is a minimum, 
\[
\L u (x_\circ) = \textrm{P.V.}\int_{\R^n} \big(u(x_\circ) - u(x_\circ+y)\big) K(dy)  \le 0.
\]
Now, if $\L u(x_\circ) < 0$, we get a contradiction. Otherwise, from \eqref{eq:nondegcond} there exists some $x_1 \in x_\circ + B_{\eps/2}(\eps \be)$ such that $u(x_1) = u(x_\circ) < 0$. In particular, $x_1\in \Omega$ and $\L u(x_1) \le 0$. 

We continue recursively: either $\L u(x_1) < 0$ and we get a contradiction, or we can find $x_2 \in x_1 + B_{\eps/2}(\eps \be)$ such that $u(x_2) = u(x_\circ) < 0$. After finitely many steps, we will get either a point $x_i\in \Omega$ such that $\L u(x_i) < 0$, or $x_i \not\in \Omega$ (since $\Omega$ is bounded), in both cases getting a contradiction.
\end{proof}

As a consequence we obtain the comparison principle for strong solutions. 
\begin{cor}[Comparison Principle for strong solutions]
\label{cor:comp_principle_G}\index{Comparison principle!General operators!Strong solutions}
Let $s\in (0, 1)$ and let $\L\in \GL$. Let $\Omega\subset \R^n$ be any bounded open set, and let $u_1, u_2 \in C_{\rm loc}^{2s+\eps}(\Omega)\cap C(\overline{\Omega})\cap L^\infty_{2s-\eps}(\R^n)$ for some $\eps>0$. Let us assume that 
\[
\left\{
\begin{array}{rcll}
\L u_1 & \ge & \L u_2 & \quad\text{in}\quad \Omega,\\
u_1 & \ge & u_2 & \quad\text{in}\quad \R^n\setminus \Omega.
\end{array}
\right.
\]
Then $u_1\ge u_2$ in $\R^n$. 
\end{cor}
\begin{proof}
This is an immediate consequence of Lemma~\ref{lem:max_principle_G}, applied to $w := u_1 - u_2$. 
\end{proof}

As a consequence of the previous result we also obtain the uniqueness of strong solutions, as in Corollary~\ref{cor:uniqueness_sqrtl}.

\subsection{Maximum principle for weak solutions} Let us now prove the maximum principle for weak solutions. 

\begin{lem}[Maximum Principle for weak solutions]
\label{lem:max_principle_G_w}\index{Maximum principle!General operators!Weak solutions}
Let $s\in (0, 1)$ and let $\L\in \GL$. Let $\Omega\subset \R^n$ be any bounded open set, and let $u$ be a weak supersolution $\L u \ge 0$ (see Definition~\ref{defi:weaksol}) with $u \ge 0$ in $\R^n\setminus \Omega$. 

 Then, $u \ge 0$ in $\R^n$. 
\end{lem}
\begin{proof}
By the definition of a weak supersolution we know that (recall \eqref{eq:bilin})
\[
\langle u, \varphi\rangle_K \ge 0,
\]
for all $\varphi \in H^s(\R^n)$ with $\varphi \equiv 0$ in $\R^n\setminus \Omega$ and $\varphi \ge 0$ in $\R^n$. Let us denote $u = u^+-u^-$ where $u^+= \max\{u, 0\}$ and $u^- = \max\{-u, 0\}$ are the positive and negative parts of $u$ respectively. Take $\varphi (x) = u^-(x)$, which by assumption satisfies $u^-(x) = 0$ in $\R^n\setminus\Omega$. Moreover, since $|a^- - b^-|\le |a-b|$ (recall Definition~\ref{defi:bilin_om}) 
\[
\begin{split}
 [u^-]_{H^s(\R^n)}^2& \asymp \langle u^-, u^-\rangle_K   = \langle u^-, u^-\rangle_{K;\Omega}    \le  \langle u, u\rangle_{K;\Omega} < \infty,
 \end{split}
\]
where we have used Lemma~\ref{lem:u_fourier} and that $u$ satisfies \eqref{eq:u_weak_sol}. 

Using now $\left(u^+(x) - u^+(z)\right)\left(u^-(x) - u^-(z)\right) \le 0$, gives
\[
\begin{split}
0 & \le - \langle u^-, u^-\rangle_K = - \iint_{\R^n\times\R^n} \big|u^-(x+y) - u^-(x)\big|^2K(dy)\, dx.
\end{split}
\]

In particular,  for a.e. $y\in {\rm supp}(K)$  and for a.e. $x\in \R^n$, $u^-(x) = u^-(x+ y)$. Since $K\not \equiv 0$ (and $K\neq \delta_0$) and $u^- = 0$ in $\R^n\setminus \Omega$, this implies $u^- \equiv 0$ in $\R^n$ and hence $u \ge 0$ in $\R^n$, as we wanted to see.
\end{proof}

We also obtain the comparison principle:

\begin{cor}[Comparison Principle for weak solutions]
\label{cor:comp_principle_G_w}\index{Comparison principle!General operators!Weak solutions}
Let $s\in (0, 1)$ and let $\L\in\GL$. Let $\Omega\subset \R^n$ be any bounded open set, and let $u_1$ and $u_2$ be such that
\[
\left\{
\begin{array}{rcll}
\L u_1 & \ge & \L u_2 & \quad\text{in}\quad \Omega,\\
u_1 & \ge & u_2 & \quad\text{in}\quad \R^n\setminus \Omega,
\end{array}
\right.
\]
in the weak sense. Then $u_1\ge u_2$ in $\R^n$. 
\end{cor}
\begin{proof}
Follows from Lemma~\ref{lem:max_principle_G_w} applied to $w := u_1 - u_2$. 
\end{proof}

Again, from the previous result we obtain the uniqueness of weak solutions.

\subsection{Maximum principle for continuous distributional solutions} 

Finally, let us prove a maximum principle for distributional solutions. 
In this case, we need to add the condition that the functions are continuous up to the boundary of the domain:
 
\begin{lem}[Maximum Principle for distributional solutions]
\label{lem:max_principle_G_d}\index{Maximum principle!General operators!Distributional solutions}
Let $s\in (0, 1)$ and   $\L\in \GL$. Let $\Omega\subset \R^n$ be any bounded open set, and   $u$ be a supersolution in the distributional sense $\L u \ge 0$ (see Definition~\ref{defi:dist}) with $u \ge 0$ in $\R^n\setminus \Omega$. Let us assume, moreover, that $u\in C(\overline{\Omega})$. Then $u \ge 0$ in~$\R^n$. 
\end{lem}
\begin{proof}
Let us argue by contradiction, and let us assume that $u < 0$ somewhere in $\Omega$. Let $-A = \min_{x\in \Omega} u(x) < 0$, and let $\Omega_{A} := \{x \in \Omega : u(x) < -A/2\}$. Then $\Omega_A \ssubset \Omega$ and $\dist(\Omega_A, \R^n\setminus \Omega) =: \delta_A > 0$ (since $u$ is continuous). 

Let $\psi\in C^\infty_c(B_1)$ be a mollifier ($\psi \ge 0$, $\int_{\R^n}\psi  = 1$), and let $\psi_\eps(x) = \eps^{-n}\psi(x/\eps)$ be its rescalings. Then, choose $\eps \le \delta_A$ such that $\min_{x\in \Omega} u_\eps \le -5A/6$, and $u_\eps \ge -3A/4$ in $\R^n\setminus\Omega_A$ where $u_\eps = u*\psi_\eps$. Then $u_\eps$ is a function such that $u_\eps \ge -3A/4$ in $\R^n\setminus\Omega_A$, $\L u_\eps \ge 0$ in $\Omega_A$ (by Lemma~\ref{lem:conv_dist_sol}), but $u \le -5A/6$ somewhere in $\Omega_A$. Since now $u_\eps$ is smooth, this contradicts the maximum principle for strong solutions, Lemma~\ref{lem:max_principle_G}, applied to $u_\eps - 3A/4$. 
\end{proof}

\begin{rem}
As we will see in Section~\ref{sec:int_reg_G}, when $f\in L^p$ and $p > \frac{n}{2s}$, distributional solutions are H\"older continuous in $\Omega$, so in that case the assumption $u\in C(\overline{\Omega})$ is a condition only on the boundary $\partial\Omega$. 
\end{rem}

\begin{rem}
The continuity assumption in $\overline{\Omega}$ is very important, otherwise there are distributional solutions like 
\[
u(x) = \big(1-|x|^2\big)^{s-1}_+\quad\text{in}\quad \R^n
\]
which satisfy  $\fls u = 0$ in $B_1$  and $u = 0$ in $\R^n\setminus B_1$, but $u\not\equiv 0$. In this case, $u$ is a distributional solution (in the relaxed sense of Remark~\ref{rem:relax_reg_dist}, $u\in L^1(\R^n)$), but $u\notin H^s(\R^n)$, so it is not a weak solution in $B_1$.

  The previous function $u$ is analogous to the solution $u = \chi_{\Omega}$ for $s = 1$, which satisfies $\Delta u = 0$  inside $\Omega$, $u = 0$ on $\partial\Omega$, but $u\not\equiv 0$ in $\Omega$. 
\end{rem}

Again we have a comparison principle as well: 
\begin{cor}[Comparison Principle for distributional solutions]
\label{cor:comp_principle_G_d}\index{Comparison principle!General operators!Distributional solutions}
Let $s\in (0, 1)$ and $\L\in\GL$. Let $\Omega\subset \R^n$ be any bounded open set, and let $u_1$ and $u_2$ be such that
\[
\left\{
\begin{array}{rcll}
\L u_1 & \ge & \L u_2 & \quad\text{in}\quad \Omega,\\
u_1 & \ge & u_2 & \quad\text{in}\quad \R^n\setminus \Omega,
\end{array}
\right.
\]
in the distributional sense. Assume, moreover, that $u_1, u_2\in C(\overline{\Omega})$. Then, $u_1\ge u_2$ in $\R^n$. 
\end{cor}
\begin{proof}
Follows as the proof of Corollary~\ref{cor:comp_principle_G_w}.
\end{proof}

\subsection{$L^\infty$ bounds}
\label{ssec:Linftybounds}\index{Linfty bounds@$L^\infty$ bounds!General operators}
As a direct consequence of the maximum principle for weak solutions, and by means of a barrier argument, we next prove an $L^\infty$ bound for solutions to the Dirichlet problem with bounded right-hand side:

\begin{lem}
\label{lem:Linftybound}
Let $s\in (0,1)$ and $\L\in \GL$. Let $\Omega\subset \R^n$ be any bounded open set, $f\in L^\infty(\Omega)$, and $g\in L^\infty(\R^n\setminus \Omega)$. 

Let $u$ be a weak solution to 
\[
\left\{
\begin{array}{rcll}
\L u & = & f & \quad\text{in}\quad \Omega,\\
u & = & g & \quad\text{in}\quad \R^n\setminus \Omega.
\end{array}
\right.
\]
Then, $u\in L^\infty(\Omega)$ with 
\[
\|u\|_{L^\infty(\Omega)} \le \|g\|_{L^\infty(\R^n\setminus \Omega)} + C \|f\|_{L^\infty(\Omega)},
\]
for some constant $C$ depending only on $n$, $s$, $\lambda$, $\Lambda$, and ${\rm diam}(\Omega)$. 
\end{lem}

In order to prove it, we use the following  simple barrier:
\begin{lem}
\label{lem:barrierLinfty}
Let $s\in (0,1)$ and $\L\in \GL$. Let $\Omega\subset \R^n$ be any bounded open set. Then, there exists a function $w\in C^\infty_c(\R^n)$ such that 
\[
\left\{
\begin{array}{rcll}
\L w & \ge & 1 & \quad\text{in}\quad \Omega,\\
w & \ge & 0 & \quad\text{in}\quad \R^n\setminus \Omega,\\
w & \le & C & \quad\text{in}\quad \Omega,
\end{array}
\right.
\]
for some constant $C$ depending only on $\lambda$, $\Lambda$, $n$, $s$, and ${\rm diam}(\Omega)$. 
\end{lem}

\begin{proof}
Let $B_R$ be any large ball such that $\Omega\subset B_R$, and let 
\[
\eta(x):= \big(M^2-|x|^2\big)_+,
\]
with $M=3R$.
Observe that, for each $x\in B_R$, we will have 
\[
2\eta(x) - \eta(x+y)-\eta(x-y) \ge |y|^2 \quad \textrm{for}\quad |y|<2R,
\]
while
\[
2\eta(x) - \eta(x+y)-\eta(x-y) \ge 0 \quad \textrm{for}\quad |y|>2R,
\]
since $|x|<R$ and $|x\pm y|>R$.
Therefore, 
\[
\L \eta (x) \ge \frac12\int_{B_{2R}} |y|^2 K(dy)\geq c>0,
\]
where we have used the lower bound \eqref{eq:Kellipt_gen_l} on ellipticity.
The constant $c$ depends only on $\lambda$, $n$, $s$, and $R$, which in turn can be taken to be (after a translation) $R = 2\,{\rm diam}(\Omega)$. 

Hence, by taking $w = \frac{1}{c}\eta$, we obtain the desired result, with $C = \frac{1}{c}$. 
\end{proof}

Thanks to the previous barrier and by the comparison principle, we obtain: 
\begin{proof}[Proof of Lemma~\ref{lem:Linftybound}]
Let $v(x) = \|g\|_{L^\infty(\R^n\setminus \Omega)} + \|f\|_{L^\infty(\Omega)} w(x)$, where $w$ is given by Lemma~\ref{lem:barrierLinfty}. Then, we clearly have $\L u \le \L v$ in $\Omega$, and $g \le v$ in $\R^n\setminus \Omega$. Thus, by the comparison principle, Corollary~\ref{cor:comp_principle_G_w}, we have $u \le v$ in $\Omega$, and in particular, $u \le \|g\|_{L^\infty(\R^n\setminus \Omega)} + C\|f\|_{L^\infty(\Omega)}$ in $\Omega$. 

Applying the same argument to $-u$, the result follows. 
\end{proof}

More generally, when $f\in L^p(\Omega)$ with $p > \frac{n}{2s}$ and $g\in L^\infty_{2s-\eps}(\R^n\setminus\Omega)$, we also have an $L^\infty$ bound for weak solutions. 

\begin{prop}[$L^\infty$ bound with $L^p$ right-hand side]
\label{prop:Lpbound}
Let $s\in (0,1)$ and $\L\in \GL$. Let $\Omega\subset \R^n$ be any bounded open set, $f\in L^p(\Omega)$ for some $ \frac{n}{2s}< p \le \infty$ with $p \ge 1$, and $g\in L^\infty_{2s-\eps}(\R^n)$ for some $\eps > 0$ satisfying \eqref{eq:u_weak_sol}. Let $u$ be a weak solution of 
\[
\left\{
\begin{array}{rcll}
\L u & = & f & \quad\text{in}\quad \Omega,\\
u & = & g & \quad\text{in}\quad \R^n\setminus \Omega.
\end{array}
\right.
\]
Then, $u\in L^\infty(\Omega)$ with 
\[
\|u\|_{L^\infty(\Omega)} \le C\left( \|g\|_{L^\infty_{2s-\eps}(\R^n)} +  \|f\|_{L^p(\Omega)}\right),
\]
for some constant $C$ depending only on $n$, $s$, $\lambda$, $\Lambda$, $p$, $\eps$, and ${\rm diam}(\Omega)$. 
\end{prop}

We first prove the following intermediate lemma:
\begin{lem}
\label{lem:asdffdsa}
Let $T > 0$ and $\beta\ge 1$. Then,
\[
\left|a_T^\beta a - b_T^\beta b\right|^2\le (1+\beta)^2 (a-b) \big(a_T^{2\beta} a - b_T^{2\beta}b\big)\quad\text{for all}\quad a, b \in \R,
\]
where we denote $a_T := \min\{|a|, T\}$ for any $a\in \R$.
\end{lem}
\begin{proof}
Let $f_\beta(t) :=  t_T^\beta t$, which satisfies $f'_\beta(t) = T^\beta$ if $|t|> T$ and $f_\beta'(t) = (\beta + 1)|t|^\beta$ if $|t| < T$. In particular
\begin{equation}
\label{eq:leminttonto}
t_T^\beta \le f'_\beta(t) \le (\beta + 1) t_T^\beta. 
\end{equation}
Using now 
\[
\begin{split}
|f_\beta(a)- f_\beta(b)|^2 & = \left(\int_a^b f'_\beta(t)\, dt \right)^2 \le |a-b|\left|\int_a^b (f'_\beta)^2(t)\, dt\right|
\end{split}
\]
Together with \eqref{eq:leminttonto} we get
\[
\begin{split}
|f_\beta(a)- f_\beta(b)|^2 & \le (\beta+1)^2|a-b|\left|\int_a^b f'_{2\beta}(t) \, dt \right|\\
& =(\beta+1)^2 |a-b||f_{2\beta}(a) - f_{2\beta}(b)|,
\end{split}
\]
which is the desired inequality. 
\end{proof}

We can now give the proof of the $L^\infty$ bound with $L^p$ right-hand side:

\begin{proof}[Proof of Proposition~\ref{prop:Lpbound}]
Let us assume, for notational simplicity only, that $K(dy) = K(y)\, dy$. Let $R > 0$ such that $\Omega\subset B_{R/2}$ (after a translation we can take $ R = {\rm diam}(\Omega)$). Let us then consider
\[
g_1:= g \chi_{\R^n\setminus B_R}\qquad\text{and}\qquad g_2 = g - g_1 = g \chi_{B_R}.
\]
In particular, $\|\L g_1 \|_{L^\infty(\Omega)}\le C \|g\|_{L^\infty_{2s-\eps}(\R^n )}$ by Lemma~\ref{lem:Lu}. If we define $u_1 = u - g_1$, then $u_1$ satisfies 
\[
\left\{
\begin{array}{rcll}
\L u_1 & = & h & \quad\text{in}\quad \Omega,\\
u_1 & = & g_2 & \quad\text{in}\quad \R^n\setminus \Omega
\end{array}
\right.
\]
in the weak sense (notice that since $g$ satisfies \eqref{eq:u_weak_sol}, we have that $g_2$ satisfies \eqref{eq:u_weak_sol} as well), where $h = f - \L g_1$ is such that $\|h\|_{L^p(\Omega)}\le \|f\|_{L^p(\Omega)}+  C \|g\|_{L^\infty_{2s-\eps}(\R^n )}$, and $\|g_2\|_{L^\infty(\R^n )}\le C \|g\|_{L_{2s-\eps}^\infty(\R^n\setminus\Omega)}$. Observe that $u = u_1$ in $\Omega$. 

Let us split $u_1 = v+w$ with $v$ a solution to 
\[
\left\{
\begin{array}{rcll}
\L v & = & h & \quad\text{in}\quad \Omega,\\
v & = & 0 & \quad\text{in}\quad \R^n\setminus \Omega.
\end{array}
\right.
\]
given by Theorem~\ref{thm:exist_weak_sol} (together with Remark~\ref{rem:allp} for the case $n \le 2s$), and $w = u - v$. By Lemma~\ref{lem:Linftybound} we already have an $L^\infty$ bound for $w$ (since $\L w = 0$ in $\Omega$ and $w = g_2$ is bounded in $\R^n\setminus \Omega$), so let us find an $L^\infty$ bound for $v$ as well. After dividing by a constant and rescaling the domain, we can assume that
\begin{equation}
\label{eq:assume_f_Om}
\|h\|_{L^p(\Omega)}\le 1\quad\text{and}\quad |\Omega|\le 1.
\end{equation}

Let $\beta \ge 0$, $T\ge 1$, and let us denote $v_T := \min\{|v|, T\}$, so that by Lemma~\ref{lem:asdffdsa} for all $x, z \in \R^n$, 
\[
\frac{\big|v_T^\beta(x) v(x) - v_T^\beta(z)v(z)\big|^2}{(1+\beta)^2}\le \left(v(x) - v(z)\right)\left(v^{2\beta}_T(x) v(x) - v_T^{2\beta}(z) v(z)\right).
\]
Thus, integrating against $K(dy)\, dx$ with $z= x+y$ we get
\[
\begin{split}
&\frac{1}{(1+\beta)^2} \int_{\R^n}\int_{\R^n}\big|v_T^\beta(x) v(x) - v_T^\beta(x+y)v(x+y)\big|^2\, K(dy)\, dx\leq \\
& ~~~  \le  \int_{\R^n}\int_{\R^n}\left(v(x) - v(x+y)\right)\left(v^{2\beta}_T(x) v(x) - v_T^{2\beta}(x+y) v(x+y)\right)K(dy)\,  dx.
\end{split}
\]
Since $v\in H^s(\R^n)$ with $v \equiv 0$ in $\R^n\setminus \Omega$, we have $v_T^{2\beta} v\in H^s(\R^n)$  with $v_T^{2\beta}v  \equiv 0 $ in $\R^n\setminus\Omega$ as well. Notice that the second term in the inequality above is the weak formulation of $\L v = h$ in $\Omega$, \eqref{eq:weak_sol_meas}, taking $\varphi = v_T^{2\beta} v$. On the other hand, the first term is comparable to $[v_T^\beta v]_{H^s(\R^n)}^2$ by Plancherel's theorem (see Lemma~\ref{lem:u_fourier}). In all, we have 
\[
[v_T^\beta v]_{H^s(\R^n)}^2 \le C(1+\beta)^2\int_{\R^n} h v_T^{2\beta} v
\]
for some constant $C$ depending only on $n$, $s$, and $\lambda$ (and independent of~$\beta$). 
By the fractional Sobolev inequality, Theorem~\ref{thm:FSI}, and H\"older's inequality with $1 = \frac{1}{p} + \frac{1}{p'}$, we have 
\[
\|v_T^\beta v\|^2_{L^q(\Omega)} \le C (1+\beta)^2 \|h\|_{L^p(\Omega)} \|v_T^{2\beta}v \|_{L^{p'}(\Omega)},
\]
where $q = \frac{2n}{n-2s}$. If we assume that the right-hand side is uniformly bounded in $T$, we can then let $T\uparrow \infty $ to obtain (recall \eqref{eq:assume_f_Om})
\[
\|v^{1+\beta}\|^2_{L^q(\Omega)} \le C (1+\beta)^2 \|v^{1+2\beta} \|_{L^{p'}(\Omega)},
\]
whenever the upper bound is finite. Denoting $\gamma = (1+\beta) q$ we can rewrite the previous inequality as follows (for a possibly different $C$):
\[
\|v\|_{L^\gamma(\Omega)} \le C^{\frac{1}{\gamma }}\gamma^{\frac{q}{\gamma}} \|v\|_{L^{\frac{2\gamma-q}{q}p'}(\Omega)}^{\frac{2\gamma-q}{2\gamma}}.
\]
By H\"older's inequality, we have
\[
\|v\|^{\frac{2\gamma-q}{2\gamma}}_{L^{\frac{2\gamma-q}{q}p'}(\Omega)}\le \|v\|^{\frac{q}{2\gamma p'}}_{L^{\frac{2\gamma}{q}p'}(\Omega)}\le \max\left\{1, \|v\|_{L^{\frac{2\gamma}{q}p'}(\Omega)}\right\},
\]
where we have also used that $q \le 2\gamma p'$. In all, for any $\gamma \ge q$, 
\[
\|v\|_{L^\gamma(\Omega)} \le C^{\frac{1}{\gamma }}\gamma^{\frac{q}{\gamma}} \max\left\{1,\|v\|_{L^{ \gamma/ \sigma}(\Omega)}\right\},
\]
where $\sigma = \frac{q}{2p'} > 1$, since $p > \frac{n}{2s}$. We apply now the previous inequality recursively, starting with $\gamma = q$: for any $m\in \N$, 
\[
\|v\|_{L^{\sigma^m q}(\Omega)} \le C^{\frac{1}{\sigma^m q}} (\sigma^m q)^{{\sigma^{-m}}} \max\left\{1,\|v\|_{L^{\sigma^{m-1}q}(\Omega)}\right\}.
\]
Observe that for $m = 0$ we already have that, since $v\in H^s(\R^n)$, by the fractional Sobolev inequality $v\in L^q(\Omega)$. After dividing by a constant we can assume $\|v\|_{L^q(\Omega)}\le 1$ and hence we get that
\[
\|v\|_{L^{\sigma^m q}(\Omega)} \le \prod_{i = 0}^m C^{\frac{1}{\sigma^i q}} (\sigma^i q)^{{\sigma^{-i}}} = \left(qC^{\frac{1}{q}}\right)^{\sum_{i = 0}^m\sigma^{-i}}\sigma^{\sum_{i = 0}^m i \sigma^{-i}}\le C,
\]
for some $C$ that is independent of $m$. Letting $m\to \infty$, we deduce
\[
\|v\|_{L^\infty(\Omega)} \le C,
\]
which is the desired result.
\end{proof}

\section{Interior regularity}

\label{sec:int_reg_G}

We now turn our attention to the interior regularity of solutions of 
\[
\L u = f \quad\text{in}\quad \Omega\subset \R^n. 
\]
This is a topic that has attracted a lot of attention in the last years; we refer to \cite{Bas09, BC09,Cozzi17,DSV2,DRSV,DK20,Fal,IS2,JX,KM,PK,RS-stable,Ser,Ser2,SS,Sil06,SiICM} for several interior regularity results under various assumptions on the kernels.

The strategy to establish interior regularity estimates for general integro-differential operators $\L$ is very different from the one we have used for $\fls$ (cf. Sections~\ref{sec:interior_regularity} and \ref{sec:fls}). 

In this section, we will establish the following interior regularity result, stating that solutions gain $2s$ derivatives when the kernels are regular (in the sense of Definition~\ref{defi:Gmu}).

\index{Interior regularity!General operators}

\begin{thm}\label{thm-interior-linear-2}
Let $s\in (0,1)$ and $\L\in \G_s(\lambda, \Lambda; \alpha)$ for some $\alpha>0$ such that $2s+\alpha\notin\N$. Let $f\in C^\alpha(B_1)$, and let $u\in L^\infty_{2s-\eps}(\R^n)$ for some $\eps > 0$ be a distributional solution of
\[
\L u = f \quad\text{in}\quad B_1.
\]

Then, $u\in C_{\rm loc}^{2s+\alpha}(B_{1})$ with
\[
\|u\|_{C^{2s+\alpha}(B_{1/2})} \le C\left(\|u\|_{L^\infty_{2s-\eps}(\R^n)} + \|f\|_{C^\alpha(B_1)}\right)
\]
for some $C$ depending only on $n$, $s$, $\alpha$, $\eps$,  $\lambda$, $\Lambda$, and $[\L]_\alpha$. 
\end{thm}

Notice that in Theorem~\ref{thm-interior-linear-2} we need   the kernel to have the same degree of regularity as the right-hand side, $[\L]_\alpha < \infty$, in order to gain $2s$ derivatives (recall \eqref{Calpha-assumption} and Definition~\ref{defi:Gmu}). Without this assumption, the estimate is false:

\begin{lem}
\label{lem:counter}
Let $s\in (0, 1)$. 
There exists $\L \in \GL$ with kernel $K(y)\, dy$ comparable to the fractional Laplacian, \eqref{eq:Kcompfls}, and a distributional solution of $\L u = 0$ in $B_1$ such that  $u\notin C^{2s+\eps}(B_{1/2})$ for any $\eps > 0$. 
\end{lem}
\begin{proof}
\label{proof:lem_counter}
 Let $v\in L^\infty(\R)$ and $\L\in \GL$ be the ones constructed in Lemma~\ref{lem:counterexamples_strong}, such that $v$ is compactly supported, $v\equiv 0$ in $(-2, 2)$, and $\L v\notin C^\eps((-\frac12, \frac12))$ for any $\eps >0$. 

Let $u$ be the solution to 
\[
\left\{
\begin{array}{rcll}
\L u & = & 0 & \quad\text{in}\quad (-1, 1)\\
u & = & v & \quad\text{in}\quad (-1, 1)^c,
\end{array}
\right.
\]
and let $w = u - v$, which satisfies $\L w\in L^\infty((-1, 1))$ and $w = 0$ in $(-1, 1)^c$. Hence, by Corollary~\ref{cor:globCsreg_loc_Lip_2} (which appears later in the book), $w\in C^\delta(\R )$ for some $\delta > 0$. On the other hand, $\L w = -\L v\notin C^\eps((-1, 1))$ for any $\eps > 0$, which by Lemma~\ref{lem:Lu_2}-\ref{it:lem_Lu2_ii} (using that $w\in C^\delta(\R )$) means that $w\notin C^{2s+\eps}((-1, 1))$ for any $\eps > 0$. Since $w = u$ in $(-1, 1)$, we are done. 
\end{proof}

When the right-hand side is in $L^p$, we prove that solutions are $C^{2s-n/p}$, with no regularity assumption on the kernel. 

\begin{thm} \label{thm-interior-linear-Lp}
Let $s\in (0,1)$ and $\L\in \GL$. Let $p \ge 1$ satisfying $\frac{n}{2s}< p \le \infty$ and $2s-\frac{n}{p}\notin\N$. Let $f\in L^p(B_1)$, and let $u\in L^\infty_{2s-\eps}(\R^n)$ for some $\eps > 0$ be a distributional solution of
\[
\L u = f \quad\text{in}\quad B_1.
\]

Then, $u\in C_{\rm loc}^{2s-n/p}(B_1)$ with 
\[
\|u\|_{C^{2s-n/p}(B_{1/2})} \le C\left(\|u\|_{L^\infty_{2s-\eps}(\R^n)} + \|f\|_{L^p(B_1)}\right)
\]
for some $C$ depending only on $n$, $\eps$, $p$, $\lambda$, and $\Lambda$.
\end{thm}

In particular, when $f\in L^\infty(B_1)$, we have that $u\in C^{2s}(B_1)$ except for $s = \frac12$, in which case we have $u\in C^{1-\eps}(B_1)$ for any $\eps > 0$.

In order to obtain an estimate like the one in Theorem~\ref{thm-interior-linear-2} but without regularity on the kernel, we must instead require global $C^\alpha$ regularity on $u$ (rather than  just boundedness).

\begin{prop}\label{prop-interior-linear}
Let $s\in (0,1)$ and $\L\in \GL$. Let $\alpha>0$ such that $2s+\alpha\notin \N$. Let $f\in C^\alpha(B_1)$, and let $u\in C^{\alpha}(\R^n)$ be a distributional solution of
\[
\L u = f \quad\text{in}\quad B_1.
\]

Then $u\in C_{\rm loc}^{2s+\alpha}(B_1)$ with 
\[
\|u\|_{C^{2s+\alpha}(B_{1/2})} \le C\left(\|u\|_{C^\alpha(\R^n)} + \|f\|_{C^\alpha(B_1)}\right)
\]
for some $C$ depending only on $n$, $s$, $\eps$, $\alpha$, $\lambda$, and $\Lambda$.
\end{prop}

\begin{rem}
\label{rem:int_weak}
Thanks to Lemma~\ref{lem:weak_distr}, the estimates in Theorems~\ref{thm-interior-linear-2} and \ref{thm-interior-linear-Lp} (for $p \ge\frac{2n}{n+2s}$ and $n > 2s$) and Proposition~\ref{prop-interior-linear}  hold for weak solutions, too. 
\end{rem}

 \begin{rem}
 In  Theorem~\ref{thm-interior-linear-2} we obtain $C^{2s+\alpha}$ interior regularity estimates by imposing $\alpha$-regularity on the kernel of the operator, whereas in Proposition~\ref{prop-interior-linear} we do so by imposing $C^\alpha$ global regularity of the solution. One could also obtain the interpolating statements, in which the $C^{2s+\alpha}$ interior regularity is attained by imposing $C^\gamma$ global regularity of the solution and $(\alpha-\gamma)$-regularity on the kernel, for some $\gamma\in (0, \alpha)$ (cf. Proposition~\ref{prop-interior-linear-x}). 
 \end{rem}

The proofs of Theorems~\ref{thm-interior-linear-2}, \ref{thm-interior-linear-Lp}, and Proposition~\ref{prop-interior-linear} that we present here are based on the ideas of \cite{FR,RS-stable,Sim}, and will follow by a contradiction and a blow-up argument, combined with a Liouville theorem for $\L$. 
In order to prove the Liouville theorem, moreover, we need first to introduce the heat kernel for an operator $\L\in \GL$ and describe some of its properties. 
We do so in the following.

\subsection{Heat kernel} \index{Heat kernel!General operators} 

The heat kernel for an operator $\L$ is the fundamental solution to the corresponding heat equation in the Euclidean space $\R^n$. 
That is,  a function $p(t, x): (0, \infty)\times \R^n\to \R_{\ge 0}$ such that, for any $\varphi \in C^\infty_c(\R^n)$, 
\[
\Phi(t, x) := [p(t, \cdot) * \varphi](x)
\]
satisfies (cf. Section~\ref{sec:heatkernel_sqrtl}): 
\begin{equation}
\label{eq:fract_heat_equation_L}
\left\{
\begin{array}{rcll}
\partial_t \Phi  + \L \Phi  & = & 0& \quad\text{in}\quad (0, \infty) \times \R^n\\
\Phi(0, x) & = &   \varphi(x)&\quad\text{for}\quad x\in \R^n,
\end{array}
\right.
\end{equation}
where  $\Phi(0, \cdot) = \lim_{t\downarrow 0}\Phi(t, \cdot)$ is understood as a uniform limit. Formally, $p(t, x)$ satisfies 
\[
\left\{
\begin{array}{rcll}
\partial_t p + \L p  & = & 0& \quad\text{in}\quad (0, \infty) \times \R^n\\
p(0,x)& = &  \delta_{0}&\quad\text{for}\quad x\in \R^n,
\end{array}
\right.
\]
where $\delta_{0}$ denotes the Dirac delta at 0. In subsection~\ref{ssec:inf_gen}  we referred to $p$ as the \emph{Markov transition function} associated to a L\'evy process with infinitesimal generator $\L$; see \eqref{eq:heat_eq_p}. 
\begin{lem}
\label{lem:est_heat_kernel}
Let $\L\in \GL$ with Fourier symbol $\A$, and let
\[
\hat{p}(t, \xi) := e^{-t\A(\xi)}.
\]
Then $p(t, x) := \mathcal{F}^{-1}(\hat{p}(t, \xi))(x)$ is the heat kernel for $\L$, that is, \eqref{eq:fract_heat_equation_L} holds. 
Moreover, for all $t > 0$,
\begin{equation}
\label{eq:pbounds}
p(t, x) \ge 0,\quad \int_{\R^n} p(t, x) \, dx= 1,\quad  \|\nabla_x p(1, \cdot)\|_{L^\infty(\R^n)} \le C_1,
\end{equation}
for some $C_1$ that depends only on $n$, $s$, $\lambda$, and $\Lambda$.
In addition, for any $\delta > 0$,
\begin{equation}
\label{eq:pbounds2}
\int_{\R^n} |x|^{2s-\delta} p(1, x) \, dx \le C_2
\end{equation}
for some $C_2$ that depends only on $n$, $s$, $\delta$, $\lambda$, and $\Lambda$.
\end{lem}

\begin{proof}
Notice that $\hat{p}(t, \xi) = e^{\Psi(\xi)}$ where $\Psi(\xi) = -t\A(\xi)$ is of the form \eqref{eq:Psi}. In particular, by Remark~\ref{rem:p_nonnegative},  we necessarily have that $p(t, x)$ is a probability measure, for each $t > 0$ (and $p(0, dx) = \delta_0(dx)$), so $p\ge 0$ and $\int_{\R^n} p(t, x)\, dx = 1$. (Alternatively, from Proposition~\ref{prop:equiv_fourier}, $e^{-t\tilde{\Lambda}|\xi|^{2s}}\le \hat{p}(t, \xi) \le e^{-t\tilde{\lambda}|\xi|^{2s}}$  and hence $\int_{\R^n} p(t, x)\, dx = \hat{p}(t, 0) = 1$.)
 
The expression \eqref{eq:fract_heat_equation_L} in the Fourier space is 
\[
\partial_t \hat{p}(t, \xi) \hat{\varphi}(\xi) + \A(\xi) \hat{p}(t, \xi)\hat{\varphi}(\xi) = 0,
\]
which holds for $\hat{p}(t, \xi) = e^{-t\A(\xi)}$. Since $\varphi\in C^\infty_c(\R^n)$, the first equation in \eqref{eq:fract_heat_equation_L} holds. On the other hand, since $\hat{p}(t, \xi) \to 1$ as $t\downarrow 0$, we have $p(t, x) \rightharpoonup \delta_{0}$ and the second equation also holds.

For the last expression in \eqref{eq:pbounds}, notice that $\hat{p}(t, \xi)\le e^{-t\tilde{\lambda}|\xi|^{2s}}$ for some $\tilde\lambda > 0$ depending only on $n$, $s$, $\lambda$, and $\Lambda$ (since $\L \in \GL$, see Proposition~\ref{prop:equiv_fourier}). In particular, $\hat{p}(1, \xi)$ has an exponential decay for $\xi\to \infty$, which implies that $p$ is smooth with 
\[
\|\nabla_x p (1, \cdot)\|_{L^\infty(\R^n)}\le C_1
\]
for some $C_1$ that depends only on $n$, $s$, $\lambda$, and $\Lambda$. 

Finally, let us show \eqref{eq:pbounds2}. To do that, let $\eta$ be a fixed radial cut-off function such that $\eta \in C^\infty(\R^n)$, $0\le \eta \le 1$, $\eta \equiv 1 $ in $B_{1/2}$ and $\eta \equiv 0$ in $\R^n\setminus B_1$. We define, for any $R\ge 1$ fixed, 
\[
\varphi(x) := (1+|x|^2)^{s-\delta/2} \eta(x/R),
\]
and for any $\rho \ge 1$, 
\[
\varphi_\rho(x) := \rho^{-2s+\delta}\varphi(\rho x) = (\rho^{-2}+|x|^2)^{s-\delta/2} \eta(\rho x/R),
\]
Then, we have
\[
\begin{split}
\|D^2\varphi_\rho \|_{L^\infty(B_2\setminus B_1)} & \le C\bigg(1+\frac{\rho}{R}\|\nabla \eta\|_{L^\infty(B_{2\rho/R}\setminus B_{\rho/R})} \\
&\qquad  \qquad +\frac{\rho^2}{R^2} \|D^2\eta\|_{L^\infty(B_{2\rho/R}\setminus B_{\rho/R})} \bigg),
\end{split}
\]
for some $C$ depending only on $n$, $s$, and $\delta$, and where we are using that $(\rho^{-2}+|x|^2)^{s-\delta/2}$ is smooth in $B_2\setminus B_1$. Now notice that 
\[
\|\nabla \eta\|_{L^\infty(B_{2\rho/R}\setminus B_{\rho/R})} =  \|D^2\eta\|_{L^\infty(B_{2\rho/R}\setminus B_{\rho/R})} = 0\quad\text{if}\quad \frac{\rho}{R}\ge 1,
\] (since $\eta$ is supported in $B_1$), and hence we actually have 
\[
\|\varphi_\rho \|_{C^2(B_2\setminus B_1)}  \le C,
\]
for some $C$ depending only on $n$, $s$, and $\delta$ (and $\eta$, but it is fixed universally). On the other hand, by definition we also have that 
\[
\|\varphi_\rho\|_{L^\infty_{2s-\delta/2}(\R^n)}\le C,
\]
for some $C$ depending only on $n$, $s$, and $\delta$. Hence, we can apply Lemma~\ref{lem:Lu} to deduce that 
\[
\|\tilde \L \varphi_\rho\|_{L^\infty(B_{7/4}\setminus B_{5/4})}\le C\quad\text{for any}\quad \tilde \L \in \GL,
\]
for some $C$ depending only on $n$, $s$, $\delta$, $\lambda$, and $\Lambda$. In particular, by the scale invariance of the class $\GL$, Remark~\ref{rem:Scale_invariance}, we have that for any $\rho \ge 1$,
\[
|\L \varphi| \le C \rho^{-\delta}\quad\text{in}\quad   B_{7\rho/4}\setminus B_{5\rho/4}. 
\]
Together with the fact that, again by Lemma~\ref{lem:Lu}, $|\L \varphi| \le C$ in $B_2$, we deduce that 
\[
\|\L \varphi\|_{L^\infty(\R^n)}\le C,
\]
for some $C$ depending only on $n$, $s$, $\delta$, $\lambda$, and $\Lambda$ (in particular, it is independent of $R$). We can now compute, by \eqref{eq:fract_heat_equation_L} together with  Remark~\ref{rem:conv_dist_sol}:
\[
\begin{split}
\big| p(1, \cdot) * \varphi  -\varphi\big| & = \left|\int_0^1 \partial_t \left[ p(t, \cdot) * \varphi \right] \, dt\right| \\
& \le \int_0^1 \left|\L \left[ p(t, \cdot) * \varphi\right] \right|\, dt\\
& = \int_0^1 \left|p(t, \cdot) * (\L \varphi) \right|\, dt\\
&  \le C \int_0^1 \int_{\R^n} p(t, y) \, dy \, dt = C, 
\end{split}
\]
where we have also used that $p(t, x)$ is a probability density for each $t \ge 0$. Hence, since $\varphi$ is radially symmetric, in particular we obtain
\[
\big(p(1, \cdot)*\varphi\big)(0) = \int_{\R^n} p(1, x)\varphi(x)\, dx \le |\varphi(0)|+ C = 1+C. 
\]
Since the constants are independent of $R$, we can now let $R\to \infty$ in the definition of $\varphi$ to obtain the desired result, with $C_2 = 1+C$. 
\end{proof}

 \begin{rem}
 \label{rem:fund_sol}
 The previous lemma, and in particular \eqref{eq:pbounds2}, gives an integral decay for the heat kernel $p(1, x)$. 
Under the extra assumption that the kernel $K$ of the operator $\L$ is absolutely continuous and comparable to the one of fractional Laplacian,
 \begin{equation}\label{rem-comp}
 K(y)\asymp |y|^{-n-2s}\quad \textrm{in}\quad \R^n,
\end{equation}  
then in fact we have a \textit{pointwise} decay for $p$, namely
 \begin{equation}\label{rem-heat-pointwise}
 p(t,x) \asymp \frac{t}{t^{\frac{n+2s}{2s}}+|x|^{n+2s}}\quad \textrm{in}\quad \R^n\times(0,\infty).
\end{equation} 
Indeed, assume for simplicity that $K$ is in addition homogeneous, \eqref{eq:Khom}-\eqref{eq:Khom2}.
Then, we have
\begin{equation}
\label{eq:homog_p}
 p(t, x) =   t^{-\frac{n}{2s}}p(1, t^{-\frac{1}{2s}}x) \quad\text{for all}\quad t>0,~~x\in \R^n.
 \end{equation}
On the other hand, for any $f\in C^\infty_c(\R^n\setminus\{0\})$ ,
\[
 \begin{split}
\int_{\R^n} f(x)K(x)\, dx & = -\L f(0) = \partial_t\big|_{t = 0} (p(t,\cdot)*f)(0) \\
& = \lim_{t\downarrow 0}\frac{1}{t} \int_{\R^n} p(t, x) f(x)\, dx,
 \end{split}
 \]
and from \eqref{eq:homog_p}, for any set $A\subset \R^n$ with $0\notin A$, 
 \[
 \int_{A} K(x)\, dx = \lim_{R \to \infty}  R^{n+2s} \int_{A} p(1, R x)\, dx.
 \]
 Now, heuristically, since \eqref{rem-comp} holds, the previous equality should imply a pointwise asymptotic bound for $p(1, x)$ of the form $p(1, x)\sim |x|^{-n-2s}$, when $|x|$ is large. 
 That is, we should have a bound of the form
\[
 p(1, x)\asymp \frac{1}{1+|x|^{n+2s}}, 
\]
which implies \eqref{rem-heat-pointwise}.
We refer to \cite{CKS} for an actual proof of \eqref{rem-heat-pointwise} for all operators $\L$ satisfying \eqref{rem-comp}. 
 
Finally, notice that, when $n > 2s$, this bound for the heat kernel also yields sharp bounds for the fundamental solution $\Gamma$ of the operator $\L$ (i.e., $\L \Gamma = \delta_0$), which can be obtained by integrating in time the heat kernel (cf. second proof of Lemma~\ref{lem:fundamental_solution}):
 \[
 \Gamma (x) = \int_0^\infty p(t, x) \asymp |x|^{2s-n},
 \]
 where we have used \eqref{rem-heat-pointwise}.
 \end{rem}

\subsection{Liouville's theorem} \index{Liouville's theorem!General operators} Let us start by stating and proving Liouville's theorem for globally bounded solutions to integro-differential equations $\L u = 0$.

\begin{thm}[Liouville's Theorem]
\label{thm:Liouville}
Let $s\in (0,1)$ and $\L\in \GL$. Let $u\in L^\infty(\R^n)$ be a distributional solution to 
\[
\L u = 0 \quad \text{in}\quad \R^n.
\]
Then, $u$ is constant. 
\end{thm}

\begin{proof} We will show that $u$ is continuous, with a  (H\"older) modulus of continuity that depends only on the ellipticity constants of $\L$. Since the class $\GL$ is invariant under scaling, this will imply that $u$ is in fact constant. Up to dividing by a constant, let us assume $|u|\le 1 $ in $\R^n$.

Let $R \ge 1$ and let us define $u_R(x) := u(Rx)$. Observe that   $\|u\|_{L^\infty(\R^n)} = \|u_R\|_{L^\infty(\R^n)}\le 1$, and that 
\[
\L_R u_R = 0\quad\text{in}\quad \R^n
\]
in the distributional sense, where $\L_R$ is an operator of the form \eqref{eq:Lu_nu} with kernel $K_R$   given by $K_R(dy) = R^{2s} K(R\, dy)$ (where $K$ is the kernel of~$\L$; see Remark~\ref{rem:Scale_invariance}), and hence $\L_R\in \GL$. In particular, thanks to Lemma~\ref{lem:est_heat_kernel}, the corresponding heat kernel, $p_R$, satisfies 
\[
\|\nabla_x p_R(1, x)\|_{L^\infty(\R^n)} \le C
\]
for some $C$ independent of $R$. 

Observe now that $u_R = p_R(1, \cdot)* u_R$. Indeed, by \eqref{eq:fract_heat_equation_L} (which holds here by density) together with Lemma~\ref{lem:conv_dist_sol} (see Remark~\ref{rem:conv_dist_sol}) we get:
\[
u_R - p_R(1, \cdot) * u_R  = -\int_0^1 \partial_t \left[ p_R(t, \cdot) * u_R \right] \, dt =  \int_0^1 \L \left[ p_R(t, \cdot) * u_R \right] \, dt = 0.
\]
Now, given $x, x'\in \R^n$  we have, for any $M \ge 0$ (denoting $p_R(z) := p_R(1, z)$),
\[
\begin{split}
\left| u_R(x) - u_R(x')\right| & = \left|\int_{\R^n}\left(p_R(x-y) - p_R(x'-y)\right) u_R(y) \, dy\right|\\
& \le C M^n |x-x'|  +  C\int_{|y|\ge M} \left(p_R(y)+p_R(x'-x+y)\right)\, dy, 
\end{split}
\]
where we have used $p_R \ge 0$, $|u_R|\le 1$ in $\R^n$, and $\|\nabla_x p_R\|_{L^\infty(\R^n)}\le C$ independently of $R$. 
Now, by \eqref{eq:pbounds2}, $\int_{B_M^c} p_R(x)\, dx \le  M^{\delta-2s} \int_{B_M^c} |x|^{2s-\delta} p_R(x)\, dx \le C M^{\delta-2s}$ for some $C$ depending only on $n$, $s$, $\delta$, $\lambda$, and $\Lambda$, but independent of the operator $\L_R$.

Fix now $M = |x-x'|^{-\frac{1}{2n}}$ so that
\[
\left| u_R(x) - u_R(x')\right| \le C \left(|x-x'|^{\frac12}  +  \|p_R\|_{L^1(B^c_M\cup B_M^c(x'-x))}\right)
\le C|x-x'|^{\gamma},
\]
with $\gamma>0$. 
Hence, rescaling back to $u$, we get
\[
\left| u(x) - u(x')\right| \le  C\left|\frac{x-x'}{R}\right|^\gamma = C \frac{|x-x'|^\gamma}{R^\gamma},
\]
for any $x, x'\in \R^n$.
Letting $R \to \infty$ we deduce $u$ is constant. 
\end{proof}

\subsection{A compactness argument} 
\label{ssec:compactness}
Let us state some useful preliminary results on the compactness argument to prove the interior regularity estimate. 

In the following statement we prove a quantitative estimate for solutions in very large balls. It can be a seen as a quantitative version of Liouville's theorem. 
\begin{prop}
\label{prop:compactness}
Let $s\in (0,1)$, $\delta > 0$, and $\L\in \GL$. Let $\alpha > 0$ with $2s+\alpha \notin\N$, and let $u\in C_{\rm loc}^{2s+\alpha}(\R^n)\cap L^\infty_{2s-\eps}(\R^n)$ for some $\eps > 0$, with $[u]_{C^{2s+\alpha}(\R^n)}\le 1$ and 
\[
\L u = f\quad\text{in}\quad B_{{1}/{\delta}}
\]
for some $f\in C^\alpha(B_{{1}/{\delta}})$ such that $[f]_{C^\alpha(B_{{1}/{\delta}})} \le \delta$. 

Then, for every $\eps_\circ > 0$ there exists $\delta_\circ >0 $ depending only on $\eps_\circ$, $n$, $s$, $\alpha$, $\eps$, $\lambda$, and $\Lambda$, such that if $\delta < \delta_\circ$,
\[
\|u-q\|_{C^\nu(B_1)}\le \eps_\circ,
\]
where $q$ is the Taylor polynomial of $u$ at 0 of degree $\nu:= \lfloor 2s+\alpha\rfloor$. 
\end{prop}

In the proof, we will use incremental quotients on functions $f:\R^n\to \R$ (see Section~\ref{sec:incremental_quotients} in Appendix~\ref{app.A}). That is, for $h\in \R^n$ we define the first order incremental quotient $D_h f$ as 
\[
D_h f (x) := \frac{f(x+h)- f(x) }{|h|}.
\]
More generally, we define the $m$-th order incremental quotient $D^m_h f$ recurrently as 
\[
D_h^m f(x) = D_h (D_h^{m-1} f(x)) = \frac{1}{|h|^m}\sum_{k = 0}^m {m \choose k} (-1)^{m-k} f(x+k h).
\]
Among other properties, we will use that if $f\in C^{m-1}(\R^n)$ and $D_h^m f$ is constant for any $h\in \R^n$, then $f$ is a polynomial of degree $m$, see Lemma~\ref{it:H10}.

\begin{proof}[Proof of Proposition~\ref{prop:compactness}]
Let us argue by contradiction,  i.e., let us assume that the statement does not hold. That is, there exists some $\eps_\circ > 0$ such that for any $k\in \N$, there are $u_k\in C_{\rm loc}^{2s+\alpha}(\R^n)\cap L^\infty_{2s-\eps}(\R^n)$ with $[u_k]_{C^{2s+\alpha}(\R^n)}\le 1$, $f_k\in C^\alpha(B_k)$ with $[f_k]_{C^\alpha(B_k)}\le \frac{1}{k}$, and $\L_k\in \GL$ such that
\[
\L_k u_k = f_k \quad\text{in}\quad B_k
\]
but 
\[
\|u_k - q_k\|_{C^\nu(B_1)}\ge \eps_\circ,
\]
where $q_k$ is the Taylor polynomial of $u_k$ at 0 of order $\nu$. If we denote $v_k := u_k - q_k$, we have 
\begin{equation}
\label{eq:vkzero}
v_k(0) = |\nabla v_k(0)| =   \dots  = |D^\nu v_k(0)| = 0, 
\end{equation}
and by assumption, up to a subsequence, we know that $v_k \to v$ in $C^\nu_{\rm loc}(\R^n)$ for some $v$ with $[v]_{C^{2s+\alpha}(\R^n)}\le 1$ (by Arzel\`a-Ascoli), and satisfying \eqref{eq:vkzero} as well as $\|v\|_{C^\nu(B_1)}\ge \eps_\circ$.

We now look for an equation satisfied by  the limiting function $v$. Let us define, for a fixed $h\in \R^n$ and $m := \lceil \alpha \rceil$,
\[
 V_k := |h|^{\nu+1}D_h^{\nu+1} u_k = |h|^{\nu+1}D_h^{\nu+1} v_k, \qquad F_k := |h|^m D_h^{m} f_k,
\]
with $\|F_k\|_{L^\infty(B_k)} \le \frac{C}{k}|h|^{\alpha}$ (see, for example, Lemma~\ref{it:H7_gen}). Observe that, by linearity and translation invariance of the operators $\L_k$, since $\nu \ge m-1$, 
\[
\L_k  V_k = |h|^{\nu+1-m}D_h^{\nu+1-m} F_k  \to 0\quad\text{locally uniformly in }\R^n. 
\]
On the other hand, and by assumption, we have  $\|V_k\|_{L^\infty(\R^n)}\le C |h|^{2s+\alpha}$ and $[ V_k]_{C^{2s+\alpha}(\R^n)}\le  C(\nu)$ (again by Lemma~\ref{it:H7_gen}). That is, up to a subsequence and thanks to Proposition~\ref{prop:stab_distr}, $V_k$ is converging as $k \to \infty$ to some bounded $V\in C_{\rm loc}^{2s+\alpha}(\R^n)$ such that
\[
\L V=0 \quad \text{in}\quad \R^n,
\]
for some $\L \in \GL$. By Liouville's theorem, Theorem~\ref{thm:Liouville}, $V$ is constant. Observe that from the $C^\nu_{\rm loc}$ convergence $v_k \to v$, we have that $V = |h|^{\nu+1}D_h^{\nu+1} v$ is constant for every $h$, and by Lemma~\ref{it:H10} we have that $v$ is a polynomial of degree $\nu+1$. But, since $[v]_{C^{2s+\alpha}(\R^n)}\le 1$ and $2s+\alpha < \nu+1$, we must have that $v$ is a polynomial of degree $\nu$. Because it satisfies \eqref{eq:vkzero}, it must be $v \equiv 0$, which is a contradiction with $\|v\|_{C^\nu(B_1)}\ge \eps_\circ>0$.
\end{proof}

We also need an analogous result for distributional solutions when the right hand side is in $L^p$:
\begin{prop}
\label{prop:compactness_Lp}
Let $s\in (0,1)$, $\delta > 0$, and $\L\in \GL$. Let $p \ge 1$ and $\frac{n}{2s}< p \le \infty$ with $2s-\frac{n}{p}\notin\N$, and let $u\in C_{\rm loc}^{2s-n/p}(\R^n)\cap L^\infty_{2s-\eps}(\R^n)$ for some $\eps > 0$, with $[u]_{C^{2s-n/p}(\R^n)}\le 1$ and 
\[
\L u = f\quad\text{in}\quad B_{ {1}/{\delta}}
\]
in the distributional sense, for some $f\in L^p(B_{ {1}/{\delta}})$ such that $\|f\|_{L^p(B_{ {1}/{\delta}})} \le \delta$. 

Then, for every $\eps_\circ > 0$ there exist $\delta_\circ >0 $ depending only on $\eps_\circ$, $n$, $s$, $p$, $\eps$, $\lambda$, and $\Lambda$, such that if $\delta < \delta_\circ$,
\[
\|u-q\|_{C^\nu(B_1)}\le \eps_\circ,
\]
where $q$ is the Taylor polynomial of $u$ at 0 of degree $\nu:= \lfloor 2s-n/p\rfloor$. 
\end{prop}
\begin{proof}
The proof is analogous to the proof of Proposition~\ref{prop:compactness}. In this case we take $m = 0$ (so $D_h^0 f = f$) and we observe that by the properties of distributional solutions, if $\L_k u_k = f_k$ in the distributional sense, then $\L_k V_k = |h|^{\nu+1}D^{\nu+1}_h f_k$ in the distributional sense as well, where $V_k = |h|^{\nu+1} D^{\nu+1}_h v_k$ is converging locally uniformly to some $V\in C^{2s-\frac{n}{p}}(\R^n)$. Now we have $D^{\nu+1}_h f_k\to 0$ locally in $L^p$ for $p \ge 1$, and this is enough to apply the stability result for distributional solutions, Proposition~\ref{prop:stab_distr}, and conclude the proof as in Proposition~\ref{prop:compactness}. 
\end{proof}

We will also need the following general lemma that provides a blow-up sequence:

\begin{lem}\label{lem-interior-blowup}
Let $\mu>0$ with $\mu\notin \mathbb Z$, let $\mathcal S: C^\mu(\R^n)\to\R_{\ge 0}$, and let $\delta>0$.
Then,
\begin{enumerate}[leftmargin=*, label=(\roman*)]
\item \label{it:lem_int_blowup_i} either we have
\[[u]_{C^\mu(B_{1/2})} \leq \delta [u]_{C^\mu(\R^n)} + C_\delta\big(\|u\|_{L^\infty(B_1)} + \mathcal S(u)\big)\]
for all $u\in C^\mu(\R^n)$, for some $C_\delta$ depending only on $\mu$, $\mathcal S$, and $\delta$,

\item \label{it:lem_int_blowup_ii} or there is a sequence $u_k\in C^\mu(\R^n)$, with 
\begin{equation}\label{lem-interior-blowup2}
\frac{\mathcal S(u_k)}{[u_k]_{C^\mu(B_{1/2})}} \longrightarrow 0,
\end{equation}
and there are $r_k\to0$, $x_k\in B_{1/2}$, such that if we define
\begin{equation}\label{lem-interior-blowup3}
v_k(x) := \frac{u_k(x_k+r_k x)}{r_k^\mu[u_k]_{C^\mu(\R^n)}},
\end{equation}
and we denote by $p_k$ the $\nu$-th order Taylor polynomial of $v_k$ at 0, then 
\[
\|v_k - p_k\|_{C^\nu(B_1)} > \frac{\delta}{2}
\]
for $k$ large enough, where $\nu=\lfloor \mu\rfloor$.
\end{enumerate}
\end{lem}

\begin{proof}
Assume \ref{it:lem_int_blowup_i} does not hold, and let $\mu=\nu+\beta$, where $\nu\in\mathbb Z$ and $\beta\in(0,1)$.
Then, there is a sequence $u_k\in C^\mu(\R^n)$ such that
\[ [u_k]_{C^\mu(B_{1/2})} \geq \delta [u_k]_{C^\mu(\R^n)} + k\big(\|u_k\|_{L^\infty(B_{1})} + \mathcal S(u_k)\big).\]
Such sequence clearly satisfies \eqref{lem-interior-blowup2} (with rate $\frac{1}{k}$).
Let $x_k,y_k\in B_{1/2}$ be such that 
\[\frac{\big|D^\nu u_k(x_k)-D^\nu u_k(y_k)\big|}{|x_k-y_k|^\beta} \geq \frac{1}{2}[u_k]_{C^\mu(B_{1/2})}.\]
Then, we claim that 
\[r_k:=|x_k-y_k|\longrightarrow 0.\]
Indeed, if $r_k^\beta\geq 4\varepsilon>0$ for all $k\in \N$, with $\varepsilon\ll 1$, then
\[\frac{1}{2}[u_k]_{C^\mu(B_{1/2})} \leq \frac{2\|u_k\|_{C^\nu(B_{1/2})}}{|x_k-y_k|^\beta} \leq \frac{\varepsilon[u_k]_{C^\mu(B_{1/2})} + C_\varepsilon\|u_k\|_{L^\infty(B_{1/2})}}{r_k^\beta},\]
where we used the interpolation inequality in Proposition~\ref{it:H9}.
Now, since $\|u_k\|_{L^\infty(B_{1/2})} \leq \frac1k [u_k]_{C^\mu(B_{1/2})}$, the previous inequality yields
\[\frac{1}{2} \leq \frac{1}{4} + \frac{C_\varepsilon}{4\eps k},\]
a contradiction if $k$ is large enough.
Thus, $r_k\to0$ as $k\to \infty$, as wanted.

Define now $v_k$ as in \eqref{lem-interior-blowup3} (which clearly satisfies $[v_k]_{C^\mu(\R^n)}=1$), and let  $p_k$ be the Taylor polynomials of order $\nu$, so that $w_k := v_k - p_k$ satisfies $w_k(0)=...=|D^\nu w_k(0)|=0$.
If we let 
\[z_k := \frac{x_k-y_k}{r_k} \in \partial B_1,\]
we have 
\[|D^\nu w_k(z_k)| = \frac{\big|D^\nu u_k(x_k)-D^\nu u_k(y_k)\big|}{r_k^\beta [u_k]_{C^\mu(\R^n)}} \geq  \frac{\frac{1}{2}[u_k]_{C^\mu(B_{1/2})}}{[u_k]_{C^\mu(\R^n)}} > \frac{\delta}{2},\]
and the lemma follows.
\end{proof}

\subsection{Proof of the interior estimates} \label{ssec:interiorproofs}Before proving Theorem~\ref{thm-interior-linear-2}, let us first show the following intermediate result. Notice that the following statement is close to our desired estimate: if we could let $\delta \downarrow 0$ and $C_\delta$ remained constant, we would prove the regularity for $u$. Alternatively, if we could take the H\"older norm on the right-hand side in $B_{1/2}$ instead of $\R^n$, we would also be done by taking $\delta$ small enough. As we will see, this statement is easier to prove and will still yield Theorem~\ref{thm-interior-linear-2} as a consequence.

\begin{prop}
\label{prop:interior-linear}
Let $s\in (0,1)$, $\L\in \GL$, and $u \in C^\infty_c(\R^n)$. Then, the following holds. 
\begin{enumerate}[leftmargin=*, label=(\roman*)]
\item  \label{it:int_lin_i} Let $\alpha>0$ with $2s+\alpha\notin \N$. Then, for any $\delta > 0$ there exists $C_\delta$ such that 
\[
[u]_{C^{2s+\alpha}(B_{1/2})} \le \delta [u]_{C^{2s+\alpha}(\R^n)} + C_\delta \left(\|u\|_{L^\infty(B_1)} + [\L u]_{C^\alpha(B_1)}\right),
\]
where $C_\delta$ depends only on $\delta$, $n$, $s$, $\alpha$, $\lambda$, and $\Lambda$. 
\item \label{it:int_lin_ii}  Let $p \ge 1$ with $\frac{n}{2s} < p \le \infty$. Then, for any $\delta > 0$ there exists $C_\delta$ such that 
\[
[u]_{C^{2s-n/p}(B_{1/2})} \le \delta [u]_{C^{2s-n/p}(\R^n)} + C_\delta \left(\|u\|_{L^\infty(B_1)} + \|\L u\|_{L^p(B_1)}\right),
\]
where $C_\delta$ depends only on $\delta$, $n$, $s$, $p$, $\lambda$, and $\Lambda$. 
\end{enumerate}
\end{prop}

\begin{proof}
Let us start with part \ref{it:int_lin_i}. 
We use Lemma~\ref{lem-interior-blowup} with $\mu = 2s+\alpha$ and
\[
\mathcal{S}(w) =  \left\{
\begin{array}{ll}
\inf_{\tilde \L \in \GL} [\tilde \L w]_{C^\alpha(B_1)}  & \text{if}\quad w\in  C_c^\infty(\R^n),
\\
\infty & \text{otherwise,}
\end{array}
\right.
\]
so that the mapping $\mathcal{S}:C^\mu(\R^n)\to \R_{\ge 0}$ depends only on $n$, $s$,  $\alpha$, $\lambda$, and~$\Lambda$. 
Thus, either Lemma~\ref{lem-interior-blowup}~\ref{it:lem_int_blowup_i} holds, in which case we are done, or there exists a sequence $u_k\in C^{\infty}_c(\R^n)$ and $\L_k\in \GL$ such that if we denote $f_k = \L_k u_k$ then
\begin{equation}
\label{eq:fktozero}
\frac{[f_k]_{C^\alpha(B_1)}}{[u_k]_{C^{2s+\alpha}(B_{1/2})}} \le \frac{2\mathcal{S}(u_k)}{[u_k]_{C^{2s+\alpha}(B_{1/2})}}\to 0,
\end{equation}
and for some $x_k \in B_{1/2}$ and $r_k \downarrow 0$,
\[
v_k(x) := \frac{u_k(x_k+r_k x)}{r_k^{2s+\alpha}[u_k]_{C^{2s+\alpha}(\R^n)}}
\]
satisfies 
\begin{equation}
\label{eq:vkcontradiction}
\|v_k - p_k\|_{C^\nu(B_1)} > \frac{\delta}{2},
\end{equation}
where $p_k$ is the $\nu$-th order Taylor polynomial of $v_k$ at 0.  Then, there exists $\tilde \L_k\in \GL$ such that 
\[
\tilde \L_k v_k (x) = \frac{f_k(x_k+r_k x)}{r_k^\alpha [u_k]_{C^{2s+\alpha}(\R^n)}}
\] (see Remark~\ref{rem:Scale_invariance}) so that, from \eqref{eq:fktozero}, 
\[
[\tilde \L_k v_k]_{C^\alpha (B_{1 /(2r_k)})} 
= \frac{[f_k(x_k+r_k\cdot)]_{C^\alpha (B_{1 /(2r_k)})}}{r_k^\alpha [u_k]_{C^{2s+\alpha}(\R^n)}}
\leq  \frac{[f_k]_{C^\alpha(B_1)}}{[u_k]_{C^{2s+\alpha}(\R^n)}}
\to 0, 
\] 
as $k\to \infty$. Since by definition $[v_k]_{C^\mu(\R^n)} = 1$,  we have that $v_k$ satisfies the hypotheses of Proposition~\ref{prop:compactness} for any fixed $\delta_\circ>0$, if $k$ is large enough. 
In particular, taking $\eps_\circ$ sufficiently small in Proposition~\ref{prop:compactness} we get a contradiction with \eqref{eq:vkcontradiction}.

Part \ref{it:int_lin_ii} follows in the same way, defining 
\[
\mathcal{S}(w) =  \left\{
\begin{array}{ll}
\inf_{\tilde \L \in \GL} \|\tilde \L w\|_{L^p(B_1)} & \text{if}\quad w\in   C_c^\infty(\R^n),
\\
\infty & \text{otherwise,}
\end{array}
\right.
\]
using Proposition~\ref{prop:compactness_Lp} instead of Proposition~\ref{prop:compactness}, and recalling that the translation and scale invariance are also true for distributional solutions. The exponent $2s-\frac{n}{p}$ comes from the scaling,  $\|f_k(x_k+r_k\cdot)\|_{L^p (B_{1 /(2r_k)})}\le r_k^{-\frac{n}{p}}\|f_k\|_{L^p(B_1)}$. 
\end{proof}

We now have all the tools to prove the interior regularity:

\begin{proof}[Proof of Theorem \ref{thm-interior-linear-2}]
We   divide the proof into three steps. 

\begin{steps}
\item \label{step:main}
 By Proposition~\ref{prop:interior-linear}-\ref{it:int_lin_i}, for any $\delta> 0$ there exists $C_\delta$ depending only on $\delta$, $n$, $s$, $\alpha$, $\lambda$, and $\Lambda$, such that
\begin{equation}
\label{eq:touse}
[u]_{C^{2s+\alpha}(B_{1/2})}\le \delta [u]_{C^{2s+\alpha}(\R^n)} + C_\delta \left(\|u\|_{L^\infty(B_1)}+[\L u ]_{C^\alpha(B_1)}\right)
\end{equation}
for any $u\in C^\infty_c(\R^n)$.

 Let $\eta\in C^\infty_c(B_3)$ such that $\eta \equiv 1 $ in $B_{2}$, and apply \eqref{eq:touse} to $u\eta$ for any $u\in C^\infty(\R^n)\cap L^\infty_{2s-\eps}(\R^n)$. That is, for any $\delta > 0$ there exists some $C_\delta$ such that 
\[
[\eta u]_{C^{2s+\alpha}(B_{1/2})}\le \delta [\eta u]_{C^{2s+\alpha}(B_3)} + C_\delta \left(\|\eta u\|_{L^\infty(B_1)}+[\L (\eta u)]_{C^\alpha(B_1)}\right),
\]
for any $u\in C^\infty(\R^n)\cap L^\infty_{2s-\eps}(\R^n)$.

Since $\L (\eta u) = \L u +  \L (u-\eta u)$ in $B_1$, and $u - \eta u \equiv 0$ in $B_{2}$, we have from Lemma~\ref{lem:Lu_2}-\ref{it:lem_Lu2_i} (since $\L\in \G_s(\lambda, \Lambda; \alpha)$)
\[
[\L (u-\eta u)]_{C^\alpha(B_1)}\le C \|u(1-\eta)\|_{L^\infty_{2s-\eps}(\R^n)}\le C \|u\|_{L^\infty_{2s-\eps}(\R^n)}. 
\]
Hence, 
\[
[\L (\eta u)]_{C^\alpha(B_1)} \le [\L u]_{C^\alpha(B_1)} + C \|u\|_{L^\infty_{2s-\eps}(\R^n)}, 
\]
and we get (taking $\delta$ smaller if necessary, thanks to Proposition~\ref{prop:A_imp}))
\begin{equation}
\label{eq:touse2}
[u]_{C^{2s+\alpha}(B_{1/2})}\le \delta [u]_{C^{2s+\alpha}(B_4)} + C_\delta \left(\|u\|_{L^\infty_{2s-\eps}(\R^n)}+[\L u]_{C^\alpha(B_1)}\right)
\end{equation}
for any $u\in C^\infty(\R^n)\cap L^\infty_{2s-\eps}(\R^n)$.

\item  \label{step:SAL} In order to deduce the desired estimate from \eqref{eq:touse2}, we proceed by a standard  interpolation and covering argument (see, for example, \cite[Lemma 2.27 or Theorem 2.20]{FR4}). Let us consider, for $\mu> 0$, $\mu\notin \N$, the following weighted norm in a domain $\Omega\subset \R^n$,
\[
[w]^*_{\mu; \Omega} := \sup_{B_{2r}(z) \subset \Omega}\left(r^\mu [w]_{C^\mu(B_r(z))}\right),
\]
 where the supremum is taken among all balls such that $B_{2r}(z)\subset \Omega$ with $z\in \Omega$ and $r > 0$, i.e., all $z\in \Omega$ and $r < \frac12 \dist(z, \R^n\setminus\Omega)$.  Observe first that, for some constant $C$ depending only on $n$ and $\mu$, 
\begin{equation}
\label{eq:willuse}
[w]_{\mu, \Omega}^* \le C \sup_{B_{2r}(z) \subset \Omega}\left(r^\mu [w]_{C^\mu(B_{r/8}(z))}\right).
\end{equation}
\begin{figure}
\centering
\makebox[\textwidth][c]{\includegraphics[scale = 1]{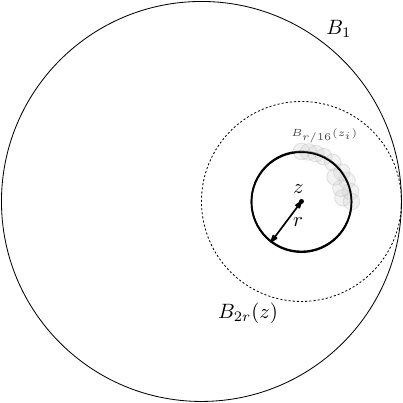}}
\caption{\label{fig:08} Covering of $B_1$ and $B_r(z)$ from \ref{step:SAL}.}
\end{figure}
Indeed,   each $B_r(z)$ with $B_{2r}(z)\subset B_1$  can be covered with $N$ smaller balls $(B_{r/16}(z_i))_{i = 1}^N$ with $z_i\in B_r(z)$ and, since $B_{r }(z_i) \subset B_1$, we have 
\[
2^{-\mu} r^\mu [w]_{C^\mu(B_{r/16}(z_i))} \le \sup_{B_{2\rho}(\bar z)\subset B_1} \rho^\mu[w]_{C^\mu(B_{\rho/8}(\bar z))}.
\]
(See Figure~\ref{fig:08}.) Combining the previous inequality  with the fact that $[w]_{C^\mu(B_{r}(z))} \le \sum_{i = 1}^N [w]_{C^\mu(B_{r/16}(z_i))} $ we obtain \eqref{eq:willuse} after taking a supremum. 

Rescaling the estimate \eqref{eq:touse2} to any ball\footnote{We are using here that if $u_r(x) := u(rx)$ then $\|u_r\|_{L^\infty_{2s-\eps}(\R^n)}\le \|u\|_{L^\infty_{2s-\eps}(\R^n)}$ for $r < 1$, and that if $\L u = f$, then there exists some $\L_r\in \GL$ such that $\L_r u_r (x) = r^{2s} f(rx)$ (see Remark~\ref{rem:Scale_invariance}), so that the same estimate applies.
} $B_{r/4}(z) \subset B_{2r}(z) \subset B_1$ we have
\[
\begin{split}
\left(\frac{r}{4}\right)^{2s+\alpha}[u]_{C^{2s+\alpha}(B_{r/8}(z))}& \le \delta r^{2s+\alpha}[u]_{C^{2s+\alpha}(B_{r}(z))} \\
& \quad + C_\delta \left(\|u\|_{L^\infty_{2s-\eps}(\R^n)}+r^{2s+\alpha} [\L u]_{C^\alpha(B_{r/2})}\right)\\
& \le \delta [u]^*_{2s+\alpha; B_1} + C_\delta \left(\|u\|_{L^\infty_{2s-\eps}(\R^n)}+ [\L u]_{C^\alpha(B_1)}\right).
\end{split}
\]
Taking the supremum among all balls $B_{2r}(z)\subset B_1$  and using \eqref{eq:willuse}, we get 
\[
\frac{1}{C} [u]^*_{2s+\alpha; B_1} \le \delta [u]^*_{2s+\alpha; B_1} + C_\delta \left(\|u\|_{L^\infty_{2s-\eps}(\R^n)}+ [\L u]_{C^\alpha(B_1)}\right).
\]
Let now $\delta> 0$ fixed such that $\delta \le \frac{1}{2C}$ to obtain 
\[
[u]^*_{2s+\alpha; B_1} \le C \left(\|u\|_{L^\infty_{2s-\eps}(\R^n)}+ [\L u]_{C^\alpha(B_1)}\right),
\]
for any $u\in C^\infty(\R^n)\cap L^\infty_{2s-\eps}(\R^n)$. 
Since $[u]_{C^{2s+\alpha}(B_{1/2})} \le 
[u]^*_{2s+\alpha; B_1}$, and $\|u\|_{C^{2s+\alpha}(B_{1/2})} \le C [u]_{C^{2s+\alpha}(B_{1/2})} +C\|u\|_{L^\infty(B_{1/2})}$ (see Proposition \ref{it:H9}), we deduce that
\[
\|u\|_{C^{2s+\alpha}(B_{1/2})} \le C \left(\|u\|_{L^\infty_{2s-\eps}(\R^n)}+ [\L u]_{C^\alpha(B_1)}\right)
\]
for all $u\in C^\infty(\R^n)\cap L^\infty_{2s-\eps}(\R^n)$.

\item   Let us now drop the smoothness assumption on $u$, and let us suppose that $u\in L^\infty_{2s-\eps}(\R^n)$ solves 
\[
\L u = f\quad \text{in}\quad B_1
\]
in the distributional sense. 
We regularize $u$ by convoluting it with a smooth mollifier. That is, let $\psi\in C^\infty_c(B_1)$, and let $\psi_\delta (x) := \delta^{-n}\psi(x/\delta)$, so that $\int_{\R^n}\psi_\delta = 1$ and $\psi_\delta\in C^\infty_c(B_\delta)$. We define 
\[
u_\delta = u*\psi_\delta.
\]
By Lemma~\ref{lem:conv_dist_sol},
\[
\L u_\delta = f *  \psi_\delta =: f_\delta\quad\text{in}\quad B_{1-\delta}
\]
in the strong sense (since $u_\delta\in C^\infty(\R^n)$). Thanks to the previous steps (after a rescaling and covering argument, to have an estimate in $B_{1-\delta}$ instead of $B_1$), we now have 
\[
\begin{split}
\|u_\delta\|_{C^{2s+\alpha}(B_{1/2})} & \le C \left(\|u_\delta\|_{L^\infty_{2s-\eps}(\R^n)}+ [f_\delta]_{C^\alpha(B_{1-\delta})}\right)\\
& \le  C \left(\|u \|_{L^\infty_{2s-\eps}(\R^n)}+ [f]_{C^\alpha(B_1)}\right).
\end{split}
\]
We are also using that, if $u\in L^\infty_{2s-\eps}(\R^n)$ and $f\in C^\alpha(B_1)$, then we have $\|u_\delta\|_{L^\infty_{2s-\eps}(\R^n)}\le  2\|u\|_{L^\infty_{2s-\eps}(\R^n)}$ and $[f_\delta]_{C^\alpha(B_{1-\delta})}\le [f]_{C^\alpha(B_1)}$. Letting $\delta \downarrow 0$, since $u$ is locally bounded, $u_\delta\to u$ pointwise almost everywhere.  Moreover, by Arzel\`a-Ascoli, $u_\delta \to u$ in $C^\nu$ norm (up to a subsequence), with $\|u\|_{C^{2s+\alpha}(B_{1/2})}
\le \limsup_{\delta\downarrow 0}\|u_\delta\|_{C^{2s+\alpha}(B_{1/2})}$. 
Hence, $u\in C^{2s+\alpha}(\overline{B_{1/2}})$ and we get the desired estimate. \qedhere
\end{steps}
\end{proof}

Now, minor modifications of the previous proof allow us to show:

\begin{proof}[Proof of Theorem~\ref{thm-interior-linear-Lp}]
The proof follows in the same way as that of Theorem~\ref{thm-interior-linear-2} (using Proposition~\ref{prop:interior-linear}-\ref{it:int_lin_ii}) modifying \ref{step:main} as follows: using the same notation as in the proof of Theorem~\ref{thm-interior-linear-2} we now have (thanks to Lemma~\ref{lem:Lu}):
\[
\begin{split}
\left\|\L  (u-\eta u) \right\|_{L^\infty(B_{1/2})} &   \le C\|u\|_{L^\infty_{2s-\eps}(\R^n)}.
\end{split}
\]
The result is now obtained in the same way as done in the proof of Theorem~\ref{thm-interior-linear-2}. 
\end{proof}

 And we also get:
\begin{proof}[Proof of Proposition~\ref{prop-interior-linear}]
The proof is the same as that of Theorem~\ref{thm-interior-linear-2}, but in \ref{step:main} we use Lemma~\ref{lem:Lu_2} part \ref{it:lem_Lu2_ii} instead of part \ref{it:lem_Lu2_i}. 
\end{proof}

\subsection{Liouville's theorem for solutions with growth}
\index{Liouville's theorem!General operators!Solutions with growth} 
As an immediate consequence of the interior estimates, we obtain Liouville's theorem (Theorem~\ref{thm:Liouville}) but for solutions that may have some growth (up to a power $2s-\eps$):

\begin{cor}[Liouville's Theorem with growth]
\label{cor:Liouville_growth}
Let $s\in (0,1)$ and $\L\in \GL$. Let $u\in L^\infty_{\beta}(\R^n)$ for some $\beta\in [0, 2s)$ be a distributional solution to 
\[
\L u = 0 \quad \text{in}\quad \R^n.
\]
Then, $u(x) = a+b\cdot x$, with $b = 0$ if $\beta < 1$.
\end{cor}
\begin{proof}
Let $u_R(x) := \frac{u(Rx)}{R^{\beta}}$ for $R \ge 1$. Then
\[
\|u_R\|_{L^\infty_{\beta}(\R^n)} = \left\|\frac{u(Rx)}{R^{\beta}+|Rx|^{\beta}}\right\|_{L^\infty(\R^n)} \le \|u\|_{L^\infty_{\beta}(\R^n)}.
\]
By Theorem~\ref{thm-interior-linear-Lp}, we have for any $\beta' < 2s$,
\[
[u_R]_{C^{\beta'}(B_{1/2})} = R^{\beta'-\beta} [u]_{C^{\beta'}(B_{R/2})} \le C  \|u\|_{L^\infty_{\beta}(\R^n)}.
\]
Choosing $\beta' > \beta$, and letting $R\to \infty$, we obtain $[u]_{C^{\beta'}(\R^n)} = 0$ and hence the desired result (if $\beta < 1$ we can choose $\beta' < 1$, and deduce that $u$ is constant). 
\end{proof}

It is also possible to prove higher-order Liouville theorems for solutions with polynomial growth. In that case, however, one needs to define a gene\-ralized notion of solution that allows functions with arbitrary (polynomial) growth; see \cite{DSV2}. See also \cite{ADEJ20, GK23, FW16} for the highest level of generality under which Liouville's theorem holds for general L\'evy operators.

\subsection{The strong maximum principle}
\label{ssec:strongmaximumprinciple}

We finish this section by proving the strong maximum principle for general operators $\L \in \GL$.

\begin{thm}
\label{thm:strong_maximum_principle}\index{Strong maximum principle}
Let $s\in (0, 1)$, let $\L \in\GL$, and let $\Omega\subset \R^n$ be any bounded domain. Given the equation
\begin{equation}
\label{eq:smp}
\left\{
\begin{array}{rcccll}
\L u & = & f & \ge &0&\quad\text{in}\quad \Omega\\
u & \ge & 0 &  & & \quad\text{in}\quad \R^n\setminus \Omega,
\end{array}
\right.
\end{equation}
let us assume that one of the following situations occurs
\begin{enumerate}[(a)]

\item either $u$ is such that \eqref{eq:u_weak_sol} holds and it satisfies \eqref{eq:smp} in the weak sense for some $f\in L^p(\Omega)$ for $p \ge\frac{2n}{n+2s}$ and $n > 2s$;

\item or $u\in L^\infty_{2s-\eps}(\R^n)\cap C(\overline{\Omega})$ for some $\eps >0 $ satisfies \eqref{eq:smp} in the distributional sense, for some $f\in L^\infty_{\rm loc}(\Omega)$. 

\end{enumerate}
Then, either $u > 0$ in $\Omega$ or $u\equiv 0$  in $\Omega$. 
\end{thm}

\begin{proof}
By the corresponding maximum principles (see Lemmas~\ref{lem:max_principle_G_w}, and \ref{lem:max_principle_G_d}) we already know that $u \ge 0$ in $\Omega$. We now divide the proof into two steps:
\begin{steps}

\item \label{it:step1_smp} Let us assume first that $f\in L^\infty_{\rm loc}(\Omega)$ and $u\in L^\infty_{2s-\eps}(\R^n)$. Then weak solutions are distributional (see Lemma~\ref{lem:weak_distr}) and thanks to the interior regularity estimates from Theorem~\ref{thm-interior-linear-Lp}, we have that $u\in C^{2s-\eps}(\Omega)$ for any $\eps > 0$.

Now, if $u\not\equiv 0$ in $\Omega$, the set $\Omega^+ :=\{u > 0\}\cap \Omega$ is open and nonempty. We want to show that $\Omega^+ = \Omega$. Suppose that this is not the case: this means that there exists some $B_r(y_\circ)\subset \Omega^+$ such that for some $x_\circ\in \partial B_r(y_\circ)$, $u(x_\circ) = 0$ and $x_\circ\in \Omega$ (see Figure~\ref{fig:03}).

\begin{figure}
\centering
\makebox[\textwidth][c]{\includegraphics[scale = 1]{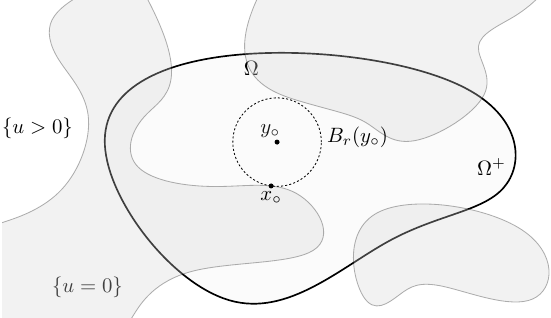}}
\caption{\label{fig:03} Graphical representation the setting in the contradiction argument from \ref{it:step1_smp}.}
\end{figure}

We consider the function 
\[
\psi(x) := (r^2-|x-y_\circ|^2)^\beta_+.
\]
By Lemma~\ref{lem-radialsubsol} (rescaled) we know that there is $\beta < 2s$ such that
\[
\L \psi \le -c \quad\text{in}\quad B_r(y_\circ) \setminus B_{r-\eta}(y_\circ) \quad\text{for some}\quad c, \eta > 0. 
\]
 Let $0<c_* = \min \{u(x) : x\in B_{r-\eta}(y_\circ)\}$. We have:
\[
\left\{
\begin{array}{rcll}
u &\ge& c_*\psi& \quad\text{in}\quad B_{r-\eta}(y_\circ)\cup B_r^c(y_\circ)\\
\L u & > & \L (c_* \psi)& \quad\text{in}\quad B_{r}(y_\circ)\setminus B_{r-\eta}(y_\circ).
\end{array}
\right.
\]
By the corresponding comparison principle (since we can choose $\beta > s-\frac12$, by Lemma~\ref{lem:distHs}, $\psi$ is a weak solutions as well) we have that $u\ge c_*\psi$ in $\R^n$. 

However, if we let $x_t := (1-t)x_\circ + ty_\circ$ for $t\in [0, 1]$, then $u(x_0) = 0$ and 
\[
u(x_t) \ge c_*\psi(x_t) = c_*\big(r^2-r^2(1-t)^2\big)^\beta \ge c_*r^{2\beta} t^\beta(2-t)^{\beta}.
\]
Since $\beta < 2s$, this contradicts the fact that $u\in C^{2s-\eps}(\Omega)$ for all $\eps > 0$. Hence, $\Omega^+ = \Omega$ and the result follows.

\item   If $u\notin L^\infty_{2s-\eps}(\R^n)$ or $f\notin L^\infty_{\rm loc}(\Omega)$ and $u$ is a weak solution, we can define for $m\in \N$,
\[
u_m (x) := \min\{u(x), m\},
\]
and consider $\bar u_m$ to be the weak solution of
\[
\left\{
\begin{array}{rcll}
\L \bar u_m & = & \min\{1, f\} & \quad\text{in}\quad \Omega,\\
\bar u_m &  = & u_m & \quad\text{in}\quad \R^n\setminus \Omega,
\end{array}
\right.
\]
given by Theorem~\ref{thm:exist_weak_sol}. Observe that we can indeed apply Theorem~\ref{thm:exist_weak_sol} since the right-hand side is now bounded and 
\[
|u_m(x) - u_m(y)|\le |u(x) - u(y)|\quad\text{for all}\quad x, y \in \R^n,
\]
so $\langle u_m, u_m\rangle_{K;\Omega}\le \langle u, u\rangle_{K;\Omega}< \infty$.

By the  maximum principle (Lemma~\ref{lem:max_principle_G_w}), $0\le \bar u_m \le u$. Now we can use the first step  for each $m\in \N$,  and deduce that either $\bar u_m > 0$ in $\Omega$ for some $m\in \N$ (in which case $u > 0$ in $\Omega$ and we are done), or $\bar u_m \equiv 0$ in $\Omega$ for all $m\in \N$ (in which case we would like to deduce that $u\equiv 0$ in $\Omega$). In such situation, $\L \bar u_m \equiv 0$ in $\Omega$ (since $\bar u_m \ge 0$ in $\R^n\setminus\Omega$) and thus $\L u = \L \bar u_m = 0$ in $\Omega$ for all $m\in \N$. In particular, $\bar u_m$ minimizes the energy \eqref{eq:energy_weak} (with $f \equiv 0$) among functions with the same boundary datum, and we have 
\[
\langle \bar u_m, \bar u_m \rangle_{K;\Omega} \le \langle u_m, u_m\rangle_{K;\Omega} \le \langle u, u\rangle_{K;\Omega} < \infty.
\]

We can now apply the monotone convergence theorem letting $m \to \infty$ (and using that $\bar u_m \equiv 0$ in $\Omega$) to deduce that 
\[
\langle u\chi_{\Omega^c}, u\chi_{\Omega^c} \rangle_{K;\Omega}\le \langle u, u\rangle_{K;\Omega}. 
\]
Since $u$ is the unique minimizer of the energy with prescribed exterior datum (by Theorem~\ref{thm:exist_weak_sol}), we obtain that $u = u \chi_{\Omega^c}$ and thus $u \equiv 0$ in $\Omega$, as we wanted to see. \qedhere 
\end{steps}
\end{proof}

\begin{rem}
Observe that, contrary to what occurs with the fractional Laplacian (see Lemma~\ref{lem:max_principle_s}), if $u$ is not strictly positive everywhere in $\Omega$ then $u \equiv 0$ in $\Omega$, instead of $u \equiv 0$ in $\R^n$. Indeed, since the kernel $K$ could vanish in open sets, there are operators $\L \in \GL$ and functions $u \ge 0$ such that $u \equiv 0$ in $\Omega$, $\L u = 0$ in $\Omega$, but $u > 0$ somewhere in $\R^n\setminus \Omega$. 
\end{rem}

\section{Equations with $x$-dependence}

\label{sec:x-dependence}

The results we have established so far give a quite complete understanding of the interior regularity of solutions to linear and \emph{translation invariant} equations of order $2s$.

The next very natural question is to understand what happens when the operators under consideration are \emph{not} translation invariant, i.e., they have $x$-dependence.
For this, we have two different classes of operators: in non-divergence form  and in  divergence form:

\begin{itemize}[leftmargin=*]
\item \index{Non-divergence-form equations} \emph{Non-divergence-form} operators are those of the form 
\[
\begin{split}
\L(u,x) & = {\rm P.V.} \int_{\R^n}\big(u(x)-u(x+y)\big)K(x,y)\, dy   \\
& = \frac12 \int_{\R^n}\big(2u(x)-u(x+y)-u(x-y)\big)K(x,y)\, dy,
\end{split}
\]
with
\[
K\geq0\qquad \textrm{and}\qquad K(x,y)=K(x,-y)\quad \textrm{for all}\quad x,y\in \R^n.
\]
They correspond (in the limiting case $s=1$) to operators of the type 
\[\L(u,x)= {\rm tr}(A(x)D^2 u)= \sum_{i,j=1}^n a_{ij}(x)\partial_{ij}u,\]
and arise naturally when studying fully nonlinear elliptic equations.
Equations involving these operators \emph{cannot} be studied by using weak solutions and energy functionals, and instead require the use of viscosity solutions (see Chapter \ref{ch:fully_nonlinear}). We can, however, establish a priori interior regularity estimates (Schauder estimates).

The kernels we consider, $K(x, y)\, dy$, can in general be measures $K(x, dy)$, that will satisfy \eqref{eq:nu_cond} for each $x\in \R^n$. In that case, the operators are given by \eqref{x-dependence-L-nu}. 

\item \index{Divergence-form equations} \emph{Divergence-form} operators, instead, are those of the form
\[
\begin{split}
\L (u,x) & = {\rm P.V.} \int_{\R^n}\big(u(x)-u(z)\big)K(x,z) \, dz\\
& = {\rm P.V.} \int_{\R^n}\big(u(x)-u(x+y)\big)K(x,x+y) \, dy,
\end{split}
\]
with 
\[
K\geq0\qquad \textrm{and}\qquad K(x,z)=K(z,x)\quad \textrm{for all}\quad x,z\in \R^n.
\]
They correspond (in the limiting case $s=1$) to operators of the type 
\[\L(u,x)= {\rm div}\big(A(x)\nabla u\big)= \sum_{i,j=1}^n \partial_i\big(a_{ij}(x)\partial_{j}u\big),\]
and arise naturally in the Calculus of Variations.
Equations involving these operators are studied by using weak solutions and integration by parts, and some results (such as the existence of solutions and maximum principles) may be established similarly to the ones we proved in this chapter.

Notice that, as before, the kernel $K(x,z)\, dz$ could actually be a general measure $K(x,dz)$.
In that case, the previous expressions are given by \eqref{divergence-form0}-\eqref{divergence-form1}. 
\end{itemize}

We refer to the works of Barrios-Figalli-Valdinoci \cite{BFV}, Imbert-Silvestre \cite{IS2}, Jin-Xiong \cite{JX}, and Serra \cite{Ser2} for Schauder estimates in non-divergence form, and to the work of Fall \cite{Fal} for a Schauder estimate in divergence form. Here, we will extend these Schauder-type estimates to a more general class of operators, following the presentation from \cite{FR5}.

\subsection{Schauder estimates for equations in non-divergence form}
\index{Non-divergence-form equations}

We consider first operators with $x$-dependence in non-divergence form, that is, of the type  
\begin{equation}\label{x-dependence-L-nu}
\begin{split}
\L(u,x)& := {\rm P.V.}\int_{\R^n}\big(u(x)-u(x+y)\big)K_x(dy) \\
& = \frac12 \int_{\R^n}\big(2u(x)-u(x+y)-u(x-y)\big)K_x(dy),
\end{split}
\end{equation}
(cf. \eqref{eq:Lu_nu}) where $(K_x)_{x\in \R^n}$ is a family of L\'evy measures satisfying \eqref{eq:nu_cond} and with ellipticity conditions \eqref{eq:Kellipt_gen_L}-\eqref{eq:Kellipt_gen_l} uniform in $x\in \R^n$. Namely, if for a given $x_\circ \in \R^n$ we denote $\L_{x_\circ}$ the translation invariant operator with L\'evy measure $K_{x_\circ}$ (i.e., the operator with ``frozen coefficients''), we consider in this section operators $\L$ of the form \eqref{x-dependence-L-nu} such that $\L_x\in \GL$ for all $x\in \R^n$.

Our goal is to prove Schauder-type estimates for such a class of operators. In order to do that, we   require a certain regularity in $x$ of the corresponding kernels. That is, we need to have, in some integral sense, ``H\"older continuous coefficients'': if $\alpha\in (0, 1]$ is fixed, we will assume  
\begin{equation}\label{x-dependence-L2}
\int_{B_{2\rho}\setminus B_\rho} \big|K_x(dy)-K_{x'}(dy)\big|  \leq M|x-x'|^\alpha \rho^{-2s},
\end{equation}
for all $x, x'\in \R^n$ and all $\rho > 0$, and for some $M > 0$.

As in Theorem \ref{thm-interior-linear-2}, we also need some regularity of each $K_x$ in the $y$-variable, namely  
\begin{equation}\label{x-dependence-L3}
\sup_{x\in\R^n} \, [K_x]_\alpha \leq M.
\end{equation}
(Recall \eqref{Calpha-assumption}.)

Under the previous assumptions, we then have the following a priori estimate (see Remark~\ref{rem:visc_x_dep} for a discussion on the a priori regularity):

\begin{thm}\label{thm-interior-linear-x}
Let $s\in (0,1)$, and let $\L$ be an operator of the form \eqref{x-dependence-L-nu} with $\L_x\in \GL$ for all $x\in \R^n$. Suppose, moreover, that $\L$ satisfies the regularity assumptions \eqref{x-dependence-L2}-\eqref{x-dependence-L3} for some $\alpha\in(0,1]$  such that $2s+\alpha\notin\N$, and $M > 0$.  

Let $f\in C^\alpha(B_1)$, and let $u\in C_{\rm loc}^{2s+\alpha}(B_1) \cap L^\infty_{2s-\eps}(\R^n)$ for some $\eps > 0$ be any solution of
\[
\L(u,x) = f \quad\text{in}\quad B_1.
\]

Then, 
\[
\|u\|_{C^{2s+\alpha}(B_{1/2})} \le C\left(\|u\|_{L^\infty_{2s-\eps}(\R^n)} + \|f\|_{C^\alpha(B_1)}\right)
\]
for some $C$ depending only on $n$, $s$, $\alpha$, $\eps$, $\lambda$, $\Lambda$, and $M$. 
\end{thm}

Notice that, as in Theorem~\ref{thm-interior-linear-2}, we need the kernel to have the same degree of regularity (in the $y$ variable) as the right-hand side in order to gain $2s$ derivatives. 
As we saw in Lemma~\ref{lem:counter}, without this assumption, the previous estimate is false.

Still, for solutions satisfying a global regularity assumption of the type $u\in C^\gamma(\R^n)$, one expects to be able to remove (or at least, weaken) such an assumption (cf. Proposition~\ref{prop-interior-linear}). This is what we do in the following result, which holds under the weaker regularity assumption (in $x$) 
\begin{equation}\label{x-dependence-L4}
\sup_{\begin{subarray}{c} [\phi]_{C^\gamma(\R^n)}\leq 1 \\ \phi(0)=0 \end{subarray}} 
\left|\int_{B_{2\rho}\setminus B_\rho} \phi\, d(K_x - K_{x'})\right| \leq M|x-x'|^\alpha\rho^{\gamma-2s},
\end{equation}
for all $x,x'\in \R^n$ and all $\rho>0$. The expression on the right-hand side of \eqref{x-dependence-L4} can be understood as some kind of weighted Kantorovich norm measuring the distance between the measures $K_x$ and $K_{x'}$ (see \cite{Han92, Han99}). Notice that, in the following, no regularity in the $y$ variable is assumed when $\alpha\leq \gamma$.

\begin{prop}\label{prop-interior-linear-x}
Let $s\in (0,1)$,  $\gamma\in (0, 2s)$, $\alpha\in (0,1]$ such that $2s+\alpha\notin \N$, and let $\L$ be an operator of the form \eqref{x-dependence-L-nu}  such that $\L_x\in \GL$ for all $x\in \R^n$. Suppose, moreover, that 
 $\L$ satisfies \eqref{x-dependence-L4} for some $M > 0$. 
 In case $\alpha>\gamma$, assume in addition that $\sup_{x\in B_1} [K_x]_{\alpha-\gamma} \leq M$.

Let $f\in C^\alpha(B_1)$, and let $u\in C^{2s+\alpha}(B_1)\cap C^{\gamma}(\R^n)$ be any solution of
\[
\L(u,x) = f \quad\text{in}\quad B_1.
\]

Then
\[
\|u\|_{C^{2s+\alpha}(B_{1/2})} \le C\left(\|f\|_{C^\alpha(B_1)} + \|u\|_{C^\gamma(\R^n)}\right)
\]
for some $C$ depending only on $n$, $s$, $\alpha$, $\gamma$, $\lambda$, $\Lambda$, and $M$.
\end{prop}

Notice that the assumption \eqref{x-dependence-L4} is scale invariant, and in case of stable operators $\GLh$ (of the form \eqref{eq:Lu_stab}) it is equivalent to 
\begin{equation}\label{x-dependence-L4_h} 
\sup_{\begin{subarray}{c} \|\phi\|_{C^\gamma(\mathbb S^{n-1})}\leq 1  \end{subarray}}
\left|\int_{\mathbb S^{n-1}} \phi\, d(\zeta_x - \zeta_{x'})\right| \leq M|x-x'|^\alpha.\end{equation}
Moreover, it is the minimal scale-invariant assumption that  ensures the property
\[\big\|\L_{x_1} w-\L_{x_2}w\big\|_{L^\infty(B_{3/4})} \leq M|x_1-x_2|^\alpha \big(\|w\|_{C^{2s+\alpha}(B_1)}+\|w\|_{C^{\gamma}(\R^n)}\big),\]
for all $w\in C^{2s+\alpha}(B_1)\cap C^{\gamma}(\R^n)$, and thus we expect it to be the minimal assumption under which these Schauder-type estimates hold.

It is interesting to notice that, for $\alpha\leq \gamma$, \eqref{x-dependence-L4} allows for completely singular L\'evy measures.
For example, one could have operators of the type (cf.~Example~\ref{ex:four_symb_ex})
\[\sum_{i=1}^n \big(-\partial_{\textbf{v}_i(x)}^2\big)^s,\]
where the directions $\textbf{v}_i$ are smooth functions of $x$ satisfying ${\rm det}(\textbf{v}_i)_i\neq0$.
  Notice also that by choosing $\alpha\le \gamma$, \eqref{x-dependence-L4_h} allows us to consider H\"older continuous functions~$\textbf{v}_i$ (with exponent $\max\{\alpha, \alpha/\gamma\}$).    Such operators are not covered, for example, by the stronger assumption \eqref{x-dependence-L2}.

\begin{rem} 
\label{rem:visc_x_dep}
The estimates we prove in this section are \textit{a priori} estimates, in the sense that we assume the solutions to be $C^{2s+\alpha}$.
This is because for the equations we are considering here, one cannot define weak nor distributional solutions in general.
Still, using the theory of viscosity solutions introduced in Chapter~\ref{ch:fully_nonlinear}, one can actually show by a regularization procedure analogous to the one in Section~\ref{sec:regularization} that the previous estimates are also valid for viscosity solutions (see, e.g., \cite{Fer23}). 
\end{rem}

\begin{proof}[Proof of Theorem \ref{thm-interior-linear-x}]
Let us divide the proof into three  steps.
\begin{steps}
\item\label{step:1nondiv} We start with an initial reduction. By a standard covering argument, we may prove the estimate in $B_{r_\circ/2}$ instead of $B_{1/2}$, where $r_\circ$ is a small fixed constant to be chosen later (depending only on $n$, $s$, $\alpha$, $\eps$, $\lambda$, $\Lambda$, and $M$).

In order to do that, we define $u_\circ(x) := u(r_\circ x)$, so that if $\L$ has kernels $K_x(dy)$, and we consider $\L^{r_\circ}$ to be the operator with kernels $K^{r_\circ}_x(dy) = r_\circ^{ 2s} K_{r_\circ x}(r_\circ \, dy)$   (which has the same ellipticity constants as $K_x$, cf. Remark~\ref{rem:Scale_invariance}), then $u_\circ$ satisfies
\[
\L^{r_\circ}(u_\circ, x) = f_\circ(x) := r_\circ^{2s} f(r_\circ x) \quad\text{in}\quad B_1.
\]

If we can now prove the desired estimate for $u_\circ$,
\[
\|u_\circ\|_{C^{2s+\alpha}(B_{1/2})} \le C\left(\|u_\circ\|_{L^\infty_{2s-\eps}(\R^n)} + \|f_\circ\|_{C^\alpha(B_1)}\right),
\]
we will be done, since 
\[
\begin{split}
r_\circ^{2s+\alpha} \|u\|_{C^{2s+\alpha}(B_{r_\circ/2})}& \le \|u_\circ\|_{C^{2s+\alpha}(B_{1/2})}\\
& \le C\left(\|u_\circ\|_{L^\infty_{2s-\eps}(\R^n)} + \|f_\circ\|_{C^\alpha(B_1)}\right)\\
&   \le C\left(\|u\|_{L^\infty_{2s-\eps}(\R^n)} + \|f\|_{C^\alpha(B_1)}\right).
\end{split}
\]
In particular, since $r_\circ$ will be fixed universally, this will yield the estimate in $B_{r_\circ/2}$, and after a covering (and rescaling) argument, in the whole $B_{1/2}$. 

The advantage of this rescaling is that, now, the new operator $\L^{r_\circ}$ satisfies a new condition \eqref{x-dependence-L2} with constants depending on $r_\circ$ (since we have ``expanded'' the space by a factor $r_\circ^{-1}$):
\begin{equation}\label{x-dependence-L2-2}
\sup_{\rho > 0}\sup_{x, x'\in \R^n} \rho^{2s}\int_{B_{2\rho}\setminus B_\rho} \frac{\big|K^{r_\circ}_x(dy)-K^{r_\circ}_{x'}(dy)\big|}{|x-x'|^\alpha}   \leq Mr^\alpha_\circ =:\delta,
\end{equation}
whereas condition \eqref{x-dependence-L3} remains the same (that is, with the same constant);
\begin{equation}\label{x-dependence-L3-2}
\sup_{x\in\R^n} \, [K^{r_\circ}_x]_\alpha \leq M.
\end{equation} 

 In all, up to replacing $u$, $\L$, and $f$,  by $u_\circ$, $\L^{r_\circ}$, and $f_\circ$, we can assume without loss of generality that condition \eqref{x-dependence-L2} holds with $M = \delta> 0$ arbitrarily small, but fixed, that will be chosen universally:
\begin{equation}\label{x-dependence-L2_delta}
\sup_{\rho > 0}\sup_{x, x'\in \R^n} \rho^{2s}\int_{B_{2\rho}\setminus B_\rho} \frac{\big|K_x(dy)-K_{x'}(dy)\big|}{|x-x'|^\alpha}   \leq \delta,
\end{equation}

\item \label{step:step2prove} Let $\L_0 \in \G_s(\lambda, \Lambda;\alpha)$ be the operator $\L_x$ corresponding to $x=0$.
Then, by the regularity estimates for translation invariant equations, Theorem~\ref{thm-interior-linear-2}, we have
\[
\|u\|_{C^{2s+\alpha}(B_{1/4})} \le C\left(\|u\|_{L^\infty_{2s-\eps}(\R^n)} + \|\L_0u\|_{C^\alpha(B_{1/2})}\right).
\]
Moreover, we also have
\[\|\L_0u\|_{C^\alpha(B_{1/2})} \leq \|f\|_{C^\alpha(B_{1/2})} + \|\L(u,\cdot) - \L_0u\|_{C^\alpha(B_{1/2})}.\]
We now claim that 
\begin{equation}\label{ngfhhghg} 
\|\L(u,\cdot) - \L_0u\|_{C^\alpha(B_{1/2})} \leq C\left(\delta\|u\|_{C^{2s+\alpha}(B_1)} + \|u\|_{L^\infty_{2s-\eps}(\R^n)}\right).
\end{equation}

For this, we define, for any $x\in B_{1/2}$ fixed, a new operator 
\[
\MR_x := \L_x - \L_0,
\]
so that $\L(u, x) - \L_0 u = \MR_x u$. Observe that $\MR_x$ has kernel $R_x(dy) := K_x(dy) - K_0(dy)$, but it does not belong to $\GL$ (since it is not positive). 

We have, nonetheless, that by assumption \eqref{x-dependence-L2_delta} applied with $x' = 0$,
\[
\rho^{2s}\int_{B_{2\rho\setminus B_\rho}} |R_x(dy)| \le \delta |x|^\alpha \le \delta\qquad\text{for all}\quad \rho > 0. 
\]
In particular, we can apply Lemma~\ref{lem:Lu} (see Remark~\ref{rem:nonneg_kernels}) and deduce 
\[
|\MR_x u(x)| \le C \delta \left( \|u\|_{C^{2s+\alpha}(B_1)} + \|u\|_{L^\infty_{2s-\eps}(\R^n)}\right),
\]
for any $x\in B_{1/2}$, which gives the $L^\infty$ bound on \eqref{ngfhhghg}. 

On the other hand, given any $x_1, x_2\in B_{1/2}$, we now want to bound the difference
\[
|\MR_{x_1} u(x_1) - \MR_{x_2} u(x_2)|\le |\MR_{x_1} u(x_1) - \MR_{x_1} u(x_2)|+|\MR_{x_1} u(x_2) - \MR_{x_2} u(x_2)|.
\]
For the first term, we  use Lemma~\ref{lem:Lu_2}-\ref{it:lem_Lu2_i} (together with Remarks~\ref{rem:nonneg_kernels} and~\ref{rem:Lu_2}) with operator $\MR_{x_1}$ fixed, to deduce
\[
|\MR_{x_1} u(x_1) - \MR_{x_1} u(x_2)| \le C|x_1-x_2|^\alpha\left(\delta \|u\|_{C^{2s+\alpha}(B_1)}+ 2M \|u\|_{L^\infty_{2s-\eps}(\R^n)}\right).
\]
We have also used here that, by assumption \eqref{x-dependence-L3}, $[R_{x_1}]_\alpha\le 2M$. 

For the second term, we can define yet another operator
\[
\tilde\MR_{x_1,x_2} := \MR_{x_1} - \MR_{x_2},
\]
which has kernel $\tilde R_{x_1,x_2} (y) := K_{x_1}(y) - K_{x_2}(y)$ and satisfies,  by \eqref{x-dependence-L2_delta},
\[
\rho^{2s}\int_{B_{2\rho}\setminus B_\rho}|\tilde R_{x_1,x_2}(dy)|  \le \delta |x_1-x_2|^\alpha. 
\]
Thus, we can  apply again Lemma~\ref{lem:Lu} and Remark~\ref{rem:nonneg_kernels} to deduce
\[
\begin{split}
|\MR_{x_1} u(x_2) - \MR_{x_2} u(x_2)| & = |\tilde \MR_{x_1,x_2} u(x_2)|\\
& \le C\delta |x_1-x_2|^\alpha (\|u\|_{C^{2s+\alpha}(B_1)}+\|u\|_{L^\infty_{2s-\eps}(\R^n)}). 
\end{split}
\]

Putting everything together, we have obtained \eqref{ngfhhghg}, where $\delta$ can still be chosen. 

\item In all, we have shown 
\[
\|u\|_{C^{2s+\alpha}(B_{1/4})} \le C\left(\delta\|u\|_{C^{2s+\alpha}(B_1)} + \|u\|_{L^\infty_{2s-\eps}(\R^n)} + \|f\|_{C^\alpha(B_{1/2})}\right)
\]
for all $u\in C^{2s+\alpha}(B_1)\cap L^\infty_{2s-\eps}(\R^n)$, where $\L(u, \cdot) = f$. We now choose $\delta$ sufficiently small (depending only on $n$, $s$, $\alpha$, $\eps$, $\lambda$, $\Lambda$, and $M$) so that we can repeat the argument in \ref{step:SAL} of the proof of Theorem~\ref{thm-interior-linear-2} (on page \pageref{step:SAL}), to deduce that (after a covering argument)
\[
\|u\|_{C^{2s+\alpha}(B_{1/2})} \le C\left(\|u\|_{L^\infty_{2s-\eps}(\R^n)} + \|f\|_{C^\alpha(B_{1})}\right),
\]
for all $u\in C^{2s+\alpha}(B_1)\cap L^\infty_{2s-\eps}(\R^n)$, where $\L(u, \cdot) = f$.

From \ref{step:1nondiv}, this proves the estimate in a universally small ball $B_{r_\circ}$, and after a further covering and rescaling, this shows the desired result. \qedhere
\end{steps}
\end{proof}

To prove Proposition  \ref{prop-interior-linear-x} we need the following.

\begin{lem}\label{lem-sdhhh}
Let $s\in(0,1)$,   $\alpha\in(0,1]$, and $\gamma\in(0,2s)$.
Let $\L_1,\L_2\in \GL$, with L\'evy measures $K_1,K_2$, be such that
\begin{equation}\label{iserbgieub}
\sup_{\begin{subarray}{c} [\phi]_{C^\gamma(\R^n)}\leq 1 \\ \phi(0)=0 \end{subarray}} 
\left|\int_{B_{2\rho}\setminus B_\rho} \phi \, d(K_1 - K_2)\right| \leq \theta \rho^{\gamma-2s}\qquad\text{for all}\quad \rho > 0,
\end{equation}
for some $\theta > 0$. 
 In case $\alpha>\gamma$, assume in addition that $[\L_i]_{\alpha-\gamma} \leq M$ for $i= 1,2$, and some $M > 0$. 

Then, if $w\in C^{2s+\alpha}(B_1)\cap C^{\gamma}(\R^n)$ we have
\begin{equation}
\label{eq:Linfty_nondivform}
\big\|(\L_1-\L_2)w\big\|_{L^\infty(B_{1/2})} \leq C 
\theta \big(\|w\|_{C^{2s+\alpha}(B_1)}+[w]_{C^{\gamma}(\R^n)}\big)
\end{equation}
as well as  
\begin{equation}
\label{eq:Calpha_nondivform}
\big[(\L_1-\L_2)w\big]_{C^\alpha(B_{1/2})}  \le 
\left\{
\begin{array}{ll}
  C 
\theta \big( \|w\|_{C^{2s+\alpha}(B_1)}+[w]_{C^{\gamma}(\R^n)}\big)& \quad\text{if}\quad \alpha \le \gamma\\[0.2cm]
 C 
\big(\theta \|w\|_{C^{2s+\alpha}(B_1)}+M [w]_{C^{\gamma}(\R^n)}\big)& \quad\text{if}\quad \alpha > \gamma.
\end{array}
\right.
\end{equation}
The constant $C$ depends only on $n$, $s$, $\alpha$, and $\gamma$.
\end{lem}

\begin{proof}
We denote 
\[
\MR := \L_1-\L_2.
\]
We divide the proof into four steps.
\begin{steps}
 \item \label{step:firsttodo}
We first prove an $L^\infty$ bound for $\MR w$, \eqref{eq:Linfty_nondivform}.
 Dividing by a constant if necessary, we may assume $\|w\|_{C^{2s+\alpha}(B_1)}+[w]_{C^\gamma(\R^n)}\leq 1$, and by taking $\alpha$ smaller if necessary, we may also assume $2s+\alpha < 2$ and $2s+\alpha\neq 1$.

Let $x_\circ\in B_{1/2}$.
Taking $\phi(y) = 2w(x_\circ)-w(x_\circ+y)-w(x_\circ-y)$ in \eqref{iserbgieub}, we find that for any  $\rho>0$ (using that $[w]_{C^\gamma(\R^n)}\le 1$), 
\begin{equation}
\label{eq:firstineqlem}
\left|\int_{B_{2\rho}\setminus B_\rho} \hspace{-0.35cm} \big(2w(x_\circ)-w(x_\circ+y)-w(x_\circ-y)\big)  (K_1 - K_2)(dy)\right| \leq 4\theta\rho^{\gamma-2s}.
\end{equation}
On the other hand, since $w\in C^{2s+\alpha}(B_1)$ with $\|w\|_{C^{2s+\alpha}(B_1)}\le 1$, we have (see  Lemma~\ref{it:H7_gen})
\[
\begin{split}
\left\|\phi\right\|_{L^\infty(B_{2\rho})}& \le C \rho^{2s+\alpha},\\
\|\phi\|_{C^{2s+\alpha}(B_{2\rho})}& \le 4,
\end{split}
\] 
for any $\rho <\frac14$. 
Applying Lemma~\ref{lem:interp_mult} to the function $\phi$ in $B_{2\rho}$ we obtain 
\[\left[ \frac{\phi}{\rho^{2s+\alpha-\gamma}} \right]_{C^\gamma(B_{2\rho}\setminus B_\rho)} \leq C,\]
to get, by \eqref{iserbgieub},
\begin{equation}
\label{eq:secineqlem}
\left|\int_{B_{2\rho}\setminus B_\rho} \hspace{-0.35cm} \big(2w(x_\circ)-w(x_\circ+y)-w(x_\circ-y)\big) (K_1 - K_2)(dy)\right| \leq C\theta\rho^{\alpha}.
\end{equation}
Using the first inequality, \eqref{eq:firstineqlem}, for $\rho=2^{k}$, $k=-1,0,1,2,3,...$, and the second one, \eqref{eq:secineqlem}, for $k=-2,-3,...$,  and summing a geometric series, we deduce that
\[
\left|\int_{\R^n} \big(2w(x_\circ)-w(x_\circ+y)-w(x_\circ-y)\big)\,  (K_1 - K_2)(dy)\right| \leq C\theta.
\]
This proves the $L^\infty$ bound for $\MR w$.
\item 
Let us now show the $C^\alpha$ bounds, \eqref{eq:Calpha_nondivform}.
For this, we argue as in the proof of Lemma~\ref{lem:Lu_2} (and Remarks~\ref{rem:Lu_2} and \ref{rem:nonneg_kernels}), and split $w=u_1+u_2$, with $u_1:=\eta w$ and $\eta\in C^\infty_c(B_1)$ such that $\eta \ge 0$, $\eta \equiv 0$ in $\R^n\setminus B_{3/4}$ and $\eta \equiv 1$ in $B_{2/3}$.
Notice that $\|u_1\|_{C^{2s+\alpha}(\R^n)}\le C\|w\|_{C^{2s+\alpha}(B_1)}$ and $[u_2]_{C^\gamma(\R^n)} \leq C[w]_{C^\gamma(\R^n)}$.

We prove first the bound for $\MR u_1$ in the case $2s+\alpha \le 2$. Dividing by a constant if necessary, we assume $\|u_1\|_{C^{2s+\alpha}(\R^n)}\leq 1$. Observe that the $L^\infty$ bound for $\MR u_1$ follows by \ref{step:firsttodo}. We now want to bound the $C^\alpha$ seminorm, and more precisely, we will bound 
\begin{equation}
\label{eq:sameasbefore222}
\big|\MR u_1(x_\circ) + \MR u_1(-x_\circ) - 2 \MR u_1(0) \big| \le C \theta r^\alpha.
\end{equation}

Let $x_\circ\in B_{1/2}$ be fixed, and $r := |x_\circ|$. We split
\[
\begin{split}
\MR u_1(x_\circ) & = \frac12 \int_{B_r}\big(2u_1(x_\circ)-u_1(x_\circ+y)-u_1(x_\circ-y)\big)(K_1-K_2)(dy)\\
& \hspace{-5mm} +  \frac12 \int_{\R^n\setminus B_r}\big(2u_1(x_\circ)-u_1(x_\circ+y)-u_1(x_\circ-y)\big)(K_1-K_2)(dy).
\end{split}
\]
Then, since $2s+\alpha \le 2$ and by interpolation as in the previous step using Lemma~\ref{lem:interp_mult} (on expressions of the form \eqref{eq:toused124}), we have  
\begin{equation}
\label{eq:theexp}
\begin{split}
\big[u_1(x_\circ+\cdot\,)+u_1(x_\circ-\cdot\,)-2u_1(x_\circ)\big]_{C^\gamma(B_{2\rho}\setminus B_\rho)} &\le  C\rho^{2s+\alpha-\gamma},\\
\big[u_1(x_\circ\pm \cdot\,)+u_1(-x_\circ\pm \cdot\,)-2u_1(\pm \,\cdot\,)\big]_{C^\gamma(\R^n)} &\le  Cr^{2s+\alpha-\gamma}.
\end{split}
\end{equation}
If we now denote $\delta^2_h v(x)$ the second order centered increments, 
\[
  \delta^2_h v(x) = \frac{v(x+h)+v(x-h)}{2} - v(x),
\]
and 
\[
\MR u_1(x_\circ) + \MR u_1(-x_\circ) - 2 \MR u_1(0)  = 2\int_{\R^n} \phi(x_\circ, y)\, (K_1-K_2)(dy)
\]
with 
\[
\phi(x_\circ, y) := \delta^2_{y} u_1(x_\circ) + \delta^2_{y} u_1(-x_\circ) -2 \delta^2_{y} u_1(0),
\]
then the expressions \eqref{eq:theexp} imply, by definition \eqref{iserbgieub},
\[
\left|\int_{B_{2\rho}\setminus B_\rho} \phi(x_\circ,y)(K_1-K_2)(dy) \right| \le C\theta \min\{\rho^\alpha, \rho^{\gamma-2s}r^{2s+\alpha-\gamma}\}.
\]
Summing first for $\rho\in (0, r)$ (and taking the first argument in the $\min$), and summing then for $\rho > r$ (and taking the second argument in the $\min$) we obtain \eqref{eq:sameasbefore222}.

\item Let us now show \eqref{eq:sameasbefore222} in the case  $2<2s+\alpha<3$. By Lemma~\ref{lem:A_imp_2}-\ref{it:A:3} and the interpolation in Lemma~\ref{lem:interp_mult} (cf. \eqref{eq:2salpha1}) we have, on the one hand,
\begin{equation}
\label{eq:2salpha1_2}
\begin{split}
& \big[u_1(x_\circ+\cdot\,)\hspace{-0.5mm}+\hspace{-0.5mm}u_1(x_\circ-\cdot\,)\hspace{-0.5mm}-\hspace{-0.5mm}2u_1(x_\circ)\hspace{-0.5mm}-\hspace{-0.5mm}u_1(\,\cdot\,)\hspace{-0.5mm}-\hspace{-0.5mm}u_1(-\,\cdot\,)\hspace{-0.5mm}+\hspace{-0.5mm}2u_1(0)\big]_{C^\gamma(B_{2\rho}\setminus B_\rho)} \leq \\
& \qquad \le C\rho^{2-\gamma} r^{2s+\alpha-2} + C \rho^{2\left(1-\frac{\gamma}{2s+\alpha}\right)} r^{(2s+\alpha-2)\left(1-\frac{\gamma}{2s+\alpha}\right)}=:I_1.
\end{split}
\end{equation}
On the other hand, we want to find an appropriate bound for 
\begin{equation}
\label{eq:2salpha2_2}
\begin{split}
& \big[u_1(x_\circ+\cdot\,)\hspace{-0.71mm}+\hspace{-0.71mm}u_1(-x_\circ+\cdot\,)\hspace{-0.71mm}-\hspace{-0.71mm}2u_1(\,\cdot\,)\hspace{-0.71mm}-\hspace{-0.71mm}u_1(x_\circ)\hspace{-0.71mm}-\hspace{-0.71mm}u_1(-x_\circ)\hspace{-0.71mm}+\hspace{-0.71mm}2u_1(0)\big]_{C^\gamma(B_{2\rho}\setminus B_\rho)} =\\
& \qquad\le 2[\delta_{x_\circ}^2 u_1(\,\cdot\,)]_{C^\gamma(B_{2\rho})}.
\end{split}
\end{equation}
We do so by separating into three possible cases according to the value of~$\gamma$:
\begin{itemize}
\item If $\gamma \le 2s+\alpha -2 < 1$, then we have 
\[
\frac{|\delta_{x_\circ}^2 u_1(y) - \delta_{x_\circ}^2 u_1(y')|}{|y-y'|^\gamma} \le C |x_\circ|^2 |y-y'|^{2s+\alpha-2-\gamma}\le Cr^2\rho^{2s+\alpha-2-\gamma},
\]
for all $y, y'\in B_{2\rho}$, where we have used that $2s+\alpha-2-\gamma \ge 0$, and $[u_1]_{C^{2s+\alpha}(\R^n)}\le 1$, together with Lemma~\ref{lem:A_imp_2}-\ref{it:A:3}.
\item If $2s+\alpha-2 < \gamma \le 1$, then we can use again Lemma~\ref{lem:A_imp_2}-\ref{it:A:3} but now with $t$ such that $t(2s+\alpha -2) + (1-t) = \gamma$, where $t\in [0, 1)$ by assumption on $\gamma$, to obtain 
\[
\frac{|\delta_{x_\circ}^2 u_1(y) - \delta_{x_\circ}^2 u_1(y')|}{|y-y'|^\gamma} \le C |x_\circ|^{2s+\alpha-\gamma} \le Cr^{2s+\alpha-\gamma}
\]
for all $y, y'\in B_{2\rho}$. 

\item Finally, if $\gamma > 1$ we use the second equation in Lemma~\ref{lem:A_imp_2}-\ref{it:A:1} applied to $\nabla u_1$ and with $t = \gamma-1$ to derive (since $[\nabla u_1]_{C^{2s+\alpha-1}(\R^n)}\le C$)
\[
\frac{|\delta_{x_\circ}^2 \nabla u_1(y) - \delta_{x_\circ}^2 \nabla u_1(y')|}{|y-y'|^{\gamma-1}} \le C |x_\circ|^{2s+\alpha-\gamma} \le Cr^{2s+\alpha-\gamma}
\]
for all $y, y'\in B_{2\rho}$. 
\end{itemize}

We thus have a bound for the expression \eqref{eq:2salpha2_2} of the form 
\begin{equation}
\label{eq:2salpha2_3}
[\delta_{x_\circ}^2 u_1(\,\cdot\,)]_{C^\gamma(B_{2\rho})}\le I_2^\gamma := \left\{\begin{array}{ll}
C r^2 \rho^{2s+\alpha-2-\gamma} & \text{if} \ \gamma\le 2s+\alpha -2\\
 C r^{2s+\alpha-\gamma} & \text{if} \ \gamma > 2s+\alpha -2. 
\end{array}
\right.
\end{equation}

In particular, \eqref{eq:2salpha1_2}-\eqref{eq:2salpha2_3} imply now 
\[
\left|\int_{B_{2\rho}\setminus B_\rho} \phi(x_\circ,y)(K_1-K_2)(dy) \right| \le C\theta \rho^{\gamma-2s}\min\{I_1,I_2^\gamma\}.
\]
We now sum as before for $\rho = 2^k$ and $k\in \Z$, separating between $\rho < r$ and $\rho > r$. That is, we consider 
\begin{equation}
\label{eq:wecando}
\big|\MR u_1(x_\circ) + \MR u_1(-x_\circ) - 2 \MR u_1(0) \big| \le C\theta \sum_{\substack{\rho = 2^k\\ \rho < r}}\rho^{\gamma-2s} I_1+ C\theta \sum_{\substack{\rho = 2^k\\ \rho \ge r}}\rho^{\gamma-2s} I_2^\gamma.
\end{equation}

For the first term in the sum, we observe that the exponents of $\rho$ are positive, since they are 
\[
2-\gamma + \gamma - 2s > 0
\]
and 
\[
2\left(1-\frac{\gamma}{2s+\alpha}\right) + \gamma -2s > 0 ,
\]
where we are using that $2s+\alpha > 2$ and $2 > 2s$. In the second term of \eqref{eq:wecando}, the exponent of $\rho$ is negative for any $\gamma$, since for $\gamma \le 2s+ \alpha-2$ it is $\alpha -2$, and for $\gamma > 2s+\alpha-2$ it is $\gamma -2s$. We can therefore perform the sum in \eqref{eq:wecando} and obtain \eqref{eq:sameasbefore222} also in this case. 

 Repeating around any point in $B_{1/2}$ and thanks to the $L^\infty$ bound for $\MR u_1$ and Lemma~\ref{it:H7}  we get 
\[
\|\MR u_1\|_{C^\alpha(B_{1/2})}\le C\theta.
\]

\item Finally, the bound for $u_2$ 
\[
[\MR u_2]_{C^\alpha(B_{1/2})}\le [\MR u_2]_{C^\gamma (B_{1/2})}\le C\theta  [w]_{C^\gamma(\R^n)}\quad\text{if}\quad \alpha \le \gamma
\]
follows directly from the expression \eqref{eq:firstineqlem} (with $w = u_2$), and summing for $\rho = 2^k$, for $k = -1,0,1,2,\dots$ (since $u_2 \equiv 0$ in $B_{1/2}$). On the other hand, the bound 
\[
[\MR u_2]_{C^\alpha(B_{1/2})} \le CM [w]_{C^\gamma(\R^n)}\quad\text{if}\quad \alpha > \gamma
\]
follows separately for $\L_1$ and $\L_2$ by using \eqref{eq:rem_Lu22} (and a rescaling and covering argument)\qedhere
\end{steps}
\end{proof}

We can now use Lemma~\ref{lem-sdhhh} to show Proposition  \ref{prop-interior-linear-x} (which follows by analogy with the proof of Theorem \ref{thm-interior-linear-x}). 

\begin{proof}[Proof of Proposition  \ref{prop-interior-linear-x}]
The proof is essentially the same as that of Theorem \ref{thm-interior-linear-x}, however \eqref{ngfhhghg} needs to be replaced by
\begin{equation}\label{ngfhhghg2} 
\|\L(u,\cdot) - \L_0u\|_{C^\alpha(B_1)} \leq C\left(\delta\|u\|_{C^{2s+\alpha}(B_2)} + [u]_{C^\gamma(\R^n)}\right).
\end{equation}
The proof of \eqref{ngfhhghg2} follows exactly as the proof of \eqref{ngfhhghg} by replacing the use of Lemmas~\ref{lem:Lu} and \ref{lem:Lu_2}-\ref{it:lem_Lu2_i} by \eqref{eq:Linfty_nondivform} and \eqref{eq:Calpha_nondivform} in Lemma~\ref{lem-sdhhh}, respectively, with $\theta = \delta$ small. The fact that $\theta$ can be taken to be small is for the exact same reason as in the proof of Theorem~\ref{thm-interior-linear-x}.
\end{proof}

\subsection{Schauder estimates for equations in divergence form} 
\index{Divergence-form equations} 

We now consider operators with $x$-dependence in divergence form, of the type 
\begin{equation}\label{divergence-form0}
\begin{split}
\L (u,x) & = {\rm P.V.} \int_{\R^n}\big(u(x)-u(z)\big)K(x,dz) \\
 &  = {\rm P.V.} \int_{\R^n}\big(u(x)-u(x+y)\big)K(x,x+dy) 
\end{split}
\end{equation}
where $(K(x, \cdot))_{x\in \R^n}$ is a family of measures  in $\R^n$ that satisfies the uniform ellipticity conditions 
\begin{equation}\label{divergence-ellipticity0}
r^{2s} \int_{B_{2r}(x)\setminus B_r(x)} K(x,dz)\le \Lambda\qquad\text{for all}\quad x\in \R^n,
\end{equation}
and 
\begin{equation}\label{divergence-ellipticity1}
r^{2s-2}\inf_{e\in \S^{n-1}}\int_{B_{r}(x)}|e\cdot (x-z)|^{2} K(x,dz)\ge \lambda\qquad\text{for all}\quad x\in \R^n,
\end{equation}
as well as symmetry in the two variables, in the sense that
\begin{equation}\label{divergence-form1}
\begin{split}
\int_{A}\int_B  K(x,dz)&\, dx  =\int_B\int_A K(x,dz)\, dx\\
& \quad \textrm{for all}\  A,B\subset \R^n\ \text{Borel, such that} \ A\cap B = \varnothing.
\end{split}
\end{equation}
Observe that, when $(K(x, \cdot))_{x\in \R^n}$ are absolutely continuous with respect to the Lebesgue measure, then \eqref{divergence-form1} reads as
\[
K(x, z) = K(z, x)\qquad \text{for a.e.}\quad  (x, z)\in \R^n\times \R^n.
\]

Equations of the type
\begin{equation}
\label{eq:div_weak}\L(u,x)=f(x)\quad\textrm{in}\quad \Omega
\end{equation}
have a natural weak formulation:
 
\begin{defi} Let $s\in (0, 1)$, and let $\L(\cdot, x)$ be of the form \eqref{divergence-form0}-\eqref{divergence-ellipticity0}-\eqref{divergence-ellipticity1}-\eqref{divergence-form1}. Let $\Omega\subset \R^n$ be any bounded domain, and let $f\in L^p(\Omega)$ for some $p \ge \frac{2n}{n+2s}$ and $n > 2s$. Let $u$ be such that
\[
\iint_{\R^n\times \R^n\setminus(\Omega^c\times \Omega^c)} \left(u(x) - u(z)\right)^2 K(x, dz)\, dx<\infty.
\]
We say that $u$ is a \emph{weak solution} of \eqref{eq:div_weak} if 
\[
\begin{split}\frac12 \int_{\R^n}  \int_{\R^n}   \big(u(x)-u(z)\big)\big(\eta(x)-\eta(z)\big)K( x,dz)& \, dx  =\int_{\R^n}  f\eta 
\end{split}\]
for all $\eta\in C^\infty_c(\Omega)$.

We say that $u$ is a \emph{weak supersolution} of \eqref{eq:div_weak} (resp. \emph{weak subsolution} of \eqref{eq:div_weak}) and we denote it $\L(u, x) \ge f(x)$ in $\Omega$ (resp. $\L(u, x) \le f(x)$ in~$\Omega$) if 
\[
\begin{split}\frac12 \int_{\R^n}  \int_{\R^n}   \big(u(x)-u(z)\big)\big(\eta(x)-\eta(z)\big)K( x,dz)& \, dx  \underset{\left(\text{resp. $\le$}\right)}{\ge} \int_{\R^n}  f\eta 
\end{split}\]
for all $\eta\in C^\infty_c(\Omega)$ with $\eta\ge 0$.
 
\end{defi}

It is important to notice that, in case of divergence-form equations \eqref{divergence-form0}, one cannot symmetrize the operator and write it in terms of  a second-order incremental quotient $2u(x)-u(x+y)-u(x-y)$.
In particular, when $s\geq\frac12$ one cannot evaluate in general $\L(u,x)$ pointwise\footnote{This also happens for operators in divergence form ${\rm div}(A(x) \nabla u)$ in the local case $s = 1$.} even for smooth functions $u\in C^\infty_c(\Omega)$.

In order to obtain Schauder-type estimates, in addition to the uniform ellipticity assumptions \eqref{divergence-ellipticity0}-\eqref{divergence-ellipticity1} we need to assume some $C^\alpha$ regularity of the kernels in the $x$-variable.
More precisely, we assume 
\begin{equation}\label{reg-x-div}
\int_{B_{2\rho}(x)\setminus B_\rho(x)} \big|K(x+h,h+dz)-K(x,dz)\big|\leq M|h|^\alpha \rho^{-2s}
\end{equation}
for all $x, h\in \R^n$,   and $\rho > 0$. Alternatively, such a condition can also be written as (cf. \eqref{x-dependence-L2} in non-divergence-form equations) 
\begin{equation}\label{reg-x-div-xx'}
\int_{B_{2\rho}\setminus B_\rho} \big|K(x,x+dy)-K(x',x'+dy)\big| \leq M|x-x'|^\alpha \rho^{-2s}
\end{equation}
for all $x, x'\in \R^n$ and all $\rho > 0$. 
For divergence form equations, together with the previous regularity in $x$ we also need to assume that the kernel, at small scales, is (quantitatively) \emph{almost even}:  
\begin{equation}\label{reg-x-div-even}
\int_{B_{2\rho}\setminus B_\rho} \big|K(x,x+dy)-K(x,x-dy)\big|  \leq M\rho^{\alpha-2s}
\end{equation}
for all $x\in \R^n$,    and $\rho > 0$. Conditions \eqref{reg-x-div}-\eqref{reg-x-div-even} are like a $C^\alpha$ regularity of the coefficients. Indeed,  when $K$ is absolutely continuous, since $K(x, z) = K(z, x)$  we have that \eqref{reg-x-div}-\eqref{reg-x-div-even}  hold automatically if one assumes the (stronger) pointwise condition
\[
\big|K(x+h, z+h) -K(x, z)\big| \le \frac{M|h|^\alpha}{|z-x|^{n+2s}},\qquad\text{for all}\quad x, z, h\in \R^n. 
\]
In general, though, our assumptions \eqref{reg-x-div}-\eqref{reg-x-div-even} allow for more singular kernels. 

On the other hand, in some cases we will need to assume regularity in the $y$-variable as well, given by 
\begin{equation}\label{reg-x-div-y}
\int_{B_{2\rho}(x) \setminus B_\rho(x) } \big|K(x, h+dz) - K(x, dz)\big|  \le M |h|^\theta\rho^{-2s-\theta},
\end{equation}
(for some $\theta\in (0, 1]$) for all $h\in B_{\rho/2}$ and for all $x\in \R^n$ and $\rho > 0$. Observe that the previous condition is equivalent to asking that, using the notation in \eqref{Calpha-assumption}, $\sup_{x\in \R^n} [K(x, x+\,\cdot\,)]_\theta\le M$.

The interior Schauder estimates for nonlocal divergence form equations are the following:

\begin{thm}\label{thm-interior-linear-x-div}
Let $s\in (0,1)$, $\alpha\in (0, 1]$, $\eps \in (0, \alpha)$, and let $\L$ be an operator of the form \eqref{divergence-form0}-\eqref{divergence-form1}, with kernels satisfying the ellipticity conditions \eqref{divergence-ellipticity0}-\eqref{divergence-ellipticity1}, and \eqref{reg-x-div}-\eqref{reg-x-div-even}
 for some  $M > 0$. 
 
Let $u\in C^{\beta}_{\rm loc}(B_1)\cap L^\infty_{2s-\eps}(\R^n)$ be a weak solution of
\begin{equation}\label{sjdgneorij}
\L(u,x) = f \quad\text{in}\quad B_1.
\end{equation}
with $f\in X$, and  
\begin{equation}\label{sjdgneorij2}
\beta:=\left\{ \begin{array}{rl}
1+\alpha & \textrm{if}\quad s>{\textstyle \frac12}, \vspace{1mm}\\
1+\alpha-\eps &  \textrm{if}\quad s={\textstyle \frac12}, \vspace{1mm} \\
2s+\alpha &  \textrm{if}\quad s<{\textstyle \frac12},
\end{array}\right.\qquad X := \left\{
\begin{array}{ll}
C^{\beta-2s}(B_1)& \quad\text{if}\quad \beta > 2s,\\
L^{\frac{n}{2s-\beta}}(B_1)& \quad\text{if}\quad \beta < 2s.
\end{array}
\right.
\end{equation}
  Assume in addition that $\beta\neq1$, $\beta\neq 2s$, and that \eqref{reg-x-div-y} holds if $\beta > 2s$, with $\theta = \beta-2s$. 
Then,
\[
\|u\|_{C^{\beta}(B_{1/2})} \le C\left(\|u\|_{L^\infty_{2s-\eps}(\R^n)} + \|f\|_{X}\right)
\]
The constant $C$ depends only on $n$, $s$, $\alpha$, $\eps$, $\lambda$, $\Lambda$, and $M$. 
\end{thm}
 
In case $s > \frac12$, we get $C^{1,\alpha}$ regularity, which coincides with the one  obtained in the local case, $s = 1$. For $s < \frac12$, since the equation has ``$C^\alpha$'' coefficients, the maximum regularity we expect to obtain is $C^{2s+\alpha}$, corresponding to the one for the non-divergence-form result.

  The strategy to prove this result follows a similar dichotomy and  will be different in cases $\beta > 2s$ and $\beta < 2s$. When $\beta > 2s$, we will treat $\L$ as a (nonsymmetric) operator in non-divergence form, and argue as in the proof of Theorem~\ref{thm-interior-linear-x}. Instead, when $\beta < 2s$, the proof will be by contradiction and blow-up, similarly to the proof of Theorem~\ref{thm-interior-linear-2}.

\begin{proof}[Proof of Theorem \ref{thm-interior-linear-x-div} in case $\beta>2s$]
When $\beta>2s$ (that is,  $\theta := \beta - 2s > 0$), the operator $\L$ can be evaluated pointwise on smooth functions~$u$, and thus it can be seen as a (nonsymmetric) equation in non-divergence form.
Using this, we can follow the strategy of the proof of Theorem~\ref{thm-interior-linear-x} above. We divide the proof into  six steps. 

\begin{steps}
\item \label{step:1div} 
As in \ref{step:1nondiv} of the proof of Theorem~\ref{thm-interior-linear-x}, we start with an initial reduction wherein, up to considering the rescalings $u_\circ(x) := u(r_\circ x)$ and $\L^{r_\circ}$ (with kernel $K^{r_\circ} (x, dz) = r_\circ^{ 2s}K(r_\circ x, r_\circ \, dz)$, where $K$ is the kernel of $\L$), we can assume that the operator $\L$ satisfies the ellipticity conditions \eqref{divergence-ellipticity0}-\eqref{divergence-ellipticity1}, has regularity in $y$ given by \eqref{reg-x-div-y} for some $M > 0$, but conditions \eqref{reg-x-div} and \eqref{reg-x-div-even} now become
\begin{equation}\label{reg-x-div-2}
\int_{B_{2\rho}(x)\setminus B_\rho(x)} \big|K(x+h,h+dz)-K(x,dz)\big|  \leq Mr_\circ^\alpha |h|^\alpha \rho^{-2s} =: \delta |h|^\alpha \rho^{-2s}
\end{equation}
for all $x, h\in \R^n$,   and $\rho > 0$; and 
\begin{equation}\label{reg-x-div-even-2}
\int_{B_{2\rho}\setminus B_\rho} \big|K(x,x+dy)-K(x,x-dy)\big|  \leq Mr_\circ^\alpha \rho^{\alpha-2s} =: \delta \rho^{\alpha -2s}
\end{equation}
for all $x\in \R^n$,    and $\rho > 0$; for some $\delta > 0$ a small universal constant to be chosen.

\item  Let us denote by $\L^e$ and $\L^o$ respectively the even and odd parts of the operator $\L$. Namely, we have 
\[
 \L^e(u, x) = {\rm P.V.} \int_{\R^n}\big(u(x)-u(x+y)\big)K^e(x,dy),
 \]
 where
 \[
  K^e(x, dy) := \frac{K(x, x+dy) + K(x, x-dy)}{2}
\]
is an even kernel, in the sense that $K^e(x, dy) = K^e(x, -dy)$ (i.e., $K^e(x, dy)$ is a symmetric measure); and 
\[
 \L^o(u, x) = {\rm P.V.} \int_{\R^n}\big(u(x)-u(x+y)\big)K^o(x,dy),
 \]
 where
 \[
  K^o(x, dy) := \frac{K(x, x+dy) - K(x, x-dy)}{2}
\]
is an odd kernel, in the sense that $K^o(x, dy) = -K^o(x, -dy)$. With these definitions, we have that $\L(u, x) = \L^e(u, x) + \L^o(u, x)$. Notice, moreover, that in the case of $\L^e$ we can symmetrize its expression as 
\[
 \L^e(u, x) = \frac12 \int_{\R^n}\big(2u(x)-u(x+y)-u(x-y)\big)K^e(x,dy),
 \]
which is now well-defined in $B_1$ (even without the principal value) because $u\in C^\beta_{\rm loc}(B_1)$ with $\beta > 2s$; cf. Lemma~\ref{lem:Lu}. In fact, the operator $\L^e$ is an operator in non-divergence form like the ones in Theorem~\ref{thm-interior-linear-x}, where the ellipticity conditions are satisfied thanks to \eqref{divergence-ellipticity0}-\eqref{divergence-ellipticity1} and linearity, $\sup_{x\in\R^n} [K^e(x, \cdot)]_\theta \le M$  holds  for $K^e$ thanks to \eqref{reg-x-div-y} and the triangle inequality, and \eqref{x-dependence-L2} holds for $K^e$ with $M = \delta$ by \eqref{reg-x-div-2} and the triangle inequality again. 

We proceed now as in the beginning of \ref{step:step2prove} in the proof of Theorem~\ref{thm-interior-linear-x}. Let $\L^e_0 \in \G_s(\lambda, \Lambda;\theta)$ be the translation invariant operator with kernel $K^e(0, y)$.
By the regularity estimates for translation invariant equations, Theorem~\ref{thm-interior-linear-2}, we have
\begin{equation}\label{ngfhhghg2222-2} 
\|u\|_{C^{\beta}(B_{1/4})} \le C\left(\|u\|_{L^\infty_{2s-\eps}(\R^n)} + \|\L^e_0u\|_{C^{\theta}(B_{1/2})}\right).
\end{equation}
Moreover, 
\begin{equation}\label{ngfhhghg2222-1} \|\L^e_0u\|_{C^\theta(B_{1/2})} \leq \|f\|_{C^\theta(B_{1/2})}+\|\L^o(u,\cdot)\|_{C^\theta(B_{1/2})} + \|\L^e(u,\cdot) - \L^e_0u\|_{C^\theta(B_{1/2})},
\end{equation}
and thanks to \eqref{ngfhhghg}  in \ref{step:step2prove} of the proof of Theorem~\ref{thm-interior-linear-x} (since $\L^e$ is now an operator in non-divergence form) we have 
\begin{equation}\label{ngfhhghg2222} 
\|\L^e(u,\cdot) - \L^e_0u\|_{C^\theta(B_{1/2})} \leq C\left(\delta\|u\|_{C^{\beta}(B_1)} + \|u\|_{L^\infty_{2s-\eps}(\R^n)}\right).
\end{equation}

It only remains to be bounded the $C^\theta$ norm of  $\L^o(u, \cdot)$. 

\item  That is, we now want to prove 
\begin{equation}\label{ngfhhghg3} 
\|\L^o(u,\cdot) \|_{C^{\theta}(B_{1/2})} \leq C\left(\bar \delta\|u\|_{C^{\beta}(B_1)} + \|u\|_{L^\infty_{2s-\eps}(\R^n)}\right),
\end{equation}
for some $\bar \delta > 0$ that is small whenever $\delta$ is small. In fact, we will show 
\begin{equation}\label{ngfhhghg3_2} 
\|\L^o(u,\cdot) \|_{C^{\theta}(B_{1/2})} \leq C\left(\delta_1 \|u\|_{C^{\beta_1}(B_1)} + \|u\|_{L^\infty_{2s-\eps}(\R^n)}\right),
\end{equation}
where
\begin{equation}\label{eq:deltabeta} 
\delta_1 := \delta^{\min\{1-\eps, \frac{2s}{\alpha}\}},\qquad \beta_1 := \min\{1, \beta\}.
\end{equation}

We will use that the operator $\L^o$ has a kernel $K^o$ that satisfies (combining upper ellipticity and \eqref{reg-x-div-even-2})
\begin{equation}\label{eq:Lo0}
\int_{B_{2\rho}\setminus B_\rho} \big|K^o(x,dy)\big| \leq C\rho^{-2s} \min\{1, \delta \rho^\alpha\}\quad\text{for all}\  \rho > 0, \, x\in \R^n, 
\end{equation}
as well as 
\begin{equation}\label{eq:Lo1}
\int_{B_{2\rho} \setminus B_\rho } \big|K^o(x+h,dy)-K^o(x,dy)\big| \leq \delta |h|^\alpha \rho^{-2s}\quad\text{for all}\  \rho > 0, \, x\in \R^n, 
\end{equation}
(by \eqref{reg-x-div-2} and the triangle inequality). 

  We start with the $L^\infty$ bound. For any $x\in B_{1/2}$ we have 
\[
\begin{split}
|\L^o(u, x)|& \le   \int_{\R^n}\left|u(x) - u(x+y)\right|\left|K^o(x, dy)\right| \\
&  \le  \|u\|_{C^{\beta_1}(B_1)} \int_{B_{1/2}}|y|^{\beta_1}\left|K^o(x, dy)\right| \\
&\quad +  C\|u\|_{L^\infty_{2s-\eps}(\R^n)} \int_{B^c_{1/2}}|y|^{2s-\eps}\left|K^o(x, dy)\right|.
\end{split}
\]
We   split each integral in dyadic balls and thanks to \eqref{eq:Lo0} and the fact that $\beta_1 + \alpha - 2s  \ge \theta > 0$ we get
\begin{equation}
\label{eq:Loinf}
|\L^o(u, x)|\le    C\left( \delta \|u\|_{C^{\beta_1}(B_1)}+\|u\|_{L^\infty_{2s-\eps}(\R^n)}\right),
\end{equation}
which gives the $L^\infty$ bound in  \eqref{ngfhhghg3_2} and \eqref{ngfhhghg3}.

The next step is to bound the $C^\theta$ seminorm. To do that, we consider $u=u_1+u_2$ where $u_1:=u\eta$ with $\eta\in C^\infty_c(B_{3/4})$, $0\le \eta \le 1$, and $\eta \equiv 1 $ in $B_{2/3}$, and bound each of the seminorms for $\L^o(u_1, x)$ and $\L^o(u_2, x)$ separately. 

\item  We focus our attention first on finding a bound for the seminorm of $\L^o(u_1, x)$, where we recall that $u_1\in C^\beta_c(B_1)$ with $\|u_1\|_{C^\beta(\R^n)}\le C \|u\|_{C^\beta(B_1)}$. Let us denote, given $\bar x\in B_{1/2}$ fixed, $\L_{\bar x}^o$ to be the translation invariant operator with kernel $K^o(\bar x, y)$ (which is not necessarily positive). We will bound, for any $x_1, x_2\in B_{1/2}$,
\begin{equation}
\label{eq:twotermsdif}
\left|\L^o(u_1, x_1) - \L^o(u_1, x_2)\right|\le |\L_{x_1}^o u_1(x_1)-\L_{x_2}^o u_1(x_1)|+|\L_{x_2}^o u_1(x_1)-\L_{x_2}^o u_1(x_2)|.
\end{equation}
For the first term, we have
\[
\begin{split}
|\L_{x_1}^o u_1(x_1)\hspace{-0.2mm}-\hspace{-0.2mm}\L_{x_2}^o u_1(x_1)|\hspace{-0.5mm}& \le\hspace{-0.5mm}\int_{\R^n} \hspace{-0.5mm}\left|u_1(x_1) - u_1(x_1+y)\right|\left|K^o(x_1, dy) \hspace{-0.5mm}- \hspace{-0.5mm}K^o(x_2, dy)\right|  \\
& = I_1 +I_2,
\end{split}
\]
where, since $u_1(x_1+y) = 0$ for $y\in \R^n\setminus B_2$, by denoting $r := |x_1-x_2|$,
\[
I_1 := |u_1(x_1)|\int_{\R^n\setminus B_2 } |K^o(x_1, dy)-K^o(x_2, dy)|  \le C \delta \|u\|_{L^\infty(B_1)}{r}^\alpha
\]
(using \eqref{eq:Lo1}), and 
\[
\begin{split}
I_2  & := \int_{B_2}\left|u_1(x_1) - u_1(x_1+y)\right|\left|K^o(x_1, dy)  -  K^o(x_2, dy)\right| \\
& \le C\|u\|_{C^{\beta_1}(B_1)} \int_{B_2}|y|^{\beta_1}\left|K^o(x_1, dy)  -  K^o(x_2, dy)\right| \\
& \le   C\|u\|_{C^{\beta_1}(B_1)} \sum_{\substack{\rho = 2^{-k}\\k \ge 0}} \rho^{\beta_1}\int_{B_{2\rho} \setminus B_{\rho}}\left|K^o(x_1, dy)  -  K^o(x_2, dy)\right|,
\end{split}
\]
(recall \eqref{eq:deltabeta}). By \eqref{eq:Lo0}-\eqref{eq:Lo1} we have  
\[
\begin{split}
I_2  & \le  C \delta\|u\|_{C^{\beta_1}(B_1)} \sum_{\substack{\rho = 2^{-k}\\k \ge 0}}  \rho^{\beta_1 - 2s} \min\{\rho^{\alpha}, {r}^\alpha\},\end{split}
\]
and we can split the sum into 
\[
\begin{split}
   {r}^\alpha \sum_{\substack{\rho = 2^{-k}\\1\ge\rho \ge {r}}}  \rho^{\beta_1 - 2s} +\sum_{\substack{\rho = 2^{-k}\\\rho \le {r}}}  \rho^{\beta_1 - 2s+\alpha},
\end{split}
\]
where the second term can be bounded by $C{r}^{\beta_1-2s+\alpha}$, since we have that ${\beta_1-2s+\alpha} > 0$;  and the first term is bounded by 
\[
 {r}^\alpha \sum_{\substack{\rho = 2^{-k}\\1\ge \rho \ge {r}}}  \rho^{\beta_1- 2s}\le \left\{
 \begin{array}{ll}
 C{r}^\alpha & \quad \text{if}\quad s < \frac12,\\
 C{r}^\alpha|\log{r}|& \quad \text{if}\quad s = \frac12,\\
  C{r}^{1+\alpha - 2s} & \quad \text{if}\quad s > \frac12.
 \end{array}
 \right.
\]
Using that ${r}^\alpha|\log{r}|\le C_\eps {r}^{\alpha -\eps}$, in all cases we have 
\[
I_2 \le C \delta \|u\|_{C^{\beta_1}(B_1)}  {r}^{\theta}.
\]
Together with the bound on $I_1$ and the fact that $\alpha \ge \theta$, we obtain 
\begin{equation}
\label{eq:prevtog}
|\L_{x_1}^o u_1(x_1)\hspace{-0.2mm}-\hspace{-0.2mm}\L_{x_2}^o u_1(x_1)| \le C\delta \|u\|_{C^{\beta_1}(B_1)} {r}^\theta. 
\end{equation}

Now, for the second term in \eqref{eq:twotermsdif} we use that, since $\beta_1 \le 1$, 
\[
|u_1(x_1) - u_1(x_1+y) - u_1(x_2) + u_1(x_2+y)|\le C\|u\|_{C^{\beta_1}(B_1)}\min\{{r}^{\beta_1}, |y|^{\beta_1}\},
\]
and thus, by \eqref{eq:Lo0}, 
\[
|\L_{x_2}^o u_1(x_1)-\L_{x_2}^o u_1(x_2)| \le C \|u\|_{C^{\beta_1}(B_1)}\sum_{\rho = 2^k} \rho^{-2s} \min\{\rho^{{\beta_1}} , {r}^{{\beta_1}}\} \min\{1, \delta \rho^\alpha\}.
\]
We split the sum into three terms according to the value of $\rho\in (0, \infty) = (0, r)\cup (r,   \delta^{-\frac{1}{\alpha}}) \cup (\delta^{-\frac{1}{\alpha}}, \infty)$ and bound it by
\begin{equation}
\label{eq:three_terms}
\delta \sum_{\substack{\rho = 2^k\\\rho \le  r}}  \rho^{-2s+\alpha+{\beta_1}} +\delta{r}^{\beta_1} \sum_{\substack{\rho = 2^k\\r \le \rho \le \delta^{-\frac{1}{\alpha}}}} \rho^{-2s+\alpha} +\sum_{\substack{\rho = 2^k\\ \delta^{-\frac{1}{\alpha}}\le\rho }} \rho^{-2s}{r}^{\beta_1}.
\end{equation}
The first term is immediately bounded by $C\delta r^{-2s+\alpha+{\beta_1}}$, since $-2s+\alpha +{\beta_1} > 0$; and the third term is bounded by $C{r}^{\beta_1} \delta^{\frac{2s}{\alpha}}$. Observe that, since $-2s+\alpha+\beta_1\ge\beta - 2s$ and $\beta_1\ge \beta - 2s$, we have that the first and third terms are bounded by $C\delta_1 r^{\beta-2s}$ (recall \eqref{eq:deltabeta}).

 For the second term, we have different values according to the relative values of $\alpha$ and $s$, as follows:
\begin{itemize}
\item If $2s < \alpha$, then the second term in \eqref{eq:three_terms} is bounded by  $C\delta{r}^{\beta_1} \delta^{\frac{2s-\alpha}{\alpha}} = C{r}^{\beta_1} \delta^{\frac{2s}{\alpha}}$. 
\item If $2s > \alpha$, then the second term in \eqref{eq:three_terms} is bounded by  $C\delta r^{-2s+\alpha+{\beta_1}}$. 
\item Finally, in case $2s = \alpha\le 1$, the second term is bounded by the factor $C\delta {r}^{\beta_1}\left(|\log r|+|\log \delta|\right)\le C_{ \eps} \delta^{1- \eps}{r}^{{\beta_1}- \eps}$. We have $\beta_1-\eps \ge \beta-2s$ as well (since $\eps$ was small).
\end{itemize}

Putting all terms together, we have shown that the sum \eqref{eq:three_terms} is bounded by $\delta_1 r^{\beta-2s}$ (recall \eqref{eq:deltabeta}) and thus we have 
\[
|\L_{x_2}^o u_1(x_1)-\L_{x_2}^o u_1(x_2)| \le C \delta_1 \|u\|_{C^{\beta_1}(B_1)}
 r^{\theta}.
\]
With \eqref{eq:prevtog}, this gives in \eqref{eq:twotermsdif}
\begin{equation}
\label{eq:twotermsdif2}
\left|\L^o(u_1, x_1) - \L^o(u_1, x_2)\right|\le  C \delta_1 \|u\|_{C^{\beta_1}(B_1)}|x_1-x_2|^{\theta},
\end{equation}
which bounds the seminorm $[\L^o(u_1, \cdot)]_{C^{\theta}(B_{1/2})}\le C\delta_1 \|u\|_{C^{\beta_1}(B_1)}$.

 \item Let us now bound the $C^\theta$ seminorm of $\L^o(u_2, x)$, where now $u_2$ satisfies that $u_2 \equiv 0$ in $B_{2/3}$ and $\|u_2\|_{L^\infty_{2s-\eps}(\R^n)}\le C\|u\|_{L^\infty_{2s-\eps}(\R^n)}$. To do it, we will use the regularity of the kernel in the $y$ variable, namely, \eqref{reg-x-div-y}. By triangle inequality and in terms of $K^o$, this condition reads as 
 \begin{equation}\label{reg-x-div-y_o}
\int_{B_{2\rho} \setminus B_\rho} \left|K^o(x, h+dy) - K^o(x, dy)\right|  \le M |h|^\theta\rho^{-2s-\theta} \quad\text{for all}\quad h\in B_{\rho/2},
\end{equation}
for all $x\in \R^n$, $\rho > 0$.

 Using the same notation as before, we again split as, for any $x_1, x_2\in B_{1/2}$,
\begin{equation}
\label{eq:twotermsdif_222}
\left|\L^o(u_2, x_1) - \L^o(u_2, x_2)\right|\le |\L_{x_1}^o u_2(x_1)-\L_{x_2}^o u_2(x_1)|+|\L_{x_2}^o u_2(x_1)-\L_{x_2}^o u_2(x_2)|.
\end{equation}
 Now, for the first  term we have, since $u_2(x_1) = u_2(x_2) = 0$, 
\[
\begin{split}
|\L_{x_1}^o u_2(x_1)-\L_{x_2}^o u_2(x_1)|&  \le \int_{\R^n\setminus B_{1/6}} |u_2(x_1+y)| \left|K^o(x_1, dy) - K^o(x_2, dy)\right| \\
& \le C \delta   \|u\|_{L^\infty_{2s-\eps}(\R^n)} |x_1-x_2|^\alpha, 
\end{split}
\]
where we have used the fact that $|u_2(z)|\le C(1+|z|^{2s-\eps}) \|u\|_{L^\infty_{2s-\eps}(\R^n)}$ together with \eqref{eq:Lo1}. 

For the second term, we have 
\[
\begin{split}
|\L_{x_2}^o u_2(x_1)\hspace{-0.3mm}-\hspace{-0.3mm}\L_{x_2}^o u_2(x_2)| & \hspace{-0.3mm}\le\hspace{-0.3mm} \int_{ B^c_{1/2}}\hspace{-1.5mm} |u_2(z)|\left|K^o(x_2, \hspace{-0.15mm}-\hspace{-0.15mm}x_1\hspace{-0.15mm}+\hspace{-0.15mm}dz) \hspace{-0.3mm}-\hspace{-0.3mm} K^o(x_2, \hspace{-0.15mm}-\hspace{-0.15mm}x_2\hspace{-0.15mm}+\hspace{-0.15mm}dz)\right|  \\
& \le C M \|u\|_{L^\infty_{2s-\eps}(\R^n)} |x_1-x_2|^\theta, 
\end{split}
\]
where we have used again the bound on $u_2$ now together with  \eqref{reg-x-div-y_o} (cf. \ref{it:step2Lu2} of the proof of Lemma~\ref{lem:Lu_2}). Since $\alpha \ge \theta$, we have shown that 
\begin{equation}
\label{eq:tihrdto}
[\L^o(u_2, \cdot)]_{C^\theta(B_{1/2})}\le C (M+\delta) \|u\|_{L^\infty_{2s-\eps}(\R^n)}. 
\end{equation}
\item Thanks to \eqref{eq:Loinf}-\eqref{eq:twotermsdif2}-\eqref{eq:tihrdto} we have now shown \eqref{ngfhhghg3_2}, and thus, \eqref{ngfhhghg3}. Together with \eqref{ngfhhghg2222-2}-\eqref{ngfhhghg2222-1}-\eqref{ngfhhghg2222} this shows  
\[
\|u\|_{C^{\beta}(B_{1/4})} \le C\left(\bar \delta\|u\|_{C^{\beta}(B_1)} + \|u\|_{L^\infty_{2s-\eps}(\R^n)} + \|f\|_{C^\theta(B_{1/2})}\right)
\]
for all $u\in C^{\beta}(B_1)\cap L^\infty_{2s-\eps}(\R^n)$, where $\L(u, \cdot) = f$, and where we can still choose $\bar \delta$ small (depending only on $n$, $s$, $\alpha$, $\eps$, $\lambda$, $\Lambda$, and $M$).  Hence, we can use the argument in \ref{step:SAL} of the proof of Theorem~\ref{thm-interior-linear-2} (on page \pageref{step:SAL}), to deduce that (after a covering argument)
\[
\|u\|_{C^{\beta}(B_{1/2})} \le C\left(\|u\|_{L^\infty_{2s-\eps}(\R^n)} + \|f\|_{C^\theta(B_{1})}\right),
\]
for all $u\in C^{\beta}(B_1)\cap L^\infty_{2s-\eps}(\R^n)$, where $\L(u, \cdot) = f$, as wanted.
\qedhere
\end{steps}
\end{proof}

In case $\beta<2s$ (i.e., $s>\frac12$ and $\beta=1+\alpha<2s$) the equation cannot be seen as a non-divergence-form equation, and thus we need a different argument.
We will proceed by a compactness argument, like the one we did for translation invariant equations.   The following proof works only when $\beta < 2s$.  We first show a quantitative Liouville-type estimate for solutions in very large balls (cf. Proposition~\ref{prop:compactness}).
 
\begin{prop}
\label{prop:compactness-x-div}
Let $s\in (\frac12,1)$, $\delta > 0$, and $\alpha \in(0,1)$ be  such that $1+\alpha<2s$.
Let $\L$ be an operator of the form \eqref{divergence-form0}-\eqref{divergence-form1}, with kernels satisfying the ellipticity conditions \eqref{divergence-ellipticity0}-\eqref{divergence-ellipticity1}.
Let $u\in C^{1+\alpha}(\R^n)\cap L^\infty_{2s-\eps}(\R^n)$ for some $\eps>0$ and $f\in L^q_{\rm loc}(\R^n)$ for some $q \ge 1$, and with $[u]_{C^{1+\alpha}(\R^n)}\le 1$ and $\|f\|_{L^q(B_{{1}/{\delta}})} \le \delta$.
Assume in addition that 
\begin{equation}\label{x-dependence-div}
\begin{split} 
& \|\nabla u\|_{L^\infty(\R^n)} \hspace{-1mm}\int_{B_{2\rho}(x)\setminus B_\rho(x)} \hspace{-1mm}\big|K(x\hspace{-0.5mm}+\hspace{-0.5mm}h,h\hspace{-0.5mm}+\hspace{-0.5mm}dz)\hspace{-0.5mm}-\hspace{-0.5mm}K(x,dz)\big|    \leq \delta |h|^\alpha \rho^{-2s}  \\
 & \|\nabla u\|_{L^\infty(\R^n)} \hspace{-1mm}\int_{B_{2\rho}\setminus B_\rho}\hspace{-1mm} \big|K(x,x+dy)\hspace{-0.5mm}-\hspace{-0.5mm}K(x, x-dy)\big|    \leq \delta \rho^{\alpha-2s}
\end{split}
\end{equation}
for all $x,h\in \R^n$, and $\rho > 0$.
Suppose also that $u$ satisfies
\[
\L(u,x) = f\quad\text{in}\quad B_{{1}/{\delta}}
\]
in the weak sense.

Then, for every $\eps_\circ > 0$ there exists $\delta_\circ >0 $ depending only on $\eps_\circ$, $n$, $s$, $\alpha$, $q$, $\eps$, $\lambda$, and $\Lambda$,  such that if $\delta < \delta_\circ$,
\[
\|u-\ell\|_{C^{1}(B_1)}\le \eps_\circ,
\]
where $\ell(x)=u(0)+\nabla u(0)\cdot x$.
\end{prop}

\begin{proof}
The proof is a modification of that of Proposition \ref{prop:compactness}.  For the sake of readability, we will assume here that 
\[
K(x, dz) = K(x, z)\ dz,
\]
so that, in particular, $K(x, z) = K(z, x)$ from \eqref{divergence-form1}.  
We divide it into three steps:
\begin{steps}
\item Assume that the statement does not hold. 
Then, there exists some $\eps_\circ > 0$ such that for any $k\in \N$, there are $u_k\in C_{\rm loc}^{1+\alpha}(\R^n)\cap L^\infty_{2s-\eps}(\R^n)$ with $[u_k]_{C^{1+\alpha}(\R^n)}\le 1$, $f_k\in L^q(B_k)$ with $\|f_k\|_{L^q(B_k)}\le \frac{1}{k}$, and $\L^{(k)}$ as in the statement such that
\[
\L^{(k)} (u_k, x) = f_k \quad\text{in}\quad B_k,
\]
in the weak sense, with
\begin{equation}
\label{eq:byassumptdiv}
\begin{split} 
&\|\nabla u_k\|_{L^\infty(\R^n)}\int_{B_{2\rho}(x)\setminus B_\rho(x)} \big|K^{(k)}(x+\bar h,z+\bar h)-K^{(k)}(x,z)\big|dz \leq \frac{1}{k} |\bar h|^\alpha \rho^{-2s}\\
& \|\nabla u_k\|_{L^\infty(\R^n)}\int_{B_{2\rho}\setminus B_\rho} \big|K^{(k)}(x,x+y)-K^{(k)}(x, x-y)\big|dy \leq \frac{1}{k} \rho^{\alpha-2s},
\end{split}
\end{equation}
for all $x, \bar h\in \R^n$, and $\rho > 0$,  but 
\[
\|u_k - \ell_k\|_{C^1(B_1)}\ge \eps_\circ.
\]

Let us define, for a fixed $h\in \R^n$,
\[
 V_k := u_k(x+h)+u_k(x-h)-2u_k(x), \quad \ F_k := f_k(x+h)+f_k(x-h)-2f_k(x),
\]
with $\|F_k\|_{L^q(B_k)} \le \frac{4}{k}$, $\|V_k\|_{C^{1+\alpha}(\R^n)} \leq C_h$.

\item Let $K_0^{(k)}(y)$ be the kernel denoting the even part of $K(x, x+y)$ at $x = 0$, i.e., $K_0^{(k)}(y) := \frac12K^{(k)}(0,y)+\frac12K^{(k)}(0,-y)$. Observe that by assumption, the operators with kernel $K_0^{(k)}$ belong to $\GL$. 

Now, for any $\eta\in C^\infty_c(\R^n)$ we have (for $k$ such that ${\rm supp}\,\eta\subset B_k$)
\[\int_{\R^n} \int_{\R^n} \hspace{-2mm} \big(V_k(x)-V_k(z)\big)\big(\eta(x)-\eta(z)\big) K^{(k)}_0(z-x)dz\, dx = 2\int_{\R^n}\hspace{-2mm}  F_k\eta+ 2E_0-E_h-E_{-h},\]
where
\[E_0:= \int_{\R^n} \int_{\R^n}  E_0(x, z)\, dz\,dx,\]
with 
\[
E_0(x, z) := \big(u_k(x)-u_k(z)\big)\big(\eta(x)-\eta(z)\big) \big(K^{(k)}(x,z) - K^{(k)}_0(z-x)\big),
\]
and the expressions for $E_{\pm h}$ are analogous, replacing $u_k(x)$, $u_k(z)$, and $K^{(k)}(x,z)$ by $u_k(x\pm h)$, $u_k(z\pm h)$,  and $K^{(k)}(x\pm h,z\pm h)$ respectively. By symmetry in the roles of $x$ and $z$ (here we use that $K(x, z) = K(z, x)$, or \eqref{divergence-form1} in the case of non-absolutely continuous measures),
\begin{equation}
\label{eq:absE0}
\begin{split}
|E_0|& \le \int_{{\rm supp}\, \eta}\int_{\R^n} |E_0
(x, z)|\, dz\, dx+\int_{\R^n}\int_{{\rm supp}\, \eta} |E_0 (x, z)|\, dz\, dx\\
& \le 2\int_{{\rm supp}\, \eta}\int_{\R^n} |E_0
(x, z)|\, dz\, dx.
\end{split}
\end{equation}

Using that $u_k$ are globally Lipschitz and that $\eta\in C^1(\R^n)$ we have
\[
\begin{split} \big|u_k(x)-u_k(z)\big|\,\big|\eta(x)-\eta(z)\big| & \leq  C\|\nabla u\|_{L^\infty(\R^n)}|x-z|\min\{1, |x-z|\},
\end{split}\]
for some constant depending only on $\eta$. Hence, from \eqref{eq:absE0} we have
\begin{equation}
\label{eq:E0bound}
|E_0|\le C \sum_{\substack{\rho = 2^j\\ j \in \Z}}  \rho \min\{1, \rho\}\int_{{\rm supp}\, \eta}I_\rho(x) \, dx.
\end{equation}
where
\[
I_\rho(x) := \|\nabla u\|_{L^\infty(\R^n)}\int_{B_{2\rho}\setminus B_\rho}
\left| K^{(k)}(x, x+y)-K_0^{(k)}(y)\right|\, dy.
\]

We split
\[
\begin{split}
\left|K^{(k)}(x, x+y)-K_0^{(k)}(y)\right|  & \le \frac12 \left|K^{(k)}(x, x+y)-K^{(k)}(0, y)\right| \\
& \quad +\frac12 \left|K^{(k)}(x, x-y)-K^{(k)}(0, -y)\right|\\
& \quad +\frac12 \left|K^{(k)}(x, x+y)-K^{(k)}(x, x-y)\right|,
\end{split}
\] 
so that 
\[
\begin{split}
I_\rho (x)&  \le \|\nabla u\|_{L^\infty(\R^n)}\int_{B_{2\rho}\setminus B_\rho}
\left| K^{(k)}(x, x+y)-K^{(k)}(0, y)\right|\, dy\\
& \quad +\frac{\|\nabla u\|_{L^\infty(\R^n)}}{2} \int_{B_{2\rho}\setminus B_\rho}
\left| K^{(k)}(x, x+y)-K^{(k)}(x, x-y)\right|\, dy.
\end{split}
\]
The first term can be bounded thanks to the first inequality in \eqref{eq:byassumptdiv} (putting $\bar h = x$ and $x = 0$), and the second term is directly bounded thanks to the second inequality in \eqref{eq:byassumptdiv}, so that we obtain 
\[
I_\rho(x) \le \frac{1}{k} \rho^{-2s}\left(|x|^\alpha +   \rho^\alpha\right). 
\]
Putting it back into the bound on $|E_0|$, \eqref{eq:E0bound}, we get
\[
 |E_0|\le \frac{C}{k } \sum_{\substack{\rho = 2^j\\ j \in \Z}}  \rho^{1-2s} \min\{1, \rho\}\int_{{\rm supp}\, \eta}\left(|x|^\alpha +   \rho^\alpha\right)\, dx.
\]
Since $\eta$ is fixed and compactly supported, the last integral is finite and bounded  by $(1+\rho^\alpha)$ (up to a constant depending only on $\eta$), so that 
\[
 |E_0|\le \frac{C}{k } \sum_{\substack{\rho = 2^j\\ j \in \Z}}  \rho^{1-2s} \min\{1, \rho\}\left(1+   \rho^\alpha\right) \le \frac{C}{k},
\]
where the last sum is finite since 
\[
\rho^{1-2s} \min\{1, \rho\}\left(1+   \rho^\alpha\right)\le\left\{
\begin{array}{ll}
2\rho^{1-2s+\alpha}&\quad\text{if} \ \rho \ge 1,\\
2\rho^{2-2s}&\quad\text{if} \ \rho < 1,
\end{array}
\right.
\]
and  $1 +\alpha < 2s < 2$.

\item We have proved $|E_0|\le \frac{C}{k}$, and the same bounds hold for $E_{\pm h}$ as well (for example, simply by considering the test functions $\eta(\cdot\pm h)$ instead of $\eta$). 
Thus, together with the fact that  $\|F_k\|_{L^q(B_k)} \le \frac{4}{k}$, we get 
\[\int_{\R^n} \int_{\R^n} \big(V_k(x)-V_k(z)\big)\big(\eta(x)-\eta(z)\big) K^{(k)}_0(z-x)dz\,dx \longrightarrow 0 \quad \textrm{as}\ k\to\infty.\]
By  Arzel\`a-Ascoli, the functions $V_k$ converge (up to a subsequence) in $C^1_{\rm loc}(\R^n)$ to a function $V\in C^{1+\alpha}(\R^n)$. 
On the other hand, as in the proof of Proposition~\ref{prop:stab_distr}, the measures $\min\{1,|y|^2\}K^{(k)}_0(dy)$ converge weakly to a limiting measure $\min\{1,|y|^2\}\bar K_0(dy)$ that will satisfy 
\[\int_{\R^n} \int_{\R^n} \big(V(x)-V(z)\big)\big(\eta(x)-\eta(z)\big) \bar K_0(z-x)dz\, dx= 0.\]
Notice  that since $V\in C^{1+\alpha}(\R^n)$, it has finite energy \eqref{eq:u_weak_sol} on compact sets. Together with the fact that the previous equality holds  for any $\eta\in C^\infty_c(\R^n)$, we have that $V$ solves $\bar \L_0 V=0$ in $\R^n$ in the weak sense (where $\bar \L_0\in \GL$ is the limiting operator with kernel $\bar K_0$), and by Liouville's theorem, Theorem~\ref{thm:Liouville} (together with Lemma~\ref{lem:weak_distr}), we get that $V$ is constant.
As in the proof of Proposition~\ref{prop:compactness}, if we define 
\[
v_k := u_k - \ell_k,
\]
then, $v_k(0) = |\nabla v_k(0)|  = 0$ with $[v_k]_{C^{1+\alpha}(\R^n)}\le 1$ and $v_k\to v$ in $C^1_{\rm loc}$, for some $v$ with $v(0) = |\nabla v(0)| = 0$ and $[v]_{C^{1+\alpha}(\R^n)}\le 1$. Since  $V_k(x) = v_k(x+h) + v_k(x-h) - 2 v_k(x)$, we get $V(x)=v(x+h)+v(x-h)-2v(x)$, which is constant (for every $h\in \R^n$ fixed). By Lemma~\ref{it:H10} we have that $v$ is a quadratic polynomial, and the condition $[v]_{C^{1+\alpha}(\R^n)}\le 1$ implies it is actually linear. Because it also satisfies $v(0) = |\nabla v(0)| = 0$, it must be $v \equiv 0$, which is a contradiction with $\|v\|_{C^1(B_1)}\ge \eps_\circ>0$. 
\qedhere
\end{steps} 
\end{proof}

Thanks to the previous Liouville-type statement, we get the following estimate, which is almost the desired result in case $\beta < 2s$:

\begin{prop}
\label{prop:interior-linear-x-div}
Let $s\in (\frac12,1)$ and $\alpha\in (0,1)$ be such that $1+\alpha<2s$, and let $q = \frac{n}{2s-1-\alpha}$. Let $\L$ be an operator of the form \eqref{divergence-form0}-\eqref{divergence-form1}, with kernels satisfying \eqref{divergence-ellipticity0}-\eqref{divergence-ellipticity1} and \eqref{reg-x-div}-\eqref{reg-x-div-even} for some $M > 0$.
Then, the following holds. 

For any $\delta > 0$ there exists $C_\delta$ such that
\[
[u]_{C^{1+\alpha}(B_{1/2})} \le \delta [u]_{C^{1+\alpha}(\R^n)} + C_\delta \left(\|u\|_{L^\infty(B_1)} +\|\nabla u\|_{L^\infty(\R^n)} + \|f\|_{L^q(B_1)}\right)
\]
for any $u \in C^{1+\alpha}_c(\R^n)$ satisfying $\L(u,x)=f$ in $B_1$ in the weak sense.
The constant $C_\delta$ depends only on $\delta$, $n$, $s$, $\alpha$,    $M$, $\lambda$, and $\Lambda$.
\end{prop}

\begin{proof}
Let us denote,  for $w\in C_c^{1+\alpha}(\R^n)$,
\[
\tilde{\mathcal{S}} (w) := \inf\left\{\|g\|_{L^q(B_1)} :  
\begin{array}{l}
\tilde \L(w,x)=g \ \text{in the weak sense, for some $\tilde\L$ of}\\
\text{the form \eqref{divergence-form0}-\eqref{divergence-form1}-\eqref{divergence-ellipticity0}-\eqref{divergence-ellipticity1}}\\\text{and satisfying \eqref{reg-x-div}-\eqref{reg-x-div-even}, with $M>0$.}
\end{array}\right\}.
\]

We use Lemma~\ref{lem-interior-blowup} with $\mu = 1+\alpha$ and
\[
\mathcal{S}(w) =  \left\{
\begin{array}{ll} \displaystyle
\tilde{\mathcal{S}} (w) + \|\nabla w\|_{L^\infty(\R^n)} & \text{if}\quad w\in  C_c^{1+\alpha}(\R^n).
\\
\infty & \text{otherwise,}
\end{array}
\right.
\]
Notice that the mapping $\mathcal{S}:C^{1+\alpha}(\R^n)\to \R_{\ge 0}$ depends only on $n$, $s$, $\alpha$,    $M$, $\lambda$, and~$\Lambda$.

Thus, either Lemma~\ref{lem-interior-blowup}~\ref{it:lem_int_blowup_i} holds, in which case we would have
\[
[u]_{C^{1+\alpha}(B_{1/2})} \le \delta [u]_{C^{1+\alpha}(\R^n)} + C_\delta \left(\|u\|_{L^\infty(B_1)} +\|\nabla u\|_{L^\infty(\R^n)} + \|f\|_{L^q(B_1)}\right),
\]
or there exists a sequence $u_k\in C^{1+\alpha}_c(\R^n)$ and $\L_{k}$ of the previous form such that $\L_k(u_k,x)=f_k$ in the weak sense,
\begin{equation}
\label{eq:fktozero-x-div}
\frac{\|f_k\|_{L^q(B_1)}+\|\nabla u_k\|_{L^\infty(\R^n)}}{[u_k]_{C^{1+\alpha}(B_{1/2})}} \le \frac{2\mathcal{S}(u_k)}{[u_k]_{C^{1+\alpha}(B_{1/2})}}\to 0,
\end{equation}
and for some $x_k \in B_{1/2}$ and $r_k \downarrow 0$,
\[
v_k(x) := \frac{u_k(x_k+r_k x)}{r_k^{1+\alpha}[u_k]_{C^{1+\alpha}(\R^n)}}
\]
satisfies 
\begin{equation}
\label{eq:vkcontradiction-x-div}
\|v_k - \ell_k\|_{C^1(B_1)} > \frac{\delta}{2},
\end{equation}
where $\ell_k$ is the 1st order Taylor polynomial of $v_k$ at 0.  
Then, by scaling, there exists an operator of the form \eqref{divergence-form0}-\eqref{divergence-form1}-\eqref{divergence-ellipticity0}-\eqref{divergence-ellipticity1}, $\tilde \L_k$, such that 
\[
\tilde \L_k(v_k,x) = r_k^{2s-1-\alpha}\frac{f_k(x_k+r_k x)}{ [u_k]_{C^{1+\alpha}(\R^n)}} =:\tilde f_k(x)
\]
in the weak sense.
More precisely, if $\L_k$ has kernel $K^k(x,z)$, then $\tilde \L_k$ has kernel $\tilde K^k(x,dz)$ given by   
\[
\tilde K^k (x,dz) := r_k^{2s} K^k (x_k+r_kx,x_k+r_k \, dz).
\]
Moreover, since $\L_k$ satisfies \eqref{reg-x-div}-\eqref{reg-x-div-even}, $\tilde\L_k$ also satisfies them with a smaller $M$, i.e., as in \ref{step:1div} of the proof of Theorem~\ref{thm-interior-linear-x-div} in case $\beta>2s$ (on page \pageref{step:1div}) we have 
\[
\int_{B_{2\rho}(x)\setminus B_\rho(x)} \big|\tilde K^k(x+h,h+dz)-\tilde K^k(x,dz)\big| \leq   Mr_k^{\alpha} |h|^\alpha \rho^{-2s},
\]
and 
\[
\int_{B_{2\rho}\setminus B_\rho} \big|\tilde K^k(x,x+dy)-\tilde K^k(x,x-dy)\big|  \leq Mr_k^\alpha \rho^{\alpha-2s}
\]
for all $x, h\in \R^n$,    and $\rho > 0$.
We also have
\[   \|\nabla v_k\|_{L^\infty(\R^n)} \leq \frac{r_k^{-\alpha}\|\nabla u_k\|_{L^\infty(\R^n)}}{[u_k]_{C^{1+\alpha}(\R^n)}},\]
so that 
\[
\begin{split}
&  \|\nabla v_k\|_{L^\infty(\R^n)} \hspace{-1mm}\int_{B_{2\rho}(x)\setminus B_\rho(x)}\hspace{-1mm} \big|\tilde K^k(x+h,h+dz)-\tilde K^k(x,dz)\big| \leq  \\
& \hspace{7.5cm}\le  \frac{M \|\nabla u_k\|_{L^\infty(\R^n)}}{[u_k]_{C^{1+\alpha}(\R^n)}} |h|^\alpha \rho^{-2s},
\\
&  \|\nabla v_k\|_{L^\infty(\R^n)}\hspace{-1mm} \int_{B_{2\rho}\setminus B_\rho}\hspace{-1mm} \big|\tilde K^k(x,x+dy)-\tilde K^k(x,x-dy)\big|  \leq \frac{M \|\nabla u_k\|_{L^\infty(\R^n)}}{[u_k]_{C^{1+\alpha}(\R^n)} \rho^{2s-\alpha}}
\end{split}
\]
for all $x, h\in \R^n$, $\rho > 0$, and with $\|\nabla u_k\|_{L^\infty(\R^n)} / [u_k]_{C^{1+\alpha}(\R^n)} \to 0$ as $k\to \infty$ (thanks to \eqref{eq:fktozero-x-div}).

Also from \eqref{eq:fktozero-x-div} and since $q = \frac{n}{2s-1-\alpha}$,
\[
\|\tilde f_k\|_{L^q (B_{1 /(2r_k)})} 
= \frac{\|f_k(x_k+r_k\,\cdot\,)\|_{L^q(B_{1 /(2r_k)})}}{r_k^{1+\alpha-2s} [u_k]_{C^{1+\alpha}(\R^n)}}
\leq  r_k^{2s-1-\alpha-\frac{n}{q}}\frac{\|f_k\|_{L^q(B_1)}}{[u_k]_{C^{1+\alpha}(\R^n)}}
\to 0, 
\] 
as $k\to \infty$. 
In all, since by definition $[v_k]_{C^{1+\alpha}(\R^n)} = 1$,  we have that $v_k$ satisfies all the hypotheses of Proposition~\ref{prop:compactness-x-div} for any fixed $\delta_\circ>0$, if $k$ is large enough. 
In particular, taking $\eps_\circ$ sufficiently small in Proposition~\ref{prop:compactness-x-div} we get a contradiction with \eqref{eq:vkcontradiction-x-div}.
\end{proof}

We can now give the final part of the proof of Theorem \ref{thm-interior-linear-x-div}:

\begin{proof}[Proof of Theorem \ref{thm-interior-linear-x-div} in case $\beta<2s$]
Notice that, since $\beta<2s$, then $\beta=1+\alpha<2s$ and $s\in(\frac12,1)$.
By Proposition~\ref{prop:interior-linear-x-div}, for any $\delta> 0$ there exists $C_\delta$ depending only on $\delta$, $n$, $s$, $\alpha$,    $M$, $\lambda$, and $\Lambda$,    such that 
\begin{equation}
\label{eq:touse-x-div}
[u]_{C^{\beta}(B_{1/2})}\le \delta [u]_{C^{\beta}(\R^n)} + C_\delta \left(\|u\|_{L^\infty(B_1)}+\|\nabla u\|_{L^\infty(\R^n)}+\|f\|_{L^q(B_1)}\right)
\end{equation}
for any $u\in C^{\beta}_c(\R^n)$ satisfying $\L(u,x)=f$ in the weak sense, where $ q := \frac{n}{2s-\beta} = \frac{n}{2s-1-\alpha}$.  

Let $\eta\in C^\infty_c(B_3)$ such that $\eta \equiv 1 $ in $B_{2}$, and consider the function $u\eta$ for $u\in C^{\beta}(\R^n)\cap L^\infty_{2s-\eps}(\R^n)$ satisfying $\L(u,x)=f$ in the weak sense.
Since $u - \eta u \equiv 0$ in $B_{2}$, we have that
\begin{equation}\label{claim-x-div}
\|\L (u-\eta u,x)\|_{L^\infty(B_1)} \le C \|u\|_{L^\infty_{2s-\eps}(\R^n)}.
\end{equation}

Hence, we have $\L (\eta u,x)=g$ in the weak sense, with 
\[
\|g\|_{L^q(B_1)} \le \|f\|_{L^q(B_1)} + C \|u\|_{L^\infty_{2s-\eps}(\R^n)}.
\]
Apply now \eqref{eq:touse-x-div} to $u\eta$, to get
\begin{equation}
\label{eq:touse2-x}
[u]_{C^{\beta}(B_{1/2})}\le   \delta [u]_{C^{\beta}(B_4)} + C_\delta \left(\|u\|_{L^\infty_{2s-\eps}(\R^n)}+\|\nabla u\|_{L^\infty(B_4)}+\|f\|_{L^q(B_1)}\right)
\end{equation}
for any $u\in C^{\beta}(B_4)\cap L^\infty_{2s-\eps}(\R^n)$. By interpolation (Proposition~\ref{it:H9}) we know that 
\[
\|\nabla u\|_{L^\infty(B_4)}\le \delta[u]_{C^\beta(B_4)} + C_\delta \|u\|_{L^\infty(B_4)}, 
\]
so that \eqref{eq:touse2-x} becomes
\begin{equation}
\label{eq:touse2-x2}
[u]_{C^{\beta}(B_{1/2})}\le   2\delta [u]_{C^{\beta}(B_4)} + C_\delta \left(\|u\|_{L^\infty_{2s-\eps}(\R^n)}+\|f\|_{L^q(B_1)}\right)
\end{equation}
for any $u\in C^{\beta}(B_4)\cap L^\infty_{2s-\eps}(\R^n)$.

Now, by the exact same interpolation argument as in \ref{step:SAL}  of the proof of Theorem \ref{thm-interior-linear-2} on page \pageref{step:SAL}, we deduce from \eqref{eq:touse2-x2} that
\[
\|u\|_{C^{\beta}(B_{1/2})} \le C \left(\|u\|_{L^\infty_{2s-\eps}(\R^n)}+ \|f\|_{L^q(B_1)}\right)
\]
for all $u\in C^{\beta}(B_1)\cap L^\infty_{2s-\eps}(\R^n)$ such that $\L(u, x) = f$ in $B_1$ in the weak sense, which proves Theorem \ref{thm-interior-linear-x-div}.
\end{proof}

\section{H\"older regularity up to the boundary}
\label{sec:bdryregularity}

As we have seen in Section~\ref{sec:int_reg_G}, solutions to $\L u = f$ in $B_1$ with  bounded right-hand side are $C^{2s}$ inside $B_1$, and are $C^\infty$ in $B_1$ whenever $f$ is $C^\infty$ and $\L\in \G_s(\lambda, \Lambda; \mu)$ for all $\mu > 0$. 
We now study the regularity \emph{up to the boundary}. Consider for example the problem 
\[
\left\{
\begin{array}{rcll}
\L u & = & f& \quad\text{in}\quad B_1\\
u & = & 0& \quad\text{in}\quad \R^n \setminus B_1. 
\end{array}
\right.
\]
What can we say about the regularity of $u$ in $\overline{B_1}$?

A simple example for $\L = \fls$ showed us that, already in dimension~1, the function 
\[
v(x) = \big(1-x^2\big)_+^s
\] 
satisfies 
\[
\fls v = q_{s}\quad\text{in}\quad B_1,
\]
for some constant $q_{s}> 0$ (recall Proposition~\ref{prop:sqrtlv11D} for $\sqrtl$ and Proposition~\ref{prop:count_ex_s} for $\fls$). 
In particular, we do not expect better regularity than $C^s$ up to the boundary, even if we have $C^\infty$ interior regularity. 
We emphasize that this is a nonlocal phenomenon, since in the local analogue typically the interior regularity determines the regularity up to the boundary (in smooth domains with smooth boundary datum).

In view of the previous example, one may conjecture that, at least in smooth domains $\Omega\subset \R^n$, one has $u\in C^s(\overline{\Omega})$. This is precisely what we prove in this section, provided $\L\in \GL$  is such that its L\'evy measure $K$ is \emph{homogeneous}, \eqref{eq:Lu_stab0} (that is, $\L\in \GLh$, recall Definition~\ref{defi:Gh}): 
\[
\left.
\begin{array}{c}
\text{$\L\in \GLh$}\\[0.25cm]
\left\{
\begin{array}{rcll}
\L u & = & f& \quad\text{in}\quad \Omega\\
u & = & 0& \quad\text{in}\quad \R^n \setminus \Omega
\end{array}
\right.
\\[0.45cm]
\text{$\Omega$ is a $C^{1,\alpha}$ domain}
\end{array}
\right\}
\quad\Longrightarrow\quad
\|u \|_{C^s(\overline{\Omega})}\le C \|f\|_{L^\infty(\Omega)}. 
\]

In other words, we will need to assume that  $\L$ is a  \emph{stable operator}, see \eqref{eq:Khom}-\eqref{eq:Khom2}; or equivalently, $\L$ is of the form \eqref{eq:Lu_stab}.
This result was first proved in \cite{RS-Cs} for the fractional Laplacian in $C^{1,1}$ domains, and later on in \cite{RS-C1} for general stable operators in $C^{1,\alpha}$ domains; see also \cite{RS-stable} and the results of Grubb \cite{Grubb,Grubb2,Gru22}.

In order to prove this result, the strategy is to first show that 
\[
|u(x)| \le  Cd^s(x),
\]
and then combine it with the interior estimates to get that $u\in C^s(\overline{\Omega})$.   
Here, and throughout the book, we denote
\begin{equation}
\label{eq:ddef}
d_\Omega(x) := {\dist}(x, \Omega^c)\qquad\text{for any}\quad x\in \R^n. 
\end{equation}
 When there is no possible confusion about the domain $\Omega$, we will simply denote $d := d_\Omega$. 

The precise result we prove is the following: 
\begin{thm}[Global regularity in $C^{1,\alpha}$ domains]
\label{thm:globCsreg}\index{Boundary regularity!General operators}
Let $s \in (0, 1)$ and let $\L \in \GLh$. Let $\alpha \in (0, 1)$, and let $\Omega$ be any bounded $C^{1,\alpha}$ domain. Let $f\in L^\infty(\Omega)$, and let $u$ be the weak solution of 
\[
\left\{
\begin{array}{rcll}
\L u &=& f& \quad\text{in}\quad \Omega\\
u &=& 0& \quad\text{in}\quad \R^n \setminus \Omega.
\end{array}
\right.
\]
Then $u \in C^s(\overline{\Omega})$ with 
\[
\|u \|_{C^s(\Omega)}\le C \|f\|_{L^\infty(\Omega)}
\]
for some $C$ depending only on $n$, $s$, $\Omega$, $\lambda$, and $\Lambda$. 
\end{thm}

It turns out that, if we want a fine description of the boundary regularity of solutions, we need the kernels $K$ to be homogeneous (namely, $\L$ to be  a stable operator, \eqref{eq:Lu_stab}). This property, that was not important for the interior regularity, becomes essential for the proof that we present here of the regularity up to the boundary\footnote{See \cite{RW24}   for the case of non-homogeneous kernels.}. The main difference comes from the following lemma, stating that the one dimensional functions $u(x) = (x\cdot \be)^s_+$ satisfy $\L u = 0$  in $\{x\cdot \be > 0\}$ when $\L$ is homogeneous (cf. Proposition~\ref{prop:onedbarrier} for $\fls$):
\begin{lem}
\label{lem:L1d}
Let $s\in (0, 1)$, and let $\L\in \GLh$. Let 
\[
u(x) = (x\cdot\be)^s_+
\]
for some $\be\in \S^{n-1}$. Then 
\[
\L u = 0\quad\text{in}\quad \{x\cdot\be > 0\}. 
\]
\end{lem}

\begin{proof}
We denote $\bar u (t) = t^s_+$ for $\bar u:\R\to\R$, so that $u(x) = \bar u(x\cdot \be)$. 
We compute $\L u (x)$ by using polar coordinates $y=r\theta$, with $\theta\in \mathbb S^{n-1}$ and $r>0$, so that (as an abuse of notation, since $K$ is homogeneous)
\[K(dy)=\frac{dr}{r^{1+2s}}\,K(d\theta)\]
(c.f. \eqref{eq:Lu_stab}). 
Then, we have
\[
\begin{split}
\L& u (x)  = \frac{1}{2} \int_{\S^{n-1}}\int_0^{\infty} \big(2u(x) - u(x+r\theta) - u(x-r\theta)\big) \frac{dr}{r^{1+2s}}K(d \theta) \\
& = \frac{1}{4} \int_{\S^{n-1}}\int_{-\infty}^{\infty} \big(2u(x) - u(x+r\theta) - u(x-r\theta)\big) \frac{dr}{|r|^{1+2s}}K(d\theta)  \\
& = \frac{1}{4} \int_{\S^{n-1}}\int_{-\infty}^{\infty} \big(2\bar u(x\cdot \be) - \bar u((x+r\theta)\cdot\be ) - \bar u((x-r\theta)\cdot \be )\big) \frac{dr}{|r|^{1+2s}}K(d \theta) \\
& = c_s\int_{\S^{n-1}} \big({(-\Delta)^s_\R}\bar u\big)(x\cdot \be+r \theta\cdot \be)\big|_{r = 0}K(d \theta)  = 0
\end{split}
\]
for $x\in \R^n$ such that $x\cdot\be > 0$, by Proposition~\ref{prop:onedbarrier}.
\end{proof}

It is important to emphasize that the previous lemma is \emph{not} true for general operators $\L \in \GL$ without the homogeneity assumption.

\subsection{The case of convex domains}
\index{Boundary regularity!General operators!Convex domains}
As a consequence of Lemma~\ref{lem:L1d}, we can use the function $(x\cdot\be)_+^s$ as a barrier from above to show that, if $u$ satisfies 
\[
\left\{
\begin{array}{rcll}
\L u & = & f& \quad\text{in}\quad \Omega\\
u & = & 0& \quad\text{in}\quad \R^n \setminus \Omega 
\end{array}
\right.
\]
  for some convex and bounded domain $\Omega\subset \R^n$, where $\L\in \GL$ is such that $K$ is homogeneous, \eqref{eq:Lu_stab0}-\eqref{eq:Lu_stab} (i.e., $\L\in \GLh$), and $f$ with $\|f\|_{L^\infty(\Omega)}\le 1$, then 
\begin{equation}
\label{eq:convex_holds}
|u(x)|\le C d^s(x)\quad\text{for}\quad x\in \R^n,
\end{equation}
for some $C$ depending only on $n$, $s$, $\Lambda$, $\lambda$, and ${\rm diam}(\Omega)$.

\begin{figure}
\centering
\makebox[\textwidth][c]{\includegraphics[scale = 1]{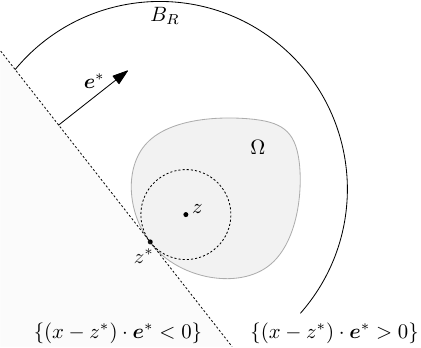}}
\caption{\label{fig:04} Setting in which we want to use $\phi$ as a barrier from above for $u$.}
\end{figure}

Indeed, for each $z\in \Omega$, let $z^*\in \partial \Omega$ such that $|z-z^*| = d(z)$, and let $\be^* =\frac{z-z^*}{|z-z^*|}\in \mathbb{S}^{n-1}$. Since $\Omega$ is convex, we know that 
$
\Omega \subset \{x\in \R^n : (x-z^*)\cdot\be^* > 0\}
$ (see Figure~\ref{fig:04}).
We define a translation of the 1-dimensional barrier from Lemma~\ref{lem:L1d},
\[
\phi(x) = \big((x-z^*)\cdot \be^*\big)^s_+.
\]
Observe that $\phi > 0$ in $\Omega$ and, by Lemma~\ref{lem:L1d}, $\L \phi = 0$ in $\{(x-z^*)\cdot \be^* > 0\}$. On the other hand, if we assume $0\in \Omega$ (after a translation), for any $x\in \Omega$, and for $R = {\rm diam}(\Omega)$,
\[
\begin{split}
-\L \big(\phi\chi_{B_R^c} \big) (x) & = \int_{\R^n}  \big((x+y-z^*\big)\cdot \be^*)^s_+\chi_{B_R^c} (x+y)\, K(dy)  \\
& \ge   \int_{2R}^\infty \int_{\mathbb{S}^{n-1}} \big(\theta\cdot \be^*\big)^s_+ K(d \theta)   \frac{dr}{r^{1+s}} \ge c_\circ > 0,
\end{split}
\]
see Remark~\ref{rem:equiv_ellipt_homog}. 
Since $\L \phi = 0$ for $x\in \Omega$ this implies $(\L \phi\chi_{B_R})(x) \ge c_\circ>0$ in $\Omega$. By defining
\[
v :=  \frac{1}{c_\circ} \phi \chi_{B_R},
\]
we have that $v \ge 0 = u$ in $\R^n \setminus \Omega$, and 
\[
\L u = f \le   1\le \L v\quad\text{in}\quad \Omega. 
\]
By the comparison principle, Corollary~\ref{cor:comp_principle_G_w}, we get that $ v \ge u $ in $\R^n$. In particular, $u(z) \le C\phi(z) = Cd^s(z)$. Repeating for every $z\in \Omega$ and replacing $u$ by $-u$, we get that \eqref{eq:convex_holds} holds for all $x\in \R^n$.

More generally, for $C^{1,\alpha}$ domains $\Omega$, we need a more appropriate barrier, adapted to the domain. Such barrier  is constructed in Appendix~\ref{app.B}; see Corollary~\ref{cor:supersol_domains}.

\subsection{Proof of the boundary regularity}
\label{ssec:boundary_reg}
Before proving Theorem~\ref{thm:globCsreg}, let us state and prove a local version of the result.

We will use the following: 
 
\begin{defi}[Regular domains]
\label{defi:varrho}
Given a $C^{k,\alpha}$ domain with $k \in \N$, and $\alpha\in (0, 1]$, we say that $\varrho$ is a $C^{k,\alpha}$-radius for $\Omega$ if, for any $x_\circ\in \partial\Omega$ with $\nu_\circ \in \S^{n-1}$ the unit normal to $\partial\Omega$, we have
\[
\varrho^{-1} \mathcal{R} \left(\partial\Omega-x_\circ\right)\cap (B'_1\times [-1,1]) = \left\{ 
\begin{array}{l}
(x', x_n)\in B'_1\times [-1,1] : x_n = \varphi(x'),\\
\text{with}\quad \|\nabla \varphi\|_{C^{k-1,\alpha}(B_1')}\le 1
\end{array}\right\},
\]
for some $\varphi\in C^{k,\alpha}(B_1')$, where $\mathcal{R}$ is any rotation such that $\mathcal{R}(\nu_\circ) = e_n$, and $B_1'\subset\R^{n-1}$ is the unit ball.
\end{defi}

Observe that if $\varrho>0$ is a $C^{k,\alpha}$-radius for $\Omega$, then $0<\varrho' <\varrho$ is also a $C^{k,\alpha}$-radius for $\Omega$.

The local version of Theorem~\ref{thm:globCsreg} is then the following:

\begin{prop}[Boundary regularity in $C^{1,\alpha}$ domains]
\label{prop:globCsreg_loc}
Let $s \in (0, 1)$ and let $\L \in \GLh$. Let $\alpha \in (0, 1)$, and let $\Omega$ be any  $C^{1,\alpha}$ domain  with  $C^{1,\alpha}$-radius $\varrho_\circ > 0$.  Let $f\in L^\infty(\Omega\cap B_1)$, and let  $u\in L^\infty_{2s-\eps}(\R^n)$ for some $\eps > 0$ be a weak solution of
\[
\left\{
\begin{array}{rcll}
\L u &=& f& \quad\text{in}\quad \Omega \cap B_1\\
u &=& 0& \quad\text{in}\quad B_1 \setminus \Omega.
\end{array}
\right.
\]
Then, $u \in C_{\rm loc}^s(B_{1})$ with
\[
\|u \|_{C^s(B_{1/2})}\le C\left( \|u\|_{L^\infty_{2s-\eps}(\R^n)}+ \|f\|_{L^\infty(\Omega\cap B_1)}\right)
\]
for some $C$ depending only on $n$, $s$, $\varrho_\circ$, $\eps$, $\lambda$, and $\Lambda$. 
\end{prop}

\begin{proof}
Let $\tilde u = u\chi_{B_2}$, and let $\L \tilde u = \tilde f$ in $\Omega\cap B_1$, with 
\[
\begin{split}
\|\tilde f\|_{L^\infty(\Omega\cap B_1)}& \le \| f\|_{L^\infty(\Omega\cap B_1)} + \|\L (u\chi_{B^c_2}) \|_{L^\infty(B_1)}\\
& \le \| f\|_{L^\infty(\Omega\cap B_1)} + C\|u\|_{L^\infty_{2s-\eps}(\R^n)}
\end{split}
\]
by Lemma~\ref{lem:Lu}. After dividing by a constant, we assume that $\|\tilde f \|_{L^\infty(\Omega\cap B_1)}\le 1$ and $\|\tilde u\|_{L^\infty(\R^n)}\le 1$, and we notice that $u = \tilde u$ in $B_1$, and $\tilde u \equiv 0$ in $\R^n\setminus B_2$.  

We divide the proof into two steps.
\begin{steps}
\begin{figure}
\centering
\makebox[\textwidth][c]{\includegraphics[scale = 1.0]{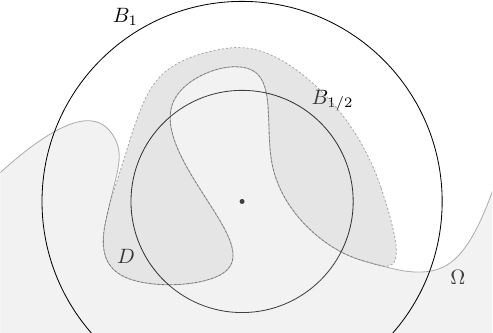}}
\caption{\label{fig:05} A possible example of the set $D$ from \ref{it:step1bdry}.}
\end{figure}
\item \label{it:step1bdry} Let  $D\subset B_1$ with
\[
B_{1/2}\setminus \Omega\subset D\subset B_{1}\setminus \Omega
\]
 be a $C^{1,\alpha}$ set with $C^{1,\alpha}$-radius $c \varrho_\circ$  (for some $c > 0$ universal); see Figure~\ref{fig:05}. Thanks to Corollary~\ref{cor:supersol_domains} applied to $\tilde \Omega := B_3\setminus D\supset \Omega\cap B_2$, there exists $\delta > 0$ (depending only on $n$, $s$, $\alpha$, $\varrho_\circ$, $\lambda$, and $\Lambda$) and $\varphi$ such that denoting $N_\delta :=  \{0 < d_{\tilde\Omega}(x) < \delta\}$,
 \[
\left\{
\begin{array}{rcll}
\L \varphi & \ge & 1& \quad\text{in}\quad N_\delta\\
  \varphi(x)& \le & \frac{1}{\delta}d^s_{\tilde\Omega}(x)   & \quad\text{for}\quad  x\in \tilde\Omega.
\end{array}
\right.
\]  
(Since $\varphi\in H^s(\R^n)$, this is also satisfied in the weak sense.) 
In particular, and since $u = 0$ in $\R^n\setminus B_2$,
 \[
 C \L \varphi \ge \L \tilde u \quad\text{in}\quad N_\delta\cap B_1\cap \Omega\qquad\text{and}\qquad C  \varphi \ge \tilde u\quad\text{in} \quad \R^n\setminus (N_\delta\cap \Omega)
 \]
 for some $C$ large enough (depending on $\delta$). 
By the comparison principle, Corollary~\ref{cor:comp_principle_G_w} , we have that $C \varphi \ge \tilde u$ in $N_\delta\cap B_1\cap \Omega$ as well. Consequently
\begin{equation}
\label{eq:boundbdry}
|\tilde u(x)| = |u(x)| \le Cd_{\tilde \Omega}(x)^s =  Cd_\Omega(x)^s\quad\text{for}\quad x\in B_{3/4}. 
\end{equation}
\item \label{it:step2bdry} We now combine \eqref{eq:boundbdry} with the interior estimates in Theorem~\ref{thm-interior-linear-Lp} to get the desired result (recall that weak solutions in $L^\infty_{2s-\eps}(\R^n)$ are distributional solutions as well, see Lemma~\ref{lem:weak_distr}).

Indeed, let $x_1, x_2\in \Omega\cap B_{1/2}$, and let $r = d_\Omega(x_1) \le d_\Omega(x_2)$. We separate now into two cases:
\begin{itemize}[leftmargin=*]
\item If $\rho:= |x_1 - x_2|\ge \frac{r}{2}$, then 
\[
\left|u(x_1) - u(x_2) \right|\le |u(x_1)|+|u(x_2)| \le Cr^s+ C(\rho+r)^s\le C \rho^s,
\]
thanks to \eqref{eq:boundbdry}. 
\item If $\rho :=|x_1- x_2|< \frac{r}{2}$ instead, then $B_r(x_1)\subset \Omega \cap B_1$ and $x_2 \in B_{r/2}(x_1)$. We consider the rescaled function $\tilde u_r(x)$ defined by 
\[
\tilde u_r(x) := \tilde u(x_1+rx), 
\]
which (thanks to \eqref{eq:scale_invariant_stable}) satisfies 
\[
\L \tilde u_r (x) = r^{2s} \tilde f(x_1 + rx) =: f_r(x) \quad\text{in}\quad B_1,
\]
 and 
$
\|\tilde u_r\|_{L^\infty(B_2)}\le C r^s 
$
from \eqref{eq:boundbdry}. Moreover, if we denote $x_*\in \partial\Omega$ such that $r = |x_1 - x_*|$, we have for any $x\in B_{1/(2r)}$,
\[
\begin{split}
|\tilde u_r(x)| & \le C d_\Omega^s(x_1+rx)\le C\left(|x_1-x_*|+r|x|\right)^s\\
& \le C\left(|x_1-x_*|^s+r^s|x|^s\right) \le C r^s(1+|x|^s). 
\end{split}
\]
Furthermore, $\|\tilde u_r\|_{L^\infty(\R^n)} = \|\tilde u\|_{L^\infty(\R^n)} \le 1$. Hence, $\tilde u_r\in L^\infty_s(\R^n)$, with $\|\tilde u_r\|_{L^\infty_s(\R^n)}\le C r^s$. Applying Theorem~\ref{thm-interior-linear-Lp} (see Remark~\ref{rem:int_weak}) with $C^s$ norm on the left-hand side, and $p = \infty$, we get
\[
[\tilde u_r]_{C^s(B_{1/2})} \le C \left(\|\tilde u_r\|_{L^\infty_s(\R^n)} + \|f_r\|_{L^\infty(B_{3/4})}\right) \le Cr^s.
\]
Since $[\tilde u_r]_{C^s(B_{1/2})} = r^s [\tilde u]_{C^s(B_{r/2}(x_1))}$, we deduce $[\tilde u]_{C^s(B_{r/2}(x_1))}\le C$ and $|\tilde u(x_1) - \tilde u(x_2)|\le C r^s$. 
\end{itemize}

In all, since $u = \tilde u$ in $B_{1/2}$, for any $x_1, x_2\in \Omega\cap B_{1/2}$ we have $| u(x_1) -   u(x_2)|\le C|x_1-x_2|^s$, which   gives the desired result.\qedhere
\end{steps}
\end{proof}

As a consequence, we have:

\begin{proof}[Proof of Theorem~\ref{thm:globCsreg}]
We cover $\Omega$ with finitely many balls $B_{1/2}(x_i)$, and apply Proposition~\ref{prop:globCsreg_loc} to each of them to get the desired result (notice that $u$ is bounded by Lemma~\ref{lem:Linftybound}).
\end{proof}

\begin{rem}
\label{rem:takedistinstead}
Theorem~\ref{thm:globCsreg} says that if $u$ is a weak solution, then it is continuous up to the boundary. From now one we can consider distributional solutions that are continuous up to the boundary and this will in particular include weak solutions (recall Lemma~\ref{lem:weak_distr}). Observe that Theorem~\ref{thm:globCsreg} and Proposition~\ref{prop:globCsreg_loc} also hold true for continuous distributional solutions (by using the comparison principle in Corollary~\ref{cor:comp_principle_G_d} instead of Corollary~\ref{cor:comp_principle_G_w}). 

This also allows us to expand our class of functions to nonenergetic solutions; see   Proposition~\ref{prop:globCsreg_loc_2} below. 
\end{rem}

\subsection{Hopf's lemma}
\index{Hopf's lemma}
We have shown that, given a $C^{1,\alpha}$ domain $\Omega$ and $\L u = f$ in $\Omega$, $u = 0$ in $\R^n\setminus \Omega$, then up to a constant we have $|u(x)| \le  d^s(x)$. We next show that, if moreover we assume $u \ge 0$ and $f\ge 0$, then we also have the lower bound $u(x) \ge cd^s(x)$ for some  $c> 0$, whenever $u\not\equiv 0$   in $\Omega$. 

\begin{prop}[Hopf's Lemma]
\label{prop:Hopf}
Let $s \in (0, 1)$ and let $\L \in \GLh$. Let $\alpha \in (0, 1)$, and let $\Omega$ be any  $C^{1,\alpha}$ domain.  Let $u\in L^\infty_{2s-\eps}(\R^n)\cap C(\overline{B_{1 }})$ for some $\eps > 0$ satisfy distributionally
\[
\left\{
\begin{array}{rcccll}
\L u & =  & f& \ge & 0& \quad\text{in}\quad \Omega \cap B_1\\
u &\ge& 0&&& \quad\text{in}\quad \R^n\setminus \Omega,
\end{array}
\right.
\]
for some $f\in L^\infty_{\rm loc}(\Omega)$. 
Then either $u \equiv 0$  in $\Omega$  or 
\[
u(x)\ge c d^s(x)\quad\text{for}\quad x \in B_{1/2}
\]
for some $c > 0$. 
\end{prop}

\begin{proof}
By the strong maximum principle (Theorem~\ref{thm:strong_maximum_principle}) we know that if $u\not\equiv 0$ in $\Omega$ then $u > 0$ in $\Omega$. Thus, let us assume $u > 0$ in $\Omega$. 

\begin{figure}
\centering
\makebox[\textwidth][c]{\includegraphics[scale = 1.0]{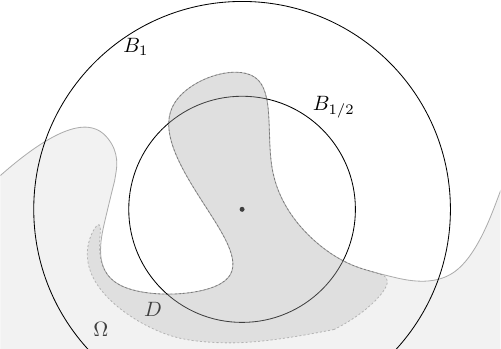}}
\caption{\label{fig:06} A possible example of the set $D$ from Hopf's lemma, Proposition~\ref{prop:Hopf}.}
\end{figure}

Let  $D\subset B_1$ with
\begin{equation}
\label{eq:from1dom}
\Omega\cap B_{1/2}\subset D \subset \Omega\cap B_{1}
\end{equation} 
 be a $C^{1,\alpha}$ domain (in particular, it has $C^{1,\alpha}$-radius $\varrho > 0$), see Figure~\ref{fig:06}. We use the subsolution $\varphi$ from Corollary~\ref{cor:subsol_domains}, which satisfies
  \begin{equation}
\label{eq:from2dom}
\left\{
\begin{array}{rcll}
\L \varphi & \le & -1& \quad\text{in}\quad N_\delta\\
 \varphi(x)& \ge & \delta d_D^s(x)  & \quad\text{for}\quad  x\in D, 
\end{array}
\right.
\end{equation}
for some $\delta > 0$, and where $N_\delta :=  \{0 < d_D(x) < \delta\}$. Let $c_*$ be defined as 
\[
c_* = \min\{u(x) : x\in D\setminus N_\delta\} > 0.
\]
Observe that $c_*>0$ since $u> 0$ in $D$ and $u$ is continuous. Then, we have that 
\[
 c_* \delta \L \varphi \le \L u\quad\text{in}\quad N_\delta\qquad\text{and}\qquad c_*\delta \varphi \le u\quad\text{in}\quad\R^n\setminus N_\delta.
\]
By the comparison principle, Corollary~\ref{cor:comp_principle_G_d} we have that $c_*\delta \varphi\le u$ in $B_1$. In particular, from \eqref{eq:from1dom} and \eqref{eq:from2dom} we get the desired result. 
\end{proof}

\begin{rem}
Proceeding as in the proof of the strong maximum principle, Theorem~\ref{thm:strong_maximum_principle}, we also have a version of Hopf's lemma for weak solutions.
\end{rem}

\begin{rem}
The strong maximum principle, Theorem~\ref{thm:strong_maximum_principle}, holds for distributional solutions with right-hand side $f\in L^1_{\rm loc}(\Omega)$. 
In this case, inside the proof one can use the unique weak solution to 
\[
\left\{
\begin{array}{rcll}
\L \bar u & = & \min\{1, f\} & \quad\text{in}\quad \Omega,\\
\bar u &  = & 0& \quad\text{in}\quad \R^n\setminus \Omega,
\end{array}
\right.
\]
as a barrier from below for $u$ when applying the comparison principle, Corollary~\ref{cor:comp_principle_G_d}. We are using here that $\bar u \in C(\overline{\Omega})$ thanks to Theorem~\ref{thm:globCsreg}.

As a consequence, in Hopf's lemma, Proposition~\ref{prop:Hopf}, we could assume instead that $f\in L^1_{\rm loc}(\Omega)$, and we would obtain the same result. 
\end{rem}

\subsection{Boundary regularity in $C^1$ and Lipschitz domains}

Let us now show how to obtain an estimate up to the boundary in more general domains ($C^1$ and Lipschitz), with operators whose kernels are not necessarily homogeneous.

In this case, solutions are not $C^s$ but only $C^\delta$ up to the boundary, for some small $\delta>0$.
On the other hand, in order to construct the barriers in Lipschitz domains,   we need a pointwise assumption on the kernels of the form
\begin{equation}
\label{kernels-Lipschitz-domains}
K(dy) = K(y)\, dy\quad\text{and}\quad \frac{\lambda}{|y|^{n+2s}} \leq K(y) \le \frac{\Lambda}{|y|^{n+2s}}\quad\text{in}\quad \R^n. 
\end{equation}
Our next estimate  reads as follows. 

\begin{prop}[Boundary regularity in $C^1$ and Lipschitz domains]
\label{prop:globCsreg_loc_Lip}\index{Boundary regularity!Operators comparable to fractional Laplacian!Lipschitz domains}\index{Boundary regularity!General operators!$C^1$ domains}
Let $s \in (0, 1)$, let $\L \in \GL$, and let $\Omega$ be any domain. Suppose, moreover, that
\begin{enumerate}[(i)]
\item\label{it:globCs_casei} either $\Omega$ is a $C^1$ domain, 
\item \label{it:globCs_caseii} or $\Omega$ is  a Lipschitz domain and $\L$ satisfies \eqref{kernels-Lipschitz-domains}.
\end{enumerate}

Let $f\in L^\infty(\Omega\cap B_1)$, and  $u\in L^\infty_{2s-\eps}(\R^n)$ for some $\eps > 0$ be any weak solution of 
\[
\left\{
\begin{array}{rcll}
\L u &=& f& \quad\text{in}\quad \Omega \cap B_1\\
u &=& 0& \quad\text{in}\quad B_1 \setminus \Omega.
\end{array}
\right.
\]
Then $u \in C_{\rm loc}^\delta(B_{1})$ with 
\[
\|u \|_{C^\delta(B_{1/2})}\le C\left( \|u\|_{L^\infty_{2s-\eps}(\R^n)}+ \|f\|_{L^\infty(\Omega\cap B_1)}\right)
\]
for some $\delta > 0$ and $C$ depending only on $n$, $s$, $\Omega$, $\eps$, $\lambda$, and $\Lambda$. 
\end{prop}

\begin{proof}
The proof is essentially the same as in Proposition~\ref{prop:globCsreg_loc}, using now the barrier constructed in Lemma~\ref{supersolution-d-eps-2} for case~\ref{it:globCs_casei}, or the barrier from  Lemma~\ref{supersolution-d-eps} in case~\ref{it:globCs_caseii}.
Namely, we take $\varphi:= C\dr_D^\eps$ for some $\eps>0$ small enough (see Definition~\ref{defi:distance}), and repeat the exact same proof of Proposition~\ref{prop:globCsreg_loc} (with exponent $\eps$ instead of $s$).

There is only one difference: the function $\varphi$ does not belong to $H^s(\R^n)$ now, and therefore we cannot directly use the comparison principle for weak solutions.
To solve this, we consider slightly larger domains $D_t\supset D$ for $t>0$, with $D_t$ being Lipschitz domains converging to $D$ as $t\to0$, and such that $\partial D_t\cap \partial\Omega=\varnothing$.
Then, the functions $\varphi_t:= C\dr_{D_t}^\eps$ are in the correct energy space (since they are smooth inside $D_t$, we have that  \eqref{eq:u_weak_sol} holds for $\varphi_t$) and thus by the comparison principle for weak solutions we deduce that $u\leq \varphi_t$ in $\Omega\cap B_{1/2}$.
Letting then $t\to0$ we get that $u\leq \varphi$, and the rest of the proof is the same.
\end{proof}
\begin{rem}
\label{rem:smallLipconst}
In fact, in part \ref{it:globCs_casei} one could instead consider Lipschitz domains with a sufficiently small Lipschitz constant; see Remark~\ref{rem:lemma_ch2_toadd}.   Namely, we say that a domain $\Omega$ has Lipschitz constant $\eta$ if for every $z\in \partial\Omega$ we have that, up to a rotation, $ \partial\Omega\cap B_\eta(z)$ is locally the graph of  a Lipschitz function with Lipschitz constant $\eta$. Then,  there exists some $\eta_\circ > 0$ depending only on $n$, $s$, $\lambda$, and $\Lambda$, such that Proposition~\ref{prop:globCsreg_loc_Lip} still holds swapping assumption \ref{it:globCs_casei} with the following:
\textit{\begin{enumerate}
\item[(i')] \label{it:globCs_casei_prime} $\Omega$ is a Lipschitz domain with Lipschitz constant $\eta \le \eta_\circ$.
\end{enumerate}}
\noindent The same holds for Corollary~\ref{cor:globCsreg_loc_Lip_2} below. 
\end{rem}

And we obtain the global regularity in this setting as a consequence. 

\begin{cor}
\label{cor:globCsreg_loc_Lip_2}
Let $s \in (0, 1)$, let $\L \in \GL$, and let $\Omega$ be any bounded domain. Suppose, moreover, that
\begin{enumerate}[(i)]
\item  either $\Omega$ is a $C^1$ domain, 
\item   or $\Omega$ is  a Lipschitz domain and $\L$ satisfies \eqref{kernels-Lipschitz-domains}.
\end{enumerate}
Let $f\in L^\infty(\Omega)$, and let $u$ be the weak solution to 
\[
\left\{
\begin{array}{rcll}
\L u &=& f& \quad\text{in}\quad \Omega\\
u &=& 0& \quad\text{in}\quad \R^n \setminus \Omega.
\end{array}
\right.
\]
Then $u \in C^\delta(\overline{\Omega})$ with 
\[
\|u \|_{C^\delta(\Omega)}\le C \|f\|_{L^\infty(\Omega)}
\]
for some $\delta > 0$ and $C$ depending only on $n$, $s$, $\Omega$, $\lambda$, and $\Lambda$. 
\end{cor}

\begin{proof}
Apply Proposition~\ref{prop:globCsreg_loc_Lip} to finitely many balls covering $\partial\Omega$, and combine with interior regularity estimates, Theorem~\ref{thm-interior-linear-Lp}. 
\end{proof}

\subsection{Existence of continuous distributional solutions}
\index{Boundary regularity!General operators!Continuous exterior datum}

As a consequence of the regularity up to the boundary for general operators from Corollary~\ref{cor:globCsreg_loc_Lip_2}, and thanks to the comparison principle for continuous distributional solutions, we obtain the general existence result for continuous distributional solutions: 

\begin{thm}
\label{thm:exist_dist_sol}
Let $s\in (0, 1)$, let $\L\in \GL$, and let $\Omega$ be any bounded $C^{1}$ domain. Let $g\in  L^\infty_{2s-\eps}(\R^n) $ for some $\eps > 0$ be continuous at all points on~$\partial \Omega$, and let $f\in L^\infty(\Omega)$. Then, there exists  a unique distributional solution $u\in C(\overline{\Omega})\cap L^\infty_{2s-\eps}(\R^n)$ to the equation 
\[
\left\{
\begin{array}{rcll}
\L u & =&  f& \quad\text{in}\quad \Omega\\
u & = & g& \quad\text{in}\quad \R^n\setminus \Omega. 
\end{array}
\right.
\]
\end{thm}
\begin{proof}
By considering $g\chi_{B_{2R}}$ instead of $g$, where $\Omega\subset B_R$, we can assume that $g$ is globally bounded (cf. the beginning of the proof of Proposition~\ref{prop:globCsreg_loc}). Thus, from now on, we will use $g\in L^\infty(\R^n)$. 

We divide the proof into three steps. 
\begin{steps}
\item Let $\bar g$ be any smooth extension of $g$ inside $\Omega$, that is, $\bar g\in L^\infty(\R^n) \cap C(\overline{\Omega})\cap C^\infty(\Omega)$ with $\bar g = g$ in $\R^n\setminus\Omega$. We can take, for example, the harmonic extension of $g$ inside $\Omega$.

Let $\delta > 0$, and let $\Omega_\delta := \{x\in \Omega : d_\Omega(x) >\delta\}$. We   consider the problem 
\[
\left\{
\begin{array}{rcll}
\L u_\delta & =&  f& \quad\text{in}\quad \Omega_\delta\\
u_\delta & = & \bar g& \quad\text{in}\quad \R^n\setminus \Omega_\delta,
\end{array}
\right.
\]
which now has a unique weak solution by Theorem~\ref{thm:exist_weak_sol}. This is because we   have that, since $\bar g$ is smooth near $\partial\Omega_\delta$,  
\[
\langle \bar g,\bar g\rangle_{K; \Omega_\delta}\le C_\delta < \infty, 
\]
for some $C_\delta$ that might blow-up as $\delta \downarrow 0$. 

Hence,  $u_\delta$ is globally bounded (by Lemma~\ref{lem:Linftybound}) and therefore it is also a distributional solution to $\L u_\delta = f$ in $\Omega_\delta$ (by Lemma~\ref{lem:weak_distr}). Thanks to the interior regularity estimates from Theorem~\ref{thm-interior-linear-Lp}, by Arzel\`a-Ascoli we also know that $u_\delta$ converges locally uniformly in $\Omega$ to some $u\in C(\Omega)\cap L^\infty(\R^n)$ with $u = g$ in $\R^n\setminus\Omega$, which by Proposition~\ref{prop:stab_distr} (where $\L_k$ is constant and equal to $\L$ for all $k$) is also a distributional solution to $\L u = f$. It only remains to be seen that $u$ is continuous at all points on $\partial \Omega$. 
\item We want to show that there exists some modulus of continuity $\omega$ independent of $\delta$ such that 
\[
|g(x) - u_\delta(y)|\le \omega(|x-y|)\quad\text{for all}\quad x\in \partial\Omega, \ y\in \Omega. 
\]
To do that, let us first build an appropriate barrier. 

Notice that the previous inequality already holds for any $y\in \Omega\setminus\Omega_\delta$ (with $\omega$ being the modulus of $\bar g$ on $\partial\Omega$). 

Let $\xi\ge 0$ be the unique weak solution, given by Theorem~\ref{thm:exist_weak_sol}, to 
\[
\left\{
\begin{array}{rcll}
\L \xi  & =&  1& \quad\text{in}\quad \Omega\\
\xi & = & 0& \quad\text{in}\quad \R^n\setminus \Omega.
\end{array}
\right.
\]
 By Corollary~\ref{cor:globCsreg_loc_Lip_2}, $\xi\in C^\gamma(\R^n)$ for some $\gamma > 0$, where $\gamma$ and the $C^\gamma$ estimates on $\xi$ depend only on $n$, $s$, $\Omega$, $\lambda$, and $\Lambda$.

Let now $x_\circ\in \partial\Omega$ fixed, and let us consider the function 
\[
\psi_{x_\circ}(x) := C\xi(x) + \frac{2|x-x_\circ|^2}{1+|x-x_\circ|^2}.
\]
Since $h_{x_\circ}(x) := \frac{2|x-x_\circ|^2}{1+|x-x_\circ|^2}$ is smooth and globally bounded, by Lemma~\ref{lem:Lu} there exists some $C_0$ universal (depending only on $n$ and $\Lambda$) such that 
\[
|\L  h_{x_\circ}|\le C_0\quad\text{in}\quad\R^n.
\]
Hence, we can choose $C$ large enough (depending on $C_0$) in the definition of $\psi_{x_\circ}$ such that 
\[
\L \psi_{x_\circ} \ge 1\quad\text{in}\quad\Omega. 
\]
Moreover, $\psi_{x_\circ}\in C^\gamma(\R^n)$, $\psi(x_\circ) = 0$, $\psi_{x_\circ} > 0$ in $\R^n\setminus\{x_\circ\}$, and $\psi_{x_\circ} \ge 1$ in $\R^n\setminus B_1(x_\circ)$. 
\item Let  $\eps > 0$, and let us define 
 \[
 w_\eps = g(x_\circ) + \eps + k_\eps \psi_{x_\circ},
 \]
where $k_\eps$ is chosen large enough (depending on $\eps$, but also on $\Omega$, $n$, $s$, $\lambda$, $\Lambda$, and on the modulus of continuity of $g$ on $\partial\Omega$) so that 
\[
 |\bar g-g(x_\circ)|\le \eps + k_\eps \psi_{x_\circ} \quad\text{in}\quad \R^n. 
\]
Hence, we have 
  \begin{equation}
  \label{eq:we1}
 w_\eps \ge \bar g \quad\text{in}\quad \R^n.
 \end{equation}
 By assuming $k_\eps \ge \|f\|_{L^\infty(\Omega)}$, we have from the fact that $\L \psi_{x_\circ} \ge 1$ in $\Omega$,
  \begin{equation}
  \label{eq:we2}
 \begin{split}
\L w_\eps \ge k_\eps \L (\psi_{x_\circ}) \ge \|f\|_{L^\infty(\Omega)} \quad\text{in}\quad \Omega.
 \end{split}
 \end{equation}
 
Thus, \eqref{eq:we1}-\eqref{eq:we2} let us apply the comparison principle for weak solutions in $\Omega_\delta$, Corollary~\ref{cor:comp_principle_G_w}, comparing $w_\eps$ and $u_\delta$, to deduce 
\[
u_\delta \le w_\eps\quad\text{in}\quad\R^n, \qquad\text{for all}\quad\delta > 0.  
\]

By continuity of $\psi_{x_\circ}$, for each $\eps > 0$ there exists some $\tilde\delta > 0$ (independent of $x_\circ$) such that $w_\eps \le g(x_\circ) +2\eps$ in $B_{\tilde \delta}(x_\circ)$.  This yields $u_\delta \le w_\eps \le g(x_\circ) +2\eps$ in $B_{\tilde \delta}(x_\circ)$, and letting $\delta\downarrow 0$ we obtain that for each $\eps > 0$ there exists some $\tilde\delta > 0$ such that
\[
u \le g(x_\circ) + 2\eps\quad\text{in}\quad B_{\tilde\delta}(x_\circ). 
\]
Repeating from below (defining $\tilde w_\eps = g(x_\circ) - \eps - k_\eps \psi_{x_\circ}$) we obtain that~$u$ is continuous at any $x_\circ\in  \partial\Omega$. Since it was continuous in $\Omega$ by interior regularity, we deduce that $u\in C(\overline{\Omega})$. Finally, the uniqueness follows from the comparison principle, Corollary~\ref{cor:comp_principle_G_d}.\qedhere
\end{steps}
\end{proof}

\subsection{Boundary regularity with unbounded right-hand side}
\index{Boundary regularity!General operators!Unbounded right-hand-side} \label{ssec:unboundedrhs}

 Let us now establish two boundary regularity results with a right-hand side that blows up as it approaches the boundary. 
They will be used in subsection~\ref{ssec:nonzerodatum} and Section~\ref{sec:higherorder} below.

  We start with solutions to an equation with right-hand side $f$ that may blow up at a rate $d^{\alpha-s}$. Notice that, $f$ belongs only to $L^p$ for $p > \frac{1}{s}$ and so, in general, one cannot talk about weak solutions in this setting (see Definition~\ref{defi:weaksol}). We show that, even if $f$ is not bounded, we can still recover the optimal $C^s$ regularity of solutions.

\begin{prop} \label{prop:globCsreg_loc_2}
Let $s \in (0, 1)$ and let $\L \in \GLh$. Let $\alpha >0$, and let $\Omega$ be any  $C^{1,\alpha}$ domain with $C^{1,\alpha}$-radius $\varrho_\circ>0$. Let $f d^{s-\alpha} \in L^\infty(\Omega\cap B_1)$, and let  $u\in C(B_1)\cap L^\infty_{2s-\eps}(\R^n)  $ for some $\eps > 0$ be a distributional solution to 
\[
\left\{
\begin{array}{rcll}
\L u &=& f& \quad\text{in}\quad \Omega \cap B_1\\
u &=& 0& \quad\text{in}\quad B_1 \setminus \Omega.
\end{array}
\right.
\]
Then $u \in C_{\rm loc}^s(B_{1})$ with 
\[
\|u \|_{C^s(B_{1/2})}\le C\left( \|u\|_{L^\infty_{2s-\eps}(\R^n)}+ \|f d^{s-\alpha}\|_{L^\infty(\Omega\cap B_1)}\right)
\]
for some $C$ depending only on $n$, $s$, $\alpha$, $\varrho_\circ$, $\eps$, $\lambda$, and $\Lambda$. 
\end{prop}

\begin{proof}
The proof is a modification of the proof of Proposition~\ref{prop:globCsreg_loc}. We highlight the only differences: 

\begin{itemize}[leftmargin=*]
\item We need to consider now (continuous) distributional solutions because otherwise the growth on the right-hand side is not compatible with the finite energy constraint of weak solutions. The interior regularity estimates were already stated for distributional solutions, whereas for the comparison principle we need to use Corollary~\ref{cor:comp_principle_G_d} instead of Corollary~\ref{cor:comp_principle_G_w}. 

\item In \ref{it:step1bdry} of the proof of Proposition~\ref{prop:globCsreg_loc} we used the barrier from Corollary~\ref{cor:supersol_domains} with right-hand side 1.   We now proceed slightly differently, to obtain again \eqref{eq:boundbdry}. Indeed, we consider $D\subset B_{7/8}$ similar in construction to the one in Proposition~\ref{prop:Hopf}: 
\[
\Omega\cap B_{5/6}\subset D\subset \Omega\cap B_{7/8},
\]
such that $D$ is a $C^{1,\alpha}$ domain. We split $u = u_1 + u_2$, where $u_1$ is the unique distributional solution in $C(\overline{D})$ (given by Theorem~\ref{thm:exist_dist_sol}) to 
\[
\left\{
\begin{array}{rcll}
\L u_1 &=&0& \quad\text{in}\quad D\\
u_1 &=& u& \quad\text{in}\quad \R^n\setminus D.
\end{array}
\right.
\]

In particular, we immediately have that $u_1\in C^s(B_{3/4})$ by Proposition~\ref{prop:globCsreg_loc} (rescaled) and Remark~\ref{rem:takedistinstead}, with
\[
\|u_1 \|_{C^s(B_{3/4})}\le C \|u\|_{L^\infty_{2s-\eps}(\R^n)}.
\]

On the other hand, $u_2$ is a distributional solution to
\[
\left\{
\begin{array}{rcll}
\L u_2 &=&f& \quad\text{in}\quad D\\
u_2 &=& 0& \quad\text{in}\quad \R^n\setminus D,
\end{array}
\right.
\]
where $fd_\Omega^{s-\alpha}\in L^\infty(\Omega\cap B_1)$, and since $D\subset \Omega$, we also have $fd_D^{s-\alpha}\in L^\infty(D)$. We can therefore use the barrier from Corollary~\ref{cor:supersol_domains} with right-hand side $d_D^{\alpha-s}$ together with the comparison principle Corollary~\ref{cor:comp_principle_G_d} to  obtain \eqref{eq:boundbdry} now too.

\item In \ref{it:step2bdry}, we proceed in the same way by observing that now 
\[
\|f_r\|_{L^\infty(B_{3/4})}\le C r^{2s}r^{ -s}\le  r^s,\]
and hence the same proof still works.  \qedhere
\end{itemize}
\end{proof}

 For the second result, we consider even more singular $f$, and we get a lower regularity of the solution. In this case, we can also deal with $C^1$ domains:

\begin{prop}
\label{prop:globCsreg_loc_3}
Let $s \in (0, 1)$ and let $\L \in \GLh$. Let $\alpha \in (0, s)$, and let $\Omega$ be any  $C^{1}$ domain. Let $f d^{2s-\alpha} \in L^\infty(\Omega\cap B_1)$, and let  $u\in C(B_1)\cap L^\infty_{2s-\eps}(\R^n)  $ for some $\eps > 0$ be a distributional solution to 
\[
\left\{
\begin{array}{rcll}
\L u &=& f& \quad\text{in}\quad \Omega \cap B_1\\
u &=& 0& \quad\text{in}\quad B_1 \setminus \Omega.
\end{array}
\right.
\]
Then $u \in C_{\rm loc}^\alpha(B_{1})$ with 
\[
\|u \|_{C^\alpha(B_{1/2})}\le C\left( \|u\|_{L^\infty_{2s-\eps}(\R^n)}+ \|f d^{2s-\alpha}\|_{L^\infty(\Omega\cap B_1)}\right)
\]
for some $C$ depending only on $n$, $s$, $\alpha$, $\Omega$, $\eps$, $\lambda$, and $\Lambda$. 
\end{prop}

\begin{proof}
The proof is a modification of the proof of Proposition~\ref{prop:globCsreg_loc_2} (which, at the same time, is modification of the proof of Proposition~\ref{prop:globCsreg_loc}). The differences with respect to the proof Proposition~\ref{prop:globCsreg_loc_2} are:

\begin{itemize}[leftmargin=*]

\item In the second bullet point, we use the barrier from Proposition~\ref{lem11_2} with right-hand side $\alpha -2s$ instead of that from Corollary~\ref{cor:supersol_domains}, to obtain in this case
\[
|\tilde u(x)| = |u(x)| \le   Cd_\Omega(x)^\alpha \quad\text{for}\quad x\in B_{3/4}.
\]

\item In the third bullet point, regarding \ref{it:step2bdry} of the proof of Proposition~\ref{prop:globCsreg_loc_2}, we proceed in the same way by observing that now 
\[
\|f_r\|_{L^\infty(B_{3/4})}\le C r^{2s}r^{ \alpha-2s}\le  r^\alpha,\]
and hence the same proof still works. \qedhere
\end{itemize}
\end{proof}
 
\begin{rem}
As in Remark~\ref{rem:smallLipconst}, one could consider instead Lipschitz domains with small Lipschitz constants. 
\end{rem}

\subsection{Boundary regularity with nonzero exterior data}
\label{ssec:nonzerodatum}
\index{Boundary regularity!General operators!H\"older exterior datum}

Up until now, we were considering problems with exterior data equal to zero (see Theorem~\ref{thm:globCsreg} and Proposition~\ref{prop:globCsreg_loc}). Thanks to the results in subsection~\ref{ssec:unboundedrhs}, we can extend them to more general exterior data. 

Given a bounded domain $\Omega\subset \R^n$, we   now consider an exterior datum $g:\R^n\setminus \Omega\to \R$ such that, for some $\alpha\in (0, \min\{1, 2s\})$ and $\eps\in (0, 2s-\alpha)$, 
\begin{equation}
\label{eq:gexterior}
 |g(x)|+\frac{|g(x) - g(y)|}{|x-y|^\alpha+|x-y|^{2s-\eps}}\le C_\circ  \qquad\text{for all}\quad x\in \partial \Omega, \ y \in \R^n\setminus \Omega. 
\end{equation}
Condition \eqref{eq:gexterior} basically says that $g$ is bounded on $\partial\Omega$, it is ``$C^\alpha$'' when comparing points on $\partial\Omega$ with points on $\R^n\setminus\Omega$ ``close'' to $\partial\Omega$, and it is $L^\infty_{2s-\eps}(\R^n)$ at infinity. For example, if $g\in C^\alpha(\R^n\setminus \Omega)$ it immediately satisfies \eqref{eq:gexterior}, whereas in general \eqref{eq:gexterior} also allows for discontinuities far away from $\partial\Omega$ and a growth up to a power $2s-\eps$ at infinity.

The regularity up to the boundary is then the following:  

\begin{prop} 
\label{prop:globCsreg_loc4}
Let $s \in (0, 1)$ and let $\L \in \GLh$. Let $\alpha \in (0, \min\{1, 2s\})$, $\eps\in (0, 2s-\alpha)$, and let $\Omega\subset \R^n$ be a bounded domain.  Let $f\in L^\infty(\Omega\cap B_1)$, let $g$ satisfy \eqref{eq:gexterior} for some $C_\circ > 0$, and let    $u\in C(B_1)$ be a distributional solution of
\[
\left\{
\begin{array}{rcll}
\L u &=& f& \quad\text{in}\quad \Omega \cap B_1\\
u &=& g& \quad\text{in}\quad B_1 \setminus \Omega.
\end{array}
\right.
\]
Then,
\begin{enumerate}[(i), leftmargin=*]
\item \label{case:i_g} If $\alpha < s$ and $\Omega$ is a $C^1$ domain, then $u \in C_{\rm loc}^\alpha(B_{1})$ with 
\[
\|u \|_{C^\alpha(B_{1/2})}\le C\left( C_\circ + \|f\|_{L^\infty(\Omega\cap B_1)} \right)
\]
for some $C$ depending only on $n$, $s$, $\alpha$, $\Omega$, $\eps$, $\lambda$, and $\Lambda$. 
\item \label{case:ii_g} If $\alpha > s$ and $\Omega$ is a $C^{1,\gamma}$ domain, then $u \in C_{\rm loc}^s(B_{1})$ with 
\[
\|u \|_{C^s(B_{1/2})}\le C\left( C_\circ + \|f\|_{L^\infty(\Omega\cap B_1)}\right)
\]
for some $C$ depending only on $n$, $s$, $\alpha$, $\Omega$, $\eps$, $\lambda$, and $\Lambda$. 
\end{enumerate}
\end{prop}
\begin{proof}
We divide the proof into three steps. 
\begin{steps}
\item Let $\bar g\in C^\infty(\Omega)\cap C^\alpha(\overline{\Omega})$ be a smooth extension of $g$ to $\Omega$ satisfying 
\begin{equation}
\label{eq:smoothgext}
|D^2 \bar g|\le C C_\circ d^{\alpha-2}\quad\text{in}\quad \Omega.
\end{equation}
For example, we can take $\bar g$ to be the solution to 
\[
\left\{
\begin{array}{rcll}
\Delta \bar g & = & 0&\quad\text{in}\quad\Omega,\\
\bar g & = & g& \quad\text{on}\quad \partial\Omega.  
\end{array}
\right.
\]
Indeed, by classical estimates for harmonic functions, we have 
\[
[\bar g]_{C^\alpha(\overline{\Omega})}\le C [g]_{C^\alpha(\partial \Omega)} \le CC_\circ,
\]
for some $C$ depending only on $n$, $\alpha$, and $\Omega$. Now, given $x_\circ\in \Omega$ with $r = d(x_\circ)$, we have 
\begin{equation}
\label{eq:bargtogether}
|D^2\bar g(x_\circ)|\le \|D^2 \bar g\|_{L^\infty(B_{r/2}(x_\circ))}\le \frac{C}{r^2} \osc_{B_r(x_\circ)}\bar g \le C C_\circ r^{\alpha -2},
\end{equation}
and thus \eqref{eq:smoothgext} holds. 

\item We now claim that 
\begin{equation}
\label{eq:thankstogggg}
|\L\bar g|\le CC_\circ d^{\alpha-2s}\quad\text{in}\quad \Omega.
\end{equation}
To prove this, we rescale the estimate from Lemma~\ref{lem:Lu} applied to $h(x) = \bar g(x_\circ + x) - \bar g(x_\circ)$ (also using $2s+\eps \le 2$ for $\eps > 0$ small enough, and that $\L$ is homogeneous), to get
\begin{equation}
\label{eq:bargtogether2}
r^{2s} \|\L \bar g \|_{L^\infty(B_{r/4}(x_\circ))}\le C\Lambda\left( r^2\|D^2 \bar g\|_{L^\infty(B_{r/2}(x_\circ))}+\|h(r\cdot)\|_{L^\infty_{2s-\eps}(\R^n)}\right).
\end{equation}

Notice now that, if $r|x| \ge 1$, 
\[
\frac{|h(rx)|}{1+|x|^{2s-\eps}} = \frac{|\bar g(x_\circ + rx) - \bar g(x_\circ)|}{r^{2s-\eps}|x|^{2s-\eps}}\frac{r^{2s-\eps}|x|^{2s-\eps}}{1+|x|^{2s-\eps}} \le 2C_\circ r^{2s-\eps},
\]
where we used \eqref{eq:gexterior} and $r^{2s-\eps}|x|^{2s-\eps}\ge r^\alpha|x|^\alpha$ (since $2s-\eps  >\alpha$). On the other hand, if $r|x| < 1$ we have 
\[
\frac{|h(rx)|}{1+|x|^{2s-\eps}} = \frac{|\bar g(x_\circ + rx) - \bar g(x_\circ)|}{r^{\alpha}|x|^{\alpha}}\frac{r^{\alpha}|x|^{\alpha}}{1+|x|^{2s-\eps}} \le 2C_\circ r^{\alpha},
\]
where we now used \eqref{eq:gexterior} and $r^{2s-\eps}|x|^{2s-\eps}< r^\alpha|x|^\alpha$, and the fact that $|x|^\alpha \le 1+|x|^{2s-\eps}$.

In all, we have that $\|h(r\cdot)\|_{L^\infty_{2s-\eps}(\R^n)}\le CC_\circ r^\alpha$ (since $r$ is bounded), and from \eqref{eq:bargtogether}-\eqref{eq:bargtogether2} we get 
\[
 \|\L \bar g \|_{L^\infty(B_{r/4}(x_\circ))}\le CC_\circ \Lambda r^{\alpha - 2s},
\]
that is, \eqref{eq:thankstogggg} holds, or $(\L \bar g) d^{\alpha-2s}\in L^\infty(\Omega)$. 

\item Finally, let $u = \bar g + w$, where $w$ satisfies 
\[
\left\{
\begin{array}{rcll}
\L w &=& f-\L\bar g& \quad\text{in}\quad \Omega \cap B_1\\
w &=& 0& \quad\text{in}\quad B_1 \setminus \Omega.
\end{array}
\right.
\]
We apply, in case \ref{case:i_g}, Proposition~\ref{prop:globCsreg_loc_3} to $w$ to get 
\[
\|w \|_{C^\alpha(B_{1/2})}\le C\left( \|w\|_{L^\infty_{2s-\eps}(\R^n)}+ \|f d^{2s-\alpha}\|_{L^\infty(\Omega\cap B_1)} + C_\circ\right)
\]
thanks to \eqref{eq:thankstogggg}. In case \ref{case:ii_g}, we apply instead Proposition~\ref{prop:globCsreg_loc_2} by writing $2s-\alpha = s - (\alpha -s)$, where $\alpha - s > 0$ to get
\[
\|w \|_{C^s(B_{1/2})}\le C\left( \|w\|_{L^\infty_{2s-\eps}(\R^n)}+ \|f d^{2s-\alpha}\|_{L^\infty(\Omega\cap B_1)} + C_\circ\right).
\]
Using that $w = u -\bar g$, and that $\bar g\in C^\alpha(\overline{\Omega})$, we get the desired result.\qedhere
\end{steps}
\end{proof}

\section{Higher-order boundary regularity}
\label{sec:higherorder}\index{Boundary regularity!General operators!Higher-order}

For second-order elliptic PDEs like
\[
\left\{
\begin{array}{rcll}
-\Delta  u  & =& f & \quad\text{in}\quad \Omega,\\
u & = & 0 & \quad\text{on}\quad \partial\Omega,
\end{array}
\right.
\]
 we know that if $f$ and $\partial\Omega$ are sufficiently smooth, then we can always obtain higher regularity of $u$ up to the boundary.
 Namely, if $f\in C^\alpha(\overline{\Omega})$ and $\partial\Omega$ is $C^{2+\alpha}$, with $\alpha\notin\mathbb N$, then $u\in C^{2+\alpha}(\overline{\Omega})$. 

In the case of integro-differential equations,  solutions to 
\begin{equation}
\label{eq:hoeq}
\left\{
\begin{array}{rcll}
\L u & =& f & \quad\text{in}\quad \Omega,\\
u & = & 0 & \quad\text{in}\quad \R^n\setminus \Omega,
\end{array}
\right.
\end{equation}
 are in general no better than $C^s(\overline{\Omega})$, even if $\L = \fls$, and $f$ and $\partial \Omega$ are~$C^\infty$. 
This is in contrast with the interior regularity: we have shown that higher regularity of $f$ and $\L$  do imply higher regularity of the solution inside the domain (see Theorem~\ref{thm-interior-linear-2}). In fact, already for $f\in L^\infty(\Omega)$, we have that the interior regularity is better than the boundary regularity, even in smooth domains!

A natural question is then: is there any higher-order boundary regularity estimate, analogous to that of the case $s = 1$?

At the moment we know that, close to the boundary, solutions to \eqref{eq:hoeq}   behave roughly as  
\[
u \asymp d^s,\qquad \text{where}\quad d(x) = \dist(x, \Omega^c). 
\]

We  wonder whether we can improve this expansion of the solution at boundary points. That is, if we denote
\[
u  = \eta \, d^s ,
\]
then thanks to the interior estimates, the $C^s$ regularity up to the boundary for $u$ can be understood as the boundedness of $\eta$. Is it now true that  higher regularity of $\partial\Omega$, $f$, and $\L$, gives higher regularity of $\eta$? 

The answer to this question is positive, and the best known results essentially show the following: for any $\beta > 1$, with $\beta, \beta\pm s\notin \N$, we have 
\[
\left\{
\begin{array}{l}
\partial\Omega\in C^\beta\\[1mm]
f\in C^{\max\{\beta-1-s, 0\}}(\overline{\Omega})\\[1mm]
\text{$K$ is homogeneous }\\[-0.3mm]
\text{and ``regular enough''}
\end{array}
\right.
\quad\Longrightarrow\qquad 
\frac{u}{d^s}\in C^{\beta-1}(\overline{\Omega});
\]
see \cite{RS-C1,RS-Cs} for $\beta<1+s$, and Grubb \cite{Grubb,Grubb2,Gru22}, Abatangelo and the second author \cite{AR20}, and Abels-Grubb \cite{AG}, for $\beta>1+s$.

In this section we will establish this result in the simplest case $\beta\in (1, 1+s)$.

We need to consider operators $\L \in \GL$ that are stable (i.e., with $K$ homogeneous, \eqref{eq:Lu_stab0}-\eqref{eq:Lu_stab}, or $\L\in \GLh$) and furthermore have absolutely continuous kernels satisfying the upper bound
\begin{equation}
\label{eq:boundabovekernel}
K(dy) = K(y)\, dy\quad\text{and}\quad K(y) \le \frac{\Lambda}{|y|^{n+2s}}\quad\text{for all}\quad y \in \R^n. 
\end{equation}

Under these conditions, we have:

\begin{thm}[Higher-order boundary regularity]
\label{thm:higher_bdry_reg}
Let $s \in (0, 1)$ and let $\L \in \GLh$  satisfy \eqref{eq:boundabovekernel}.  Let $\alpha \in (0, s)$, and $\Omega$ be any bounded $C^{1,\alpha}$ domain. Let $f\in L^\infty(\Omega)$, and let $u$ be any weak solution of
\[
\left\{
\begin{array}{rcll}
\L u &=& f& \quad\text{in}\quad \Omega\\
u &=& 0& \quad\text{in}\quad \R^n \setminus \Omega.
\end{array}
\right.
\]
Then $u/{d^s} \in C^\alpha(\overline{\Omega})$ with 
\[
\left\|\frac{u}{d^s} \right\|_{C^\alpha(\overline\Omega)}\le C \|f\|_{L^\infty(\Omega)},
\]
for some $C$ depending only on $n$, $s$, $\alpha$, $\Omega$, $\lambda$, and $\Lambda$. 
\end{thm}

When $s=1$, this result is equivalent to $u\in C^{1,\alpha}(\overline\Omega)$, which is the optimal boundary regularity in $C^{1,\alpha}$ domains.

The proof of Theorem \ref{thm:higher_bdry_reg} will be done via a contradiction and compactness argument.
For this, we first need a classification result for solutions in a half-space.

\subsection{Liouville's theorem in the half-space}

In the following, we prove a Liouville theorem in the half-space. 
We denote $\R^n_+ := \{x\in \R^n : x_n > 0\}$ and $\R^n_- = \R^n \setminus \R^n_+$.

\begin{thm}[Liouville theorem in the half-space]
\label{thm:Liouville_half_space}\index{Liouville's theorem!General operators!Half-space}
Let $s \in (0, 1)$ and let $\L \in \GLh$. Let $u\in C(\R^n)\cap L^\infty_{2s-\eps}(\R^n)$ for some $\eps > 0$ be any distributional solution of 
\[
\left\{
\begin{array}{rcll}
 \L u & = & 0 \quad\text{in}\quad \R^n_+,\\
  u & = & 0 \quad\text{in}\quad \R^n_-.
\end{array}
\right.
\]
Then 
\[
u(x) = \kappa  (x_n)^s_+
\]
for some $\kappa \in \R$. 
\end{thm}

\begin{proof}
Let $\mu = 2s-\eps$, and notice that  the function
\[
v_R(x) := R^{-\mu} u(Rx)\quad\text{for}\quad R \ge 1
\]
 satisfies 
\[
|v_R(x)| \le \|u\|_{L^\infty_\mu(\R^n)} R^{-\mu} \big(1+R^\mu|x|^\mu\big)\le  \|u\|_{L^\infty_\mu(\R^n)}\big(1+|x|^\mu\big),
\]
so that $\|v_R\|_{L^\infty_\mu(\R^n)}\le \|u\|_{L^\infty_\mu(\R^n)}$. 
On the other hand, by homogeneity  (see \eqref{eq:scale_invariant_stable}) we know that $\L v_R = \L\, u $ in $\R^n_+$. Hence, we can apply the boundary estimates from Proposition~\ref{prop:globCsreg_loc} (see Remark~\ref{rem:takedistinstead} or Proposition~\ref{prop:globCsreg_loc_2}) to deduce that
\[
[u]_{C^s(B_{R/2})} = R^{\mu-s} [v_R]_{C^s(B_{1/2})} \le C R^{\mu-s} \|u\|_{L^\infty_\mu(\R^n)}\quad\text{for}\quad R\ge 1.
\]

Let now $h \in \R^n$ such that $h_n = 0$. Let us define
\[
w(x) := \frac{u(x+h)-u(x)}{|h|^s}. 
\]
Then, by the previous consideration and if $|h|\le R/4$ we know 
\[
\|w\|_{L^\infty(B_R)} \le C \|u\|_{L^\infty_\mu(\R^n)} R^{\mu - s}\quad\text{for}\quad R\ge 1. 
\]
Observe, also, that $\L w = 0$ in $\R^n_+$ and $w = 0$ in $\R^n_-$.  In particular, we can repeat the previous argument with $R^s w$ to deduce
\[
[w]_{C^s(B_{R/2})} \le C R^{\mu-2s} \|u\|_{L^\infty_\mu(\R^n)}\quad\text{for}\quad R\ge 1.
\]
Letting $R\to \infty$, we deduce $w \equiv 0$ in $\R^n$ (since $w = 0$ in $\R^n_-$) and hence $u(x+h_\circ) = u(x)$ for any $h_\circ$ such that $(h_\circ)_n = 0$. In particular, $u$ depends only on the $x_n$ variable, and we have 
\[
u(x) = \bar u(x_n).
\]
By Lemma~\ref{lem:L1d}, $\bar u$ satisfies
\[
\left\{
\begin{array}{rcll}
 {(-\Delta)^s_\R} \bar u & = & 0 \quad\text{in}\quad (0,\infty),\\
  \bar u & = & 0 \quad\text{in}\quad (-\infty, 0].
\end{array}
\right.
\]
Since $|\bar u(t)| \le C(1+t^{2s-\eps})$, we are  done by the classification of one-dimensional solutions  from Theorem~\ref{thm:class_1d_fls}.
\end{proof}

We refer to \cite[Theorem 3.10]{AR20} for a higher-order version of this Liouville-type result in a half-space.

\subsection{Expansion around  boundary points}

In order to prove Theorem~\ref{thm:higher_bdry_reg}, the main point is to show the following expansion of the solution at  boundary points. 

\begin{prop}
\label{prop:expansion_bdry}
Let $s\in (0, 1)$, and let $\L\in \GLh$ satisfy  \eqref{eq:boundabovekernel}. Let $\alpha\in (0, s)$,  and let $\Omega\subset \R^n$ be any domain whose boundary is a $C^{1,\alpha}$ graph in $B_{1}$, with $C^{1,\alpha}$ norm bounded by $1$, and $0\in \partial\Omega$.  Let $u\in L^\infty_{2s-\eps}(\R^n)\cap C(B_1)$ for some $\eps > 0$ be a distributional solution to 
\[
\left\{
\begin{array}{rcll}
\L u & = & f& \quad\text{in}\quad \Omega\cap B_{1}\\
u & = & 0& \quad\text{in}\quad B_{1}\setminus \Omega,
\end{array}
\right.
\]
for some $f\in L^\infty(\Omega)$. Then, there exists $q_0\in \R$ such that
\[
\left\|u - q_0 d^s\right\|_{L^\infty(B_r)}\le C \left(\|f\|_{L^\infty(\Omega\cap B_1)}+ \|u\|_{L^\infty_{2s-\eps}(\R^n)}\right) r^{s+\alpha},
\]
 for all $r\in (0, 1)$.
The constant $C$ depends only on $n$, $s$, $\eps$, $\alpha$, $\lambda$, and $\Lambda$. 
\end{prop}

The proof of Proposition~\ref{prop:expansion_bdry} is somewhat similar to the interior regularity case (see subsection~\ref{ssec:compactness}), starting with a quantitative version of the Liouville-type result for solutions in very large balls and domains very close to a half-space:

\begin{prop}
\label{prop:compactness_bdry}
Let $s\in (0, 1)$, $\delta > 0$, and let $\L\in \GLh$ satisfy  \eqref{eq:boundabovekernel}. Let $\alpha\in (0, s)$,  and let $\Omega\subset \R^n$ be a domain whose boundary is a $C^{1,\alpha}$ graph in $B_{1/\delta}$, with $C^{1,\alpha}$ norm bounded by $\delta$, and $0\in \partial\Omega$.

Let $u\in L^\infty_{2s-\eps}(\R^n)\cap C(B_{1/\delta})$ for some $\eps > 0$ with $\|u\|_{L^\infty_{2s-\eps}(\R^n)}\le 1$ be a  distributional solution to 
\[
\left\{
\begin{array}{rcll}
\L u & = & f& \quad\text{in}\quad \Omega\cap B_{{1}/{\delta}}\\
u & = & 0& \quad\text{in}\quad B_{{1}/{\delta}}\setminus \Omega,
\end{array}
\right.
\]
for some $f$ with $\|f d^{ s-\alpha}\|_{L^\infty(\Omega)}\le \delta$. 

Then, for any $\eps_\circ > 0$ there exists $\delta_\circ > 0$ depending only on $\eps_\circ$, $n$, $s$, $\eps$, $\alpha$, $\lambda$, and $\Lambda$, such that if $\delta  <\delta_\circ$, 
\[
\left\|u-c_\circ d^s\right\|_{L^\infty(B_1)}\le \eps_\circ
\]
for some $c_\circ\in \R$. Moreover, we may take
\[
c_\circ = \frac{\int_{B_1} ud^s}{\int_{B_1}d^{2s}}. 
\]

\end{prop}

\begin{proof}
Let us argue by contradiction and let us assume that the statement does not hold. That is, there exists some $\eps_\circ > 0$ such that, for any $k\in \N$, there are:
\begin{itemize}
\item $\Omega_k\subset \R^n$ domains with boundary given by a $C^{1,\alpha}$ graph in $B_k$, with $C^{1,\alpha}$ norm bounded by $\frac{1}{k}$, and $0\in \partial\Omega_k$
\item $u_k\in L^\infty_{2s-\eps}(\R^n)\cap C(B_{k})$ with $\|u_k\|_{L^\infty_{2s-\eps}(\R^n)}\le 1$,
\item $f_k$ with $\|f_k d_k^{s-\alpha}\|_{L^\infty(\Omega_k)}\le  \frac{1}{k}$, where $d_k := d_{\Omega_k}$ denotes the distance to $\R^n\setminus \Omega_k$,
\item $\L_k\in \GLh$ satisfying  \eqref{eq:boundabovekernel},
\end{itemize} 
and they are such that 
\[
\left\{
\begin{array}{rcll}
\L_k u_k & = & f_k& \quad\text{in}\quad \Omega_k\cap B_{k}\\
u_k & = & 0& \quad\text{in}\quad B_{k}\setminus \Omega_k,
\end{array}
\right.
\]
in the distributional sense, with
\begin{equation}
\label{eq:contr_bdry}
\|u_k - c_k d_k^s\|_{L^\infty(B_1)} > \eps_\circ\quad\text{for all}\quad k\in \N,\quad c_k = \frac{\int_{B_1} u_kd_k^s}{\int_{B_1}d_k^{2s}},
\end{equation}
for some $\eps_\circ > 0$. 
From the estimates  of Proposition~\ref{prop:globCsreg_loc_2} in a ball $B_R$ with $R\ge 1$ fixed we know that, for $k$ large enough, 
\[
\|u_k \|_{C^{s}(B_{R/2})}\le C(R)\left(\|f_k d_k^{s-\alpha}\|_{L^\infty(\Omega_k\cap B_R)} + \|u_k\|_{L^\infty_{2s-\eps}(\R^n)}\right) \le C(R),
\]
and so $u_k$ converges locally uniformly in $\R^n$ to some $u_\infty\in L^\infty_{2s-\eps}(\R^n)\cap C(\R^n)$ with $\|u_\infty\|_{L^\infty_{2s-\eps}(\R^n)}\le 1$. 
Moreover, after a rotation we have that $\partial \Omega_k \to \partial \R^n_+$ (locally uniformly as graphs, or in the Hausdorff distance, up to a subsequence), and so 
\[u_\infty = 0 \quad \textrm{in}\quad \R^n_-.\] 
Notice, also, that $c_k \to c_\infty$, where
\[
c_\infty = \frac{\int_{B_1} u_\infty d_\infty^s}{\int_{B_1}d_\infty^{2s}}
\] 
and $d_\infty(x) = (x_n)_+$. 

On the other hand, $f_k\to 0$ as $k\to \infty$ locally uniformly in $\R^n_+$, and we can apply the stability of distributional solutions, Proposition~\ref{prop:stab_distr}, to deduce that 
\[\L_\infty u_\infty = 0\quad \textrm{in}\quad \R^{n}_+\]
for some $\L_\infty\in \GL$.  
Observe also that since $\L_k\in \GL$ are of the form \eqref{eq:Lu_stab} and satisfying \eqref{eq:boundabovekernel}, the same holds true for $\L_\infty$.

We can therefore apply the Liouville theorem in the half-space, Theorem~\ref{thm:Liouville_half_space}, to deduce that $u_\infty(x) = c_\infty (x_n)^s_+$, which contradicts \eqref{eq:contr_bdry} if $k$ is large enough.
\end{proof}

\subsubsection{The regularized distance}  In order to proceed with the proof of Proposition~\ref{prop:expansion_bdry}, we need to introduce the concept of \emph{regularized distance}, which will be necessary at exactly one point in the proof: given a domain $\Omega\subset \R^n$, we will work with powers of the distance to the boundary, which in general are not smooth functions (not even close to $\Omega$, unless $\Omega$ itself is smooth). In particular, in general it is not possible to evaluate (in the strong sense) $\L (d_\Omega)$.

\begin{defi}[Regularized distance]
\label{defi:distance}\index{Regularized distance}
Let $\Omega\subset \R^n$ be a $C^{\beta}$ domain with $\beta > 1$ and $\beta\notin\N$. We say that 
\[
\dr_\Omega:\R^n\to \R_{\ge 0}
\]
is a \emph{regularized distance} if $\dr_\Omega\in C^\infty(\Omega)\cap C^\beta(\overline{\Omega})$, there exists a constant $C$ such that 
\begin{equation}
\label{eq:compOmega}
d_\Omega \le \dr_\Omega \le C d_\Omega \quad\text{in}\quad \R^n,
\end{equation}
and 
\begin{equation}
\label{eq:compOmega2}
|D^k \dr_\Omega| \le C_k \dr_\Omega^{\beta-k} \quad\text{in}\quad \Omega,\quad\text{for all}\quad k \ge \beta
\end{equation}
for some $C_k$ depending only on $\Omega$, $k$, and $\beta$.

If $\Omega$ is  Lipschitz domain, we say that $\dr_\Omega:\R^n\to \R_{\ge 0}$ is a \emph{regularized distance} if $\dr_\Omega\in C^\infty(\Omega)\cap C^{0,1}(\overline{\Omega})$, and \eqref{eq:compOmega}-\eqref{eq:compOmega2} hold for $\beta = 1$.
\end{defi}

\begin{rem}\label{rem:quotient} Since $d_\Omega$ is $C^\beta$ close to $\partial\Omega$, $\dr_\Omega$ is $C^\beta(\overline{\Omega})$, and both $d_\Omega$ and $\dr_\Omega$ are vanishing linearly on $\partial\Omega$,   we have that
\[
\frac{d_\Omega}{\dr_\Omega} \in C^{\beta-1}(\overline{N_\delta}),
\]
where $N_\delta := \{0 < d_\Omega < \delta\}$ and $\delta$ depends only on the $C^\beta$-radius of $\Omega$ (see Definition~\eqref{defi:varrho}). 
Indeed, after flattening the boundary it is enough to show it when $\Omega = \{x_n > 0\}$ locally. In this case, we can write 
\[
\frac{u}{x_n} = \int_0^1 u_{x_n}(x', tx_n)\, dt \in C^{\beta-1}
\]
for $u = d_\Omega, \dr_\Omega$, since $u$ is $C^\beta$ in $\Omega$, and $u_{x_n}\in C^{\beta-1}$. Hence, since furthermore $u/{x_n}$ is uniformly positive and bounded in $\Omega$, we get that 
\[
\frac{d_\Omega}{\dr_\Omega}  =\frac{d_\Omega}{x_n}\frac{x_n}{\dr_\Omega} \in C^{\beta-1}(\overline{N_\delta}),
\]
as claimed.
\end{rem}

\begin{rem}\label{rem:dregisenough}
Thanks to Remark~\ref{rem:quotient}, in  Proposition~\ref{prop:expansion_bdry} we can equivalently take the regularized distance $\dr_\Omega$ in place of $d_\Omega$. Indeed, since $w := {\dr_\Omega}/{d_\Omega}\in C^\alpha(B_{r_\circ})$ for some $r_\circ$ universal (by Remark~\ref{rem:quotient}), and  $1\le w \le C$ in $\Omega$, we also have that $w^s\in C^\alpha(B_{r_\circ})$ and
\[
\left\|\frac{\dr^s_\Omega}{d^s_\Omega} - w(0)\right\|_{L^\infty(B_r)}\le Cr^{\alpha}\quad\Rightarrow \quad \left\|{\dr^s_\Omega}- w(0){d^s_\Omega} \right\|_{L^\infty(B_r)}\le Cr^{s+\alpha},
\]
for any $r\in (0, r_\circ)$.  
By the triangle inequality we now have that the statement holds for $\dr_\Omega$ if and only if it holds for $d_\Omega$. 
\end{rem}

For $\beta > 1$ and $\beta\notin\N$, any $C^\beta$ (or Lipschitz) domain has a regularized distance (see Lemma~\ref{lem:distance} in Appendix~\ref{app.B}).

\subsubsection{Proof of the expansion around boundary points} 

We can next proceed with the proof of Proposition~\ref{prop:expansion_bdry}. 

\begin{proof}[Proof of Proposition~\ref{prop:expansion_bdry}]
 We will show it for $\dr_\Omega$ instead of $d_\Omega$, which is sufficient thanks to Remark~\ref{rem:dregisenough}. 

We divide the proof into five steps. It will eventually follow by contradiction, with two preliminary steps.

\begin{steps}
\item \label{step:1comp} We start with the following claim: 
\begin{claim*}
Let $\mu > s$, let $v$ be such that $\|v\|_{L^\infty(B_1)}\le 1$ and vanishes on a domain $B_1\setminus \Omega$ with Lipschitz constant bounded by 1, and let $\dr_\Omega$ denote the regularized distance (according to Definition~\ref{defi:distance} and Lemma~\ref{lem:distance}) corresponding to $\Omega$. 
If for any $r \in (0, 1)$ there exists $q_r\in \R$ such that 
\[
\|v - q_{r} \dr_\Omega^s\|_{L^\infty(B_r)} \le C_1 r^\mu,
\]
then
\begin{equation}
\label{eq:1circ}
\|v - c_\circ \dr_\Omega^s\|_{L^\infty(B_r)}\le C_\circ r^\mu\quad\text{for all}\quad r \in (0, 1),
\end{equation}
for some $c_\circ$ and $C_\circ$ bounded by $C(C_1+1)$, where $C$ depends only on $n$ and~$s$.
\end{claim*}
To prove this, observe first that, since $\|\dr_\Omega^s \|_{L^\infty(B_r)}\le c r^s$ for some dimensional constant $c$, we have for any $r \in \left(0, \frac12\right]$ and $\beta\in [1, 2]$,
\begin{equation}
\label{eq:limitmus}
\begin{split}
|q_{\beta r}-q_{r}| & \le C r^{-s}\|q_{\beta r}\dr_\Omega^s-q_{r}\dr_\Omega^s\|_{L^\infty(B_r)}\\
& \le C r^{-s}\|v-q_{ r}\dr_\Omega^s\|_{L^\infty(B_r)}+C r^{-s}\|v-q_{ \beta r}\dr_\Omega^s\|_{L^\infty(B_{\beta r})}\\
& \le C C_1 r^{\mu-s}. 
\end{split}
\end{equation}
On the other hand, since $\|v\|_{L^\infty(B_1)}\le 1$, we have that 
\[
|q_{r}|\le C(C_1+1)\quad\text{for all}\quad r\in {\textstyle \left[\frac12, 1\right]}.
\]
Together with \eqref{eq:limitmus} and the fact that $\mu > s$, this implies the existence of a limit 
\[
 q_0:=\lim_{r\downarrow 0} q_{ r},
\]
which satisfies
\[
|q_0 - q_{ r}|\le \sum_{j = 0}^\infty |q_{ 2^{-j}r}  - q_{ 2^{-j-1}r}| \le CC_1 r^{\mu-s}\sum_{j = 0}^\infty 2^{-j(\mu-s)}\le C C_1 r^{\mu-s}
\]
for any $r\in (0, 1)$. In particular, 
\[
|q_0|\le C(C_1 + 1),
\]
and, for any $r\in (0, 1)$,
\[
\begin{split}
\|v - q_0 \dr_\Omega^s\|_{L^\infty(B_r)} & \le\|v - q_{r} \dr_\Omega^s\|_{L^\infty(B_r)} +|q_{r}-q_0 |\|\dr_\Omega^s\|_{L^\infty(B_r)} \\
& \le C(C_1 +1)r^\mu.
\end{split}
\]
That is, \eqref{eq:1circ} holds with $c_\circ = q_0$ and $C_\circ = C(C_1+1)$.

\item   \label{step:2comp} We now prove the following:
\begin{claim*}
Let  $\mu > s$. Let $\{u_m\}_{m\in \N}$ with $\|u_m\|_{L^\infty(B_1)}\le 1$ for all $m\in \N$ be a family of functions such that each $u_m$ vanishes on a domain $B_1\setminus \Omega_m$ with Lipschitz constant bounded by 1. Let $\dr_m = \dr_{\Omega_m}$ denote the regularized distance to $B_1\setminus \Omega_m$. 

If there exist $(q_{m, r})_{m\in \N, r\in (0, 1)}$ such that 
\[
\sup_{m \in \N}\sup_{r\in(0,1)} r^{-\mu} \|u_m - q_{m, r}\dr_m^s\|_{L^\infty(B_r)} = \infty,
\]
then there exist sequences $r_k \downarrow 0$ and $m_k \to \infty$ such that 
\begin{equation}
\label{eq:forlater}
r_k^{-\mu}\|u_{m_k} - q_k \dr_{m_k}^s\|_{L^\infty(B_{r_k})} \to \infty\quad\text{as}\quad k\to \infty,
\end{equation}
and
\begin{equation}
\label{eq:def_vk}
v_k(x) := \frac{\left(u_{m_k} - q_k \dr_{m_k}^s\right)(r_k x)}{\|u_{m_k} - q_k \dr_{m_k}^s\|_{L^\infty(B_{r_k})}} 
\end{equation}
satisfies $\|v_k\|_{L^\infty(B_R)} \le C R^\mu$ for any $1\le R \le \frac{1}{2r_k}$, for some $C$ that depends only on $n$.
\end{claim*}

In order to show the claim, let us define
\[
Q_{m, r}:= \|u_{m} - q_{m, r} \dr_{m}^s\|_{L^\infty(B_{r})}
\quad\text{and}\quad
\theta(\rho) := \sup_{m\in \N}\, \sup_{1>r > \rho} \, \rho^{-\mu}Q_{m, \rho},
\]
so that $\theta(\rho)$ is nonincreasing, finite for any $\rho > 0$, and by assumption $\theta(\rho)\uparrow \infty$ as $\rho\downarrow 0$. In particular, there exist sequences $m_k\in \N$ and $r_k\downarrow 0$ such that
\begin{equation}
\label{eq:forlater2}
r_k^{-\mu}Q_{m_k, r_k} \ge \frac12\theta(r_k) \to \infty\quad\text{as}\quad k\to \infty.
\end{equation}
On the other hand, as in \eqref{eq:limitmus} we have
\[
|q_{m, 2r}-q_{m, r}|\le Cr^{-s}(Q_{m, r}+Q_{m, 2r})\le Cr^{\mu-s}\theta(r)
\]
and, for $R = 2^N$ for some $N\in \N$, 
\[
|q_{m, Rr}-q_{m, r}|\le C\sum_{j = 0}^{N-1} (2^j r)^{\mu-s}\theta(2^jr)\le   C R^{\mu-s}r^{\mu-s}\theta(r)
\]
(and similarly for any $R\le \frac{1}{2r}$).
Thanks to this, if we define $v_k$ as \eqref{eq:def_vk} we have, for any $1 \le R \le \frac{1}{2r_k}$,
\[
\begin{split}
\|v_k\|_{L^\infty(B_R)} & = Q_{m_k, r_k}^{-1}\|u_{m_k} - q_{m_k, r_k} d_{m_k}^s\|_{L^\infty(B_{Rr_k})} \\
& \le Q_{m_k, r_k}^{-1}\left( Q_{m_k, Rr_k} + |q_{m_k, Rr_k} - q_{m_k, r_k} |\|d_{m_k}^s\|_{L^\infty(B_{Rr_k})} \right)\\
& \le  Q_{m_k, r_k}^{-1} \theta(Rr_k)R^\mu r_k^\mu + CQ_{m_k, r_k}^{-1} R^{\mu-s}r_k^{\mu-s}\theta(r_k) (Rr_k)^s\\
& \le 2 \theta(Rr_k) [\theta(r_k)]^{-1}R^\mu  + C R^\mu r_k^\mu Q_{m_k, r_k}^{-1}\theta(r_k)\\
& \le C R^\mu,
\end{split}
\]
for some $C$ that depends only on $n$, where we have also used the monotonicity of $\theta$ and \eqref{eq:forlater2}.

\item After the two preliminary steps, let us start with the body of the proof. Up to dividing by a constant, we assume that 
\[
\|f\|_{L^\infty(\Omega\cap B_1)}+ \|u\|_{L^\infty_{2s-\eps}(\R^n)}\le 1.
\]
We now argue by contradiction, and assume that the statement does not hold. In particular, there exist sequences of 
\begin{itemize}
\item $\Omega_m\subset \R^n$ domains with boundary given by a $C^{1,\alpha}$ graph in $B_1$ and $C^{1,\alpha}$ norm bounded by 1, with $0\in \partial\Omega_m$, and where we denote by $\dr_m := \dr_{\Omega_m}$ the regularized distance corresponding to $\Omega_m$,
\item $u_m\in L^\infty_{2s-\eps}(\R^n)\cap C({B_1})$ with $\|u_m\|_{L^\infty_{2s-\eps}(\R^n)}\le 1$,
\item $f_m\in L^\infty(\Omega_m)$ with $\|f_m\|_{L^\infty(\Omega)}\le 1$,
\item $\L_m \in \GLh$ operators  satisfying \eqref{eq:boundabovekernel}, 
\end{itemize}
such that 
\[
\left\{
\begin{array}{rcll}
\L_m u_m & = & f_m& \quad\text{in}\quad \Omega_m\cap B_{1}\\
u_m & = & 0& \quad\text{in}\quad B_{1}\setminus \Omega,
\end{array}
\right.
\]
and
\[
\sup_{r\in (0, 1)} r^{-s-\alpha} \|u_m - c_m \dr_m^s\|_{L^\infty(B_r)} \ge m \quad\text{for any}\quad c_m\in \R. 
\]
By \ref{step:1comp} (with $\mu = s+\alpha$) we have that for any family $(q_{m, r})_{m\in \N, r\in (0, 1)}$,
\[
\sup_{m \in \N}\sup_{r\in(0,1)} r^{-s+\alpha} \|u_m - q_{m, r}\dr_m^s\|_{L^\infty(B_r)} = \infty.
\]
We choose our $q_{m,r}$ as:
\[
q_{m, r} := \frac{\int_{B_{r}} u_{m} \dr_{m}^s}{\int_{B_{r}} \dr_{m}^{2s}}.
\]

By \ref{step:2comp},  and using the notation there, there exist sequences $m_k\to \infty$ and $r_k\downarrow 0$ such that if we define $v_k$ as in \eqref{eq:def_vk} with the $q_{m, r}$ above,
\[
v_k(x) =  \frac{\left(u_{m_k} - q_k \dr_{m_k}^s\right)(r_k x)}{Q_k},\quad \text{with}\quad q_k =  \frac{\int_{B_{r_k}} u_{m_k} \dr_{m_k}^s}{\int_{B_{r_k}} \dr_{m_k}^{2s}},
\]
and
\[
Q_k := Q_{m_k,r_k} = \|u_{m_k} - q_k \dr_{m_k}^s\|_{L^\infty(B_{r_k})},
\]
then 
\[
\|v_k\|_{L^\infty(B_R)}\le C R^{s+\alpha}  \quad\text{for}\quad 1 \le R\le \frac{1}{2r_k} 
\]
and
\[
\|v_k\|_{L^\infty(B_1)} = 1. 
\]

Notice that, since $|u_{m_k}|\le C \dr_m^s$ by Proposition~\ref{prop:globCsreg_loc}, we have $|q_k|\le C $ for some $C$ depending only on $n$, $s$, $\lambda$, $\Lambda$, $\alpha$, and $\eps$

Furthermore, by assumption we have $|u_{m_k}(x)|\le 1+|x|^{2s-\eps}$ in $\R^n$, and so we can also bound $v_k$ outside of $B_{1/(2r_k)}$ by 
\[
|v_k(x)|\le Q_k^{-1}\left(1+r_k^{2s-\eps}|x|^{2s-\eps} + |q_k| r_k^s|x|^s\right) \le C  Q_k^{-1} r_k^{2s-\eps}|x|^{2s-\eps}
\]
for all $x\in \R^n\setminus B_{1/(2r_k)}$. Assuming, without loss of generality, that $\eps < s- \alpha$, then $Q_k^{-1} r_k^{2s-\eps}\to 0$ as $k\to \infty$ by \eqref{eq:forlater}, and we obtain 
\begin{equation}
\label{eq:vkboundmu}
\|v_k\|_{L^\infty_{2s-\eps}(\R^n)}\le C
\end{equation}
for some $C$ depending only on $n$, $s$, $\lambda$, $\Lambda$, $\alpha$, and $\eps$.

\item Let us now see what equation each $v_k$ satisfies. If we denote 
\[
\tilde\Omega_k := \frac{1}{r_k}\Omega_{m_k},\quad\text{and}\quad \tilde \dr_k(x) := \frac{1}{r_k} \dr_{m_k}(r_k x),
\]
 we have 
\[
\left\{
\begin{array}{rcll}
\L_{m_k} v_k & = & \tilde f_k(x)& \quad\text{in}\quad \tilde\Omega_k\cap B_{1/r_k},\\
v_k & = & 0& \quad\text{in} \quad B_{1/r_k}\setminus \tilde\Omega_k,
\end{array}
\right.
\]
where 
\[
\tilde f_k(x):= Q_k^{-1}r_k^{2s} \left(f_{m_k}(r_kx) - q_k (\L_{m_k} \dr^s_{m_k})(r_kx)\right). 
\]
This is the point where we use the fact that we have the regularized distance and not the standard distance. Indeed, we need to make sense of the term $\L_{m_k} \dr^s_{m_k}$, which is only well defined if $\dr_{m_k}$ is smooth enough. Then, we can use Proposition~\ref{lem11} (and the fact that $|f_{m_k}|\le 1$) to obtain a bound for~$\tilde f_k$,
\[
|\tilde f_k \tilde \dr_k^{s-\alpha}|(x) \le Q_k^{-1}r_k^{2s}|\tilde \dr_k^{s-\alpha}(x)| + CQ_k^{-1}r_k^{2s}|q_k| |\dr^{\alpha-s}_{m_k} (r_kx)| |\tilde \dr_k^{s-\alpha}(x)|.
\]
Since $\dr_{m_k} (r_kx) = r_k \tilde \dr_k (x)$ and $\tilde \dr_k(x) \le \frac{1}{r_k}$ in $B_{1/r_k}$, we get 
\[
|\tilde f_k \tilde d_k^{s-\alpha}| \le CQ^{-1}_k r_k^{s+\alpha}\quad \textrm{in}\quad  \tilde\Omega_k \cap B_{1/r_k},\]
where we have also used that $|q_k|$ is bounded. 

In all, 
\begin{equation}
\label{eq:rhstozero}
\big\|\tilde f_k \tilde \dr_k^{s-\alpha}\big\|_{L^\infty(\tilde\Omega_k\cap B_{1/r_k})} \le CQ^{-1}_k r_k^{s+\alpha} \to 0\quad\text{as}\quad k \to \infty, 
\end{equation}
by \eqref{eq:forlater} (recall $\mu = s+\alpha$ and the definition of $Q_k$). 

\item 
Let us denote 
\[
\begin{split}
d_{m_k}(x) &:= \dist(x, \R^n\setminus \Omega_{m_k}),\\
\tilde d_{k}(x)& := \dist(x, \R^n\setminus {r^{-1}_k} \Omega_{m_k}),
 \end{split}
\] 
so that $\tilde d_k (x) = \frac{1}{r_k} d_{m_k}(r_k x)$. By  Remark~\ref{rem:dregisenough}, we have that
\[
\big\|{\dr_{m_k}} - w_k(0){d_{m_k}}\big\|_{L^\infty(B_r)}\le C r^{1+\alpha},
\]
for any $r\in (0, r_\circ)$, with $w_k := {\dr_{m_k}}/{d_{m_k}}\in C^\alpha(B_1)$. In terms of $\tilde d_k$ and $\tilde\dr_k$ we have 
\begin{equation}
\label{eq:convdist}
\big\|{\tilde \dr_{k}} - w_k(0){\tilde d_{k}}\big\|_{L^\infty(B_1)}\le C r_k^{\alpha},
\end{equation}
that is, $\tilde\dr_k$ converges locally uniformly to $w_k(0)\tilde d_k$, as $k\to \infty$. 

On the other hand, thanks to \eqref{eq:vkboundmu}-\eqref{eq:rhstozero}, up to dividing by a universal constant depending only on $n$ and $s$, we have that  $v_k$ satisfies the hypotheses of Proposition~\ref{prop:compactness_bdry} for any $\delta > 0$ provided that $k$ is large enough. In particular, for any $\eps_\circ>0$ there exists some $k$ large enough such that 
\[
\big\|v_k - c_k \tilde d_k^s\big\|_{L^\infty(B_1)}\le \eps_\circ,\quad\text{for some}\quad 
c_k \in \R.
\]
Combining this with \eqref{eq:convdist}, we get
\begin{equation}
\label{eq:integrating-inequality}
\left\|v_k - \frac{c_k}{w_k(0)} \tilde \dr_k^s\right\|_{L^\infty(B_1)} \le 2\eps_\circ,
\end{equation}
if $k$ is large enough. Notice, however, that by definition of $v_k$, $\dr_k$, and $q_k$,
\[
\int_{B_1} v_k \tilde \dr^s_k = Q_k^{-1}r_k^{-1-n}\left( \int_{B_{r_k}} u_{m_k}  \dr_{m_k}^s- q_k \int_{B_{r_k}} \dr_{m_k}^{2s}\right) = 0, 
\]
so that integrating \eqref{eq:integrating-inequality} against $\tilde \dr_k^s$ in $B_1$ we get
$
|c_k| \le C \eps_\circ
$
 for some $C$ independent of $k$ (using that $\int_{B_1} \tilde d_k^{2s}$ and $w_k^{-1}(0)$ are uniformly bounded below). Thus
\[
\|v_k\|_{L^\infty(B_1)}\le C\eps_\circ,
\]
for $k$ large enough. 
However,  by definition of $v_k$  we have $\|v_k\|_{L^\infty(B_1)} = 1$, which is a contradiction if $\eps_\circ < 1/C$. This completes the proof. 
\qedhere
\end{steps}
\end{proof}

\subsection{Proof of the higher-order  boundary regularity}

Let us finish with the proof of the higher-order boundary regularity result, Theorem~\ref{thm:higher_bdry_reg}. As in the proof of the boundary regularity, we start with a local version of the estimate, that has interest on its own: 
\begin{prop}
\label{prop:local_bdry}
Let $s\in (0, 1)$, and let $\L\in \GLh$ satisfy  \eqref{eq:boundabovekernel}. Let $\alpha\in (0, s)$,  and let $\Omega\subset \R^n$ be any domain whose boundary is a $C^{1,\alpha}$ graph in $B_{1}$, with $C^{1,\alpha}$ norm bounded by $1$. 

Let $u\in L^\infty_{2s-\eps}(\R^n)\cap C(B_1)$ for some $\eps > 0$  be any distributional solution to 
\[
\left\{
\begin{array}{rcll}
\L u & = & f& \quad\text{in}\quad \Omega\cap B_{1}\\
u & = & 0& \quad\text{in}\quad B_{1}\setminus \Omega,
\end{array}
\right.
\]
for some $f\in L^\infty(\Omega)$. Then, $u/{d^s}\in C^\alpha_{\rm loc}(B_1\cap \overline{\Omega})$ with
\[
\left\|\frac{u}{d^s}\right\|_{C^\alpha( \overline{\Omega}\cap B_{1/2})}\le C \left(\|f\|_{L^\infty(\Omega\cap B_1)}+  \|u\|_{L^\infty_{2s-\eps}(\R^n)}\right) 
\]
for some $C$ depending only on $n$, $s$, $\eps$,  $\alpha$, $\lambda$, and $\Lambda$. 
\end{prop}
 \begin{proof}
 We assume $ \|f\|_{L^\infty(\Omega\cap B_1)}+  \|u\|_{L^\infty_{2s-\eps}(\R^n)} \le  1$ after dividing by a constant if necessary. We split the proof into two steps.

 \begin{steps}
 \item \label{it:step1_1111} From Proposition~\ref{prop:expansion_bdry} we know that for each $z\in \partial\Omega\cap B_{1/2}$ there exists some $q = q(z)$   such that  
\[
 \|u - q(z) d^s\|_{L^\infty(B_r(z))}\le C r^{s+\alpha}\quad\text{for all}\quad r\in (0, {\textstyle \frac12}).
\]
Equivalently, thanks to Remark~\ref{rem:dregisenough}, we have  for some $q'(z)$,
 \begin{equation}
 \label{eq:from1111_2}
 \|u - q'(z) \dr^s\|_{L^\infty(B_r(z))}\le C r^{s+\alpha}\quad\text{for all}\quad r\in (0, {\textstyle \frac12}).
 \end{equation}
Let now $x_\circ \in \Omega\cap B_{1/2}$ with $r_\circ := d(x_\circ) \le \rho_\circ$ for some small universal $\rho_\circ$, and let $z_\circ\in \partial \Omega\cap B_{3/4}$ be a projection of $x_\circ$ towards $\partial\Omega$ (in particular, $|x_\circ-z_\circ| = r_\circ$). 
 
 We define the function
 \[
 v (x) := u(x_\circ + r_\circ x) - q'(z_\circ) \dr^s(x_\circ + r_\circ x). 
 \]
 Then, we have (by \eqref{eq:from1111_2})
 \[
 \|v\|_{L^\infty(B_R)}\le C (Rr_\circ)^{s+\alpha}\quad\text{for all}\quad 1\le R \le cr_\circ^{-1}. 
 \]
 On the other hand, from $\|u\|_{L^\infty_{2s-\eps}(\R^n)}\le 1$,  
 \[
 \|v\|_{L^\infty(B_R)}\le C (Rr_\circ)^{2s-\eps}\quad\text{for all}\quad R \ge cr_\circ^{-1}. 
 \]
  In all, assuming without loss of generality that $\eps < \alpha - s$, we have 
 \[
 \|v\|_{L^\infty_{2s-\eps}(\R^n)}\le C r_\circ^{s+\alpha}. 
 \]
 
 Observe, also, that thanks to Proposition~\ref{lem11}, for all $x\in B_{3/4}$,
 \[
 \begin{split}
 |\L v(x)|& \le r_\circ^{2s} f(x_\circ+r_\circ x) + C|q'(z_\circ)|r_\circ^{2s} \dr^{\alpha-s}(x_\circ + r_\circ x)
 \\
 & \le r_\circ^{2s} + C r_\circ^{s+\alpha} \le C r_\circ^{s+\alpha},
 \end{split}
 \]
 where we have  used that $q'(z_\circ)$ is bounded and that, for $x\in B_{3/4}$, the function $\dr(x_\circ+r_\circ x)$ is comparable to $r_\circ$. 
 Hence, by $C^\alpha$ interior estimates with bounded right-hand side  (Theorem~\ref{thm-interior-linear-Lp} rescaled to balls $B_{1/2}$ and $B_{3/4}$) we have:
 \[
 [v]_{C^\alpha(B_{1/2})}\le C\left(\|v\|_{L^\infty_{2s-\eps}(\R^n)} +  \|\L v \|_{L^\infty(B_{3/4})} \right) \le C r_\circ^{s+\alpha}. 
 \]
 In terms of $u$ this is
 \begin{equation}
 \label{eq:from2222_2}
 [u - q'(z_\circ) \dr^s]_{C^\alpha(B_{r_\circ/2}(x_\circ))}\le C r_\circ^s. 
 \end{equation}
 
 \item We now want to obtain a bound for the H\"older norm of ${u}/{d^s}$ from \eqref{eq:from2222_2}. To do that, two preliminary observations are in order. 
 
 We will use that, on the one hand, for any $u_1, u_2\in C^\alpha(D)$ we have $u_1u_2\in C^\alpha(D)$ with 
 \begin{equation}
 \label{eq:boundprodholder}
 [u_1 u_2]_{C^\alpha(D)} \le \|u_1\|_{L^\infty(D)}[u_2]_{C^\alpha(D)}+[u_1]_{C^\alpha(D)}\|u_2\|_{L^\infty(D)},
 \end{equation}
 (see \eqref{eq:APP_In1} in Proposition~\ref{prop:A_imp}). 
 
 On the other hand, we will also use that, since
\[
 \begin{split}
&  \|\dr^{-s}\|_{L^\infty(B_{r_\circ/2}(x_\circ)) } \le C r_\circ^{-s},\\
 & \|\nabla (\dr^{-s})\|_{L^\infty(B_{r_\circ/2}(x_\circ))}   = C \|\dr^{-s-1}\nabla \dr\|_{L^\infty(B_{r_\circ/2}(x_\circ))}\le Cr_\circ^{-s-1},
 \end{split}
\]
we have that by interpolation 
 \begin{equation}
 \label{eq:boundprodholder_s}
\begin{split}
  [\dr^{-s}]_{C^\alpha(B_{r_\circ/2}(x_\circ))}&  = \sup_{x, y\in B_{r_\circ/2}(x_\circ)} \frac{|\dr^{-s}(x) - \dr^{-s}(y)|^\alpha}{|x-y|^\alpha}|\dr^{-s}(x) - \dr^{-s}(y)|^{1-\alpha} \\
  & \le   C[\dr^{-s}]^\alpha_{C^{0,1}(B_{r_\circ/2}(x_\circ))}  \|\dr^{-s}\|_{L^\infty(B_{r_\circ/2}(x_\circ)) }^{1-\alpha}\le Cr_\circ^{-s-\alpha}.
  \end{split}
 \end{equation}
 Combining \eqref{eq:boundprodholder} with \eqref{eq:boundprodholder_s} in \eqref{eq:from2222_2} (also recalling \eqref{eq:from1111_2}) we get
 \[
 \begin{split}
 \left[\frac{u}{\dr^s}\right]_{C^\alpha(B_{r_\circ/2}(x_\circ))} & =  \left[\frac{u}{\dr^s}-q'(z_\circ)\right]_{C^\alpha(B_{r_\circ/2}(x_\circ))}\\
 & \le \|u - q'(z_\circ) \dr^s\|_{L^\infty(B_{r_\circ/2}(x_\circ))} \left[\dr^{-s}\right]_{C^\alpha(B_{r_\circ/2}(x_\circ))} \\
 & \quad +[u - q'(z_\circ) \dr^s]_{C^\alpha(B_{r_\circ/2}(x_\circ))}\|\dr^{-s}\|_{L^\infty(B_{r_\circ/2}(x_\circ))}\\
 & \le C r_\circ^{s+\alpha}r_\circ^{-s-\alpha} + C r_\circ^{s}r_\circ^{-s} = C.
 \end{split}
 \]
 That is, we have obtained a universal bound for $\left[{u}/{\dr^s}\right]_{C^\alpha(B_{r_\circ/2}(x_\circ))} $ for any $r_\circ\le \rho_\circ$. Since $u$ has interior $C^{2s}$ regularity and $\dr^s$ is smooth and uniformly positive away from the  boundary, the previous estimate also holds for $\rho_\circ \le r_\circ \le \frac12$:
 \begin{equation}
 \label{eq:frompreveq_2}
  \left[\frac{u}{\dr^s}\right]_{C^\alpha(B_{r_\circ/2}(x_\circ))} \le C \quad\text{for all} \ \ x_\circ\in  B_{1/2}\ \ \text{with}\ \ B_{r_\circ(x_\circ)}\subset \Omega\cap B_1.
 \end{equation}
 
 Using again \eqref{eq:boundprodholder}, since $u/\dr^s$ and $\dr^s / d^s$ are $C^\alpha$ in $\Omega$ (see Remark~\ref{rem:dregisenough}) we get from \eqref{eq:frompreveq_2},
\[
  \left[\frac{u}{d^s}\right]_{C^\alpha(B_{r_\circ/2}(x_\circ))} \le C \quad\text{for all} \ \ x_\circ\in  B_{1/2}\ \ \text{with}\ \ B_{r_\circ(x_\circ)}\subset \Omega\cap B_1.
\]
 This is enough to conclude that $u/d^s\in C^\alpha(B_{1/2}\cap \overline{\Omega})$ with a bound of the form  $\|u/d^s\|_{C^\alpha(B_{1/2}\cap \overline{\Omega})}\le C$ (see, e.g., Lemma~\ref{lem:Lipholder} in Appendix~\ref{app.A}). \qedhere
 \end{steps}
 \end{proof}

From the previous result, we finally obtain the following:  
 
\begin{proof}[Proof of Theorem~\ref{thm:higher_bdry_reg}]
We cover $\Omega$ with finitely many balls $B_{1/2}(x_i)$, and apply Proposition~\ref{prop:local_bdry} to each of them to get the desired result.
\end{proof}

\section{Further results and open problems}

In this chapter, we have obtained a quite complete understanding of the regularity of solutions to linear nonlocal equations of order $2s$. 
Still, there are many interesting research directions that we did not discuss yet; we briefly do it next.

\subsection{General equations in divergence form}
\label{sect-div-open}

Under appropriate regularity assumptions on $K(x,z)$ in the $x$-variable, we  know that solutions of divergence-form equations 
\begin{equation}\label{divergence-eq}
\L(u,x)=0 \quad \textrm{in}\quad \Omega\subset \R^n,
\end{equation}
with operators as in \eqref{divergence-form0}-\eqref{divergence-form1},  satisfy similar regularity estimates to the ones we established in Theorems \ref{thm-interior-linear-2} and \ref{thm-interior-linear-Lp}; see Section~\ref{sec:x-dependence} and \cite{FR5,Fal}.  

Recall that, in this case, the ellipticity conditions should read as
\begin{equation}\label{divergence-ellipticity0B}
r^{2s} \int_{B_{2r}(x)\setminus B_r(x)} K(x,dz) \le \Lambda
\end{equation}
and 
\begin{equation}\label{divergence-ellipticity1B}
\lambda \le 
r^{2s-2}\inf_{e\in \S^{n-1}}\int_{B_{\Lambda r}(x)\setminus B_r(x)}|e\cdot (x-z)|^{2} K(x,dz)
\end{equation}
for all $r>0$ and $x\in \R^n$, with $\Lambda\geq \lambda>0$.

A quite different problem concerns the regularity of solutions of divergence-form equations \eqref{divergence-eq}-\eqref{divergence-form0}-\eqref{divergence-form1}, with ellipticity assumptions \eqref{divergence-ellipticity0B}-\eqref{divergence-ellipticity1B}, but with \emph{no} regularity assumption in $x$.
These are called equations with \emph{bounded measurable coefficients}.

The main regularity question in this context is to understand whether solutions to these equations are H\"older continuous or not\footnote{See also \cite{KMS, KNS22} for regularity results in $L^p$ spaces.}.
The first results in this direction were established by Kassmann \cite{Kas09} (see also Chen-Kumagai \cite{CH03}, Bass-Levin \cite{BL02b}, and Caffarelli-Chan-Vasseur \cite{CCV}), who proved that, under the stronger assumption 
\begin{equation}\label{divergence-unif}
 0<\frac{\bar\lambda}{|x-z|^{n+2s}} \leq K(x,z)\leq \frac{\bar\Lambda}{|x-z|^{n+2s}},
\end{equation}
any solution of  \eqref{divergence-eq}-\eqref{divergence-form0}-\eqref{divergence-form1} in $B_1$ satisfies
\begin{equation}\label{divergence-holder}
\|u\|_{C^{0,\alpha}(B_{1/2})} \leq C\|u\|_{L^\infty(\R^n)},
\end{equation}
for some positive constants $\alpha$ and $C$ depending only on $n$, $s$, $\bar\lambda$ and $\bar \Lambda$.

Notice that the techniques we used in this chapter in order to establish regularity estimates do not apply in case of equations with bounded measurable coefficients, and completely different methods are required.
The main estimate in \cite{Kas09} is proved by developing a nonlocal version of the classical Moser iteration.

After the results in \cite{Kas09,CH03,BL02b,CCV}, the natural question was then to ask whether  \eqref{divergence-holder} holds for more general kernels, not satisfying \eqref{divergence-unif}.
The best known result in this direction follows from the works of Chaker-Silvestre \cite{CS20}, Dyda-Kassmann \cite{DK20}, and Imbert-Silvetre \cite{IS}, and states that \eqref{divergence-holder} holds for any solution of \eqref{divergence-eq}-\eqref{divergence-form0}-\eqref{divergence-form1}, with ellipticity assumptions \eqref{divergence-ellipticity0B}-\eqref{divergence-ellipticity1B}, as long as the kernel $K(x,z)$ satisfies the additional assumption
\[
\big|\{z\in B: K(x,z)\geq \mu |x-z|^{-n-2s} \}\big|\geq \theta|B|
\]
for any ball $B\subset \R^n$ and $x\in B$, with $\theta\in(0,1)$ and $\mu>0$.
Notice that this assumption is much weaker than \eqref{divergence-unif}, but still leaves the following:

\vspace{2mm}

\noindent\textbf{Open question 2.1}:
\textit{Does the H\"older estimate \eqref{divergence-holder} hold for all solutions of \eqref{divergence-eq}-\eqref{divergence-form0}-\eqref{divergence-form1} in $B_1$, under the  general ellipticity assumptions \eqref{divergence-ellipticity0B}-\eqref{divergence-ellipticity1B}?}

\vspace{2mm}

Thanks to the results in \cite{IS}, the question essentially reduces to showing that, under the ellipticity assumptions \eqref{divergence-ellipticity0B}-\eqref{divergence-ellipticity1B}, the coercivity estimate
\[ \int_{\R^n} \int_{\R^n} \frac{|u(x)-u(z)|^2}{|x-z|^{n+2s}}\,dzdx \leq C\int_{\R^n} \int_{\R^n} |u(x)-u(z)|^2 K(x,z) dzdx \]
holds for all $u\in H^s(\R^n)$, with $C$ depending only on the constants $\lambda,\Lambda$ appearing in \eqref{divergence-ellipticity0B}-\eqref{divergence-ellipticity1B}.

This is a challenging and very natural question, which arises as well in the context of the Boltzmann equation; see \cite{CS20} for more details.

\subsection{Harnack's inequality}

In Chapter \ref{ch:fract_Lapl} we saw that the square root of Laplacian satisfies a Harnack inequality; see Proposition~\ref{prop:Harnack_sqrt}.
It is then natural to wonder if more general operators $\L \in \GL$ --- or more general operators of the form \eqref{divergence-form0}-\eqref{divergence-form1}-\eqref{divergence-ellipticity0B}-\eqref{divergence-ellipticity1B} --- satisfy a similar Harnack inequality \cite{DKP14,DKP16,Cozzi17}.

It turns out that the condition
\[
 0<\frac{\bar\lambda}{|y|^{n+2s}} \leq K(y)\leq \frac{\bar\Lambda}{|y|^{n+2s}}
\]
is enough for the Harnack inequality to hold:
\begin{equation}\label{Harnack}
\left\{\begin{array}{rcll}
\L u &=& 0 & \quad\textrm{in}\quad B_1 \\
u &\geq& 0 & \quad\textrm{in}\quad \R^n 
\end{array}\right.
\qquad \Longrightarrow\qquad
\sup_{B_{1/2}} u \leq C\inf_{B_{1/2}} u,
\end{equation}
with $C$ depending only on $n$, $s$, $\bar \lambda$, $\bar \Lambda$ (see Corollary~\ref{cor-Harnack}).

However, there are (linear and translation invariant) operators $\L\in\GL$ for which such Harnack inequality fails.
A simple example is given by $\L = (-\partial_{x_1x_1}^2)^s+(-\partial_{x_2x_2}^2)^s$ in $\R^2$; see \cite{BC09,BS2}.
This raises the following:

\vspace{2mm}

\noindent\textbf{Open question 2.2}:
\textit{Can one characterize those operators $\L$ for which the Harnack inequality \eqref{Harnack} holds?}

\vspace{2mm}

In case of stable operators, i.e., when $K$ is homogeneous, this question was completely solved by Bogdan and Sztonyk; see \cite{BS2,BS}.
On the other hand, for general operators $\L \in \GL$, a characterization of those operators satisfying a \emph{parabolic} Harnack inequality follows from the results of Chen, Kumagai, and Wang \cite{CKW20}.

\subsection{Boundary regularity estimates in $L^p$ spaces}

Concerning the boundary regularity of solutions, here we have studied solutions $u$ to the Dirichlet problem 
\begin{equation}\label{Dir-pb-Ch2}
\left\{
\begin{array}{rcll}
\L u &=& f& \quad\text{in}\quad \Omega\\
u &=& 0& \quad\text{in}\quad \R^n \setminus \Omega,
\end{array}
\right.
\end{equation}
for operators $\L\in\GL$ with homogeneous kernels $K$, i.e. $\L\in \GLh$, and 
with $f\in L^\infty(\Omega)$ or $f\in C^\alpha(\overline\Omega)$.
In this case, the best known results essentially show the following: for any $\beta > 1$, with $\beta, \beta\pm s\notin \N$, we have 
\[
\left\{
\begin{array}{l}
\partial\Omega\in C^\beta\\
f\in C^{\max\{\beta-1-s, 0\}}(\overline{\Omega})\\
\text{$K|_{\mathbb S^{n-1}}$ is ``regular enough''}
\end{array}
\right.
\quad\Longrightarrow\qquad 
\frac{u}{d^s}\in C^{\beta-1}(\overline{\Omega}). 
\]
This sharp higher regularity result was first proved in the case $\beta = \infty$ by Grubb \cite{Grubb, Grubb2} (see also Hormander \cite{Hor65,Hor07}), then by the second author and Serra \cite{RS-C1} in case $\beta\in(1,1+s)$, and finally by Abatangelo and the second author in \cite{AR20}, for all $\beta>1+s$; see also Abels-Grubb \cite{AG}.

In case $f\in L^p(\Omega)$, the best known results have been established by Abels and Grubb \cite{AG}, where they proved fine estimates in Bessel-potential type $H^s_q$ spaces.
As a particular case of their results, one has:
\[
\left\{
\begin{array}{l}
\partial\Omega\in C^2\\
f\in L^p({\Omega}),\quad ps>1\\
K|_{\mathbb S^{n-1}}\in C^\infty(\mathbb S^{n-1})
\end{array}
\right.
\quad\Longrightarrow\qquad 
\left.\frac{u}{d^s}\right|_{\partial\Omega}\in W^{s-\frac1p,\,p}(\partial{\Omega}),
\]
where $u/d^s|_{\partial\Omega}$ should be understood in a trace sense; see \cite[Theorem 4.5]{AG}.
Here, $W^{\alpha,p}$ denotes the fractional Sobolev space; \eqref{Walphap}.   See also \cite{DR23} for another related result in weighted Sobolev spaces for general stable operators in $C^{1,\alpha}$ domains.

An interesting question that remains completely open in this context is to establish $L^p$-based regularity estimates for $u/d^s$ in Lipschitz domains.

\vspace{2mm}

\noindent\textbf{Open question 2.3}:
\textit{When $\Omega$ is Lipschitz and $\L$ is a stable operator, can one prove that solutions to \eqref{Dir-pb-Ch2} satisfy $u/d^s\in L^p(\partial\Omega)$?}

\vspace{2mm}

This is not known even for the fractional Laplacian, $\L=(-\Delta)^s$, in which case the most general global results say that $u$ belongs to the  Besov space $B_{2, \infty}^{\min\{2s, s+1/2\}}(\Omega)$ when $f\in L^2(\Omega)$, \cite{BN23}.

\subsection{Non-symmetric operators}

The operators $\L\in \GL$ that we have considered in this Chapter are all \emph{symmetric}, in the sense that $\int_{\R^n} u\,\L v = \int_{\R^n} \L u\,v$ for all $u,v\in C^\infty_c(\R^n)$.
In terms of the kernel $K$ this is equivalent to the symmetry condition $K(y)=K(-y)$, which is what allowed us to write the operator \eqref{eq:L_def} as \eqref{eq:Lu_nu}.

The natural generalization of $\GL$ to the context of non-symmetric kernels is the class of operators of the type
\[\L u(x)= \left\{ \begin{array}{ll}
\displaystyle 
\int_{\R^n} \big(u(x)-u(x+y)\big)K(dy) 
 & \quad \textrm{if}\ s\in(0,{\textstyle \frac12})    \\
\displaystyle  
{\rm P.V.}\int_{\R^n} \big(u(x)-u(x+y)\big)K(dy) +b\cdot \nabla u(x) 
  & \quad \textrm{if}\ s={\textstyle \frac12} \\
\displaystyle 
\int_{\R^n} \big(u(x)-u(x+y)+\nabla u(x)\cdot y\big)K(dy)
    & \quad \textrm{if}\ s\in({\textstyle \frac12},1),
\end{array}
\right.\]
where $b\in \R^n$, and the kernel $K\geq0$ satisfies the additional cancellation property
\[\int_{B_{2r}\setminus B_r} y\,K(dy)=0\quad \textrm{in case}\ s={\textstyle \frac12}.\]
The corresponding uniform ellipticity assumptions are
\[r^{2s}\int_{B_{2r}\setminus B_r} K(dy) \leq \Lambda\qquad \textrm{for all}\quad r>0\]
and
\[0<\lambda \leq r^{2s-2}\int_{B_{2r}\setminus B_r} |e\cdot y|^2K(dy)\qquad \textrm{for all}\quad r>0 \quad \textrm{and}\quad e\in \S^{n-1}.\]

These operators have Fourier symbols of the form $\mathcal A+i\mathcal B$, with
\[\mathcal A(\xi)\asymp |\xi|^{2s}\qquad \textrm{and}\qquad |\mathcal B(\xi)|\lesssim |\xi|^{2s}\]
for all $\xi\in \R^n$. 
Thanks to this, the proofs of the interior regularity results we presented in Section \ref{sec:int_reg_G} still hold for this more general class of operators; see \cite[Theorem~3.8]{DRSV}.

On the other hand, the boundary regularity of solutions that we presented in Sections \ref{sec:bdryregularity}  and \ref{sec:higherorder} used very strongly the fact that the kernels under consideration were {symmetric} (and homogeneous).
Namely, a crucial ingredient in those proofs is the fact that the function $(x\cdot e)_+^s$ is a 1D solution of $\L w=0$ in $\{x\cdot e>0\}$, for any $e\in \S^{n-1}$.

For non-symmetric operators, the boundary regularity was developed in \cite[Theorem~1.2]{DRSV}, where Dipierro, the second author, Serra, and Valdinoci developed a boundary regularity theory which extends the one for symmetric operators that we presented in Sections \ref{sec:bdryregularity}  and \ref{sec:higherorder} here.
A key observation is that, for any non-symmetric operator $\L$ with homogeneous kernel, and for any $e\in \S^{n-1}$, there exists an explicit exponent 
\[\gamma(e,\L):= s+ \frac{1}{\pi}\arctan \left(\frac{\mathcal B(e)}{\mathcal A(e)}\right),\qquad \gamma(e,\L) \in (0,2s)\cap (2s-1,1), \]
for which the function $(x\cdot e)_+^{\gamma(e,\L)}$ is a 1D solution of $\L w=0$ in $\{x\cdot e>0\}$.

An interesting problem that remains open is the following:

\vspace{2mm}

\noindent\textbf{Open question 2.4}:
\textit{Develop a higher-order boundary regularity theory for nonsymmetric operators, extending the one in \cite{AR20} for symmetric operators.}

\vspace{2mm}

 This would have applications to obstacle problems for non-symmetric operators, too; see subsection~\ref{ssec:high_reg_fb}.

%% file: chap3.tex
%
%
%

\chapter{Fully nonlinear equations}
\label{ch:fully_nonlinear}

In this chapter we study fully nonlinear equations of the type
\begin{equation}\label{ch3-0}
\mathcal I u = 0 \quad \textrm{in}\quad \Omega\subset\R^n,
\end{equation}
where $\mathcal I$ is an operator of the form
\begin{equation}\label{ch3-1}
\mathcal I w := \inf_{\gamma\in \Gamma} \big\{ -\L_\gamma w\big\}\qquad \textrm{or}\qquad 
\mathcal I w := \inf_{b\in \mathcal B}\sup_{a\in\mathcal A} \big\{ -\L_{ab} w\big\},
\end{equation}
and $\L_\gamma$ or $\L_{ab}$ are linear integro-differential elliptic operators of order $2s$.

When $s=1$, these correspond to fully nonlinear elliptic PDE,
\[F(D^2 u) = 0 \quad \textrm{in}\quad \Omega\subset\R^n,\]
and their study has been a major research direction since the second half of the 20th century, with many important contributions by Nirenberg, Krylov, Safonov, Evans, Caffarelli, and many others; see \cite{CC, NTV14, FR4} for more details.

For integro-differential operators, the regularity theory for fully nonlinear equations of the type \eqref{ch3-0}-\eqref{ch3-1} was developed in a series of famous papers by Caffarelli and Silvestre \cite{CS,CS2,CS3}, and many more results have been established since then; see \cite{CD12,CD16,C-LK,DZ,GS,GS2,JX,KM,Kriv,Mou,RS-Duke,RS-bdryH,Ser,Ser2,SS,Hui}.

Here we will present first the basic results concerning the existence and uniqueness of (viscosity) solutions for these equations, to then establish some of the main known regularity results.

\section{Preliminaries}

\index{Fully nonlinear integro-differential equations}

We next present a very heuristic motivation for the equations \eqref{ch3-0}-\eqref{ch3-1}, and we refer to \cite{Lions1,Lions2} for a rigorous derivation in the case of second-order elliptic equations.

Let $\Gamma=\{1,...,N\}$, and $\{X_t^i\}_{i\in\Gamma}$ be a finite collection of L\'evy processes, with infinitesimal generators $\L_i$. 
Given a bounded smooth domain $\Omega\subset \R^n$, and a payoff function $g\in C(\R^n\setminus \Omega)$, we saw in subsection~\ref{subsect-2.1.4} that if we consider the process $x+X_t^i$ (i.e., a ``particle'' starts at $x\in \Omega$ and moves randomly according to $X_t^i$), then the \emph{expected payoff} of the first hitting time~$\tau_x$ (the first time the particle falls outside $\Omega$), denoted $u_i(x)$, is given by the solution of
\[
\left\{
\begin{array}{rcll}
\L_i u_i & = & 0& \quad\text{in}\quad \Omega\\
u_i & = & g& \quad\text{in}\quad \R^n\setminus \Omega.
\end{array}
\right.
\]

Here, we consider instead the following \emph{control} problem: suppose that at any time $t$ we can choose to change the parameter $i\in \Gamma$ so that  the particle that has reached the position $x$ at time $t$  must then move according to $X_t^i$ during a very small time interval $(t,t+\delta)$. 
We get   a process $X_t^{\gamma}$, where $\gamma:\Omega\to \Gamma$ may actually change from point to point.

Our goal is to minimize\footnote{Notice that, if we wanted to \emph{maximize} the payoff, we could simply consider $-g$ instead of $g$.} the expected payoff, by making the optimal choice of $i$ at each $x\in \Omega$.
The minimum expected payoff is then given by 
\[
u(x) := \inf_{\substack{\text{all possible choices}\\ \text{of } \gamma:\Omega\to \Gamma}} \E[g(X_{\tau_x}^{\gamma})]. 
\] 
Observe  that, since we can always choose to continue with a constant $\gamma\equiv i$, we have $u(x) \le \E\left[u(x+X_t^i)\right]$ for every $x\in \Omega$ (and every $i\in \Gamma$), and from \eqref{eq:Ptdef}-\eqref{eq:def-infinit-generat} we get $-\L_i u \ge 0$ in $\Omega$ for all $i\in \Gamma$. 
On the other hand, assuming that $u$ is regular enough, one can show that for short time intervals we will have $u(x) = \E\left[u(x+X_t^i)\right]+o(t)$ for some $i\in \Gamma$, so that (again by \eqref{eq:Ptdef}-\eqref{eq:def-infinit-generat}) at each $x\in \Omega$ we  have that $\L_i u(x) = 0$ for some $i\in \Gamma$. 
The previous expressions can be put into a single (nonlinear) equation as 
\[
\min_{i\in \Gamma}\big\{-\L_i u\big\}  = 0.
\]

One can repeat the same argument for any general family of stochastic processes $\{X_t^\gamma\}_{\gamma\in\Gamma}$, and the value function $u$ will then solve the fully nonlinear equation 
\begin{equation}
\label{eq:Dir_pb_fully}
\left\{
\begin{array}{rcll}
\I u & = & 0&\quad\text{in}\quad \Omega,\\
u & = & g&\quad \text{in}\quad \R^n\setminus \Omega,
\end{array}
\right.
\end{equation}
where 
\[\mathcal I w := \inf_{\gamma\in \Gamma} \big\{ -\L_\gamma w\big\}.\]

Similar considerations with zero-sum two-player stochastic games\footnote{Namely, two players with conflicting interests govern the evolution of the particle by choosing from different sets of indices.} lead to the nonlinear operators 
\[
\I w := \inf_{b\in \B}\,\sup_{a\in \A}\, \big\{- \L_{ab} w\big\}.
\]
Finally, when running costs are considered we get the problem \eqref{eq:Dir_pb_fully} with 
\[
\begin{array}{c}
\I u = \inf_{\gamma\in \Gamma}\big\{-\L_\gamma u + c_\gamma\big\}\\[0.1cm]
\text{\small (stochastic control)}
\end{array}
\quad\begin{array}{c}
\text{or}\\
\phantom{1}
\end{array}\quad
\begin{array}{c}
\I u = \inf_{b\in \B}\sup_{a\in\A}\big\{-\L_{ab} u + c_{ab}\big\}\\[0.1cm]
\text{\small (zero-sum games)}
\end{array}.
\]

As explained in more detail below, here we will be interested in the case where $\L_\gamma$ are {integro-differential} operators of the form
\[
\L_\gamma u(x) = {\rm P.V.} \int_{\R^n}\big(u(x)-u(x+y)\big)K_\gamma(y)\,dy,
\]
with $K_\gamma(y)\asymp |y|^{-n-2s}$.

\begin{rem}[Second-order equations]
When $\L_\gamma u$ or $\L_{ab}$ are second-order uniformly elliptic operators of the form
\[
\L u (x) = -\sum_{i,j = 1}^n a_{ij}(x) \partial_{ij} u(x),\qquad\lambda{\rm Id} \le (a_{ij})_{ij} \le \Lambda{\rm Id},
\]
then $\I u$ is a (local) fully nonlinear uniformly elliptic operator $F(D^2 u) = 0$. 
When $\I u = \inf_{\gamma\in \Gamma}\left\{-\L_\gamma u + c_\gamma\right\}$ then $F$ is, in addition, a concave function. 

Notice that, in this case,  the stochastic processes $X_t^\gamma$ have continuous paths, and thus only the value of $g$ on $\partial\Omega$ matters in \eqref{eq:Dir_pb_fully}. 
\end{rem}

\subsection{Extremal operators}

\index{Extremal operators}

An important observation is the following:
\begin{lem}
\label{lem:32}
Assume that all $\L_\gamma$ and $\L_{ab}$ belong to a class of linear operators~$\LL$; and $c_\gamma$ and $c_{ab}$ are functions of $x$.  Let 
\begin{equation}
\label{eq:Iu_form}
\begin{array}{c}
\I (u, x) = \inf_{\gamma\in \Gamma}\big\{-\L_\gamma u(x) + c_\gamma(x)\big\}
\\[0.2cm]
\text{or}\\[0.2cm]
\I (u, x)= \inf_{b\in \B}\sup_{a\in \A}\big\{-\L_{ab} u(x) + c_{ab}(x)\big\}.
\end{array}
\end{equation}
Then, 
\[
\inf_{\L \in \LL} \, \big\{-\L v (x) \big\} \le \I (u + v, x) - \I (u, x) \le \sup_{\L\in \LL}\,\big\{- \L v(x)\big\},
\]
whenever all the terms are pointwise well-defined. 
\end{lem}
\begin{proof}
We have
\[
\begin{split}
 \inf_{\L\in \LL} \left\{-\L v (x)\right\}
 & \le \inf_{\gamma\in \Gamma}\left\{-\L_\gamma v(x)\right\}\\
& \hspace{-1cm}\le \inf_{\gamma\in \Gamma}\left\{-\L_\gamma u(x) - \L_\gamma v(x) + c_\gamma(x) \right\} - \inf_{\gamma\in \Gamma}\left\{-\L_\gamma u (x) + c_\gamma (x)\right\}\\
& \hspace{-1cm}\le \sup_{\gamma\in \Gamma} \left\{-\L_\gamma v (x)\right\}
 \le \sup_{\L\in \LL}\left\{-\L v(x)\right\}. 
\end{split}
\]
In the case $\L_{ab}$ we have instead 
\[
-\L_{ab} (u + v) +c_{ab}  \le -\L_{ab} u + c_{ab} + \sup_{\L\in \LL} \left\{-\L v\right\},
\]
and taking $\inf_{b\in \B}\sup_{a\in \A}$ we get 
\[
\I (u + v, x) \le \I (u, x) + \sup_{\L \in \LL} \{-\L v(x)\}. 
\]
The other inequality follows in a similar way. 
\end{proof}

Given a class $\LL$ of linear L\'evy  operators \eqref{eq:L_def}-\eqref{eq:cond_nu}, we define the \emph{extremal operators}:
\begin{equation}
\label{eq:extremal_def}
\M^+_{\LL} u = \sup_{\L \in \LL}\big\{- \L u \big\}\qquad\text{and}\qquad 
\M^-_{\LL} u = \inf_{\L\in \LL} \big\{-\L u\big\}. 
\end{equation}
Thanks to the previous lemma, if $\I u $ is of the form \eqref{eq:Iu_form}, then 
\begin{equation}
\label{eq:extremal_eq} \index{Ellipticity!Viscosity solutions}
\M_\LL^-\left(u-v\right) \le \I u - \I v \le \M_\LL^+\left(u-v\right).
\end{equation}
We say that $\I$ is \emph{elliptic with respect to the class} $\LL$ when \eqref{eq:extremal_eq} holds,   for any $u, v\in C^2_b(\R^n)$. 

Notice that $\M^+_{\LL}$ and $\M^-_{\LL}$ are themselves fully nonlinear operators, and that they are elliptic with respect to $\LL$ (by Lemma~\ref{lem:32}). 

\begin{rem}
When $\LL$ is the class of second order uniformly elliptic operators with ellipticity constants $\lambda$ and $\Lambda$, the previous operators are called \emph{Pucci operators}, \cite{CC, FR4}:
\[
\begin{split}
\pm \M^\pm_{\Lambda, \lambda} (D^2 u ) & = \sup_{\lambda{\rm Id}\le A\le \Lambda{\rm Id}} \big(\pm{\rm tr}(A D^2 u)\big).
\end{split}
\]
\end{rem}

\begin{rem}
 It is interesting to notice that, when $\LL$ is the class of L\'evy-type operators \eqref{eq:levytypeoperator},   any  operator $\I$ satisfying \eqref{eq:extremal_eq} (or more generally, \eqref{eq:GMP}) must be of the form \eqref{eq:Iu_form}, with $c_{ab}$ constant; see \cite{GS2, GS3}. See \cite{GS2, GS3} for more general results for operators with $x$ dependence. 
\end{rem}

\begin{rem}
 In \eqref{eq:extremal_def}, we take $-\L$ instead of $\L$ because traditionally this is the sign taken in the context of non-divergence-form equations.
 Notice that the linear operators $\L$ have the same sign as $-\Delta$, while the fully nonlinear operators $\I$ have the same sign as $\Delta$.
\end{rem}

\begin{rem}
In \eqref{eq:Iu_form}, or in Definitions~\ref{defi:LL} and \ref{defi:II} below, we could have chosen the class of operators $\I$ to be translation invariant (that is, with $c_{ab}$ constants), since we are already taking the linear operators $\L_{ab}$ to be independent of $x$. However, as we will see, allowing the lower-order terms $c_{ab}$ to depend on $x$ is convenient when studying the interior regularity of solutions to $\I u = 0$ (even in the case where $\I$ is translation invariant). 
\end{rem}

\subsection{The class of operators}

In this chapter,   we focus our attention on the class of linear operators $\L$ with kernels comparable to the one of fractional Laplacian (see \eqref{eq:abs_cont_K}-\eqref{eq:Kcompfls} in subsection~\ref{ssec:comparable}). 

Namely, we consider operators of the form 
\begin{equation}
\label{eq:Lu1}
\begin{split}
\L u(x) & = {\rm P.V.} \int_{\R^n}\big(u(x)-u(x+y)\big)K(y) \, dy \\
& = \frac12 \int_{\R^n}\big(2u(x)-u(x+y)-u(x-y)\big)K(y)\,dy,
\end{split}
\end{equation}
with
\begin{equation}
\label{eq:Kint1}
K(y) = K(-y)\qquad\text{in}\quad \R^n,
\end{equation}
and
\begin{equation}
\label{eq:compdef}
0< \frac{\lambda}{|y|^{n+2s}}\le K(y) \le \frac{\Lambda}{|y|^{n+2s}}\qquad\text{in}\quad   \R^n,
\end{equation}
and we define the classes $\LL_s(\lambda, \Lambda)$ and $\LL_s(\lambda, \Lambda;\mu)$ for  $\mu > 0$ as follows (cf.~Definitions~\ref{defi:G} and \ref{defi:Gmu}):
\begin{defi}
\label{defi:LL}
Let $s\in (0, 1)$, $0< \lambda \le \Lambda$, and $\mu > 0$ with $\mu\notin\N$. We define
\[
\LL_s(\lambda, \Lambda) :=\big\{\L : 
\text{\eqref{eq:Lu1}-\eqref{eq:Kint1}-\eqref{eq:compdef} holds}
 \big\},
\]
and 
\[
\LL_s(\lambda, \Lambda; \mu) :=\big\{\L \in \LLL : [\L]_{C^\mu} < \infty\big\},
\]
where we have denoted, for $\L\in \LLL$ with kernel $K$, 
\[
[\L]_{C^\mu} := \sup_{\rho > 0}\,\rho^{n+2s+\mu} [K]_{C^\mu(B_\rho^c)}.
\]
 We also denote $\LL_s(\lambda, \Lambda; 0) := \LLL$.
\end{defi}
\begin{rem}
\label{rem:equiv_seminorms}
Arguing as in \eqref{eq:prop_kernel} or \eqref{eq:prop_kernel4} one can see that the norm $[\L]_{C^\mu}$ is equivalent to considering instead semi-norms in $B_{2\rho}\setminus B_{\rho}$:
\[
[\L]_{C^\mu} \asymp \sup_{\rho > 0}\,\rho^{n+2s+\mu} [K]_{C^\mu(B_{2\rho}\setminus B_\rho)}.
\]
\end{rem}

Notice that the class $\LLL$ is \emph{scale invariant}. 
That is, for any $\L \in \LLL$ there is $\L_r\in \LLL$ such that
\begin{equation}
\label{eq:scaleinvariance_comp} \index{Scale invariance!Operators comparable to fractional Laplacian}
\big(\L_r u(r\, \cdot\,)\big)(x) = r^{2s}(\L u)(rx)
\end{equation}
for every $u$ such that $\L u(rx)$ is well-defined. 
More precisely, if $\L$ has kernel~$K$, then the kernel $K_r$ of $\L_r$ is given by $K_r(y) := r^{n+2s} K(ry)$. 

We then define the classes of fully nonlinear integro-differential operators $\II_s(\lambda, \Lambda)$ and $\II_s(\lambda, \Lambda; \mu)$ as follows, in which we consider operators of the form:
\begin{equation}
\label{eq:Iofheform00}
\I (u, x) = \inf_{b\in \B}\sup_{a\in \A}\big\{-\L_{ab} u(x) + c_{ab}(x)\big\}\quad\text{with}\quad \L_{ab}\in \LLL.
\end{equation}

\begin{defi}
\label{defi:II}
Let $s\in (0, 1)$ and $0< \lambda \le \Lambda$. We define
\[
\II_s(\lambda, \Lambda) :=\left\{\I : \begin{array}{l}
\I \text{ is of the form \eqref{eq:Iofheform00}}, ~~\I(0, \,\cdot\,)\in L^\infty(\R^n),\\
\text{and}~(c_{ab}(x))_{ab}~\text{are equicontinuous in $\R^n$}
\end{array} \right\},
\]
that is, there exists some $\sigma(\I):[0, +\infty)\to [0, +\infty)$ continuous, nondecreasing, with $\sigma(\I)(0) = 0$, such that $|c_{ab}(x)-c_{ab}(y)|\le \sigma(\I)(|x-y|)$ for all $x, y\in \R^n$, and $(a, b)\in \A\times\B$. 

We also define, given $\mu > 0$,
\[
\II_s(\lambda, \Lambda; \mu) :=\big\{\I \in \IIL: [\I]_{C^\mu}<+\infty 
 \big\}.
\]
where for $\I \in \III$ we denote
\[
[\I]_{C^\mu} := \sup_{(a,b)\in \A\times\B} [\L_{ab}]_{C^\mu}. 
\]
In particular, in the expression \eqref{eq:Iofheform00} we have that $\L_{ab}\in \LL_s(\lambda, \Lambda; \mu)$ for all $(a, b)\in \A\times \B$. When $\mu = 0$, we denote furthermore $\II_s(\lambda, \Lambda; 0) := \III$.
\end{defi}

 The  extremal operators corresponding to the class $\LLL$ have a relatively simple closed expression: 
\begin{equation}
\label{eq:Mpexplicit}
\begin{split}
\mathcal{M}_{\LL_s(\lambda, \Lambda)}^+u(x) = &\frac12  \int_{\R^n}\bigg\{  \Lambda\big(u(x+y)+u(x-y)-2u(x)\big)_+ \\
&\qquad\quad  - \lambda\big(u(x+y)+u(x-y)-2u(x)\big)_-\bigg\}\frac{dy}{|y|^{n+2s}},
\end{split}
\end{equation}
and 
\begin{equation}
\label{eq:Mpexplicit2}
\begin{split}
\mathcal{M}_{\LL_s(\lambda, \Lambda)}^-u(x) = &\frac12  \int_{\R^n}\bigg\{  \lambda\big(u(x+y)+u(x-y)-2u(x)\big)_+ \\
& \qquad - \Lambda\big(u(x+y)+u(x-y)-2u(x)\big)_-\bigg\}\frac{dy}{|y|^{n+2s}}.
\end{split}
\end{equation}

Throughout this chapter we will denote 
\begin{equation}
\label{eq:MMpm}
\Mp := \MMp \qquad\text{and}\qquad \Mm := \MMm.
\end{equation}
From the definitions, we immediately have the following properties:
\begin{enumerate}[(a)]
\item \index{Translation invariance!Extremal operators} $\Mpm$ are \emph{translation invariant}, i.e., 
\[\Mpm \big(a u(\,\cdot\, + b)\big) (x) = a (\Mpm u)(x+b)\]
for any $a\ge 0$ and $b\in \R^n$, 
\item \index{Rotation invariance!Extremal operators} $\Mpm$ are \emph{rotation invariant}, i.e.,  $\Mpm \big(u(O \,\cdot\,)\big) (x) = (\Mpm u)(O x)$ for any orthogonal transformation $O\in \mathcal{O}(n)$, 
\item \index{Scale invariance!Extremal operators} $\Mpm$ are \emph{scale invariant} of order $2s$, i.e., 
\[\Mpm \big(u(r\,\cdot\,)\big)(x) = |r|^{2s}(\Mpm u)(rx)\]
for any $r\in \R$, 
\item $\Mp (-u) = -\Mm u$,
\item $\Mp(u+v) \le \Mp u + \Mp v$, 
\item $\Mm(u+v)\ge \Mm u + \Mm v$,
\end{enumerate}
whenever these expressions are well-defined (in particular, also for viscosity solutions, see Section~\ref{sec:viscosity}). Thus, the extremal operators are translation, rotation, and scale invariant. 

\subsection{Strong solutions} 

\index{Strong solutions!Fully nonlinear equations}

Let us start with some basic properties of the classes $\LLL$ and $\III$.

In the following lemma, we say that $u\in C^{2s+\eps}_r(x_\circ)$ for some $r, \eps > 0$ with $2s+\eps< 2$ and $2s+\eps \neq 1$ if
\begin{equation}
\label{eq:locestCr}
[u]_{C^{2s+\eps}_r(x_\circ)} := \sup_{h\in B_r} \frac{|u(x_\circ + h)+u(x_\circ-h)-2u(x_\circ)|}{|h|^{2s+\eps}}\ <\infty.
\end{equation}
If  $2s+\eps = 2$, we say that  $u\in C_r^{1,1}(x_\circ)$.

\begin{lem}
\label{lem:Lu_LL}
Let $s\in (0, 1)$, let $\eps>0$ such that $2s+\eps\le 2$, and let $\L\in \LL_s(\lambda, \Lambda)$. 
Then, for any $u \in C^{2s+\eps}_r(0)\cap L^1_{\omega_s}(\R^n)$ for some $r \in (0, 1)$,  we have that $\L u (0)$ is well-defined and 
\[
|\L u (0)|\le C\Lambda \left( r^{\eps} [u]_{C^{2s+\eps}_r(0)}+ r^{-2s} |u(0)| + r^{-n-2s}\|u\|_{L^1_{\omega_s}(\R^n)}\right)
\]
for some $C$ that depends only on $n$, $s$, and $\eps$. 

Moreover, for any $u \in C^{2s+\eps}(B_1)\cap L^1_{\omega_s}(\R^n)$, we have $\L u \in C(B_1)$ with a modulus of continuity in $B_{1/2}$ that depends only on $u$, $n$, $s$, and $\Lambda$ (in particular, it is independent of $\L$). 
\end{lem}
\begin{proof}
The first part follows as in the proof of Lemmas~\ref{lem:laplu} or \ref{lem:Lu}, by bounding:
\[
\begin{split}
\big|\L u (0)\big| & \le \Lambda [u]_{C^{2s+\eps}_r(0)}\int_{B_r}\frac{dy}{|y|^{n-\eps}}+2\Lambda \int_{B_r^c} \frac{|u(0)|}{|y|^{n+2s}}\,dy + 2\Lambda \int_{B_r^c} \frac{|u(y)|}{|y|^{n+2s}}\, dy\\
& \le C\Lambda \left( r^{\eps} [u]_{C^{2s+\eps}_r(0)}+r^{-2s} |u(0)| +  r^{-n-2s} \|u\|_{L^1_{\omega_s}(\R^n)}\right).
\end{split}
\]

For the second part, we proceed as in the proof of Lemma~\ref{lem:Lu_2} by considering a cut-off function $\eta\in C^\infty_c(\R^n)$ such that $\eta \ge 0$, $\eta \equiv 0$ in $\R^n\setminus B_{3/4}$ and $\eta \equiv 1$ in $B_{2/3}$, and
\[
u_1 := \eta u\qquad\text{and}\qquad u_2 := (1-\eta) u,
\]
so that $u = u_1 + u_2$, with $u_1$ being compactly supported in $B_{3/4}$ and $u_2$ satisfying that $u_2 \equiv 0$ in $B_{2/3}$.

Then, exactly as in Lemma~\ref{lem:Lu_2} we have 
\begin{equation}
\label{eq:Lu1eps}
\|\L u_1\|_{C^\eps(B_{1/2})}\le C
\end{equation}
for some $C$ depending only on $n$, $s$, $u$, and $\Lambda$. On the other hand, proceeding similarly to \ref{it:step2Lu2} of the proof of Lemma~\ref{lem:Lu_2}, for $x_1, x_2\in B_{1/2}$ with $z_\circ := x_1-x_2$ and $|z_\circ|\le \frac{1}{32}$ we have 
\[
\big|\L u_2(x_1) - \L u_2 (x_2) \big|= 
  \left|\int_{B_{3/5}^c} \left(u_2(z_\circ+y) - u_2(y) \right) K(y-x_2)\, dy\right|,
\]
where we are using that $u_2$ vanishes in $B_{2/3}$. Now, from the upper bound in \eqref{eq:compdef} and since $|y-x_2|\ge c|y|$ in $B_{3/5}^c$ (because $x_2\in B_{1/2}$) we have 
\[
\big|\L u_2(x_1) - \L u_2 (x_2) \big|\le  
  C\int_{\R^n} \left| u_2(z_\circ+y) - u_2(y) \right| w(y) \, dy,
\]
where we have denoted $w(y) := (1+ |y|)^{-n-2s} $. By the triangle inequality, 
\[
\begin{split}
\big|\L u_2(x_1) - \L u_2 (x_2) \big|& \le  
  C\big\|(w u_2) (z_\circ+\,\cdot\,) - w u_2\big\|_{L^1(\R^n)}\\
  & \quad  + C\int_{\R^n} \left| w(z_\circ+y) - w(y) \right| u_2(z_\circ+y)\, dy,
  \end{split}
\]
for some $C$ depending only on $n$, $s$, $\Lambda$, and $u$. Observe that now both terms above go to zero as $|z_\circ|\downarrow 0$. Indeed, the first term is just the continuity of translations for $L^1$ functions (since $w u_2\in L^1(\R^n)$), while for the second term we notice that, since $w$ is globally Lipschitz, we have $\left| w(z_\circ+y) - w(y) \right|\le C\min\{|z_\circ|, w(z_\circ+y)\}$ and it also goes to zero as $|z_\circ|\downarrow 0$ by the dominated convergence theorem.

Hence, $\L u_2$ has a modulus of continuity that depends only on $n$, $s$, $\Lambda$, and $u$. Together with \eqref{eq:Lu1eps} this concludes the proof. 
\end{proof}
\begin{rem}
\label{rem:onesidedcondition}
In fact, the previous proof says something even more general. That is, if $u\in L^1_{\omega_s}(\R^n)$ is such that
\[
{u(x)+u(-x)-2u(0)}\le C_u |x|^{2s+\eps}\quad\text{in}\quad B_r,
\]
then we can evaluate $\L u(0)$ pointwise, with a one-sided bound
\[
\L u (0)\ge -C\Lambda \left( r^{\eps} C_u+ r^{-2s} |u(0)| + r^{-n-2s}\|u\|_{L^1_{\omega_s}(\R^n)}\right) > -\infty.
\]
Notice, however, that under this assumption it could be that $\L u(0)$ takes the value $+\infty$. Similarly, if ${u(x)+u(-x)-2u(0)}\ge -C_u |x|^{2s+\eps}$ in $B_r$ instead, the value $\L u(0)$ is again well-defined, but it may be $-\infty$. 
\end{rem}

We also have a higher regularity result, as the one in Lemma~\ref{lem:Lu_2}:
\begin{lem}
\label{lem:Lu_2_v}
Let $s\in (0, 1)$, let $\L\in \LL_s(\lambda, \Lambda; \alpha)$ for some $\alpha>0$ with $\alpha\not\in\mathbb N$. 
Then, for any $u\in C^{2s+\alpha}(B_1)\cap L^1_{\omega_s}(\R^n)$, we have  $\L u\in C_{\rm loc}^\alpha(B_{1})$ with
\[
\|\L u\|_{C^\alpha (B_{1/2})} \le C   \left(\|u\|_{C^{2s+\alpha}(B_1)} + \|u\|_{L^{1}_{\omega_s}(\R^n)}\right) ,
\]
for some $C$ depending only on $n$, $s$, $\Lambda$; $\alpha$, and $[\L]_{C^\alpha}$.
\end{lem}
\begin{proof}
The proof follows along the lines of Lemma~\ref{lem:Lu_2}. Using the notation therein, the only difference comes in the very last step, where we need to bound $\|\L u_2\|_{C^\alpha(B_{1/2})}$ in equation~\eqref{eq:weneedtoboundlaststep}. In this case, using that $[\L]_{C^\alpha}<  \infty$ (recall Definition~\ref{defi:LL}), we can bound it by 
\[
\left|D^k \L u_2(x_1) - D^k \L u_2 (x_2) \right|\le C [K]_{C^\alpha}  \|u\|_{L^1_{\omega_s}(\R^n)}  |x_1-x_2|^{\alpha-\lfloor \alpha\rfloor},
\]
and this is enough to conclude the proof.
\end{proof}

As a consequence of the previous lemmas, we also obtain a similar result for fully nonlinear operators:

\begin{cor}
\label{cor:Iu_II}
Let $s\in (0, 1)$, let $\eps>0$ such that $2s+\eps\le 2$, and let $\I\in \II_s(\lambda, \Lambda)$. 
Then, for any $u \in C^{2s+\eps}_r(0)\cap L^1_{\omega_s}(\R^n)$ for some $r \in (0, 1)$, we have that $\I( u, 0)$ is well-defined and 
\[
|\I (u,0)|\le C\Lambda \left( r^\eps [u]_{C^{2s+\eps}_r(0)}+ r^{-2s} |u(0)| + r^{-n-2s}\|u\|_{L^1_{\omega_s}(\R^n)}\right)+|\I(0, 0)|
\]
for some $C$ that depends only on $n$, $s$, and $\eps$.  

Moreover, if $u \in C^{2s+\eps}(B_1)\cap L^1_{\omega_s}(\R^n)$ we have
\[
\|\I (u, \cdot) \|_{L^\infty(B_{1/2})}\le C\Lambda \left( \|u\|_{C^{2s+\eps}(B_1)}+\|u\|_{L^1_{\omega_s}(\R^n)}\right) + \|\I(0, \cdot)\|_{L^\infty(B_{1/2})}
\]
and $\I u \in C(B_1)$ with a modulus of continuity in $B_{1/2}$ that depends only on $u$, $n$, $s$, $\Lambda$, and $\sigma(\I)$ (recall Definition~\ref{defi:II}).
\end{cor}

\begin{proof}
It follows from Lemma~\ref{lem:Lu_LL}, the definition of $\I$,  and the fact that the infimum and supremum of equicontinuous functions has the same modulus of continuity. 
\end{proof}
\begin{rem}
\label{rem:onesidedcondition2}
As in Remark~\ref{rem:onesidedcondition}, it is still possible to have a well-defined value (including $\pm\infty$) for $\I(u, 0)$ with a one-sided pointwise condition for $u$ at 0.  
\end{rem}

And the higher regularity: 

\begin{cor}
\label{cor:Iu_II_2}
Let $s\in (0, 1)$, let $\I\in \II_s(\lambda, \Lambda; \alpha)$ for some $\alpha\in (0, 1)$ (see Definition~\ref{defi:II}) be of the form
\[
\I (u, x) = \inf_{b\in \B}\sup_{a\in \A}\big\{-\L_{ab} u(x) + c_{ab}(x)\big\},\qquad \L_{ab}\in \LL_s(\lambda, \Lambda; \alpha),
\]
and satisfying 
\[
\sup_{(a, b)\in \A\times\B} [c_{ab}]_{C^\alpha(\R^n)}\le C_\circ.
\]
Then, for any $u\in C^{2s+\alpha}(B_1)\cap L^1_{\omega_s}(\R^n)$, we have  $\I u\in C_{\rm loc}^\alpha(B_{1})$ and
\[
\|\I (u, \cdot)\|_{C^\alpha (B_{1/2})} \le C   \left(\|u\|_{C^{2s+\alpha}(B_1)} + \|u\|_{L^{1}_{\omega_s}(\R^n)}\right)+ \|\I(0, \cdot)\|_{L^\infty(B_{1/2})}+C_\circ,
\]
with $C$ depending only on $n$, $s$, $\Lambda$, $\alpha$, and $[\I]_{C^\alpha}$.
\end{cor}
\begin{proof}
It follows from Lemma~\ref{lem:Lu_2_v}, Corollary~\ref{cor:Iu_II}, and again from the definition of $\I$  since the infimum and supremum of equicontinuous functions has the same modulus of continuity.
\end{proof}
\begin{rem}
Here, and in Lemma~\ref{lem:Lu_2_v} above, we also have the analogue of Remark~\ref{rem:Lu_2} in this context. 
\end{rem}

We also obtain that the ellipticity condition from Lemma~\ref{lem:32} holds as soon as $u$ and $v$ are $C^{2s+\eps}$ at a point:

\begin{cor}
\label{cor:elliptC11}
Let $s\in (0, 1)$,  let $\I\in \II_s(\lambda, \Lambda)$, and let $\Mpm$ be given by \eqref{eq:MMpm}. Then, for any $u, v \in C^{2s+\eps}_r(x_\circ)\cap L^1_{\omega_s}(\R^n)$ for some $\eps, r>0$, we have 
\[
\Mm(u-v)(x_\circ) \le \I(u, x_\circ) - \I(v, x_\circ) \le \Mp(u-v)(x_\circ).
\]
\end{cor}
\begin{proof}
The proof follows as in Lemma~\ref{lem:32}, where now all the terms are well-defined thanks to Lemma~\ref{lem:Lu_LL} and Corollary~\ref{cor:Iu_II}.
\end{proof}

\section{Viscosity solutions}

\index{Viscosity solutions}\label{sec:viscosity}

In this section we turn our attention to the existence and uniqueness of solutions. In the previous chapter, in subsection~\ref{ssec:existence_weak} we showed the existence of solutions to (linear) integro-differential elliptic equations in a variational way (i.e., by minimizing an energy functional), and in order to do that we introduced the notion of \emph{weak solution}. In this chapter, however, the same method does not apply: in general, fully nonlinear equations do not have a variational formulation (in particular, they do not come from an energy functional). To construct solutions, therefore, we will need to rely on other methods, and in this case, it will be through a defining factor in elliptic problems: the comparison principle. We will use Perron's method to prove existence of solutions, and in order to do that we need a new generalized notion of solution: \emph{viscosity solutions}. 

Notice that if the new notion were too restrictive, then we might be able to prove   a comparison principle (and uniqueness of solutions) but it might be difficult to prove existence. Otherwise, if we relaxed the notion of solution too much, we might be able to prove existence, but within such a general class that uniqueness might not be possible. Hence, one has to find the right balance. In this case, such a balance is found with the notion of viscosity solution in Definition~\ref{defi:viscosity} below, which allows the proofs of existence, uniqueness, and stability of solutions. Originally introduced by Crandall and Lions in 1983 in the study of first-order equations, \cite{CL83}, in 1988 Jensen showed that the concept is also well-posed for second order elliptic equations, \cite{Jen88}, and since then it has become prevalent in the analysis of elliptic problems. In the context of integro-differential equations, the basic theory of viscosity solutions was developed by Caffarelli and Silvestre in \cite{CS}.

This section is based partly on \cite{CS, CS3, S-viscosity, FR4}; see also \cite{Mou} and \cite{BCI,BCI2,BI}.

\subsection{Definition and basic properties} 
In the following, we recall the notion of semi-continuity:

\begin{defi} \index{Semi-continuity}
A function $u$ is lower semi-continuous at a point $x_\circ\in \R^n$ if $u(x_\circ) > -\infty$ and 
\[
\liminf_{x\to x_\circ} u(x) \ge u(x_\circ).
\]
Similarly, it is upper semi-continuous at $x_\circ \in \R^n$ if $u (x_\circ) < \infty$ and 
\[
\limsup_{x\to x_\circ } u(x) \le u(x_\circ).
\]
Of course, $u$ is continuous at $x_\circ$ if and only if it is both upper and lower semi-continuous.

Given any set $D\subset \R^n$, we say that $u$ is lower (resp. upper) semi-continuous in $D$, and we denote it  $u\in {\rm LSC}(D)$ (resp. $u\in {\rm USC}(D)$) if for any $x_\circ\in D$,
\[
\liminf_{D\ni x\to x_\circ} u(x) \ge u(x_\circ)>-\infty,\quad\left(\text{resp.}\quad \limsup_{D\ni x\to x_\circ}  u(x) \le u(x_\circ)<\infty\right).
\]
In particular, if $D$ is open, then this is equivalent to $u$ being lower (resp. upper) semi-continuous at every $x_\circ\in D$. Any lower (resp. upper) semi-continuous function  defined on a compact set always attains its minimum (resp. maximum). 
\end{defi}

We can now give the definition of viscosity sub- or supersolution to an equation of the type
\begin{equation}
\label{eq:defivisc}
\I (u, x) = f(x)
\end{equation}
 by evaluating the corresponding operator on smooth functions touching from above or below: 
\begin{defi}[Viscosity solutions] 
\label{defi:viscosity}
Let $s\in (0, 1)$, let $\I\in \II_s(\lambda, \Lambda)$, let $\Omega\subset \R^n$ be any open set, and let $f\in C(\Omega)$.
\begin{itemize}[leftmargin=1cm]
\item We say that $u\in {\rm USC}({\Omega})\cap L^1_{\omega_s}(\R^n)$ is a \emph{viscosity subsolution} to \eqref{eq:defivisc} in $\Omega$, and we denote   
\[\I (u, x) \ge f(x)\quad  \textrm{in}\quad \Omega,\] 
if for any $x\in \Omega$ and any neighborhood of $x$ in $\Omega$, $N_x\subset \Omega$, and for any test function $\phi\in L^1_{\omega_s}(\R^n)$ such that $\phi\in C^2(N_x)$, $\phi(x) = u(x)$, and $\phi \ge u$ in all of $\R^n$,   we have $\I (\phi,x) \ge f(x)$. 

\item We say that $u\in {\rm LSC}(\Omega)\cap L^1_{\omega_s}(\R^n)$ is a \emph{viscosity supersolution} to \eqref{eq:defivisc} in $\Omega$, and we denote   
\[\I (u, x) \le f(x)\quad  \textrm{in}\quad \Omega,\]
if for any $x\in \Omega$ and any neighborhood of $x$ in $\Omega$, $N_x\subset \Omega$, and for   any   test function $\phi\in L^1_{\omega_s}(\R^n)$ such that $\phi\in C^2(N_x)$, $\phi(x) = u(x)$, and $\phi \le u$ in all of $\R^n$,   we have $\I (\phi,x) \le f(x)$. 
\item We say that $u\in C({\Omega})\cap L^1_{\omega_s}(\R^n)$ is a \emph{viscosity solution} to \eqref{eq:defivisc} in $\Omega$, and we denote $\I(u, x) = f(x)$ in $\Omega$, if it is both a viscosity subsolution and supersolution. 
\end{itemize}
\end{defi}

Observe that the notion of viscosity solution only requires $u$ to be continuous in $ {\Omega}$. In particular, there might be points $x\in \Omega$ at which there is no function $\phi\in C^2$  touching $u$ at $x_\circ$ (from above and/or below); this is allowed by the previous definition. 

\begin{rem}
\label{rem:max_min}
Notice that, if  $u_1$ and $u_2$ are viscosity subsolutions in $\Omega$, say $\I(u_i, x)\ge f(x)$, then $\bar u = \max\{u_1, u_2\}$ is a viscosity subsolution in $\Omega$ as well, $\I(\bar u, x) \ge f(x)$. Indeed, given any function touching $\bar u$ from above, then it is touching either $u_1$ or $u_2$ from above at the same point. Similarly, if $v_1$ and $v_2$ are viscosity supersolutions in $\Omega$, then $\ubar{v} := \min\{v_1, v_2\}$ is a viscosity supersolution in $\Omega$ too.
\end{rem}

Let us start by showing that viscosity solutions are equivalent to strong solutions when the functions are regular enough:

\begin{lem}
\label{lem:equiv_defi} \index{Equivalent notions!Strong and viscosity}
Let $s\in (0, 1)$,  let $\I\in \II_s(\lambda, \Lambda)$, and let $\Omega\subset \R^n$ be any open set. Let $u \in C_{\rm loc}^{2s+\eps}(\Omega)\cap L^1_{\omega_s}(\R^n)$ for some $\eps > 0$. Then, $u$ satisfies $\I (u, x) = f(x)$ in $\Omega$ for some $f\in C(\Omega)$ in the strong sense if and only if it satisfies it in the viscosity sense. 
\end{lem}

\begin{proof} 
Let us first show that if $u$ is a strong solution, it is a viscosity solution. Observe that $\I(u, x) = f(x)$ is pointwise well-defined in the strong sense and $f\in C(\Omega)$, by Corollary~\ref{cor:Iu_II}.  Take $\phi$ to be a test function as in the definition of subsolution in Definition~\ref{defi:viscosity}:  $\phi\ge u$ in $\R^n$, $\phi\in C^2$ around $x$, and $\phi(x) = u(x)$. Then $w:= \phi-u\in C^{2s+\eps}_r(x)$ for some $r > 0$ (recall~\eqref{eq:locestCr}), and for any $\L \in \LL_s(\lambda, \Lambda)$ we can compute (using $w(x) = 0$)
\[
\L w (x) = -\int_{\R^n} \big(w(x+y)+w(x-y)\big)|y|^{-n-2s}\,dy\le 0,
\]
which is well-defined. In particular, by Lemma~\ref{lem:Lu_LL}, for any $\L \in \LL_s(\lambda, \Lambda$), 
\[
\L \phi(x) \le \L u(x) <+\infty,
\]
and therefore, $\I(\phi,x) \le \I (u,x) = f(x)$. Since $\phi$ was arbitrary, $u$ is   a viscosity subsolution. We can similarly check that $u$ is also a viscosity supersolution, and thus it is a viscosity solution to the equation $\I( u, x) = f(x)$ in $\Omega$. 

Conversely, let us suppose that $u$ is a viscosity solution. Notice that, by Corollary~\ref{cor:Iu_II}, we already know that $\I (u, x) = \tilde f(x)$, for some $\tilde f \in C(\Omega)$. We want to show $\tilde f = f$. 

Let now $x_\circ \in \Omega$ be arbitrary. After a translation, rescaling, and addition of a constant, we assume 
\[
x_\circ = 0, \quad B_1\subset \Omega, \quad u(0) = 0. 
\]

Moreover, if $2s+\eps > 1$ we can also subtract a hyperplane and assume furthermore  $\nabla u(0) = 0$.
In all cases we have, since $u\in C^{2s+\eps}(\Omega)$,
\begin{equation}
\label{eq:ubounded2se}
|u(x)| \le C|x|^{2s+\eps}\quad \text{in}\quad B_{1/2}. 
\end{equation}

Let $\delta > 0$, and let $c_\delta\ge0$ the smallest value for which 
\[
\frac{1}{\delta}|x|^2+c_\delta \ge u\quad\text{in}\quad B_{1/2}. 
\]
In particular, there exists some $x_\delta\in \overline{B_{1/2}}$ such that $\frac{1}{\delta} |x_\delta|^2 +c_\delta= u(x_\delta)$. From \eqref{eq:ubounded2se}
\[
\frac{1}{\delta}|x_\delta|^2+c_\delta \le C|x_\delta|^{2s+\eps}\quad\Rightarrow\quad |x_\delta|^{2-2s-\eps} \le C\delta.
\]

Let us denote $r_\delta := C \delta^{\frac{1}{2-2s-\eps}} = |x_\delta|$, and let us consider the test function
\[
\phi_\delta(x) = \left\{
\begin{array}{ll}
\frac{1}{\delta}|x|^2 + c_\delta& \quad\text{if}\quad |x|<2r_\delta,\\
u(x) & \quad\text{if}\quad |x|\ge 2r_\delta,
\end{array}
\right.
\]
which satisfies $\phi_\delta \ge u$ in $\R^n$, $\phi_\delta(x_\delta) = u(x_\delta)$, and $\phi_\delta$ is smooth around 0. Thus, on the one hand we know (since $u$ is a viscosity solution to $\I (u, x) = f(x)$) that 
\[
\I (\phi_\delta,x_\delta) \ge f(x_\delta). 
\]
And on the other hand, thanks to Lemma~\ref{lem:Lu_LL} and Corollary~\ref{cor:elliptC11},  denoting $v_\delta := \phi_\delta-u\ge 0$ and $\Mp$ given by \eqref{eq:MMpm},
\[
\begin{split}
\I(\phi_\delta, x_\delta)  -\I (u,x_\delta) & \le \Mp v_\delta(x_\delta)\\
& \le C \left( r_\delta^{\eps}[v_\delta]_{C_{2r_\delta}^{2s+\eps}(x_\delta)} + r_\delta^{-n-2s} \|v_\delta\|_{L^1(B_{2r_\delta})}\right).
\end{split}
\]
Observe that, from \eqref{eq:ubounded2se} and by definition of $\phi$, 
\[
\begin{split}
0\le v_\delta(x)& \le 
 C\left(|x|^{2s+\eps}+\frac{1}{\delta}|x|^2\right)\chi_{|x|\le 2r_\delta}\le C r_\delta^{2s+\eps},
\end{split}
\]
and in particular, $\|v_\delta\|_{L^1(B_{2r_\delta})} \le r_\delta^{n+2s+\eps}$.

Moreover, 
\[
[v_\delta]_{C_{2r_\delta}^{2s+\eps}(x_\delta)} \le [u]_{C_{2r_\delta}^{2s+\eps}(x_\delta)}  + \frac{1}{\delta}[|x|^2]_{C_{2r_\delta}^{2s+\eps}(x_\delta)}.
\]
Now, we know that $[u]_{C_{2r_\delta}^{2s+\eps}(x_\delta)}  < \infty$ by assumption on $u$ and
\[
[|x|^2]_{C_{2r_\delta}^{2s+\eps}(x_\delta)} =  r_\delta^{2-2s-\eps} [|x|^2]_{C_{2}^{2s+\eps} (x_\delta/r_\delta)} \le C r_\delta^{2-2s-\eps}
\]
Thus, using that $\delta = cr_\delta^{2-2s-\eps}$, we obtain
\[
[v_\delta]_{C_{2r_\delta}^{2s+\eps/2}(x_\delta)} \le C.
\]
In all, we have 
\[
\I(\phi_\delta, x_\delta) - \I (u,x_\delta)  \le C r_\delta^\eps,
\]
that is, 
\[
f(x_\delta) \le C r_\delta^\eps + \tilde f(x_\delta).
\]
Since $f$ and $\tilde f$ are continuous, we can let $\delta\downarrow 0$ (so that $r_\delta \downarrow 0$ and $x_\delta\to x_\circ$) to deduce
\[
f(x_\circ) \le \tilde f(x_\circ). 
\]

Repeating the same procedure with test functions from below, we obtain $f(x_\circ) \ge \tilde f(x_\circ)$. Hence, $\I (u, x) = f(x)$ in the strong sense, as  wanted. 
\end{proof}
\begin{rem}
\label{rem:visc_C2seps}
The proof of the previous lemma actually shows that, in the definition of viscosity sub- and supersolution, Definition~\ref{defi:viscosity}, we can equivalently consider functions $\phi$ that are $C^{2s+\eps}$ in a neighbourhood of $x$, instead of $C^2$. 
\end{rem}

As a consequence of Corollary~\ref{cor:Iu_II} we also obtain the following result:

\begin{lem}
\label{lem:C11punt1}
Let $s\in (0, 1)$,  let $\I\in \II_s(\lambda, \Lambda)$, and let $\Omega\subset \R^n$ open. Let $u\in {\rm LSC}(\Omega)\cap L^1_{\omega_s}(\R^n)$ such that $\I (u, x) \le 0$ in $\Omega$ in the viscosity sense. If $x_\circ \in \Omega$ and $\phi\in C^{1,1}_r(x_\circ)\cap L^1_{\omega_s}(\R^n)$ for some $r > 0$ is such that $\phi\le u$ in $\R^n$ and $\phi(x_\circ) = u(x_\circ)$, then $\I(\phi, x_\circ) \le 0$. Moreover, $\I(u, x_\circ) \le 0$ pointwise. 
\end{lem}
\begin{proof}
Notice that $\I(\phi, x_\circ)$ is defined classically thanks to Corollary~\ref{cor:Iu_II}. Moreover, since $\phi$ is $C^{1,1}$ at $x_\circ$, there exists a quadratic polynomial $p$ such that $\phi \ge p$ in $B_r(x_\circ)$, $\phi(x_\circ) = p(x_\circ)$, and $\nabla \phi(x_\circ) = \nabla p(x_\circ)$. In particular, for some $C > 0$ we have
\begin{equation}
\label{eq:quadratic_approx}
|\phi(x) - p(x)|\le C |x-x_\circ|^2\quad\text{in}\quad B_r(x_\circ).
\end{equation}

Let us define, for any $0 < \rho < r$,
\[
\phi_\rho(x) = \left\{
\begin{array}{ll}
p&\quad\text{in}\quad B_\rho(x_\circ)\\
\phi&\quad\text{in}\quad \R^n\setminus B_\rho(x_\circ).
\end{array}\right.
\]
Then $\phi_\rho$ is an admissible test function, and therefore $\I(\phi_\rho, x_\circ) \le 0$. On the other hand, by Corollary~\ref{cor:elliptC11}
\[
\I(\phi, x_\circ) - \I(\phi_\rho, x_\circ) \le \Mp((\phi-p)\chi_{B_\rho}).
\]
Finally, thanks to Corollary~\ref{cor:Iu_II}, since $[\phi-p]_{C^{1,1}_\rho(x_\circ)}< \infty$ and $\|\phi-p\|_{L^1_{\omega_s}(B_\rho)}\le \rho^{n+2}$ (by \eqref{eq:quadratic_approx}), we get (recall \eqref{eq:MMpm})
\[
\Mp\big((\phi-p)\chi_{B_\rho}\big)\le C \rho^{2-2s}\downarrow 0\quad\text{as}\quad \rho\downarrow 0. 
\]
That is, passing to the limit, $\I(\phi, x_\circ) \le 0$ as we wanted to see. 

Finally, observe that since $\phi$ touches $u$ at $x_\circ$ from below, we have
\[
u(x_\circ+x) + u(x_\circ-x) -2u(x_\circ) \ge \phi(x_\circ+x)+\phi(x_\circ-x)-2\phi(x_\circ) \ge -C_\phi|x|^{2s+\eps}
\] in some $B_r$ for $r > 0$, and therefore the value of $\I(u, x_\circ)$ is well-defined (though it could be $-\infty$), see Remarks~\ref{rem:onesidedcondition} and \ref{rem:onesidedcondition2}. In particular, since $\I(\phi, x_\circ)\le 0$, we have $\I(u, x_\circ) \le 0$.  
\end{proof}
\begin{rem}
\label{rem:C11}
Notice that the previous result implies that, in the definition of viscosity sub- and supersolution, Definition~\ref{defi:viscosity}, we can equivalently take functions $\phi$ that are pointwise $C^{1,1}$ at $x$, $\phi\in C^{1,1}_r(x)$ for some $r> 0$. Even more generally, as in the proof of Lemma~\ref{lem:equiv_defi}, we could equivalently take test functions that are pointwise $C^{2s+\eps}$ at $x$ (see Remark~\ref{rem:visc_C2seps} as well).
\end{rem}

\subsection{Stability} \index{Stability!Viscosity solutions}

One of the most important properties of any notion of solution is their potential stability under appropriate limits. For example, in Proposition~\ref{prop:stab_distr} we   saw that distributional solutions (to linear equations) are stable under $L^1$ limits. Here, we prove that viscosity solutions (to fully nonlinear equations) are stable under uniform limits and, more generally, that viscosity sub- and supersolution are stable under half-relaxed limits.

\begin{defi}[Half-relaxed limits] \label{defi:half_relaxed}  Let $\Omega\subset\R^n$ open, and let $(u_k)_{k\in \N}$ be a sequence such that $\inf_k u_k(x)$ is locally bounded from below in~$\Omega$. We say that $u$ is the \emph{(lower) half-relaxed limit} of $u_k$ in $\Omega$, and we denote it   $u_k\tosd u$ in $\Omega$, if 
\begin{equation}
\label{eq:halfrelaxedpointwise} \index{Half-relaxed limit}
u(x) = \liminfs{k\to\infty} u_k(x) := \inf\left\{\liminf_{k\to\infty} u_k(x_k) : x_k \to x\right\}
\end{equation}
for every $x\in \Omega$.

On the other hand, if $u_k$ is such that $\sup_k u_k(x)$ is locally bounded from above in $\Omega$, we say that $u$ is the \emph{(upper)  half-relaxed limit} of $u_k$, and we denote it $u_k\tosu u$ in~$\Omega$, if 
\[
u(x) = \limsups{k\to\infty} u_k(x) := \sup\left\{\limsup_{k\to\infty} u_k(x_k) : x_k \to x\right\}
\]
for every $x\in \Omega$.
\end{defi}

By definition, the lower half-relaxed limit is always \emph{lower semi-continuous} in $\Omega$, and the upper half-relaxed limit is always \emph{upper semi-continuous} in~$\Omega$. Notice, also, that if $u_k\tosu  u$ in $\Omega$, then $-u_k\tosd {-u}$ in $\Omega$.  Finally, if a sequence of continuous functions converges locally uniformly, it converges in the half-relaxed way described above. 

The notion of half-relaxed limits has the following two important properties:
\begin{lem}
\label{lem:min_K}
Let $\Omega\subset \R^n$ be open, $u_k \tosd u$ in $\Omega$, and let $K \ssubset \Omega$  compact. Then, for every $\eps > 0$ there exists some $k_\circ\in \N$ such that 
\[
u_k \ge \min_K u - \eps\quad\text{for all}\quad k\ge k_\circ. 
\]
\end{lem}
\begin{proof}
Since $u$ is lower semi-continuous, the minimum is achieved in $K$ compact. Now, arguing by contradiction, let us suppose that there are $k_j\in \N$ with $k_j \to \infty$ as $j \to \infty$, and $x_j\in K$ such that $u_{k_j}(x_j) < \min_K u - \eps_\circ$ for some $\eps_\circ$ fixed. Up to a subsequence, we have $x_j \to x_*\in K$, and hence 
\[
\min_K u - \eps_\circ\ge \liminf_{j\to \infty} u_{k_j}(x_j) \ge u(x_*),
\]
a contradiction. 
\end{proof}

We also have the following: 

\begin{lem}
\label{lem:loc_minima}
Let $\Omega\subset \R^n$ be open, and let $u_k \tosd u$ in $\Omega$ for some $u_k \in {\rm LSC}(\Omega)$. If $u$ has  a strict local minimum at $x_\circ\in \Omega$ (which is the minimum in $\overline{B_r(x_\circ)}\subset \Omega$ for some $r > 0$) then there exists a sequence of indices $k_j\to \infty$ and points $x_j \to x_\circ$ as $j\to \infty$ such that $u_{k_j}(x_j)\to u(x_\circ)$ and $u_{k_j}$ has a local minimum at $x_j$ (which is the minimum in the same ball $\overline{B_r(x_\circ)}$).
\end{lem}
\begin{proof}
Up to a translation and up to adding a constant, let us assume $x_\circ = 0$ and $u(0) = 0$ (notice that the local minimum of $u$ is finite, since it is lower semi-continuous and $u_k$ are uniformly bounded below). By assumption, there exists $r > 0$ such that $u(x) > 0$ for all $0 < |x|\le r$. In particular, for every $0 < \rho < r$ there exists some $\eps > 0$ such that 
\[
\min_{K_\rho} u \ge \eps,
\]
where we have denoted $K_\rho := \overline{B_r\setminus B_\rho}$. By Lemma~\ref{lem:min_K} applied with $K_\rho$, for $k$ large enough we have 
\[
u_k \ge \frac{\eps}{2}\quad\text{in}\quad K_\rho.
\]

Now, by definition of half-relaxed limit, there exist sequences of indices $k_j\to \infty$ and of points $y_j\to 0$ as $j\to \infty$ such that $u_{k_j}(y_j)\to 0$. Let $x_j\in \overline{B_r}$ be the point where the minimum of $u_{k_j}$ in $\overline{B_r}$ is attained (which exists because $u_k$ is lower semi-continuous). In particular, $u_{k_j}(x_j) \le u_{k_j}(y_j) \to 0$, that is, $u_{k_j}(x_j) \le \eps/4$ for $j$ large enough. Since $u_{k_j} \ge \eps/2$ in $K_\rho$ (again, for $j$ large enough), this implies that $x_j \in B_\rho$. That is, $u_{k_j}$ attains its minimum in $\overline{B_r}$, inside $B_\rho$. By repeating this argument choosing smaller $\rho> 0$, we can extract a subsequence $k_m:= k_{j_m}$ to get the desired result. Notice that $x_j \to 0$ and  $u_{k_j}(x_j) \to 0$ since $u_k\tosd u$. This completes the proof.
\end{proof}
 
Notice that by taking $-u_k$ and $-u$ in the previous lemma we obtain the corresponding result for upper semi-continuous functions (involving maxima and local maxima)  

On the other hand, we also need the notion of convergence of operators $\I(u, x)$: 
\begin{defi}[Weak convergence of operators]
\label{defi:conv_I}
Let $s\in (0, 1)$. Let $(\I_k)_{k\in \N}$ be a sequence of operators with $\I_k\in \III$ and let $\I\in \III$. We say that $\I_k$ weakly converges to $\I$ in $\Omega$, and we denote it 
\[
\I_k\rightharpoonup \I\quad\text{in}\quad \Omega,
\] if for every $x_\circ\in \Omega$ and every function $v\in L^1_{\omega_s}(\R^n)$ such that $v\in C^2(B_r(x_\circ))$ with $B_r(x_\circ)  \subset \Omega$, we have $\I_k(v, x) \to \I(v, x)$ uniformly in $B_{r/2}(x_\circ)$. 
\end{defi}

In all, we now have the ingredients to establish the stability result:
 
\begin{prop}[Stability of viscosity supersolutions]
\label{prop:stab_super}
Let $s\in (0, 1)$. Let $\Omega\subset \R^n$, and let assume that for every $k\in \N$ we have
\begin{enumerate}[leftmargin=*, label=(\roman*)]
\item  $\I_k, \I \in \III$, with $\I_k\rightharpoonup \I$ in $\Omega$  as $k\to \infty$,
\item   $u_k, u\in {\rm LSC}(\Omega)$, with $u_k\tosd u$ in $\Omega$  and $u_k\to u$ in $L^1_{\omega_s}(\R^n)$ as $k\to \infty$,
\item   $f_k, f\in C(\Omega)$, with $f_k\to f$ locally uniformly in $\Omega$ as $k\to \infty$,
\item $\I_k(u_k, x) \le f_k(x)$ in $\Omega$ in the viscosity sense.
\end{enumerate}
Then, $\I(u, x) \le f(x)$ in $\Omega$ in the viscosity sense. 
\end{prop}
\begin{proof}
Let us fix $x_\circ \in \Omega$, and let $\phi$ be any function touching $u$ from below at $x_\circ$ such that $\phi \le u$ in $\R^n$ and $\phi$ is $C^2$ in a neighborhood $B_r(x_\circ)\subset \Omega$ of~$x_\circ$, with $r > 0$. In the following, we may assume without loss of generality that $\phi = u$ in $\R^n\setminus \overline{B_r(x_\circ)}$. We want to show that $\I(\phi, x_\circ)\le f(x_\circ)$. Notice, first, that we can also assume that $u-\phi$ has a strict local minimum at $x_\circ$: otherwise take $\phi_\eps = \phi-\eps|x-x_\circ|^2$ in $B_r(x_\circ)$ instead (and $\phi_\eps = \phi$ in $B_r^c(x_\circ)$). Then $u-\phi_\eps$ has a strict local minimum at $x_\circ$, and if we showed the result for $\phi_\eps$ we would have $\I(\phi_\eps, x_\circ)\le f(x_\circ)$ for any $\eps > 0$. We can then use Lemma~\ref{lem:32} or Corollary~\ref{cor:elliptC11} (since $\phi$ and $\phi_\eps$ are $C^2$ around $x_\circ$): 
\[
\I(\phi, x_\circ) \le \I(\phi_\eps, x_\circ) + \eps \Mp\big(|\,\cdot\,-x_\circ|^2\chi_{B_r(x_\circ)}\big)(x_\circ)\le  f(x_\circ) + C \eps
\]
where $\Mp$ is given by \eqref{eq:MMpm}, and we have bounded the last term for some $C$ independent of $\eps$ by Lemma~\ref{lem:Lu_LL}. Letting  $\eps\downarrow 0$ we would have $\I(\phi, x_\circ) \le f(x_\circ)$ as well. Hence, we can assume that $u-\phi$ has a strict local minimum at~$x_\circ$. 

Let $\phi_k$ be a sequence of test functions defined for each $u_k$ in the following way
\[
\phi_k := \left\{\begin{array}{ll}
u_k &\quad \text{in}\quad  \R^n\setminus \overline{B_r(x_\circ)}\\
\phi + c_k& \quad \text{in}\quad \overline{B_r(x_\circ)},
\end{array}
\right.
\]
where
\[
c_k := \max\big\{c \in \R : \phi + c\le u_k \quad\text{in}\quad \overline{B_r(x_\circ)}\big\}. 
\]

We  observe that, since $u_k\tosd  u$ in $\Omega$, up to a subsequence we have $c_k\to 0$ and there exists $x_k\to x_\circ$ such that $\phi_k$ touches $u_k$ at $x_k$ from below. Indeed, since the minimum of $u-\phi$ in $\overline{B_r(x_\circ)}$ is 0 (attained at $x_\circ$ strictly) and $u_k-\phi \tosd u-\phi$, by Lemma~\ref{lem:loc_minima} we have that for some $x_k \to x_\circ$,
\[
u_k-\phi \ge u_k(x_k)- \phi(x_k) =: c_k \to 0\quad\text{in}\quad\overline{B_r(x_\circ)}.
\]

Hence, $\phi_k$ is a test function for $u_k$ and since $u_k$ is a supersolution, we have $\I_k(\phi_k, x_k)\le f_k(x_k)$, which is defined classically. We would now like to pass this inequality to the limit. In order to do it, let us compute first, for any $x\in B_{r/4}(x_\circ)$ and $\L \in \LLL$,
\[
\begin{split}
|\L(\phi - \phi_k)(x)|& \le\int_{\R^n\setminus B_{r/2}} |\phi_k(x+y) - \phi(x+y)| K(y)\, dy + C_r |\phi_k(x) - \phi(x)|\\
& \le \Lambda\int_{\R^n\setminus B_{r/2}(x)} |\phi_k(z) - \phi(z)| |z-x|^{-n-2s}\, dy + C_r |c_k|\\
& \le C_r\left(\|u_k-u\|_{L^1_{\omega_s}(\R^n)} + |c_k|\right),
\end{split}
\]
for some constant $C$ that depends on $n$, $s$, $\Lambda$, and $r$. Hence, 
since both $\phi_k$ and $\phi$ are $C^2$ in $B_{r/2}(x_\circ)$, by ellipticity (Corollary~\ref{cor:elliptC11}) we have
\[
\begin{split}
\|\I_k(\phi_k, \cdot) - \I_k(\phi, \cdot)\|_{L^\infty(B_{r/4}(x_\circ))}& \le \left\|\sup_{\L\in \LLL} |\L (\phi - \phi_k)|\right\|_{L^\infty(B_{r/2})}\\
& \le C_r\left(\|u_k-u\|_{L^1_{\omega_s}(\R^n)} + |c_k|\right)\downarrow 0,
\end{split}
\]
as $k\to \infty$, by assumption. On the other hand, we also know that, from the weak convergence of $\I_k$ towards $\I$, 
\[
\|\I_k(\phi, \cdot) - \I(\phi, \cdot)\|_{L^\infty(B_{r/4}(x_\circ))}\downarrow 0\quad\text{as}\quad k \to \infty. 
\]

In all, by the triangle inequality we obtain that $\I_k(\phi_k, x)\to \I(\phi, x)$ uniformly in $B_{r/4}(x_\circ)$. In particular, since each $\I_k(\phi_k, x)$ is continuous, the limit $\I(\phi, x)$ is continuous as well and again by the triangle inequality
\[
|\I_k(\phi_k, x_k) - \I(\phi, x_\circ)|\le |\I_k(\phi_k, x_k) - \I(\phi, x_k)| + |\I(\phi, x_k) - \I(\phi, x_\circ)|\downarrow 0
\]
as $k\to \infty$. Hence, since $\I_k(\phi_k, x_k)\le f_k(x_k)$, and $f_k$ are continuous functions converging uniformly to $f$, we get $\I(\phi, x_\circ)\le f(x_\circ)$ and the result is proved. 
\end{proof}
\begin{rem}[Stability of viscosity subsolutions]
\label{rem:stab_subsol}
By considering $-u_k$ and $-u$ in Proposition~\ref{prop:stab_super} we obtain the stability of viscosity subsolutions, where now $(u_k)_{k\in \N}$ and $u \in {\rm USC}(\Omega)$ are such that $u_k\tosu u$ in $\Omega$. 
\end{rem}

\subsection{The comparison principle}

\index{Comparison principle!Viscosity solutions!Continuous}

Using the stability of viscosity solutions, we are now ready to prove  the comparison principle: 

\begin{thm}[Comparison principle]
\label{thm:comparison_viscosity}
Let $s\in (0, 1)$, let $\I\in \III$, and let $\Omega\subset\R^n$   be any bounded open set. Let $u\in {\rm  LSC}(\overline{\Omega})\cap L^1_{\omega_s}(\R^n)$ and $v\in {\rm  USC}(\overline{\Omega})\cap L^1_{\omega_s}(\R^n)$ be such that $u \ge v$ in $\R^n\setminus \Omega$ and 
\[
\I(u, x) \le f(x) \quad\text{and}\quad  \I (v, x)\ge f(x)\quad\text{in}\quad \Omega
\]
in the viscosity sense, for some $f\in C(\Omega)$. Then $u \ge v$ in $\R^n$.
\end{thm}

 In order to prove it, we first show that  the ellipticity conditions from Lemma~\ref{lem:32} (or Corollary~\ref{cor:elliptC11}) also hold true in the viscosity sense: 

\begin{prop}
\label{prop:viscosity_ellipticity}  \index{Ellipticity!Viscosity solutions}
Let $s\in (0, 1)$, let $\Omega\subset \R^n$ be any open set, and let $\I \in \III$. Let $u\in {\rm LSC}({\Omega})\cap L^1_{\omega_s}(\R^n)$ and $v\in {\rm USC}({\Omega})\cap L^1_{\omega_s}(\R^n)$, and let $f, g\in C({\Omega})$. Assume that 
\[
\I (u, x) \le f(x) \quad \text{and}\quad \I (v, x) \ge g(x)\qquad \text{in}\quad \Omega
\]
in the viscosity sense. Then,   we have
\[
\Mm(u -v) \le f - g  \quad\text{in}\quad \Omega
\]
in the viscosity sense,   where $\Mm$ is given by \eqref{eq:MMpm}.
\end{prop}

The previous result follows by a now classical idea of Jensen,  \cite{Jen88}, of regularizing a semi-continuous function through its sup- or inf-convolution. In the following, given an open set $\Omega$, we fix some $D \ssubset\Omega$, with $D$ open and bounded.

Now, given $\eps > 0$ and $u\in {\rm LSC}({\Omega})$, we define the inf-convolution $u_\eps$ in~$\overline{D}$ as 
\begin{equation}
\label{eq:u_eps} \index{inf-convolution}
u_\eps(x) := \inf_{\overline{D}} \left(u(z) + \frac{|x-z|^2}{\eps}\right),\quad\text{for any}\quad x\in \overline{D},
\end{equation}
 and $u_\eps = u$ in $\R^n\setminus D$. 
If $v\in {\rm USC}({\Omega})$, we define the sup-convolution $v^\eps$ in~$\overline{D}$ as
\[\index{sup-convolution}
v^\eps(x) := \sup_{\overline{D}} \left(v(z) - \frac{|x-z|^2}{\eps}\right)\quad\text{for any}\quad x\in \overline{D},
\]
 and $v^\eps = v$ in $\R^n\setminus D$. 
Notice that, immediately by definition, $u_\eps$ and $v^\eps$ are Lipschitz in $\overline{D}$ (being the supremum and infimum of Lipschitz functions) and $u_\eps \le u$ and $v^\eps\ge v$. Notice, also, that $u_\eps = -(-u)^\eps$.

Finally, we also have that $u_\eps$ is semiconcave, in the sense that $u_\eps - C_\eps|x|^2$ is concave for some constant $C_\eps$.

We then have the following:
\begin{lem}
\label{lem:ueps_conv}
Let $D$  and $\Omega$ be open, with $D\ssubset \Omega$. Let $u\in {\rm LSC}(\Omega)$, and  let $u_\eps$  be defined as above. Then $u_\eps\uparrow u$ pointwise in $\R^n$, $u_\eps$ is semiconcave in $\overline{D}$, and $u_\eps\tosd u$ in $D$ as $\eps\downarrow 0$. 
\end{lem}
\begin{proof}
We know that $u_\eps = u $ in $\R^n\setminus D$, and $u_{\eps'}(x) \le u_\eps(x) \le u(x)$ in $\R^n$, for $\eps\le \eps'$. Let now $x_\circ\in D$ be fixed. If $u(x_\circ) = +\infty$, then $\liminf_{x\to x_\circ} u(x) = \infty$ and $\liminf_{\eps\downarrow 0} u_\eps(x_\eps) = \infty$ for any $x_\eps\to x_\circ$ as $\eps\downarrow 0$, and $u_\eps\tosd u$ at $x_\circ$ (in the sense \eqref{eq:halfrelaxedpointwise} from Definition~\ref{defi:half_relaxed}). We assume therefore $u(x_\circ) < +\infty$. 

Since $u\in {\rm LSC}(\Omega)$ and $\overline{D}\subset \Omega$, by definition of $u_\eps(x)$ for any $x\in D$ there exists $z_\eps(x)\in \overline{D}$ such that 
\[
u_\eps(x) = u(z_\eps(x)) + \frac{|x-z_\eps(x)|^2}{\eps}
\]
(lower semi-continuous functions always achieve their minimum in a compact set). In particular, since $u(x_\circ)< +\infty$ and $u$ is bounded below in $\overline{D}$, we obtain $x_\eps \to x_\circ$ as $\eps \downarrow 0$, where we have denoted $x_\eps := z_\eps(x_\circ)$. Moreover, we also have 
\[
u(x_\eps) -u(x_\circ)\le u_\eps(x_\circ)-u(x_\circ) \le 0
\]
and taking the $\liminf$ in $\eps\downarrow 0$ (again by lower semi-continuity of $u$), 
\[
0 = \liminf_{\eps\downarrow 0} u(x_\eps) -u(x_\circ)\le \liminf_{\eps\downarrow 0} u_\eps(x_\circ)-u(x_\circ) \le 0.
\]
Hence, since $u_\eps(x_\circ) \le u(x_\circ)$, we get  $\lim_{\eps\downarrow 0} u_\eps(x_\circ) = u(x_\circ)$ and $u_\eps\uparrow u$ pointwise in $\R^n$.

Let us now show $u_\eps \tosd u$ in $D$.  Let $x_k \to x_\circ$ and $\eps_k \downarrow 0$ be fixed sequences. If $\liminf_k u_{\eps_k}(x_k) = +\infty$   we are done. Assume $\liminf_k u_{\eps_k}(x_k) < +\infty$. Then, since
\[
u(x_k) \ge u_{\eps_k}(x_k) = u(z_{\eps_k}(x_k)) + \frac{|x_k - z_{\eps_k}(x_k)|^2}{\eps_k},
\]
we have 
\[
\liminf_{k\to\infty} |x_k - z_{\eps_k}(x_k)|^2 \le \liminf_{k\to \infty} \eps_k \big(u(x_k) - u(z_{\eps_k}(x_k)\big) = 0,
\]
where we have used that $\liminf_k u_{\eps_k}(x_k) < +\infty$ and that $u$ is bounded below in $\overline{D}$. Thus, there is a subsequence of $z_{\eps_k}(x_k)$ converging to $x_\circ$, and since $u$ is lower semi-continuous we get
\[
\liminf_{k\to \infty } u(x_k) \ge \liminf_{k\to \infty} u_{\eps_k} (x_k) \ge \liminf_{k\to\infty} u(z_{\eps_k}(x_k)) \ge u(x_\circ).
\]
Combined with the fact that $u_\eps(x_\circ)\uparrow u(x_\circ)$, we get that $u_\eps\tosd u$ at $x_\circ\in D$ (in the sense \eqref{eq:halfrelaxedpointwise}). Since $D$ is arbitrary, $u_\eps \tosd u$ in $D$. 

Finally, let again $x_\circ\in \overline{D}$ be fixed and let $x_\eps = z_\eps(x_\circ)\in \overline{D}$. Then, 
\[
u_\eps(x) \le u(x_\eps) + \frac{|x-x_\eps|^2}{\eps}\quad\text{for any}\quad x\in \overline{D},
\]
with equality at $x = x_\circ$. That is, there is a paraboloid of opening $\frac{2}{\eps}$ touching $u_\eps$ from above at $x = x_\circ$. Since $x_\circ$ is arbitrary, we get that $u_\eps$ can be touched from above at any point in $\overline{D}$ by a paraboloid of opening $\frac{2}{\eps}$; that is, $u_\eps$ is semiconcave. Alternatively, 
\[
u_\eps(x) - \frac{|x|^2}{\eps}\le u(x_\eps) + \frac{|x-x_\eps|^2 - |x|^2}{\eps} = u(x_\eps) +\frac{|x_\eps|^2 - 2x\cdot x_\eps}{\eps}
\]
with equality at $x = x_\circ$, and so $u_\eps(x)- \frac{|x|^2}{\eps}$ is below a tangent line, and hence it is concave.
\end{proof}

Thanks to the result above, we have the following proposition, saying that not only $u_\eps$ is a good lower semi-continuous approximation of $u$, it is also a good approximation as a supersolution:
\begin{prop}
\label{prop:ueps_conv}
Let $s\in (0, 1)$, let $\Omega\subset \R^n$ be any open set, and $\overline{D}\subset \Omega$. Let $\I \in \III$, let $u\in {\rm LSC}(\Omega)\cap L^1_{\omega_s}(\R^n)\cap L^\infty(D)$, and let $u_\eps$ be defined by \eqref{eq:u_eps}. 
Let us suppose that $\I (u, x) \le 0$ in $D$ in the viscosity sense. 

Then, for every $D'\ssubset D$ and every $\eps > 0$ small enough there exists $\delta_\eps$ such that $\I( u_\eps, x) \le   \delta_\eps$ in the viscosity sense in $D'$, with $\delta_\eps \downarrow 0$ as $\eps \downarrow 0$.
\end{prop}
\begin{proof}
We use the same notation as in Lemma~\ref{lem:ueps_conv}, and we fix some $x_\circ\in D'$ and denote $x_\eps = z_\eps(x_\circ)\in \overline{D}$. That is, 
\[
u(x_\circ) \ge u_\eps(x_\circ) = u(x_\eps) + \frac{1}{\eps}|x_\circ-x_\eps|^2\quad\Rightarrow \quad |x_\circ - x_\eps|^2 \le \eps \osc_{\overline{D}} u,
\]
where $\osc_{\overline{D}} u < \infty$ since $u\in L^\infty (D)$.  That is, $x_\eps \to x_\circ$ as $\eps \downarrow 0$, at a rate independent of $x_\circ$. 

Let now $\phi$ be any test function touching $u_\eps$ from below at $x_\circ$, i.e., $\phi \le u$ and  $\phi(x_\circ) = u_\eps(x_\circ)$. We want to show that $\I(\phi, x_\circ) \le \delta_\eps$.

We have,  for any $x\in D'$, 
\[
\phi(x+x_\circ -x_\eps) \le u_\eps(x + x_\circ - x_\eps) \le u(x) + \frac{1}{\eps}|x_\circ - x_\eps|^2,
\]
where we assume $\eps$ is small enough such that $x+x_\circ - x_\eps \in D$. Observe that the previous inequalities are all equalities at $x = x_\eps$, and therefore, the function 
\[
\phi(x+x_\circ -x_\eps) -\frac{1}{\eps}|x_\circ - x_\eps|^2
\]
is a test function for $u$ touching from below at $x_\eps$. In particular, 
\[
\I\big(\phi(\,\cdot\,+x_\circ-x_\eps), x_\eps\big) = \inf_{b\in \B} \sup_{a\in \A}\big\{(\L_{ab} \phi)(x_\circ) + c_{ab}(x_\eps)\big\}\le 0.
\]
Hence, using that $c_{ab}(x_\eps)-c_{ab}(x_\circ) \ge - \omega(|x_\eps-x_\circ|)$ for all $(a, b)\in \A\times \B$ for some modulus of continuity $\omega$ independent of $(a, b)$, we get 
\[
\I (\phi, x_\circ) \le \omega(|x_\eps-x_\circ|) \le \omega \big( {(\eps \,\osc_{\overline{D}} u)^{1/2}}\big)=: \delta_\eps,
\]
as we wanted to see. 
\end{proof}

\begin{rem}\label{rem:changeroles}By changing the sign of $u$, the statements of Lemma~\ref{lem:ueps_conv} and Proposition~\ref{prop:ueps_conv} have their analogues for the sup-convolution $u^\eps$. Namely, we have that $u^\eps\downarrow u$ pointwise in $\R^n$, $u^\eps$ is semiconvex in $\overline{D}$, and $u^\eps\tosu u$ in $D$ as $\eps \downarrow 0$. Moreover, if $\I (u, x) \ge 0$ in $D$, then $\I (u^\eps, x) \ge -\delta_\eps$ in $D'$, with $\delta_\eps\downarrow 0$ as $\eps \downarrow 0$. 
\end{rem}

We can now show Proposition~\ref{prop:viscosity_ellipticity}:
\begin{proof}[Proof of Proposition~\ref{prop:viscosity_ellipticity}]
Let us fix any $D'\ssubset D\ssubset \Omega$, and let us show $\Mm(u-v) \le f-g$ in $D'$. Since $D'$ is arbitrary, this will imply $\Mm(u-v) \le f-g$ in $\Omega$. We divide the proof into two steps:

\begin{steps}
\item \label{it:ellipti1} We assume first that $u, v\in L^\infty_{\rm loc}(\Omega)$. By Proposition~\ref{prop:ueps_conv} (see also Remark~\ref{rem:changeroles}) we have\footnote{We apply Proposition~\ref{prop:ueps_conv} to the operator $\I_g(u, x) = \I(u, x) - g$ (and $\I_f(u, x) = \I(u, x) - f$), where we are considering the new $c_{ab}'(x) = c_{ab}(x) - f(x)$, that are still equicontinuous in $\overline{D}$.} $\I (u_\eps, x) \le f + \delta_\eps$ and $\I (u^\eps, x) \ge g -\delta_\eps$ in $D'$ for some $\delta_\eps\downarrow 0$ as $\eps\downarrow 0$. Let $\phi$ be a test function touching $u_\eps - v^\eps$ from below at some point $x_\circ \in D'$. Then, since $u_\eps$ and $-v^\eps$ are  semiconcave (by Lemma~\ref{lem:ueps_conv} and Remark~\ref{rem:changeroles}), and $u_\eps - v^\eps$ can be  touched from below by a paraboloid at~$x_\circ$, this implies that both $u_\eps$ and $v^\eps$ are $C^{1,1}$ at $x_\circ$, that is, $u_\eps, v^\eps\in C^{1,1}_{r}(x_\circ)$ for some $r > 0$. By Corollary~\ref{cor:elliptC11}, 
\[
\Mm(u_\eps - v^\eps)(x_\circ) \le \I(u_\eps, x_\circ) -\I(v^\eps, x_\circ)\le f-g+2\delta_\eps.
\]
(We also used here Lemma~\ref{lem:C11punt1} to take $u_\eps$ and $v^\eps$ as admissible test functions at $x_\circ$, see also Remark~\ref{rem:C11}). Since $\phi \le u_\eps-v^\eps$, this implies $\Mm \phi (x_\circ) \le f-g+2\delta_\eps$ and thus $\Mm (u_\eps-v^\eps)  (x_\circ) \le f-g+2\delta_\eps$ in $D'$ in the viscosity sense. 

We  use the stability property from Proposition~\ref{prop:stab_super} (and Remark~\ref{rem:stab_subsol}) together with Lemma~\ref{lem:ueps_conv} and Proposition~\ref{prop:ueps_conv} to deduce that 
\[
\Mm (u -v)  \le f-g\quad\text{in}\quad D'
\]
in the viscosity sense. Since $D'$ was arbitrary, we get the desired result  whenever $u, v\in L^\infty_{\rm loc}(\Omega)$.

\item Let us now prove the proposition for any $u\in {\rm LSC}({\Omega})\cap L^1_{\omega_s}(\R^n)$ and $v\in {\rm USC}({\Omega})\cap L^1_{\omega_s}(\R^n)$. For this,  let $D\ssubset \Omega$, and let us fix any $\eta\in C^\infty(\R^n)\cap L^\infty(\R^n)$ such that $\Mm \eta  \ge 1$ in $D$\footnote{This is satisfied, for example, by any test function $\eta\in C^\infty_c(\R^n)$ such that $0 \le \eta\le M$ and $\eta \equiv M$ in $D$, if $M$ is large enough.}. Then, we have that 
\[
\left\{
\begin{array}{rcll}
\I(-C\eta, x) & \le&  \I(0, x)  -C& \quad\text{in}\quad D,\\
\I(C\eta, x) & \ge & \I(0, x)+  C& \quad\text{in}\quad D.
\end{array}
\right.
\]  
We take $C = \|\I(0, x)\|_{L^\infty(\Omega)} + \|f\|_{L^\infty(D)}$, so that $\I(-C\eta, x) \le -\|f\|_{L^\infty(D)}$ and $\I(C\eta, x) \ge \|f\|_{L^\infty(D)}$ in $D$. Then, by definition of viscosity sub- and supersolution (see Remark~\ref{rem:max_min}), we have that if we denote $u_\ell := \min\{u, \ell-C\eta\}$ and $v^\ell := \max\{v, -\ell + C\eta\}$ for any $\ell > 0$,
\[
\I(u_\ell, x) \le f\quad\text{and}\quad \I(v^\ell, x) \ge g\quad\text{in}\quad D. 
\] 
Moreover, now $u_\ell$ and $v^\ell$ are bounded, so that by \ref{it:ellipti1}  we deduce
\[
\Mm(u_\ell - v^\ell)\le f - g \quad\text{in}\quad D. 
\]
Finally, we let $\ell\to \infty$, so that we can apply the stability in Proposition~\ref{prop:stab_super} and deduce 
\[
\Mm(u - v)\le f - g \quad\text{in}\quad D. 
\]
Since $D$ was arbitrary, we now get the desired result. \qedhere
\end{steps}
\end{proof}

The following lemma is the minimum principle, and says that supersolutions in a domain attain their minimum in the exterior of such domain. In order to have it, we need to impose that the upper semi-continuity  holds up to the boundary.

\begin{lem}
\label{lem:fromabove} \index{Maximum principle!Viscosity solutions!Continuous}
Let $s\in (0, 1)$, and let $\Omega\subset \R^n$ be any bounded open set. Let  $\Mm $ be given by \eqref{eq:MMpm}, and let $u\in {\rm LSC}(\overline{\Omega})\cap L^1_{\omega_s}(\R^n)$  be such that 
\[
\left\{
\begin{array}{rcll}
 \Mm u & \le & 0 & \quad \text{in}\quad \Omega,\\
  u& \ge & 0 & \quad \text{in}\quad \R^n\setminus \Omega,
\end{array}
\right.
\]
in the viscosity sense. Then, $u \ge 0$ in $\R^n$. 
\end{lem}
\begin{proof} Let us suppose that the conclusion is not true, and the function $u$ attains its negative minimum at some $x_\circ\in \overline{\Omega}$. Notice, first, that by assumption $x_\circ\in \Omega$, since $u \ge 0$ on $\partial \Omega$. Then, we immediately have that the function 
\[
\varphi(x) := u(x_\circ) \chi_{\Omega}(x)
\]
is an admissible test function from below for $u$, so that by the viscosity definition  we should have $\Mm \varphi(x_\circ) \le 0$. However, since $\Omega$ is bounded, we can directly compute $\Mm \varphi(x_\circ)  > 0$ ($x_\circ$ is a global minimum), a contradiction. 
\end{proof}
\begin{rem}
Notice that the condition that we need on $u$ (which is implied by the hypotheses) is that for any $x_\circ\in \Omega$,
\[
\liminf_{\Omega\ni x \to x_\circ} u \ge 0.
\]
\end{rem}

Finally, we can show the comparison principle for viscosity solutions:
\begin{proof}[Proof of Theorem~\ref{thm:comparison_viscosity}]
By Proposition~\ref{prop:viscosity_ellipticity} we have $\Mm (u-v)\le 0$ in~$\Omega$, and by assumption $u- v \ge 0$  in $\R^n\setminus \Omega$. Hence, by Lemma~\ref{lem:fromabove} we deduce 
\[
u - v \ge   0\quad\text{in}\quad \R^n,
\]
which is the desired result. 
\end{proof}

Notice that as a consequence of the comparison principle we have the uniqueness of viscosity solutions:

\begin{cor}[Uniqueness of continuous viscosity solutions]
\label{cor:uniqueness_viscosity}
 \index{Uniqueness!Viscosity solutions!Continuous}
Let $s\in (0, 1)$, let $\I\in \III$, and let $\Omega\subset\R^n$ be any bounded open set. Let $u_1, u_2 \in C(\overline{\Omega})\cap L^1_{\omega_s}(\R^n)$ be such that $u_1 = u_2$ in $\R^n\setminus \Omega$ and 
\[
\I(u_i, x) = f(x) \quad\text{in}\quad   \Omega,\quad  \text{for}\quad i = 1,2,
\]
for some $f\in C({\Omega})$, in the viscosity sense. Then $u = v$ in $\R^n$.
\end{cor}
 \begin{proof}
 We apply the comparison principle, Theorem~\ref{thm:comparison_viscosity}, to obtain both $u \ge v$ and $u \le v$ in $\R^n$. 
 \end{proof}
 
Finally, as in subsection~\ref{ssec:Linftybounds}, the comparison principle also yields the $L^\infty$ bound for solutions with bounded exterior data: 
\begin{cor}
\label{cor:Linftybound_visc} \index{Linfty bounds@$L^\infty$ bounds!Viscosity solutions!Continuous}
Let $s\in (0,1)$ and   $\Mpm $ be given by \eqref{eq:MMpm}. Let $\Omega\subset \R^n$ be any bounded open set, and $g\in L^\infty(\R^n\setminus \Omega)$. 

Let $u\in C(\overline{\Omega})\cap L^1_{\omega_s}(\R^n)$ be a viscosity solution of 
\begin{equation}
\label{eq:sol_u_g}
\left\{
\begin{array}{rcll}
\Mp u & \ge & -C_\circ & \quad\text{in}\quad \Omega,\\
\Mm u & \le & C_\circ & \quad\text{in}\quad \Omega,\\
u & = & g & \quad\text{in}\quad \R^n\setminus \Omega.
\end{array}
\right.
\end{equation}
Then, 
\[
\|u\|_{L^\infty(\Omega)} \le \|g\|_{L^\infty(\R^n\setminus \Omega)} + C C_\circ,
\]
for some constant $C$ depending only on $n$, $s$, $\lambda$, $\Lambda$, and ${\rm diam}(\Omega)$. 
\end{cor}
\begin{proof}
The proof follows as that in Lemma~\ref{lem:Linftybound}. Let $w\in C^\infty_c(\R^n)$ be the barrier from Lemma~\ref{lem:barrierLinfty}. Then, $v := \|g\|_{L^\infty(\R^n\setminus \Omega)} + C_\circ w$  satisfies 
\[
\left\{
\begin{array}{rcll}
\Mp v & \le & -C_\circ & \quad\text{in}\quad \Omega,\\
v & \ge & \|g\|_{L^\infty(\R^n\setminus \Omega)} & \quad\text{in}\quad \R^n\setminus \Omega,\\
v & \le & \|g\|_{L^\infty(\R^n\setminus \Omega)}+CC_\circ  & \quad\text{in}\quad \Omega,
\end{array}
\right.
\]
for some constant $C>0$ depending only on $n$, $s$, $\lambda$, $\Lambda$, and ${\rm diam}(\Omega)$. We have $\Mp u \ge \Mp v$ in $\Omega$ and $v \ge u$ in $\R^n\setminus \Omega$. By the comparison principle, Theorem~\ref{thm:comparison_viscosity} (applied with operator $\I(u, x) := \Mp u(x) + C_\circ$) we deduce
$
v\ge u 
$
in $\R^n$.  That is, 
\[
u \le \|g\|_{L^\infty(\R^n\setminus \Omega)}+CC_\circ \quad\text{in}\quad \Omega.
\]
By taking $-u$ instead of $u$, we prove the other inequality to obtain the desired result.
\end{proof}

In fact, in order to get an $L^\infty$ bound in $\Omega$, it is enough to impose that $g$ is bounded near $\Omega$; see Corollary~\ref{cor:Linftybound_visc_2_2} below.

\subsection{Comparison principle without boundary continuity}
 
In the comparison principle above, Theorem~\ref{thm:comparison_viscosity}, we had to impose for the functions $u$ and $v$ to be lower (resp. uppper) semi-continuous up to the boundary of $\overline{\Omega}$. For (local) elliptic PDEs, this is a necessary assumption, since otherwise we could have discontinuous solutions like $u(x) = 1$ in $\Omega$ with $u(x) = 0$ on $\partial\Omega$. Interestingly, it turns out that for nonlocal equations, this assumption is often not necessary, and that we can prove the comparison principle even for solutions that are \textit{discontinuous} on $\partial \Omega$. For this, we need some extra assumptions on the domain $\Omega$. More precisely, we will consider domains   $\Omega\subset \R^n$ satisfying:
\begin{equation}
\label{eq:cond_Om}
\text{for all $r  > 0$ and $z\in \partial\Omega$ there is a ball $B_{\kappa r}(x_{r, z})\subset \Omega^c \cap B_r(z)$,}
\end{equation}
for some $\kappa > 0$.  Notice that any bounded Lipschitz domain satisfies \eqref{eq:cond_Om}, but the assumption is actually much more general, and includes quite rough domains. (This is often called the \emph{exterior corkscrew condition} \cite{JK82}.) 

In this context, we have the comparison principle for bounded solutions: 

\begin{thm}[Comparison principle without boundary continuity]
\label{thm:comparison_viscosity_2} \index{Comparison principle!Viscosity solutions!Bounded}
Let $s\in (0, 1)$, let $\I\in \III$, and let $\Omega\subset\R^n$ be any bounded  open set satisfying \eqref{eq:cond_Om} for some $\kappa > 0$. Let $u\in {\rm  LSC}({\Omega})\cap L^1_{\omega_s}(\R^n)$ with $\inf_{x\in \Omega} u > -\infty$ and $v\in {\rm  USC}({\Omega})\cap L^1_{\omega_s}(\R^n)$ with $\sup_{x\in \Omega} v  < \infty$ be such that $u \ge v$ in $\R^n\setminus \Omega$ and 
\[
\I(u, x) \le f(x) \quad\text{and}\quad  \I (v, x) \ge f(x) \quad\text{for}\quad x\in \Omega
\]
in the viscosity sense, for some $f\in C(\Omega)$. Then $u \ge v$ in $\R^n$.
\end{thm}
Notice that neither $u$ and $v$ nor $f$ are assumed to be (semi-) continuous up to the boundary. In order to prove this result, it is enough to show that the minimum principle holds in this context: 
\begin{lem}
\label{lem:fromabove_2} \index{Maximum principle!Viscosity solutions!Bounded}
Let $s\in (0, 1)$, and let $\Omega\subset \R^n$ be any bounded open set satisfying \eqref{eq:cond_Om} for some $\kappa > 0$. Let  $\Mm $ be given by \eqref{eq:MMpm}, and let $u\in {\rm LSC}({\Omega})\cap L^1_{\omega_s}(\R^n)$ with $\inf_{x\in \Omega} u > -\infty$  be such that 
\[
\left\{
\begin{array}{rcll}
 \Mm u & \le & 0 & \quad \text{in}\quad \Omega,\\
  u& \ge & 0 & \quad \text{in}\quad \R^n\setminus \Omega,
\end{array}
\right.
\]
in the viscosity sense. Then, $u \ge 0$ in $\R^n$. 
\end{lem}

\begin{figure}
\centering
\makebox[\textwidth][c]{\includegraphics[scale = 1]{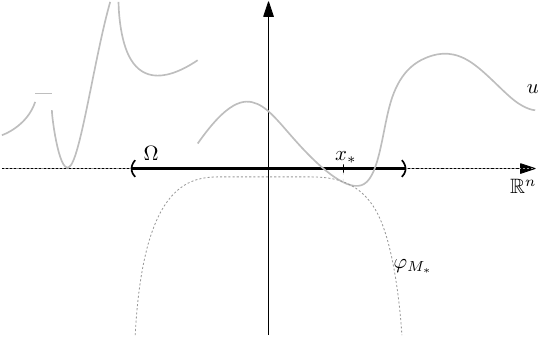}}
\caption{\label{fig:09_2} The barrier $\varphi_{M, \eps}$ slidden from below   blows-up on the boundary of $\Omega$.}
\end{figure}

\begin{proof}
The proof follows by sliding an appropriate (singular) barrier from below. We take 
\[
\varphi(x) :=  \dr_\Omega^{-\eps_\circ/2}(x)\chi_{\Omega}(x)
\]
given by Lemma~\ref{barrier-singular-eps}, which satisfies 
\[
\Mp \varphi \le -c_* < 0\quad\text{in}\quad \Omega. 
\]
If we now define, for $M > 0$, 
 $
 \varphi_{M} := -M\varphi,
 $
 then it satisfies $\Mm \varphi_{M} \ge M c_* > 0$ in~$\Omega$. We   use  $\varphi_{M }$ as a barrier from below. Let us denote $M_*$ as
\[
 M_* := \inf\{ M > 0 : \varphi_{M} \le u\quad\text{in}\quad \R^n\}. 
 \]
(See Figure~\ref{fig:09_2}.) Observe first that $M_*$ is well-defined, since the infimum is taken on a nonempty set: $u\ge 0 = \varphi_M$ in $\R^n\setminus \Omega$ for every $M \ge 0$; and since $\varphi \ge c_\circ > 0$ in $\Omega$ for some $c_\circ$, and $\inf_{x\in \Omega} u >-\infty$  by assumption, we have that if $M$ is large enough then $u \ge \varphi_M$ in $\Omega$. 

Let us suppose now that $M_*  >0$. Since  $u\in {\rm LSC}({\Omega})$  is bounded below, we have that the minimum of $u-\varphi_{M_*}$ (which is zero) is attained at some $x_*\in \Omega$ (observe that $\varphi$ blows-up on $\partial\Omega$). Hence, $\varphi_{M_*}$ is an admissible test function for $u$, and  from the viscosity condition we have $\Mm \varphi_{M_*}(x_*) \le 0$, which contradicts $\Mm \varphi_{M_*} >0$ in $\Omega$. Thus $M_* = 0$ and we get $u \ge 0$ in $\R^n$, as wanted. 
\end{proof}

As a consequence, we immediately get:
\begin{proof}[Proof of Theorem~\ref{thm:comparison_viscosity_2}]
By Proposition~\ref{prop:viscosity_ellipticity} we have $\Mm (u-v)\le 0$ in $\Omega$, and by assumption $u- v \ge 0$  in $\R^n\setminus \Omega$. Hence, by Lemma~\ref{lem:fromabove_2} we deduce 
$
u - v \ge   0 
$ in $\R^n$, 
which is the desired result. 
\end{proof}

We also get the following two corollaries; the first one on the uniqueness of (bounded) solutions: 

\begin{cor}[Uniqueness of bounded viscosity solutions]
\label{cor:uniqueness_viscosity_2} \index{Uniqueness!Viscosity solutions!Bounded}
Let $s\in (0, 1)$,  let $\I\in \III$, and let $\Omega\subset\R^n$ be any bounded open set satisfying \eqref{eq:cond_Om} for some $\kappa > 0$. Let $u_1, u_2 \in {  C}({\Omega})\cap L^1_{\omega_s}(\R^n)\cap L^\infty(\Omega)$ be such that $u_1 = u_2$ in $\R^n\setminus \Omega$ and 
\[
\I(u_i, x) = f(x) \quad\text{in}\quad   \Omega,\quad \text{for}\quad i = 1,2,
\]
for some $f\in C({\Omega})$, in the viscosity sense. Then $u = v$ in $\R^n$.
\end{cor}
 \begin{proof}
 We apply the comparison principle, Theorem~\ref{thm:comparison_viscosity_2}, to obtain both $u \ge v$ and $u \le v$ in $\R^n$. 
 \end{proof}
 And the second one, which is a corollary on the $L^\infty$ bound with bounded exterior datum (near the boundary): 
 \begin{cor}
\label{cor:Linftybound_visc_2_2} \index{Linfty bounds@$L^\infty$ bounds!Viscosity solutions!Bounded}
Let $s\in (0,1)$ and $\Mpm$ be given by \eqref{eq:MMpm}. Let $\Omega\subset \R^n$  be any bounded open set satisfying \eqref{eq:cond_Om} for some $\kappa > 0$, and let $g\in L^1_{\omega_s}(\R^n\setminus\Omega)$ be such that 
\[
|g(x)|\le C_g\quad\text{in}\quad \{x\in \R^n\setminus\Omega : \dist(x, \Omega) \le \rho\},
\]
holds for some $C_g \ge 0$ and $\rho > 0$. 

Let $u\in C( {\Omega})\cap L^1_{\omega_s}(\R^n)\cap L^\infty(\Omega)$ be a viscosity solution of \eqref{eq:sol_u_g}.  Then, 
\[
\|u\|_{L^\infty(\Omega)} \le  C\left(\|g\|_{L^1_{\omega_s}(\R^n\setminus \Omega)} + C_\circ\right)+ C_g,
\]
for some constant $C$ depending only on $n$, $s$, $\lambda$, $\Lambda$, $\rho$, and ${\rm diam}(\Omega)$. 
\end{cor}
\begin{proof}
Let us denote $\Omega_\rho := \{x\in \R^n : \dist(x, \Omega) \le \rho/2\}$. Let $\L\in \LLL$, and notice that for any $x\in\Omega$
\[
\L\chi_{\Omega_\rho}(x) = \int_{\R^n\setminus (\Omega_\rho - x)} K(y)\, dy,
\]
and thus, if we assume $\Omega\subset B_R$ for some $R \le 2{\rm diam}(\Omega)$ (after a translation) we have
\[
\Mp \chi_{\Omega_\rho}  \le -\lambda\int_{B_{2R}^c}|y|^{-n-2s}\, dy = -C_1\quad\text{in}\quad \Omega,
\]
for some $C_1$ that depends only on $n$, $s$, $\lambda$, and ${\rm diam}(\Omega)$. 

On the other hand, thanks to Lemma~\ref{lem:Lu_LL} we know that $\Mp (g\chi_{\Omega_\rho^c})\le C_2\|g\|_{L^1_{\omega_s}(\R^n)}$ in $\Omega$, for some $C_2$ depending only on $n$, $s$, $\Lambda$, and $\rho$. Thus, if we consider the function $v := g\chi_{\Omega_\rho^c} + C_* \chi_{\Omega_\rho}$ we have by the previous inequalities
\[
\Mp v \le  C_2\|g\|_{L^1_{\omega_s}(\R^n)} - C_*C_1 \le -C_\circ,
\]
where we have chosen $C_*\ge C_1^{-1}\left(C_\circ+C_2\|g\|_{L^1_{\omega_s}(\R^n\setminus\Omega)}\right)$. If we furthermore impose $C_* \ge C_g$ we can apply the comparison principle, Theorem~\ref{thm:comparison_viscosity_2}, to $v$ and $u$ and deduce 
\[
u \le v \le C_*\quad\text{in}\quad \Omega. 
\]
By replacing the role of $u$ by $-u$, we obtain the desired result. 
\end{proof}

\subsection{Existence of viscosity solutions} 
\label{subsection3.2.5} 

\index{Viscosity solutions!Existence}

Having proved the comparison principle for viscosity solutions, we now have the tools required to establish their existence. 

In the following, we say that $g\in L^1_{\omega_s}(\R^n)$ is bounded near $\partial\Omega$ with respect to $\R^n\setminus\Omega$ if there exists $\rho > 0$ and $C_g > 0$ such that 
\begin{equation}
\label{eq:bounded_from_outside}
|g|\le C_g \quad\text{in}\quad \left\{x\in \R^n\setminus \Omega : 0<\dist(x, \Omega) <\rho\right\}. 
\end{equation}

We say that $g\in L^1_{\omega_s}(\R^n)$ is continuous on $\partial\Omega$ with respect to $\R^n\setminus\Omega$ if there exists a modulus of continuity $\sigma:[0, \infty]\to [0, \infty]$  such that
\begin{equation}
\label{eq:cont_from_outside}
|g(x) - g(y)|\le \sigma(|x-y|)\quad\text{for all}\quad x\in \partial\Omega, ~~y\in \R^n\setminus\Omega. 
\end{equation}

Notice that if $g$ is continuous on $\partial\Omega$, then it is   bounded near $\partial\Omega$. 

The theorem we want to prove is the following: 
 
\begin{thm}[Existence of viscosity solutions]
\label{thm:existence_visc}
Let $s\in (0, 1)$, and let $\Omega$ be any bounded open set satisfying
\begin{equation}
\label{eq:cond_Om0}
\text{for all $r  > 0$ and $z\in \partial\Omega$ there is a ball $B_{\kappa r}(x_{r, z})\subset \Omega^c \cap B_r(z)$,}
\end{equation}
for some $\kappa > 0$. Let $\I \in \III$, and let $g\in L^1_{\omega_s}(\R^n)$ be bounded near $\partial\Omega$ with respect to $\R^n\setminus \Omega$, in the sense of \eqref{eq:bounded_from_outside}. Then, there exists a unique viscosity solution $u\in C({\Omega})\cap L^1_{\omega_s}(\R^n)\cap L^\infty(\Omega)$ of
\[
\left\{
\begin{array}{rcll}
\I(u, x) & = & 0& \quad\text{in}\quad \Omega\\
u & = & g& \quad\text{in}\quad \R^n\setminus \Omega. 
\end{array}
\right.
\]
Moreover, if $g$ is continuous on $\partial\Omega$ with respect to $\R^n\setminus \Omega$, in the sense of \eqref{eq:cont_from_outside}, then $u\in C(\overline{{\Omega}})$ as well. 
 \end{thm}

Notice that the assumption on the domain is very mild. It includes any bounded Lipschitz domain $\Omega\subset \R^n$, but also much rougher domains, and even sets whose boundary has Hausdorff dimension strictly greater than $n-1$, like Koch's snowflake. It is quite surprising (and new) to obtain such a general existence and uniqueness result for fully nonlinear equations.  
 
 The proof of the previous result follows by Perron's method, since we already have the comparison principle in Theorem~\ref{thm:comparison_viscosity}. We will show it first for globally bounded functions, and then we use it  to obtain the result for an exterior datum in $L^1_{\omega_s}(\R^n)$.
  
 We start with the following lemma, that says that the infimum of a family of supersolutions is a supersolution as well: 
 \begin{lem}
 \label{lem:inf_is_super}
Let $s\in (0, 1)$, let $\I \in \III$, and let $D\subset \R^n$. Let $(u_a)_{a\in \A}$ be a family of supersolutions, $u_a\in {\rm LSC}( {D})\cap L^1_{\omega_s}(\R^n)$   uniformly bounded from below in $D$, such that 
 \[
 \I(u_a, x) \le 0\quad\text{in}\quad \mathring{D},\quad\text{for all}\quad a\in \A
 \]
 in the viscosity sense.  Let 
 \[
 u(x) := \inf_{a\in \A} u_a(x),
 \]
 and let us consider its lower semi-continuous envelope in $D$, 
 \[
 u_*(x) := \inf\left\{\liminf_{k\to \infty} u(x_k) : D \ni x_k\to x\right\},
 \]
 with $u_* = u$ in $\R^n\setminus D$. 
 Then $u_*\in {\rm LSC}(D) \cap  L^1_{\omega_s}(\R^n)$ is a viscosity supersolution, $\I(u_*, x) \le 0$ in $\mathring{D}$. 
 \end{lem}
 \begin{proof}
By definition we have $u_*\in {\rm LSC}( D)$. Let any $x_\circ\in \mathring{D}$, and let $\phi$ be any test function such that $\phi\in C^2(B_r(x_\circ))$ for some $r > 0$,  $\phi \le u_*$ in $\R^n$, and $\phi(x_\circ) = u_*(x_\circ)$. As in the beginning of the proof of Proposition~\ref{prop:stab_super} we can assume that $u_*-\phi$ has a strict local minimum at $x_\circ$.

 By definition of $u_*$, there exists a sequence $(u_{a_k})_{k\in \N}$ and $x_k\to x_\circ$ such that $\lim_{k\to \infty} u_{a_k}(x_k) = u_*(x_\circ)$. We can define $\tilde u_*$ to be the lower half-relaxed limit of $u_{a_k}$, 
 \[
 \tilde u_* := \liminfs{k\to \infty} u_{a_k} \in {\rm LSC}( D),
 \]
 and by assumption we have $u_*\le \tilde u_*$ with $u_*(x_\circ) = \tilde u_*(x_\circ)$. Thus, the function $\tilde u_*-\phi$ has a strict local minimum at $x_\circ$, and by Lemma~\ref{lem:loc_minima} we can find indices $k_j\to \infty$ as $j\to \infty$ and points $x_j \to x_\circ$ such that $u_{a_{k_j}}(x_j) \to u_*(x_\circ)$ and $u_{a_{k_j}}-\phi$ has a local minimum at $x_j$. Since $u_{a_{k_j}}$ is a viscosity supersolution, we must have 
 \[
 \I(\phi, x_j) \le 0.
 \] 
 Therefore, since $x_j \to x_\circ$ and by continuity of $\I(\phi, x)$ around $x_\circ$ (by Corollary~\ref{cor:Iu_II}) we get $\I(\phi, x_\circ) \le 0$. That is, $u_*$ is a viscosity supersolution. 
 \end{proof}

  We then have the following: 
  \begin{prop}
 \label{prop:existence_main_bdd}  \index{Perron's method}
Let $s\in (0, 1)$, and let $\Omega$ be any bounded  open set satisfying \eqref{eq:cond_Om0} for some $\kappa > 0$. Let $\I \in \III$, and let $g\in L^\infty(\R^n)$. Then, there exists a unique viscosity solution $u\in C({\Omega})\cap L^\infty(\R^n)$ to
\[
\left\{
\begin{array}{rcll}
\I(u, x) & = & 0& \quad\text{in}\quad \Omega\\
u & = & g& \quad\text{in}\quad \R^n\setminus \Omega. 
\end{array}
\right.
\]
 \end{prop}
 \begin{proof}
 The uniqueness directly follows from Corollary~\ref{cor:uniqueness_viscosity_2}. The proof now follows by Perron's method. We  divide it into two steps: 
 \begin{steps}
 \item  Let 
 \[
 \mathcal{S} := \left\{v \in {\rm LSC}({\Omega})\cap L^\infty(\R^n): \begin{array}{l} \I(v, x)\le 0 \ \text{ in   $\Omega$  in the viscosity sense} \\
 \hspace{0.93cm} v \ge g \ \text{ in}\ \R^n\setminus\Omega\end{array}\right\},
 \]
 and let us define 
 \[
 u(x) := \inf_{v\in \mathcal{S}} v(x). 
 \]
 Let us assume $\Omega\subset B_R$ for some $R > 0$, and consider $\eta\in C^\infty_c(B_{ 2R})$ with $0 \le \eta \le 1$ and $\eta \equiv 1$ in $B_{R}$. Then, $\Mp \eta \le - c < 0$ in $B_R$ and the function 
 \[
 v_\circ := \|g\|_{L^\infty(\R^n\setminus\Omega)} + c^{-1}\|\I(0, x)\|_{L^\infty(\R^n)}\eta\quad\text{is such that}\quad v_\circ \in \mathcal{S}.
 \] 
 In particular, $\mathcal{S}$ is nonempty. Moreover, since $-v_\circ\le g$ in $\R^n\setminus\Omega$, and 
 \[
 \I(-v_\circ, x) \ge \Mm(-v_\circ)+\I(0, x) \ge  -\Mp(v_\circ)-\|\I(0, x)\|_{L^\infty(\Omega)}\ge 0 \quad\text{in}\quad \Omega,
 \]
 by the comparison principle (Theorem~\ref{thm:comparison_viscosity_2}) all elements in $\mathcal{S}$ are bounded below by $-v_\circ$, and therefore, $u$ is globally bounded. 
 We define its lower semi-continuous envelope in $ {\Omega}$, 
 \[
 u_*(x) := \inf\left\{\liminf_{k\to \infty} u(x_k) :  {\Omega}\ni x_k\to x\right\},
 \]
 with $u_* = u$ in $\R^n\setminus {\Omega}$. By Lemma~\ref{lem:inf_is_super} we already know that
 \[
 \I(u_*, x) \le 0\quad\text{in}\quad \Omega. 
 \]
Moreover, since we can change the value of $v_\circ$ around one point outside of $\R^n\setminus {\Omega}$, we also immediately get that $u = u_* = g$ in $\R^n\setminus\overline{\Omega}$. Alternatively, the function $g\chi_{\Omega^c}+C\chi_{\bar\Omega}$ for $C$ large enough belongs to $\mathcal{S}$, and therefore, $u = u_* = g$ in $\R^n\setminus\overline{\Omega}$.

In all, $u_*\in \mathcal{S}$ and (since $u_*\le u$) $u_* \equiv u$ in $\R^n$.

\item \label{step:prevh00} Let us now show that $u$ is a subsolution as well. To do so, let us consider its upper semi-continuous envelope in ${\Omega}$,
  \[
 u^*(x) := \sup\left\{\limsup_{k\to \infty} u(x_k) : {\Omega}\ni x_k\to x\right\},
 \]
 with $u^* = u$ in $\R^n\setminus{\Omega}$, and let us show $\I(u^*, x)  \ge 0$ in $\Omega$.
 
 Assume now by contradiction that $u^*$ is not a subsolution, that is, there exist some $x_\circ\in \Omega$ and some test function $\phi\in C^2(B_r(x_\circ))$ for $r > 0$ such that $\phi \ge u^*$ in $\R^n$, $\phi(x_\circ) = u^*(x_\circ)$ but $\I(\phi, x_\circ) < 0$. Arguing as in the first part of the proof of Proposition~\ref{prop:stab_super} (by taking $\phi +\eps|x-x_\circ|^2$ around $x_\circ$ for some small $\eps > 0$, if necessary) we can assume that $u^*-\phi$ has a strict local maximum at $x_\circ$, that is, $u^* > \phi$ in $B_r(x_\circ)\setminus\{x_\circ\}$. Notice, also, that by continuity of $\I(\phi, x)$ around $x_\circ$ (see Corollary~\ref{cor:Iu_II}) we have that $\I(\phi, x) <0$ in $B_\rho(x_\circ)$ for some $\rho > 0$ small.
 
 Let us now consider $\phi_\delta := \phi -\delta$ for some $\delta > 0$. Since $\phi(x) > u^*(x) \ge u(x)$ for $x\in B_r(x_\circ)\setminus\{x_\circ\}$, we have that for $\delta> 0 $ small enough, $\phi_\delta > u$ in $B_r(x_\circ)\setminus B_\rho(x_\circ)$ as well. Let us define 
 \[
 u_\delta = \left\{
 \begin{array}{ll} 
 \min \{u, \phi_\delta\} & \quad\text{in}\quad B_\rho(x_\circ),\\
 u  & \quad\text{in}\quad B^c_\rho(x_\circ).
 \end{array}
 \right.
 \]
 Notice that $u_\delta$ is a supersolution, since it coincides with $u$ in $B_\rho^c(x_\circ)$, and is the infimum of two supersolutions in $B_\rho(x_\circ)$ (recall Remark~\ref{rem:max_min}). This means that $u_\delta\in \mathcal{S}$, and therefore $u_\delta \ge u$. In particular, we have that $\phi-\delta \ge u$ in $B_r(x_\circ)$, and thus $\phi(x_\circ) - \delta \ge u^*(x_\circ)$, a contradiction. That is, $u^*$ is a subsolution. 
 
 But then,   by the comparison principle (Theorem~\ref{thm:comparison_viscosity_2}), since $u$ is a supersolution, $u^*$ is bounded above, and $u = u^* = g$ in $\R^n\setminus \Omega$, we get that $u^* \le u $ in $\Omega$, which means that $u = u^*$. Therefore, $u$ is continuous in ${\Omega}$, and it is both a sub- and supersolution. This concludes the proof. 
 \qedhere
 \end{steps}
 \end{proof}

 In order to prove the existence of viscosity solutions that are continuous  up to the boundary, we will assume that the  bounded domain $\Omega\subset \R^n$ is such that
 \begin{equation}
 \label{eq:domaincondition}
 \begin{split}
 &\text{for every $x_\circ\in \partial\Omega$, there exists a function $\psi_+\in C^2(\Omega)\cap C(\R^n)$}\\
 &\text{such that $\psi_+(x_\circ) = 0$, $\psi_+ > 0$ in $\R^n\setminus \{x_\circ\}$, $\Mp \psi_+ \le -1$ in $\Omega$,}\\
&  \text{and $\psi_+\ge 1$ in $\R^n \setminus B_{1}(x_\circ)$.}
 \end{split}
 \end{equation}

As we will see below (in Lemma~\ref{lem:lip_satisfies}), condition \eqref{eq:domaincondition} is implied by condition \eqref{eq:cond_Om0}. 
For the existence of bounded viscosity solutions in Proposition~\ref{prop:existence_main_bdd}, the requirement on the domain was due to the comparison principle, which was proved by constructing a supersolution that is \emph{singular} on $\partial\Omega$. 
When we want to construct solutions that are continuous, instead, the limiting factor is the existence of a supersolution that \emph{vanishes} on $\partial\Omega$, as in \eqref{eq:domaincondition}.

 In this context, we have:

 \begin{prop}
 \label{prop:existence_main}
Let $s\in (0, 1)$, and let $\Omega$ be any bounded open set satisfying \eqref{eq:domaincondition} for some $\kappa > 0$. Let $\I \in \III$, and let $g\in L^\infty(\R^n)$ be continuous on $\partial\Omega$ with respect to $\R^n\setminus \Omega$, in the sense of \eqref{eq:cont_from_outside}. 
Then, there exists a unique viscosity solution $u\in C(\overline{\Omega})\cap L^\infty(\R^n)$ to
\[
\left\{
\begin{array}{rcll}
\I(u, x) & = & 0& \quad\text{in}\quad \Omega\\
u & = & g& \quad\text{in}\quad \R^n\setminus \Omega. 
\end{array}
\right.
\]
 \end{prop}
 
 \begin{proof} The proof goes in parallel to that of Proposition~\ref{prop:existence_main_bdd}. 
 The uniqueness directly follows from Corollary~\ref{cor:uniqueness_viscosity}. We divide it into three steps: 
 \begin{steps}
 \item  Let 
 \[
 \mathcal{S} := \left\{v \in {\rm LSC}(\overline{\Omega})\cap L^1_{\omega_s}(\R^n): \begin{array}{l} \I(v, x)\le 0 \ \text{ in   $\Omega$  in the viscosity sense} \\
  \hspace{0.93cm}v \ge g \ \text{ in}\ \R^n\setminus\Omega\end{array}\right\},
 \]
 and let us define 
 \[
 u(x) := \inf_{v\in \mathcal{S}} v(x). 
 \]
 As in the proof of Proposition~\ref{prop:existence_main_bdd}, $\mathcal{S}$ is nonempty and $u$ is globally bounded. 
 We define its lower semi-continuous envelope in $\overline{\Omega}$, 
 \[
 u_*(x) := \inf\left\{\liminf_{k\to \infty} u(x_k) : \overline{\Omega}\ni x_k\to x\right\},
 \]
 with $u_* = u$ in $\R^n\setminus\overline{\Omega}$, which by  Lemma~\ref{lem:inf_is_super} satisfies
 \[
 \I(u_*, x) \le 0\quad\text{in}\quad \Omega. 
 \]
 Notice, also, that as before we have $u = u_* = g$ in $\R^n\setminus\overline{\Omega}$. 
 
 \item \label{step:prevh} Let us now show that $u_* = g$ on~$\partial\Omega$ and that $u_*$ is continuous on $\partial \Omega$, that is, for every $x_\circ\in \partial\Omega$ and $x_k \to x_\circ$, then $\liminf_{k\to \infty} u_*(x_k) = \limsup_{k \to \infty} u_*(x_k) = g(x_\circ)$. 
 
  Let $x_\circ\in \partial\Omega$,  $\eps > 0$, and let us define 
 \[
 w_\eps^\pm = g(x_\circ) \pm \eps \pm k_\eps \psi_+,
 \]
 where $\psi_+$ is the function from the condition \eqref{eq:domaincondition} at $x_\circ$ and $k_\eps$ is chosen large enough (depending on $\eps$, but also on $g$, $\Omega$, and $\psi_+$) so that 
 \begin{equation}
 \label{eq:gchosen}
 |g-g(x_\circ)|\le \eps + k_\eps \psi_+ \quad\text{in}\quad \R^n\setminus \Omega. 
 \end{equation}
Hence, we have 
  \[
\begin{split}
 w_\eps^+& \ge g \quad\text{in}\quad \R^n\setminus \Omega, \\
  w_\eps^-& \le g \quad\text{in}\quad \R^n\setminus \Omega.
  \end{split}
 \] 
 By assuming $k_\eps \ge \|\I(0, x)\|_{L^\infty(\R^n)}$, we have from the fact that $\Mp \psi_+ \le  -1$ in $\Omega$,
 \[
 \begin{split}
 \I(w_\eps^-, x) \ge \I(0, x) -k_\eps \Mp (\psi_+) \ge \I(0, x) + k_\eps\ge 0\quad\text{in}\quad \Omega,\\
\I(w_\eps^+, x) \le \I(0, x) +k_\eps \Mp (\psi_+) \le \I(0, x) - k_\eps\le 0\quad\text{in}\quad \Omega.
 \end{split}
 \]
 In particular, $w_\eps^+\in \mathcal{S}$. By continuity of $\psi_+$, for each $\eps > 0$ there exists some $\delta > 0$ such that $w_\eps^+ \le g(x_\circ) +2\eps$ in $B_\delta(x_\circ)$. This yields $u\le w_\eps^+\le g(x_\circ) +2\eps$ in $B_\delta(x_\circ)$, so that if $x_k\to x_\circ$, then 
 \[
 \limsup_{k\to \infty}u(x_k) \le g(x_\circ) + 2\eps. 
 \]

On the other hand, by comparison principle all elements in $\mathcal{S}$ are above $w_\eps^-$ for any $\eps > 0$. Again, by continuity of $\psi_+$, for each $\eps > 0$ there exists some $\delta > 0$ such that $w_\eps^- \ge g(x_\circ) -2\eps$ in $B_\delta(x_\circ)$, and therefore $u_*\ge w_\eps^-\ge  g(x_\circ) -2\eps$ in $B_\delta(x_\circ)$, so that if $x_k\to x_\circ$, 
 \[
 \liminf_{k\to \infty} u(x_k) \ge g(x_\circ) - 2\eps. 
 \]
 Since $\eps > 0$ is arbitrary, we have that if $ x_k\to x_\circ$ then 
 \[
 \lim_{k\to \infty} u_*(x_k ) = g(x_\circ).
 \]
 Hence $u_* = g$ on $\partial\Omega$ as well,  and $u_*$ is continuous on $\partial\Omega$. Therefore, $u_*\in \mathcal{S}$ and (since $u_*\le u$) $u_* \equiv u$ in $\R^n$. In particular, $u\in {\rm LSC}(\overline{\Omega})$ and $\I(u, x) \le 0$ in $\Omega$.  
 \item It remains to see that $u$ is a subsolution as well, by defining first its upper semi-continuous envelope in $\overline{\Omega}$,
  \[
 u^*(x) := \sup\left\{\limsup_{k\to \infty} u(x_k) : \overline{\Omega}\ni x_k\to x\right\},
 \]
 with $u^* = u$ in $\R^n\setminus\overline{\Omega}$, and let us show $\I(u^*, x)  \ge 0$ in $\Omega$.
 
 Observe  that, since by \ref{step:prevh} $u$ is continuous on $\partial\Omega$, we have $u^* = g$ in the whole $\R^n\setminus\Omega$. Now, arguing by contradiction as in \ref{step:prevh00} from Proposition~\ref{prop:existence_main_bdd}, we obtain that $u^*$ is a subsolution. 
 
 But then, again by comparison principle (Theorem~\ref{thm:comparison_viscosity}), since $u$ is a supersolution, and $u = u^* = g$ in $\R^n\setminus \Omega$, we get that $u^* \le u $ in $\Omega$, which means that $u = u^*$. Therefore, $u$ is continuous in $\overline{\Omega}$, and is both a sub- and supersolution. This concludes the proof. 
 \qedhere
 \end{steps}
 \end{proof}

In order to prove our main result, we will also need to show that the domains satisfying \eqref{eq:cond_Om} are such that \eqref{eq:domaincondition} holds:

\begin{lem}
\label{lem:lip_satisfies}
Let $s\in (0, 1)$, and let $\Omega\subset \R^n$ be any bounded open set satisfying \eqref{eq:cond_Om0} for some $\kappa > 0$. Then $\Omega$ satisfies condition \eqref{eq:domaincondition}, and we can take a $\psi_+$ that has a  H\"older  modulus of continuity which is independent of the point $x_\circ\in \partial\Omega$   and $\psi_+(x) \ge \frac12 |x-x_\circ|^2$ in $B_1$.
\end{lem}
\begin{proof}
Let $x_\circ\in \partial\Omega$, and let $\varphi_1 := \dr_\Omega^\eps$ be the function from Lemma~\ref{supersolution-d-eps}. Then, $\varphi_1\ge 0$ satisfies
\[
\Mp \varphi_1 \le -c\quad\text{in}\quad\Omega
\] 
by Lemma~\ref{supersolution-d-eps}, and $\varphi_1\in C^2(\Omega)\cap C(\R^n)$. Let $\varphi_2\in C^\infty(\R^n)$ with $0\le \varphi_2\le 1$ such that $\varphi_2 \equiv 1$ in $\R^n\setminus B_1$, $\varphi_2(0) = 0$, $\varphi(x) > 0$ for $x\neq 0$, and $\varphi(x) \ge \frac12 |x|^2$ in $B_1$. In particular, $\Mp \varphi_2$ is globally bounded by a constant $C$ and if we take 
\[
\psi_+ := C_1\varphi_1+\varphi_2(\,\cdot\, + x_\circ)
\]
then 
\[
\Mp \psi_+ \le - C_1c + C \quad\text{in}\quad\Omega.
\]
By choosing $C_1$ large enough, $\psi_+$ satisfies all the conditions in \eqref{eq:domaincondition}. 
\end{proof}
Hence, we have all the ingredients to show the existence of viscosity solutions in Lipschitz (and more general) domains:
\begin{proof}[Proof of Theorem~\ref{thm:existence_visc}]
Uniqueness follows from Corollary~\ref{cor:uniqueness_viscosity_2}. 

If $g$ is globally bounded, then we are done by Propositions~\ref{prop:existence_main_bdd} or \ref{prop:existence_main} and Lemma~\ref{lem:lip_satisfies}. Otherwise, since $g$ is  bounded in a $\rho$-neighborhood of $\partial\Omega$, we denote
\[
\Omega_\rho := \left\{x\in \R^n\setminus \Omega : \dist(x, \partial\Omega) \le \rho/2\right\},
\]
and we consider $g_\rho = g\chi_{\Omega_\rho}$, and $\bar g_\rho := g-g_\rho$, which is defined in $\R^n$ (extended by zero inside $\Omega$). If $\I$ is of the form 
\[
\I (u, x) = \inf_{b\in \B}\sup_{a\in \A}\big\{-\L_{ab} u(x) + c_{ab}(x)\big\}, 
\]
we define $\I_\rho$ as
\[
\I_\rho (v, x) = \inf_{b\in \B}\sup_{a\in \A}\big\{-\L_{ab} v(x) - \L_{ab} \bar g_\rho(x)+ c_{ab}(x)\big\},\qquad \L_{ab}\in \LLL.
\]
Notice that by Lemma~\ref{lem:Lu_LL}, $\L_{ab} \bar g_\rho$ is bounded in $\Omega_\rho$ and continuous, with a modulus of continuity independent of $\L_{ab}$. Thus,  $\I_\rho\in \III$ as well, with 
\begin{equation}
\label{eq:boundIrho}
\|\I_\rho(0, x)\|_{L^\infty(\R^n)}\le \|\I(0, x)\|_{L^\infty(\R^n)}+C \|g\|_{L^1_{\omega_s}(\R^n)}, 
\end{equation}
for some $C$ depending only on $n$, $s$, $\Lambda$, and $\rho$. 
Let now $\bar u$ be the (unique) solution to 
\[
\left\{
\begin{array}{rcll}
\I_\rho(\bar u, x) & = & 0& \quad\text{in}\quad \Omega\\
\bar u & = & g_\rho& \quad\text{in}\quad \R^n\setminus \Omega, 
\end{array}
\right.
\]
which exists by Proposition~\ref{prop:existence_main_bdd}, since $g_\rho$ is bounded near $\partial\Omega$, and it is continuous on $\partial\Omega$ if $g$ is continuous on $\partial\Omega$, by Proposition~\ref{prop:existence_main} and Lemma~\ref{lem:lip_satisfies}. Then, the function $u = \bar u + \bar g_\rho$ satisfies, by construction of $\I_\rho$ and $g_\rho$,
\[
\left\{
\begin{array}{rcll}
\I(u, x) & = & 0& \quad\text{in}\quad \Omega\\
u & = & g& \quad\text{in}\quad \R^n\setminus \Omega.
\end{array}
\right.
\]
 This proves the desired result.
\end{proof}

Finally, we mention that, as a consequence of the proof of existence in Proposition~\ref{prop:existence_main}, we actually have the following boundary continuity estimate (used for the boundary regularity in Corollary~\ref{cor:bdry_reg_int} later on):

\begin{cor}
\label{cor:boundaryregularity_visc}
Let $s\in (0, 1)$, and let $\Omega$ be any bounded open set satisfying \eqref{eq:cond_Om0}. Let $\I \in \III$, and let $g\in L^1_{\omega_s}(\R^n)$ be continuous on $\partial\Omega$ with respect to $\R^n\setminus \Omega$ (with modulus $\sigma$, in the sense of \eqref{eq:cont_from_outside}). Then, the unique viscosity solution $u\in C(\overline{\Omega})\cap L^1_{\omega_s}(\R^n)$ of
\[
\left\{
\begin{array}{rcll}
\I(u, x) & = & 0& \quad\text{in}\quad \Omega\\
u & = & g& \quad\text{in}\quad \R^n\setminus \Omega
\end{array}
\right.
\]
satisfies
\[
|u(x) - u(y)|\le \omega(|x-y|)\quad\text{for all}\quad x\in \partial\Omega,\, y \in \R^n,
\]
for some modulus of continuity $\omega$ that depends only on $n$, $s$, $\lambda$, $\Lambda$, $\sigma$, $\|g\|_{L^1_{\omega_s}(\R^n)}$, $\|g\|_{L^\infty(\partial\Omega)}$, $\|\I(0, x)\|_{L^\infty(\R^n)}$, and $\Omega$, but is independent of $u$. Moreover, if $\sigma$ is $\alpha$-H\"older continuous, then $\omega$ is $\bar \alpha$-H\"older  continuous as well, where $\bar \alpha$ depends only on $\alpha$, $\Omega$, $n$, $s$, $\lambda$, and $\Lambda$. 
\end{cor}
\begin{proof}
Let us first assume that $g$ is bounded. Recall that in \ref{step:prevh} of the proof of Proposition~\ref{prop:existence_main} we showed that for any $\eps$ and $x_\circ$ there exists $\delta > 0$ such that 
\[
g(x_\circ) -2\eps \le u(x_\circ) \le g(x_\circ) +2\eps,
\]
and $\delta$ depends only on the choice of $k_\eps$ and the modulus of continuity of $\psi_+$. By Lemma~\ref{lem:lip_satisfies}, $\psi_+$ is H\"older continuous with some exponent $\beta > 0$, depending only on $\Omega$ (and $n$, $s$, $\lambda$, and $\Lambda$). On the other hand, $k_\eps$ in the proof of Proposition~\ref{prop:existence_main} is chosen so that \eqref{eq:gchosen} holds, so that its value is determined  only by $\eps$, the modulus of continuity of $g$ around $x_\circ$ (hence,~$\sigma$), $\psi_+$, and an upper bound for $\|g\|_{L^\infty(\R^n)}$ and $\|\I(0, x)\|_{L^\infty(\R^n)}$. This shows the existence of $\omega$, with an added dependence on $\|g\|_{L^\infty(\R^n)}$. 

Notice, moreover, that when $\sigma$ is $\alpha$-H\"older, since $\psi_+(x)\ge \frac12|x-x_\circ|^2$ (see  Lemma~\ref{lem:lip_satisfies}), in order to ensure that \eqref{eq:gchosen} holds it is enough to choose $k_\eps \ge C \eps^{-2/\alpha}$. Hence, we can choose $\delta = C \eps^{\frac{2+\alpha}{\alpha\beta}}$ and $u$ is H\"older continuous, with exponent $\frac{\alpha\beta}{2+\alpha}$ (thus,  $\omega$ is H\"older).

Finally, for the general case (and in order to remove the last dependence on $\|g\|_{L^\infty(\R^n)}$), we just notice that by the proof of Theorem~\ref{thm:existence_visc} the modulus $\omega$ now depends only  on $\eps$, $\sigma$, $\Omega$, and an upper bound for $\|g_\rho\|_{L^\infty(\R^n)}$ and $\|\I_\rho(0, x)\|_{L^\infty(\R^n)}$. Note, however, that $\|g_\rho\|_{L^\infty(\R^n)}$ is controlled by both $\sigma$ and $\|g \|_{L^\infty(\partial\Omega)}$, whereas $\|\I_\rho(0, x)\|_{L^\infty(\R^n)}$ is controlled (see \eqref{eq:boundIrho}) by $\|\I(0, x)\|_{L^\infty(\R^n)}$ and $C\|g\|_{L^1_{\omega_s}(\R^n\setminus\Omega)}$. This shows the existence of the modulus $\omega$ in the general case, and thus  we are done.  
\end{proof}

\section{Harnack's inequality and H\"older estimates}

\index{Krylov-Safonov}

In this section, we prove Harnack's inequality and H\"older estimates for integro-differential elliptic equations in non-divergence form with bounded measurable coefficients.
This is the nonlocal analogue of the Krylov--Safonov theorem for second-order elliptic PDE \cite{KS79,CC}.

Namely, if $\Mpm$ are given by \eqref{eq:MMpm}, we will study solutions to equations of the form 
\begin{equation}
\label{eq:bdd_meas_coeff}
\left\{
\begin{array}{rcll}
\Mp u &\geq &-C_\circ&\quad \textrm{in} \quad  B_1,\\
\Mm u &\leq &C_\circ&\quad \textrm{in} \quad  B_1,
\end{array}
\right.
\end{equation}
for some $C_\circ\ge 0$. 
The expression \eqref{eq:bdd_meas_coeff} is satisfied, for example, by any solution to $\I(u, x) = 0$ in $B_1$ for some $\I\in \III$ (with $C_\circ = \|\I(0, x)\|_{L^\infty(B_1)}$). 
Moreover, if $\J\in \III$ is translation invariant and $v$ solves 
\[\J v = 0 \quad\textrm{in}\quad B_1,\]
 then the incremental quotients of $v$ satisfy
\[
\mathcal M^+ \left(\frac{v(x+h)-v(x)}{|h|}\right) \ge  0 \ge \mathcal M^- \left(\frac{v(x+h)-v(x)}{|h|}\right) \quad \textrm{in}\quad B_{1-|h|},
\]
which is \eqref{eq:bdd_meas_coeff} with $C_\circ = 0$; see \eqref{eq:extremal_eq} and Proposition \ref{prop:viscosity_ellipticity}.
As a consequence, in case that $v\in C^1$, we can take $h\to0$ to find that derivatives of $v$ also satisfy \eqref{eq:bdd_meas_coeff}.

When $u$ satisfies \eqref{eq:bdd_meas_coeff}, it is said that $u$ satisfies a  \emph{non-divergence-form equation with bounded measurable coefficients}. 
This is because, if $u\in C^2(B_1)$ is a strong solution to \eqref{eq:bdd_meas_coeff}, one can show that (by definition of $\mathcal M^\pm$) for each $x\in B_1$ there exists some $\L_x\in \LLL$ such that 
\[
\L_x u(x) = f(x),
\]
with $\|f\|_{L^\infty(B_1)}\le C_\circ$. Here, the operators $\L_x$ are of the form 
\begin{equation}
\label{eq:Lx1}
\L_x u(x) = \frac12 \int_{\R^n}\big(2u(x)-u(x+y)-u(x-y)\big)K(x, y)\,dy,
\end{equation}
where the kernel $K(x, y)$ satisfies
\begin{equation}
\label{eq:Lx2}
0< \frac{\lambda}{|y|^{n+2s}}\le K(x, y) \le \frac{\Lambda}{|y|^{n+2s}}\qquad\text{for all}\quad   x, y\in \R^n,
\end{equation}
but   has no regularity  on the $x$ variable (apart from being uniformly bounded in $x$ for each $y$ fixed, between $\lambda|y|^{-n-2s}$ and $\Lambda|y|^{-n-2s}$). Conversely, if $u\in C^2(B_1)$ is such that it satisfies 
\[
\L_x u(x) = f(x)\quad\text{in}\quad B_1,
\]
for some $f$ bounded and $\L_x$ of the form \eqref{eq:Lx1}-\eqref{eq:Lx2}, then it clearly satisfies \eqref{eq:bdd_meas_coeff} pointwise in $B_1$, with $C_\circ = \|f\|_{L^\infty(B_1)}$ (again, by definition of $\mathcal M^\pm$).

\subsection{Weak Harnack inequality for supersolutions}

The weak Harnack inequality is a key tool in the study of non-divergence-form equations.
For second-order PDE ($s=1$), this is the key step in the Krylov--Safonov theorem, and its proof is actually quite delicate; see \cite{CC}.

Quite surprisingly, in case of integro-differential operators with kernels satisfying \eqref{eq:Lx2}, its proof is much easier, and the conclusion is even stronger\footnote{More precisely, the term on the left-hand side of the estimate in Theorem~\ref{half-Harnack-sup} is a global $L^1_{\omega_s}$ norm of the solution, while in the local case this would be a local $L^1$ norm in $B_{1/2}$. Of course, this comes at a price, which is that the constant $C$ blows-up as $s\uparrow 1$.}.
The following proof originates from \cite{Sil06}.

\begin{thm} \label{half-Harnack-sup} \index{Half-Harnack!Fully nonlinear!Supersolutions}
Let $s\in(0,1)$, and let $\Mm$ be given by \eqref{eq:MMpm}. Assume that $u\in {\rm LSC}(B_1)\cap L^1_{\omega_s}(\R^n)$ satisfies
\[
\left\{
\begin{array}{rcll}
\Mm u &\leq &C_\circ&\quad \textrm{in} \quad  B_1\\
u &\geq &0&\quad \textrm{in} \quad \R^n.
\end{array}
\right.
\]
in the viscosity sense, for some~$C_\circ\geq0$. Then,
\[\|u\|_{L^1_{\omega_s}(\R^n)}\leq C\left(\inf_{B_{1/2}}u+C_\circ\right).\]
The constant $C$ depends only on $n$, $s$, $\lambda$,  and $\Lambda$. 
\end{thm}

\begin{proof}
Let $\eta\in C^\infty_c(B_{3/4})$ be such that $0\leq \eta\leq 1$ and $\eta\equiv1$ in $B_{1/2}$. Let 
\begin{equation}\label{ufgeruerug2}
t := \max\{\tau \ge 0 : u \ge \tau \eta\quad\text{in}\ \R^n\} \le \inf_{B_{1/2}} u.
\end{equation}
Since $u$ and $\eta$ are continuous in $B_1$, there exists $x_\circ\in B_{3/4}$ such that $u(x_\circ)=t\eta(x_\circ)$. That is, $u$ can be touched from below by a smooth function, and in particular we can evaluate $\Mm u$ pointwise at $x_\circ$ (see Remark~\ref{rem:onesidedcondition2}). By ellipticity  (recall Proposition~\ref{prop:viscosity_ellipticity}) we have
\begin{equation}\label{ufgeruerug}
\Mm (u-t\eta)(x_\circ) \leq \Mm u(x_\circ)-t\,\Mm\eta (x_\circ) \leq C_\circ+C_1t,
\end{equation}
 for some~$C_1>0$. Furthermore, since $u-t\eta\geq0$ in $\R^n$ and $(u-t\eta)(x_\circ)=0$,  we can also bound pointwise
\[\begin{split}
\Mm (u-t\eta)(x_\circ)& \geq \lambda\int_{\R^n}\frac{u(z)-t\eta(z)}{|x_\circ-z|^{n+2s}}dz
\\
& \geq c\int_{\R^n}\frac{u(z)-t\eta(z)}{1+|z|^{n+2s}}\,dz
\\
& \geq c\int_{\R^n}\frac{u(z)}{1+|z|^{n+2s}}dz-C_2t,\end{split}\]
for some~$c$, $C_2>0$.
From this and \eqref{ufgeruerug2}-\eqref{ufgeruerug}, we obtain that
\[ \inf_{B_{1/2}} u\geq t\geq -c_1C_\circ+c_2\int_{\R^n}\frac{u(z)}{1+|z|^{n+2s}}dz,\]
for some~$c_1$, $c_2>0$, as desired.
\end{proof}

\subsection{Iteration and H\"older estimate}

Let us now show how to iterate the weak Harnack inequality (Theorem \ref{half-Harnack-sup}) to get a H\"older estimate for solutions to nonlocal equations.

Notice that one needs to be careful when doing this since, contrary to what happens for second-order PDE, in the present setting the weak Harnack inequality requires $u\geq0$ in the full space.

\begin{thm}\label{C^alpha-bmc} \index{Holder estimate@H\"older estimate!Viscosity solutions!Bounded}
Let $s\in(0,1)$ and let $\Mpm $ be given by \eqref{eq:MMpm}.

Let $u\in C(B_1)\cap L^\infty(\R^n)$   be any viscosity solution to  a non-divergence-form equation with bounded measurable coefficients, i.e.,
\[\left\{
\begin{array}{rcll}
\Mp u &\geq & -C_\circ &\quad \textrm{in}\quad B_1\\
\Mm u &\leq & C_\circ &\quad \textrm{in}\quad B_1,
\end{array}
\right.\]
for some~$C_\circ\geq0$.
Then $u\in C^\gamma_{\rm loc}(B_1)$ with 
\[\|u\|_{C^{\gamma}(B_{1/2})}\leq C\left(\|u\|_{L^\infty(\R^n)}+C_\circ\right),\]
where $C$ and $\gamma>0$ depend only on $n$, $s$,  $\lambda$, and $\Lambda$.
\end{thm}

The proof of this result will be based on the following version of the weak Harnack inequality.

\begin{lem}\label{weak-Harnack-v2}
Let $s\in(0,1)$, $\alpha \in (0, 2s)$, and let $\Mm$ be given by \eqref{eq:MMpm}.
There exist $\eps_\circ> 0$ depending only on $n$, $s$, $\lambda$, and $\Lambda$,  and $R_\circ>1$ depending only on $n$, $s$, $\alpha$, $\lambda$, and $\Lambda$, such that the following holds for any $\eps \leq \eps_\circ$ and $R \geq R_\circ$.

Let $u\in C(B_2)\cap L^1_{\omega_s}(\R^n)$ satisfy
\[
\left\{
\begin{array}{rcll}
\Mm u &\leq& \varepsilon&\quad \textrm{in}\quad B_2,\\
u &\geq& 0&\quad \textrm{in}\quad B_R,
\end{array}
\right.
\]
in the viscosity sense, and 
\[
\int_{B_R} \frac{u(x)}{1+|x|^{n+2s}}\,dx \ge \frac{|B_1|}{2}.
\]
Assume in addition that $u$ has controlled growth,
\[|u(x)|\leq 2+|x|^{\alpha}\quad \textrm{for}\quad |x|\geq R.\]

Then, 
\[\inf_{B_1}u \geq \theta>0,\]
with $\theta$ depending only on $n$, $s$, $\lambda$, and $\Lambda$.
\end{lem}
\begin{proof}
We have, by Proposition~\ref{prop:viscosity_ellipticity},
\[
\Mm(u\chi_{B_R}) = \Mm(u - u\chi_{B_R^c}) \le \Mm u - \Mm(u\chi_{B_R^c}).
\]

Observe that, for any $x\in B_2$, if $R > 3$
\[
\Mp (u\chi_{B_R^c})(x) \le C\int_{B_R^c}\frac{|u(y)|}{|x-y|^{n+2s}}\, dy  \le C_\alpha \int_{B_R^c}|y|^{\alpha-n-2s} dy\le C_\alpha R^{\alpha-2s},
\]
for some $C_\alpha$ that depends only on $n$, $s$, $\alpha$, and $\Lambda$, and hence
\[
\Mm(u\chi_{B_R}) \leq \Mm u + \Mp(u\chi_{B_R^c})
 \le \eps_\circ + C_\alpha R^{\alpha -2s}\quad\text{in}\quad B_2. 
\]
We then apply Theorem \ref{half-Harnack-sup} (after a rescaling by a factor 2) to the function $u\chi_{B_R}$, to obtain, 
\[
\inf_{B_1} u \ge \frac{1}{C} \int_{B_R} \frac{u(x)}{1+|x|^{n+2s}}\,dx - \eps_\circ  -C_\alpha R^{\alpha -2s}\ge \frac{|B_1|}{2C}  - \eps_\circ  -C_\alpha R_\circ^{\alpha -2s}.
\]
Choose $\eps_\circ\le \frac{|B_1|}{8C}$  and  $R_\circ>3$ large enough so that $C_\alpha R_\circ^{\alpha-2s}\le \frac{|B_1|}{8C} $, to get the result with $\theta = \frac{|B_1|}{4C}$. 
\end{proof}

As a consequence of the previous lemma, we can prove an oscillation decay result.

\begin{lem}\label{oscillat-decay2}
Let $s\in(0,1)$, $\alpha\in (0, 2s)$, and let $\Mpm$ be given by \eqref{eq:MMpm}. Let $R_\circ>1$ and $\varepsilon_\circ>0$ be given by Lemma \ref{weak-Harnack-v2}. 
Then, there exists $R\geq R_\circ$ large enough depending only on $n$, $s$, $\alpha$, $\lambda$, and $\Lambda$, such that the following statement holds. 

Let $u\in C(B_R)\cap L^1_{\omega_s}(\R^n)$ satisfy
\[
\left\{
\begin{array}{rcll}
\Mp u &\geq& -\varepsilon_\circ&\quad \textrm{in}\quad B_R,\\
\Mm u &\leq& \eps_\circ &\quad \textrm{in}\quad B_R,
\end{array}
\right.
\]
with
\[|u|\leq 1\quad \textrm{in}\quad B_R\]
and
\[|u(x)|\leq 1+|x|^{\alpha}\quad \textrm{for}\quad |x|\geq R.\]

Then, there exists $\gamma>0$ such that
\[\osc_{B_{R^{-k}}} u \leq CR^{-k\gamma}\qquad\text{for all}\quad k \in \N,\]
for some constants $\gamma$ and $C$ depending only on $n$, $s$, $\alpha$, $\lambda$, and $\Lambda$
\end{lem}

\begin{proof}
We split the proof into two steps.

\begin{steps}
\item \label{it:step1eitheror}
We first show that 
\begin{equation}
\label{eq:eitherorlemma}
\text{either   \quad $-1\leq u\leq 1-2\theta$ \ in \ $B_1$\qquad or \quad $2\theta-1\leq u\leq 1$ \ in \ $B_1$.}
\end{equation}
with $\theta>0$ depending only on $n$, $s$, $\lambda$, and $\Lambda$.

For this, we may assume $\int_{B_R}\frac{u(x)}{1+|x|^{n+2s}}\,dx \leq 0$ (otherwise take $u\mapsto -u$), and define
\[w:=1-u.\]
Then, we have
\[
\left\{
\begin{array}{rcll}
\Mp w &\geq& -\varepsilon_\circ&\quad \textrm{in}\quad B_R\\
\Mm w &\leq& \eps_\circ &\quad \textrm{in}\quad B_R\\
w & \ge & 0 & \quad \textrm{in}\quad B_R,
\end{array}
\right.
\]
and
\[ \int_{B_R} \frac{w(x)}{1+|x|^{n+2s}}\,dx \geq \int_{B_R} \frac{dx}{1+|x|^{n+2s}} \geq \frac{|B_1|}{2}.\]
Moreover, 
\[|w(x)|\leq 2+|x|^\alpha\quad \textrm{for}\quad |x|\geq R.\]
Hence, by Lemma \ref{weak-Harnack-v2}, since $R\geq R_\circ$ we deduce that
\[\inf_{B_1} w \geq 2\theta>0,\]
or equivalently, $u\leq 1-2\theta$ in $B_1$, for some $\theta$ depending only on $n$, $s$, $\lambda$, and~$\Lambda$ (half the value obtained from Lemma~\ref{weak-Harnack-v2}). If we had taken $-u$ instead of $u$, we would have obtained $u\geq 2\theta-1$ in $B_1$, thus giving \eqref{eq:eitherorlemma}. That is, 
\[
\osc_{B_1} u = 2(1-\theta). 
\]

\item 
We next show how to iterate \ref{it:step1eitheror}, \eqref{eq:eitherorlemma}, in order to get the desired result.
Notice that since $|u\pm\theta| \leq 1-\theta$ in $B_1$, the function
\[u_R(x):=\frac{u\left(\frac{x}{R}\right) \pm \theta}{1-\theta}\]
satisfies
\[|u_R|\leq 1\quad \textrm{in}\quad B_R,\]
and solves the equation
\[\left\{\begin{array}{rcccll}
\Mp u_R &\geq & -\frac{R^{-2s} \varepsilon_\circ}{1-\theta} &\geq &-\varepsilon_\circ & \quad \textrm{in}\quad B_R, \\
\Mm u_R& \leq & \frac{R^{-2s} \varepsilon_\circ}{1-\theta} &\leq& \varepsilon_\circ & \quad \textrm{in}\quad B_R,
\end{array}\right.\]
provided that $R$ is large enough. 
In addition, outside $B_R$ it satisfies
\[|u_R(x)| \leq \left\{\begin{array}{ll}
\displaystyle{\frac{1+\theta}{1-\theta}} & \quad \textrm{if}\quad R\leq |x|\leq R^2, \\[0.35cm]
{\displaystyle{\frac{\left|\frac{x}{R}\right|^\alpha + \theta}{1-\theta}}} & \quad \textrm{if}\quad |x|\geq R^2,
\end{array}
\right.
\]
and hence
\[
|u_R(x)|  \leq \left|\frac{x}{R}\right|^\alpha \frac{1+\theta}{1-\theta}\leq |x|^\alpha\quad\text{for all}\quad |x|\ge R,\]
if $R$ is large enough so that $R^\alpha \geq \frac{1+\theta}{1-\theta}$.

Thus, $u_R$ satisfies again the hypotheses of \ref{it:step1eitheror}, and therefore
\[
|u_R\pm \theta|\le 1-\theta\quad\text{in}\quad B_1,
\]
that is, 
\[
\osc_{B_{R^{-1}}} u \le 2 (1-\theta)^2.
\]
Iterating the procedure $k $ times we get
\[ \osc_{B_{R^{-k}}} u \leq 2(1-\theta)^k \leq 2R^{-k\gamma} \quad\text{for all}\quad k\in \N,\]
where $\gamma>0$ is chosen so that $1-\theta= R^{-\gamma}$.\qedhere
\end{steps}
\end{proof}

We can finally deduce the desired H\"older  estimate.

\begin{proof}[Proof of Theorem \ref{C^alpha-bmc}]
Dividing $u$ by a constant if necessary we may assume $\|u\|_{L^\infty(\R^n)}\leq 1$ and $C_\circ\leq \varepsilon_\circ$, where $\eps_\circ$ is given by Lemma~\ref{weak-Harnack-v2}.
Then, for any $x_\circ\in B_{1/2}$ we define
\[\tilde u(x):= u\left({\textstyle \frac{x_\circ+x}{2R}}\right),\]
where $R>1$ is given by Lemma \ref{oscillat-decay2}.
The function $\tilde u$ satisfies the hypotheses of Lemma \ref{oscillat-decay2}, and therefore 
\[\osc_{B_{R^{-k}}} \tilde u \leq C_1 R^{-k\gamma} \quad \textrm{for all} \quad k\geq0.\]
This yields 
\[\|u-u(x_\circ)\|_{L^\infty(B_{R^{-k}}(x_\circ))} \leq \osc_{B_{R^{-k}}(x_\circ)} u \le \osc_{B_{R^{-k+2}}} \tilde u \leq C_2 R^{-k\gamma},\] 
which implies
\[\big|u(x)-u(x_\circ)\big| \leq C_3|x-x_\circ|^\gamma\]
for all $x\in B_1$, and we are done.
\end{proof}

In the next subsection we will prove that we can actually replace the $L^\infty (\R^n)$ norm of $u$ by its $L^1_{\omega_s}(\R^n)$ norm.

\subsection{$L^\infty$ bound for subsolutions}

We next prove the other half of Harnack's inequality, which reads as follows.
This was first proved in \cite{CS2}, and the proof we present here is from \cite{DRSV}.

\begin{thm}\label{half-Harnack-sub} \index{Half-Harnack!Fully nonlinear!Subsolutions} \index{Linfty bounds@$L^\infty$ bounds!Fully nonlinear}
Let $s\in(0,1)$ and let $\Mp$ be given by \eqref{eq:MMpm}. Assume that $u\in {\rm USC}(B_1)\cap L^1_{\omega_s}(\R^n)$ satisfies
\[\Mp u\geq -C_\circ\quad \textrm{in}\quad B_1\]
in the viscosity sense, for some~$C_\circ\geq0$.
Then,
\[\sup_{B_{1/2}}u\leq C\left(\|u\|_{L^1_{\omega_s}(\R^n)}+C_\circ\right).\]
The constant $C$ depends only on $n$, $s$, $\lambda$, and $\Lambda$.
\end{thm}

\begin{proof} 
We divide the proof into three steps. 
\begin{steps}
\item Let us start with a series of reductions. We may first assume that $C_\circ=0$: otherwise we consider $\tilde u:=u-C_1C_0\eta$ for some $C_1  > 0$, with $\eta\in C^\infty_c(B_2)$, $\eta\equiv1$ in $B_1$, $0\le \eta\le 1$ in $\R^n$.
Then, we have that $\Mp \eta \leq -c<0$ in $B_1$, and hence 
\[
\Mp\tilde u \geq -C_\circ -C_1 C_\circ \Mp \eta\geq 0\quad\text{in}\quad B_1, 
\]
provided that $C_1$ is large enough. If we now show the result for $\tilde u$, we would have
\[\sup_{B_{1/2}}u\le \sup_{B_{1/2}}\tilde u + C_1C_\circ \leq C\left(\|u\|_{L^1_{\omega_s}(\R^n)}+C_1C_\circ \|\eta\|_{L^1_{\omega_s}(\R^n)}+C_\circ\right),\]
which also gives the result for $u$ (since $\|\eta\|_{L^1_{\omega_s}(\R^n)}$ is bounded universally). 

We may further assume that~$u\geq0$; otherwise we can consider $u_+:=\max\{u,0\}$ instead of $u$, which satisfies $\Mp u_+ \ge 0$ in the viscosity sense in~$\R^n$ (see Remark~\ref{rem:max_min}). Proving the result for $u_+$ would yield
\[\sup_{B_{1/2}}u \le \sup_{B_{1/2}}u_+ \leq C\left(\|u_+\|_{L^1_{\omega_s}(\R^n)}+C_\circ\right)\le C\left(\|u\|_{L^1_{\omega_s}(\R^n)}+C_\circ\right),\]
which is the result for $u$ as well. 

Finally, observe that, after dividing $u$ by a constant,
we may assume in addition that 
\begin{equation}\label{soigjer}
\|u\|_{L^1_{\omega_s}(\R^n)} = \int_{\R^n} \frac{u(x)}{1+|x|^{n+2s}}\,dx=1.
\end{equation}

In this setting, it suffices to prove that if $u \in{\rm  USC}(B_1)\cap L^1_{\omega_s}(\R^n)$ satisfies
\begin{equation}
\label{eq:u_USC_satisfies}
\left\{
\begin{array}{rcll}
\Mp u & \ge & 0 & \quad \text{in}\quad B_1,\\
u & \ge & 0 & \quad\text{in}\quad \R^n,
\end{array}
\right.\qquad\text{and}\qquad \|u\|_{L^1_{\omega_s}(\R^n)}= 1,
\end{equation}
then
\[
u(0) \leq C
\]
for some constant $C$ depending only on $n$, $s$, $\lambda$, and $\Lambda$. 
Once this is proved, we simply apply this to every point in $B_{1/2}$ (after scaling and translating), and the result follows.

\item \label{it:step2claim} Let us show the following claim:

\noindent \textbf{Claim}. There exists $\delta$ and $c_\circ$ depending only on $n$, $s$, $\lambda$, and $\Lambda$, such that if $u$ satisfies \eqref{eq:u_USC_satisfies}, $u(x_\circ)\ge M$ for $x_\circ\in B_{1/2}$, and $(c_\circ M)^{-1/n} <\frac12$, then 
\[\sup_{B_{{r_\circ}}(x_\circ)} u > (1+\delta) M, \quad \textrm{for} \quad {r_\circ}:= (c_\circ M)^{-1/n} <{\textstyle \frac12}.\]

To prove the claim,   we argue by contradiction. That is, let us suppose that
\[\sup_{B_{{r_\circ}}(x_\circ)} u \le (1+\delta)M.\]
We will reach the contradiction by using Theorem~\ref{half-Harnack-sup} for an auxiliary function.
 
 Consider
\[v(x) := (1+\delta)M - u\big(x_\circ+({r_\circ}/2) x\big).\]
Then, $v_\pm\geq0$ everywhere, and
\begin{equation}\label{alsoriyhur}
v_+(x)\equiv (1+\delta)M - u\big(x_\circ+({r_\circ}/2)x\big) \quad{\mbox{ for all}}\quad x\in B_2.
\end{equation}
Since $\Mp u\geq 0$ in $B_1$ by assumption, and using Proposition~\ref{prop:viscosity_ellipticity},
\begin{equation}
\label{eq:vintheviscositysense}
\Mm v_+  - \Mp v_-\le \Mm v = -\Mp \big(u(x_\circ+(r_\circ/2)\,\cdot\,)\big) \le 0\quad\text{in}\quad B_1
\end{equation}
in the viscosity sense. On the other hand, we have,  for every~$x\in B_1$ (notice that $v_- \equiv 0$ in $B_2$),
\[\begin{split}
\Mp v_- (x)& \leq  \Lambda \int_{\R^n\setminus B_2}v_-(y)\frac{dy}{|x-y|^{n+2s}}\\
  & \leq  \Lambda \int_{\R^n\setminus B_2}u(x_\circ + ({r_\circ}/2) x)\frac{dy}{|x-y|^{n+2s}} \quad \textrm{in}\quad B_1,
\end{split}\]
where we have used that $v_-(x) \le u(x_\circ+({r_\circ}/2) x)$ since $u\ge 0$ in $\R^n$. 
After changing variables we then have
\begin{equation}\begin{split}\label{oet4ighty387}
\Mp v_- (x)
&\leq \frac{\Lambda\, 2^n}{r_\circ^n} \int_{{\R^n\setminus B_{{r_\circ}}(x_\circ)}} 
\frac{u(z)}{\left|x-\frac2{{r_\circ}}(z-x_\circ)\right|^{n+2s}}\,dz \quad \textrm{in}\quad B_1.
\end{split}
\end{equation}
Let us see that, for every~$z\in \R^n\setminus B_{{r_\circ}}(x_\circ)$,
\begin{equation}\label{bv46b754vb4}
\left|x-\frac2{{r_\circ}}(z-x_\circ)\right|\geq \frac2{{r_\circ}}|z-x_\circ|-|x|\ge \frac14 (1+|z|).
\end{equation}
Indeed, if~$z\in B_2\setminus B_{{r_\circ}}(x_\circ)$ then
\[ \frac2{{r_\circ}}|z-x_\circ|-|x|\geq 2-|x|
\geq\frac12+\frac{|z|}4\geq \frac14(1+|z|).
\]
If instead~$z\in\R^n\setminus B_2$,
then we have~$|z-x_\circ|\ge \frac34 |z|$ and therefore
\[
 \frac2{{r_\circ}}|z-x_\circ|-|x|\geq
\frac3{2{r_\circ}}|z|-\frac12|z|\geq \frac52|z|\geq 1+|z|,
\]
so that \eqref{bv46b754vb4} always holds. 

As a consequence of~\eqref{bv46b754vb4} (recall \eqref{oet4ighty387} and $\|u\|_{L^1_{\omega_s}(\R^n)}= 1$),
we obtain that for any $x\in B_1$,
\[
\Mp v_- (x)\le  \frac{\Lambda 2^n}{r_\circ^n} \int_{{\R^n\setminus B_{{r_\circ}}(x_\circ)}} 
\frac{u(z)}{\left|x-\frac2{{r_\circ}}(z-x_\circ)\right|^{n+2s}}\,dz\leq  c_\circ M {\Lambda} C_1,
\]
for some~$C_1>0$ depending only on $n$ and $s$. That is, from \eqref{eq:vintheviscositysense},
\[
\Mm v_+ \le c_\circ M\Lambda C_1\quad\text{in}\quad B_1,
\]
in the viscosity sense.

Notice that $v_+(0)=(1+\delta)M-u(x_\circ)<\delta M$.
Using now Theorem~\ref{half-Harnack-sup}, we deduce  
\[
\ave_{B_2} v_+(x) \,dx\leq C\big(v_+(0)+c_\circ M\Lambda C_1\big)\leq \frac{M}{2},
\]
provided that  $\delta$ and $c_\circ$ are sufficiently small, depending only on $n$, $s$, $\lambda$, and~$\Lambda$. 
Equivalently, recalling~\eqref{alsoriyhur}
and employing the change of variable~$y:=x_\circ+({r_\circ}/2)x$, this means that
\[\ave_{B_{{r_\circ}}(x_\circ)} \big((1+\delta)M - u(y)\big)\,dy \leq \frac{M}{2},\]
which gives 
$$ \ave_{B_{{r_\circ}}(x_\circ)} u(y)\,dy\geq (1+\delta)M
-\frac{M}2\geq \frac{M}2.$$
Since~$B_{{r_\circ}}(x_\circ) \subset B_1$ and   $1+|x|^{n+2s} \le 2$  in $B_1$, we finally obtain that for some $c$ depending only on $n$, 
\[
\begin{split}
\frac{M}{2} & \leq c{r_\circ}^{-n}\int_{B_{{r_\circ}}(x_\circ)} u(x)\,dx \\
& \leq  2 c{r_\circ}^{-n} \int_{\R^n} \frac{u(x)}{1+|x|^{n+2s}}\,dx = 2c_\circ c M,\end{split}\]
where we used \eqref{soigjer} and the definition of~${r_\circ}$.
By choosing $c_\circ$ small enough depending only on $n$, $s$, $\lambda$, and $\Lambda$, we reach a contradiction, and the claim is proved.

\item We now use the previous claim to finish the proof of the theorem.
Namely, we will show that if $u(0)>M_0$, with $M_0$ sufficiently large depending only on $n$, $s$, $\lambda$, and $\Lambda$, then this leads to $\sup_{B_{1/2}} u=\infty$, which is a contradiction.

Indeed, let $r_0 := (c_\circ M_0)^{-1/n}< \frac12$, with $c_\circ$ given by the claim in \ref{it:step2claim}.
Then,  there exists~$z_1\in B_{r_\circ}$ such that 
\[u(z_1)>(1+\delta)M_\circ =:M_1.\]
Applying iteratively the claim, we find $z_k$, $M_k$ and $r_k$ satisfying
\[z_{k+1}\in B_{r_{k}}{(z_k)},\qquad r_k= (c_\circ M_k)^{-1/n}, \qquad
u(z_k)> M_k := (1+\delta)M_{k-1}.\]
Since $M_k=(1+\delta)^k M_0 \to \infty$ as~$k\to\infty$, and 
\[|z_k| \leq \sum_{j=0}^{k-1} r_j \leq CM_0^{-1/n}\sum_{j=0}^{k-1} (1+\delta)^{-j/n} < \frac12,\]
if $M_0$ is large enough depending only on $n$, $s$, $\lambda$, and $\Lambda$, we have 
\[
\sup_{B_{1/2}} u\ge u(z_k) > (1+\delta)^k M_0 \to \infty\quad\text{as}\quad k\to \infty,
\]
reaching a contradiction (recall $u\in {\rm USC}(B_1)$). Thus $u(0)$ is bounded by a constant depending only on $n$, $s$, $\lambda$, and $\Lambda$, and we are done.\qedhere
\end{steps}
\end{proof}

As a consequence of this $L^\infty$ bound we can now improve  Theorem~\ref{C^alpha-bmc} to allow solutions in $L^1_{\omega_s}(\R^n)$:

\begin{thm}\label{C^alpha-bmc_2} \index{Holder estimate@H\"older estimate!Viscosity solutions!L1@$L^1_{\omega_s}$}
Let $s\in(0,1)$ and let $\Mpm$ be given by \eqref{eq:MMpm}. Let $u\in C(B_1)\cap L_{\omega_s}^1(\R^n)$ be any viscosity solution   to a non-divergence-form equation with bounded measurable coefficients, i.e.,
\[\left\{
\begin{array}{rcll}
\Mp u &\geq & -C_\circ &\quad \textrm{in}\quad B_1,\\
\Mm u &\leq & C_\circ &\quad \textrm{in}\quad B_1,
\end{array}
\right.\]
for some~$C_\circ\geq0$.
Then, $u\in C^\gamma_{\rm loc}(B_1)$ with 
\[\|u\|_{C^{\gamma}(B_{1/2})}\leq C\left( \|u\|_{L_{\omega_s}^1(\R^n) }+C_\circ\right),\]
where $C$ and $\gamma>0$ depend only on $n$, $s$, $\lambda$, and $\Lambda$.
\end{thm}
\begin{proof}
Let us define $w_1 := \eta u$ and $w_2 = u-w_1$, where $0\le \eta \le 1$ is a cut-off function such that $\eta \equiv 1$ in $B_{3/4}$ and $ \eta \equiv 0$ in $\R^n\setminus B_{4/5}$. Then, we know that by ellipticity
\[\left\{
\begin{array}{rcll}
\Mp w_1 &\geq & -C_\circ - \Mp w_2 &\quad \textrm{in}\quad B_{2/3}\\
\Mm w_1 &\leq & C_\circ - \Mm w_2&\quad \textrm{in}\quad B_{2/3},
\end{array}
\right.\]
in the viscosity sense.

Since $w_2 = 0$ in $B_{3/4}$, by Corollary~\ref{cor:Iu_II} we know that
\[
\|\Mpm w_2\|_{L^\infty(B_{2/3})} \le C \|w_2\|_{L^1_{\omega_s}(\R^n)}\le C \|u\|_{L^1_{\omega_s}(\R^n)}.
\]
and therefore, $w_1$ satisfies
\[\left\{
\begin{array}{rcll}
\Mp w_1 &\geq & -C_\circ - C \|u\|_{L^1_{\omega_s}(\R^n)}&\quad \textrm{in}\ B_{2/3}\\
\Mm w_1 &\leq & C_\circ + C \|u\|_{L^1_{\omega_s}(\R^n)}&\quad \textrm{in}\ B_{2/3}
\end{array}
\right.\]
in the viscosity sense. 
On the other hand, by Theorem~\ref{half-Harnack-sub} (applied to $u$ and $-u$) and a covering argument, we know that
\[
\|w_1\|_{L^\infty(\R^n)} \le \|u\|_{L^\infty(B_{4/5})}\le C \left(\|u \|_{L^1_{\omega_s}(\R^n)}+C_\circ\right).
\]
We can therefore use Theorem~\ref{C^alpha-bmc}  with $w_1$ to obtain (after another covering argument),
\[
\|u\|_{C^\gamma(B_{1/2})} = \|w_1\|_{C^\gamma(B_{1/2})} \le C \left(\|u \|_{L^1_{\omega_s}(\R^n)}+C_\circ\right),
\]
as we wanted to see. 
\end{proof}

\subsection{Harnack's inequality}

As an immediate consequence of the two half Harnack inequalities proved above, the Harnack inequality follows.

\begin{cor}[Harnack's inequality] \label{cor-Harnack} \index{Harnack's inequality!Fully nonlinear}
Let $s\in(0,1)$ and let $\Mpm$ be given by \eqref{eq:MMpm}. Let $u\in C(B_1)\cap L^1_{\omega_s}(\R^n)$ be  any viscosity solution to a non-divergence-form equation with bounded measurable coefficients, i.e.
\[\left\{
\begin{array}{rcll}
\Mp u &\geq & -C_\circ &\quad \textrm{in}\quad B_1\\
\Mm u &\leq & C_\circ &\quad \textrm{in}\quad B_1,
\end{array}
\right.\]
for some~$C_\circ\geq0$.
Assume in addition that $u\geq0$ in $\R^n$.
Then,
\[\sup_{B_{1/2}}u\leq C\left(\inf_{B_{1/2}} u+C_\circ\right),\]
where $C$ depends only on $n$, $s$, $\lambda$, and $\Lambda$.
\end{cor}

\begin{proof}
The result follows from Theorems \ref{half-Harnack-sup} and \ref{half-Harnack-sub}.
\end{proof}

Notice that, thanks to Theorems~\ref{half-Harnack-sup} and \ref{half-Harnack-sub}, we have that if $C_\circ = 0$ then both $\sup_{B_{1/2}} u$ and $\inf_{B_{1/2}} u$ are comparable to $\|u\|_{L^1_{\omega_s}(\R^n)}$. 

Moreover, if we only assume that $u\ge 0$ in $B_1$ (instead of $\R^n$) in the previous setting, we get instead 
\[\sup_{B_{1/2}}u\leq C\left(\inf_{B_{1/2}} u+C_\circ+\int_{\R^n}\frac{u_-(y)}{1+|y|^{n+2s}}\, dy\right),\]
where $u_-(y) = \max\{0, -u(y)\}$ is the negative part of $u$. This follows by applying Corollary~\ref{cor-Harnack} to $u_+ := u+u_-$ in $B_{3/4}$.

\subsection{Boundary regularity}

\index{Boundary regularity!Fully nonlinear}

Combining Corollary~\ref{cor:boundaryregularity_visc} with interior estimates, we obtain uniform continuity of $u$ in $\Omega$. The following proof is from \cite{CS3}:

\begin{cor}
\label{cor:bdry_reg_int}
Let $s\in (0, 1)$, and let $\Omega$ be any bounded open set satisfying \eqref{eq:cond_Om0}. Let $g\in L^1_{\omega_s}(\R^n)$ be continuous on $\partial\Omega$ with respect to $\R^n\setminus \Omega$  with modulus $\sigma$, in the sense of \eqref{eq:cont_from_outside}, and let $\Mpm$ be given by \eqref{eq:MMpm}. Let $u\in C(\Omega)\cap L^\infty(\Omega)\cap L^1_{\omega_s}(\R^n)$ be any viscosity solution of 
\[
\left\{
\begin{array}{rcll}
\Mp u & \ge & -C_\circ & \quad\text{in}\quad \Omega,\\
\Mm u & \le & C_\circ & \quad\text{in}\quad \Omega,\\
u & = & g & \quad\text{in}\quad \R^n\setminus \Omega.
\end{array}
\right.
\]
Then, $u\in C(\overline{\Omega})$ with 
\[
|u(x) - u(y)|\le \omega(|x-y|)\quad\text{for all}\quad x, y \in \overline{\Omega},
\]
for some modulus of continuity $\omega$ that depends only on $n$, $s$, $\lambda$, $\Lambda$, $\sigma$, $\|g\|_{L^1_{\omega_s}(\R^n)}$, $\|g\|_{L^\infty(\partial\Omega)}$, $C_\circ$, and $\Omega$, but is independent of $u$. Moreover, if $\sigma$ is $\alpha$-H\"older continuous, then $\omega$ is $\bar \alpha$-H\"older continuous as well, where $\bar \alpha$ depends only on $\alpha$, $\Omega$, $n$, $s$, $\lambda$, and $\Lambda$. 
\end{cor}
\begin{proof}
As in the proof of Theorem~\ref{thm:existence_visc} or Corollary~\ref{cor:boundaryregularity_visc}, we may assume that $g$ is bounded. The dependence on $\|g\|_{L^\infty(\R^n)}$ can then be replaced by a dependence on $\sigma$, $\|g\|_{L^1_{\omega_s}(\R^n)}$, and $\|g\|_{L^\infty(\partial\Omega)}$.

Let $w_\pm \in C(\overline{\Omega})\cap L^1_{\omega_s}(\R^n)$ be the unique viscosity solution  to 
\[
\left\{
\begin{array}{rcll}
\Mpm w_\pm & = & \mp C_\circ & \quad\text{in}\quad \Omega,\\
w_\pm & = & g & \quad\text{in}\quad \R^n\setminus \Omega,
\end{array}
\right.
\]
given by Theorem~\ref{thm:existence_visc} (applied with operator $\I(w, x)  := \Mpm w \pm C_\circ$). Then, by comparison principle, Theorem~\ref{thm:comparison_viscosity_2}, we  have that
\begin{equation}
\label{eq:thankstowpm}
w_- \le u \le w_+ \quad\text{in}\quad \R^n.
\end{equation}
Let now $x, y \in \overline{\Omega}$, and let $\bar x\in \partial\Omega$ be such that $\dist(x, \partial\Omega) = |x-\bar x| =: 2r$. Let us also denote $\rho := |x-y|$, so that we want to show $|u(x) - u(y)|\le \omega(\rho)$ for some $\omega$.  By Corollary~\ref{cor:boundaryregularity_visc} applied to both $w_+$ and $w_-$, we already know that this is true if $x\in \partial\Omega$ (thanks to \eqref{eq:thankstowpm}, since $w_- = w_+ = u$ on $\partial \Omega$), for some modulus that, as an abuse of notation, we still denote $\sigma$ (which is H\"older if $\sigma$ was H\"older). In particular, $u\in C(\overline{\Omega})$. 

If $ \rho > r/2$ we have
\[
|u(x) - u(y)|\le |u(x) - u(\bar x)|+|u(y) -u(\bar x)| \le \sigma(2r) +\sigma(\rho+2r)\le 2\sigma(5\rho)  
\]
and we can take $\omega(t) =2\sigma(5t)$ in this case. 

On the other hand, if $ \rho < r/2$ we let $\tilde u$ be the truncation
\[
\tilde u(x) = \left\{\begin{array}{ll}
u(\bar x) +\sigma(4r)& \quad\text{if}\quad u > u(\bar x) +\sigma(4r),\\
u(\bar x) -\sigma(4r)& \quad\text{if}\quad u < u(\bar x) +\sigma(4r),\\
u& \quad\text{otherwise}.
\end{array}
\right.
\]
For any $z\in \R^n$ we have 
\[
\left|u(z) - \tilde u(z)\right|\le \min\left\{\big(\sigma(|z-x|+r)-\sigma(4r)\big)_+, 2\|u\|_{L^\infty(\R^n)}\right\}.
\]
We consider the rescalings 
\[
u_r(z) := u(x+rz)\quad\text{and}\quad \bar u_r(z) := \bar u(x+rz),
\]
so that,  
\[
|u_r(z) - \tilde u_r(z)|\le \min\left\{\left[\sigma(r|z|+r)-\sigma(4r)\right]_+, 2\|u\|_{L^\infty(\R^n)}\right\}.
\]
Observe that in $B_{2}$ we have that $u_r = \bar u_r$. Thus,   we can bound $\Mp (u_r - \bar u_r)$ in $B_1$ by 
\[
|\Mp (u_r - \bar u_r)|\le C  \int_{  B_1^c}  \min\left\{\left[\sigma(r|y+\,\cdot\,|+r)-\sigma(4r)\right]_+, \|u\|_{L^\infty(\R^n)}\right\}\, \frac{dy}{|y|^{n+2s}}.
\]
In particular, 
\[
\|\Mp (u_r - \tilde u_r)\|_{L^\infty(B_1)}\le \mu(r),
\]
for some $\mu(r)\downarrow 0$ as $r\downarrow 0$, at a rate that depends only on $n$, $s$, $\Lambda$, and and upper bound for $\|u\|_{L^\infty(\R^n)}$. Observe, though, that by Corollary~\ref{cor:Linftybound_visc}, $\|u\|_{L^\infty(\R^n)}$ is controlled by $\|g\|_{L^\infty(\R^n\setminus \Omega)}$ and $C_\circ$. Moreover, if $\sigma$ is $\alpha$-H\"older with $\alpha < 2s$, then $\mu(r)$ is also $\alpha$-H\"older (that is, $\mu(r) \le C r^\alpha$ for some $ C > 0$). 

By ellipticity, Proposition~\ref{prop:viscosity_ellipticity}, we therefore have
\[
\Mp \tilde u_r \ge \Mp u_r - \Mp ( u_r-\tilde u_r)\ge -C_\circ r^{2s}  - \mu(r)\quad\text{in}\quad B_1. 
\]
We similarly obtain that
\[
\Mm \tilde u_r \le C_\circ r^{2s} + \mu(r),
\]
and, by the interior estimates in Theorem~\ref{C^alpha-bmc} applied to $\tilde u_r - u(\bar x)$, we get
\[
|u(x) - u(y)| = \big|\tilde u_r(0) - \tilde u_r((x-y)/r)\big|\le C \big(\sigma(4r) + \mu(r)\big)\left(\frac{\rho}{r}\right)^\alpha=:\nu(r) \left(\frac{\rho}{r}\right)^\alpha, 
\]
with $\rho < r/2$. Notice that $\nu(r)$ is a modulus of continuity, and by  making $\nu(r)$ larger if necessary, we may assume that $\nu(r) r^{-\alpha}$ is decreasing. 

Since $\rho < r/2$ we have
\[
|u(x) - u(y)|\le \nu(r) \left(\frac{\rho}{r}\right)^\alpha \le 2^\alpha \nu(2\rho) \downarrow 0\quad\text{as}\quad\rho\downarrow 0,
\]
as we wanted to see.  Finally,  if $\sigma$ is H\"older, then $\nu$ can be taken to be H\"older, and so $u$ (and hence, $\omega$) is H\"older continuous. 
\end{proof}

\section{Approximation of viscosity solutions}
\label{sec:regularization}

In order to prove the interior regularity results of Section~\ref{sec:int_reg}, and similarly to what happened in the linear case, we will first establish a priori estimates. Then, to recover estimates for general viscosity solutions, we will  need a result on the approximation of a viscosity solution by strong (or smooth) solutions. This is precisely the goal of this section. 

We want to find a regularization of a solution to    $\I(u,x) = 0$ in $B_1$. In the case of weak or distributional solutions to linear equations, this was directly accomplished in Lemma~\ref{lem:conv_dist_sol} by means of the convolution against a smooth mollifier. For viscosity solutions to fully nonlinear equations, the process is mcuh more involved, and we will also need to define regularized versions of $\I$. The results in this section follow closely the proofs of \cite{Fer23}, which are based partly on \cite{CS3, Kriv, Ser2}. 
  
Given an operator $\I\in \III$ of the form 
\begin{equation}
\label{eq:ILab}
\I  (u, x) = \inf_{b\in \B}\sup_{a\in \A}\big\{-\L_{ab}  u(x) + c_{ab}(x)\big\},\qquad \L_{ab} \in \LL_s(\lambda, \Lambda),
\end{equation}
we say that it has modulus $\sigma$ if $c_{ab}(x)$ are equicontinuous with modulus $\sigma$ (recall Definition~\ref{defi:II}). We will prove the following:

  \begin{thm}
\label{thm:cinftysol00}
Let $s\in (0, 1)$, and let $\I \in \III$. Let $u\in C(B_1)\cap L^1_{\omega_s}(\R^n)$ be any viscosity solution of
\[
\I (u, x) = f(x)\quad\text{in}\quad B_1
\]
for some $f\in C(B_1)$. 

Then, there exist  sequences of functions $u^{(\eps)}, f_\eps\in C^\infty_c(\R^n)$ such that 
\[
\begin{split}
 u^{(\eps)}\to u&\quad\text{uniformly in $B_{3/4}$ and in $L^1_{\omega_s}(\R^n)$},\\
 f_\eps\to f&\quad\text{uniformly in $B_{3/4}$,}
  \end{split}
 \]
and a sequence of operators $\I_\eps\in \III$ with 
\[
\I_\eps (0, x) \to \I(0, x) \quad\text{uniformly in $B_{3/4}$},\]
as $\eps \downarrow 0$, such that
\[
\I_\eps(u_\eps, x) = f_\eps(x) \quad\text{in}\quad B_{3/4}.
\]
Moreover, if $\I$ is translation invariant (resp. concave), then $\I_\eps$ are translation invariant (resp. concave). Finally, if $\I$ has modulus $\sigma$, then $\I_\eps$ have modulus $\sigma$ as well.  
\end{thm}
  \begin{rem}
The new operators $\I_\eps$ are $C^\infty$, in the sense that for any $w\in C^\infty_c(\R^n)$ we have $\I_\eps(w, x) \in C^\infty(\R^n)$, with vanishing derivatives at infinity; see the proof of Proposition~\ref{prop:cinfty1}.  Moreover, if $\I\in \II_s(\lambda, \Lambda; \theta)$ for some $\theta\in (0, 1)$, then $\I_\eps\in \II_s(\lambda, \Lambda; \theta)$ with $[\I_\eps]_{C^\theta}\le C[\I]_{C^\theta}$, with $C$ depending only on $n$, $s$, $\lambda$, $\Lambda$, and $\theta$; see Remark~\ref{rem:reg_inherited}. 
\end{rem}

  We first establish a simpler approximation result in the case of regular kernels, which will be used for the interior regularity later on. After that, we use the strategies developed in the first part to prove that viscosity solutions can always be approximated by strong $C^{2s+\delta}$ solutions, and using this we finally prove Theorem~\ref{thm:cinftysol00}.

\subsection{A first approximation result}
In the following, we consider operators of the form  
\begin{equation}
\label{eq:Iepsisoftheform0}
\I (u, x) = \inf_{b\in \B}\sup_{a\in \A}\big\{-\L_{ab} u(x) + c_{ab}(x)\big\},\qquad \L_{ab}\in \LL_s(\lambda, \Lambda; \theta),
\end{equation}
from which we define its regularized version as
\begin{equation}
\label{eq:Iepsisoftheform}
\I_\eps (u, x) = \inf_{b\in \B}\sup_{a\in \A}\left\{-\L_{ab}^{(\eps)} u(x) + c_{ab}(x)\right\},\qquad \L_{ab}^{(\eps)}\in \LL_s(\lambda, \Lambda; \theta),
\end{equation}
where   the $c_{ab}$ are the same as above. The first regularization or approximation result (and the one we will use in Section~\ref{sec:int_reg}) is then the following:

\begin{prop}
\label{prop:regularization}
Let $s\in (0, 1)$, and let $\I \in \II_s(\lambda, \Lambda; \theta)$ with $\theta\in [0, 1)$ be of the form \eqref{eq:Iepsisoftheform0}. Let us assume, moreover, that 
\[
\sup_{(a, b)\in \A\times\B} [c_{ab}]_{C^\theta(\R^n)}<\infty,
\]
where we denote $[\, \cdot\,]_{C^0} = {\rm osc}(\, \cdot\, )$. 

Let $u\in C(B_1)\cap L^1_{\omega_s}(\R^n)$ be any viscosity solution of
\[
\I (u, x) = 0\quad\text{in}\quad B_1. 
\]
Then, there exist: a sequence of functions,
\[\qquad \qquad u^{(\eps)}\in C^{2s+\theta}_{\rm loc}(B_{3/4})\cap {C(B_1)}\cap L^1_{\omega_s}(\R^n)\qquad \textrm{if}\quad 2s+\theta\notin\N,\]
or $u^{(\eps)}\in C_{\rm loc}^{1-\delta}(B_{3/4})\cap {C(B_1)}\cap L^1_{\omega_s}(\R^n)$ for any $\delta > 0$ if $\theta = 0$ and $s = \frac12$; and  a sequence of operators $\I_\eps\in \II_s(\lambda, \Lambda; \theta)$ of the form \eqref{eq:Iepsisoftheform} and satisfying $[\I_\eps]_{C^\theta}\le C_1 [\I]_{C^\theta}$ if  $\theta > 0$, with $C_1$ depending only on $n$, $s$, $\lambda$, $\Lambda$, and $\theta$,
such that, 
\[
\left\{
\begin{array}{rcll}
\I_\eps(u^{(\eps)}, x) & = & 0 & \quad\text{in}\quad B_{3/4}\\
u^{(\eps)} & = & u & \quad\text{in}\quad \R^n\setminus B_{3/4},
\end{array}
\right.
\]
and
\[
u^{(\eps)}\to u\quad\text{locally uniformly in $B_{3/4}$}.
\]
 Moreover, we have 
\begin{equation}
\label{eq:bound_reg}
\|u^{(\eps)}\|_{L^\infty(B_{3/4})}+\|u^{(\eps)}\|_{L^1_{\omega_s}(\R^n)}\le C \left(\| u\|_{L^1_{\omega_s}(\R^n)} + \|\I(0, x)\|_{L^\infty(B_{3/4})}\right)
\end{equation}
for some $C$ depending only on $n$, $s$, $\lambda$, and $\Lambda$. 
\end{prop}

Let us start with the construction of $\I_\eps$.  Let $\psi:[0, \infty) \to [0, \infty)$ be a given fixed cut-off function such that 
\[
\left\{\begin{array}{l}
\psi\in C^\infty_c([0, \infty)), \\
\psi = 1 \quad\text{in}\quad [0, 1/2],\\ 
\psi = 0 \quad\text{in}\quad [1, \infty),\\
\text{$\psi$ is monotone nonincreasing}.
\end{array}
\right.
\]
 Given $\L\in \LLL$ with kernel $K$ (which satisfies \eqref{eq:compdef}) and $\eps > 0$, we define $\L^{(\eps)}$ to be the operator that has kernel $K_\eps$ given by 
 \begin{equation}
 \label{eq:Kepsdef}
 K_\eps(y) = \big(1-\psi(|y|/\eps) \big) K(y) + \psi(|y|/\eps)|y|^{-n-2s}.
 \end{equation}
 Notice that with this definition we still have $\L^{(\eps)}\in \LLL$ (see Figure~\ref{fig:12}). Moreover, we have:
 
 \begin{figure}
\centering
\makebox[\textwidth][c]{\includegraphics[scale = 1]{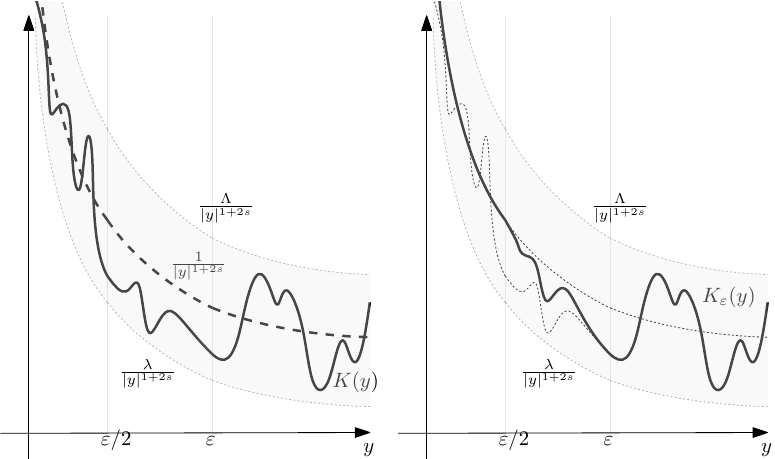}}
\caption{\label{fig:12} The kernels $K$ and $K_\eps$ given by \eqref{eq:Kepsdef}.}
\end{figure}
 
 \begin{lem}
 \label{lem:Lepsalpha} Let $s\in (0, 1)$. If $\L\in \LL_s(\lambda, \Lambda; \theta)$ for some $\theta\in (0, 1)$, then $\L^{(\eps)}\in \LL_s(\lambda, \Lambda; \theta)$ as well, with 
\[
 \big[\L^{(\eps)}\big]_{C^\theta} \le C    [\L]_{C^\theta}  
\]
 for some constant $C$ depending only on $n$, $s$, $\theta$, $\lambda$, and $\Lambda$. 
 \end{lem}
 \begin{proof}
 Indeed, notice that using \eqref{eq:APP_In1} on the definition of $K_\eps$ we get
 \[
 \begin{split}
 [K_\eps]_{C^\theta(B_\eps\setminus B_{\eps/2})} & \le  [K]_{C^\theta(B_\eps\setminus B_{\eps/2})} + (\Lambda+1)\eps^{-n-2s}  \big[\psi(|y|/\eps)\big]_{C^\theta(B_\eps\setminus B_{\eps/2})} \\
&  \qquad+ \big[|y|^{-n-2s}\big]_{C^\theta(B_\eps\setminus B_{\eps/2})}\\
& \le \frac{C}{\eps^{n+2s+\theta}}\left([\L]_{C^\theta}+  [\psi]_{C^\theta(B_1\setminus B_{1/2})} + \big[|y|^{-n-2s}\big]_{C^\theta(B_1\setminus B_{1/2})}\right),
 \end{split}
 \]
 which, since $\psi$ is fixed, directly implies 
 \[
  [K_\eps]_{C^\theta(B_\eps\setminus B_{\eps/2})} \le C([\L]_{C^\theta}+1)\eps^{-n-2s-\theta}. 
 \]
 On the other hand, we already know that 
 \[
   [K_\eps]_{C^\theta(B_\rho\setminus B_{\rho/2})} \le C[\L]_{C^\theta}\rho^{-n-2s-\theta}\quad\text{for any}\quad \rho \ge 2\eps,
 \]
 since  $K_\eps = K$ there, and 
 \[
   [K_\eps]_{C^\theta(B_\rho\setminus B_{\rho/2})} \le C \rho^{-n-2s-\theta}\quad\text{for any}\quad \rho \le \eps/2,
 \]
 since $K_\eps (y)= |y|^{-n-2s}$ there. Thus, we have 
  \[
  [K_\eps]_{C^\theta(B_\rho\setminus B_{\rho/2})} \le C([\L]_{C^\theta}+1)\rho^{-n-2s-\theta}\quad\text{for all}\quad \rho > 0.
 \]
 Finally, since $[\L]_{C^\theta}\ge c > 0$ for all $\L \in \LL_s(\lambda, \Lambda; \theta)$, the result follows. 
 \end{proof}

Let now $\I \in \II_s(\lambda, \Lambda; \theta)$ for some $\theta\in [0, 1)$, i.e., of the form \eqref{eq:Iepsisoftheform0}. We define $\I_\eps$ as \eqref{eq:Iepsisoftheform} with $\L_{ab}^{(\eps)}$ given by \eqref{eq:Kepsdef}. By Lemma~\ref{lem:Lepsalpha} we immediately have that 
\[
\I_\eps\in \II_s(\lambda, \Lambda; \theta)
\]
 as well, with $[\I_\eps]_{C^\theta} \le C [\I]_{C^\theta}$  if $\theta >0$. Furthermore, $\I_\eps$ weakly converges to $\I$ as $\eps\downarrow 0$:

 \begin{lem} 
 \label{lem:weakconv}
Let $s\in (0, 1)$, and let $\I, \I_\eps\in \II_s(\lambda, \Lambda)$ be as above. Then 
 \[
 \I_\eps\rightharpoonup \I\quad\text{in}\quad \R^n,\quad \text{as}\quad \eps\downarrow 0,
 \]
 in the sense of Definition~\ref{defi:conv_I}. 
 \end{lem}
 \begin{proof}
 Let $x_\circ\in \R^n$, and let $v\in L^1_{\omega_s}(\R^n)$ such that it is $C^2$ in $B_r(x_\circ)$ for some $r > 0$. Let us compute, for any $x\in B_{r/2}(x_\circ)$ and $\L\in \II_s(\lambda, \Lambda)$ with kernel $K$,
 \[
 \L^{(\eps)} v(x) - \L v(x) =\frac12\int_{\R^n} \big(2v(x)-v(x+y)-v(x-y)\big)\big(K_\eps(y)-K(y)\big)\,dy.
 \]
 Since $K_\eps(y)-K(y) = \psi(|y|/\eps)\left(|y|^{-n-2s}-K(y)\right)$ we can bound the right-hand side by 
 \[
\big| \L^{(\eps)} v(x) - \L v(x)\big|\le C \int_{B_\eps} \big|2v(x)-v(x+y)-v(x-y)\big||y|^{-n-2s}\,dy.
 \]
 Now, using that $v$ is $C^2$ in $B_r(x_\circ)$ and taking $\eps < r/4$ we get
  \[
  \begin{split}
\big| \L^{(\eps)} v(x) - \L v(x)\big|& \le C \|v\|_{C^2(B_{3r/4}(x_\circ))}\int_{B_\eps} |y|^{-n-2s+2}\,dy\\
& \le C \|v\|_{C^2(B_{3r/4}(x_\circ))}\eps^{2-2s}.
\end{split}
 \]
 Hence, we can bound 
\[
\begin{split}
\I_\eps (v, x) & = \inf_{b\in \B}\sup_{a\in \A}\left\{-\L_{ab}^{(\eps)} u(x) + c_{ab}(x)\right\}\\
& \le \inf_{b\in \B}\sup_{a\in \A}\left\{-\L_{ab} u(x) + c_{ab}(x)\right\} + C \|v\|_{C^2(B_{3r/4}(x_\circ))}\eps^{2-2s}\\
& = \I(v, x) + C \|v\|_{C^2(B_{3r/4}(x_\circ))}\eps^{2-2s}.
\end{split}
\] 
On the other hand, we also get similarly, 
\[
\I_\eps (v, x) \ge \I(v, x) - C \|v\|_{C^2(B_{3r/4}(x_\circ))}\eps^{2-2s}.
\]
In all, we have proved that 
\[
\|\I_\eps(v, \cdot) - \I(v, \cdot)\|_{L^\infty(B_{r/2}(x_\circ))}\le C \|v\|_{C^2(B_{3r/4}(x_\circ))}\eps^{2-2s}\downarrow 0
\]
as $\eps\downarrow 0$, and hence $\I_\eps\rightharpoonup\I$.  
\end{proof}

We want to use the previous operators $\I_\eps$ to construct a series of regular solutions approximating a given solution. That is, let $\I \in \II_s(\lambda, \Lambda; \theta)$ for some $\theta\in [0, 1)$, and let $u\in C(B_1)\cap L^1_{\omega_s}(\R^n)$ be such that 
\begin{equation}
\label{eq:estrelletax2}
\I (u, x) = 0\quad\text{in}\quad B_1. 
\end{equation}
We then define $u^{(\eps)}$ to be the unique solution, given by Theorem~\ref{thm:existence_visc}, to 
\begin{equation}
\label{eq:uepsqualdef}
\left\{
\begin{array}{rcll}
\I_\eps(u^{(\eps)}, x) & = & 0 & \quad \text{in}\quad B_{3/4}\\
u^{(\eps)} & = & u & \quad \text{in}\quad \R^n\setminus B_{3/4}. 
\end{array}
\right.
\end{equation}

\begin{lem}
\label{lem:ueps1}
Let $s\in (0, 1)$ and $\I \in \III$. 
Let $u\in C(B_1)\cap L^1_{\omega_s}(\R^n)$ be any solution of \eqref{eq:estrelletax2}, and let $u^{(\eps)}\in C(B_{3/4})\cap L^\infty(B_{3/4})\cap L^1_{\omega_s}(\R^n)$ be the unique solution of \eqref{eq:uepsqualdef}. Let  $\gamma > 0$ be given by Theorem~\ref{C^alpha-bmc_2}. Then 
\[
\| u^{(\eps)}\|_{L^\infty(B_{3/4})}+\|u^{(\eps)}\|_{C^\gamma(B_{1/2})}\le C \left(\| u\|_{L^1_{\omega_s}(\R^n)} + \|\I(0, x)\|_{L^\infty(B_{3/4})}\right),
\]
 for some $C$ depending only on $n$, $s$, $\lambda$, and $\Lambda$.
\end{lem} 
\begin{proof}
 Observe that,  
\[
\Mp u \ge \I(u, x) - \I(0, x) \ge \Mm u 
\]
and therefore, from Theorem~\ref{C^alpha-bmc_2} (or Theorem~\ref{half-Harnack-sub} applied to $u$ and $-u$),  
\begin{equation}
\label{eq:asdfophsdfadsf}
\|u\|_{L^\infty(B_{7/8})}\le C\left( \|u\|_{L^1_{\omega_s}(\R^n)}+ \|\I(0, x)\|_{L^\infty(B_1)}\right).
\end{equation}
Let $\eta\in C^\infty(B_1)$ with $0\le \eta\le 1$ be a cut-off function such that $\eta \equiv 0$ in $B_{6/7}^c$ and $\eta \equiv 1$ in $B_{5/6}$. Then, the function $\eta u^{(\eps)}$ satisfies 
\[
\tilde \I_\eps(\eta u^{(\eps)}, x) = 0\quad\text{in}\quad B_{3/4},
\]
where if $\I_\eps$ is of the form \eqref{eq:Iepsisoftheform}, then $\tilde \I_\eps$ is (also using that $(1-\eta)u = (1-\eta)u^{(\eps)}$): 
\[
\tilde \I_\eps (v, x) := \inf_{b\in \B}\sup_{a\in \A}\left\{-\L_{ab}^{(\eps)} v(x) + c_{ab}(x)-\L_{ab}^{(\eps)} [(1-\eta)u](x)\right\}.
\]
In particular, by Corollary~\ref{cor:Linftybound_visc}  together with the ellipticity condition in Proposition~\ref{prop:viscosity_ellipticity}, we get
\begin{equation}
\label{eq:asdfophadsf}
\|\eta u^{(\eps)}\|_{L^\infty(B_{3/4})}\le C \left(\|\eta u\|_{L^\infty(B^c_{3/4})} + \|\tilde \I_\eps(0, x)\|_{L^\infty(B_{3/4})}\right).
\end{equation}
Now, by Lemma~\ref{lem:Lu_LL} we know that for any $\L \in \LLL$,
\[
\big\|\L[(1-\eta)u]\big\|_{L^\infty(B_{3/4})} \le C\|u\|_{L^1_{\omega_s}(\R^n)}
\]
for some $C$ depending only on $n$, $s$, $\lambda$, and $\Lambda$. Thus, 
\begin{equation}
\label{eq:asdfophadsf2}
\|\tilde \I_\eps (0, x) - \I (0, x)\|_{L^\infty(B_{3/4})} \le C\|u\|_{L^1_{\omega_s}(\R^n)}.
\end{equation}
Since we also have $\eta u^{(\eps)} = u^{(\eps)}$ in $B_{3/4}$ and $\|\eta u\|_{L^\infty(B_{3/4}^c)} \le \|u\|_{L^\infty(B_{7/8})}$, from \eqref{eq:asdfophadsf}-\eqref{eq:asdfophadsf2}-\eqref{eq:asdfophsdfadsf} we obtain 
\[
\| u^{(\eps)}\|_{L^\infty(B_{3/4})}\le C \left(\| u\|_{L^1_{\omega_s}(\R^n)} + \|\I(0, x)\|_{L^\infty(B_{3/4})}\right).
\]
In particular, since $u = u^{(\eps)}$ in $\R^n\setminus B_{3/4}$, this also implies 
\begin{equation}
\label{eq:asdfophadsf3}
\| u^{(\eps)}\|_{L^1_{\omega_s}(\R^n)}\le C \left(\| u\|_{L^1_{\omega_s}(\R^n)} + \|\I(0, x)\|_{L^\infty(B_{3/4})}\right).
\end{equation}
We finally apply Theorem~\ref{C^alpha-bmc_2} to $u^{(\eps)}$ as above, to deduce, together with \eqref{eq:asdfophadsf2}-\eqref{eq:asdfophadsf3}, 
\[
\| u^{(\eps)}\|_{C^\gamma(B_{1/2})}\le C \left(\| u\|_{L^1_{\omega_s}(\R^n)} + \|\I(0, x)\|_{L^\infty(B_{3/4})}\right),
\]
which completes the proof.
\end{proof}

We now want to show that the solutions $u^{(\eps)}$ are qualitatively regular (that is, strong solutions) in the interior of $B_{3/4}$. In order to do it, we use the structure of the operator $\I_\eps$, which  behaves like a fractional Laplacian. Thus, we need the interior estimates for the fractional Laplacian in the case of viscosity solutions: 

\begin{prop}[Interior estimates for viscosity solutions of $\fls$]
\label{prop:viscosity_fls} \index{Interior regularity!Fractional Laplacian!Viscosity solutions}
Let $s\in (0, 1)$, and let $f\in C^\theta(B_1)$ for some $\theta\in [0, 1)$. Let $u\in C(B_1) \cap L^1_{\omega_s}(\R^n)$ satisfy
\[
\fls u = f \quad \text{in}\quad B_1
\]
in the viscosity sense. Then, if $2s+\theta\notin \N$,  $u\in C_{\rm loc}^{2s+\theta}(B_1)$ with 
\[
\|u\|_{C^{2s+\theta}(B_{1/2})}\le C \left(\|u\|_{L^1_{\omega_s}(\R^n)} + \|f\|_{C^\theta(B_1)}\right),
\]
for some $C$ depending only on $n$, $s$, and $\theta$. If $\theta = 0$ and $s = \frac12$, then $u\in C^{1-\delta}(B_1)$ for any $\delta > 0$. 
\end{prop}
\begin{proof}
Let us define $v_1$ and $v_2$ to be the solutions of
\[
\left\{
\begin{array}{rcll}
\fls v_1 & = & f& \quad\text{in}\quad B_{3/4},\\
v_1 & = & 0& \quad\text{in}\quad \R^n\setminus B_{3/4},
\end{array}
\right.
\]
and
\[
\left\{
\begin{array}{rcll}
\fls v_2 & = & 0& \quad\text{in}\quad B_{3/4},\\
v_2 & = & u& \quad\text{in}\quad \R^n\setminus B_{3/4}.
\end{array}
\right.
\]
Here,  $v_1$ is the unique weak solution to its problem, given by Theorem~\ref{thm:exist_weak_sol}, which satisfies (by Theorems~\ref{thm:globCsreg}, \ref{thm-interior-linear-Lp}, and \ref{thm-interior-linear-2}) $v_1\in C(\overline{B_{3/4}})\cap C_{\rm loc}^{2s+\theta}(B_{3/4})$ with interior estimates. On the other hand, $v_2\in C(\overline{B_{3/4}})\cap C^\infty(B_{3/4})$ is given by the Poisson kernel representation, Proposition~\ref{prop:poisson_kernel_ball_s}. Let us define
\[
v := v_1+v_2 \in C(\overline{B_{3/4}})\cap C_{\rm loc}^{2s+\theta}(B_{3/4}),
\]
and we assume $\theta > 0$. By Lemma~\ref{lem:equiv_defi}, $v$ satisfies 
 \[
\left\{
\begin{array}{rcll}
\fls v & = & f& \quad\text{in}\quad B_{3/4},\\
v & = & u& \quad\text{in}\quad \R^n\setminus B_{3/4},
\end{array}
\right.
\]
in the viscosity sense, too. By uniqueness of viscosity solutions, Corollary~\ref{cor:uniqueness_viscosity}, $v = u$ in $\R^n$, and by the a priori interior regularity estimates for the fractional Laplacian, Theorem~\ref{thm:interior_regularity_fls} (since $v\in C_{\rm loc}^{2s+\theta}(B_{3/4})$), we have 
\[
\|u\|_{C^{2s+\theta}(B_{1/2})}\le C \left(\|u\|_{L^1_{\omega_s}(\R^n)} + \|f\|_{C^\theta(B_1)}\right),
\]
as we wanted to see. 

Finally, if $\theta = 0$, we can regularize $v_1$ first taking the convolution against a smooth function (by Lemma~\ref{lem:conv_dist_sol}, see also Remark~\ref{rem:relax_reg_dist}), to obtain a sequence of approximate smooth solutions $v_1^\delta+v_2\to v_1+v_2$ locally uniformly as $\delta \downarrow 1$ in $B_1$ (which are also viscosity solutions, by Lemma~\ref{lem:equiv_defi}). Then, proceeding as before, we are done  thanks to the stability of viscosity solutions,  Proposition~\ref{prop:stab_super}. 
\end{proof}

The following is the qualitative result on the regularity of $u^{(\eps)}$:

\begin{lem}
\label{lem:regularization}
Let $s\in (0, 1)$. Let $u^{(\eps)}$ be defined as above, \eqref{eq:uepsqualdef}, for a fixed $\I \in \II_s(\lambda, \Lambda; \theta)$ with $\theta\in [0, 1)$ of the form \eqref{eq:Iepsisoftheform0}. Let us assume, moreover, that 
\[
\sup_{(a, b)\in \A\times\B} [c_{ab}]_{C^\theta(\R^n)}\le C_\circ <\infty,
\]
where we denote $[\, \cdot\,]_{C^0} = {\rm osc}(\, \cdot\, )$.  Then, if $2s+\theta\notin\N$, $u^{(\eps)}\in C^{2s+\theta}_{\rm loc}(B_{3/4})$. If $\theta = 0$ and $s  = \frac12$, we have $u^{(\eps)}\in C_{\rm loc}^{1-\delta}(B_{3/4})$ for any $\delta > 0$.
\end{lem}
\begin{proof}
For the sake of readability, let us denote $v = u^{(\eps)}$. Notice that, by Lemma~\ref{lem:ueps1} and a covering argument, we already know that $v$ is $C^\gamma$ inside~$B_{3/4}$. 

 We now express the operator $\I_\eps$ as follows:
\begin{equation}
\label{eq:Iepssplit}
\begin{split}
\I_\eps(v, x) & = -c^{-1}_{n,s}\fls v(x) +\inf_{b\in \B}\sup_{a\in \A}\left\{\tilde \L^{(\eps)}_{ab} v (x) + c_{ab}(x)\right\}\\
& =  -c^{-1}_{n,s}\fls v(x) + f_\eps(x),
\end{split}
\end{equation}
where we have denoted, for each $\tilde \L = \tilde \L^{(\eps)}_{ab}$,
\[
\begin{split}
\tilde \L v (x) & = \left(c^{-1}_{n,s}\fls-\L_{ab}^{(\eps)}\right) v(x)\\
&  = \frac12 \int_{\R^n}\big(2v(x) - v(x+y)-v(x-y)\big) \tilde K_\eps(y)\, dy\\
& = \int_{B_{\eps/2}^c}\big(v(x) - v(x+y)\big) \tilde K_\eps(y)\, dy,
\end{split}
\]
with
\[
\tilde K_\eps(y) = \big(1-\psi(|y|/\eps)\big) \left(|y|^{-n-2s}-K_{ab}(y)\right)\in L^1(\R^n) ,
\]
where $c_{n, s}$ is the constant in \eqref{eq:cns}, and $K_{ab}$ is the kernel of the operator $\L_{ab}$ (recall \eqref{eq:Iepsisoftheform0}-\eqref{eq:Iepsisoftheform}). In particular, 
\begin{equation}
\label{eq:fromexp}
\tilde \L v (x) = C_\eps v(x) - \int_{B_{\eps/2}^c(x)} v(z) \tilde K_\eps(z-x)\, dz.
\end{equation}
Observe that, since $\I \in \II_s(\lambda, \Lambda; \theta)$ (see Definition~\ref{defi:LL}), for all $0\le \mu \le \theta$,
\[
\big|\tilde K_\eps(x_\circ+h) -\tilde K_\eps(x_\circ)\big| \le C_\eps |x_\circ|^{-n-2s-\mu} |h|^\mu\quad\text{for any}\ x_\circ\in B_{\eps/2}^c,\ h \in B_{\eps/4}.
\]
Let now $x\in B_{3/4}$ fixed, and let 
\[
\rho = \min\left\{\frac{\eps}{4}, \frac12 \dist(x, \partial B_{3/4})\right\} = \min\left\{\frac{\eps}{4},\frac12 \left(\frac34 - |x|\right)\right\}.\] 
Then, from \eqref{eq:fromexp} we can bound
\[
\begin{split}
[\tilde \L v]_{C^\mu(B_{\rho}(x))} & \le C_\eps [ v]_{C^\mu(B_{\rho}(x))} +C_\eps \int_{B_{\eps/2}^c(x)} v(z) |z-x|^{-n-2s-\mu}\, dz \\
& \le C_\eps \left( [ v]_{C^\mu(B_{\rho}(x))} +\|v\|_{L^1_{\omega_s}(\R^n)}\right),
\end{split}
\]
for any $\tilde \L = \tilde \L^{(\eps)}_{ab}$, where $C_\eps$ is independent of $(a, b)\in \A\times\B$.  In \eqref{eq:Iepssplit} we can therefore bound the H\"older semi-norms of $f_\eps$ (being the $\inf\sup$ of H\"older functions\footnote{The $\inf\sup$ of $C^\mu$ functions is $C^\mu$, whenever $\mu <1$. If $\mu \ge 1$, then it is at most Lipschitz ($C^{0,1}$). This is why this proof does not obtain higher regularity even if $\theta\gg 1$.}) as
\[
[f_\eps]_{C^\mu(B_{\rho}(x))}\le C_\eps \left( [ v]_{C^\mu(B_{\rho}(x))} +\|v\|_{L^1_{\omega_s}(\R^n)} + C_\circ\right).
\]
Thus, we obtain that 
\begin{equation}
\label{eq:vimplication}
v\in C_{\rm loc}^\mu(B_{3/4})\quad\text{and}\quad 0\le \mu \le\theta \quad \Longrightarrow \quad f_\eps\in C_{\rm loc}^\mu(B_{3/4}),
\end{equation}
in a qualitative way.

Moreover, $v$ satisfies, by assumption
\[
\fls v = f_\eps\quad\text{in}\quad B_{3/4}. 
\]
Hence, we can now use interior estimates for viscosity solutions with  the fractional Laplacian, Proposition~\ref{prop:viscosity_fls}, together with a bootstrap argument to conclude: we already know that $v\in C^\gamma(B_{3/4})$, hence by \eqref{eq:vimplication} we have $f_\eps\in C^\gamma(B_{3/4})$ and by the interior estimates in Proposition~\ref{prop:viscosity_fls} $v\in C^{2s+\min\{\gamma, \theta\}}(B_{3/4})$. If $\theta > \gamma$, we can iteratively repeat this until $\gamma+ ms > \theta$ for some $m\in \N$, at which point we have to stop when we reach $C^\theta$ regularity of $f_\eps$. A final application of interior estimates implies $C^{2s+\theta}$  regularity of $v$. If $\theta = 0$, we only apply the iteration once. 
\end{proof}

We can finally prove the regularization result:
\begin{proof}[Proof of Proposition~\ref{prop:regularization}]\label{proof:thmregularization}
We construct $\I_\eps$ and $u^{(\eps)}$ as \eqref{eq:Iepsisoftheform} and \eqref{eq:uepsqualdef}. Then, Lemma~\ref{lem:weakconv} gives the weak convergence of $\I_\eps$ to  $\I$, and Lemma~\ref{lem:ueps1} and a covering argument give  the locally uniform convergence of $u^{(\eps)}$ in $B_{3/4}$ (by Arzel\`a-Ascoli, up to taking subsequences), towards some function $\tilde u  \in C(B_{3/4})\cap L^\infty(B_{3/4})$.  Moreover, since $u^{(\eps)}$ is uniformly bounded in $B_{3/4}$, it converges to $\tilde u$ in $L^1_{\omega_s}(\R^n)$ as well, by the dominated convergence theorem and the locally uniform convergence. Hence we are in a situation where we can apply Proposition~\ref{prop:stab_super} (and Remark~\ref{rem:stab_subsol}) to deduce that $\tilde u\in C(B_{3/4})\cap L^1_{\omega_s}(\R^n)\cap L^\infty(B_{3/4})$ satisfies
\[
\left\{
\begin{array}{rcll}
\I(\tilde u, x) & = & 0 & \quad\text{in}\quad B_{3/4}\\
\tilde u & = & u & \quad\text{in}\quad \R^n\setminus B_{3/4}.
\end{array}
\right.
\]
By the uniqueness of bounded viscosity solutions, Corollary~\ref{cor:uniqueness_viscosity_2}, we have $\tilde u = u$, and moreover $u\in C(B_1)$  (by Theorem~\ref{thm:existence_visc}).

The bound of the $L^1_{\omega_s}(\R^n)$ norm of $u^{(\eps)}$ is a consequence of Lemma~\ref{lem:ueps1}, and finally the interior regularity is due to Lemma~\ref{lem:regularization}. This completes the proof. 
\end{proof}

\subsection{Approximation by strong solutions} \label{ssec:classicalapprox}

Proposition~\ref{prop:regularization} suffices for our purposes in subsections \ref{ssec:C1alpha} and \ref{ssec:C1alpha2}, on the interior  regularity of fully nonlinear elliptic equations. It gives an approximating sequence to a viscosity solution by \emph{smoother} solutions, which in the case $\theta > 0$ are strong.  Let us now show that, with a bit more of work, also in the  most general case $\theta = 0$ we can consider  strong solutions as the approximating sequence: 

\index{Approximation by strong solutions}
\begin{prop}
\label{prop:regularization_2}
Let $s\in (0, 1)$, and let $\I \in \II_s(\lambda, \Lambda)$. Let $u\in C(B_1)\cap L^1_{\omega_s}(\R^n)$ be any viscosity solution of
\[
\I (u, x) = 0\quad\text{in}\quad B_1. 
\]
Then, there exist $\delta > 0$, a sequence of functions,
\[  C^{2s+\delta}_{\rm loc}(B_{3/4}) \cap {C^\delta_c(\R^n)} \ni u^{(\eps)} \to u \  \text{locally uniformly in $B_{3/4}$ and  in $L^1_{\omega_s}(\R^n)$},\]
 and  a sequence of operators $\hat\I_\eps\in \II_s(\lambda, \Lambda)$ of the form \eqref{eq:Iepsisoftheform2},  such that
\[
\begin{split}
&\hat \I_\eps(u^{(\eps)}, x)   =  0  \quad\text{in } B_{3/4}\\ 
& \hat\I_\eps\rightharpoonup \I \quad\text{in the sense of Definition~\ref{defi:conv_I}},
\end{split}
\]
as $\eps\downarrow 0$. Moreover, we have 
\begin{equation}
\label{eq:bound_reg_2}
\|u^{(\eps)}\|_{L^\infty(B_{3/4})}+\|u^{(\eps)}\|_{L^1_{\omega_s}(\R^n)}\le C \left(\| u\|_{L^1_{\omega_s}(\R^n)} + \|\I(0, x)\|_{L^\infty(B_{3/4})}+\sigma(\eps)\right)
\end{equation}
for some $C$ depending only on $n$, $s$, $\lambda$, and $\Lambda$, and where $\sigma(\eps)\downarrow 0$ as $\eps\downarrow 0$. 
\end{prop}

In order to prove it, we proceed following a similar strategy to the one before. Now, however, we need to regularize the $c_{ab}(x)$ in the definition of~$\I$, as well as the value of $u$ outside of $B_{3/4}$. We will do that by means of a convolution. 

Let us fix a mollifier $\varphi$ such that
\begin{equation}
\label{eq:mollifiervarphi}
\text{$\varphi\in C^\infty_c(B_1)
$ is radial, with $\varphi \ge 0$ and $\textstyle{\int_{B_1}}\varphi = 1$,}
\end{equation}
and we consider the rescalings
\begin{equation}
\label{eq:mollifiervarphi2}
\varphi_\eps(x) := \frac{1}{\eps^n}\varphi\left(\frac{x}{\eps}\right)\in C^\infty_c(B_\eps). 
\end{equation}

We   define $\hat \I_\eps$ analogously to \eqref{eq:Iepsisoftheform} but also regularizing the terms $c_{ab}(x)$. Lemma~\ref{lem:weakconv} still holds in this case:
 \begin{lem} 
 \label{lem:weakconv2}
Let $s\in (0, 1)$, and for any $\I$ of the form \eqref{eq:ILab} we consider 
\begin{equation}
\label{eq:Iepsisoftheform2}
\hat \I_\eps (u, x) := \inf_{b\in \B}\sup_{a\in \A}\left\{-\L_{ab}^{(\eps)} u(x) + c^{(\eps)}_{ab}(x)\right\},\qquad \L_{ab}^{(\eps)}\in \LL_s(\lambda, \Lambda),
\end{equation}
where $\L_{ab}^{(\eps)}$ are the corresponding operators to $\L_{ab}$ but with kernel given by \eqref{eq:Kepsdef}, and where $c_{ab}^{(\eps)}(x) := (\varphi_\eps * c_{ab})(x)$.  Then 
 \[
\hat \I_\eps\rightharpoonup \I\quad\text{in}\quad \R^n,\quad \text{as}\quad \eps\downarrow 0,
 \]
 in the sense of Definition~\ref{defi:conv_I}. 
 \end{lem}
 \begin{proof}
The proof is exactly the same as that of Lemma~\ref{lem:weakconv}, where we now use that since $c_{ab}(x)$ are equicontinuous, $c_{ab}^{(\eps)}(x)$ converges locally uniformly  to $c_{ab}(x)$ as $\eps\downarrow 0$ independently of $(a, b)\in \A\times \B$ (that is, depending only on the modulus $\sigma$ from Definition~\ref{defi:II}). 
 \end{proof}
 
 If $\I \in \II_s(\lambda, \Lambda)$ and $u\in C(B_1)\cap L^1_{\omega_s}(\R^n)$ is a viscosity solution to
\begin{equation}
\label{eq:estrelletax22}
\I (u, x) = 0\quad\text{in}\quad B_1,
\end{equation}
we   define our new functions $u^{(\eps)}$ to be the unique solution, given by Theorem~\ref{thm:existence_visc}, to 
\begin{equation}
\label{eq:uepsqualdef2}
\left\{
\begin{array}{rcll}
\hat \I_\eps(u^{(\eps)}, x) & = & 0 & \quad \text{in}\quad B_{3/4}\\
u^{(\eps)} & = & (u\chi_{B_{1/\eps}}) * \varphi_\eps & \quad \text{in}\quad \R^n\setminus B_{3/4}. 
\end{array}
\right.
\end{equation}

In doing so, the following analogue of 	Lemma~\ref{lem:ueps1} also holds now:

\begin{lem}
\label{lem:ueps12}
Let $s\in (0, 1)$ and $\I \in \III$. Let $u\in C(B_1)\cap L^1_{\omega_s}(\R^n)$ be any viscosity solution of \eqref{eq:estrelletax22}, and let $u^{(\eps)}\in C(B_{3/4})\cap  L^\infty(B_{3/4})\cap L^1_{\omega_s}(\R^n) $ be the unique solution of \eqref{eq:uepsqualdef2}. Let  $\gamma > 0$ be given by Theorem~\ref{C^alpha-bmc_2}. Then 
\[
\| u^{(\eps)}\|_{L^\infty(B_{3/4})}+\|u^{(\eps)}\|_{C^\gamma(B_{1/2})}\le C \left(\| u\|_{L^1_{\omega_s}(\R^n)} + \|\I(0, x)\|_{L^\infty(B_{3/4})}+\sigma(\eps)\right),
\]
 for some $C$ depending only on $n$, $s$, $\lambda$, and $\Lambda$, and where $\sigma$ is the modulus of continuity associated to $c_{ab}$ in the definition of $\I$ (see Definition~\ref{defi:II}). 
\end{lem} 
\begin{proof}
The proof is the same as that of Lemma~\ref{lem:ueps1},  using that 
\[
\|(u\chi_{B_{1/\eps}}) * \varphi_\eps\|_{L^1_{\omega_s}(\R^n)}\le 2 \|u \chi_{B_{1/\eps}}\|_{L^1_{\omega_s}(\R^n)}\le 2 \|u\|_{L^1_{\omega_s}(\R^n)}. 
\]
The main difference is the appearance of $\sigma(\eps)$ on the right-hand side of the estimate. This is because we now define 
\[
\tilde \I_\eps (v, x) := \inf_{b\in \B}\sup_{a\in \A}\left\{-\L_{ab}^{(\eps)} v(x) + c^{(\eps)}_{ab}(x)-\L_{ab}^{(\eps)} [(1-\eta)u](x)\right\},
\]
and therefore the bound \eqref{eq:asdfophadsf2} becomes
\[
\|\tilde \I_\eps (0, x) - \I (0, x)\|_{L^\infty(B_{3/4})} \le C\|u\|_{L^1_{\omega_s}(\R^n)} + \sup_{(a, b)\in \A\times \B } \|c_{ab}^{(\eps)} - c_{ab}\|_{L^\infty(B_{3/4})}.
\]
From the definition of $\sigma$ and recalling that $c_{ab}^{(\eps)}$ is a convolution of $c_{ab}$, this gives the desired estimate. 
\end{proof}

By regularizing the boundary datum we  have now improved the regularity of $u^{(\eps)}$ with respect to the previous case, Lemma~\ref{lem:regularization}: 

\begin{lem}
\label{lem:regularization2}
Let $s\in (0, 1)$ and let  $\I \in \II_s(\lambda, \Lambda)$. Let  $u\in C(B_1)\cap L^1_{\omega_s}(\R^n)$ be any viscosity solution to $\I(u, x) = 0$ in $B_1$, and let $u^{(\eps)}$ be defined by \eqref{eq:uepsqualdef2}.  Then, there exists $\delta  >0$ independent of $\eps > 0$ such that $u^{(\eps)}\in C^{2s+\delta}_{\rm loc}(B_{3/4})\cap C_c^\delta(\R^n)$.
\end{lem}
\begin{proof}
For the sake of readability, we denote $v = u^{(\eps)}$. Observe that
\[
\|\nabla ((u\chi_{B_{1/\eps}}) * \varphi_\eps)\|_{L^\infty(\R^n)}\le C_\eps,
\]
for some $C_\eps$ that might blow-up as $\eps\downarrow 0$. This is enough to deduce, from Corollary~\ref{cor:bdry_reg_int}, that there exists some $\delta > 0$ (independent of $\eps> 0$) such that $v \in C^\delta_c(\R^n)$.

As in Lemma~\ref{lem:regularization}, we rewrite the operator $\hat \I_\eps$ as 
\[
\hat\I_\eps(v, x)  =  -c^{-1}_{n,s}\fls v(x) + f_\eps(x),
\]
where 
\[
f_\eps(x) := \inf_{b\in \B}\sup_{a\in \A}\left\{\tilde \L^{(\eps)}_{ab} v (x) + c^{(\eps)}_{ab}(x)\right\}
\]
and, for each $\tilde \L = \tilde\L^{(\eps)}_{ab}$, 
\[
\tilde \L v (x) = C_\eps v(x) - \int_{B_{\eps/2}^c} v(z+x) \tilde K_\eps(z)\, dz,
\]
with
\[
\tilde K_\eps(y) = \big(1-\psi(|y|/\eps)\big) \left(|y|^{-n-2s}-K_{ab}(y)\right)\in L^1(\R^n).
\]
In particular, since $v\in C_c^\delta(\R^n)$, for any $x\in B_{3/4}$ and $h\in \R^n$ small, 
\[
\left|\tilde \L v(x+h) - \tilde \L v (x)\right|\le  C_\eps |h|^\delta + \int_{B_{\eps/2}^c} C_\eps |h|^\delta \tilde K_\eps(z) \, dz \le C_\eps|h|^\delta. 
\]
Together with the fact that $c_{ab}^{(\eps)}\in C^\infty$, we get  
\[
f_\eps(x) \in C^\delta_{\rm loc}(B_{3/4}) 
\]
for some $\delta > 0$ independent of $\eps$. By the interior estimates for viscosity solutions with the fractional Laplacian, Proposition~\ref{prop:viscosity_fls}, we deduce $v\in C^{2s+\delta}_{\rm loc}(B_{3/4})$, as wanted. 
\end{proof}

We can finally prove Proposition~\ref{prop:regularization_2}:
\begin{proof}[Proof of Proposition~\ref{prop:regularization_2}]
We proceed as in the proof of Proposition~\ref{prop:regularization} on page~\pageref{proof:thmregularization}, with the corresponding changes in this new situation. 

We construct $\hat \I_\eps$ and $u^{(\eps)}$ as \eqref{eq:Iepsisoftheform2} and \eqref{eq:uepsqualdef2}, and  Lemma~\ref{lem:weakconv2} gives the weak convergence of $\hat \I_\eps$ to  $\I$, while Lemma~\ref{lem:ueps12} and a covering argument give  the locally uniform convergence  in $B_{3/4}$ and the convergence in $L^1_{\omega_s}(\R^n)$ of $u^{(\eps)}$ to some $\tilde u  \in C(B_{3/4})\cap L^\infty(B_{3/4})\cap L^1_{\omega_s}(\R^n)$.  Proposition~\ref{prop:stab_super} (and Remark~\ref{rem:stab_subsol}) now imply that $\tilde u$ satisfies
\[
\left\{
\begin{array}{rcll}
\I(\tilde u, x) & = & 0 & \quad\text{in}\quad B_{3/4}\\
\tilde u & = & u & \quad\text{in}\quad \R^n\setminus B_{3/4},
\end{array}
\right.
\]
and by uniqueness (Corollary~\ref{cor:uniqueness_viscosity_2}), we have $\tilde u = u$, and $u\in C(B_1)$  (by Theorem~\ref{thm:existence_visc}).  The bound on the $L^1_{\omega_s}(\R^n)$ norm of $u^{(\eps)}$ is a consequence of Lemma~\ref{lem:ueps12}, and finally its qualitative interior regularity is due to Lemma~\ref{lem:regularization2}. This completes the proof. 
\end{proof}

\subsection{Equivalence between viscosity and distributional solutions} As a consequence of Proposition~\ref{prop:regularization_2} we obtain that, in the \emph{linear} case (taking operators $\L\in \LLL$, recall Definition~\ref{defi:LL}), the notions of viscosity (Definition~\ref{defi:viscosity}) and distributional (Definition~\ref{defi:dist}) solution are equivalent:

\begin{lem}
\label{lem:visc_dist}
Let $s\in (0, 1)$, $\L \in \LLL$, $u\in L^\infty_{2s-\eps}(\R^n)$, and $f \in C(B_1)$. Then, $u$ solves $\L u = f$ in $B_1$ in the distributional sense if and only if it does so in the  viscosity sense.
\end{lem}
\begin{proof}
If $u$ is a distributional solution, it is continuous by Theorem~\ref{thm-interior-linear-Lp}, and we can regularize it and consider 
\[
u_\eps := u * \varphi_\eps, 
\]
where $\varphi_\eps$ is given by \eqref{eq:mollifiervarphi}-\eqref{eq:mollifiervarphi2}. Then $u_\eps$ satisfies 
\[
\L u_\eps = f_\eps \quad\text{in}\quad B_{1-\eps}
\]
in the strong sense (see Lemma~\ref{lem:conv_dist_sol}), and therefore, in the viscosity sense as well (by Lemma~\ref{lem:equiv_defi}). Taking the limit $\eps \downarrow 0$, by Proposition~\ref{prop:stab_super} $u$ is a viscosity solution to $\L u = f$ in $B_1$.

Conversely, if $u\in C(B_1)$ is a viscosity solution to the equation, by Proposition~\ref{prop:regularization_2} it can be approximated by  strong  solutions (and therefore, distributional solutions, see Lemma~\ref{cor:dist_strong_smooth}) $u_\eps\to u$  to an equation of the form 
\[
\hat{\L}_\eps u_\eps = f_\eps \quad\text{in}\quad {B_{3/4}},
\]
with a sequence of explicit operators $\hat \L_\eps$. 

Then, the limit $\eps\downarrow 0$ is a distributional solution to $\L_\infty u = f$ in $B_{3/4}$, where by construction $\L_\infty = \L$. A covering argument, yields that $\L u = f $ in $B_1$ in the distributional sense. 
\end{proof}

\subsection{Approximation by $C^\infty$ solutions}

Our next goal is to finally prove that we can actually approximate any viscosity solution by $C^\infty_c(\R^n)$ solutions. We consider $\I$  of the form \eqref{eq:ILab},  where $c_{ab}$ are equicontinuous with modulus $\sigma$, and we want to show Theorem~\ref{thm:cinftysol00}.

In order do it, we will combine the approximation by strong solutions in Proposition~\ref{prop:regularization_2} with the next result.

\begin{prop}
\label{prop:cinfty1}
Let $s\in (0, 1)$, and let $\I\in \III$ of the form \eqref{eq:ILab} with modulus $\sigma$. Let $u\in C^{2s+\delta}_{\rm loc}(B_1)\cap C^\delta_c(\R^n)$ be any   solution of 
\[
\I(u, x) = f(x) \quad\text{in}\quad B_1
\] 
for some $f\in C(B_1)$  and $\delta > 0$.  Let $(\varphi_\eps)_{\eps > 0}$ be given by  \eqref{eq:mollifiervarphi}-\eqref{eq:mollifiervarphi2}.

Then, there exists some $\I_\eps\in \III$ of the form \eqref{eq:ILab} with modulus~$\sigma$ such that the sequence $u_\eps := u * \varphi_\eps\in C^\infty_c(\R^n)$ satisfies 
\[
\I_\eps(u_\eps, x) = f_\eps(x)\quad\text{in}\quad B_1
\]
for some $f_\eps\in C^\infty(B_1)$ such that 
\[
f_\eps \to  f\quad\text{uniformly in $B_{3/4}$  as $\eps\downarrow 0$.}
\]
Moreover, 
\[
\I_\eps(0, x) \to  \I(0, x) \quad\text{uniformly in $B_{3/4}$  as $\eps\downarrow 0$.}
\]
\end{prop}

  Before proving Proposition~\ref{prop:cinfty1}, we consider the following lemma on the representation of Lipschitz functions: 

\begin{lem}
\label{lem:lip_rep}
Let $f\in C^1(\R^n)$ with $\nabla f\in L^\infty(\R^n)$. Then
\[
f(x) = \inf_{z\in \R^n}\sup_{v\in E\subset \R^n} \left\{v\cdot x - v\cdot z + f(z)\right\},
\]
where $ E :=\nabla f(\R^n)\subset \R^n$. 
\end{lem}
 \begin{figure}
\centering
\makebox[\textwidth][c]{\includegraphics[scale = 1]{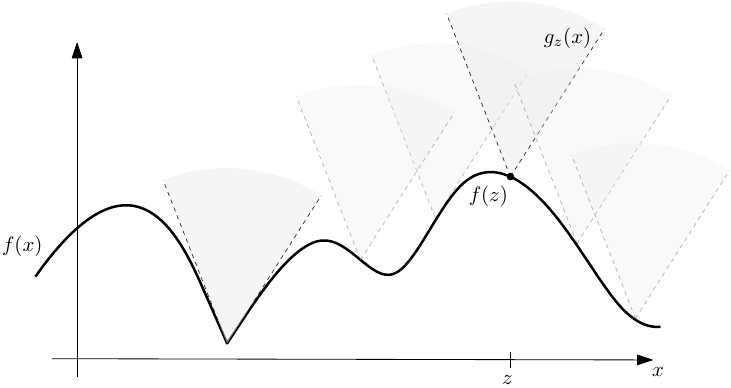}}
\caption{\label{fig:16}Any Lipschitz function $f$ can be expressed as the infimum of convex functions $g_z(x)$.}
\end{figure}
\begin{proof}
Let, for each $z\in \R^n$ fixed, 
\[
g_z(x) := \sup_{v\in E\subset \R^n}  v\cdot (x -  z) + f(z).
\]
We then have
\[
f(x)  = (x-z)\cdot \int_0^1 \nabla f((1-t)z+ tx)\, dt + f(z)  \le g_z(x).
\]
Thus, 
\[
f(x) \le \inf_{z\in \R^n} g_z(x),
\]
(see Figure~\ref{fig:16}) and since $f(z) = g_z(z)$ for all $z\in \R^n$, we are done. 
\end{proof}

We can now prove the approximation of strong solutions by smooth solutions:

\begin{proof}[Proof of Proposition~\ref{prop:cinfty1}]
We divide the proof into four steps. For the sake of readability, we assume $f = 0$. The general case follows analogously by taking $\I(\cdot, x) - f(x)$.
\begin{steps}
\item We define $c_{ab}^{(\eps)} := c_{ab}*\varphi_\eps\in C^\infty(\R^n)$ and  consider
\[
\hat{\I}_\eps (v, x) := \inf_{b\in \B}\sup_{a\in \A}\left\{-\L_{ab}  v(x) + c^{(\eps)}_{ab}(x)\right\}.
\]
Notice that $\L_{ab} u_\eps \in C^\delta_{\rm loc}(B_1)$ (see Lemma~\ref{lem:Lu_2}-\ref{it:lem_Lu2_ii}) with local uniform (in $a$, $b$, and $\eps$) estimates in $B_1$, as well as $\L_{ab} u_\eps \in C^\infty(\R^n)$ (locally uniformly in $a$ and $b$, but not in $\eps$) with vanishing derivatives at infinity. Since $c_{ab}$ are equicontinuous, the family $\L_{ab}  u_\eps(x) + c^{(\eps)}_{ab}(x)$ is locally equicontinuous in $B_1$. In particular, there exists a modulus of continuity\footnote{Actually, $\omega(r) = Cr^\delta + \sigma(r)$.} $\omega$  such that $\L_{ab}  u_\eps(x) + c^{(\eps)}_{ab}(x)$ is continuous with modulus $\omega$ in $B_{3/4}$, for all $(a, b)\in \A\times \B$ and $\eps \ge 0$.  

Hence, in fact, $(\hat{\I}_\eps(u_\eps, x))_{\eps\ge 0}$ is locally equicontinuous in $B_1$, and 
\begin{equation}
\label{eq:hatieps}
  \hat{\I}_\eps(u_\eps, x)\to 0\quad\text{locally uniformly in $B_1$},
\end{equation}
 (recall $f \equiv 0$) as well as 
\begin{equation}
\label{eq:hatieps2}
  \hat{\I}_\eps(0, x)\to \I(0, x) \quad\text{locally uniformly in $B_1$}.
\end{equation}

\item We now consider, for any $\eps > 0$ fixed, a finite collection of points $G_\eps :=\{y_1,\dots,y_{N_\eps}\}$ with $y_i \in B_{3/4}$ for $1\le i\le N_\eps$ such that $\dist(z, G_\eps) \le \zeta$ for all $z\in B_{3/4}$, where $\zeta = \zeta(\eps)$ is chosen small enough so that $\omega(\zeta)\le \eps/4$  (where $\omega$ is the modulus of continuity of the previous step).

We   want to take a finite redefinition of $\hat{\I}_\eps$ such that its value at $u_\eps$ and $0$ is not altered too much. Namely, for any $y_i\in G_\eps$, we consider $b_i, b_{N_\eps+i}\in \B$ such that if we define 
\[
\begin{split}
\mathcal{G}_i(v, x) & := \sup_{a\in \A}\left\{-\L_{ab_i}  v(x) + c^{(\eps)}_{ab_i}(x)\right\}\quad\text{for}\quad 1 \le i \le 2N_\eps
\end{split}
\]
then 
\[
\begin{array}{rcccll}
 0&\le &\mathcal{G}_i(u_\eps, y_i) - \hat{\I}_\eps(u_\eps, y_i)  & \le& {\eps}/{4},&\\[0.1cm]
 0&\le  &\mathcal{G}_{N_\eps +i}(0, y_i) - \hat{\I}_\eps(0, y_i)  & \le&  {\eps}/{4} & \quad \text{for} \quad 1\le i \le N_\eps.
\end{array}
\]
 Together with the fact that $\mathcal{G}_i(v, x) \ge \hat{\I}_\eps(v, x)$ in $\R^n$ for all $1 \le i \le 2N_\eps$, and from the choice of $\zeta$, we have 
\begin{equation}
\label{eq:tocombg_i}
\begin{array}{rcccll}
  0&\le&{\displaystyle \inf_{1\le i \le 2N_\eps}} \mathcal{G}_i(u_\eps, x) - \hat{\I}_\eps (u_\eps, x)& \le  {\eps}/{2}&\quad\text{in}\quad B_{3/4},\\[0.2cm]
 0&\le & {\displaystyle \inf_{1\le i \le 2N_\eps} }\mathcal{G}_i(0,  x) - \hat{\I}_\eps (0, x)& \le  {\eps}/{2}&\quad\text{in}\quad B_{3/4}. 
\end{array}
\end{equation}
Similarly, for each $1\le i \le 2N_\eps$ fixed, and for any $y_j\in G_\eps$ we consider $a_{ij}, a_{i, N_\eps+j}\in \A$ such that 
\[
\begin{split}
\left|-\L_{a_{ij} b_i} u_\eps(y_j) + c_{a_{i,j} b_i}^{(\eps)}(y_j) -\mathcal{G}_i(u_\eps, y_j)\right|& \le  {\eps}/{4},\\
\left|c_{a_{i,N_\eps + j} b_i}^{(\eps)}(y_j) - \mathcal{G}_i(0, y_j)\right|& \le  {\eps}/{4}\quad\text{for}\quad 1 \le j \le N_\eps.
\end{split}
\]
In particular, again by the choice of $\zeta$ above, we have that  
\[
\begin{split}
\left|\sup_{1\le j \le 2N_\eps} \left\{-\L_{a_{ij} b_i} u_\eps(x) + c_{a_{ij} b_i}^{(\eps)}(x)\right\} - \mathcal{G}_i(u_\eps, x)\right|& \le  {\eps}/{2}\quad\text{in}\quad B_{3/4},\\
\left|\sup_{1\le j \le 2N_\eps}   c_{a_{ij} b_i}^{(\eps)}(x) - \mathcal{G}_i(0, x)\right|& \le  {\eps}/{2}\quad\text{in}\quad B_{3/4}.
\end{split}
\]

Combined with \eqref{eq:tocombg_i} we get 
\[
\begin{split}
\left| \inf_{1\le i \le 2N_\eps} \sup_{1\le j \le 2N_\eps} \left\{-\L_{a_{ij} b_i} u_\eps(x) + c_{a_{ij} b_i}^{(\eps)}(x)\right\} - \hat{\I}_\eps (u_\eps, x)\right|& \le {\eps}\quad\text{in}\quad B_{3/4},\\
\left| \inf_{1\le i \le 2N_\eps} \sup_{1\le j \le 2N_\eps} \left\{ c_{a_{ij} b_i}^{(\eps)}(x)\right\} - \hat{\I}_\eps (0, x)\right|& \le {\eps}\quad\text{in}\quad B_{3/4}.
\end{split}
\]
Thus, we can define 
\[
\I_\eps^*(v, x ) := \inf_{1\le i \le 2N_\eps} \sup_{1\le j \le2 N_\eps} \left\{-\tilde \L_{ij} v(x) + \tilde c_{ij}^{(\eps)}(x)\right\}
\]
where
\[
\tilde\L_{ij} := \L_{a_{ij} b_i}\in   \LL_s(\lambda, \Lambda)\quad\text{and}\quad \tilde c_{ij}^{(\eps)}= c_{a_{ij} b_i}^{(\eps)}\quad\text{for all}\quad 1\le i,j\le 2N_\eps
\]
and we have that 
\begin{equation}
\label{eq:Iepfftohat}
\begin{split}
\big|\I_\eps^*(u_\eps, x) -\hat{\I}_\eps(u_\eps, x)\big|&\le \eps\quad\text{in}\quad B_{3/4},\\
\big|\I_\eps^*(0, x) -\hat{\I}_\eps(0, x)\big|&\le \eps\quad\text{in}\quad B_{3/4}.
\end{split}
\end{equation}
The key  difference now is that $\I^*_\eps$ is a \emph{finite} $\inf\sup$.

\item Let us denote, for the sake of readability, $N := 2N_\eps$. We define $F_\eps:\R^{N\times N} \to \R$ as 
\[
F_\eps(\{x_{ij}\}_{1\le i,j\le N}) = F_\eps 
\begin{pmatrix}
x_{11} & x_{12} & \dots & x_{1 N}\\
x_{21} & x_{22} & \dots & x_{2 N}\\
\vdots & \vdots & \ddots & \vdots\\
x_{N 1} & x_{N 2} & \dots & x_{N N}
\end{pmatrix} = \inf_{1\le i \le N} \sup_{1 \le j \le N} x_{ij},
\]
so that 
\begin{equation}
\label{eq:Iepsf}
\I_\eps^*(v, x) = F_\eps\left(\left\{-\tilde \L_{ij} v(x) + \tilde c_{ij}^{(\eps)}(x)\right\}_{1\le i,j\le N}\right).
\end{equation}

Then, $F_\eps$ is a piecewise linear function with $|\nabla F_\eps| = 1$ almost everywhere and such that for a.e. $x\in \R^{N\times N}$, $\nabla F_\eps(x) \in  \{\be_{ij}\}_{1\le i,j\le N}$, where $\be_{ij} \in \R^{N\times N}$ is the matrix with $(\be_{ij})_{ij} = 1$ and $(\be_{ij})_{k\ell} = 0$ for all $(k, \ell)\neq (i, j)$. 

In particular, by considering a regularization $F_\eps^r := F_\eps * \varphi_\eps$, where $\varphi_\eps\in C^\infty_c(B_\eps)$ with $B_\eps\in \R^{N\times N}$ (see \eqref{eq:mollifiervarphi}-\eqref{eq:mollifiervarphi2}) we have that $F_\eps^r \in C^\infty(\R^{N\times N})$ with 
\[
{\rm Grad}( F_\eps^r) := \bigcup_{x\in \R^{N\times N}} \nabla F_\eps^r(x)\subset \partial \, {\rm Conv}\big( \{\be_{ij}\}_{1\le i,j\le N}\big),
\]
where ${\rm Conv}(A)$ denotes the convex hull of $A\in \R^{N\times N}$. Since $|\nabla F_\eps |\le 1$, 
\begin{equation}
\label{eq:FepsFepsr}
\|F_\eps - F_\eps^r\|_{L^\infty(\R^{N\times N})}\le \eps,
\end{equation}
and   we can write it as 
\[
F_\eps^r(x) = \inf_{z\in \R^{N\times N}}\sup_{M\in {\rm Grad}(F_\eps^r)} \big\{M \cdot x - M\cdot z+F^r_\eps(z)\big\}.
\]
(This representation  is valid for any Lipschitz function, see Lemma~\ref{lem:lip_rep}.)  We then define
\begin{equation}
\label{eq:Iepsvfollows}
\I_\eps(v, x) := F_\eps^r\left(\left\{-\tilde \L_{ij} v(x) + \tilde c_{ij}^{(\eps)}(x)\right\}_{1 \le i,j\le N}\right),
\end{equation}
so that\footnote{If $f \not\equiv 0$, we would have now $\I_\eps(v, x)  - (f*\varphi_\eps)(x)$ as a regularized version of $\I(v, x) - f(x)$, since $\sum_{i, j} M_{ij} = 1$.}
\begin{equation}\label{eq:repIeps}
\begin{split}
\I_\eps(v, x) & = \inf_{z\in \R^{N\times N}}\sup_{M\in {\rm Grad}(F_\eps^r)} \left\{ \sum_{i,j =1}^{N} \left(-M_{ij} \tilde \L_{ij} v(x) + M_{ij} \tilde c_{ij}^{(\eps)}(x)\right) +C_{M,z}^\eps\right\}\\
& \hspace{-1.4cm}=\hspace{-1mm} \inf_{z\in \R^{N\times N}}\hspace{-0.5mm}\sup_{M\in {\rm Grad}(F_\eps^r)} \hspace{-0.5mm}  \left\{\hspace{-0.5mm} - \hspace{-0.5mm}\left(\sum_{i,j =1}^{N} M_{ij} \tilde \L_{ij} \right) v(x) + \left( \sum_{i, j = 1}^N M_{ij} \tilde c_{ij}^{(\eps)}(x) +C_{M,z}^\eps\right)\hspace{-0.5mm} \right\},
\end{split} 
\end{equation}
where 
\[
C_{M, z}^\eps := F_\eps^r(z) - M\cdot z.
\]
In particular, since $\sum_{i,j=1}^{N} M_{ij} = 1$, $M_{ij}\ge 0$, and $\LLL$ is convex, we have that $\I_\eps\in \III$ with
\[
\I  (v, x) = \inf_{b\in \hat  \B}\sup_{a\in \hat  \A}\left\{-\hat  \L^{(\eps)}_{ab}  v(x) + \hat c^{(\eps)}_{ab}(x)\right\},\qquad \hat  \L^{(\eps)}_{ab} \in \LL_s(\lambda, \Lambda),
\]  
and where $\hat  c^{(\eps)}_{ab}$ are equicontinuous with modulus $\sigma$  (the same as for $c_{ab}$). 

\item To finish, we notice that by the chain rule, since $\tilde \L_{ij} u_\eps, \tilde c_{ij}^{(\eps)}\in C^\infty(\R^n)$, it follows from \eqref{eq:Iepsvfollows} that $\I_\eps(u_\eps, x) \in C^\infty(\R^n)$. 

Moreover, thanks to \eqref{eq:FepsFepsr}-\eqref{eq:Iepsf} together with \eqref{eq:Iepfftohat} and \eqref{eq:hatieps}-\eqref{eq:hatieps2}, we have 
\[
\begin{array}{ll}
\I_\eps(u_\eps, x)  \to 0&\quad\text{uniformly in $B_{3/4}$}\\
\I_\eps(0, x)  \to \I(0, x) &\quad\text{uniformly in $B_{3/4}$}.
\end{array}
\]
This completes the proof.\qedhere
\end{steps}
\end{proof}

With this, we can complete the approximation result by $C^\infty_c$ solutions: 

\begin{proof}[Proof of Theorem~\ref{thm:cinftysol00}]
By defining the operator $\J(\cdot, x) := \I(\cdot, x) - f(x)$, we consider first the sequence of functions $u^{(\eps)}$ from Proposition~\ref{prop:regularization_2} applied with operator $\J$ in $B_{5/6}$ (after a scaling argument), so $u^{(\eps)}\in C^{2s+\delta}(B_{5/6})\cap C^\delta_c(\R^n)$. Notice that this also generates a sequence of operators $\hat \J_\eps(\cdot, x) = \hat \I_\eps(\cdot, x) - (f*\varphi_\eps)(x)$, so that if $\I$ had modulus $\sigma$, so do $\hat \I_\eps$  (see \eqref{eq:Iepsisoftheform2}).

Each $u^{(\eps)}$ can then be regularized by applying  Proposition~\ref{prop:cinfty1} (rescaled to $B_{5/6}$), which together with a diagonal argument yields the desired result.  
\end{proof}

\begin{rem}
\label{rem:reg_inherited}
In  Theorem~\ref{thm:cinftysol00} we have that, in fact, $f_\eps = f*\varphi_\eps$. Furthermore, notice that from the proof of  Proposition~\ref{prop:cinfty1}, and more precisely, from the representation \eqref{eq:repIeps} together with Lemma~\ref{lem:Lepsalpha}, we have that if $\I\in \II_s(\lambda, \Lambda; \theta)$ for some $\theta\in (0, 1)$, then $\I_\eps\in \II_s(\lambda, \Lambda; \theta)$ as well, with $[\I_\eps]_{C^\theta}\le C[\I]_{C^\theta}$, and $C$ depending only on $n$, $s$, $\lambda$, $\Lambda$, and $\theta$. Finally, also from \eqref{eq:repIeps}, if $\I$ is of the form \eqref{eq:ILab}, and $\I_\eps$ is of the form 
\[
\I_\eps  (u, x) = \inf_{b'\in \B_\eps}\sup_{a'\in \A_\eps}\big\{-\L_{a'b'}^{(\eps)}  u(x) + c_{a'b'}^{(\eps)}(x)\big\},\qquad \L_{a'b'} \in \LL_s(\lambda, \Lambda),
\]
then for any $(a', b') \in \A_\eps\times \B_\eps$, 
\[
[c^{(\eps)}_{a'b'}]_{C^\mu(\R^n)}\le  \sup_{(a, b)\in \A\times \B} [c_{ab}]_{C^\mu(\R^n)},
\]
for $\mu > 0$, $\mu\notin\N$.
\end{rem}

\section{Interior regularity results}\label{sec:int_reg}
\index{Interior regularity!Fully nonlinear}

  After having proved the existence, uniqueness, and approximation of (viscosity) solutions to fully nonlinear equations of the type $\I(u, x) = 0$ in $\Omega\subset \R^n$, we now turn our attention to their regularity. For this, we will use the Harnack inequalities and the H\"older estimates that we established in the previous section. 

For second-order (uniformly elliptic) fully nonlinear PDE of the form 
\begin{equation}
\label{eq:Fd2uint}
F(D^2 u ) = 0\quad\text{in}\quad \Omega\subset \R^n
\end{equation}
(i.e., when $ s = 1$),  there are two main interior regularity results   (see e.g. \cite{CC, FR4}): 
\begin{itemize}[leftmargin=1cm]

\item By the Krylov--Safonov theorem, solutions to fully nonlinear equations \eqref{eq:Fd2uint} are $C^{1,\alpha}$ in $\Omega$, for some small $\alpha > 0$; \cite{KS79}.

\item By the Evans--Krylov theorem, if $F$ is \emph{concave} (or convex) then solutions to \eqref{eq:Fd2uint} are $C^{2,\alpha}$ in $\Omega$, for some $\alpha > 0$; \cite{Eva82, Kry82}. 
In particular, they are strong solutions and the equation \eqref{eq:Fd2uint} holds pointwise. 

\end{itemize}

Here, we will establish the analogous results in the nonlocal setting. 
These nonlocal regularity estimates were first established by Caffarelli and Silvestre in \cite{CS, CS2}, and later refined by Kriventsov \cite{Kriv} and Serra \cite{Ser,Ser2}; see also \cite{S-nonlocal}.

We will start by showing the $C^{1,\alpha}$ regularity of solutions to fully nonlinear equations, for some $\alpha > 0$: 
\[
\left\{
\begin{array}{l}
\I(u, x) = 0 \quad\text{in}\quad \Omega, \\[0.2cm]
{\small \I (u, x) =  \inf_{b}\sup_a \left\{-\L_{ab} u + c_{ab}\right\}}\\[0.1cm]
{\small \text{with }\L_{ab}\in \LL_s(\lambda, \Lambda; \theta)}
\end{array}
\right.
\quad\Longrightarrow\quad u \in C^{1,\alpha}(\Omega)
\]
provided that $\theta \ge \max\{1+\alpha - 2s, 0\}$. We will then show $C^{2s+\alpha}$ regularity of solutions when the operator $\I$ is concave or convex: 
\[
\left\{
\begin{array}{l}
\I(u, x) = 0 \quad\text{in}\quad \Omega, \\[0.2cm]
{\small \I (u, x) =  \inf_{a} \left\{-\L_{a} u + c_{a}\right\}}\\[0.1cm]
{\small \text{with }\L_{a}\in \LL_s(\lambda, \Lambda; \alpha)}
\end{array}
\right.
\quad\Longrightarrow\quad u \in C^{2s+\alpha}(\Omega).
\]
In particular, solutions are strong and the equation $\I(u, x) = 0$  holds pointwise in $\Omega$. We refer to Theorems~\ref{C1alpha} and \ref{EvansKrylov} below for the precise statements.

\subsection{$C^{1,\alpha}$ regularity}\label{ssec:C1alpha}

\index{Interior regularity!Fully nonlinear!C1alpha@$C^{1,\alpha}$ estimates}

\index{Krylov-Safonov}

In the following theorem, given $\theta \ge 0$ we consider operators of the form 
\begin{equation}
\label{eq:Iform11}
\I (u, x) = \inf_{b\in \B}\sup_{a\in \A}\big\{-\L_{ab} u(x) + c_{ab}(x)\big\},\qquad \L_{ab}\in \LL_s(\lambda, \Lambda; \theta),
\end{equation}
and we prove the $C^{1,\alpha}$ regularity of solutions for \emph{all} $s\in(0,1)$.
In particular, when $s\leq \frac12$,  this is enough to conclude that solutions are classical or strong.

\begin{thm}\label{C1alpha}
Let $s\in(0,1)$, $\gamma>0$ be given by Theorem~\ref{C^alpha-bmc_2}, and $\alpha\in(0,\gamma)$. Let $\theta:=(1+\alpha-2s)_+$, and let $\I \in \II_s(\lambda, \Lambda; \theta)$ be of the form \eqref{eq:Iform11} and satisfying 
\begin{equation}
\label{eq:Iform22}
\sup_{(a, b)\in \A\times\B} [c_{ab}]_{C^\theta(\R^n)}\le C_\circ,
\end{equation}
where we denote $[\,\cdot\,]_{C^0} = \osc (\,\cdot\,)$. 

 Let $u\in C(B_1)\cap L^1_{\omega_s}(\R^n)$ be a viscosity solution of 
\[\I (u, x)=0\quad \textrm{in}\quad B_1.\]
 Then, $u\in C_{\rm loc}^{1,\alpha}(B_1)$ and 
\begin{equation}\label{sdnvwieuh} 
\|u\|_{C^{1+\alpha}(B_{1/2})} \leq 
C\left( \|u\|_{L^1_{\omega_s}(\R^n)}+\I(0, 0) +C_\circ \right),
\end{equation}
for some constant $C$ depending only on $n$, $s$, $\alpha$, $\lambda$, $\Lambda$, and $[\I]_{C^\theta}$ if $\theta > 0$.
\end{thm}

When the operator $\I$ is local (i.e., $s = 1$), this $C^{1,\alpha}$ estimate follows by applying iteratively the H\"older estimate from Theorem~\ref{C^alpha-bmc_2} to the incremental quotients $\frac{u(x+h)-u(x)}{|h|^\beta}$, improving an exponent $\gamma > 0$ at each step (see, for example, \cite[Chapter 4]{FR4}). For nonlocal equations, however, one has to take care of the tails of the functions, and this is why the proof of this result is more involved. 

The proof we present here is based on a blow-up argument, very similar to the ones we saw in  Section~\ref{sec:int_reg_G}.
For this, we will need the following  Liouville-type theorem (cf. Proposition~\ref{prop:compactness}), where we now consider operators of the form  
\begin{equation}
\label{eq:Ioftheform}
\I (u, x) = \inf_{b\in \B}\sup_{a\in \A}\big\{-\L_{ab} u(x) + c_{ab}(x)\big\},\qquad \L_{ab}\in \LLL.
\end{equation}

\index{Liouville's theorem!Fully nonlinear}
\begin{prop}\label{KrylovSafonov-Liouv}
Let $s\in(0,1)$, $\delta > 0$, $\gamma>0$ be given by Theorem~\ref{C^alpha-bmc_2}, and let $\alpha\in(0,\gamma)$. Let $\I\in \III$ of the form \eqref{eq:Ioftheform} and satisfying 
\[
\sup_{(a, b)\in \A\times\B} [c_{ab}]_{C^\theta(\R^n)}\le \delta\quad\text{for some}\quad \theta\in [0, 1),
\]
where we denote $[\,\cdot\,]_{C^0} = \osc (\,\cdot\,)$.

Assume that $u\in C_{\rm loc}^{1,\alpha}(\R^n)\cap L^1_{\omega_s}(\R^n)$ satisfies
$[u]_{C^{1+\alpha}(\R^n)} \leq 1$
and
\[\I (u, x)=0\quad \textrm{in}\ B_{ {1}/{\delta}}\]
in the viscosity sense.

Then, for every $\varepsilon_\circ>0$, there exists $\delta_\circ > 0$ depending only on $\eps_\circ$, $n$, $s$, $\alpha$, $\theta$, $\lambda$, and $\Lambda$, such that if $\delta < \delta_\circ$,
\[\|u-\ell\|_{C^1(B_1)} \leq \varepsilon_\circ,\]
where $\ell(x) = u(0) +x\cdot \nabla u(0) $. 
\end{prop}
 
\begin{proof}
We split the proof into two steps.

\begin{steps}
\item Let us argue by contradiction and let us assume that the statement does not hold. That is, there exists some $\eps_\circ > 0$ such that for any $k\in \N$, there are $u_k\in C_{\rm loc}^{1+\alpha}(\R^n)\cap L^1_{\omega_s}(\R^n)$ with $[u_k]_{C^{1+\alpha}(\R^n)}\le 1$,  and $\I_k\in \II_s(\lambda, \Lambda)$ of the form 
\[
\I_k (u, x) = \inf_{b\in \B_k}\sup_{a\in \A_k}\left\{-\L^{(k)}_{ab} u(x) + c^{(k)}_{ab}(x)\right\},\qquad \L^{(k)}_{ab}\in \LLL,
\]
with 
\begin{equation}
\label{eq:cthetacondition}
\sup_{(a, b)\in \A_k\times\B_k} [c^{(k)}_{ab}]_{C^\theta(\R^n)}\le \frac{1}{k},
\end{equation}
and such that
\[
\I_k (u_k, x) = 0\quad\text{in}\quad B_k
\]
but 
\[
\|u_k - \ell_k\|_{C^1(B_1)}\ge \eps_\circ,
\]
where $\ell_k(x) = u_k(0) + \nabla u_k(0)\cdot x$. If we denote $v_k := u_k - \ell_k$, we have 
\begin{equation}
\label{eq:vkzero_s}
v_k(0) = |\nabla v_k(0)| =  0,
\end{equation}
and 
\begin{equation}
\label{eq:vkzero_s2}
\|v_k\|_{C^1(B_1)}\ge \eps_\circ\quad\text{and}\quad [v_k]_{C^{1+\alpha}(\R^n)}\le 1.
\end{equation}
Up to a subsequence, we know that $v_k \to v$ in $C^1_{\rm loc}(\R^n)$ for some $v$ with $[v]_{C^{1+\alpha}(\R^n)}\le 1$ (by Arzel\`a-Ascoli), and satisfying \eqref{eq:vkzero_s} as well.

In particular, this bound implies that 
\[|v(x)|\leq |x|^{1+\alpha}\qquad \textrm{and}\qquad |\nabla v(x)|\leq |x|^{\alpha} \qquad\text{in}\quad \R^n.\]

We would like to evaluate $\I(v, x)$. 
However, when $1+\alpha\geq2s$, the function $v$ does not necessarily belong to $L^1_{\omega_s}(\R^n)$, so we cannot evaluate $\I(v, x)$.
Still, we would like to show that $v$ solves $\I (v, x)=0$ \emph{in some sense}.
This can be done as follows  (cf. the proof of Proposition~\ref{prop:compactness}).

First  notice that, since $\I_k (u_k, x)=0$ in $B_{k}$  with $\I_k$ satisfying \eqref{eq:cthetacondition}, we have that for any $h\in B_1$   (using also the ellipticity, Proposition~\ref{prop:viscosity_ellipticity})
\[
\begin{split}
\Mp \bigl(v_k(x+h)-v_k(x)\bigr) & = \Mp \bigl(u_k(x+h)-u_k(x)+h\cdot \nabla u_k(0)\bigr) \\
& = \Mp \bigl(u_k(x+h)-u_k(x)\bigr) \\
& \geq \I_k(u_k(\cdot+h), x) \\
& \geq \I_k(u_k, x+h)-\frac{1}{k}|h|^\theta = -\frac{1}{k}|h|^\theta\quad \textrm{in}\quad B_{k-1},
\end{split}
\]
and, similarly, 
\[\Mm \bigl(v_k(x+h)-v_k(x)\bigr) \leq  \frac{1}{k} |h|^\theta\quad \textrm{in}\quad B_{k-1}.\]
Since 
\[\big\|v_k(x+h)-v_k(x)\big\|_{L^\infty(B_\rho)} \leq C\|\nabla v_k\|_{L^\infty(B_{\rho+1})} \leq  C\rho^{\alpha}\]
for all $\rho>0$,   we have $\|v_k(x+h)-v_k(x)\|_{L^1_{\omega_s}(\R^n)}\leq C$ uniformly in $k$ and thus $v_k(\cdot - h) - v_k$ is converging to $v(\cdot - h) - v$ in $L^1_{\omega_s}(\R^n)$ by the dominated convergence theorem.
Therefore, we can apply the stability result from Proposition~\ref{prop:stab_super} (see also Remark~\ref{rem:stab_subsol})  and  pass to the limit the previous inequalities to get
\begin{equation}\label{owsgfvb0}
\Mp \bigl(v(x+h)-v(x)\bigr) \geq 0 \geq \Mm \bigl(v(x+h)-v(x)\bigr) \quad \textrm{in}\quad \R^n
\end{equation}
for all $h\in B_1$.

This tells us that the incremental quotients of the limiting function $v$ solve an equation with bounded measurable coefficients.

\item 
We next prove that any function $u$ satisfying \eqref{eq:vkzero_s2}-\eqref{owsgfvb0} must be affine, and thus we get a contradiction with \eqref{eq:vkzero_s2}, since $v(0)=|\nabla v(0)|=0$, \eqref{eq:vkzero_s}.
For this, the idea is to apply the estimate from Theorem \ref{C^alpha-bmc} in large balls $B_R$, with $R\to\infty$.

Indeed, given $h\in B_1$ let us define the function 
\[w(x):=\frac{v(x+h)-v(x)}{|h|}.\]
Then, by \eqref{owsgfvb0} we have
\[\Mp w \geq 0 \geq \Mm w \quad \textrm{in}\quad \R^n\]
and since $|\nabla v(x)|\leq |x|^\alpha$,
\[|w(x)|\leq 1+|x|^{\alpha} \quad \textrm{in}\quad \R^n\]
In particular, by Theorem~\ref{C^alpha-bmc}, we have
\[[w]_{C^{\gamma}(B_1)} \leq C.\]
Moreover, applying the same result to the rescaled function $w_\rho(x)=w(\rho x)/\rho^{\alpha}$, we get
\[[w]_{C^{\gamma}(B_\rho)} \leq C\rho^{\alpha-\gamma}.\]
Finally,   since $\alpha<\gamma$, letting $\rho\to\infty$ we deduce that $w\equiv {\rm constant}$ in $\R^n$, which yields that $v$ is affine  (see Lemma~\ref{it:H10}). 
This gives a contradiction with \eqref{eq:vkzero_s}-\eqref{eq:vkzero_s2}, and thus the result is proved.
\qedhere
\end{steps}
\end{proof}

We can now give the proof of the $C^{1,\alpha}$ estimates (cf. the proofs of Theorems~\ref{thm-interior-linear-2} and \ref{thm-interior-linear-Lp}, and Proposition~\ref{prop-interior-linear} in subsection~\ref{ssec:interiorproofs}):

\begin{proof}[Proof of Theorem \ref{C1alpha}] 
Let us split the proof into three steps
\begin{steps}
\item \label{it:step1c1alpha}
Given an operator $\J \in \III$ of the form \eqref{eq:cthetacondition}, let us denote 
\[
c_\theta(\J) := \sup_{(a, b)\in \A\times\B} [c_{ab}]_{C^\theta(\R^n)}.
\]
 We first claim that for any $\delta>0$ and any $u \in C^{1+\alpha}(\R^n)$ we have 
\begin{equation}\label{iasjddjdj}
[u]_{C^{1+\alpha}(B_{1/2})} \leq \delta [u]_{C^{1+\alpha}(\R^n)} + C_\delta\big(\|u\|_{L^\infty(B_1)} + \mathcal{S}(u)\big),
\end{equation}
with $C_\delta$ depending only on $n$, $s$, $\alpha$, $\delta$, $\lambda$, and $\Lambda$, where
\[
\mathcal{S}(u) := \inf\left\{c_\theta(\J) : \begin{array}{l} \J(u, x) = 0~~\text{in}~~B_1~~\text{in the viscosity sense,}\\
\text{ for some $\J\in \III$} \end{array}
\right\},
\]
and we set $\mathcal{S}(u) = \infty$ if $u\notin C(B_1)\cap L^1_{\omega_s}(\R^n)$ or if the set is empty.
Indeed, let us apply Lemma~\ref{lem-interior-blowup} with\footnote{It is important to notice that $\mathcal S$ depends only on $n$, $s$, $\alpha$, and $\Lambda$. This is what gives the dependence of the constant $C_\delta$.} $\mathcal{S}$ as above. Then, either \eqref{iasjddjdj} holds, or we have a sequence $u_k \in C^{1+\alpha}(\R^n)$ and a sequence of operators $\I_k\in \IIL$ of the form 
\[
\I_k (u, x) = \inf_{b\in \B_k}\sup_{a\in \A_k}\left\{-\L^{(k)}_{ab} u(x) + c^{(k)}_{ab}(x)\right\},\qquad \L^{(k)}_{ab}\in \LLL,
\]
such that $\I_k(u_k, x) = 0$ in $B_1$ in the viscosity sense and 
\[\frac{c_\theta(\I_k)}{[u_k]_{C^{1+\alpha}(B_{1/2})}}\le \frac{2\mathcal{S}(u_k)}{[u_k]_{C^{1+\alpha}(B_{1/2})}}  \longrightarrow 0,\]
and there are $r_k\to0$, $x_k\in B_{1/2}$, for which the rescaled functions
\[
v_k(x):= \frac{u_k(x_k+r_k x)}{r_k^{1+\alpha} [u_k]_{C^{1+\alpha}(\R^n)}}
\]
satisfy $[v_k]_{C^{1+\alpha}(\R^n)}=1$ and 
\begin{equation}
\label{eq:vkellkcontradict}
\|v_k - \ell_k\|_{C^1(B_1)} > \frac{\delta}{2},
\end{equation}
where $\ell_k$ is the first order Taylor expansion of $v_k$ at 0. 
Moreover, the functions $v_k$ solve
\[
\tilde \I_k (v_k, x) = 0\quad\text{in}\quad B_{1/(2r_k)},
\]
 where 
\[
\tilde \I_k (v, x) = \inf_{b\in \B_k}\sup_{a\in \A_k}\left\{-(\L^{(k)}_{ab})_{r_k} v(x) + r_k^{2s-1-\alpha}[u_k]^{-1}_{C^{1+\alpha}(\R^n)}c^{(k)}_{ab}(x_k+r_k x)\right\},
\]
and $\big(\L^{(k)}_{ab}\big)_{r_k}\in \LLL$ are the corresponding rescalings from \eqref{eq:scaleinvariance_comp} such that 
\[
\big(\L^{(k)}_{ab}\big)_{r_k} \big(u(x_k+r_kx)\big) = r_k^{2s}\big(\L^{(k)}_{ab} u\big)(x_k+r_kx).
\]
In particular, we have
\[
c_\theta(\tilde \I_k) = \frac{r_k^{2s+\theta-1-\alpha} c_\theta(\I_k)}{[u_k]_{C^{1+\alpha}(\R^n)}} \longrightarrow 0\quad\text{as}\quad k \to \infty.\]
This means that, if $k$ is large enough, then $v_k$ satisfies the assumptions of Proposition \ref{KrylovSafonov-Liouv}, and therefore we have
\[\|v_k-\ell_k\|_{C^1(B_1)} \leq \frac{\delta}{4}.\]
This contradicts \eqref{eq:vkellkcontradict}, and thus \eqref{iasjddjdj} is proved.

\item \label{it:step2C1alpha}
We next show that for any $u\in C_{\rm loc}^{1+\alpha}(B_1)\cap L^1_{\omega_s}(\R^n)$ such that
\[
\I(u, x) = 0\quad\text{in}\quad B_1,
\]
for some $\I\in \II_s(\lambda, \Lambda; \theta)$ of the form \eqref{eq:Iform11}-\eqref{eq:Iform22}, we have \eqref{sdnvwieuh}. We proceed as in the proof of Theorem~\ref{thm-interior-linear-2} on page~\pageref{step:main}. 

Let $\eta\in C^\infty_c(B_1)$ be such that $\eta \equiv 1 $ in $B_{1/2}$, and apply \eqref{iasjddjdj} to $u\eta$. That is, for any $\delta > 0$ there exists some $C_\delta$ such that 
\[
[u]_{C^{1+\alpha}(B_{1/2})}\le \delta [\eta u]_{C^{1+\alpha}(B_1)} + C_\delta \left(\|u\|_{L^\infty(B_1)}+\mathcal{S} (\eta u)\right).
\]
Now, since $\I(u, x) = 0$ in $B_1$, we have that 
\[
\tilde \I(\eta u, x) = 0\quad\text{in}\quad B_1,
\]
where $\tilde \I\in \II_s(\lambda, \Lambda; \theta)$ is of the form 
\[
\tilde \I (v, x) = \inf_{b\in \B}\sup_{a\in \A}\big\{-\L_{ab} v(x) + \tilde c_{ab}(x)\big\},\qquad \L_{ab}\in \LL_s(\lambda, \Lambda; \theta),
\]
and $\tilde c_{ab}(x)$ can be expressed in terms of $c_{ab}$ from the definition of $\I$ as
\[
\tilde c_{ab}(x) = c_{ab}(x) + \L_{ab}(\eta u - u)(x). 
\]
In particular, by Lemma~\ref{lem:Lu_2_v} and since $\eta u - u \equiv 0$ in $B_2$, we have
\[
[\tilde c_{ab}]_{C^\theta(B_1)} \le [c_{ab}]_{C^\theta(B_1)} +C\|(\eta-1) u\|_{L^1_{\omega_s}(\R^n)},
\]
and taking a $C^\theta$ extension to $\R^n$ we obtain 
\[
\mathcal{S}(\eta u) \le c_\theta(\tilde \I) \le c_\theta(\I) +C\|(\eta-1)u\|_{L^1_{\omega_s}(\R^n)}\le C_\circ + C\|u\|_{L^1_{\omega_s}(\R^n)}.
\]
Thus,  up to making $\delta$ smaller (see Proposition~\ref{prop:A_imp}) we get 
\[
[u]_{C^{1+\alpha}(B_{1/2})}\le \delta [ u]_{C^{1+\alpha}(B_1)} + C_\delta \left(\|u\|_{L^\infty(B_1)}+\|u\|_{L^1_{\omega_s}(\R^n)}+C_\circ \right),
\]
which, as in \ref{step:SAL} in the proof of Theorem~\ref{thm-interior-linear-2}, implies
\[
[u]_{C^{1+\alpha}(B_{1/2})}\le C \left(\|u\|_{L^\infty(B_1)}+\|u\|_{L^1_{\omega_s}(\R^n)}+C_\circ \right).
\]
Now by a rescaling and covering argument we can actually write the bound in $B_{1/2}$ and $B_{3/4}$ as 
\[
[u]_{C^{1+\alpha}(B_{1/2})}\le C \left(\|u\|_{L^\infty(B_{3/4})}+\|u\|_{L^1_{\omega_s}(\R^n)}+C_\circ \right).
\]
We finish by noticing that
\[
\Mp u \ge \I(u, x) -\I(0, x) = -\I(0, x) \ge\Mm u 
\]
in the viscosity sense, and hence by Theorem~\ref{C^alpha-bmc_2} (again, after a covering argument) we get
\[
\|u\|_{L^\infty(B_{3/4})}\le C \left(\|u\|_{L^1_{\omega_s}(\R^n)}+\|\I(0, x)\|_{L^\infty(B_1)}\right).
\]
Since $\|\I(0, x)\|_{L^\infty(B_1)}\le |\I(0, 0)| + C_\circ$ this implies 
\begin{equation}
\label{eq:apriori_visc}
[u]_{C^{1+\alpha}(B_{1/2})}\le C \left(\|u\|_{L^1_{\omega_s}(\R^n)}+|\I(0, 0)|+C_\circ \right)
\end{equation}
for any $u\in C^{1+\alpha}(B_1)\cap L^1_{\omega_s}(\R^n)$, as we wanted.

\item \label{it:step2_C1alpha}
We finally show by approximation that \eqref{sdnvwieuh} holds for any viscosity solution $u\in  C(B_1)\cap L^1_{\omega_s}(\R^n)$ to
$
\I(u, x) = 0
$
in $B_1$, with $\I\in \II_s(\lambda, \Lambda; \theta)$ of the form \eqref{eq:Iform11}-\eqref{eq:Iform22}. Indeed, by Proposition~\ref{prop:regularization} we can find a sequence $u^{(\eps)}\in C^{2s+\theta}_{\rm loc}(B_{3/4})$ such that
\[
\left\{
\begin{array}{rcll}
\I_\eps(u^{(\eps)}, x) & = & 0 & \quad\text{in}\quad B_{3/4}\\
u^{(\eps)} & = & u & \quad\text{in}\quad \R^n\setminus B_{3/4},
\end{array}
\right.
\]
for some sequence of operators $\I_\eps$ that satisfy the same hypotheses as $\I$ (up to universal constants), and 
\[
u^{(\eps)}\to u\quad\text{locally uniformly in $B_{3/4}$}.
\]
In particular, since $\theta = (1+\alpha -2s)_+$, we have that $2s+\theta \ge 1+\alpha$ and we can apply the a priori estimates from \eqref{eq:apriori_visc} to deduce
\[
[u^{(\eps)}]_{C^{1+\alpha}(B_{1/2})}\le C \left(\|u^{(\eps)}\|_{L^1_{\omega_s}(\R^n)}+|\I(0, 0)|+C_\circ \right)
\]
(where we also used that, from Proposition~\ref{prop:regularization}, $\I_\eps(0,0) = \I(0, 0)$). Moreover, thanks to \eqref{eq:bound_reg} we get 
\[
\|u^{(\eps)}\|_{C^{1+\alpha}(B_{1/2})}\le C \left(\|u\|_{L^1_{\omega_s}(\R^n)}+|\I(0, 0)|+C_\circ \right).
\]
Since the right-hand side is independent of $\eps$, and $u^{(\eps)}$ converges locally uniformly to $u$ in $B_{3/4}$, the limit has the same bound, 
\[
\|u\|_{C^{1+\alpha}(B_{1/2})}\le C \left(\|u\|_{L^1_{\omega_s}(\R^n)}+|\I(0, 0)|+C_\circ \right),
\]
as we wanted to see. 
\qedhere
 \end{steps}
\end{proof}

\subsection{$C^{2s+\alpha}$ regularity for concave equations}\label{ssec:C1alpha2}

 \index{Interior regularity!Fully nonlinear!Concave equations}
 
 \index{Evans-Krylov}

We will now consider operators of the form 
\begin{equation}
\label{eq:concaveI}
\I (u, x) = \inf_{a\in \A}\big\{-\L_{a} u(x) + c_{a}(x)\big\},\qquad \L_{a}\in \LL_s(\lambda, \Lambda; \alpha),
\end{equation}
for some $\alpha\in [0, 1)$. Observe that, with this definition, the operators satisfy
\begin{equation}
\label{eq:concavity}
\I\big(t u_1 + (1-t) u_2,x\big) \ge t \I(u_1, x)+(1-t)\I(u_2, x)\quad\text{for}\quad t\in (0, 1),
\end{equation}
that is, they are \emph{concave}. 

The next result shows that solutions are $C^{2s+\alpha}$ for some small $\alpha>0$, and in particular they are strong solutions.
This is the nonlocal analogue of the celebrated Evans--Krylov theorem.

\begin{thm}\label{EvansKrylov}
 Let $s\in(\frac12,1)$, and $0 <\lambda \le\Lambda$. There exists $\alpha_\circ> 0$ depending only on $n$, $s$, $\lambda$, and $\Lambda$,   such that the following statement holds:
 
 Let $\alpha \in (0, \alpha_\circ)$, and let $\I \in \II_s(\lambda, \Lambda; \alpha)$ be of the form \eqref{eq:concaveI} and satisfying
\[
\sup_{a\in \A}\, [c_a]_{C^\alpha(\R^n)}\le C_\circ.
\]

Assume that $u\in C(B_1)\cap L^1_{\omega_s}(\R^n)$ satisfies
\[\I (u, x)=0\quad \textrm{in}\ B_1\]
in the viscosity sense. Then, $u\in C^{2s+\alpha}_{\rm loc}(B_1)$ and 
\[\|u\|_{C^{2s+\alpha}(B_{1/2})} \leq C\left(\|u\|_{L^1_{\omega_s}(\R^n)}+|\I(0, 0)|+C_\circ\right),\]
for some constant $C$ depending only on $n$, $s$, $\alpha$, $\lambda$, and $\Lambda$.
\end{thm}

Notice that we prove the $C^{2s+\alpha}$ estimate for $s > \frac12$; otherwise the $C^{1,\alpha}$ estimate from Theorem~\ref{C1alpha} is stronger.

The proof of this result will also be based on a blow-up and compactness argument.
However, in this case the Liouville-type theorem we need to prove is highly nontrivial (it was proved in \cite{CS2,Ser2}). 
As in the $C^{1,\alpha}$ regularity, such Liouville theorem holds for operators with no regularity assumption on the kernels, 
\begin{equation}
\label{eq:Ioftheform_conc}
\I (u, x) = \inf_{a\in \A}\big\{-\L_{a} u(x) + c_{a}(x)\big\},\qquad \L_{a}\in \LLL.
\end{equation}
\begin{prop}
\label{prop:liouville_concave} \index{Liouville's theorem!Fully nonlinear!Concave equations}
 Let $s\in(\frac12,1)$ and $0<\lambda \le \Lambda$. There exists $\alpha_\circ> 0$ depending only on $n$, $s$, $\lambda$, and $\Lambda$,   such that the following statement holds: 
 
Let $\delta > 0$, and $\alpha\in (0, \alpha_\circ)$ with that $2s+\alpha\notin\N$. Let $\I\in \III$ be a concave operator of the form \eqref{eq:Ioftheform_conc} and satisfying 
\[
\sup_{ a \in \A} \, [c_{a}]_{C^\alpha(\R^n)}\le \delta.
\]

Assume that $u\in C_{\rm loc}^{2s+\alpha}(\R^n)\cap L^1_{\omega_s}(\R^n)$ satisfies
$[u]_{C^{2s+\alpha}(\R^n)} \leq 1$
and
\[\I (u, x)=0\quad \textrm{in}\ B_{ {1}/{\delta}}\]
in the viscosity sense.

Then, for every $\varepsilon_\circ>0$, there exists $\delta_\circ > 0$ depending only on $\eps_\circ$, $n$, $s$, $\alpha$, $\lambda$, and $\Lambda$, such that if $\delta < \delta_\circ$,
\[\|u-p\|_{C^\nu(B_1)} \leq \varepsilon_\circ,\]
where $p$ is the  Taylor polynomial of $u$ at 0 of order $\nu := \lfloor 2s+\alpha\rfloor$. 
\end{prop}
\begin{proof}
The value of $\alpha_\circ$ will be fixed later in the proof. We divide the proof into five steps. 
\begin{steps}
\item Let us argue by contradiction and let us assume that the statement does not hold for a given $\alpha> 0$. That is, there exists some $\eps_\circ > 0$ such that for any $k\in \N$, there are $u_k\in C_{\rm loc}^{2s+\alpha}(\R^n)\cap L^1_{\omega_s}(\R^n)$ with $[u_k]_{C^{2s+\alpha}(\R^n)}\le 1$,  and $\I_k\in \II_s(\lambda, \Lambda)$ of the form 
\begin{equation}
\label{eq:cthetacondition_conc0}
\I_k (u, x) = \inf_{a\in \A_k}\left\{-\L^{(k)}_{a} u(x) + c^{(k)}_{a}(x)\right\},\qquad \L^{(k)}_{a}\in \LLL,
\end{equation}
with 
\begin{equation}
\label{eq:cthetacondition_conc}
\sup_{a\in \A_k} [c^{(k)}_{a}]_{C^\alpha(\R^n)}\le \frac{1}{k},
\end{equation}
and such that
\[
\I_k (u_k, x) = 0\quad\text{in}\quad B_k
\]
but 
\[
\|u_k - p_k\|_{C^\nu(B_1)}\ge \eps_\circ,
\]
where $p_k$ is the  Taylor polynomial of $u_k$ at 0 of order $\nu$ (i.e., either linear or quadratic). If we denote $v_k := u_k - p_k$, we have 
\begin{equation}
\label{eq:vkzero_s_conc}
v_k (0) = |\nabla v_k (0)| = |D^\nu v_k (0)| = 0, 
\end{equation}
and 
\[
\|v_k\|_{C^\nu(B_1)}\ge \eps_\circ\quad\text{and}\quad [v_k]_{C^{2s+\alpha}(\R^n)}\le 1.
\]
For any $\alpha'< \alpha$, up to a subsequence we know that $v_k \to v$ in $C^{2s+\alpha'}_{\rm loc}(\R^n)$, for some $v$ with 
\begin{equation}
\label{owsgfbb}
\|v \|_{C^\nu(B_1)}\ge \eps_\circ\quad\text{and}\quad [v]_{C^{2s+\alpha}(\R^n)}\le 1,
\end{equation}
  and satisfying \eqref{eq:vkzero_s_conc} as well, 
\begin{equation}
\label{eq:vkzero_s_conc0}
v(0) = |\nabla v(0)| = |D^\nu v (0)| = 0.
\end{equation}
 
In particular, this bound implies that 
\begin{equation}
\label{eq:incr_quotients}
|v(x)|\leq |x|^{2s+\alpha}\qquad \textrm{and}\qquad |\nabla v(x)|\leq |x|^{2s+\alpha-1} \qquad\text{in}\quad \R^n. 
\end{equation}

Now, as in the proof of Proposition~\ref{KrylovSafonov-Liouv}, the function $v$ does not necessarily belong to $L^1_{\omega_s}(\R^n)$, so we cannot  evaluate $\I(v, x)$. We proceed by taking incremental quotients instead, which thanks to \eqref{eq:incr_quotients} belong to $L^1_{\omega_s}(\R^n)$. 
Thus,  proceeding exactly as in the proof of Proposition \ref{KrylovSafonov-Liouv}, for any $h\in B_1$ we have
\begin{equation}\label{owsgfvb}
\Mp\bigl(v(x+h)-v(x)\bigr) \geq 0 \geq \Mm\bigl(v(x+h)-v(x)\bigr) \quad \textrm{in}\quad \R^n.
\end{equation}
This somehow tells us that the limiting function $u$ solves a fully nonlinear equation. However we still have not  used the fact that $\I$ is a \textit{concave} operator.
For this, notice that for any nonnegative $\mu\in L^1(\R^n)$ with compact support in $B_\rho$ and $\int_{\R^n}\mu=1$, we have by concavity of $\I_k$, \eqref{eq:concavity} (using Jensen's inequality for integrals, since all the terms are well-defined)
\[
\I_k \left(\int_{\R^n} u_k(\,\cdot\,+h)d\mu(h), x\right)\ge \int_{\R^n} \I_k \big(u_k(\,\cdot\, +h), x\big)d\mu(h)\quad\text{in}\quad B_{k-\rho}.
\]
Notice, also, that from \eqref{eq:cthetacondition_conc0}-\eqref{eq:cthetacondition_conc}, we know 
\[
\I_k \big(u_k(\,\cdot\, +h), x\big)\ge \I_k (u_k, x) - \frac{1}{k}|h|^\alpha = - \frac{1}{k}|h|^\alpha\quad\text{in}\quad B_{k-|h|}.
\]
On the other hand, by ellipticity (Proposition \ref{prop:viscosity_ellipticity}) we have 
\[
\begin{split}
\Mp\left(\int_{\R^n} u_k(x+h)\,d\mu(h) - u_k(x) \right)&  \ge \I_k \left(\int_{\R^n} u_k(\,\cdot\,+h)d\mu(h), x\right) \\
&  \ge - \frac{1}{k}\int_{\R^n}|h|^\alpha \, d\mu(h) \quad\text{in}\quad B_{k-\rho}.
\end{split}
\]
Finally, notice that since $2s > 1$ this implies 
\[
\begin{split}
\Mp\left(\int_{\R^n} v_k(x+h)\,d\mu(h) - v_k(x) \right)\ge - \frac{1}{k}\int_{\R^n}|h|^\alpha \, d\mu(h) \quad\text{in}\quad B_{k-\rho}
\end{split}
\]
as well.

Now, exactly as in the proof of Proposition~\ref{KrylovSafonov-Liouv}, we have that the functions $\int_{\R^n} v_k(x+h)d\mu(h) - v_k(x)$ are uniformly bounded in $L^1_{\omega_s}(\R^n)$ and converge locally uniformly and in $L^1_{\omega_s}(\R^n)$ to $\int_{\R^n} v(x+h)d\mu(h) - v(x)$. We can therefore apply the stability of viscosity solutions under uniform limits from Proposition~\ref{prop:stab_super} to pass the previous inequality to the limit, to get
\[\Mp \left(\int_{\R^n} v(x+h)d\mu(h) - v(x)\right)\geq 0 \quad \textrm{in}\quad \R^n\]
for any nonnegative $\mu\in L^1(\R^n)$ with compact support and $\int_{\R^n}\mu=1$.
Taking symmetric functions $\mu(h)=\mu(-h)$, this is equivalent to
\begin{equation}\label{owsgfvv}
\Mp\left(\int_{\R^n} \big(v(x+h)+v(x-h)-2v(x)\big)d\mu(h)\right)\geq 0 \qquad \textrm{in}\quad \R^n.
\end{equation}

\item 
We want to prove that any function $v$ satisfying \eqref{owsgfbb}-\eqref{owsgfvb}-\eqref{owsgfvv} must be a polynomial of degree $\nu$, thus reaching a contradiction with \eqref{owsgfbb}-\eqref{eq:vkzero_s_conc0}.

For this, our goal will be to prove that there exists some $\alpha_\circ > 0$ (that will be independent of $\alpha$), such that 
\begin{equation}\label{clailm}
[v]_{C^{2s+\alpha_\circ}(B_\rho)} \leq C\rho^{\alpha-\alpha_\circ}
\end{equation}
for all $\rho\geq1$. In particular, when $\alpha < \alpha_\circ$ we will get our desired result by letting $\rho\to \infty$.

Let us define
\[\delta^2_h v(x):= \frac{v(x+h)+v(x-h)}{2} - v(x),\]
and 
\[P(x):= \int_{\R^n} \big(\delta^2_h v(x) - \delta^2_h v(0)\big)_+\frac{dh}{|h|^{n+2s}}, \]
\[N(x):= \int_{\R^n} \big(\delta^2_h v(x) - \delta^2_h v(0)\big)_-\frac{dh}{|h|^{n+2s}}. \]

When $2s+\alpha<2$, since $[u]_{C^{2s+\alpha}(\R^n)}\leq 1$, we have 
\[\big|\delta^2_h v(x) - \delta^2_h v(0)\big| \leq 2 \min\left\{|h|^{2s+\alpha},\,|h|^{2s+\alpha-1}|x|\right\}\]
for all $x,h\in\R^n$.
When $2s+\alpha>2$, instead, we have that  
\[\big|\delta^2_h v(x) - \delta^2_h v(0)\big| \leq \min\left\{|h|^2|x|^{2s+\alpha-2},\,|h|^{2s+\alpha-1}|x|\right\}.\]
(See Lemma~\ref{lem:A_imp_2}.)  In both cases, a simple computation in polar coordinates gives that
\begin{equation}
\label{eq:frompn}
|P(x)|+|N(x)|\leq C_1 |x|^\alpha \quad \textrm{in}\quad \R^n,
\end{equation}
with $C_1$ depending only on $n$, $s$, $\alpha$, $\lambda$, and $\Lambda$.
Moreover, dividing $v$ by a constant if necessary, we may assume $C_1 =1$, 
\begin{equation}
\label{eq:frompn0}
|P(x)|+|N(x)|\leq |x|^\alpha \quad \textrm{in}\quad \R^n.
\end{equation}

We next want to show that 
\begin{equation}\label{pqpwpwe}
0 \leq P \leq 2^{-k\alpha_\circ} \quad \textrm{in} \quad B_{2^{-k}}
\end{equation}
for all $k\in \mathbb Z$, for some $\alpha_\circ>0$ depending only on $n$, $s$, $\lambda$, $\Lambda$ (but not on $\alpha$). We only need to prove it for $k > 0$, since for $k \le 0$ we already know it by \eqref{eq:frompn0} (since we will have $\alpha < \alpha_\circ$). It suffices to show that
\begin{equation}\label{pqpwpwf}
P\leq 1-\theta\quad \textrm{in}\quad B_{1/2},
\end{equation}
for some $\theta>0$ depending only on $n$, $s$, $\lambda$, and $\Lambda$ (in particular, independent of $\alpha$).
Once we have this, \eqref{pqpwpwe} follows for all $k\geq1$ by iteration (for some small $\alpha_\circ>0$ that depends on $\theta$, namely $1-\theta = 2^{-\alpha_\circ}$), and the same bound for $N$ is completely analogous.

To show \eqref{pqpwpwf}, let $x_\circ\in B_{1/2}$ be such that $P(x_\circ)=\max_{B_{1/2}}P$, and let 
\[U := \big\{h\in \R^n : \delta^2_h v(x_\circ) > \delta^2_h v(0)\big\} = -U.\]
In particular, we have
\[ P(x_\circ) = \int_U \big(\delta^2_h v(x_\circ) - \delta^2_h v(0)\big)\frac{dh}{|h|^{n+2s}},\]
\[ N(x_\circ) = \int_{U^c} \big(\delta^2_h v(0) - \delta^2_h v(x_\circ)\big)\frac{dh}{|h|^{n+2s}}.\]
Let us define
\[ 
\begin{split}
w(x) &:= \int_U \big(\delta^2_h v(x) - \delta^2_h v(0)\big)\frac{dh}{|h|^{n+2s}},\\
\bar w(x)& := \int_{U^c} \big(\delta^2_h v(x) - \delta^2_h v(0)\big)\frac{dh}{|h|^{n+2s}},
\end{split}\]
and notice that  
\[  w \le P \quad\text{in}\quad \R^n, \quad  P\leq 1 \quad \textrm{in}\quad B_1, \quad \textrm{and} \quad w(x_\circ)=P(x_\circ).\]
Let $\mu>0$ to be chosen later, and define the set
\[D:= \{x\in B_1 : w\geq 1-\mu\}.\]

Our next goal will be to prove the following: 

\noindent\textbf{Claim}. There is a constant $\eta>0$ depending only on $n$, $s$, $\lambda$, and $\Lambda$, for which we have
\[|D|\leq (1-\eta)|B_1|.\]

\item Before proving the claim, let us first observe that $w$ and $\bar w$ are subsolutions. Indeed, by \eqref{owsgfvv} we have
\[\Mp \left( \int_{\R^n} \delta^2_h v(x) d\mu(h)\right) \geq 0\quad \textrm{in} \quad \R^n, \]
for any $\mu \in L^1(\R^n)$ symmetric, with compact support, $\mu \ge 0$. We now want to let  $\mu \rightharpoonup |h|^{-n-2s}\chi_{U^c}$. If we define $\mu_\eps := |h|^{-n-2s}\chi_{U^c}(h)\chi_{B_{1/\eps}\setminus B_\eps}(h)\, dh$ (which is a symmetric measure) and $\mu_0 := |h|^{-n-2s}\chi_{U^c}$, then 
\[
\int_{\R^n}\delta^2_h v(x) \big(d\mu_\eps(h)-d\mu_0(h)\big)\to 0
\]
locally uniformly in $\R^n$, using that $[v]_{C^{2s+\alpha}(\R^n)}\le 1$. Moreover, 
\[
\int_{\R^n}\delta^2_h v(x) d\mu_\eps(h) \le  C_h |x|,
\]
and therefore they are uniformly in $L^\infty_{2s-\eps'}(\R^n)$ (for $\eps' < 2s-1$). Hence, by the dominated convergence theorem, the function $\int_{\R^n}\delta^2_h v(x) d\mu_\eps(h)$ converges to $\int_{\R^n}\delta^2_h v(x) d\mu_0(h)$ in $L^1_{\omega_s}(\R^n)$ as well. We can therefore apply the stability of subsolutions (see Remark~\ref{rem:stab_subsol}) to deduce that 
\[
\Mp w \ge 0 \quad\text{and}\quad \Mp \bar w \ge 0\quad\text{in}\quad \R^n.
\]

Observe, also, that if we define $\varphi(x) := v(x+x_\circ) - v(x)$ for some $x_\circ\in \R^n$ fixed, then by \eqref{owsgfvb} we know that
\[
\Mp \varphi \ge 0\quad\text{in}\quad \R^n. 
\]
Since $\varphi\in C^{2s+\alpha}(\R^n)$, we can use the explicit formula \eqref{eq:Mpexplicit} for $\Mp$, to deduce that 
\[
\int_{\R^n} \left(\Lambda (\delta^2_h \varphi(x))_+ - \lambda (\delta^2_h \varphi(x))_- \right) |h|^{-n-2s}\, dh \ge 0. 
\]
Evaluating at $x = 0$ this implies $\Lambda P(x) - \lambda N(x) \ge 0$, that is
\[
\frac{\Lambda}{\lambda} P(x) \ge N(x).
\]
By using the other inequality, $\Mm \varphi\le 0$ in $\R^n$, we deduce
\begin{equation}
\label{eq:PNineq}
\frac{\lambda}{\Lambda} P (x) \le N(x) \le \frac{\Lambda}{\lambda} P(x)\quad\text{in}\quad \R^n. 
\end{equation}

\item Let us now  prove the claim by contradiction. Assume that $|D|\geq (1-\eta)|B_1|$, with $\eta>0$ to be chosen later.
Notice that this means that $w\geq 1-\mu$ in most of $B_1$, and in particular
\[0 \leq P-w \leq \mu \quad \textrm{in}\quad D.\]
By definition, we have $P-N=w+\bar w$. Together with \eqref{eq:PNineq}, this implies 
\[\bar w =(P-w)-N \leq \mu-\frac{\lambda}{\Lambda} P \leq \mu-\frac{\lambda}{\Lambda} (1-\mu)=-c_\circ \quad \textrm{in}\quad D,\]
where we may choose for example $\mu=\frac{\lambda}{4\Lambda}$ and $c_\circ=2\mu$.

We now  use the $L^\infty$ bound for subsolutions on $\bar w$, Theorem \ref{half-Harnack-sub}, to finish the proof of the claim.
Let $r_\circ>0$ small enough, and consider the function
\[w_\circ (x) := \big(\bar w(r_\circ x)+c_\circ\big)_+,\]
which  still satisfies $\Mp w_\circ \geq 0$ in $\R^n$ in the viscosity sense (see Remark~\ref{rem:max_min}).

Notice that 
\[w_\circ \equiv0 \quad \textrm{in}\quad {\textstyle\frac{1}{r_\circ}}D\subset B_{1/r_\circ},\]
and 
\[\|w_\circ\|_{L^\infty(B_\rho)} \leq \|P\|_{L^\infty(B_{r_\circ \rho})} + c_\circ 
\leq (r_\circ \rho)^{\alpha_\circ}+c_\circ\]
for all $\rho\geq 1/r_\circ$, from \eqref{eq:frompn}. Therefore, since $|B_{1/r_\circ} \setminus ({\textstyle\frac{1}{r_\circ}}D)| \leq \eta|B_{1/r_\circ}| $, we   have that  
\[\begin{split}
\int_{\R^n} \frac{w_\circ(x)}{1+|x|^{n+2s}}\, dx & \leq (1+c_\circ) \eta|B_{1/r_\circ}|+ \int_{B^c_{1/r_\circ}} \frac{(r_\circ |x|)^{\alpha_\circ}+c_\circ}{|x|^{n+2s}}\,dx \\
& \leq C(\eta r_\circ^{-n} + r_\circ^{2s}).
\end{split}\]
Hence, by Theorem  \ref{half-Harnack-sub}, we deduce that 
\[\sup_{B_{1/2}} w_\circ \leq C\eta r_\circ^{-n} + Cr_\circ^{2s}, \]
with $C$ depending only on $n$, $s$, $\lambda$, and $\Lambda$.
Choosing now $r_\circ$ small enough, and then $\eta$ small enough, we deduce that
\[w_\circ(0) \leq \frac{c_\circ}{2},\]
which yields $\bar w(0)<0$, a contradiction since $\bar w(0)=0$ by definition.
Hence, the claim is proved.

\item Note now that what we proved in the claim is equivalent to
\begin{equation}\label{claim-bis2}
\big|\{x\in B_1 : w< 1-\mu\}\big| \geq \eta|B_1|.
\end{equation}
Since $w$ is a subsolution,
\[\Mm (1-w)_+ \leq \Mp (1-w)_- \quad \textrm{in}\quad \R^n.\]
Moreover, since $1 - w \ge 0$ in $B_1$ and $w \le P\le |x|^{\alpha_\circ}$ in $\R^n\setminus B_1$ (by \eqref{eq:frompn0}), we have that for any $x\in B_{3/4}$, 
\[
\Mp (1-w)_-(x) \le \Lambda \int_{\R^n\setminus B_1} \frac{(1-w(z))_- }{|x-z|^{n+2s}}\, dz \le C \Lambda \int_{\R^n\setminus B_1} \frac{(1-|z|^{\alpha_\circ})_- }{1+|z|^{n+2s}}\, dz.
\]
Observe that the last term goes to zero as $\alpha_\circ\downarrow 0$, so that combined with the previous inequality, for any $\delta_\circ > 0$ we can find $\alpha_\circ > 0$ such that 
\[
\Mm (1-w)_+ \le \delta_\circ\quad\text{in}\quad B_{3/4}. 
\]
We   apply Theorem \ref{half-Harnack-sup}, to obtain that
\[\inf_{B_{1/2}} (1-w)_+ \geq c\|(1-w)_+\|_{L^1_{\omega_s}(\R^n)} - \delta_\circ \geq \frac{c\eta\mu}{2} =: \theta,\]
where we used \eqref{claim-bis2} and we have chosen $\delta_\circ \le \frac{c\eta\mu}{2}$. 
This implies
\[w(x_\circ)\leq 1-\theta\quad \textrm{in}\quad B_{1/2}\]
and, since $P\leq P(x_\circ)=w(x_\circ)$ in $B_{1/2}$,   \eqref{pqpwpwf} follows.
Thus, we have proved \eqref{pqpwpwe}, or equivalently
\[P(x)\leq C |x|^{\alpha_\circ}\quad \textrm{in}\quad B_1.\]
Moreover, the same bound holds for $N(x)$.

Note that, for $\tau \in \R^n$,
\[-c_{n,s}\big(P(\tau)-N(\tau)\big) = (-\Delta)^s\big(v(\,\cdot\, +\tau)-v\big)(0). \]
Furthermore, the point 0 in the definition of $P$ and $N$ may be replaced by any other point $z\in B_{1/2}$, and hence we have proved that
\[\big|(-\Delta)^s\big(v(\cdot +\tau)-v\big)(z)\big| \leq C|\tau|^{\alpha_\circ}\quad \textrm{for all}\quad z\in B_{1/2}. \]
This means that 
\[\left|(-\Delta)^s\left(\frac{v(x+\tau)-v(x)}{|\tau|^{\alpha_\circ}}\right)\right| \leq C \quad \textrm{in}\quad B_{1/2}. \]

On the other hand, notice also that from \eqref{eq:incr_quotients} we know  
\[
\frac{|v(x+\tau)-v(x)|}{|\tau|^{\alpha_\circ}} \le |\tau|^{1-\alpha_\circ}(1+|x|^{2s+\alpha-1})\quad\text{in}\quad \R^n,
\]
and hence, $\frac{v(x+\tau)-v(x)}{|\tau|^{\alpha_\circ}}\in L^\infty_{2s-\eps}(\R^n)$ for $\eps < 1-\alpha$ and independently of $\tau$. By Theorem~\ref{thm-interior-linear-Lp}, we deduce that
\[\left\|\frac{v(x+\tau)-v(x)}{|\tau|^{\alpha_\circ}}\right\|_{C^{2s}(B_{1/4})} \leq C,\]
and therefore (see Lemma~\ref{it:H8})
\[\|v\|_{C^{2s+\alpha_\circ}(B_{1/4})} \leq C.\]

The whole argument in the previous steps can now be applied to every scale $\rho\geq1$, i.e., to the rescaled functions 
\[v_\rho(x):= \frac{v(\rho x)}{\rho^{2s+\alpha}},\]
which satisfy the same assumptions as $v$.
Doing so, we find
\[[v_\rho]_{C^{2s+\alpha_\circ}(B_{1/4})} \leq C.\]
Rescaling back to $v$, we get
\[[v]_{C^{2s+\alpha_\circ}(B_{\rho/4})} \leq C\rho^{\alpha-\alpha_\circ}\quad \textrm{for all} \quad \rho\geq1,\]
and letting $\rho\to\infty$ we deduce that $v$ must be a polynomial of degree $\nu$. Together with \eqref{eq:vkzero_s_conc0}, this implies $v \equiv 0$, a contradiction with \eqref{owsgfbb}.
Hence the proposition is proved. \qedhere
\end{steps}
\end{proof}

Once we have a Liouville-type theorem, the proof of the interior regularity follows as in the proof of Theorem~\ref{C1alpha}:

\begin{proof}[Proof of Theorem~\ref{EvansKrylov}]
The proof follows exactly as   the proof of Theorem~\ref{C1alpha} on page~\pageref{it:step1c1alpha}, thanks to the Liouville-type theorem, Proposition~\ref{prop:liouville_concave}: 

Using the same notation as in there, we first show that for any $\delta>0$  and $u\in C^{2s+\alpha}(\R^n)$, 
\begin{equation}\label{iasjddjdj2}
[u]_{C^{2s+\alpha}(B_{1/2})} \leq \delta [u]_{C^{2s+\alpha}(\R^n)} + C_\delta\big(\|u\|_{L^\infty(B_1)} + \mathcal{S}(u)\big),
\end{equation}
with $C_\delta$ depending only on $n$, $s$, $\alpha$, $\delta$, $\lambda$, and $\Lambda$, where
\[
\mathcal{S}(u) := \inf\big\{c_\alpha(\J) : \J(u, x) = 0~~\text{in}~~B_1~\text{ for some $\J\in \III$ concave}
\big\},
\]
and $\mathcal{S}(u) = \infty$ if $u\notin C(B_1)\cap  L^1_{\omega_s}(\R^n)$ or if the set is empty.  By Lemma~\ref{lem-interior-blowup}   either \eqref{iasjddjdj2} holds, or we have sequences $u_k \in C^{2s+\alpha}(\R^n)$,  $r_k\to0$, and $x_k\in B_{1/2}$, for which the rescaled functions
\[
v_k(x):= \frac{u_k(x_k+r_k x)}{r_k^{2s+\alpha} [u_k]_{C^{2s+\alpha}(\R^n)}}
\]
satisfy $[v_k]_{C^{2s+\alpha}(\R^n)}=1$,
$
\|v_k - p_k\|_{C^\nu(B_1)} > \frac{\delta}{2},
$ 
(where $p_k$ is the $\nu$-th order Taylor expansion of $v_k$ at 0, $\nu = \lfloor 2s+\alpha\rfloor$), and 
\[
\tilde \I_k (v_k, x) = 0\quad\text{in}\quad B_{1/(2r_k)}
\]
for some $\tilde \I_k\in \III$ concave such that $c_\alpha(\tilde \I_k) \rightarrow 0$ as $k \to \infty$; in particular, contradicting Proposition~\ref{prop:liouville_concave}, and thus proving  \eqref{iasjddjdj2}. 

Once we have \eqref{iasjddjdj2} we can directly show 
\begin{equation}
\label{eq:apriori_visc2}
[u]_{C^{2s+\alpha}(B_{1/2})}\le C \left(\|u\|_{L^1_{\omega_s}(\R^n)}+|\I(0, 0)|+C_\circ \right)
\end{equation}
for any $u\in C^{2s+\alpha}(B_{1})\cap L^1_{\omega_s}(\R^n)$ with $\I(u, x) = 0$ in $B_1$, proceeding as in \ref{it:step2C1alpha} of the proof of Theorem~\ref{C1alpha} (cf. proof of Theorem~\ref{thm-interior-linear-2} on page~\pageref{step:main}). Finally, thanks to the approximation result, Proposition~\ref{prop:regularization}, we can again approximate any function $u\in C(B_1)\cap L^1_{\omega_s}(\R^n)$  by functions $u^{(\eps)}\in C_{\rm loc}^{2s+\alpha}(B_{3/4})$ such that
\[
\left\{
\begin{array}{rcll}
\I_\eps(u^{(\eps)}, x) & = & 0 & \quad\text{in}\quad B_{3/4}\\
u^{(\eps)} & = & u & \quad\text{in}\quad \R^n\setminus B_{3/4},
\end{array}
\right.
\]
for some sequence of operators $\I_\eps$ that satisfy the same hypotheses as $\I$ (up to universal constants), and $u^{(\eps)}\to u$ locally uniformly in $B_{3/4}$. Thus, as in \ref{it:step2_C1alpha}, we obtain \eqref{eq:apriori_visc2} for any viscosity solution $u\in C(B_1)\cap L^1_{\omega_s}(\R^n)$. 
\end{proof}

\section{Further results and open problems}

In this chapter we have established the main known interior regularity results for solutions to fully nonlinear integro-differential equations of order $2s$.

In particular, two of the main results yield that solutions to \emph{concave} equations 
\begin{equation}\label{question3.1}
\inf_{a\in \A}\big\{-\L_a u+c_a\big\} =0\quad\textrm{in}\quad B_1
\end{equation}
are $C^{\max\{1,2s\}+\alpha}$ in $B_1$ for all $s\in (0,1)$.
A first question that arises immediately is the following:

\vspace{2mm}

\noindent\textbf{Open question 3.1}:
\textit{If the kernels of $\L_a$ are regular enough, can one prove that solutions to \eqref{question3.1} are actually $C^{1+s+\alpha}$? Or even $C^{2+\alpha}$ for every $s\in(\frac12,1)$?}

\vspace{2mm}

By analogy with the case $s=1$, we do not expect solutions to these equations to be more regular than $C^{1+2s}$.

\subsection{Uniform estimates as $s\uparrow1$}

An important feature of the regularity theory developed by Caffarelli and Silvestre in \cite{CS,CS2,CS3} is that the constants in all the estimates that they establish do \emph{not} blow-up as $s\uparrow1$.
In particular, their proof yields as a limiting case the classical estimates of Krylov--Safonov and Evans--Krylov.

This is especially relevant in case of the weak Harnack inequality (Theorem \ref{half-Harnack-sup}), where the short and simple proof that we presented here is purely nonlocal and does not work for $s=1$, while the proof in \cite{CS} is much more delicate, as it must include somehow the proof of the Krylov--Safonov theorem.

If we substitute our Theorem \ref{half-Harnack-sup} by the weak Harnack inequality in \cite{CS}, then the rest of the proofs in this chapter can be easily adapted to yield constants that are uniform as $s\uparrow1$.

\subsection{Regularity estimates in $L^p$ spaces}

For second-order fully nonlinear equations of the form
\begin{equation}\label{question3.2}
\inf_{a\in \A}\big\{-\L_a u+c_a\big\} = f\quad\textrm{in}\quad B_1,
\end{equation}
with $s=1$, a celebrated result of Caffarelli \cite{CafW2p} establishes that if $f\in L^p(B_1)$  with $p\geq n$, then $u\in W^{2,p}_{\rm loc}(B_1)$.

For nonlocal equations, however, almost nothing is known in this direction:

\vspace{2mm}

\noindent\textbf{Open question 3.2}:
\textit{Assume that $f\in L^p$ in \eqref{question3.2}, with $\L_a\in \mathfrak L_s(\lambda,\Lambda)$. 
Can one show that $u\in W^{2s,p}$, when $p$ is large enough?}

\vspace{2mm}

With the current techniques, this problem seems out of reach, and hence completely new ideas are probably needed.

In case of second-order uniformly elliptic equations, an essential tool in order to establish $W^{2,p}$ estimates for \eqref{question3.2} is the so-called ABP estimate:
\[\left\{\begin{array}{rcl}
\sum_{i,j}a_{ij}(x)\partial_{ij}u &\geq& f \quad\textrm{in}\quad \Omega \\
u&\leq& 0 \quad \textrm{on}\quad \Omega
\end{array}\right.\qquad \Longrightarrow\qquad \sup_\Omega u \leq C\|f\|_{L^n}.
\]

Unfortunately, a nonlocal version of the ABP estimate is only known for $f\in L^\infty$ \cite{CS} or under strong structural hypothesis on the operators \cite{GS}.
Even the following basic question remains completely open.

\vspace{2mm}

\noindent\textbf{Open question 3.3}:
\textit{Assume $\mathcal M^+ u \geq f$ in $B_1$ and $u\leq 0$ in $\R^n\setminus B_1$, with $|f|\leq 1$.
Can one show that $\sup_{B_1} u \leq C\|f\|_{L^p}^\theta$, for some $\theta>0$ and $p<\infty$?}

\vspace{2mm}

Recall that $\mathcal M^+$ was defined in \eqref{eq:MMpm}.
Notice that in case of second-order PDE (i.e., $s=1$) this holds for $p=n$ and $\theta=1$ --- and the assumption $\|f\|_{L^\infty}\leq 1$ is not necessary.
For $s\in(0,1)$, the main result in  \cite{GS} establishes that the corresponding inequality holds for $p=n$ and $\theta=s$, for operators with kernels of the particular form $K(x,y)=y^TA(x)y/|y|^{n+2s+2}$.

\subsection{Boundary regularity}

Throughout this chapter we mainly studied the \emph{interior} regularity of solutions for fully nonlinear operators of the form
\[
\I u :=\inf_{a\in \A}\sup_{b\in\B}\big\{-\L_{ab} u+c_{ab}\big\}.
\]
It is then natural to wonder what can be said about the \emph{boundary} regularity of solutions to Dirichlet problems of the form
\[\left\{
\begin{array}{rcll}
\I u &=& f &\quad\textrm{in}\quad \Omega \\
u&=&0&\quad\textrm{in}\quad \R^n\setminus\Omega.
\end{array}\right.\]

The boundary regularity theory for fully nonlinear equations was developed by the second author and Serra in \cite{RS-Duke}, where we proved that, if the kernels of the operators $\L_{ab}$ are \emph{homogeneous}, i.e., 
\[K_{ab}(y)=\frac{K_{ab}(y/|y|)}{|y|^{n+2s}}\]
(as in in Sections \ref{sec:bdryregularity} and \ref{sec:higherorder}), then the following holds:
\[
\left\{
\begin{array}{l}
\partial\Omega\in C^{2,\alpha}\\
f\in C^{1+\alpha-s}(\overline{\Omega})\\
\text{$K_{ab}|_{\mathbb S^{n-1}}$ are ``regular enough''}
\end{array}
\right.
\quad\Longrightarrow\qquad 
\frac{u}{d^s}\in C^{1+\alpha}(\overline{\Omega}),
\]
for some small $\alpha>0$.

The exponent $\alpha$ comes from a related estimate for equations with bounded measurable coefficients with homogeneous kernels; see \cite[Proposition~1.1]{RS-Duke}.

Furthermore, it turns out that the homogeneity assumption is necessary, and for non-homogeneous kernels solutions are not even comparable to $d^s$ near the boundary.
For example, even for the extremal operators $\mathcal M^\pm$ in dimension $n=1$, there exist two exponents $0<\beta_1<s<\beta_2$ such that the functions $(x_+)^{\beta_i}$ satisfy
\[\mathcal M^+ (x_+)^{\beta_1}=0\qquad \textrm{and}\qquad \mathcal M^- (x_+)^{\beta_2}=0\qquad \textrm{in} \quad \R_+.\]
We refer to \cite[Section 2]{RS-Duke} for more details.

\subsection{More general kernels}

In Chapter \ref{ch:2} we saw that, in case of linear translation invariant equations, the interior regularity theory can be developed for \emph{any} $\L\in\GL$, i.e., with nonnegative symmetric kernels satisfying
\[
r^{2s} \int_{B_{2r}\setminus B_r} K(dy) \le \Lambda\qquad\text{for all}\quad r  >0,
\]
\[
0 < \lambda \le r^{2s-2}\inf_{e\in \S^{n-1}}\int_{B_{\Lambda r}\setminus B_r}|e\cdot y|^{2} K(dy)\qquad\text{for all}\quad r  >0.
\]
Notice that here $K$ could even be a measure, not necessarily absolutely continuous, as it happens for example when $\L=(-\partial_{x_1x_1}^2)^s+\cdots+(-\partial_{x_nx_n}^2)^s$.

One can then consider equations with bounded measurable coefficients of the type
\begin{equation}
\label{open-q3.50}
\L_x u = 0\quad \textrm{in}\quad B_1,
\end{equation}
for linear $x$-dependent operators with kernels satisfying
\begin{equation}
\label{open-q3.5}
r^{2s} \int_{B_{2r}\setminus B_r} K(x,dy) \le \Lambda\qquad\text{for all}\quad r  >0,
\end{equation}
\begin{equation}
\label{open-q3.5b}
0 < \lambda \le r^{2s-2}\inf_{e\in \S^{n-1}}\int_{B_{\Lambda r}\setminus B_r}|e\cdot y|^{2} K(x,dy)\qquad\text{for all}\quad r  >0.
\end{equation}
No regularity in $x$ is assumed.

An outstanding problem in this context is then the following:

\vspace{2mm}

\noindent\textbf{Open question 3.4}:
\textit{Can one prove a H\"older estimate  $\|u\|_{C^\gamma(B_{1/2})}\leq C\|u\|_{L^\infty(\R^n)}$ for solutions to general equations \eqref{open-q3.50}-\eqref{open-q3.5}-\eqref{open-q3.5b}?}

\vspace{2mm}

The best known result in this direction is due to Schwab and Silvestre \cite{SS}, who established such a H\"older estimate under the additional assumption
\[
\big|\{y\in B_{2r}\setminus B_r: K(x,y)\geq \lambda |y|^{-n-2s} \}\big|\geq \mu|B_{2r}\setminus B_r|
\]
for any $r\in(0,1)$, for some $\mu>0$.
This assumption is much weaker than 
\[\frac{\lambda}{|y|^{n+2s}}\leq K(x,y)\leq \frac{\Lambda}{|y|^{n+2s}},\]
but still leaves the above question completely open.

On the other hand, we refer to \cite{FR5} for a H\"older estimate in case of operators $\L$ with general kernels \eqref{open-q3.5}-\eqref{open-q3.5b}, under the extra assumption that they have ``small oscillation'' in the $x$-variable.

It is interesting to notice the similarity with respect to the case of divergence-form operators, in which an analogous question is completely open as well; see subsection \ref{sect-div-open}.

%% file: chap4.tex
%
%
%

\chapter{Obstacle problems}
\label{ch:obst_pb}

\index{Obstacle problem}
 
In this chapter we study obstacle problems of the type
\begin{equation}\label{obstacle-motivation1}
\min\big\{\L v,\,v-\varphi\big\} = 0 \quad\textrm{in}\quad \Omega\subset\R^n
\end{equation}
for integro-differential operators $\L$ of order $2s$.

When $\L=-\Delta$ (corresponding to $s=1$) this classical free boundary problem is quite well understood, starting with the groundbreaking work of Caffarelli in the 1970s \cite{Caf}; see \cite{PSU, FR4}.

For nonlocal operators, though, the regularity of solutions and the structure of free boundaries turns out to be much more complicated (even in case $\L=\sqrt{-\Delta}$), and there are still several open problems in this context.
The regularity theory for solutions and free boundaries was first developed for $\sqrt{-\Delta}$ and $(-\Delta)^s$, and more recently for more general integro-differential operators $\L$ of order $2s$; see  \cite{ACS,GP,FS},  \cite{S-obst,CSS}, and \cite{CRS,AR20,FRS}.

\section{Motivation}
\label{sec:motivations}

Obstacle problems for integro-differential operators appear in quite different settings.
We briefly describe some of them here, and refer to \cite{CT,CF,CDM,Serfaty,DL,F-survey} for more details.

\subsection{Optimal stopping}
 \index{Optimal stopping problem}
The first (and most classical) motivation to study this kind of obstacle problems comes from probability, in the so-called optimal stopping problem.

In this context, one considers the following stochastic control model.
We have a L\'evy process $X_t$ in $\R^n$ and some given payoff function $\varphi:\R^n\rightarrow\R$.
One wants to find the \emph{optimal stopping} strategy to maximize the expected value of $\varphi$ at the end point.
If $\L$ is the infinitesimal generator of the process~$X_t$ (c.f. subsections \ref{ssec:inf_gen} and \ref{subsect-2.1.4}), then it turns out that the value function~$v(x)$ (i.e., the maximum expected payoff we can obtain starting at~$x$) solves the following problem
\begin{equation} \label{obstacle2-motivation}
\left\{
\begin{array}{rcll}
v&\geq&\varphi & \quad \textrm{in}\quad \R^n,\\
\L v&\geq&0 & \quad \textrm{in}\quad \R^n,\\
\L v&=&0 & \quad \textrm{in}\quad \{v>\varphi\}.
\end{array}
\right.
\end{equation}
This means that the value function in any optimal stopping problem solves an obstacle problem of the type \eqref{obstacle-motivation1}.

In the context of mathematical finance, this type of problem appears as a pricing model for  American options \cite{CT}, where the function $\varphi$ represents the option's payoff, and the contact set $\{v=\varphi\}$ is the exercise region.
Notice that, in this context, the most important unknown to understand is precisely the contact set, i.e., one wants to find and/or understand the two regions $\{v=\varphi\}$ (in which we should exercise the option) and  $\{v>\varphi\}$ (in which we should wait and not exercise the option yet).
The free boundary is the separating interface between these two regions.

We refer to \cite[Chapter 6]{Eva13} for a description of the model in the case of Brownian motion, and to the book \cite{CT} for an exhaustive discussion in the case of jump processes; see also \cite{Merton,CF}.

\subsection{Interacting particle systems}

\index{Interacting particle systems}
A completely different motivation to study obstacle problems for integro-differential operators comes from the study of interacting particle systems.

Indeed, many different phenomena in physics, biology, or material science, give rise to models with interacting particles or individuals.
In such a context, the 2D mathematical model   is usually the following; see \cite{BCLR,Serfaty}.
We are given a repulsive-attractive interaction potential $W\in L^1_{\rm loc}(\R^2)$, and its associated interaction energy
\[E[\mu]:=\frac12\int_{\R^2}\int_{\R^2}W(x-y)d\mu(x)d\mu(y),\]
where $\mu$ is any probability measure in $\R^2$.

The potential $W$ is repulsive when the particles or individuals are very close, and attractive when they are far from each other.
A typical assumption is that near the origin we have
\begin{equation}\label{W}
W(z) \asymp \frac{1}{|z|^{\beta}}\quad \textrm{for}\quad z\sim 0,
\end{equation}
for some $\beta\in(0,2)$.
Moreover, it would usually grow at infinity, say $W(z)\asymp |z|^{\gamma}$ for $z\gg 1$.

An important question to be understood is that of the \emph{regularity of minimizers}, i.e., the regularity properties of the measures $\mu_0$ which minimize the interaction energy $E[\mu]$.
It turns out that the minimizer $\mu_0$ is given by $\mu_0=\L v$, where $v(x):=\int_{\R^2} W(x-y)d\mu_0(y)$ satisfies the {obstacle problem} \eqref{obstacle-motivation1}
for a certain operator $\L$ and a certain obstacle $\varphi$ that depend on $W$.
When $W$ satisfies \eqref{W} (as well as some extra conditions), such operator~$\L$ turns out to be an \emph{integro-differential} operator with a kernel
\[K(z)\asymp \frac{1}{|z|^{n+2s}}\quad \textrm{for}\quad z\sim 0,\]
with $n=2$ and $2s=2-\beta\in (0,2)$.

Therefore, understanding the regularity of minimizers $\mu_0$ of the interaction energy $E[\mu]$ is equivalent to understanding the regularity of solutions and free boundaries in obstacle problems for integro-differential operators.
In this setting, the contact set $\{v=\varphi\}$ is the support of the minimizer $\mu_0$.
We refer to \cite{CDM} and \cite{Serfaty} for more details on this topic.

\subsection{The thin obstacle problem}

\index{Thin obstacle problem}

In the particular case $\L=\sqrt{-\Delta}$, thanks to the extension property from Section \ref{sec:harm_extension}, the obstacle problem \eqref{obstacle2-motivation} for this operator is equivalent to the \emph{Signorini problem}, also known as the \emph{thin obstacle problem}.

This is a classical free boundary problem which dates back to 1933 and arises in a model of linear elasticity \cite{Sig59}. 
The problem consists in finding the elastic equilibrium configuration of a 3D elastic body, resting on a rigid frictionless surface and subject only to its mass forces.
The problem leads to a system of variational inequalities for the displacement vector in $\R^3$, which can be transformed to a scalar function $v$ that solves \eqref{obstacle-motivation1} for $\L=\sqrt{-\Delta}$; see \cite{CDV22, RS-thin,F-survey} for more details.

On the other hand, the Signorini problem gained further attention in the seventies due to its connection to mechanics and biology, where it models the process of \emph{osmosis} in the study of {semipermeable membranes}.
We refer to the classical book of Duvaut and Lions \cite{DL} for more details about these models.

\section{Basic properties of solutions}

Our goal is to study the regularity properties of solutions and free boundaries in the following class of obstacle problems: given $s\in(0,1)$ and $n>2s$ we consider solutions of
\begin{equation}\label{obst-pb}
\min\big\{\L v,\,v-\varphi\big\} = 0 \quad\textrm{in}\quad \R^n,
\end{equation} 
\begin{equation}\label{obst-pb2}
\lim_{|x|\to\infty} v(x)=0,
\end{equation}
with $\varphi\in C^\infty_c(\R^n)$, and where $\L$ is a L\'evy operator of the form
\begin{equation}\label{obst-op1}
\begin{split}
\L u(x) & = {\rm P.V.}\int_{\R^n}\big(u(x)-u(x+y)\big)K(y)dy \\
& = \frac12\int_{\R^n}\big(2u(x)-u(x+y)-u(x-y)\big)K(y)dy ,
\end{split}
\end{equation}
with 
\begin{equation}\label{obst-op2}
K(y)=K(-y) \qquad\textrm{and}\qquad \frac{\lambda}{|y|^{n+2s}} \leq K(y) \leq \frac{\Lambda}{|y|^{n+2s}},
\end{equation}
and $0<\lambda\leq \Lambda$. As in Chapter~\ref{ch:fully_nonlinear}, we denote by $\LLL$ the set of operators for which \eqref{obst-op1}-\eqref{obst-op2} hold (see Definition~\ref{defi:LL}). 
We will moreover assume that $\L$ is a stable operator, that is, 
\begin{equation}\label{obst-op3}
\textrm{$K$ is homogeneous, i.e.,} \quad
K(y) = \frac{K(y/|y|)}{|y|^{n+2s}}.
\end{equation}

\begin{defi} \label{defi:LLh} Let $s\in (0, 1)$ and  $0< \lambda\le \Lambda$. We define
\[
\LLh :=\big\{ \L \in \LL : \text{its kernel $K$ satisfies \eqref{obst-op3}}\big\}. 
\]
\end{defi}

This is the class of stable operators whose kernels are comparable to the one of the fractional Laplacian.

Notice that \eqref{obst-pb} (or, equivalently, \eqref{obstacle2-motivation}) is a \emph{free boundary problem}, in the sense that there are two unknowns: the solution $v$, and the set $\{v>\varphi\}$ where the equation $\L v=0$ holds.
The boundary of such (a priori unknown) set is the so-called free boundary $\partial\{v>\varphi\}$.

The main questions in this context are the following:
\begin{itemize}
\item What is the optimal regularity of solutions $v$ of \eqref{obst-pb}?

\item What can we say about the structure and regularity of their free boundaries?
\end{itemize}
These are the questions that we tackle in this chapter.

\subsection{Existence and uniqueness}

\index{Obstacle problem!Existence and uniqueness}
Solutions to the obstacle problem \eqref{obst-pb}-\eqref{obst-pb2} can be constructed in (at least) two different ways.
On the one hand, one can minimize the energy functional
\[\mathcal E(v) = \int_{\R^n}\int_{\R^n} \big|v(x)-v(y)\big|^2 K(x-y)\,dy\,dx\]
among all functions $v\in H^s(\R^n)$ that satisfy $v\geq \varphi$ in $\R^n$.
The minimizer $v\in H^s(\R^n)$ is then unique, it solves \eqref{obst-pb} in the weak sense, and  it decays to $0$ as $|x|\to\infty$ provided that $n>2s$. (See \cite{FR4, PSU} in the local case, and \cite{S-obst} for the case of the fractional Laplacian.)

On the other hand, one can also consider the infimum of all supersolutions $w \in {\rm LSC}(\R^n)$ that satisfy $w\geq\varphi$ (and $w\geq0$) in $\R^n$.
As in subsection~\ref{subsection3.2.5}, it turns out that such infimum $v$ is actually a (continuous) supersolution itself, and by minimality it   satisfies \eqref{obst-pb}.

The two ways to construct solutions are anyway a posteriori equivalent (see Theorem~\ref{thm:existence-obst} and Remark~\ref{rem:visstrong} below), and here for convenience we will proceed with the latter construction. 
This is precisely what we do in the following.

\begin{thm}\label{thm:existence-obst}
Let $s\in(0,1)$, $n>2s$, and let $\L \in \LLL$. Given $\varphi\in C_c(\R^n)$ there exists a unique viscosity solution $v\in C(\R^n)$ of
\begin{equation}
\label{eq:vsatisfies}
\left\{\begin{array}{rcll}
\L v & \geq & 0 & \ \textrm{in}\quad \R^n, \\
\L v & =&  0 & \ \textrm{in}\quad \{v>\varphi\}, \\
v & \geq & \varphi & \ \textrm{in}\quad \R^n,\\
 \lim_{|x|\to\infty}v(x)& = & 0.&
\end{array}\right.
\end{equation}
 Moreover, $v$ has the same modulus of continuity as $\varphi$, and if $\varphi\in {\rm Lip}(\R^n)$ then $v\in C^{2s+\eps}_{\rm loc}(\{v > \varphi\})\cap {\rm Lip}(\R^n)$ for some $\eps > 0$.
 
Finally, if $w\in {\rm LSC}(\R^n)\cap L^1_{\omega_s}(\R^n)$ satisfies 
\begin{equation}
\label{eq:wsatisfies}
\left\{
\begin{array}{rcll}
\L w & \geq & 0 & \ \textrm{in}\quad \R^n, \\
w & \geq & \varphi & \ \textrm{in}\quad \R^n ,\\
\liminf_{|x|\to\infty}w(x)&\geq & 0,&
\end{array}
\right.
\end{equation}
in the viscosity sense, then $w\geq v$ in $\R^n$. 
\end{thm}

\begin{proof}
We proceed by constructing the solution, in analogy with the construction in the existence of solutions to fully nonlinear integro-differential equations in Theorem~\ref{thm:existence_visc}. We divide the proof into five steps. 
\begin{steps}
\item 
Let us define $v$ as the infimum of all viscosity supersolutions that are above the obstacle and have nonnegative limits at infinity: 
\[
v(x) := \inf_{w\in \mathcal{S}} w(x),
\]
where 
 \[
 \mathcal{S} := \left\{w\in {\rm LSC}(\R^n)\cap L^1_{\omega_s}(\R^n): \begin{array}{l} \L w \ge 0 \ \text{ in $\R^n$ in the viscosity sense}, \\
  w\ge \varphi\ \text{in}\ \R^n, \ \liminf_{|x|\to\infty}w(x) \geq   0 \end{array}\right\}.
 \]
 
 The constant function equal to $\|\varphi\|_{L^\infty(\R^n)}$ belongs to $\mathcal{S}$, and so $v\le \|\varphi\|_{L^\infty(\R^n)}$ in $\R^n$.  On the other hand,   any function $w\in \mathcal{S}$ satisfies either $w> 0$ or $w\equiv 0$ in $\R^n$. This is a consequence of the fact that if  $w$ achieves its global minimum, it is necessarily constant (since $\L w\ge 0$ and $\L$ has a strictly positive kernel). 
 
 Indeed, observe first that $w\ge 0$ in $\R^n\setminus {\rm supp}(\varphi)$ by construction. Since~$w$ achieves its minimum in any compact set, if its minimum in ${\rm supp}(\varphi)$ was negative, it would be a global minimum and   $w$ would be constant and negative, a contradiction. Hence, $w\ge 0$ in $\R^n$, and using again that it cannot achieve its global minimum, we have that either $w\equiv 0$ in $\R^n$ or $w > 0$ in $\R^n$. If $\varphi$ is positive somewhere, we must have $w > 0$ in $\R^n$. 
 
 In particular, $0\le v \le \|\varphi\|_{L^\infty(\R^n)}$ in $\R^n$. Let us define
 \[
 v_*(x) := \inf\left\{\liminf_{k\to \infty} v(x_K) : \R^n\ni x_k \to x\right\}, 
 \]
 so that, by Lemma~\ref{lem:inf_is_super}, we have $\L v_* \le 0$ in the viscosity sense and $v_* \ge \varphi $ in $\R^n$ as well, therefore showing that $v_*\in \mathcal{S}$, and $ v = v_*\in {\rm LSC}(\R^n)$. As before, we immediately have that either $v\equiv 0$ in $\R^n$, or if $\varphi$ is positive somewhere, $v > 0$ in $\R^n$. 
 
 Observe  also  that, by construction, if $w$ satisfies \eqref{eq:wsatisfies}, then $w \ge v$. Let us prove the remaining properties of $v$. 
 
 \item We show now that $\lim_{|x|\to \infty} v(x) = 0$ and that $v$ is continuous. 
 
 To show the first part, we need a barrier from above. In this case, we can take the fundamental solution for the operator $\L$, which is a positive function $0<\Gamma\in {\rm LSC}(\R^n)\cap L^1_{\omega_s}(\R^n)$ such that $\L \Gamma \ge 0$ in $\R^n$ in the viscosity sense and 
 \[
 \Gamma(x) \le \frac{C}{|x|^{n-2s}}\quad\text{for}\quad |x|\gg 1
 \]
 (see Remark~\ref{rem:fund_sol}). 
 Thus, $M \Gamma \in \mathcal{S}$ (for $M$ large enough such that $M\Gamma \ge \varphi$) and $v(x) \le C|x|^{2s-n}$ for $|x|\gg 1$, in particular giving that  $\lim_{|x|\to \infty} v(x) = 0$.
 
 On the other hand, let $\sigma$ denote the modulus of continuity of $\varphi$:
 \[
 |\varphi(x) - \varphi(y)|\le \sigma(|x-y|)\quad \textrm{for all}\quad x,y\in\R^n. 
 \]
 We define, for any $h \in \R^n$,
 \[
 v_h(x) := v(x+h) + \sigma(|h|). 
 \]
 Then $v_h\in {\rm LSC}(\R^n)\cap L^1_{\omega_s}(\R^n)$, $\L v_h \ge 0$, $\liminf_{|x|\to \infty} v_h(x) \ge 0$, and 
 \[
 v_h(x) \ge \varphi(x+h) + \sigma(|h|) \ge \varphi(x).
 \]
That is, $v_h \in \mathcal{S}$ and $v\le v_h$ for any $h\in \R^n$ implies\footnote{We use $v(x) \le v_h(x)$ and $v(x+h) \le v_{-h}(x+h) = v(x) +\sigma(|h|)$.}
\[
-\sigma(|h|)\le v(x+h) - v(x) \le \sigma(|h|),
\]
so $v$ is continuous with modulus of continuity $\sigma$.

\item In order to see that $\L v = 0$ in $\{v > \varphi\}$, we proceed as in \ref{step:prevh00} of the proof of Proposition~\ref{prop:existence_main_bdd}. Indeed, assume that it is not true, and that there is a point $x_\circ\in \{v > \varphi\}$ and a test function $\eta$ that is $C^2$ in $B_r(x_\circ)$ for $r > 0$ and $\eta > v$ in $\R^n\setminus\{x_\circ\}$ with $\eta(x_\circ) = v(x_\circ)$ but $\L \eta(x_\circ) > 0$. By continuity of $\L \eta$ around $x_\circ$ (see Lemma~\ref{lem:Lu_LL}) we have that $\L \eta (x) <0$ in $B_\rho(x_\circ)$ for some $\rho > 0$ small.
 
We consider $\eta_\delta := \eta -\delta$ for some $\delta > 0$ small enough such that $\eta-\delta \ge \varphi$ in $B_r(x_\circ)$ (which exists by continuity of $v$ and $\varphi$). Since $\eta > v $ in $B_r(x_\circ)\setminus\{x_\circ\}$, we have that for $\delta> 0 $ small enough, $\eta_\delta > v$ in $B_r(x_\circ)\setminus B_\rho(x_\circ)$ as well. Let us define 
 \[
 v_\delta = \left\{
 \begin{array}{ll} 
 \min \{v, \eta_\delta\} & \quad\text{in}\quad B_\rho(x_\circ),\\
 v  & \quad\text{in}\quad B^c_\rho(x_\circ).
 \end{array}
 \right.
 \]
Then $v_\delta$ is a supersolution, since it coincides with $v$ in $B_\rho^c(x_\circ)$, and is the infimum of two supersolutions in $B_\rho(x_\circ)$ (recall Remark~\ref{rem:max_min}). Moreover, by construction we also have $v_\delta \ge \varphi$. This means that $v_\delta\in \mathcal{S}$, and therefore $v_\delta \ge v$. In particular,   $\eta-\delta \ge v$ in $B_r(x_\circ)$, and thus $\eta(x_\circ) - \delta \ge v(x_\circ)$, a contradiction. That is, $\L v = 0$ in $\{v > \varphi\}$, and $v$ satisfies \eqref{eq:vsatisfies}.

\item Let us now see that if $\varphi\in {\rm Lip}(\R^n)$ then $v\in C^{2s+\eps}(\{v > \varphi\})$ for some $\eps > 0$. Since $v$ has the same modulus of continuity as $\varphi$, if we define for any $h\in B_1$, $D_h v:= \frac{v(x+h) - v(x)}{|h|}$, we have $|D_h v|\le C$ for some $C$ independent of $h$, and by linearity $\L D_h v = 0$ in the viscosity sense in $\Omega_h := \left\{x : \dist(x, \{v = \varphi\}) > |h|\right\}$. We can now separate two cases: 

If $s \le \frac12$, we can use the H\"older estimates in Theorem~\ref{C^alpha-bmc} to deduce that $D_h v \in C^\gamma$ in $\Omega_h$, and hence by Lemma~\ref{it:H8} we obtain $v\in C^{1+\gamma}$ in $\{v>\varphi\}$. 

If $s  > \frac12$ instead, we can use Theorem~\ref{C1alpha} with $\alpha  = \min\left\{\frac{\gamma}{2}, 2s-1\right\}$, so that $\theta = 0$ there and $D_h v\in C^{1,\alpha}$, giving $v\in C^{2,\alpha}$ in $\{v > \varphi\}$ by Lemma~\ref{it:H8} again. In all cases, we have $v \in C^{2s+\eps}$ in $\{v >\varphi\}$, and therefore $v$ is a strong solution there.

\item It only remains to show uniqueness. 
Let   $w\in C(\R^n)$ be any function that  also satisfies \eqref{eq:vsatisfies}, and let us see that $v = w$. 

By construction, $v\le w$. Let $\bar u := w-v \ge 0$. Then $\bar u\in C(\R^n)\cap L^\infty(\R^n)$ and, whenever $\bar u > 0$, we know that $\L \bar u = \L w - \L v \le 0$ in the viscosity sense  (since $w > v\ge\varphi$ and $\L w = 0$ in $\{w > \varphi\}$). Thus, $\bar u$ satisfies
\[
\left\{
\begin{array}{rcll}
\L \bar u & \leq & 0 & \ \textrm{in}\quad \{\bar u  > 0\}, \\
\bar u & \geq & 0 & \ \textrm{in}\quad \R^n ,\\
\lim_{|x|\to\infty}\bar u(x)&= & 0.&
\end{array}
\right.\]
If $\bar u$ is not identically zero in $\R^n$, it has a global positive maximum at $x_\circ$, i.e., $u(x_\circ)>0$ and $u(x_\circ)\geq u$ in $\R^n$.
But then, by evaluating the operator $\L$ at this point we then get $\L u(x_\circ)>0$, a contradiction. 
Thus, $\bar u\equiv0$, and we are done. \qedhere
 \end{steps}
\end{proof}

\subsection{Semiconvexity}
\index{Obstacle problem!Semiconvexity}

A key tool in the study of obstacle problems is the following (semi)convexity property of solutions.

\begin{lem}[Semiconvexity]\label{semiconvex}
Let $s\in (0, 1)$, $\L \in \LLL$, $\varphi\in C^{1,1}_c(\R^n)$, and $v$ be the solution to \eqref{obst-pb}-\eqref{obst-pb2}.
Then,
\begin{enumerate}[leftmargin=*, label=(\roman*)]
\item \label{it:vsemiconvex} $v$ is semiconvex, with \[\partial_{ee}v\ge -\|\varphi\|_{C^{1,1}(\R^n)}\quad \textrm{for all}\quad  e\in \S^{n-1},\]
in the sense that $v(x)+C|x|^2$ is convex, for $C=\frac{1}{2n}\|\varphi\|_{C^{1,1}(\R^n)}$.

\item \label{it:vLipschitz} $v$ is Lipschitz, with \[\|v\|_{{\rm Lip}(\R^n)}\leq \|\varphi\|_{{\rm Lip}(\R^n)}.\]
\end{enumerate}
\end{lem}

\begin{proof}
For part \ref{it:vsemiconvex}, by Theorem~\ref{thm:existence-obst}, $v$ is the least supersolution which is above the obstacle $\varphi$ and is nonnegative at infinity.

Thus, for any given $h\in \R^n$ we may take
\[w(x)=\frac{v(x+h)+v(x-h)}{2}+C|h|^2.\]
This function clearly satisfies $\L w\geq0$ in $\R^n$, and also $w\geq\varphi$ in $\R^n$ for $C=\|\varphi\|_{C^{1,1}(\R^n)}$.
Hence, we have $w\geq v$ in $\R^n$, and therefore
\[\frac{v(x+h)+v(x-h)-2v(x)}{|h|^2}\geq -C.\]
Since $C$ is independent of $h\in B_1$, we get $\partial_{ee}v\geq -C$ for all $e\in \S^{n-1}$.

Part \ref{it:vLipschitz} follows directly from Theorem~\ref{thm:existence-obst}.
\end{proof}

\begin{rem}\label{rem:visstrong}
Thanks to Lemma~\ref{semiconvex}-\ref{it:vsemiconvex} we have that, if $\varphi\in C^{1,1}_c(\R^n)$, then $v$ is semiconvex and therefore we can in fact evaluate $\L v$   at any point in $\R^n$  (thanks to Remark~\ref{rem:onesidedcondition}), so that the expression in \eqref{eq:vsatisfies} is well-defined pointwise.
\end{rem}

\section{Boundary Harnack in Lipschitz (and more general) domains}

Boundary Harnack inequalities play a key role in many free boundary pro\-blems, and especially in obstacle-type problems.

The goal of this section is to prove the following boundary Harnack principle in Lipschitz (and more general) domains, a crucial tool in order to establish the regularity of solutions and free boundaries later on.

\begin{figure}
\centering
\makebox[\textwidth][c]{\includegraphics[scale = 1]{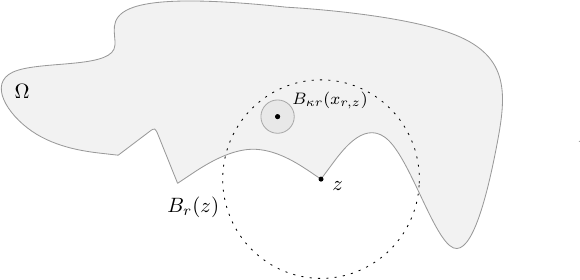}}
\caption{\label{fig:17} We consider domains satisfying assumption \eqref{Omega-Ch4}.}
\end{figure}

The only assumption that we will impose on the domain  $\Omega$ is actually the following (see Figure~\ref{fig:17}):
\begin{equation}\label{Omega-Ch4}
\textrm{for any}\quad z\in \partial\Omega,\quad \textrm{there is a ball}\quad B_{\kappa r}(x_{r,z}) \subset B_r(z)\cap \Omega
\end{equation}
for all $r\in(0,1)$, for some $\kappa>0$.  
This is often  referred as the \emph{interior corkscrew condition} \cite{JK82}.
Notice that it suffices to check \eqref{Omega-Ch4} for $r\in(0,r_\circ)$, with $r_\circ>0$, by taking $\kappa>0$ smaller if necessary.

\begin{thm}\label{thm-Lip}\index{Boundary Harnack inequality}
Let $s\in(0,1)$, let $\L\in \LLL$, and let $\Omega\subset\R^n$ be any open set satisfying \eqref{Omega-Ch4} for all $r\in(0,1)$, for some $\kappa>0$.
Then, there is $\delta>0$, depending only on $n$, $s$, $\kappa$, $\lambda$, and $\Lambda$, such that the following statement holds.

Let $u_1,u_2\in C(B_1)$ be viscosity solutions of
\begin{equation}\label{pb}
\left\{ \begin{array}{rcll}
|\L u_i| &\leq &\delta&\textrm{in }B_1\cap \Omega\\
u_i&=&0&\textrm{in }B_1\setminus\Omega
\end{array}\right.\end{equation}
satisfying
\begin{equation}\label{u-is-nonneg}
u_i\geq0\quad\mbox{in}\quad \R^n, \qquad \int_{\R^n}\frac{u_i(x)}{1+|x|^{n+2s}}\,dx=1.
\end{equation}
Then, there is $\alpha\in(0,1)$ such that
\[ \left\|\frac{u_1}{u_2}\right\|_{C^{0,\alpha}(\overline\Omega\cap B_{1/2})}\leq C.\]
The constants $\alpha$ and $C$ depend only on $n$, $s$, $\kappa$, $\lambda$, and $\Lambda$. 
\end{thm}

We will actually first prove the following boundary Harnack principle in arbitrary open sets.

\begin{thm}\label{thm-main}
Let $s\in(0,1)$, let $\L\in \LLL$, and let $\Omega\subset\R^n$ be any open set.
Let $x_\circ \in B_{1/2}$ and $\varrho>0$ be such that $B_{2\varrho}(x_\circ)\subset \Omega\cap B_{1/2}$.
Then, there exists $\delta>0$, depending only on $n$, $s$, $\varrho$, $\lambda$, and $\Lambda$, such that the following statement holds.

Let $u_1,u_2\in C(B_1)$ be viscosity solutions of \eqref{pb} satisfying \eqref{u-is-nonneg}.
Then,
\[ C^{-1}u_2\leq u_1\leq C\,u_2\quad\textrm{in}\quad B_{1/2}.\]
The constant $C$ depends only on $n$, $s$, $\varrho$, $\lambda$, and $\Lambda$.
\end{thm}

Once we   establish  Theorem \ref{thm-main}, the idea to prove Theorem~\ref{thm-Lip} is the same as in Theorem \ref{C^alpha-bmc}: A Harnack-type inequality implies a $C^\alpha$ estimate (for some small $\alpha>0$) thanks to a suitable iteration procedure.
In case of the interior Harnack, this iteration was simpler, and was carried out in Lemma \ref{oscillat-decay2}.
In the present setting of the boundary Harnack, the iteration will be a bit more delicate; see Proposition \ref{prop-iter} below.

The first boundary Harnack principle for nonlocal elliptic operators was established by Bogdan \cite{Bog97}, who proved it for the fractional Laplacian in Lipschitz domains.
Later, his result was extended to arbitrary open sets by Song and Wu \cite{SW99}; see also Bogdan-Kulczycki-Kwasnicki \cite{BKK08}, and more recently Bogdan-Kumagai-Kwasnicki \cite{BKK15} established the boundary Harnack principle in general open sets for a wide class of Markov processes with jumps.
The proof of Theorems \ref{thm-Lip} and \ref{thm-main} that we present here is based on \cite{RS-bdryH}, and can be easily extended to $x$-dependent operators with kernels $K(x,y)\asymp |y|^{-n-2s}$.

\subsection{Proof of the boundary Harnack principle in open sets}

To prove Theorem \ref{thm-main}, we need the following.

\begin{lem}\label{lem-use}
Let $s\in(0,1)$, let $\L\in \LLL$, and let $\Omega\subset\R^n$ be any open set.
Let $x_\circ\in B_{1/2}$ and $\varrho>0$ be such that $B_{2\varrho}(x_\circ)\subset \Omega\cap B_{1/2}$.

Let $C_\circ\geq0$, and $u\in C(B_1)$ be a viscosity solution of
\[\left\{ \begin{array}{rcll}
|\L u| &\leq&C_\circ&\textrm{in }B_1\cap \Omega\\
u&=&0&\textrm{in }B_1\setminus\Omega.
\end{array}\right.\]
Assume in addition that $u\geq0$ in $\R^n$.
Then,
\[\sup_{B_{3/4}}u\leq C\left(\inf_{B_\varrho(x_\circ)}u+C_\circ \right),\]
with $C$ depending only on $n$, $s$, $\varrho$, $\lambda$, and $\Lambda$. 
\end{lem}

\begin{proof}
Since $u\geq0$ in $B_1$ and $\L u\leq C_\circ$ in $B_1\cap \{u>0\}$,   it follows from the definition of viscosity solutions (see Remark~\ref{rem:max_min}) that $\L u\leq C_\circ$ in all of $B_1$ in the viscosity sense.
Thus, by Theorem \ref{half-Harnack-sub} (and a standard covering argument) we have
\[\sup_{B_{3/4}}u\leq C\left(\int_{\R^n}\frac{u(x)}{1+|x|^{n+2s}}\,dx+C_\circ\right).\]
Now, using Theorem~\ref{half-Harnack-sup} in the ball $B_{2\varrho}(x_\circ)$, we find
\[\int_{\R^n}\frac{u(x)}{1+|x|^{n+2s}}\,dx\leq C\left(\inf_{B_\varrho(x_\circ)} u+C_\circ\right).\]
Combining the previous estimates, the lemma follows.
\end{proof}

We next give the proof of the boundary Harnack principle in arbitrary open sets: 

\begin{proof}[Proof of Theorem \ref{thm-main}]
First, as in Lemma \ref{lem-use}, by \eqref{u-is-nonneg} we have
\begin{equation}\label{ineq0}
u_i\leq C\quad \textrm{in}\ B_{3/4}
\end{equation}
and
\begin{equation}\label{ineq0ghj}
u_i\geq c>0\quad\textrm{in}\ B_\varrho(x_\circ),
\end{equation}
provided that $\delta>0$ is small enough.
Notice that $c$ depends only on $n$, $s$, $\lambda$, $\Lambda$, and $\varrho$, but not on $\Omega$.

Let now $\xi\in C^\infty_c(B_{2/3})$ be such that $0\leq \xi\leq 1$ and $\xi\equiv1$ in $B_{1/2}$, and let $\eta\in C^\infty_c(B_\varrho(x_\circ))$ such that $0\leq\eta\leq 1$ in $B_\varrho(x_\circ)$ and $\eta=1$ in $B_{\varrho/2}(x_\circ)$.
We define
\[w:=u_1\chi_{B_{3/4}}+C_1(\xi-1)+C_2\eta.\]
Then, thanks to \eqref{ineq0}, if $C_1$ is chosen large enough, we will have
\[w\leq 0\quad\textrm{in}\quad \R^n\setminus B_{2/3}.\]
Moreover, taking now $C_2$ large enough,
\[\begin{split}
\L w& = \L u_1-\L (u_1\chi_{\R^n\setminus B_{3/4}})+C_1\L \xi+C_2\L \eta \\
&\leq \delta+C+CC_1-cC_2  \leq -1\qquad\qquad
\textrm{in}\quad \Omega \cap B_{2/3}\setminus B_\varrho(x_\circ).\end{split}\]
Here we used that $\L u_1\leq \delta$ in $\Omega\cap B_1$, that 
\[\L (u_1\chi_{\R^n\setminus B_{3/4}})\geq -C\int_{\R^n}\frac{u_1(x)}{1+|x|^{n+2s}}\, dx\geq -C \quad \textrm{in}\quad B_{2/3},\]
 that $\L \xi\leq C$, and that 
\[\L \eta\leq -c<0\quad  \textrm{in}\quad B_1\setminus B_\varrho(x_\circ).\]
As a consequence, for any $C_3\leq \delta^{-1}$ we get 
\[\L (w-C_3u_2)\leq -1+C_3\delta\leq 0\quad \textrm{in}\ \Omega\cap B_{2/3}\setminus B_\varrho(x_\circ).\]

\begin{figure}
\centering
\makebox[\textwidth][c]{\includegraphics[scale = 1]{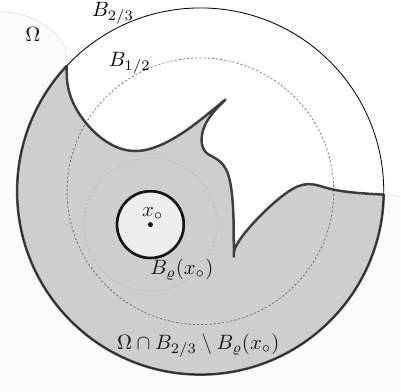}}
\caption{\label{fig:11} The domain where we apply the comparison principle in the proof of Theorem~\ref{thm-main}.}
\end{figure}

On the other hand, since $w\leq C$ in $B_\varrho(x_\circ)$ and $u_2\geq c>0$ in $B_\varrho(x_\circ)$, we have
\[w\leq C_3u_2\quad\textrm{in}\ B_\varrho(x_\circ)\]
for some constant $C_3\le \delta^{-1}$ provided that we have taken $\delta$ small enough. Notice, also, that $w \le 0$ in $B_{2/3}\setminus \Omega$ (since $u_1 = 0$ in $B_1\setminus \Omega$ and $\xi\le 1$) so that, by the comparison principle in $\Omega\cap B_{2/3}\setminus B_\varrho(x_\circ)$ (see Figure~\ref{fig:11} for a depiction of this domain) we find 
\[w\leq C_3u_2 \quad \textrm{in all of}\quad \R^n.\]

In particular, since $w\equiv u_1$ in $B_{1/2}\setminus B_\varrho(x_\circ)$, this yields
\[u_1\leq C_3u_2\quad \textrm{in}\ B_{1/2}\setminus B_\varrho(x_\circ).\]
Since $u_1$ and $u_2$ are comparable in $B_\varrho(x_\circ)$ (see \eqref{ineq0}-\eqref{ineq0ghj}), we deduce
\[u_1\leq C u_2\quad \textrm{in}\quad B_{1/2},\]
maybe with a larger constant $C$. 
Changing the roles of $u_1$ and $u_2$, we obtain the desired result. 
\end{proof}

\subsection{Iteration and proof of Theorem \ref{thm-Lip}}

We next proceed with the proof of Theorem \ref{thm-Lip}.

For this, we first need a lower bound for positive solutions $u$ in  domains satisfying \eqref{Omega-Ch4}, namely $u\geq cd^{2s-\gamma}$ in $\Omega$ for some small $\gamma>0$.

\begin{lem}\label{lem-Lip-dom}
Let $s\in(0,1)$, $\kappa>0$, and let $\L\in \LLL$. Let $\Omega\subset\R^n$ be any open set satisfying \eqref{Omega-Ch4} for all $r\in(0,2)$ and for some $\kappa > 0$.
Then, there exist $\delta>0$, $\gamma\in(0,1)$, and $c_\circ >0$ depending only on $n$, $s$, $\kappa$, $\lambda$, and $\Lambda$, such that the following statement holds.

Let $u$ be a viscosity solution of
\[\left\{ \begin{array}{rcll}
|\L u| &\leq&\delta &\textrm{in }\Omega \cap B_2\\
u&=&0&\textrm{in }B_2\setminus\Omega.
\end{array}\right.\]
Assume, in addition, that $u\geq0$ in $\R^n$, and $\|u\|_{L^\infty(B_1)}\geq 1$. Then, 
\[u\geq c_\circ d^{2s-\gamma}>0\quad \textrm{in}\quad \Omega\cap B_1.\]
\end{lem}

\begin{proof}
 Let $\dr_D$ be the regularized distance given by Lemma \ref{lem:distance}, with 
\[D:= (\Omega \cap B_1)\cup \bigcup_{\substack{ z\in \partial\Omega\cap B_1,\\ r\in(0,1)}} B_{\kappa r}(x_{r,z}).\]
Notice that $D\subset B_2$ still satisfies \eqref{Omega-Ch4}, and that $\Omega\cap B_1 \subset D$.
Then, by Lemma \ref{subsolution-2s-eps}, we have that
\[\L (\dr_D^{2s-\gamma}) \leq -1\quad \textrm{in}\quad D\cap \{\dr_D<\rho\} \]
for some small $\rho>0$.
Therefore, we have
\[\L u\geq -\delta \geq \L (\delta \dr_D^{2s-\gamma}) \quad \textrm{in}\quad D\cap \{\dr_D<\rho\}.\]

Now, since $\|u\|_{L^\infty(B_1)}\geq 1$, it follows from Lemma \ref{lem-use} that 
\[u\geq c_\rho>0 \quad \textrm{in}\quad D_\rho:=D\cap \{\dr_D\geq \rho\}.\]
This means in particular that $u\geq \delta \dr_D^{2s-\gamma}$ in $D_\rho$ if $\delta$ is small enough.
Since $u\geq0$ in $D^c$, it follows from the comparison principle that 
\[u\geq \delta \dr_D^{2s-\gamma} \quad \textrm{in}\quad D,\]
which yields the result.
\end{proof}

As a consequence, we find the following:

\begin{cor}\label{cor-Lip-dom-growth}
Let $s\in(0,1)$, $R>1$, $C_\circ >0$, and let $\L\in \LLL$. Let $\Omega\subset\R^n$ be any open set satisfying \eqref{Omega-Ch4} for all $r\in(0,R)$ and some $\kappa > 0$, and let $0\in \partial\Omega$. Then, there exist $\delta>0$, $\gamma_\circ\in(0,1)$, $R_\circ>1$, and $C_1$, depending only on $n$, $s$, $C_\circ$, $\kappa$, $\lambda$, and $\Lambda$, such that the following statement holds.

Let $R\geq R_\circ$, $\gamma\in(0,\gamma_\circ)$, and let $u$ be any viscosity solution of
\[\left\{ \begin{array}{rcll}
|\L u| &\leq& C_\circ R^{-\gamma} &\textrm{in}\quad\Omega\cap B_{2R}\\
u&=&0&\textrm{in}\quad B_{2R}\setminus\Omega.
\end{array}\right.\]
Assume, in addition, that $u\geq0$ in $\R^n$ and $u(x_{1,0}) \le 1$, where $x_{1,0}$ is given by  \eqref{Omega-Ch4}.
Then, 
\[\|u\|_{L^\infty(B_{\rho})} \leq C_1 \rho^{2s-\gamma} \quad \textrm{for all}\quad \rho\in[1,R].\]
\end{cor}

\begin{proof}
Assume we have $\rho\in [1, R]$ such that $\|u\|_{L^\infty(B_{\rho})} > C_1 \rho^{2s-\gamma}$, for $C_1$ to be chosen.
Then, we define
\[u_\rho(x):= \frac{u(\rho x)}{\|u\|_{L^\infty(B_\rho)}},\]
which satisfies
\[|\L_\rho u_\rho|\leq \frac{C_\circ}{C_1}\rho^{\gamma}R^{-\gamma}\le \frac{C_\circ}{C_1} \quad \textrm{in}\quad B_{2R/\rho}\cap (\rho^{-1}\Omega),\]
for some $\L_\rho\in \LLL$ (see \ref{eq:scaleinvariance_comp}).
Therefore, if $C_1$ is large enough, we have $C_\circ/C_1<\delta$ and we can use   Lemma \ref{lem-Lip-dom}  on the function $u_\rho$  to get that 
\[u(\rho x) \geq c_\circ \rho^{\gamma-2s}\|u\|_{L^\infty(B_\rho)}d_\Omega^{2s-\gamma}(\rho x), \]
and hence, taking $\rho x = x_{r\rho, 0}$ (we recall that $0\in \partial\Omega$, so $x_{r\rho, 0}$  given by \eqref{Omega-Ch4} is such that $B_{\kappa r\rho }(x_{r\rho, 0})\subset B_{r\rho} \cap \Omega$),
\[u(x_{r\rho,0}) \geq c_1 r^{2s-\gamma}\|u\|_{L^\infty(B_\rho)} \]
for any $r\in(0,1)$.
In particular, taking $r=1/\rho$, we find
\[ \|u\|_{L^\infty(B_\rho)} \leq \frac{1}{c_1} \rho^{2s-\gamma},\]
a contradiction if $C_1c_1\geq1$.
\end{proof}

We can now prove the main step towards the boundary Harnack inequality in domains satisfying \eqref{Omega-Ch4} (in particular, in Lipschitz domains).

\begin{prop}\label{prop-iter}
Let $s\in(0,1)$, $\kappa>0$,  and let $\L\in \LLL$.
Then, for any $\varepsilon>0$ there exist $R>1$ and $\delta>0$, depending only on $n$, $s$, $\kappa$, $\varepsilon$, $\lambda$, and~$\Lambda$, such that the following statement holds.

Let $\Omega\subset\R^n$ be any open set with $0\in\partial\Omega$ and satisfying \eqref{Omega-Ch4} for all $r\in(0,R)$ and $\kappa > 0$.
Let $u_1,u_2\in C(B_1)$ be viscosity solutions of
\[
\left\{ \begin{array}{rcll}
|\L u_i| &\leq &\delta&\textrm{in }B_{2R}\cap \Omega\\
u_i&=&0&\textrm{in }B_{2R}\setminus\Omega
\end{array}\right.
\]
satisfying $u_i\geq0$ in $B_{2R}$,  $u_i(x_{1,0}) = 1$, where $x_{1,0}$ is given by  \eqref{Omega-Ch4}, and
\[|u_i(x)|\leq 1+|x|^{2s-\varepsilon}\quad \textrm{in}\quad \R^n\setminus B_R.\]

Then, 
\[\osc_{R^{-k}}\frac{u_1}{u_2}\leq CR^{-k\alpha}\qquad\text{for all}\quad k \in \N_0,\]
with $C$ and $\alpha>0$ depending only on $n$, $s$, $\varepsilon$, $\kappa$, $\lambda$, and $\Lambda$.
\end{prop}

\begin{proof}
We split the proof into three steps.

\begin{steps}
\item \label{it:step1:propiter} We first prove the case $k=0$.
More precisely, we will prove that 
\begin{equation}
\label{eq:mMineq}
0<m \leq \frac{u_1}{u_2} \leq M\quad \textrm{in}\quad B_2
\end{equation}
for some positive constants $m<1$ and $M>1$.
For this, we want to apply Theorem \ref{thm-main} to the functions $u_i\chi_{B_{2R}}$.

Indeed, on the one hand notice that, exactly as in Lemma \ref{lem-use}, from $u_i(x_{1,0})=1$ we deduce that 
\[
C^{-1}\le \|u_i\|_{L^\infty(B_1)}\le C\qquad\text{and}\qquad C^{-1}\le \int_{\R^n}\frac{u_i(x)}{1+|x|^{n+2s}}\,dx\le C.
\]
On the other hand, thanks to the growth assumption on $u_i$, we will have that $|\L (u_i\chi_{\R^n\setminus B_{2R}})|<CR^{-\eps}$ in $B_R$, and therefore 
\[
|\L (u_i\chi_{B_{2R}})|<\delta+CR^{-\eps}\quad\text{in}\quad B_R\setminus\Omega.
\]
Taking $\delta>0$ small enough and $R>1$ large enough, it follows from Theorem~\ref{thm-main} that 
\[
C^{-1}u_2\leq u_1\leq Cu_2\quad\text{in}\quad B_2
\] for some constant $C$, as claimed.

Notice also that, thanks to Corollary \ref{cor-Lip-dom-growth}, we have 
\begin{equation}
\label{eq:growthui}
|u_i(x)|\leq C_1|x|^{2s-\eps}\quad\text{in}\quad B_R\setminus B_1,\end{equation}
 provided that $\delta \leq R^{-\eps}$.

\item \label{it:step2:propiter} We next show how to apply iteratively  \ref{it:step1:propiter}.
Let $\theta>0$ small (to be chosen later), and let us consider the functions
\[u_1^{(R)}(x) := \frac{(u_1-\theta u_2)(x/R)}{(u_1-\theta u_2)(z_R)},\qquad 
u_2^{(R)}(x) := \frac{(u_2-\theta u_1)(x/R)}{(u_2-\theta u_1)(z_R)},\]
where $z_R:=x_{1/R,0}$ is given by \eqref{Omega-Ch4}, 
and $\theta$ is small enough so that $u_i^{(R)}$ satisfy (recall \eqref{eq:mMineq})
\[
u_i^{(R)}=0\quad\text{in}\quad B_{2R^2}\setminus (R\Omega)\qquad\text{and}\qquad u_i^{(R)} \geq0\quad\text{in}\quad B_{2R}.\]
Moreover, by definition we also have $u_i^{(R)}(Rz_R)=1$.

Now, by Lemma \ref{lem-Lip-dom} and condition \eqref{Omega-Ch4}, we have 
\[u_i(x_{r,0}) \geq c_\circ r^{2s-\gamma} \quad \textrm{for}\quad r\in(0,1),\]
for some $\gamma>0$.
In particular, since $\|u_i\|_{L^\infty(B_{1/R})}\leq C$, if $\theta$ is small enough (namely, $\theta\ll R^{\gamma-2s}$) we will have 
\[(u_1-\theta u_2)(z_R) \geq c_1R^{\gamma-2s},\]
and the same holds for $(u_2-\theta u_1)(z_R)$.
This implies that, for $R$ large, 
\[\big|\L_R u_i^{(R)}\big| \leq \frac{(1+\theta)\delta R^{-2s}}{c_1R^{\gamma-2s}} \leq \delta\quad \textrm{in}\quad B_{2R^2}\cap (R\Omega),\]
for some $\L_R\in \LLL$ (see  \eqref{eq:scaleinvariance_comp}).
Moreover, by the growth of $u_i$ in~\eqref{eq:growthui},
\[
\big|u_1^{(R)}(x)\big| \leq 
\frac{C_1(1+\theta)|x/R|^{2s-\eps}}{c_1R^{\gamma-2s}} 
\leq |x|^{2s-\varepsilon} \quad\textrm{for} \quad R\leq |x|\leq R^2,
\] 
where we used that we can choose $\varepsilon>0$ small and $R$ large so that $R^{-\gamma}\ll R^{-\varepsilon}$.
On the other hand, by the growth assumption on $u_i$ in $\R^n\setminus B_R$, we also have 
\[
|u_1^{(R)}(x)| \leq 
\frac{C|x/R|^{2s-\varepsilon}}{c_1R^{\gamma-2s}} \leq |x|^{2s-\varepsilon}  \quad \textrm{for}\quad |x|\geq R^2,
\] 
and the same bounds hold for $u_2^{(R)}$.

Thus, the functions $u_i^{(R)}$ satisfy the same assumptions as $u_i$, and using the case $k=0$ from  \ref{it:step1:propiter} we deduce
\[0<m \leq \frac{u_1^{(R)}}{u_2^{(R)}} \leq M\quad \textrm{in}\quad B_1.\]
Rescaling back to $u_i$, this means that 
\[m \frac{(u_1-\theta u_2)(z_R)}{(u_2-\theta u_1)(z_R)}  \leq \frac{u_1-\theta u_2}{u_2-\theta u_1} \leq M \frac{(u_1-\theta u_2)(z_R)}{(u_2-\theta u_1)(z_R)}\quad \textrm{in}\quad B_{1/R}.\]
After some manipulations, the upper bound is equivalent to
\[\frac{u_1/u_1(z_R)}{u_2/u_2(z_R)} \leq
M - \theta (M-1)\frac{\frac{u_2(z_R)}{u_1(z_R)}+M\frac{u_1(z_R)}{u_2(z_R)} - (M+1)\theta}{1-\theta^2M+\theta(M-1)\frac{u_1(z_R)}{u_2(z_R)}},
\]
if $\theta$ is small enough such that the denominator in the right-hand side is positive, and an analogous expression holds for the lower bound.

Now, using that 
\[C^{-1}u_2(z_R)\leq u_1(z_R) \leq Cu_2(z_R),\]
and taking $\theta$ smaller if necessary, this yields
\[ m + c(1-m)\theta \leq \frac{u_1/u_1(z_R)}{u_2/u_2(z_R)} \leq M-c(M-1)\theta \quad \textrm{in}\quad B_{1/R},\]
for some small $c>0$.

\item  
Fix $\theta>0$ small enough so that the previous inequalities hold (and $c\theta<1$), and let $\alpha>0$ be such that
\[R^{-\alpha}:=1-c\theta<1.\]
(Notice that $\alpha$ is small, since we were requiring $\theta\ll R^{\gamma-2s}$.) Define
\[M_k:=1+(M-1)R^{-\alpha k},\qquad m_k:=1-(1-m)R^{-\alpha k}.\]
Then, we clearly have $M_0=M$, $m_0=m$, and 
\[M_{k+1}=M_k - c\theta(M_k-1),\qquad m_{k+1}=m_k + c\theta(1-m_k).\]
Thus, iterating \ref{it:step2:propiter} as many times as necessary, we get
\[ m_k \leq \frac{u_1/u_1(z_{k})}{u_2/u_2(z_{k})} \leq M_k \quad \textrm{in}\quad B_{R^{-k}},\]
where $z_{k} := x_{R^{-k-1}, 0}$, and hence, if we denote 
\[
S_k := \frac{u_1(z_{k})}{u_2(z_{k})},
\]
we have
\begin{equation}
\label{eqkihvasodihasd}
\big(1-C R^{-k\alpha}\big) S_k
\leq \frac{u_1}{u_2} \leq 
\big(1+CR^{-k\alpha}\big) S_k
\quad \textrm{in}\quad B_{R^{-k}}. 
\end{equation}
In particular, evaluating it at $z_{k+1}$, 
\[\left|S_{k+1}-S_k \right| \leq CR^{-k\alpha} S_k.\]
That is, since $S_0 \le C$, 
\[
S_k \le \prod_{i = 0}^k \left(1+CR^{-i\alpha}\right) \le e^{C \sum_{i = 0}^k R^{-\alpha k }} \le C,
\]
Hence, we actually get 
\[
|S_{k+1} - S_k|\le C R^{-k\alpha}
\]
and summing a geometric series again we obtain 
\[\left|S_k- c_\circ\right| \leq CR^{-k\alpha},\]
for some $c_\circ>0$. Plugged into \eqref{eqkihvasodihasd} this implies 
\[ \left|\frac{u_1}{u_2} - c_\circ\right| \leq CR^{-k\alpha} \quad \textrm{in}\quad B_{R^{-k}},\]
which is the desired result.\qedhere
\end{steps}
\end{proof}

We finally give the:

\begin{proof}[Proof of Theorem \ref{thm-Lip}]
We split the proof into two steps.

\begin{steps}
\item \label{it:step1:Lip} We first prove that
\begin{equation} \label{dfg-step1}
\left\|\frac{u_1}{u_2} -c_z\right\|_{L^\infty(B_r(z))} \leq C_1r^\alpha \quad \textrm{for}\quad r\in(0,r_\circ),\end{equation}
for any $z\in \partial \Omega\cap B_{1/2}$ and for some $r_\circ>0$.

Up to a translation and rotation, we may assume $z=0\in\partial\Omega$.
Thus, by Lemma \ref{lem-Lip-dom} we have that $u_i(x_{r,0})\geq cr^{2s-\gamma}$ for $r\in(0,\frac12)$ for some small $c,\gamma>0$, with $x_{r, 0}$ given by \eqref{Omega-Ch4}.
Moreover, as in Lemma \ref{lem-use}, we have that $\|u_i\|_{L^\infty(B_{3/4})}\leq C$.

Let $R>1$ be given by Proposition \ref{prop-iter}, and let $\varrho_\circ\geq 4R$.
Define
\[\tilde u_i(x) := \frac{(u_i\chi_{B_{3/4}})(x/\varrho_\circ)}{u_i(x_{1/\varrho_\circ,0})} .\]
Then, we clearly have $\tilde u_i(\varrho_\circ x_{1/\varrho_\circ,0})=1$ and since $|\L (u_i\chi_{B_{3/4}^c})|\leq C$ in $B_{1/2}$ and $u_i(x_{1/\varrho_\circ,0})\geq c\varrho_\circ^{\gamma-2s}$ (by Lemma~\ref{lem-Lip-dom}), we have, for some $\L_{\varrho_\circ} \in \LLL$ (see  \eqref{eq:scaleinvariance_comp}),
\[
|\L_{\varrho_\circ} \tilde u_i|\leq \delta+C\varrho_\circ^{-\gamma}\quad\text{in}\quad B_{\varrho_\circ/2}\cap \Omega.
\]
Therefore, taking $\varrho_\circ$ large enough and $\delta$ small enough (depending on $\varrho_\circ$), by Corollary \ref{cor-Lip-dom-growth}, we have
\[\|\tilde u_i\|_{L^\infty(B_\rho)} \leq C_1\rho^{2s-\gamma}\quad \textrm{for}\quad \rho\in(1,\varrho_\circ/4).\]
Moreover, since $\|\tilde u_i\|_{L^\infty(\R^n)} \leq C\varrho_\circ^{2s-\gamma}$, the previous bound holds for all $\rho\geq1$.

This means that these functions $\tilde u_i$ satisfy the assumptions of Proposition~\ref{prop-iter}, and hence
\[ \left\|\frac{\tilde u_1}{\tilde u_2} - c_0\right\|_{L^\infty(B_r)} \leq C_\circ r^\alpha \quad \textrm{for}\quad r\in(0,{\textstyle \frac12}),\]
for some $c_0>0$.
Rescaling back to $u_i$, we find \eqref{dfg-step1}.

\item  We now prove the H\"older estimate for $u_1/u_2$ by combining \eqref{dfg-step1} with interior estimates.

Namely, we want to prove that for any $x,y\in \Omega\cap B_{1/2}$ we have
\[\left|\left(\frac{u_1}{u_2}\right)(x) - \left(\frac{u_1}{u_2}\right)(y)\right| \leq C|x-y|^{\alpha'},\]
for some small $\alpha'>0$.
Fix $x_\circ\in \Omega\cap B_{1/2}$, and define 
\[d=d(x_\circ).\]

We first notice that the function $u_i(x_\circ+dx)$ satisfies the equation in~$B_1$ and is uniformly bounded (in $x_\circ$) in $L^1_{w_s}(\R^n)$, so by interior estimates (e.g., Theorem~\ref{C^alpha-bmc_2}) we have 
\[ [u_i]_{C^\alpha(B_{d/2}(x_\circ))} \leq Cd^{-\alpha}.\]
Since $u_2\geq c_\circ d^{2s-\gamma}$ in $B_{d/2}(x_\circ)$ (by Lemma \ref{lem-Lip-dom}) we deduce that 
\[
\frac{\left|u^{-1}_2(x_1)-u^{-1}_2(x_2)\right|}{|x_1-x_2|^\alpha} = \frac{\left|u_2(x_1)-u_2(x_2)\right|}{|u_2(x_1)u_2(x_2)||x_1-x_2|^\alpha} \le C \frac{\left|u_2(x_1)-u_2(x_2)\right|}{d^{4s-2\gamma}|x_1-x_2|^\alpha} 
\]
for any $x_1, x_2\in B_{d/2}(x_\circ)$, so that 
\[ \big[u_2^{-1}\big]_{C^\alpha(B_{d/2}(x_\circ))} \leq Cd^{2\gamma-\alpha-4s}.\]
Then, using \eqref{eq:APP_In1}, this yields.
\[
\left[\frac{u_1}{u_2}\right]_{C^{\alpha}(B_{d/2}(x_\circ))} \le C(d^{2\gamma-\alpha-4s}+d^{\gamma-\alpha -2s})
\]
and so (recall $\gamma < 2s$),
\[\left|\left(\frac{u_1}{u_2}\right)(x_\circ) - \left(\frac{u_1}{u_2}\right)(y)\right| \leq C|x_\circ-y|^\alpha d^{2\gamma-\alpha-4s}
\]
for any $y\in B_{d/2}(x_\circ)$.
In particular, choosing 
\[\alpha\leq 2 \gamma,\qquad  |x_\circ-y|\leq d^\theta,\qquad \theta\alpha>4s,\] 
we deduce that
\[\left|\left(\frac{u_1}{u_2}\right)(x_\circ) - \left(\frac{u_1}{u_2}\right)(y)\right| \leq C|x_\circ-y|^{\alpha-4s/\theta}\qquad\quad\textrm{if}\quad |x_\circ-y|\leq d^\theta.
\]

On the other hand, we want a similar bound for $|x_\circ-y|> d^\theta$.
For this, let $z\in \partial \Omega$ be such that $d(x_\circ)=|z-x_\circ|$.
Then, since $|y-z|\leq d+|x_\circ-y|$, by \eqref{dfg-step1} we find 
\[\begin{split}
\left|\left(\frac{u_1}{u_2}\right)(x_\circ) - \left(\frac{u_1}{u_2}\right)(y)\right| & \leq C\big(d+|x_\circ-y|\big)^\alpha \\
& \leq C|x_\circ-y|^{ \alpha/\theta}\qquad  \qquad
\textrm{if} \quad |x_\circ-y| \geq d^\theta.
\end{split}\]
Combining the previous inequalities, we deduce
\[\left|\left(\frac{u_1}{u_2}\right)(x_\circ) - \left(\frac{u_1}{u_2}\right)(y)\right| \leq C|x_\circ-y|^{\alpha'}, \quad \quad \alpha':=\min\{\alpha-4s/\theta,\,\alpha/\theta\}>0,\]
for all $x_\circ,y\in \Omega\cap B_{1/2}$, and the theorem is proved.\qedhere
\end{steps}
\end{proof}

\section{Regularity of free boundaries near regular points}

The goal of this section is to establish the main free boundary regularity results for the obstacle problem, \eqref{obst-pb}.

The first and main result is the following, which gives a dichotomy between regular and degenerate points, and establishes the $C^{1,\gamma}$ regularity of free boundaries near regular points. Notice that now we require the kernels to be homogeneous (or the operators, to be stable operators). 

\index{Free boundary!Regularity}\index{Free boundary}
\begin{thm} \label{thm:FB} 
Let $s\in(0,1)$, $\alpha\in (0, \min\{s, 1-s\})$, and let $\L\in \LLh $. 
There exists $\gamma > 0$ depending only on $n$, $s$, $\lambda$, and $\Lambda$, such that the following holds. 

Let $\varphi\in C^{2,\vartheta}_c(\R^n)$, with $\vartheta>\max\{2s-1,0\}$, and let $v$ be the solution to the obstacle problem \eqref{obst-pb}-\eqref{obst-pb2}. 

Then, near any free boundary point $x_\circ\in \partial\{v>\varphi\}$ we have
\begin{enumerate}[leftmargin=1.5cm, label=(\roman*)]
\item \label{it:i_fb} either 
\[v(x)-\varphi(x) = a_\circ \big((x-x_\circ)\cdot \nu\big)_+^{1+s} + O(|x-x_\circ|^{1+s+\gamma})\]
for some $a_\circ>0$ and $\nu\in \mathbb S^{n-1}$,

\item \label{it:ii_fb} or
\[v(x)-\varphi(x) = O(|x-x_\circ|^{1+s+\alpha}).\]
\end{enumerate}
Moreover, the set of (regular) points $x_\circ$ satisfying (i) is an open subset of the free boundary, and it is locally a $C^{1,\gamma}$ manifold.
\end{thm}

\index{Regular points}\index{Degenerate points}Free boundary points satisfying \ref{it:i_fb} are called \textbf{regular points}, while those satisfying \ref{it:ii_fb} are \textbf{degenerate points}.

This result was first established in \cite{ACS} for the operator $\L =\sqrt{-\Delta}$, then in \cite{CSS} in case $\L =(-\Delta)^s$ for all $s\in(0,1)$, and finally in \cite{CRS} for general operators of the form \eqref{obst-op1}-\eqref{obst-op2}-\eqref{obst-op3}; see also \cite{FR2,FRS} for the case of non-symmetric operators.

On the other hand, we have the following higher regularity result, which yields that free boundaries are actually $C^\infty$ near regular points.

\index{Free boundary!Higher regularity}
\begin{thm} \label{thm:FB2}
Let $s$, $\L $, $\varphi$, and $v$ be as in Theorem \ref{thm:FB}.
Assume in addition 
\[\varphi\in C^\infty_c(\R^n)\qquad \textrm{and}\qquad K|_{\mathbb S^{n-1}} \in C^\infty(\S^{n-1}).\]
Then, the set of regular free boundary points is locally a $C^\infty$ manifold.
\end{thm}

This result was first established in \cite{KPS,DS} for $\L =\sqrt{-\Delta}$, then in \cite{KRS,JN} in case $\L =(-\Delta)^s$ for all $s\in(0,1)$, and finally in \cite{AR20} for general operators of the form \eqref{obst-op1}-\eqref{obst-op2}-\eqref{obst-op3}.
Its proof is completely independent from that of Theorem \ref{thm:FB}, as we will see later on.

On the other hand, \emph{after} we prove Theorem \ref{thm:FB} we will establish the optimal $C^{1+s}$ regularity of solutions.
More precisely, combining a quantitative version of the free boundary regularity result with an iterative argument from \cite{CFR,FRS}, we will prove that 
\[\|v\|_{C^{1,s}(\R^n)} \leq C\|\varphi\|_{C^{2,\vartheta}(\R^n)},\]
with $C$ depending only on $n$, $s$, $\lambda$ and $\Lambda$.
This will be done in Section~\ref{sec-optimal-reg}.

As said before, the main goal of this section is to prove Theorem \ref{thm:FB}.
We first give an overview of the whole argument, to then provide a detailed proof in several steps.

\subsection{Sketch of the proof}

To establish Theorem \ref{thm:FB}, we consider 
\[u:=v-\varphi,\]
which satisfies
\[
u\geq0\qquad \textrm{and} \qquad D^2 u \ge -C_\circ \,{\rm{Id}}\quad \textrm{in}\quad \R^n,
\]
and
\begin{equation} \label{rgnn2}
\L u=f\quad \textrm{in}\quad \{u>0\} \qquad \textrm{and} \qquad \L u\geq f \quad \textrm{in} \quad \R^n,
\end{equation}
with 
\[
\|f\|_{{\rm Lip}(\R^n)} \leq C_\circ\qquad \textrm{and}\qquad \|u\|_{{\rm Lip}(\R^n)} \le C_\circ.
\]
Moreover, dividing $u$ by a constant if necessary we may assume $C_\circ=1$, and up to a translation we may assume $x_\circ=0\in \partial\{u>0\}$.

The proof then goes as follows: assume \ref{it:ii_fb} does \emph{not} hold at $x_\circ=0$, and let us show that \ref{it:i_fb} must then hold.

For this, we need a \emph{blow-up argument}: we consider the rescaled solutions 
\[u_r(x):= \frac{u(rx)}{\|u\|_{L^\infty(B_r)}},\]
and then try to take limits $r\to 0$ to get 
\[u_r\longrightarrow u_0\]
for some global solution $u_0$.
After that, we would like to classify all possible blow-ups $u_0$ (hopefully prove that they must all be of the form $u_0(x)=\kappa(x\cdot e)_+^{1+s}$ for some $\kappa>0$ and $e\in \S^{n-1}$), to then infer new information about our original solution $u$ near $0\in \partial\{u>0\}$.

This strategy, which is common in many free boundary problems (and originates in the classical theory of minimal surfaces), encounters new difficulties here, and thus requires several new ideas.

First, we choose carefully a sequence $r_k\to 0$ which, thanks to the assumption that \ref{it:ii_fb} does \emph{not} hold, allows us to prove that 
\[\|u\|_{L^\infty(B_{r_k})} \gtrsim r_k^{1+s+\alpha},\]
and more importantly that
\[|u_{r_k}(x)| \leq C\big(1+|x|^{1+s+\alpha}\big) \quad \textrm{in}\quad \R^n\]
and 
\begin{equation}\label{moderate-growth}
|\nabla u_{r_k}(x)| \leq C\big(1+|x|^{s+\alpha}\big) \quad \textrm{in}\quad \R^n,
\end{equation}
with a constant $C$ which is independent of $k$.
This is an important observation, and is proved in Lemma \ref{lem-growth-ell}.

Thanks to this, we will have that 
\[D^2 u_{r_k} \gtrsim -\frac{r_k^2}{r_k^{1+s+\alpha}}{\rm Id} \longrightarrow 0 \]
as $r_k\to 0$, which already gives enough \emph{compactness} on the sequence $u_{r_k}$   to have 
\[ u_{r_k} \longrightarrow u_0\quad\text{locally uniformly in $\R^n$}.\]
Moreover, the blow-up $u_0$ is \emph{convex}, and satisfies 
\begin{equation}\label{moderate-growth-0}|\nabla u_0(x)| \leq C\big(1+|x|^{s+\alpha}\big) \quad \textrm{in}\quad \R^n.
\end{equation}

The next step is to pass the equation to the limit, to show that $u_0$ is (in some sense) a global solution to the obstacle problem for the operator~$\L $.
This is a priori not obvious, since the equations in \eqref{rgnn2} do not make sense for a function $u_0$ that grows too much.
To solve this, we write \eqref{rgnn2} in terms of the gradient $\nabla u_{r_k}$, which has a moderate growth \eqref{moderate-growth} (since $\alpha<s$).
Namely, we have that\footnote{This is because, since $\L u=f$ in $\{u>0\}$ and $\L u\geq f$ everywhere,  we have $\L \big((u(x+h)-u(x)\big)\geq f(x+h)-f(x)$ in $\{u>0\}$.}
\[\L (\nabla u_{r_k})=\bar f_k\quad \textrm{in}\quad \{u_{r_k}>0\} \qquad \textrm{and} \qquad \L (D_h u_{r_k})\geq \bar f_{k,h} \quad \textrm{in} \quad \{u_{r_k}>0\},\]
for some $\bar f_k\to0$ and $\bar f_{k,h}\to0$, 
where
\[D_h w(x) := \frac{w(x+h)-w(x)}{|h|}.\]

This formulation is very good for our purposes, since we can then pass to the limit $r_k\to0$ and prove that the blow-up $u_0$ solves
\[\L (\nabla u_0)=0\quad \textrm{in}\quad \{u_0>0\} \qquad \textrm{and} \qquad \L (D_h u_0)\geq 0 \quad \textrm{in} \quad \{u_0>0\}.\]
This, combined with the convexity of $u_0$ and the growth \eqref{moderate-growth-0}, turns out to be enough for us to completely classify blow-ups and prove that  
\[u_0(x)=\kappa(x\cdot e)_+^{1+s},\]
see Figure~\ref{fig:13} and Proposition \ref{prop-classification-ell} below.

\begin{figure}
\centering
\makebox[\textwidth][c]{\includegraphics[scale = 1]{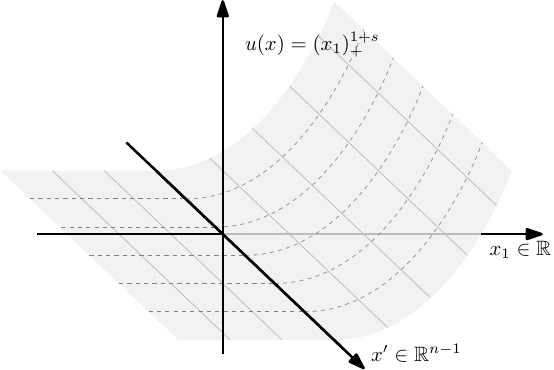}}
\caption{\label{fig:13} Blow-ups at regular points look like half-space solutions with growth $1+s$.}
\end{figure}

The next step is to transfer such information to the original solution $u$ near $0$.
As in the classical obstacle problem (corresponding to $s=1$ and $\L =-\Delta$), the key idea is to prove that 
\begin{equation}\label{positive-der}
\partial_\tau u_{r_k} \geq0\quad \textrm{in}\quad B_1
\end{equation}
for any $\tau\in \S^{n-1}$ with $\tau\cdot e\geq\frac12$, and for $r_k>0$ small enough.
This implies that the free boundary $\partial\{u_{r_k}>0\}$ is Lipschitz in $B_1$, and an appropriate application of the boundary Harnack in Lipschitz domains (Theorem \ref{thm-Lip}) yields the $C^{1,\gamma}$ regularity of the free boundary, as wanted.

The proof of \eqref{positive-der}, however, requires an extra ingredient.
Namely, it is convenient to have $C^1$ convergence of $u_{r_k}$ to the blow-up $u_0$, and this does not follow from the previous arguments.
We solve this by the following simple observation: the classification of blow-ups, combined with the contradiction and compactness arguments from Chapter~\ref{ch:2}, yields (almost-optimal) $C^{1,\mu}$ regularity estimates for solutions of \eqref{obst-pb}, for any $\mu<s$.

The optimal $C^{1,s}$ regularity estimates will be proved later on, in Section~\ref{sec-optimal-reg}.
It is actually interesting to notice that, unlike in many other free boundary problems, in order to establish the optimal regularity of solutions here we really need first the regularity of free boundaries near regular points.

\subsection{Classification of blow-ups}

\index{Classification of blow-ups}
Our first goal towards the proof of Theorem \ref{thm:FB}  will be to prove the classification of blow-ups.

\begin{prop}\label{prop-classification-ell}
Let $s \in(0,1)$, $\L\in \LLh$, and $\alpha\in (0, \min\{s, 1-s\})$.
Let $u_\circ\in {\rm Lip}_{\rm loc}(\R^n)$ be any function satisfying:

\begin{itemize}[leftmargin=0.8cm]
\item   $u_\circ$ is nonnegative and convex in $\R^n$:
\[u_\circ\geq0\quad \textrm{and} \quad D^2 u_\circ \ge 0\quad \textrm{in}\quad \R^n,\quad \textrm{with}\quad 0\in \partial\{u_\circ>0\}.\]

\item  $u_\circ\in C^1(\{u_\circ > 0\})$ solves, in the viscosity sense,
\[\L (\nabla u_\circ)=0\quad \textrm{in}\quad \{u_\circ>0\} \quad \textrm{and} \quad \L (D_h u_\circ)\geq 0 \quad \textrm{in} \quad \{u_\circ>0\},\]
for any $h\in \R^n$, where 
\[D_h u_\circ(x)={\textstyle \frac{u_\circ(x+h)-u_\circ(x)}{|h|}}.\]

\item  $u_\circ$ has a controlled growth at infinity:
\[\|\nabla u_\circ\|_{L^\infty(B_R)} \le R^{s+\alpha} \quad \textrm{for all}\quad R\ge 1.\]
\end{itemize}
Then, 
\[u_\circ(x)= \kappa (x\cdot e)_+^{1+s},\]
for some $\kappa\geq0$ and $e\in \mathbb S^{n-1}$. 
\end{prop}

\begin{rem}
The condition $\alpha<s$ is needed in order to have $\nabla u_\circ \in L^1_{w_s}(\R^n)$, so that we can evaluate $\L (\nabla u_\circ)$.
On the other hand, we also need $1+s+\alpha<2$ (i.e., $\alpha<1-s$) in order to avoid other solutions with quadratic growth, like $u_\circ(x)= x^TAx$ with $A\geq0$.
They correspond to degenerate points.
\end{rem}

To establish Proposition \ref{prop-classification-ell}, we need the following important consequence of the boundary Harnack principle.

\begin{prop}\label{prop-uniq-positive-sols}
Let $s\in(0,1)$ and $\Omega\subset\R^n$ be any  open set with $0\in \partial\Omega$.
Assume that $\Omega$ is unbounded and satisfies
\[B_{\varrho R}(-Re) \subset \Omega \quad \textrm{for some}\quad e\in \S^{n-1} \quad \textrm{and all}\quad R\geq1,\]
for some $\varrho > 0$.  Let $\L\in \LLL$, and $v_1,v_2\in C(\R^n)\cap L^1_{\omega_s}(\R^n)$ be any two  viscosity solutions of
\[
\left\{ \begin{array}{rcll}
\L v_i &=&0&\quad \textrm{in}\quad \Omega\\
v_i&=&0&\quad \textrm{in}\quad \R^n\setminus\Omega\\
v_i&\geq&0&\quad \textrm{in}\quad \R^n,
\end{array}\right.
\]
with $v_2\not\equiv0$.
 Then,
\[ v_1 \equiv \kappa\,v_2\quad\textrm{in}\quad \R^n\]
for some $\kappa\geq0$.
\end{prop}

\begin{proof}
The result will follow by applying Theorem \ref{thm-main} at scale $R$, and letting $R\to\infty$, as shown next.

Notice first that by the strong maximum principle we have either $v_i>0$ in $\Omega$ or $v_i\equiv0$ in $\R^n$.
Thus, we may assume $v_1>0$ in $\Omega$; otherwise there is nothing to prove.

\begin{steps}
\item We claim that, if $v_1\not\equiv0$, then there exists a constant $M>0$ (depending on $v_1$ and $v_2$) such that 
\begin{equation}\label{askjg}
M^{-1}v_1(x) \leq v_2(x) \leq Mv_1(x)\quad \textrm{in}\quad \R^n.
\end{equation}

To prove this, for any $R\geq2$ we consider the rescaled functions
\[v_1^{(R)}(x):= \frac{v_1(Rx)}{C_1^{(R)}},\qquad v_2^{(R)}(x):= \frac{v_2(Rx)}{C_2^{(R)}},\]
where the constants $C_i^{(R)}>0$ are chosen so that $\|v_i^{(R)}\|_{L^1_{w_s}(\R^n)}=1$.

Then, the functions $v_i^{(R)}$ satisfy the hypotheses of Theorem \ref{thm-main} in $B_1$ (with $\delta = 0$), and therefore we have
\[C^{-1}v_1^{(R)} \leq v_2^{(R)} \leq Cv_1^{(R)} \quad \textrm{in}\quad B_{1/2},\]
with $C$ independent of $R$.
Rescaling back to $v_1$ and $v_2$, this means that 
\[\frac{C^{-1}v_1(x)}{C_1^{(R)}} \leq \frac{v_2(x)}{C_2^{(R)}} \leq \frac{Cv_1(x)}{C_1^{(R)}}\quad \textrm{in}\quad B_{R/2}\]
for all $R\geq2$.
Evaluating at a fixed point $x_\circ\in \Omega\cap B_1$, we find 
\[\frac{C^{-1}v_1(x_\circ)}{v_2(x_\circ)} \leq \frac{C_1^{(R)}}{C_2^{(R)}} \leq \frac{Cv_1(x_\circ)}{v_2(x_\circ)},\]
and hence the quotient $C_1^{(R)}/C_2^{(R)}$ is uniformly   positive and uniformly bounded as $R\to \infty$.
In particular, we deduce that  
\[M^{-1}v_1(x) \leq v_2(x) \leq Mv_1(x)\quad \textrm{in}\quad B_{R/2},\]
with $M>0$ independent of $R$.
Letting $R\to \infty$, we find \eqref{askjg}.

\item  We now use  \eqref{askjg} to establish the result.
Indeed, let 
\[\kappa^* := \sup\{ \kappa>0 : v_1\geq \kappa v_2\quad \textrm{in}\quad \R^n\},\]
and define 
\[v_3:= v_1-\kappa^* v_2 \geq0.\]
Then, if $v_3\not\equiv0$, the function $v_3$ satisfies the same assumptions as $v_2$, and therefore with the exact same argument as before, we have
\[
N^{-1}v_3(x) \leq v_2(x) \leq Nv_3(x)\quad \textrm{in}\quad \R^n\]
for some $N>0$.
This would yield that $v_1\geq (\kappa^*+N^{-1})v_2$ in $\R^n$, which contradicts the definition of $\kappa^*$.
Therefore, it must be 
\[v_3\equiv0\quad \textrm{in}\quad \R^n,\]
and hence $v_1\equiv \kappa^* v_2$, as wanted.\qedhere
\end{steps}
\end{proof}

Another ingredient in the proof of Proposition \ref{prop-classification-ell} is the following.

\begin{lem}\label{lem-solution-across}
Let $s\in(0,1)$, $\L \in \LLL$, $e\in \mathbb S^{n-1}$, and $\Gamma\subset \{x\cdot e=0\}$.
Assume $w\in {\rm Lip}_{\rm loc}(\R^n)$ is a viscosity solution of
\[\L w \le 0 \quad \textrm{in}\quad \R^n\setminus \Gamma,\]
Then, $\L w\le 0$ in $\R^n$.
\end{lem}

\begin{proof}
For any $\varepsilon>0$ we consider the function $w_\varepsilon := w-\varepsilon\phi$, where $\phi$ is given by Lemma~\ref{lem-supersol1D}.
We claim that, by the choice of $\phi$ (and since $w$ is locally Lipschitz), we have $\L w_\varepsilon \leq C\varepsilon$ in $\R^n$ in the viscosity sense.

\begin{figure}
\centering
\makebox[\textwidth][c]{\includegraphics[scale = 1]{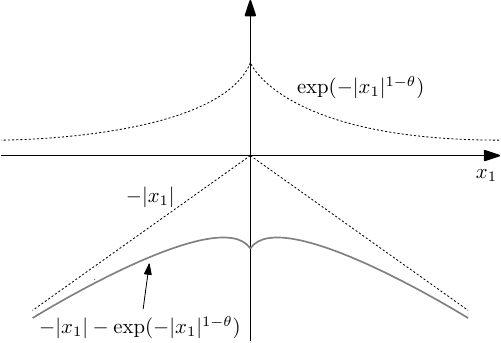}}
\caption{\label{fig:14} If $w$ is Lipschitz (in this picture, $w = -|x_1|$), the function $w-\eps\phi$ always has a positive wedge on $x_1 = 0$, that is, a wedge that cannot be touched from above by  a smooth ($C^1$) function.}
\end{figure}

Indeed, assume that a test function $\eta\in C^2$ touches $w_\varepsilon$ from above at $x_\circ\in \R^n$.
Since $w$ is Lipschitz,  by definition of $\phi$ we have that $w_\varepsilon$ has a ``positive wedge'' on $\{x\cdot e=0\}$, and therefore $x_\circ\notin \{x\cdot e=0\}$ (see Figure~\ref{fig:14}).
Then, we have $\L \eta(x_\circ) = \L w(x_\circ)-\varepsilon \L \phi(x_\circ)\leq C\varepsilon$.
This holds for every test function $\eta\in C^2$, so we deduce that $\L w_\varepsilon \leq C\varepsilon$ in $\R^n$ in the viscosity sense.
Since $w=\sup_{\varepsilon>0} w_\varepsilon$, we obtain that $\L w\leq0$ in $\R^n$.
\end{proof}

The following elementary lemma will allow us to take appropriate blow-down sequences.

\begin{lem}\label{lem-blow-down}
Let $\mu>0$, and let $w\in {\rm Lip}_{\rm loc}(\R^n)$ be any function satisfying 
\[\qquad\qquad \|\nabla w\|_{L^\infty(B_R)} \leq R^\mu \qquad \textrm{for all}\quad R\geq1. \]
Then, if $\nabla w\not\equiv0$, there exists a sequence $R_m\to\infty$ such that the rescaled functions
\[w_m(x):= \frac{w(R_m x)}{R_m \|\nabla w\|_{L^\infty(B_{R_m})}}\]
satisfy $\|\nabla w_m\|_{L^\infty(B_1)}=1$ and
\[\|\nabla w_m\|_{L^\infty(B_R)} \leq 2 R^\mu \qquad \textrm{for all}\quad R\geq1. \]
\end{lem}

\begin{proof}
Given any $m\in \mathbb N$ there must be $R_m\geq m$ such that 
\[\frac{\|\nabla w\|_{L^\infty(B_{R_m})}}{R_m^\mu} \geq \frac12\sup_{R\geq m} \frac{\|\nabla w\|_{L^\infty(B_R)}}{R^\mu} \geq \frac12\sup_{R\geq R_m} \frac{\|\nabla w\|_{L^\infty(B_R)}}{R^\mu} .\]
Therefore, we have a sequence $R_m\to\infty$ such that 
\[\frac{\|\nabla w\|_{L^\infty(B_{R_m})}}{R_m^\mu} \geq \frac12\frac{\|\nabla w\|_{L^\infty(B_{R})}}{R^\mu} \qquad \textrm{for all}\quad R\geq R_m.\]
Thus, the functions $w_m$ defined above satisfy
\[\|\nabla w_m\|_{L^\infty(B_R)} = \frac{\|\nabla w\|_{L^\infty(B_{R_mR})}}{\|\nabla w\|_{L^\infty(B_{R_m})}} \leq 2R^\mu\]
for any $R\geq1$, and we are done.
\end{proof}

We will also use the following basic lemma about convex functions:

\begin{lem}\label{lem-convex-functions}
Let $\mu>0$, and let $u_m\in {\rm Lip}_{\rm loc}(\R^n)$ be any sequence of convex functions, with 
\[\|\nabla u_m\|_{L^\infty(B_1)} =1\quad \textrm{and}\quad 
\|\nabla u_m\|_{L^\infty(B_R)}\leq 2R^{\mu}\quad \textrm{for all}\quad R\geq1.\]
Then, up to a subsequence, the functions $u_m$ converge locally uniformly to a convex function $u_\infty$ that satisfies 
\[\|\nabla u_\infty\|_{L^\infty(B_2)} \geq 1\quad \textrm{and}\quad 
\|\nabla u_\infty\|_{L^\infty(B_R)}\leq 2R^{\mu}\quad \textrm{for all}\quad R\geq1.\]
\end{lem}

\begin{proof}
The local uniform convergence of $u_m$ to $u_\infty$ follows from the uniform Lipschitz bounds for $u_m$ in balls $B_R$ and the Arzel\`a-Ascoli Theorem.
Moreover, the bound for $\nabla u_\infty$ in $B_R$ follows by uniform convergence (and the corresponding bound for $u_m$).
It therefore only remains to prove that $\|\nabla u_\infty\|_{L^\infty(B_2)} \geq 1$.
For this, notice that since $\|\nabla u_m\|_{L^\infty(B_1)} =1$, by convexity in $B_2$ it follows that $|u_m(x_m) - u_m(y_m)|\ge 1$ for some $x_m, y_m\in B_2$ with $|x_m - y_m| = 1$, and by uniform convergence  we deduce $|u_\infty(x_\infty)-u_\infty(y_\infty)| \geq 1$ for some $x_\infty, y_\infty\in B_2$ with $|x_\infty - y_\infty| = 1$. But then, using the mean value theorem, we get that $\|\nabla u_\infty\|_{L^\infty(B_2)} \geq 1$, as wanted.  
\end{proof}

Finally, we will also use the following property of (possibly unbounded) convex sets.

\begin{lem}\label{lem-blow-down-convex-sets}
Let $\Omega\subset\R^n$ be any convex set with $0\in\partial\Omega$.
Then, the rescaled sets $\Omega_R:= \frac1R\Omega$ converge as $R\to\infty$ to a convex cone $\Omega_\infty$.
\end{lem}

\begin{proof}
Notice that, for any $R_1\leq R_2$, by convexity we have $\Omega_{R_2} \subset \Omega_{R_1}$.
In particular, the limit is given by $\Omega_\infty := \bigcap_{R>0} \Omega_R$, and such set clearly satisfies $\rho\Omega_\infty=\Omega_\infty$ for any $\rho>0$.
This means that it is a cone, and the lemma is proved.
\end{proof}

We can now give the proof of the classification of blow-ups:

\begin{proof}[Proof of Proposition \ref{prop-classification-ell}]
Let us assume $u_\circ\not\equiv 0$, otherwise we are done. First, notice that the set $\{u_\circ=0\}$ is convex.
We separate the proof into two cases:

\vspace{2mm}

\noindent \emph{Case 1}.
Assume that the convex set $\{u_\circ=0\}$ contains a closed convex cone $\Sigma$ with nonempty interior (and with $0\in \Sigma$).
Then, there are $n$ independent directions $e_i\in \mathbb S^{n-1}$, $i=1,...,n$, such that $-e_i\in \Sigma$ (see Figure~\ref{fig:15}). 

\begin{figure}
\centering
\makebox[\textwidth][c]{\includegraphics[scale = 1]{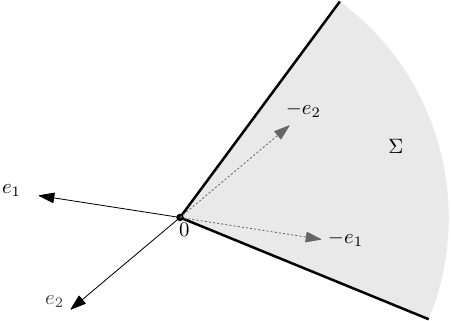}}
\caption{\label{fig:15} A cone $\Sigma$ with nonempty interior has $n$ independent directions with $-e_i\in \Sigma$, for $1\le i \le n$.}
\end{figure}

By convexity, we will have 
\begin{equation}\label{pweojrng}
v_i:=\partial_{e_i} u_\circ \geq0\quad \textrm{in}\quad\R^n.
\end{equation}
Indeed, we know that  $\partial^2_{e_i e_i} u_\circ\geq0$ and hence $\partial_{e_i}u_\circ$ is monotone in the $e_i$ direction.
Since, for any $x\in \R^n$, we will have $x-te_i\in \Sigma$ for all $t>0$ large enough, the function $(\partial_{e_i} u_\circ)(x-t e_i)$ is monotone (nonincreasing) in $t$ and converges to $0$ as $t\to+\infty$.
In particular, we have $\partial_{e_i} u_\circ(x)\geq0$, for any $x\in \R^n$, i.e., \eqref {pweojrng} holds.
Moreover, not all $v_i$ can be identically zero, and thus we assume $v_n\not\equiv0$.

We observe that $v_i$ are continuous functions, since (by convexity) $\{u_\circ>0\}$ is a Lipschitz domain and we can therefore use  Theorem~\ref{thm:existence_visc}.

Then, we can apply\footnote{The set $\Omega=\{u_\circ>0\}$ is the complement of a convex set, and in particular it satisfies the hypothesis of Theorem \ref{prop-uniq-positive-sols}.} Theorem \ref{prop-uniq-positive-sols} above to the functions $v_i$ and $v_n$, to deduce that 
\[v_i\equiv \kappa_i v_n,\quad \textrm{for}\quad i=1,...,n-1.\]
This means that, for all $i=1,...,n-1$ we have 
\[\partial_{e_i-\kappa_i e_n} u_\circ \equiv0\quad \textrm{in}\quad \R^n,\]
and thus $u_\circ$ is invariant in $n-1$ directions.
In other words,  $u_\circ$ is a (monotone) one-dimensional function,
\[u_\circ(x) = U(x\cdot e)\]
for some $e\in \S^{n-1}$, with $U\in C^1(\R)$.
In particular, we have $\{u_\circ>0\}=\{x\cdot e>0\}$.

Now, thanks to Lemma~\ref{lem:L1d_2}, the one-dimensional function $w=\partial_e u_\circ$ is continuous, solves the equation $\fls w=0$ in the viscosity sense  in a half-space $\{x\cdot e>0\}$ (and hence, in the classical or strong sense; see Lemma~\ref{lem:visc_dist} and Theorem~\ref{thm-interior-linear-2}), and $w=0$ in $\{x\cdot e\leq 0\}$. Combined with the growth condition
\[\|w\|_{L^\infty(B_R)} \leq R^{s+\alpha},\]
this implies $w(t) = a(x\cdot e)_+^s$, thanks to the classification of 1D solutions in Theorem~\ref{thm:class_1d_fls}. 
Hence, 
\[u_\circ(x) = \frac{a}{1+s}(x\cdot e)_+^{1+s},\]
as wanted.

\vspace{2mm}

\noindent \emph{Case 2}.
Assume that the convex set $\{u_\circ=0\}$ does \emph{not} contain any convex cone with nonempty interior.

By Lemma \ref{lem-blow-down} we can find a sequence $R_m\to\infty$ such that
\[u_m(x):=\frac{u_\circ(R_mx)}{R_m\|\nabla u_\circ\|_{L^\infty(B_{R_m})}}\]
satisfies 
\[\|\nabla u_m\|_{L^\infty(B_1)}=1\quad \textrm{and}\quad \|\nabla u_m\|_{L^\infty(B_R)}\leq 2R^{s+\alpha}\quad\textrm{for all}\quad R\geq1.\]
Moreover, we have $\L (D_h u_m)\ge 0$ in $\{u_m>0\}=\frac{1}{R_m}\{u_\circ>0\}$.  

By convexity of the functions $u_m$ (see Lemma \ref{lem-convex-functions}), up to a subsequence they converge locally uniformly to a function $u_\infty$ that satisfies  
\[\|\nabla u_\infty\|_{L^\infty(B_2)}\geq1\quad \textrm{and}\quad \|\nabla u_\infty\|_{L^\infty(B_R)}\leq 2R^{s+\alpha}\quad\textrm{for all}\quad R\geq1.\]

Now, since $\{u_\circ=0\}$ is convex, it follows from Lemma \ref{lem-blow-down-convex-sets} that the ``blow-down'' sequence $\frac{1}{R_m}\{u_\circ=0\}$ converges to a convex cone $\Gamma$.
Since, by assumption, the convex set $\{u_\circ=0\}$ does not contain any cone with nonempty interior, we have that $\Gamma$ must be a convex cone with empty interior (thus contained in  a hyperplane) and 
\[\L (D_h u_\infty)\ge 0\quad \textrm{in}\quad \R^n\setminus \Gamma.\]
Since $D_h u_\infty$ is locally Lipschitz,  Lemma \ref{lem-solution-across} implies $\L (D_h u_\infty)\ge 0$ in $\R^n$ for any $h$, and taking $-h$ instead of $h$ we obtain
\[
\L (D_h u_\infty)= 0\quad\text{in}\quad \R^n, 
\]
for all $h\in \R^n$.

Finally, the growth control on $\nabla u_\infty$ implies that $\|D_h u_\infty\|_{L^\infty(B_R)} \leq CR^{s+\alpha}$. Observe that for every fixed $h$, the function $D_h u_\infty$ is locally Lip\-schitz and its gradient has a controlled growth at infinity. By interior estimates (Theorem~\ref{C^alpha-bmc} applied to incremental quotients of $D_h u_\infty$) we deduce that $D_h u_\infty$ is $C^{1+\alpha}$. Moreover, when $2s> 1$, we actually have (applying instead  Theorem~\ref{C1alpha}) that $D_h u_\infty$ is $C^{2s+\alpha}$ (cf. Remark~\ref{rem:visstrong}). Thus, in all cases we have that $D_h u_\infty$ is a classical solution. Therefore, by the Liouville-type theorem with growth, Corollary \ref{cor:Liouville_growth}, we deduce that $D_h u_\infty$ is constant for every $h\in \R^n$ (recall $s+\alpha < 1$).

This means that $u_\infty(x)=a\cdot x+b$, which together with  $u_\infty(0) =0$ and $u_\infty\geq0$ implies $u_\infty\equiv 0$, contradicting $\|\nabla u_\infty\|_{L^\infty(B_2)}\geq1$.
Thus, Case~2 cannot happen, and the proposition is proved. 
\end{proof}

\subsection{Almost-optimal regularity of solutions}

\index{Obstacle problem!Almost-optimal regularity}
Before studying the regularity of the free boundary, we next show that the classification of blow-ups implies the almost-optimal regularity of solutions. 
 
\begin{prop}\label{cor-almost-optimal-ell}
Let $s \in(0,1)$, $\L\in\LLh$, and $\alpha\in (0, \min\{s, 1-s\})$.

Let $u\in {\rm Lip}_{\rm loc}(\R^n)\cap L^1_{\omega_s}(\R^n)$ with $u(0) = 0$ satisfy in the viscosity sense: 
\[
\left\{
\begin{array}{rcll}
u & \ge & 0&\quad\text{in}\quad B_2,\\
D^2 u & \ge & -{\rm Id}& \quad\text{in}\quad B_2,
\end{array}
\right.\quad
\left\{
\begin{array}{rcll}
\L u & = & f&\quad\text{in}\quad \{u > 0\}\cap B_2,\\
\L u  & \ge & f& \quad\text{in}\quad B_2,
\end{array}
\right.
\]
for some $f\in {\rm Lip}(B_2)$ with $|\nabla f|\leq 1$ in $B_2$. Suppose, moreover, that
\[
\|\nabla u\|_{L^\infty(B_R)}\le R^{s+\alpha}\quad\text{for all}\quad R \ge 1. 
\] 
Then, we have
\[\|u\|_{C^{1+s-\eps}(B_1)} \leq C_\eps\qquad \textrm{for any}\quad \eps>0,\]
with $C_\eps$ depending only on $n$, $s$, $\alpha$, $\eps$, $\lambda$, and $\Lambda$.   In particular, $u$ is a strong solution.
\end{prop}

For this, we need the following lemma about blow-up sequences.

\begin{lem}\label{lem-growth-ell}
Let $\mu>0$, and assume $w_k:\R^n\to\R^\ell$ for $k\in \N$ are a family of functions such that   $\sup_{k\in \N} \|w_k\|_{L^\infty(\R^n)} < +\infty$, and 
\[\sup_{k\in \N}\,\sup_{r>0} \,\frac{\|w_k\|_{L^\infty(B_{r})}}{r^\mu}=\infty.\]

Then, there are subsequences $w_{k_m}$ and $r_m\to0$ such that $\|w_{k_m}\|_{L^\infty(B_{r_m})} \geq c r_m^\mu$ for some $c > 0$ independent of $m$, and for which the rescaled functions 
\[\tilde w_m(x):= \frac{w_{k_m}(r_m x)}{\|w_{k_m}\|_{L^\infty(B_{r_m})}}\]
satisfy
\[\big|\tilde w_m(x)\big| \leq C\big(1+|x|^\mu\big)\quad \textrm{in}\quad \R^n,\]
with $C=2$.
\end{lem}

\begin{proof}
For every $m\in \mathbb N$ let $k_m$ and $r_m \geq \frac1m$ be such that 
\[ \frac{\|w_{k_m}\|_{L^\infty(B_{r_m})} }{r_m^{\mu} }
\geq \frac12 \sup_k \sup_{r\geq \frac1m} r^{-\mu} \|w_k\|_{L^\infty(B_r)}
\geq \frac12 \sup_k \sup_{r\geq r_m} r^{-\mu} \|w_k\|_{L^\infty(B_r)},\]
and in particular, $r_m\downarrow 0$, since $w_k$ are globally uniformly bounded.  This means that, by construction of $r_m$ and $k_m$, we have
\begin{equation}
\label{eq:puttingRrm}
r_m^{-\mu}\|w_{k_m}\|_{L^\infty(B_{r_m})} \geq \frac12 r^{-\mu}\sup_{k\in \N} \|w_{k}\|_{L^\infty(B_{r})} \quad \text{for all}\quad  r\geq r_m,\ m\in \N,
\end{equation}
and so taking $r  = 1$ we get $\|w_{k_m}\|_{L^\infty(B_{r_m})} \geq c r_m^\mu$, where we have denoted
\[
c := \frac12 \sup_{k\in \N} \|w_k\|_{L^\infty(B_1)} > 0,
\]
which is independent of $m$. On the other hand, for any $R\geq1$ we have
\[\|\tilde w_m\|_{L^\infty(B_R)} = \frac{\|w_{k_m}\|_{L^\infty(B_{R r_m})}}{\|w_{k_m}\|_{L^\infty(B_{r_m})}} \leq 2R^\mu \]
by putting $r = R r_m$ and $k = k_m$ in \eqref{eq:puttingRrm}, and we are done.
\end{proof}

To prove the almost-optimal regularity, let us show the following compactness result for viscosity solutions with linear operators $\L\in \LLL$ (analogous to the result in Proposition~\ref{prop:stab_distr} for distributional solutions):
\begin{lem}
\label{lem:compactness_viscosity} \index{Compactness of linear operators!Viscosity solutions}
Let $(\L_k)_{k\in \N}$ be a sequence of operators with $\L_k \in \LLL$ for all $k\in \N$. Then, there exists a subsequence $k_j \to \infty$ and $\L_\infty\in \LLL$ such that 
\[
\L_{k_j}\rightharpoonup \L_\infty\quad\text{in}\quad \R^n, \quad \text{as}\quad j\to \infty,
\]
in the sense of Definition~\ref{defi:conv_I}. Moreover, if $\L_k$ satisfy \eqref{obst-op3}, then $\L_\infty$ satisfies \eqref{obst-op3} as well. 
\end{lem}
\begin{proof}
Let $u\in L^1_{\omega_s}(\R^n)$ be fixed, and after a translation and rescaling, let us assume that $u\in C^2(B_1)$. We want to show that  $\L_{k_j} u\to \L_\infty u$ uniformly in $B_{1/2}$, for some subsequence $k_j$ and some operator $\L_\infty\in \LLL$. 

We   proceed similarly to Proposition~\ref{prop:stab_distr}, and if $\L_k$ has kernel $K_k$, we consider the sequence of absolutely continuous measures 
\[
 \mu_k(dy) := \min\{1, |y|^2\}K_k(y)\, dy,
\]
which by Prokhorov's theorem (and the growth of the kernels), Theorem~\ref{thm:Prokhorov}, weakly converges (up to a subsequence, $k_j$) to some absolutely continuous measure $\nu(dy) = \nu_\infty(y)\, dy$. We then define $\L_\infty$ as the operator with kernel 
\[
K_\infty(y) := \frac{\nu_\infty(y)}{\min\{1, |y|^2\}},
\]
which satisfies $\L_\infty\in \LLL$ (and also $\L_\infty\in \LLh$ if $\L_k\in \LLh$). Let us now see that $\L_{k_j} u \to \L_\infty u$ uniformly in $B_{1/2}$. 

Indeed, suppose that this is not the case. In particular, there exists a sequence $x_m \in B_{1/2}$ and a subsequence $\L_m := \L_{k_{j_m}}$ (or $\mu_m := \mu_{k_{j_m}}$) such that 
\[
|\L_m u(x_m) - \L_\infty u(x_m)| \ge \eps_\circ > 0\quad\text{for all}\quad m\in \N,
\]
for some $\eps_\circ > 0$. That is, 
\[
  \left|\int_{\R^n} \frac{2u(x_m) - u(x_m+y) - u(x_m-y)}{\min\{1, |y|^2\}} \left(\mu_{m}(dy) - \nu(dy)\right)\right| \ge \eps_\circ > 0,
\]
for all $m\in \N$. After taking a subsequence, we assume $x_m\to x_\infty \in \overline{B_{1/2}}$ as $m\to \infty$. Now, if  $m$ is large enough and since $u\in C^2(B_1)$, the previous bound implies (using $\mu_k \rightharpoonup \nu$ up to subsequences)
\[
  I_m := \left|\int_{\R^n\setminus{B_{1/4}}} u(x_m+y)\left(K_{m}(y) - K_\infty(y)\right)\, dy \right| \ge \frac{\eps_\circ}{4} > 0.
\]
However, we also have 
\[
\begin{split}
I_m &\le  2\Lambda \int_{\R^n\setminus{B_{1/4}}} |u(x_m+y)-u(x_\infty+y)||y|^{-n-2s} dy \\
&  \qquad  + \left|\int_{\R^n\setminus{B_{1/4}}} u(x_\infty+y)\left(K_{m}(y) - K_\infty(y)\right)\, dy \right|,
\end{split}
\]
where both terms go to zero: the first one being translations in $L^1_{\omega_s}(\R^n)$, and the second one by the weak convergence of $K_m$ to $K_\infty$ in $B_{1/4}^c$. This is a contradiction, and therefore, $\L_{k_j}\rightharpoonup \L_\infty$, as we wanted to see. 
\end{proof}

We can now establish the almost-optimal regularity of solutions.

\begin{proof}[Proof of Proposition \ref{cor-almost-optimal-ell}] Let $\eta\in C^\infty_c(B_4)$ with $\eta \ge 0$ be a smooth cut-off function such that $\eta \equiv 1$ in $B_3$. Then, up to dividing by a constant and replacing $u$ by $w := u \eta$, we can assume that $\nabla u$ is globally bounded. 

Indeed, $w\in L^\infty(\R^n)$ with $w = u$ in $B_2$, and it satisfies 
\[
\left\{
\begin{array}{rcll}
w & \ge & 0&\quad\text{in}\quad B_2,\\
D^2 w & \ge & -{\rm Id}& \quad\text{in}\quad B_2,
\end{array}
\right.\quad
\left\{
\begin{array}{rcll}
\L w & = & \tilde f&\quad\text{in}\quad \{w > 0\}\cap B_2,\\
\L w  & \ge & \tilde f& \quad\text{in}\quad B_2,
\end{array}
\right.
\]
where $\tilde f := f - \L (u (1-\eta))$. Since $u$ is Lipschitz, we have that $u (1-\eta)$ is also Lipschitz with bounds
\[
\|\nabla \big(u(1-\eta)\big)\|_{L^\infty(B_R)} \le C R^{s+\alpha}\quad\text{for all}\quad R \ge 1,
\]
for some $C$ depending only on $\eta$. In particular, by taking incremental quotients we see that $ \L \big(D_h (u (1-\eta))\big)$ is bounded in $B_2$, by some constant $C$ depending only on $n$, $s$, $\alpha$, $\lambda$, and $\Lambda$, independently of $h$. Hence, $\tilde f$ is Lipschitz,
\[
|\nabla \tilde f|\le C\quad\text{in}\quad B_2. 
\]
Up to dividing $w$ by $C$, we have that $w$ satisfies the same assumptions as $u$ in $B_2$, and it is globally bounded. Thus, without loss of generality, let us assume that $\nabla u$ is globally bounded by a universal fixed constant. 

Let $\mu<s$.
We prove that at every free boundary point $x_\circ\in \partial\{u>0\}\cap B_1$ we have
\begin{equation} \label{oskdn}
|\nabla u(x)|\leq C|x-x_\circ|^\mu\quad\text{for}\quad x\in \R^n,
\end{equation}
with $C$ depending only on $n$, $s$, $\mu$, $\lambda$, and $\Lambda$. 

Indeed, assume by contradiction that such estimate fails.
Then, we can find sequences $u_k$, $\L_k$, and $f_k$, satisfying the assumptions (with $\nabla u_k$ globally bounded uniformly in $k$), with $0\in \partial\{u_k>0\}$, and such that 
\[\sup_k\sup_{r>0} \frac{\|\nabla u_k\|_{L^\infty(B_{r})}}{r^\mu}=\infty.\]
Then, by Lemma \ref{lem-growth-ell}, there are sequences $k_m$ and $r_m\to0$ such that the functions
\[\tilde u_m(x) := \frac{u_{k_m}(r_m x)}{r_m\|\nabla u_{k_m}\|_{L^\infty(B_{r_m})}},\qquad 
\nabla \tilde u_m(x) := \frac{\nabla u_{k_m}(r_m x)}{\|\nabla u_{k_m}\|_{L^\infty(B_{r_m})}},\]
satisfy $\|\nabla \tilde u_m\|_{L^\infty(B_1)}=1$ and
\[|\nabla \tilde u_m(x)| \leq C\big(1+|x|^\mu\big)\quad\textrm{in}\quad \R^n.\]

Moreover, we will also have, since $\|\nabla u_{k_m}\|_{L^\infty(B_{r_m})} \geq c r_m^\mu$ for some $c> 0$ independent of $m$, 
\begin{equation}
\label{eq:isgloballyconvex}
D^2 \tilde u_m \geq - c^{-1} r_m^{1-\mu}{\rm Id} \longrightarrow 0\quad\text{in}\quad B_{2/r_m},
\end{equation}
as well as 
\[
\left\{
\begin{array}{rcll}
\L_m \tilde u_m & =&  \tilde f_m& \quad\text{in}\quad \{ \tilde u_m>0\}\cap B_{2/r_m},\\
\L_m \tilde u_m & \ge&    \tilde f_m& \quad\text{in}\quad B_{2/r_m},
\end{array}
\right.
\]
with 
\[
\tilde f_m(x) := r_m^{2s-1} \frac{f_{k_m}(r_m x)}{\|\nabla u_{k_m}\|_{L^\infty(B_{r_m})}},
\]
such that $|\nabla \tilde f_m|\leq c^{-1} r_m^{2s-\mu}\to 0$ in $B_{2/r_m}$. 
In particular, this implies
\[\L_m(D_h \tilde u_m)\geq -c^{-1} r_m^{2s-\mu} \quad \textrm{in}\quad \{\tilde u_m>0\}\cap B_{1/r_m}.\]

By the control on the gradient, a subsequence of the functions $\tilde u_m$ converges locally uniformly in $\R^n$ to a limiting function $\tilde u_\circ$, which is globally convex by \eqref{eq:isgloballyconvex}. Observe, also, that 
\[\L_m(D_h \tilde u_m)\leq c^{-1} r_m^{2s-\mu} \quad \textrm{in}\quad \{x : \dist(x, \{\tilde u_m=0\}) > |h|\} \cap B_{1/r_m},\]
and so, letting $h\downarrow 0$ and thanks to the interior estimates in Theorem~\ref{C^alpha-bmc} together with Lemma~\ref{it:H8}, we deduce that $\tilde u_m$ is $C^{1,\alpha}$ in $\{\tilde u_m>0\}$. By Lemma~\ref{lem:compactness_viscosity} together with the stability of viscosity solutions (see Proposition \ref{prop:stab_super}), we actually have 
\[
\L_\infty (\nabla \tilde u_\circ) = 0 \quad\text{in}\quad \{\tilde u_\circ>0\},
\]
for some $\L_\infty\in \LLh$, so the first two hypotheses of Proposition \ref{prop-classification-ell} are satisfied. Furthermore, we have the growth
\[\|\nabla \tilde u_\circ\|_{L^\infty(B_R)} \leq CR^\mu\]
for all $R\geq1$.
In addition, by Lemma~\ref{lem-convex-functions}  we have $\|\nabla \tilde u_\circ\|_{L^\infty(B_2)}\geq 1$.
Thanks to Proposition \ref{prop-classification-ell}, we deduce that $\tilde u_\circ=\kappa(x\cdot e)_+^{1+s}$ with $\kappa \neq 0$, a contradiction since $\mu<s$.
Thus, \eqref{oskdn} is proved.

Finally,  \eqref{oskdn} implies that $|\nabla u|\leq Cd^\mu$, where $d(x)={\rm dist}(x,\{u=0\})$. Notice, also, that 
\[
\L D_h u = D_h f \quad\text{in}\quad \{x : \dist(x, \{u = 0\}\cap B_2)\ge |h|\}
\]
in the viscosity sense, and $D_h f\in L^\infty$, so that   applying Theorems~\ref{C^alpha-bmc} or  \ref{C1alpha} to the incremental quotients   $D_h u$ we deduce that $u$ is $C^{2s+\alpha}$ in $\{u > 0\}$ (cf. Remark~\ref{rem:visstrong}), and thus,   it is a classical solution. In particular, exactly as in \ref{it:step2bdry} in the proof of Proposition~\ref{prop:globCsreg_loc}, this combined with interior regularity estimates for distributional or classical solutions, Theorem~\ref{thm-interior-linear-Lp}, yields the desired result.
\end{proof}

\subsection{Regularity of the free boundary}
 \index{Free boundary!Regularity}
The next step is to show that the free boundary is $C^{1,\gamma}$ near nondegenerate points.
For this, we need the following\footnote{It is interesting to notice that this result fails when $s=1$, i.e., in case $\L=-\Delta$.}.

\begin{lem} \label{lem-ell-positivity}
Let $s\in(0,1)$, $\L \in \LLL$, $\alpha\in(0,s)$, $\varrho_\circ\ge1$, and $c_\circ>0$.
Then, for any $R_\circ\geq1$ large enough and $\eps_\circ>0$ small enough, depending only on $n$, $s$, $\lambda$, $\Lambda$, $\varrho_\circ$, $c_\circ$, and~$\alpha$, the following holds.

Assume that $E\subset\R^n$ is closed, and $v\in C(\R^n)$ satisfies (in the viscosity sense)
\[
\left\{
\begin{array}{rcll}
\L v&\geq& -\eps_\circ&\quad \textrm{in}\quad B_{R_\circ}\setminus E,\\
v&\equiv & 0&\quad \textrm{in}\quad B_{R_\circ}\cap E,\\
v &\geq &-\eps_\circ& \quad \textrm{in}\quad B_{R_\circ},
\end{array}
\right.
\]
and 
\[
 \int_{B_{\varrho_\circ}}v_+ \geq c_\circ>0
 \qquad\text{and}\qquad |v(x)|\leq 1+|x|^{s+\alpha}\quad \textrm{in}\quad \R^n\setminus B_{R_\circ}.\]
Then, $v\geq0$ in $B_{\varrho_\circ}$.
\end{lem}

\begin{proof}
Let $w:=v\chi_{B_{R_\circ}}$, and notice that, since 
\[\big|\L \big(v(1-\chi_{B_{R_\circ}})\big)(x)\big| \leq C\int_{B_{R_\circ/2}^c}\frac{1+|x+y|^{s+\alpha}}{|y|^{n+2s}}\,dy \leq CR_\circ^{\alpha-s}\]
for any $x\in B_{R_\circ/2}$,   we have
\[
\left\{
\begin{array}{rcll}
\L w&\geq& -\eps_\circ-CR_\circ^{\alpha-s}&\quad \textrm{in}\quad B_{R_\circ/2}\setminus E,\\
w&\equiv & 0&\quad \textrm{in}\quad B_{R_\circ}\cap E,\\
w &\geq &-\eps_\circ& \quad \textrm{in}\quad B_{R_\circ},
\end{array}
\right.
\quad\text{and}\quad \int_{B_{\varrho_\circ}}w_+ \geq c_\circ>0.
\]

Let $\psi \in C^\infty_c(B_{2})$ be some radial bump function with $\psi \ge 0$ and satisfying $\psi\equiv1$ in $B_1$.
Let us define, for any $t>0$, 
\[ \psi_t (x)  = - \varepsilon_\circ - t +  \varepsilon_\circ \psi(x/\varrho_\circ).\]

If the conclusion of the lemma does not hold then $\psi_t$ touches $w$ from below at $z\in  B_{2\varrho_\circ}$ for some $t>0$.
Since $\psi_t\le -t$ in all of $\R^n$ we have that $w(z) = \psi_t(z) <0$ and hence $z$ belongs to $B_{2\varrho_\circ} \setminus E$. 
Since $\psi_t$ is $C^2$ and touches $w$ from below at $z$,  by Remark~\ref{rem:onesidedcondition} the operator $\L$ can be evaluated for $w$ at the point~$z$.

On  the one hand, since $(w-\psi_t)(z)=0$, $w\ge \psi_t$, and $\psi_t\leq0$ in all of $\R^n$, we have
\[\begin{split}
 \L (w-\psi_t)(z) & \le -\lambda\int_{\R^n} (w-\psi_t)(z+y) |y|^{-n-2s}dy \\
 & \le  -c \varrho_\circ^{-n-2s} \int_{B_{\varrho_\circ}} w_+  \\
 & \le -c\varrho_\circ^{-n-2s} c_\circ,
 \end{split}\]
 for some $c>0$ depending only on $n$, $s$, and $\lambda$. 
On the other hand
\[\L (w-\psi_t)(z)\ge \L w(z) - |\L  \psi_t(z)| \geq -\varepsilon_\circ-CR_\circ^{\alpha-s} - C\varepsilon_\circ \varrho_\circ^{-2s}.\]
We obtain a contradiction by taking $R_\circ$ large enough and $\varepsilon_\circ$ small enough (depending only on $n$, $s$, $\lambda$, $\Lambda$, $\varrho_\circ$, $c_\circ$, and~$\alpha$).  
\end{proof}

In order to prove later   the optimal regularity of solutions, it will be crucial to have the following result, which follows from the classification of blow-ups.
 
\begin{prop}\label{thm-main-ell}
Let $s \in(0,1)$, $\alpha\in (0, \min\{s, 1-s\})$, $\L\in \LLh$, and let $\varepsilon_\circ>0$ and $R_\circ>1$.
Then, there exists $\eta>0$ depending only on $n$, $s$, $\alpha$, $\eps_\circ$, $R_\circ$, $\lambda$, and $\Lambda$, such that the following holds.

Let $u\in {\rm Lip}_{\rm loc}(\R^n)\cap L^1_{w_s}(\R^n)$ satisfy:

\vspace{1mm}

$\bullet$\  $u$ is nonnegative and almost-convex:  
\begin{equation}
\label{eq:u(1)}u\geq0\quad \textrm{and} \quad D^2 u \ge - \eta \,{\rm{Id}}\quad \textrm{in}\quad \R^n,\quad \textrm{with} \quad 0\in \partial\{u>0\}.
\end{equation}

$\bullet$\  $u$ solves the obstacle problem with a small right-hand side:
\begin{equation}
\label{eq:u(2)}\L u=f~~~~\textrm{ in }~~~~ \{u>0\} \quad \textrm{and} \quad \L u\geq f ~~~~ \textrm{ in } ~~~~ \R^n, \quad \textrm{with } ~~~~ |\nabla f|\leq \eta.
\end{equation}

$\bullet$\  $u$ has a controlled growth at infinity:
\begin{equation}
\label{eq:u(3)}\|\nabla u\|_{L^\infty(B_R)} \le R^{s+\alpha} \quad \textrm{for all}\quad R\ge 1.
\end{equation}

\vspace{1mm}

\noindent Then, we have
\[
\big\|u- \kappa(x\cdot e)_+^{1+s}\big\|_{{\rm Lip} (B_{R_\circ})} \le \eps_\circ.
\]
for some $e\in \mathbb S^{n-1}$ and $\kappa\geq0$.
\end{prop}

\begin{proof}
Assume by contradiction that there is no $\eta>0$ for which the result holds.
Then, we have a sequence $\eta_k\to0$, and sequences of operators $\L_k$, functions $f_k$, and solutions $u_k$, satisfying the hypotheses with $\eta = \eta_k$ but such that
\[\big\|u_k-\kappa(x\cdot e)_+^{1+s}\big\|_{{\rm Lip}(B_{R_\circ})} \geq\eps_\circ\]
for any $e\in \mathbb S^{n-1}$ and any $\kappa\geq0$.
But then, by Proposition~\ref{cor-almost-optimal-ell}, up to a subsequence, $u_k$ converges in $C^{1+s-\eps}$ norm in compact sets to a limiting function~$u\in C^{1+s-\eps}$.
By   Proposition~\ref{prop:stab_distr} (or by Lemma~\ref{lem:compactness_viscosity} and Proposition~\ref{prop:stab_super}), there is an operator $\L\in \LLL $  such that $u$ satisfies the same conditions but with $\eta=0$ (since it is a strong solution, it is both a distributional and viscosity solution).
More precisely, we have $\L_k(u_k(x+h)-u_k(x))\geq f_k(x+h)-f_k(x)$ in $\{u_k>0\}$, and thus in the limit we get $\L(D_h u)\geq0$ in $\{u>0\}$.
Similarly we have $\L(D_{h} u)\leq0$ in $\{u(\cdot+h)>0\}$, and hence by taking $h\to0$ we deduce $\L(\nabla u)=0$ in $\{u>0\}$.
But then by Proposition~\ref{prop-classification-ell} it follows that $u(x)=\kappa(x\cdot e)_+^{1+s}$, a contradiction.
\end{proof}

We can now show the $C^{1,\gamma}$ regularity of free boundaries.

\begin{prop}\label{lemm-main-ell-2}
Let $s \in(0,1)$, $\alpha\in (0, \min\{s, 1-s\})$, $\L\in \LLh$, and $\kappa_\circ > 0$. Then, there exist $\eps  > 0$, $R_\circ > 1$, and $\eta > 0$, depending only on $n$, $s$, $\alpha$, $\kappa_\circ$, $\lambda$, and $\Lambda$, such that the following statement holds.

Let $u\in {\rm Lip}_{\rm loc}(\R^n)\cap L^1_{\omega_s}(\R^n)$ satisfy \eqref{eq:u(1)}-\eqref{eq:u(2)}-\eqref{eq:u(3)}, and assume that 
\[
\big\|u- \kappa (x\cdot e)_+^{1+s}\big\|_{{\rm Lip} (B_{R_\circ})} \le \eps,
\]
for some $\kappa\geq \kappa_\circ>0$ and $e\in \mathbb S^{n-1}$.

Then, the free boundary $\partial \{u>0\}$ is a $C^{1,\gamma}$ graph in $B_{1/2}$, and  moreover
\[\|\nabla u\|_{C^s(B_{1/2})}\leq C\]
and 
\[\|\nabla u/d^s\|_{C^\gamma(\overline{\{u>0\}} \cap B_{1/2})} \leq C,\]
with $C$   depending only on $n$, $s$, $\alpha$, $\kappa_\circ$, $\lambda$, and $\Lambda$; and   $\gamma>0$ depending only on $n$, $s$,   $\kappa_\circ$, $\lambda$, and $\Lambda$.
\end{prop}

\begin{proof}
 Let $u_\circ(x) := \kappa (x\cdot e)_+^{1+s}$. By assumption, for any direction $e'\in \mathbb S^{n-1}$ such that $e'\cdot e\geq\frac{1}{\sqrt{2}}$ we have
\[|\partial_{e'} u-\partial_{e'} u_\circ|\leq \eps\quad\textrm{in}\quad B_{R_\circ},\]
and 
\[\partial_{e'} u_\circ \geq0 \quad \textrm{in}\quad \R^n \qquad \textrm{and} \qquad \partial_{e'} u_\circ \geq c_1\kappa \quad \textrm{in}\quad \{x\cdot e\geq {\textstyle\frac{1}{\sqrt{2}}}\}.\]
Moreover, since $\L(D_h u)\geq D_h f\geq -\eta$ in $\{u>0\}$ and $\L(D_h u)\leq D_h f\leq \eta$ in $\{u(\cdot+h)>0\}$, then letting $h\to0$ (and using Proposition~\ref{prop:stab_super} or Proposition~\ref{prop:stab_distr}) we find $|\L(\nabla u)|\leq \eta$ in $\{u>0\}$.

Thus, if $\eps$ is small, we have that $w:=\partial_{e'} u$ and $E:=\{u=0\}$ satisfy
\[
\left\{
\begin{array}{rcll}
|\L w|&\leq& \eta &\quad \textrm{in}\quad B_{R_\circ}\setminus E,\\
w&\equiv & 0&\quad \textrm{in}\quad B_{R_\circ}\cap E,\\
w &\geq &-\eps& \quad \textrm{in}\quad B_{R_\circ},
\end{array}
\right.
\]
in the viscosity sense, and 
\[
 w\ge c_2 \kappa \quad \textrm{in}\quad \{x\cdot e\geq {\textstyle\frac{1}{\sqrt{2}}}\}\cap B_{R_\circ}\quad\Longrightarrow\quad  \int_{B_{\varrho_\circ}}w_+ \geq c c_2\kappa\varrho_\circ^n>0,\]
for some $\varrho_\circ$ to be chosen later, together with
\[|w(x)|\leq |x|^{s+\alpha}\quad \textrm{in}\quad \R^n\setminus B_{R_\circ}.\]
Then, choosing $R_\circ$ large enough, it follows from Lemma \ref{lem-ell-positivity} (taking $\eta$ and $\eps$ small enough, now depending on $\varrho_\circ$ as well) that $w\geq 0$ in $B_{\varrho_\circ}$, i.e. 
\[\partial_{e'} u\geq0 \quad \textrm{in}\quad B_{\varrho_\circ}\]
for all $e'\in \mathbb S^{n-1}$ such that $e'\cdot e\geq\frac{1}{\sqrt{2}}$.
This means that the free boundary $\partial\{u>0\}$ is a Lipschitz graph in $B_{\varrho_\circ}$, with Lipschitz constant bounded by~1. Indeed, let $x_\circ \in B_{\varrho_\circ}\cap \partial\{u >0\}$ be any free boundary point in $B_{\varrho_\circ}$, and let 
\[\Theta:=\bigl\{\tau\in \mathbb{S}^{n-1}: \tau\cdot e>{\textstyle\frac{1}{\sqrt{2}}}\bigr\},\]
\[\Sigma_\pm :=\bigl\{x\in B_{\varrho_\circ}: x=x_\circ \pm t\tau,\ {\rm with}\ \tau\in \Theta,\ t>0\bigr\},\]
see Figure~\ref{fig.29}.

\begin{figure}
\includegraphics[scale = 1]{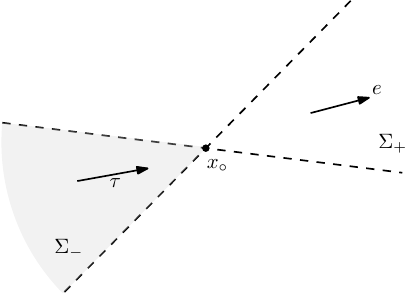}
\caption{Representation of $\Sigma_-$ and $\Sigma_+$.}
\label{fig.29}
\end{figure}

We claim that 
\begin{equation} \label{ch4-FB-Lip-Q}
\left\{
\begin{array}{rcll}
u &=&0& \ \textrm{in} \quad\Sigma_-,\\
u &>&0& \ \textrm{in} \quad\Sigma_+.
\end{array}
\right.
\end{equation}
On the one hand, since $u (x_\circ )=0$, it follows from the monotonicity property $\partial_\tau u \geq0$ (and the nonnegativity of $u $) that $u (x_\circ -t\tau)=0$ for all $t>0$ and $\tau\in\Theta$.
In particular, there cannot be any free boundary point in $\Sigma_-$.

On the other hand, by the same argument, if $u (x_1)=0$ for some $x_1\in \Sigma_+$ then we would have $u =0$ in $\bigl\{x\in B_{\varrho_\circ}: x=x_1-t\tau,\ {\rm with}\ \tau\in \Theta,\ t>0\bigr\}\ni x_\circ $, and in particular $x_\circ $ would not be a free boundary point.
Thus, $u (x_1)>0$ for all $x_1\in \Sigma_+$, and \eqref{ch4-FB-Lip-Q} is proved.

Notice that \eqref{ch4-FB-Lip-Q} yields that the free boundary $\partial\{u >0\}\cap B_{\varrho_\circ}$ satisfies both the interior and exterior cone condition (with an explicit cone), and thus it is Lipschitz (with constant 1).

Finally, if $\varrho_\circ>1$ is large enough (depending only on $n$, $s$, $\alpha$, $\lambda$, and $\lambda$), then  the functions $(\partial_{e'} u)\chi_{B_{\varrho_\circ}}$ and $(\partial_{e} u) \chi_{B_{\varrho_\circ}}$ satisfy the assumptions of  the boundary Harnack in Theorem \ref{thm-Lip} (up to dividing by a constant depending only on $n$, $s$, $\kappa_\circ$, $\lambda$, and $\Lambda$), and we deduce that
\[\left\|\frac{\partial_{e'} u}{\partial_{e} u}\right\|_{C^\gamma(\overline{\{u>0\}}\cap B_{1/2})} \leq C,\]
for some $\gamma > 0$ depending only on $n$, $s$, $\kappa_\circ$, $\lambda$, and $\Lambda$.  
Assuming $e=e_n$ and taking $e'=e_n+e_i$, for $i=1,...,n-1$, we find that the quotients $\partial_{i} u/\partial_n u$ are $C^\gamma$ in $\overline{\{u>0\}}\cap B_{1/2}$.

Now, notice that the normal vector to the level set $\{u=t\}$ is given by 
\[\nu^i = \frac{\partial_i u}{|\nabla u|} = \frac{\partial_i u/\partial_n u}{1+\sqrt{\sum_{j=1}^{n-1}(\partial_j u/\partial_n u)^2 }},\qquad i=1,...,n.\]
We have proved that this is a $C^\gamma(\overline{\{u>0\}}\cap B_{1/2})$ function, and therefore the level sets $\{u=t\}$ for $t>0$ are uniformly $C^{1,\gamma}$ graphs in $B_{1/2}$.
Taking $t\to0$, we deduce that the free boundary $\partial\{u>0\}$ is a $C^{1,\gamma}$ graph in $B_{1/2}$, as wanted.

Finally, the $C^s$ estimate for $\nabla u$ follows from Proposition~\ref{prop:globCsreg_loc}, and the $C^\gamma$ estimate for $\nabla u/d^s$ follows from Proposition \ref{prop:local_bdry}.
\end{proof}

Combining the previous results, we can finally give the proof of the free boundary regularity around regular points:

\begin{proof}[Proof of Theorem \ref{thm:FB}]
First, by Lemma~\ref{semiconvex}~\ref{it:vsemiconvex}, the solution $v$ is semiconvex.
Thus, up to dividing it by a constant, the function 
\[u:= v-\varphi\]
satisfies 
\[
\left\{
\begin{array}{rcll}
u & \ge & 0&\quad\text{in}\quad \R^n,\\
D^2 u & \ge & -{\rm Id}& \quad\text{in}\quad \R^n,
\end{array}
\right.\quad
\left\{
\begin{array}{rcll}
\L u & = & f&\quad\text{in}\quad \{u > 0\},\\
\L u  & \ge & f& \quad\text{in}\quad \R^n,
\end{array}
\right.
\]
 with $|\nabla f|\leq 1$ in $\R^n$ (we are also using here that $\L \varphi\in {\rm Lip}(\R^n)$, since $\varphi\in C^{2,\vartheta}_c(\R^n)$ with $\vartheta > \max\{0, 2s-1\}$; see Lemma~\ref{lem:Lu_2}).

Let  $x_\circ\in \partial\{u>0\}$, and assume that \ref{it:ii_fb} does \emph{not} hold at $x_\circ$.
Up to a translation we may assume $x_\circ=0$, and hence we have
\[\sup_{r>0} \frac{\|u\|_{L^\infty(B_r)}}{r^{1+s+\alpha}} = \infty.\]
Moreover, since $\|u\|_{L^\infty(B_r)} \leq r\|\nabla u\|_{L^\infty(B_r)}$, 
\[\sup_{r>0} \frac{\|\nabla u\|_{L^\infty(B_r)}}{r^{s+\alpha}} = \infty.\]

Then, by Lemma \ref{lem-growth-ell} (applied to a single function $\nabla u$, which is globally bounded, since $u$ is globally Lipschitz) there is a sequence of radii $r_m\to0$ such that $\|\nabla u\|_{L^\infty(B_{r_m})} \geq c r_m^{s+\alpha}$ for some $c > 0$  and such that the functions
\[\tilde u_m(x) := \frac{u(r_m x)}{r_m\|\nabla u\|_{L^\infty(B_{r_m})}},\qquad \nabla \tilde u_m(x)= \frac{\nabla u(r_m x)}{\|\nabla u\|_{L^\infty(B_{r_m})}},\]
satisfy $\|\nabla \tilde u_m\|_{L^\infty(B_1)}=1$ and 
\[\big|\nabla \tilde u_m(x)\big| \leq C\big(1+|x|^{s+\alpha}\big)\quad \textrm{in}\quad \R^n.\]
Moreover, we also have 
\[D^2 \tilde u_m = r_m\,\frac{D^2 u(r_m \cdot))}{\|\nabla u_m\|_{L^\infty(B_{r_m})}} \geq -c^{-1} r_m^{1-s-\alpha}{\rm Id}\longrightarrow 0\quad \textrm{in}\quad \R^n,\]
as well as 
\[
\left\{
\begin{array}{rcll}
\L \tilde u_m & = & f_m & \quad\text{in}\quad \{\tilde u_m > 0\}\\
\L \tilde u_m & \ge & f_m & \quad\text{in}\quad \R^n,
\end{array}
\right.
\quad\text{with}\quad f_m (x) := r_m^{2s-1}\frac{f(r_m x)}{\|\nabla u\|_{L^\infty(B_{r_m})}},
\]
which satisfies  $|\nabla f_m|\leq c^{-1} r_m^{s-\alpha}\to0$ in $\R^n$.

Thus, taking $m$ large enough, we can combine\footnote{Choose $\eps$, $R_\circ$, and $\eta$ from Proposition~\ref{lemm-main-ell-2}, which can then be applied by taking $\eta > 0$ smaller if necessary thanks to Proposition~\ref{thm-main-ell}.} Propositions \ref{thm-main-ell} and  \ref{lemm-main-ell-2} (with $\kappa_\circ = \frac12$, since $\|\nabla \tilde u_m\|_{L^\infty(B_1)}=1$) to deduce that the free boundary $\partial\{\tilde u_m>0\}$ is $C^{1,\gamma}$ in $B_{1/2}$, and $\nabla \tilde u_m/d_m^s$ is $C^\gamma$ up to the free boundary in $B_{1/2}$, where $d_m(x)={\rm dist}(x,\{\tilde u_m=0\})$, and $\gamma$ depends only on $n$, $s$, $\lambda$, and $\Lambda$.

Rescaling back to $u$, this means that the free boundary is $C^{1,\gamma}$ in a neighborhood of the origin, and that $\nabla u/d^s$ is $C^\gamma$ up to the free boundary in a neighborhood of the origin.
In particular, we have  
\[\big|\nabla u(x) - b_\circ d^s(x)\big| \leq C|x|^{s+\gamma}\]
for some $b_\circ\in \R^n$, given by $\lim_{x\to 0}\frac{\nabla u(x)}{d^s(x)}$, which yields
\[\big|u(x) - a_\circ d^{1+s}(x)\big| \leq C|x|^{1+s+\gamma}.\]
Since the free boundary is $C^{1,\gamma}$ near the origin,  
\[\big|d^{1+s}(x) - (x\cdot \nu)_+^{1+s}\big| \leq C|x|^{1+s+\gamma},\]
and the result follows.
\end{proof}

\subsection{Higher regularity of free boundaries}\label{ssec:high_reg_fb}
\index{Free boundary!Higher regularity}

Our next goal is to prove Theorem \ref{thm:FB2}.
The proof is based on the following. 

\begin{thm}\label{thm:u1overu2}
Let $s\in (0, 1)$  and let~$\beta>1$ be such that~$\beta,\beta\pm s\not\in\N$.
Let~$\Omega\subseteq\R^n$ be any bounded $C^\beta$ domain, and let~$\L \in \LLh$, with
$K|_{\S^{n-1}}\in C^{2\beta+1}(\S^{n-1})$.
Let~$f_1,f_2\in {C^{\beta-s}(\overline{\Omega}\cap B_1)}$, and $u_1,u_2\in L^\infty(\R^n)\cap C(B_1)$ be solutions of 
\begin{align*}
\left\lbrace\begin{aligned}
\L u_i &= f_i & & \hbox{in}\quad\Omega\cap B_1 \\
u_i &= 0 & & \hbox{in}\quad B_1\setminus\Omega,
\end{aligned}\right.
\end{align*}
with $u_2\geq c_1d^s$ in~$B_1$ for some $c_1>0$.
Then,
\[
\frac{u_1}{u_2} \in C^\beta(\overline\Omega\cap B_{1/2}).
\]
\end{thm}

We will not prove Theorem \ref{thm:u1overu2} here.
Its proof uses a higher-order version of the ideas we used to prove Proposition \ref{prop:local_bdry}; see \cite{AR20} for more details.

We next show how to use Theorem \ref{thm:u1overu2} in order to prove the higher regularity of free boundaries:

\begin{proof}[Proof of Theorem \ref{thm:FB2}]
Notice first that $v\in C^1(\R^n)$.
Let $x_\circ\in\partial\{v>\varphi\}$ be any regular point. 
By Theorem \ref{thm:FB}, there exists $r>0$ such that $\partial\{v>\varphi\}\cap B_r(x_\circ) \in  C^{\beta}$ for some~$\beta>1$.

Let us define $u=v-\varphi$
which solves
\begin{align}\label{obstacle-v}
\left\{\begin{array}{rcll}
\L u &=& f &  \quad \text{in}\quad  \{ u > 0 \} \\
u &\geq &0 &  \quad \text{in}\quad \R^n,
\end{array}\right.
\end{align}
where $f=-\L \varphi\in C^\infty(\R^n)$ (by Lemma~\ref{lem:Lu_2}). 
Note that $u\in C^1(\R^n)$ so that, 
for any $i\in\{1,\ldots,n\}$, we can differentiate \eqref{obstacle-v} to get
\begin{align}\label{obstacle-v'}
\left\lbrace\begin{aligned}
\L \big(\partial_i u\big) &= f_i & & \text{in } \{ u > 0 \}\cap B_r(x_\circ) \\
\partial_i u &= 0 & & \text{in } B_r(x_\circ)\setminus \{ u > 0 \}
\end{aligned}\right.
\end{align}
with $f_i:=\partial_i f \in C^\infty(\R^n)$. 
Suppose now, without loss of generality, 
that $e_n$ is normal to $\partial\{v>\varphi\}$ at $x_\circ$. 
Since at $x_\circ$ (and all points near $x_\circ$) we have that the free boundary is $C^\beta$ in $B_r(x_\circ)$, with $\beta>1$, it follows from the expansion in Theorem \ref{thm:FB}-\ref{it:i_fb} that 
\begin{align*}
\partial_n u \geq c_1 d^s
\qquad \text{in } \{ u > 0 \}\cap B_r(x_\circ)
\end{align*}
for some $c_1>0$.

We are therefore under the assumptions of Theorem \ref{thm:u1overu2} and we deduce that
\begin{align*}
\frac{\partial_i u}{\partial_n u} \in C^\beta\big(\overline{\{u>0\}}\cap B_r(x_\circ)\big),
\end{align*}
for any~$i\in\{1,\ldots,n-1\}$.

Now, notice that the normal vector $\nu(x)$ to the level set $\{u=t\}$ for $t>0$ and $u(x)=t$ is given by 
\begin{align*} 
\nu^i(x)= \frac{\partial_i u}{|\nabla u|}(x) = \frac{\partial_i u/\partial_n u}{\sqrt{\sum_{j=1}^{n-1}(\partial_j u/\partial_n u)^2 +1}},\qquad i=1,...,n.
\end{align*}
Therefore, denoting $\Omega=\{u>0\}$ we deduce that in $B_r(x_\circ)$ we have
\begin{align*}
\partial\Omega\in C^\beta\quad \Longrightarrow \quad \frac{\partial_i u}{\partial_n u} \in C^\beta \quad \Longrightarrow 
\quad \nu\in C^\beta \quad \Longrightarrow \quad \partial\Omega\in C^{\beta+1}.
\end{align*}

Bootstrapping this argument and recalling that $\Omega=\{v>\varphi\}$,  we find that $\partial\{v>\varphi\}\cap B_r(x_\circ)\in C^\infty$, as wanted.
\end{proof}

\section{Optimal regularity estimates}
\label{sec-optimal-reg}
\index{Obstacle problem!Optimal regularity}

After establishing the regularity of free boundaries (with the quantitative estimates from Propositions \ref{thm-main-ell} and \ref{lemm-main-ell-2}), we are now ready to prove the following optimal regularity estimates.
For this, we follow the approach from \cite{FRS,CFR}.

\begin{thm} \label{thm:OR}
Let $s\in(0,1)$ and $\L \in \LLh$.
Let $\varphi\in C^{2,\vartheta}_c(\R^n)$ with $\vartheta > \max\{0, 2s-1\}$, and let $v$ be the solution to the obstacle problem \eqref{obst-pb}-\eqref{obst-pb2}.
Then, 
\[\|v\|_{C^{1+s}(\R^n)} \leq C\|\varphi\|_{C^{2,\vartheta}(\R^n)},\]
with $C$ depending only on $n$, $s$, $\lambda$ and $\Lambda$.
\end{thm}

\begin{proof}
Let $u=v-\varphi$ and $\|\varphi\|_{C^{2,\vartheta}(\R^n)} = C_\circ$.
We first prove that 
\begin{equation}\label{iwjbgfi}
|\nabla u(x)| \leq CC_\circ|x-x_\circ|^s,
\end{equation}
for any free boundary point $x_\circ\in \partial\{u>0\}$.
Dividing by a constant if necessary, and up to a translation, we may assume $C_\circ\leq 1$ and $x_\circ=0$.

 Let $\eps_\circ$ and $R_\circ$ be given by Proposition~\ref{lemm-main-ell-2} (with $\kappa_\circ = \frac12$ and $\alpha = 0$), and let $\eta$ be the minimum between the constant  in Proposition~\ref{thm-main-ell} (with $\alpha = 0$ and these $\eps_\circ$ and $R_\circ$), and the one given in Proposition~\ref{lemm-main-ell-2} (with $\kappa_\circ  =\frac12$ and $\alpha  =0$).
Notice that $\eta$ depends only on $n$, $s$, $\lambda$, and $\Lambda$. 

Now, we consider 
\[w:=\eta\,u.\]
Then, $w$ satisfies (see Lemma~\ref{semiconvex})
\begin{equation}\label{parrot1}
w\geq0\quad \textrm{and} \quad D^2 w \ge - \eta{\rm{Id}}\quad \textrm{in}\quad \R^n,\quad \textrm{with}\quad 0\in \partial\{w>0\},
\end{equation}
\begin{equation}\label{parrot2}
\left\{
\begin{array}{rcll}
\L w&=&f&\quad \textrm{in}\quad \{w>0\} \\
\L w&\geq& f& \quad \textrm{in} \quad \R^n,
\end{array}
\right. \quad \textrm{with} \quad  |\nabla f|\leq \eta,
\end{equation}
and
\begin{equation}\label{parrot3}\|\nabla w\|_{L^\infty(\R^n)} \le \eta.
\end{equation}
We   want to apply Proposition \ref{thm-main-ell} appropriately  to get the desired estimate.

To prove it, consider the set of $r>0$ for which the following inequality does \emph{not} hold
\[\|\nabla w\|_{L^\infty(B_r)} \leq r^s.\]
If such set is empty, then \eqref{iwjbgfi} holds, so there is nothing to prove.
Otherwise, take $r_1>0$ its supremum (which exists, since $|\nabla u|$ is globally bounded), and observe that (recall that $\nabla u$ is continuous) 
\[\|\nabla w\|_{L^\infty(B_{r_1})} = r_1^s,\qquad \|\nabla w\|_{L^\infty(B_{r})} \leq r^s, \quad \text{for all}\quad r>r_1.\]
Hence, the rescaled function 
\[w_{r_1}(x) := \frac{w(r_1 x)}{r_1^{1+s}}\]
satisfies $\|\nabla w_{r_1}\|_{L^\infty(B_1)} = 1$ and
\[\|\nabla w_{r_1}\|_{L^\infty(B_R)}=
\frac{\|\nabla w\|_{L^\infty(B_{Rr_1})}}{r_1^s} \leq R^s \quad \text{for all}\quad R>1.\]
We can therefore apply Proposition \ref{thm-main-ell} to deduce that 
\[
\big\|w_{r_1}- \kappa(x\cdot e)_+^{1+s}\big\|_{{\rm Lip} (B_{R_\circ})} \le \eps_\circ,
\]
with $\kappa\geq 1-C  \eps\ge\frac12$ if $\eps$ is small (since $\|\nabla w_{r_1}\|_{L^\infty(B_1)} = 1$). 
Thus, by Proposition~\ref{lemm-main-ell-2} we obtain
\[|\nabla w_{r_1}| \leq C|x|^{s}\quad \textrm{in}\quad B_1,\]
with $C$ depending only on $n$, $s$, $\lambda$, and $\Lambda$. 
Rescaling back to $w$, this means that 
\[\|\nabla w\|_{L^\infty(B_r)} \leq Cr^s\quad \textrm{for all}\quad r\in(0,r_1).\]
Thus, by definition of $r_1$, the same bound holds for all $r>0$, and therefore \eqref{iwjbgfi} follows.

Finally, combining the inequality in \eqref{iwjbgfi} with interior regularity estimates, Theorem~\ref{prop-interior-linear} (exactly as in \ref{it:step2bdry} in the proof of Proposition~\ref{prop:globCsreg_loc}), the result follows.
\end{proof}

Actually, following the ideas of the previous proof, we can also establish that the constants in Theorem \ref{thm:FB}-\ref{it:i_fb} and \ref{it:ii_fb} do \emph{not} depend on the point~$x_\circ$.

\begin{prop} \label{prop-FB-uniform-const}
Let $s$, $\alpha$, $\L $, $\varphi$, $\vartheta$, $v$, and $\nu$ be as in Theorem \ref{thm:FB}, with  $\|\varphi\|_{C^{2,\vartheta}(\R^n)} \leq C_\circ$.
Then, for any regular free boundary point $x_\circ$ we have
\begin{equation}\label{parrot4}
\left|v(x)-\varphi(x)-a_{x_\circ}\big((x-x_\circ)\cdot \nu\big)_+^{1+s}\right| \leq CC_\circ|x-x_\circ|^{1+s+\gamma},
\end{equation}
with $a_{x_\circ}>0$, while for any singular point $x_\circ$ we have
\begin{equation}\label{parrot5}
\left|v(x)-\varphi(x)\right| \leq CC_\circ|x-x_\circ|^{1+s+\alpha}.
\end{equation}
The constants $\gamma>0$ and $C$ depend only on $n$, $s$, $\lambda$, and $\Lambda$.
\end{prop} 

\begin{proof}
Let $u=v-\varphi$.
Dividing by a constant if necessary, and up to a translation, we may assume $C_\circ\leq 1$ and $x_\circ=0$.
 Let us choose $\eps_\circ$, $R_\circ$, and   $\eta$ as in the proof of Theorem~\ref{thm:OR}, and define $w:=\eta \,u$.
Then, $w$ satisfies \eqref{parrot1}-\eqref{parrot2}-\eqref{parrot3}.

As before, in Theorem~\ref{thm:OR}, we have two cases.
Assume first the following inequality holds for all $r\in(0,1)$
\[\|\nabla w\|_{L^\infty(B_r)} \leq r^{s+\alpha}.\]
Then,
\[\left|  u(x)\right| \leq C|x|^{1+s+\alpha},\]
and in particular \eqref{parrot5} holds, so there is nothing to prove.

Otherwise, take 
\[r_1:=\sup\left\{r>0 : \|\nabla w\|_{L^\infty(B_r)} > r^{s+\alpha}\right\},\]
which satisfies $r_1\in(0,1)$. 
Then,
\[\|\nabla w\|_{L^\infty(B_{r_1})} = r_1^{s+\alpha},\qquad \|\nabla w\|_{L^\infty(B_{r})} \leq r^{s+\alpha} \quad \text{for all}\quad  r>r_1.\]
Hence, the rescaled function 
\[w_1(x) := \frac{w(r_1 x)}{r_1^{1+s+\alpha}}\]
satisfies
\[\|\nabla w_{r_1}\|_{L^\infty(B_1)} = 1,\qquad \|\nabla w_{r_1}\|_{L^\infty(B_R)} \leq R^{s+\alpha}\quad \text{for all}\quad  R>1.\]
We can therefore apply Proposition \ref{thm-main-ell} to deduce that 
\[
\big\|w_{r_1}- \kappa(x\cdot e)_+^{1+s}\big\|_{{\rm Lip} (B_{R_\circ})} \le \eps_\circ.
\]
Thus, by Proposition~\ref{lemm-main-ell-2}, we obtain
\[\|\nabla w_{r_1}/d^s\|_{C^\gamma(\overline{\{w_{r_1}>0\}}\cap B_{1/2})} \leq C,\]
and thus
\[\big\|w_{r_1}- \kappa(x\cdot e)_+^{1+s}\big\|_{{\rm Lip} (B_r)} \leq Cr^{s+\gamma}\quad \textrm{for}\quad r\in(0,1).\]

Rescaling back to $w$, this yields
\[\big\|w-\kappa r_1^\alpha(x\cdot e)_+^{1+s}\big\|_{L^\infty(B_r)} \leq Cr^{1+s+\gamma}\quad \textrm{for all}\quad r\in(0,r_1).\]
By definition of $r_1$ the same bound holds for all $r\in (0, 1)$ (assuming, without loss of generality, that $\alpha \ge \gamma$), and therefore \eqref{parrot4} follows.
\end{proof}

As a consequence of the previous quantitative estimate, we find the following:

\begin{cor} \label{cor-1+s+alpha-estimate-obst}
Let $s$, $\L $, $\varphi$, $\vartheta$, $v$, and $\nu$ be as in Theorem \ref{thm:FB}, with  $\|\varphi\|_{C^{2,\vartheta}(\R^n)} \leq C_\circ$, and let $u:=v-\varphi$.
Then, 
\[\|u/d^{1+s}\|_{C^\gamma(\R^n)} + \|\nabla u/d^s\|_{C^\gamma(\R^n)} \leq CC_\circ,\]
with $C$ and $\gamma>0$ depending only on $n$, $s$, $\lambda$, and $\Lambda$. 
\end{cor}

\begin{proof}
The result follows by combining the estimates from Proposition~\ref{prop-FB-uniform-const} with interior regularity estimates, exactly as in the proof of Proposition~\ref{prop:local_bdry}.
\end{proof}

In particular, the result in Corollary \ref{cor-1+s+alpha-estimate-obst} means that the limit
\[\lim_{\{u>0\}\ni x\to x_\circ} \frac{u}{d^{1+s}} = a_{x_\circ}\]
exists for any free boundary point $x_\circ\in \partial\{u>0\}$, and it is (H\"older) continuous in $x_\circ$.
Moreover, regular points can be characterized as those at which $a_{x_\circ}>0$, while degenerate points are those where $a_{x_\circ}=0$.

\section{Further results and open problems}

The main results we have proved in this chapter for the obstacle problem \eqref{obst-pb}-\eqref{obst-pb2} with operators of the form \eqref{obst-op1}-\eqref{obst-op2}-\eqref{obst-op3} may be summarized as follows: 
\begin{itemize}
\item Solutions are $C^{1+s}(\R^n)$, and this is optimal.

\item The free boundary splits into a set of \emph{regular points}, which is relatively open and $C^\infty$, and a closed set of \emph{degenerate points}.
\end{itemize}
We refer to Theorem \ref{thm:FB} for the dichotomy characterizing regular and degenerate points.

This gives already quite strong information about solutions and free boundaries (it is somehow the analogue of Caffarelli's celebrated result for the case of the Laplacian, \cite{Caf}), but at the same time it raises several interesting questions.

\subsection{Local estimates}

Most of the proofs we presented here are local in nature, in the sense that we only need the equation \eqref{obst-pb} to hold in a ball $B_1$ and then one can obtain local estimates in $B_{1/2}$, independently of the boundary data.

However, not all of them are local: the semiconvexity of solutions (see Lemma~\ref{semiconvex}) relies on a maximum principle type argument that really uses the fact that the problem is global.
The same proof can be carried out in bounded domains, but only when one assumes that the contact set $\{v=\varphi\}$ is at positive distance from the fixed boundary $\partial\Omega$.

To establish \emph{local} semiconvexity estimates for solutions of 
\begin{equation}\label{local-obst-4}
\min\{\L v,\,v-\varphi\}=0\quad \textrm{in}\quad B_1
\end{equation}
is much more delicate, and it was done only recently\footnote{In case of the fractional Laplacian $(-\Delta)^s$ one can use the extension property in order to prove such a semiconvexity estimate; see \cite{ACS,CV-semiconvexity,FJobst}.} by the second author, Torres-Latorre, and Weidner in \cite{RTW23}.
As a consequence, we also established optimal regularity estimates
\begin{equation}\label{local-estimate-ob}
\|v\|_{C^{1+s}(B_{1/2})} \leq C\left(\|\varphi\|_{C^\beta(B_1)} + \|v\|_{L^\infty(\R^n)} \right),
\end{equation}
for solutions of \eqref{local-obst-4}, with $\beta>1+2s$ and $K|_{\S^{n-1}} \in {\rm Lip}(\S^{n-1})$.
This raises the following:

\vspace{2mm}

\noindent\textbf{Open question 4.1}:
\textit{Does the estimate \eqref{local-estimate-ob} hold for any $\beta>1+s$?}

\vspace{2mm}

Related to this, we also have the following:

\vspace{2mm}

\noindent\textbf{Open question 4.2}:
\textit{Can one develop a similar regularity theory for operators with variable coefficients $\L(v,x)$?}

\vspace{2mm}

\subsection{More general kernels}

Many of the proofs we presented in this Chapter used quite strongly the fact that the operators $\L$ under consideration satisfy
\[0<\frac{\lambda}{|y|^{n+2s}} \leq K(y) \leq \frac{\lambda}{|y|^{n+2s}}.\] 

A natural question is then to understand whether one can develop the same regularity theory for obstacle problems driven by general (stable) operators $\L$ in the class $\GLh$ (recall Definition \ref{defi:Gh}):
\vspace{2mm}

\noindent\textbf{Open question 4.3}:
\textit{Does the regularity of solutions and free boundaries (Theorems \ref{thm:FB} and \ref{thm:OR}) hold for any stable operator $\L \in \GLh$?}

\vspace{2mm}

This question has been investigated in the recent work \cite{RW23}, where it has been answered positively under the extra assumption
\[K|_{\mathbb S^{n-1}} \in L^p(\mathbb S^{n-1})\quad \textrm{for some}\quad p>\frac{n}{2s}.\]
Moreover, thanks to the results in \cite{RW23}, it turns out that the previous question would follow for \textit{any} $\L\in\GLh$ if one can prove the uniqueness of positive $\L$-harmonic functions in cones, as stated next:

\vspace{2mm}

\noindent\textbf{Open question 4.4}:
\textit{Let $\L\in \GLh$ and $\Sigma\subset \R^n$ be any closed convex cone.
Let $w_1,w_2\in C(\R^n)$ be any \emph{positive} solutions of $\L w_i=0$ in $\Sigma^c$, with $w_i\equiv0$ in $\Sigma$.
Prove that $w_1\equiv \kappa w_2$ for some $\kappa>0$.}

\vspace{2mm}

Recall that for $\L \in \LLh$ we deduced this result (Proposition \ref{prop-uniq-positive-sols}) from the boundary Harnack inequality (Theorem \ref{thm-main}). 
However, for general operators $\L\in \GLh$ even the interior Harnack inequality may fail, and thus one would need a different approach in order to establish the uniqueness of positive solutions in cones.

\subsection{Degenerate points}

As in many other free boundary problems, an important and delicate question is to understand the set of singular (or, in this case, degenerate) points.
In case of the square root of the Laplacian $\sqrt{-\Delta}$, some of the best known results in this direction may be summarized as follows:

\begin{itemize}

\item If $\varphi$ is analytic, then the set of degenerate  free boundary points is contained in a countable union of $(n-1)$-dimensional $C^1$ manifolds, possibly except for a set of Hausdorff dimension at most $n-2$ \cite{FS,GP,CSV,FRS,FRS-obst,SY,FrSe}. 

\item There exist solutions whose free boundaries consist of {only} degenerate points.
Still, for \emph{almost every} solution, the set of degenerate points has Hausdorff dimension strictly less than $n-3$ \cite{FR3,FT}.

\end{itemize}

We refer to \cite{FS,BFR,GR,FR3} for similar results for the fractional Laplacian $(-\Delta)^s$ for any $s\in(0,1)$.

The proofs of \emph{all} these results are strongly based on the extension property of the fractional Laplacian, and in particular, that thanks to it one may use several monotonicity formulas for such operator.
These monotonicity formulas do not seem to exist for more general operators of the form \eqref{obst-op1}-\eqref{obst-op2}-\eqref{obst-op3}.
Thus, this raises then the following natural (and probably difficult) questions in this more general setting:

\vspace{2mm}

\noindent\textbf{Open question 4.5}:
\textit{Can one prove any structure or regularity result for the set of degenerate points?}

\vspace{2mm}

\noindent\textbf{Open question 4.6}:
\textit{Is it true that for almost every solution the set of degenerate points is small?}

\vspace{2mm}

As said above, we expect these questions to be quite challenging.

%% file: appendixA.tex
%
%
%

\chapter{Some properties of H\"older spaces}
\label{app.A}

In this appendix, we state some useful properties about H\"older spaces. Recall that, given $\alpha\in(0,1]$, the H\"older space $C^{0,\alpha}(\Omega)$ is the set of functions $u\in C(\Omega)$ such that
\[ [u]_{C^{0,\alpha}({\Omega})} := \sup_{\substack{x,y\in\Omega\\x\neq y}}\frac{\bigl|u(x)-u(y)\bigr|}{|x-y|^\alpha}<\infty.\]
The H\"older norm is
\[\|u\|_{C^{0,\alpha}(\Omega)}:=\|u\|_{L^\infty(\Omega)}+[u]_{C^{0,\alpha}({\Omega})}.\]
When $\alpha=1$, this is the usual space of Lipschitz functions.

More generally, we define: 
\begin{defi}[H\"older spaces]
\label{defi:Holderspace}\index{H\"older spaces@Holder spaces}
Given $k\in\mathbb N$ and $\alpha\in(0,1]$, the space $C^{k,\alpha}(\Omega)$ is the set of functions $u\in C^k(\Omega)$ such that the following norm is finite
\begin{align*}
\|u\|_{C^{k,\alpha}(\Omega)}& :=\sum_{j=1}^k\|D^ju\|_{L^\infty(\Omega)}+\sup_{\substack{x,y\in\Omega\\x\neq y}}\frac{\bigl|D^ku(x)-D^ku(y)\bigr|}{|x-y|^\alpha}\\
& =  \|u\|_{C^k(\Omega)} + [D^k u]_{C^{0,\alpha}({\Omega})}.
\end{align*}
Here, we can take any equivalent definition for the norms of $D^j u$, for example $|D^j u| = \sum_{|a| = j} |\partial^a u|$, where $a = (a_1,\dots,a_n)\in \N_0^n$ is a multi-index. When $\beta>0$ is \emph{not} an integer, we denote $C^\beta(\Omega):=C^{k,\alpha}(\Omega)$, where $\beta=k+\alpha$, with $k\in\mathbb N$, and $\alpha\in(0,1)$. We similarly denote $\|u\|_{C^\beta(\Omega)} := \|u\|_{C^{k,\alpha}({\Omega})} $ and $[u]_{C^\beta(\Omega)} := [D^k u]_{C^\alpha(\Omega)} $.

Finally, we denote $C^{k,\alpha}_{\rm loc}(\Omega)$ the set of functions $u\in C^k(\Omega)$ such that $u\in C^{k,\alpha}(K)$ for any $K\ssubset \Omega$. 
\end{defi}

\section{Some useful lemmas} Let us start by presenting some lemmas on properties of H\"older functions. 
We refer to \cite[Appendix A]{FR4} for the corresponding proofs of Lemmas \ref{it:H7}, and \ref{it:H8}.

The first one is on second order incremental quotients.

\begin{lem}
\label{it:H7} 
Assume that $\beta\in(0,2]$, $\beta\neq 1$, $\|u\|_{L^\infty(B_1)}\leq C_{\circ}$, and
\[
\sup_{\substack{h\in B_1\\ x\in\overline{B_{1-|h|}}}}\frac{\bigl|u(x+h)+u(x-h)-2u(x)\bigr|}{|h|^{\beta}}\leq C_{\circ}.
\]
Then, $u\in C^{\beta}(\overline{B_1})$ ($u\in C^{1, 1}(\overline{B_1})$ if $\beta = 2$) and $\|u\|_{C^{\beta}(\overline{B_1})}\leq CC_{\circ}$, with $C$ depending only on $n$ and $\alpha$. 
\end{lem}

And the second one is on norms of first order incremental quotients:

\begin{lem}
\label{it:H8} Assume that $\alpha\in (0,1]$, $\beta > 0$, $\|u\|_{L^\infty(B_1)}\leq C_{\circ}$,  and that for every $h\in B_1$ we have
\[
\left\|\frac{u(x+h)-u(x)}{|h|^\alpha}\right\|_{C^\beta(B_{1-|h|})}\leq C_{\circ},
\]
with $C_{\circ}$ independent of $h$.
Assume in addition that $\alpha+\beta\notin\N$.
Then, $u\in C^{\alpha+\beta}(\overline{B_1})$ and $\|u\|_{C^{\alpha+\beta}(\overline{B_1})}\leq CC_{\circ}$, with $C$ depending only on $n$, $\alpha$, and $\beta$.
\end{lem}

Let us now prove a useful characterization of H\"older spaces in Lipschitz domains. 
We define first what we mean by \emph{Lipschitz domain} (cf. Definition~\ref{defi:varrho} for more regular domains).

\begin{defi}[Lipschitz domains]
\label{defi:varrho_Lip}\index{Lipschitz domains}
Let $\Omega\subset \R^n$. We say that $\Omega$ is a Lipschitz domain with Lipschitz constant $L>0$ and $C^{0,1}$-radius $\varrho$  if, for any $x_\circ\in \partial\Omega$, we have
\[
\varrho^{-1} \mathcal{R} \left(\partial\Omega-x_\circ\right)\cap \mathcal{C}_L = \left\{ 
\begin{array}{l}
(x', x_n)\in \mathcal{C}_L : x_n = \varphi(x'),\\
\ \text{with}\  \|\nabla \varphi\|_{L^\infty(B_1')}\le L
\end{array}\right\},
\]
for some rotation $\mathcal{R}$ and some $\varphi\in C^{0,1}(B_1')$, where $\mathcal{C}_L:= B_1'\times[-L-1, L+1]$ and $B_1'\subset\R^{n-1}$ is the unit ball.  
\end{defi}

\begin{lem}
\label{lem:Lipholder}
Let $\alpha\in (0, 1]$, and let $\Omega$ be any bounded Lipschitz domain with $C^{0,1}$-radius $\varrho$ and Lipschitz constant $L$. 
Let us suppose that, for some $u\in L^\infty(\Omega)$,
\[
\osc_{B_r(x)}   u\le C_\circ r^\alpha\quad\text{for any}\quad B_{2r}(x)\subset \Omega.
\]
Then $u\in C^{0,\alpha}(\overline{\Omega})$ and 
\[
[u]_{C^{0,\alpha}(\Omega)} \le C C_\circ
\]
for some $C$ depending only on $n$, $L$, $\varrho$, ${\rm diam}(\Omega)$, and $\alpha$. 
\end{lem}
\begin{proof}
Let $x_1, x_2\in \Omega$. Let us start by assuming that $r:= |x_1-x_2|\le \frac{\varrho}{4}$. If $\dist(x_1, \partial\Omega) >\varrho$, we are done, since $B_{2r}(x_1)\subset \Omega$ and $|u(x_1) - u(x_2)|\le \osc_{B_r(x_1)}\le C_\circ r^\alpha$ by assumption. Hence, we restrict ourselves to the case where $x_1$ and $x_2$ are both at distance at most $\varrho$ from the boundary.

After a  rotation and a $\varrho$-scaling, we may assume that  $x_1, x_2\in  \{(x', x_n)\in \mathcal{C}_L : x_n > \varphi(x')\}\subset \Omega$ for some $\varphi\in C^{0, 1}(B_1')$ with $\varphi(0) = 0$ and Lipschitz constant $L$, where $\mathcal{C}_L := B_1'\times(-L-1, L+1)$. After a Lipschitz diffeomorphism we   assume furthermore that $\varphi(x') \equiv 0$ in $B_1'$, and so $\Omega \cap \mathcal{C}_L = \{x_n > 0\}\cap \mathcal{C}_L$ (the hypothesis on the oscillation of $u$ is the same after a covering argument, with a possibly larger $C_\circ$ by a multiplicative factor depending only on the diffeomorphism). We finally assume, after a translation, that $x_1 = (0, \delta_1)$ and $x_2 = (x_2', \delta_2)$, where $|x_2'|\le \frac14$ and $\delta_1, \delta_2 > 0$. We refer to Figure~\ref{fig:07} for a graphical representation of the setting described. 

\begin{figure}
\centering
\makebox[\textwidth][c]{\includegraphics[scale = 1]{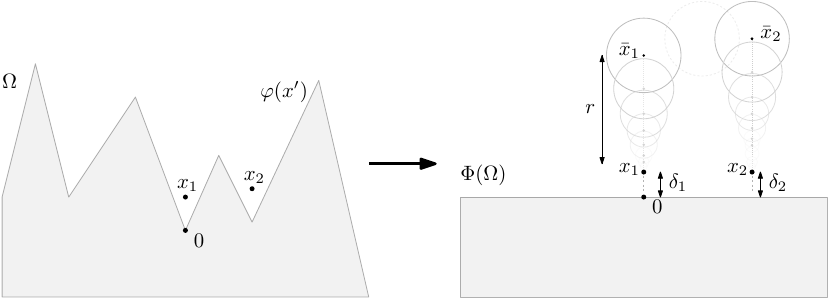}}
\caption{\label{fig:07} Graphical representation of the setting of the proof of Lemma~\ref{lem:Lipholder}. We flatten the domain by means of a Lipschitz diffeomorphism.}
\end{figure}

Notice that $B_{\delta_1}(x_1)\subset \Omega$. If $\delta_1 > r$ and $\delta_2 > r$ we are done by assumption, since there will be a path joining them inside $\Omega$ that can be covered by finitely many balls of radius $\frac{r}{2}$ entirely contained in $\Omega$. Let us suppose then that $\delta_1 \le r$, and let us define $\bar x_1 := (0, \delta_1+r)$.

We consider the path $x_{k} = (1-2^{-k})x_1 + 2^{-k}\bar x_1$, that satisfies $|x_k - x_{k+1}| =  2^{-k-1}r$ and  $(x_k)_n = \delta_1 + 2^{-k}r$. This implies that $x_{k+1}\in B_{2^{-k-1} r}(x_k)$ and $B_{2^{-k}r}(x_k)\subset \Omega$ and we can therefore apply the hypothesis of the lemma to deduce
\[
|u(x_k) - u(x_{k+1})|\le C C_\circ 2^{-(k+1)\alpha} r^\alpha, 
\]
for some constant $C$ that depends on $L$ (due to the Lipschitz diffeomorphism). By summing a geometric series we obtain that 
\[
|u(x_1) - u(\bar x_1)|\le \sum_{k \ge 0} |u(x_k) - u(x_{k+1})|\le C C_\circ r^\alpha.
\]
If $\delta_2  \le r$, we can do the same with $\bar x_2 := (x_2', \delta_2+r)$ to get in this case
\[
|u(x_2) - u(\bar x_2)|\le C C_\circ r^\alpha.
\]
Now, the points $\bar x_1$ and $\bar x_2$ satisfy that $(\bar x_1)_n > r$ and $(\bar x_2)_n > r$, and so $|u(\bar x_1) - u (\bar x_2)|\le C C_\circ r^\alpha$ as before. In all, we have 
\[
|u(x_1) - u(x_2)|\le |u(x_1) - u(\bar x_1)|+|u(\bar x_1) - u(\bar x_2)|+|u(x_2) - u(\bar x_2)| \le C C_\circ r^\alpha,
\]
as we wanted to see. 

Finally, let us suppose that $r > \frac{\varrho}{4}$. By assumption, the set 
\[
\Omega_{\varrho/4} := \{x\in \Omega : \dist(x, \partial\Omega > \varrho/4\}
\]
satisfies that any two points $y_1, y_2\in \Omega_{\varrho/4}$ can be connected by a path fully contained in $\Omega_{\varrho/4}$ of length depending only on $\varrho$, $L$, and ${\rm diam}(\Omega)$. If $\dist(x_i, \partial\Omega) \le \varrho/4$, let $\bar x_i$ be defined as above: that is, $\dist(\bar x_i, \partial\Omega) \ge \varrho/4$ and 
\[
|u(\bar x_i)  - u(x_i)|\le C C_\circ \varrho^\alpha. 
\]
We can then connect $\bar x_1$ and $\bar x_2$ by finitely many balls of radius $\varrho/8$ and center in $\Omega_{\varrho/4}$, each of which satisfies the hypotheses of the lemma, to deduce that 
\[
|u(\bar x_1)  - u(\bar x_2)|\le C C_\circ \varrho^\alpha\le C C_\circ r^\alpha,
\]
where $C$ depends only on the number of balls connecting $\bar x_1$ and $\bar x_2$. Since there is a finite path between the two points, $C$ depends only on $\varrho$, $L$, and ${\rm diam}(\Omega)$. In particular, we have
\[
|u(x_1) - u(x_2)|\le |u(x_1) - u(\bar x_1)|+|u(\bar x_1) - u(\bar x_2)|+|u(x_2) - u(\bar x_2)| \le C C_\circ r^\alpha,
\]
as we wanted to see. 
\end{proof}

\section{Incremental  quotients}\index{Incremental quotients}\label{sec:incremental_quotients}
Given $f:\R^n\to \R$, $h\in \R^n$, and $m\in \N$, the $m$-th order incremental or difference quotient $D_h^m f$ is
\[
D_h^m f(x) = 
\left\{
\begin{array}{ll}
\frac{1}{|h|}(f(x+h)-f(x))& \quad\text{if}\quad m = 1,\\
D_h (D_h^{m-1} f(x))& \quad\text{if}\quad m\ge 2.
\end{array}
\right.
\]
\begin{lem}
\label{it:H10}
Let $f\in C^{m-1}(\R^n)$ be such that $D_h^m f$ is constant for any $h\in \R^n$ (maybe depending on $h$). Then, $f$ is a polynomial of degree $m$.
\end{lem}
\begin{proof}
For $m = 1$, the assumption implies $D_h f(x) = D_h f(x-h)$, and we obtain $f(x+h) +f(x-h) = 2f(x)$ for all $x, h\in \R^n$. This is the equality case in Jensen's inequality, which for continuous functions holds only when $f$ is affine. By induction, if $D_h^m f$ is constant, then $D_h^{m-1} f$ is affine  and $D_h^{m-1} f'$ is constant, so $f'$ is a polynomial of degree $m-1$. 
\end{proof}

We also have the following basic property on the boundedness of incremental quotients for H\"older functions:

\begin{lem}
\label{it:H7_gen}
Let $\alpha\in(0,1)$, $k \in \N$, $\Omega$   convex, and let $u\in  C^{k-1}(\overline\Omega)$ be such that $[u]_{C^{k-1,\alpha}(\overline\Omega)}\leq C_{\circ}$. Then,
\[
\sup_{\substack{h\in \R^n \\ \dist(x, \partial \Omega) \ge k|h|}} \frac{|D_h^k u|}{|h|^{\alpha}}\leq C C_{\circ}.
\]
For some $C$ depending only on $n$, $\alpha$, and $k$. 
\end{lem}

\begin{proof}
For $k = 1$ this is the definition of $C^{0,\alpha}$. Let us show the case $k= 2$. By the mean value theorem 
\[
u(x+h) - u(x) = h \cdot\nabla u(x+t h)
\]
for some $t\in (0, 1)$, and therefore, 
\begin{equation}
\label{eq:tayl_form}
|u(x+h)-u(x) - h \cdot\nabla u(x)|\le |h||\nabla u(x+th) - \nabla u(x)|\le C_\circ|h|^{1+\alpha}. 
\end{equation}
Taking $-h$ instead of $h$, and adding up we obtain 
\begin{equation}
\label{eq:proc_as_in}
|u(x+h)+u(x-h)-2u(x)|\le 2C_\circ |h|^{1+\alpha},
\end{equation} 
which is what we wanted (up to replacing $x$ by $x+h$).   

For the general case, we proceed similarly. Indeed, by the Taylor theorem we know that, denoting $h = \hat{h}|h|$ for some $\hat{h}\in \mathbb{S}^{n-1}$,
\[
\begin{split}
u(x_\circ+h)-u(x_\circ)   -& \dots  -\frac{|h|^{k-2}}{(k-2)!} \partial_{\hat{h}}^{(k-2)} u(x_\circ) -\frac{|h|^{k-1}}{(k-1)!} \partial_{\hat{h}}^{(k-1)} u(x_\circ)= \\
& = \frac{|h|^{k-1}}{(k-1)!} \partial_{\hat{h}}^{(k-1)} u(x_\circ+th)-\frac{|h|^{k-1}}{(k-1)!} \partial_{\hat{h}}^{(k-1)} u(x_\circ)  
\end{split}
\]
for some $t\in (0, 1)$, and therefore,
\[
|u(x_\circ+h)-p(x_\circ, h)|\le C |h|^{k-1+\alpha},
\]
for any $h$, where $p(x_\circ, y)$ is the Taylor polynomial (in $y$) of $u$ at $x_\circ$,  of order $k-1$. In particular, since $D_h^k p(x_\circ, 0) = 0$ (where the increments are in $y$, and since $p$ is a polynomial of order $k-1$), we deduce the desired result from the triangle inequality. 
\end{proof}

If we denote by $\delta^2_h u(x)$ the second order centered increments, 
\[
  \delta^2_h u(x) = \frac{u(x+h)+u(x-h)}{2} - u(x),
\]
we then have the following technical lemma: 

\begin{lem}
\label{lem:A_imp_2}
Let $\alpha\in (0, 1]$. Then, we have the following:
\begin{enumerate}[leftmargin=*, label=(\roman*)]
\item \label{it:A:1} If $u\in C_{\rm loc}^{1,\alpha}(\R^n)$ with $[u]_{C^{1,\alpha}(\R^n)}\le 1$ then 
\[
\begin{split}
|u(x+h) - u(x) - u(h) + u(0)|& \le \min\{|h||x|^\alpha, |h|^\alpha|x|\}\\ & \le |h|^{t+\alpha(1-t)}|x|^{\alpha t + 1-t}
\end{split}
\]
and
\[
\begin{split}
\big| \delta^2_h u(x) -  \delta^2_h u(0)\big|& \le \min\{2|h|^{1+\alpha},  |h|^\alpha|x|, |h||x|^{\alpha}\}\\
& \le 2 |h|^{(1+\alpha) (1-t)+t\alpha}|x|^t
\end{split}
\]
for all $x, h\in \R^n$ and $t\in [0,1]$.

\item \label{it:A:3}  If $u\in C_{\rm loc}^{2,\alpha}(\R^n)$ with $[u]_{C^{2,\alpha}(\R^n)}\le 1$ then 
\[
\big| \delta^2_h u(x) -  \delta^2_h u(0)\big|\le \min\left\{ |h|^{2}|x|^\alpha, |h|^{1+\alpha}|x|\right\}\le |h|^{2t+(1+\alpha)(1-t)}|x|^{t\alpha + 1 - t}
\]
for all $x, h\in \R^n$ and $t\in [0, 1]$. 
\end{enumerate}
\end{lem}

\begin{proof}
Fixing $x, h\in \R^n$, let us denote 
\[
w(y) := u(x+y) - u(y).
\]
Then, the left-hand side in the first expression of \ref{it:A:1} is given by
\[
|w(h) - w(0)| \le |h||\nabla w(\bar h)|\quad\text{where}\quad \bar h = \bar t h \quad\text{for some} \ t\in [0, 1], 
\]
by the mean value formula. Notice, moreover, that 
\[
|\nabla w(\bar h)| = |\nabla u(x+\bar h) - \nabla u(\bar h)|\le |x|^\alpha
\]
 since $[u]_{C^{1, \alpha}(\R^n)}\le 1$. By the symmetry of the roles of $x$ and $h$, this gives the first inequality  in \ref{it:A:1}. We also use that, for any $a, b\ge 0$ and $t\in [0,1]$, 
 \begin{equation}
 \label{eq:minabineq}
 \min\{a, b\}\le a^t b^{1-t}.
 \end{equation}
 
 For the second one, we know that $| \delta_h^2 u(x)| \le |h|^{1+\alpha}$ and $| \delta_h^2 u(0)| \le |h|^{1+\alpha}$   by Lemma~\ref{it:H7_gen}, \eqref{eq:proc_as_in}.  On the other hand, notice also that by the triangle inequality 
 \[
 \begin{split}
 \big| \delta^2_h u(x) -  \delta^2_h u(0)\big|& \le \frac12 |u(x+h) - u(x) - u(h) + u(0)|\\
 & \quad + \frac12 |u(x-h) - u(x) - u(-h) + u(0)|,
  \end{split}
 \]
 so that the result follows by applying the first inequality (and \eqref{eq:minabineq}). 
 
Let us now show \ref{it:A:3}. Observe that 
 \[
 w(y) = \int_0^1 x\cdot \nabla u(y+tx)\, dt. 
 \]
 Then, 
 \[
 \big| \delta^2_h u(x) -  \delta^2_h u(0)\big| = \left| \delta^2_h w (0) \right| = \left|\int_0^1 x\cdot ( \delta^2_h \nabla u)(tx)\, dt\right|. 
 \]
 Now, since $\nabla u \in C^{1,\alpha}$, we have that $\big|  ( \delta^2_h \nabla u)(z)\big| \le  |h|^{1+\alpha}$ for every $z, h\in \R^n$, and thus  
 \[
 \big| \delta^2_h u(x) -  \delta^2_h u(0)\big| \le |x||h|^{1+\alpha}. 
 \]
On the other hand, by the mean value theorem for second derivatives we know that
\[
 \delta_h^2 w(0) = \frac12 h\cdot D^2 w(\bar h) h,\quad\text{where}\quad \bar h =  \bar t h,\quad\text{for some} \  \bar t\in (-1, 1). 
\]
Recalling that $w(y) = u(x+y)-u(y)$ and using that $[u]_{C^{2,\alpha}(\R^n)}\le 1$, we find
\[
\big|h\cdot \big(D^2 u(x+\bar h)-D^2 u(\bar h)\big) h\big|\le |h|^2|x|^\alpha,
\]
and therefore  the result follows (again, also using \eqref{eq:minabineq}).
 \end{proof}

\section{Interpolation inequalities}\index{Interpolation inequalities}

In the following, we present a well known interpolation inequality in H\"older spaces (see \cite[Lemma~6.35]{GT}). We give a simple proof for the convenience of the reader. 

In this section, we say that $\Omega$ is a $C^\gamma$ domain if $\Omega$ is a $C^{k,\alpha}$ domain for some $k\in \N$, $\alpha\in(0, 1]$, and $\gamma = k+\alpha$.

\begin{prop}\label{it:H9} Let $\gamma >  \beta > 0$, and let $\Omega$ be either a bounded $C^{\max\{\gamma, 1\}}$ domain, or $\Omega = \R^n$. Then, for any $\eps > 0$ and $u\in C^\gamma(\Omega)$ we have 
\[
\|u\|_{C^\beta(\Omega)}\le   C_\eps \|u\|_{L^\infty(\Omega)}+\eps[u]_{C^\gamma(\Omega)},
\]
for some $C_\eps$ depending only on $n$, $\eps$, $\gamma$, $\beta$, and $\Omega$. 
\end{prop}

\begin{proof}
We split the proof into five steps.

\begin{steps}
\item \label{step:1p}
Let us assume $\Omega = \R^n$. We first establish the case $\gamma\leq 1$, i.e.,
\begin{equation}
\label{eq:gt2}
[  u]_{C^\beta(\R^n)}\le C_\eps\|u\|_{L^\infty(\R^n)} + \eps [  u]_{C^\gamma(\R^n)},\qquad 0< \beta < \gamma \le 1. 
\end{equation}
Indeed, we have
\[
\frac{|u(x) - u(y)|}{|x-y|^\beta}\le |u(x)-u(y)|^{1-\frac{\beta}{\gamma}} \left(\frac{|u(x) - u(y)|}{|x-y|^\gamma}\right)^{\frac{\beta}{\gamma}}\quad\text{for any}\quad x, y\in \R^n,
\]
and therefore
\[
[u]_{C^\beta(\R^n)}\le C \|u\|_{L^\infty(\R^n)}^{1-\frac{\beta}{\gamma}} [u]_{C^\gamma(\R^n)}^{\frac{\beta}{\gamma}}\le C_\eps \|u\|_{L^\infty(\R^n)} + \eps [u]_{C^\gamma(\R^n)}.
\]
We have used the inequality  $t^{1-\theta} \le C_\eps + \eps t$ for any fixed $\theta \in (0, 1)$, where $C_\eps$ depends only on $\theta$ and $\eps>0$. 

\item
We next show the case $\beta = 1$ and $\gamma \in (1, 2]$, that is, 
\begin{equation}
\label{eq:gt1}
\|\nabla u\|_{L^\infty(\R^n)}\le C_\eps\|u\|_{L^\infty(\R^n)} + \eps [\nabla u]_{C^\alpha(\R^n)},\qquad 0 < \alpha \le 1.
\end{equation}

For any $x\in \R^n$, $j\in \{1, \dots, n\}$, and $\rho \in  (0, 1)$, by the mean value theorem there exists $\bar \rho \in (0, \rho)$  such that
\[
\partial_j u (x+\bar\rho e_j) = \frac{u(x+\rho e_j) - u(x)}{\rho}. 
\]
Then, we have
\[
|\partial_j u(x)|\le |\partial_j  u(x+\bar \rho e_j) - \partial_j u(x) | + \frac{| u(x+\rho e_j) - u(x)|}{\rho},
\]
and thus
\[
|\partial_j u(x)| \le \rho^\alpha [\partial_j u]_{C^\alpha(\R^n)} + 2\rho^{-1}\| u\|_{L^\infty(\R^n)}. 
\]
Summing over $1\le j\le n$, we obtain 
\[
|\nabla u(x)|\le C \rho^\alpha[\nabla u]_{C^\alpha(\R^n)} + C \rho^{-1}\|u\|_{L^\infty(\R^n)}. 
\]
By taking $\rho$ small enough such that $C \rho^\alpha < \eps$, we have shown \eqref{eq:gt1}.

\item
 We next show, for any $k\in \N$, the case $\beta=k$ and $\gamma=k+1$, 
  \begin{equation}
\label{eq:gt3}
\|u\|_{C^k(\R^n)}\le C_\eps\|u\|_{L^\infty(\R^n)} + \eps \|u\|_{C^{k+1}(\R^n)}.
\end{equation}
We proceed by induction. 
The case $k = 0$ is trivial, and $k = 1$ is a consequence of \eqref{eq:gt1} with $\alpha = 1$. Let us now suppose that it holds for some $k\in \N$, and we prove it for $k+1$. 
From \eqref{eq:gt1} we have
\[
\|u\|_{C^{k+1}(\R^n)}\le C_{\eps_0}\|u\|_{C^{k} (\R^n)} + \eps_0 \|u\|_{C^{k+2}(\R^n)}.
\] 
By induction  hypothesis, we deduce
\[
\|u\|_{C^{k+1}(\R^n)}\le C_{\eps_0}C_{\eps_1}\|u\|_{L^\infty(\R^n)} +C_{\eps_0}\eps_1 \|u\|_{C^{k+1} (\R^n)} + \eps_0 \|u\|_{C^{k+2}(\R^n)}.
\]
Choosing $\eps_0 = \frac{\eps}{2}$ and $\eps_1$ small enough so that $1-C_{\eps_0}\eps_1 \ge \frac12$, we get \eqref{eq:gt3}.  

\item
Finally, let us show how to use \eqref{eq:gt2}-\eqref{eq:gt1}-\eqref{eq:gt3} to get the desired result.
We know from \eqref{eq:gt1} that for any $\alpha \in (0, 1]$, 
\[
\|u \|_{C^k(\R^n)} \le C_{\eps_0} \|u\|_{C^{k-1}(\R^n)} + \eps_0 [u]_{C^{k+\alpha}(\R^n)}.
\]
Using \eqref{eq:gt3}, 
\[
\|u \|_{C^k(\R^n)} \le C_{\eps_0}C_{\eps_1} \|u\|_{L^\infty(\R^n)} +C_{\eps_0} \eps_1 \|u \|_{C^k(\R^n)} + \eps_0 [u]_{C^{k+\alpha}(\R^n)},
\]
so that again  choosing $\eps_0 = \frac{\eps}{2}$ and $\eps_1$ small enough so that $1-C_{\eps_0}\eps_1 \ge \frac12$, we get (also using $\|u\|_{C^{k'}(\R^n)}\le \|u\|_{C^k(\R^n)}$ for $k' \le k$)
  \begin{equation}
\label{eq:gt4}
\|u\|_{C^k(\R^n)}\le C_\eps\|u\|_{L^\infty(\R^n)} + \eps [u]_{C^{\gamma}(\R^n)}\quad\text{for any}\quad k <\gamma.
\end{equation}
On the other hand, if we denote $\lfloor \beta \rfloor = k$, thanks to \eqref{eq:gt2} we have  
\[
[u]_{C^{\beta}(\R^n)}\le C_{\eps_0} \|u\|_{C^k(\R^n)} + \eps_0 [u]_{C^{\min\{\gamma, k+1\}}(\R^n)}
\]
 which arguing as above (by means of \eqref{eq:gt3} or \eqref{eq:gt4}) completes the proof.  
 \item 
 To finish, we use  the extension theorem to prove our desired result in a given  bounded $C^\gamma$ domain $\Omega$.
 
Observe that \ref{step:1p} follows in the exact same way as before, restricting ourselves to $\Omega$. This proves the analogue of \eqref{eq:gt2} in $\Omega$. In order to obtain the analogue of \eqref{eq:gt1}, we use the extension theorem as follows. 
  
 We know that, given  $u\in C^\gamma(\Omega)$, there exists  $\bar u\in C^{\gamma}(\R^n)$, called the $C^{\gamma}$ extension of $u$ towards $\R^n$, which satisfies $u = \bar u$ in $\Omega$ and $\|\bar u\|_{C^\nu(\R^n)}\le C \|u\|_{C^\nu(\Omega)}$ for any $\nu \le \gamma$, for some $C$ depending only on $n$, $\gamma$, and $\Omega$ (see Lemma~\ref{lem:extensionthm}). In particular, we can apply \eqref{eq:gt1} to $\bar u$ and obtain  
\[
\begin{split}
\|\nabla u\|_{L^\infty(\Omega)}\le \|\nabla \bar u\|_{L^\infty(\R^n)}& \le C_\eps\|\bar u\|_{L^\infty(\R^n)} + \eps [\nabla \bar u]_{C^\alpha(\R^n)}\\
& \le CC_\eps\| u\|_{L^\infty(\Omega)} + C\eps \|u\|_{C^{1, \alpha}(\Omega)},
\end{split}
\]
 for $0 < \alpha \le 1$. That is, 
\[
(1-C\eps) \|\nabla u\|_{L^\infty(\Omega)}\le C(C_\eps+\eps)\| u\|_{L^\infty(\Omega)} + C\eps [\nabla u]_{C^{\alpha}(\Omega)}.
\]
Up to making $\eps$ smaller, we recover the analogue of \eqref{eq:gt1} in a domain $\Omega$. The proof now finishes in the exact same way as before.\qedhere
 \end{steps}
\end{proof}

In the previous proof, we have used the following version of the extension theorem, to prove the result in bounded domains. 

\begin{lem}[Extension Theorem]
\label{lem:extensionthm}
Let $\Omega\subset \R^n$ be a bounded $C^\gamma$ domain (with $\gamma \ge 1$) and let $u\in C^\gamma(\overline{\Omega})$. Then, there exists $\bar u\in C^\gamma_c(\R^n)$ such that $\bar u = u$ in $\Omega$ and 
\[
\|\bar u\|_{C^{\nu}(\R^n)}\le C \|u\|_{C^{\nu}(\Omega)} \quad\text{for all}\quad 0\le \nu \le \gamma,
\]
for some $C$ depending only on $n$, $\gamma$, and $\Omega$. 
\end{lem}

\begin{proof}
This result follows from the proof of \cite[Lemma 6.37]{GT}. 
\end{proof}

We can also use the extension theorem to derive, from the previous estimate in Proposition~\ref{it:H9}, a multiplicative version of the interpolation inequality:

\begin{lem}\label{lem:interp_mult}
Let $\gamma >\beta > 0$,  let $r \in (0, \infty]$, and let $u\in C^{\gamma}(B_r)$. Then 
\[
[u]_{C^\beta(B_r)}\le C\|u\|^{1-\frac{\beta}{\gamma}}_{L^\infty(B_r)}[u]^{\frac{\beta}{\gamma}}_{C^\gamma(B_r)} +Cr^{-\beta} \|u\|_{L^\infty(B_r)},
\]
for some $C$ depending only on $n$, $\gamma$, and $\beta$. 
\end{lem}

\begin{proof}
Let us prove first the result for $r = \infty$. In order to do that, we apply Proposition~\ref{it:H9} with $\eps = 1$ and $\Omega = \R^n$ to the function $u_\rho(x) = u(\rho x)$, to obtain
\[
\begin{split}
\rho^\beta [u]_{C^\beta(\R^n)} = [u_\rho]_{C^\beta(\R^n)}& \le   C \|u_\rho\|_{L^\infty(\R^n)}+ [u_\rho]_{C^\gamma(\R^n)} \\ & = C \|u\|_{L^\infty(\R^n)}+\rho^\gamma [u]_{C^\gamma(\R^n)}.
\end{split}
\]
By choosing $\rho = \|u\|^{\frac{1}{\gamma}}_{L^\infty(\R^n)}[u]^{-\frac{1}{\gamma}}_{C^\gamma(\R^n)}$ we get the desired result in the case $r = \infty$. 

Let us now do the case $r \in (0, \infty)$. Notice first that it is enough to prove it for $r = 1$. The result for general $r$ then follows by scaling, by applying the result for $r = 1$ to the function $u(rx)$. 

Let $\bar u$ the $C^\gamma$ extension of $u$ to $\R^n$ given by Lemma~\ref{lem:extensionthm} with $\Omega = B_1$. Now, since $\bar u$ is globally $C^\gamma$, by the case $r = \infty$ and the extension theorem we already know that 
\begin{equation}
\label{eq:endingproof}
[u]_{C^\beta(B_1)} = [\bar u]_{C^\beta(\R^n)} \le C \|\bar u\|_{L^\infty(\R^n)}^{1-\frac{\beta}{\gamma}} [\bar u]_{C^\gamma(\R^n)}^{\frac{\beta}{\gamma}} \le C \|u\|_{L^\infty(B_1)}^{1-\frac{\beta}{\gamma}} \| u \|_{C^\gamma(B_1)}^{\frac{\beta}{\gamma}}.
\end{equation}
Thanks to the interpolation result in Proposition~\ref{it:H9} we also have 
\[
\| u \|_{C^\gamma(B_1)} \le C \| u \|_{L^\infty(B_1)} + [ u ]_{C^\gamma(B_1)},
\]
which, used in \eqref{eq:endingproof} together with the fact that $(a+ b)^p \le a^p + b^p$ for any $a, b > 0$ and $p\in (0, 1)$, gives the desired result. 
\end{proof}

As a consequence of the interpolation inequality, we also have the following bound on the semi-norm of the product of functions: 

\begin{prop}
\label{prop:A_imp}
  Let $u, v\in C^{k, \alpha}(B_1)$ for some $k\in \N_0$ and $\alpha\in(0, 1)$. Then $uv\in C^{k,\alpha}(B_1)$ with 
\[
[uv]_{C^{k+\alpha}(B_1)} \le C \left([u]_{C^{k+\alpha}(B_1)} \|v\|_{C^k(B_1)} + \|u\|_{L^\infty(B_1)}\|v\|_{C^{k+\alpha}(B_1)}\right) 
\]
for some $C$ depending only on $n$, $k$, and $\alpha$. 
\end{prop}
\begin{proof}
Let us suppose first that $k = 0$. Let now $\Omega\subset \R^n$ be any open set, so that for any $x, y\in \Omega$
 \[
 \frac{|(u v)(x) - (u v)(y)|}{|x-y|^\alpha} \le\frac{|u(x) - u (y)|}{|x-y|^\alpha}|v(y)| +  |u(x)|\frac{|v(x) - v(y)|}{|x-y|^\alpha}, 
 \]
 and therefore 
  \begin{equation}
  \label{eq:APP_In1}
  [u v]_{C^\alpha(\Omega)} \le [u]_{C^\alpha(\Omega)}\|v\|_{L^\infty(\Omega)}+\|u\|_{L^\infty(\Omega)}[v]_{C^\alpha(\Omega)}.
  \end{equation}
  Taking $\Omega = B_1$ we get the case $k = 0$. Let us now suppose $k \ge 1$. By the chain rule, we have
\[
\partial^a(uv) = \sum_{b \le a} c_{ab}\partial^b u \partial^{a-b} v,
\]
 where $a,b \in (\N_0)^n$ are multi-indices, and $c_{ab}$ are constants. 
 
 Thus, by the triangle inequality and the case $k = 0$ we have
\[
\begin{split}
[u v]_{C^{k+\alpha}(B_1)}& \le C\sum_{ {|a| = k,~~   b\le a}} [\partial^b u \partial^{a-b} v]_{C^\alpha(B_1)}\\
& \hspace{-1cm}\le C\sum_{i = 0}^k\left([D^i u ]_{C^\alpha(B_1)}\|D^{k-i} v\|_{L^\infty(B_1)}+\|D^i u\|_{L^\infty(B_1)}[D^{k-i}v]_{C^\alpha(B_1)}\right),
\end{split}
\]
for some $C$ depending only on $k$. Using now the interpolation inequality, Proposition~\ref{it:H9},  
\[
\|u\|_{C^{j+\alpha}(B_1)} \le [u]_{C^{k+\alpha}(B_1)} +  C \|u\|_{L^\infty(B_1)},
\]
  for $j\in \N_0$,  and $j\le  k $, we  obtain the desired result. 
\end{proof}

%% file: appendixB.tex
%
%
%

\chapter{Construction of barriers}
\label{app.B}

\addtocontents{toc}{\protect\setcounter{tocdepth}{1}}

 In this appendix, we construct some useful barriers for nonlocal operators in different kinds of domains.

An important ingredient in these constructions will be the following:

\begin{lem}[\cite{Lieberman}]
\label{lem:distance}\index{Regularized distance}
Let $\Omega\subset \R^n$ be any open set, and let $d_\Omega(x)={\rm dist}(x,\Omega^c)$.
Then, there exists a function $\dr_\Omega \in C^\infty(\Omega)\cap {\rm Lip}(\R^n)$ satisfying
\begin{equation}
\label{eq:compddr}
d_\Omega\leq \dr_\Omega \leq C d_\Omega\quad \textrm{in}\quad \R^n,
\end{equation}
and 
\[\big|D^k\dr_\Omega\big| \leq C_kd_\Omega^{1-k}\quad \textrm{in}\quad \Omega\]
for all $k\in\mathbb N$, with $C$ and $C_k$ depending only on $n$ and $k$.

If, in addition, $\Omega$ is a bounded $C^{\beta}$ domain with $\beta > 1$ and $\beta\notin\N$, then $\dr_\Omega\in C^\beta(\overline{\Omega})$ and  
\[\big|D^k\dr_\Omega\big| \leq C_{k,\Omega}d_\Omega^{\beta-k}\quad \textrm{in}\quad \Omega\]
for all $k\in\mathbb N$, $k>\beta$.

 Finally, if $\Omega$ is a  bounded  $C^1$ domain, then $\dr_\Omega\in C^1(\overline{\Omega})$ and there exists a modulus of continuity $\omega$ depending only on $\Omega$ such that
\[
\begin{split}
 \big|D^2 \dr_\Omega\big|&\le \omega(d_\Omega) d_\Omega^{-1}\quad\text{in}\quad \Omega
\end{split}
\]
\end{lem}
 
When there is no possible confusion about what the domain $\Omega$ is, we will simply denote 
\[d:=d_\Omega \qquad \textrm{and}\qquad \dr := \dr_\Omega.\]

Notice that, in case of bounded $C^\beta$ domains with $\beta>1$ and $\beta\notin\mathbb N$, a regularized distance can be constructed by taking the solution of $\Delta \dr=1$ in~$\Omega$ with $\dr = 0$ on $\partial\Omega$, for example.

For arbitrary open sets $\Omega\subset\R^n$, the function constructed in \cite{Lieberman} is given by the equation
\[\dr(x) = \int_{B_1} d\big(x-{\textstyle\frac12}\dr(x)z\big)\phi(z)dz, \]
where $\phi\in C^\infty_c(B_1)$ is a nonnegative function with $\int_{B_1}\phi=1$.

\section{One-dimensional and radial barriers}

We first construct  one-dimensional homogeneous subsolutions.

\begin{lem}\label{lem-subsol1D}
Let $s\in(0,1)$, $e\in \mathbb{S}^{n-1}$, and   $\L\in \LLL$ --- that is, of the form \eqref{eq:Lu1}-\eqref{eq:Kint1}-\eqref{eq:compdef}. Then, there exists $\beta_\circ\in (0,2s)\cap (0,1)$ such that, for any $\beta\in (0,2s)$ with $\beta\geq\beta_\circ$, we have
 \[\L|x\cdot e|^\beta\leq -c|x\cdot e|^{\beta-2s}\quad\textrm{in}\quad \R^n\]
in the viscosity sense. The constants $c>0$ and $\beta_\circ$ depend only on $n$, $s$, $\lambda$, and $\Lambda$.
\end{lem}

\begin{proof}
We will use the extremal operators (recall \eqref{eq:extremal_def} and \eqref{eq:MMpm})
\[\Mp w:= \sup_{\L\in \LLL} \big\{-\L w\big\},\qquad
 \Mm w:= \inf_{\L\in \LLL} \big\{-\L w\big\},\]
which satisfy (since the class $\LLL$ is scale invariant of order $2s$)
\[\Mpm w_r(x) = r^{2s}\big(\Mpm w\big)(rx).\]
In particular, if $w$ is positively homogeneous of degree $\beta$, then 
\[\Mpm w(x) = \lambda^{\beta-2s}\big(\Mpm\big) w(x/\lambda).\]
 Notice also that the class $\LLL$ is rotation invariance, so it suffices to prove the result for $e=e_n$.

Thanks to the above discussion, we know that 
 \[\Mm |x_n|^\beta = c_\beta |x_n|^{\beta-2s}\quad \textrm{in}\quad \{|x_n|>0\}\]
for $\beta\in(0,2s)$, where $c_\beta\in \R$ is given by 
 \[c_\beta = \left.\Mm|x_n|^\beta\right|_{x=e_n}.\]
Moreover, it is easy to see that $c_\beta\to+\infty$ as $\beta\to 2s$.
On the other hand, when $s>\frac12$, since the function $|x_n|_+^\beta$ is convex for $\beta\geq1$, we have $c_\beta>0$ for $\beta\geq1$.
Since $c_\beta$ is continuous with respect to $\beta$, for any $s\in(0,1)$ there is $\beta_\circ\in (0,2s)\cap (0,1)$ such that $c_\beta\geq c>0$ for all $\beta\geq \beta_\circ$.

By definition of $\Mm$, this implies that for any operator $\L\in \LLL$ and all $\beta\geq \beta_\circ$ we have 
\[\L|x_n|^\beta\leq -\Mm |x_n|^\beta \leq  -c|x_n|^{\beta-2s} \quad \textrm{in}\quad \{|x_n|>0\},\]
with $c>0$. Therefore, since on $\{x_n=0\}$ the function has a minimum,  
\[\L|x_n|^\beta \le -\Mm |x_n|^\beta \leq - c |x_n|^{\beta-2s}\quad \textrm{in}\quad \R^n\]
in the viscosity sense, as we wanted to see.
\end{proof}

As a consequence, we next show the following.

\begin{lem}\label{lem-supersol1D}
Let $s\in(0,1)$, $e\in \mathbb S^{n-1}$, and   $\L\in \LLL$ --- that is, of the form \eqref{eq:Lu1}-\eqref{eq:Kint1}-\eqref{eq:compdef}.
Then, there exists $\theta >0$ such that 
\[\phi(x):= \exp\big(-|x\cdot e|^{1-\theta}\big)\]
satisfies 
 \[\L\phi\geq -C\quad\textrm{in}\quad \R^n\]
in the viscosity sense.
The constants $C$ and $\theta$ depend only on $n$, $s$, $\lambda$, and~$\Lambda$. 
\end{lem}

\begin{proof}
Thanks to Lemma \ref{lem-subsol1D}, for some $\theta > 0$ we have
\[\L|x\cdot e|^{ 1-\theta}\leq -c|x\cdot e|^{ 1-\theta-2s}\leq0 \quad \textrm{in}\quad \R^n,\]
with $c>0$.
In particular, since the difference between $\phi(x)$ and $-|x\cdot e|^{1-\theta}$ is $C^{2s+\delta}$ (for $1-\theta > s$),   the function $\phi$ satisfies $\L\phi \geq -C$ in $\R^n$ (by Lemma~\ref{lem:Lu}), as wanted.
\end{proof}

\begin{rem}
The previous results, Lemmas~\ref{lem-subsol1D} and \ref{lem-supersol1D}, work as well with essentially the same proof for any operator $\L \in \GL$, provided that one gives the corresponding definitions of $\M_{\GL}^\pm$ and viscosity super- and subsolution in this context. 
\end{rem}

We will also need the following radial subsolution. 
Recall that $\mathfrak G_s(\lambda,\Lambda)$ is the class of operators given by Definition \ref{defi:G}.

\begin{lem}\label{lem-radialsubsol}
Let $s\in(0,1)$, $\L\in \mathfrak G_s(\lambda,\Lambda)$.
Then, there exists $\beta\in (s,2s)$ such that 
 \[\L(1-|x|^2)_+^\beta\leq -1\quad\textrm{in}\quad B_1\setminus B_{1-\eta},\]
where   $\beta$  and $\eta>0$ depend only on $n$, $s$, $\lambda$, and $\Lambda$.
\end{lem}

\begin{proof}
Let $\eta>0$ to be fixed later, and let $x_\circ\in B_1\setminus B_{1-\eta}$.
After a rotation we may assume $x_\circ=(1-r)e_n$, with $r\in (0,\eta)$.

Consider the rescaled function
\[w(x):=\frac{(1-|x_\circ+rx|^2)_+^\beta}{r^{\beta}} = \big(2(x_n-1)-r|x-e_n|^2\big)^\beta_+,\]
which satisfies 
\[\|w\|_{C^2(B_{1/2})} \leq C \quad \textrm{and}\quad w(x)\geq (x_n-1)_+^\beta \quad \textrm{in}\quad \{x_n-1>r|x-e_n|^2\}.\]
Thus, we have, for any $\tilde \L \in \GL$ with kernel $K$,
\[\tilde \L w(0) \leq C - \int_{\{x_n-1>r|x-e_n|^2\}} (y_n-1)_+^\beta K(y)dy.\]
Now, as $\eta\to0$ (recall $r\in (0,\eta)$) we have
\[\int_{\{x_n-1>r|x-e_n|^2\}} (y_n-1)_+^\beta K(y)dy \longrightarrow 
\int_{\{x_n>1\}} (y_n-1)_+^\beta K(y)dy,\]
and by the ellipticity and the symmetry condition on $\tilde \L$, one can show that this converges to $+\infty$ as $\beta\to 2s$ (depending only on $n$, $s$, $\lambda$, and $\Lambda$).
Thus,
\[\int_{\{x_n-1>r|x-e_n|^2\}} (y_n-1)_+^\beta K(y)dy \longrightarrow +\infty\quad\textrm{as}\quad \eta\to0,\quad \beta\to 2s,\]
and hence
\[\tilde \L w(0)  \leq - 1 \]
for some $\beta\in(s,2s)$ and for all $r\in(0,\eta)$, provided that $\eta>0$ is small enough.
Rescaling back, and taking $\tilde \L = \L_r$ the rescaled version of $\L$ (see Remark~\ref{rem:Scale_invariance}), we deduce that 
\[\L (1-|x|^2)_+^\beta|_{x=x_\circ} = r^{\beta-2s} \L_r w(0) \leq - 1,\]
and since $x_\circ\in B_1\setminus B_{1-\eta}$ was arbitrary, we are done.
\end{proof}

Finally, we will also need the following three results, which hold for stable operators $\L$. The first one says that, for one dimensional functions, in order to compute $\L u$ it is enough to compute its fractional Laplacian:

\begin{lem}
 \label{lem:L1d_2}
Let $s\in (0, 1)$, and let $\L\in \GLh$. Let $u \in L^\infty_{2s-\eps}(\R^n)\cap C^{2s+\eps}(B_r(x))$ for some $\eps > 0$ and $r > 0$ be a one-dimensional function, $u(x) = \bar u(x\cdot\be)$ for some $\bar u:\R\to \R$ and $\be\in \S^{n-1}$. Then
\[
\L u (x)  =  c_\be \big((-\Delta)^s_\R \bar u\big)(x\cdot \be)
\]
for some $c_\be$ depending only on $n$, $s$, $\L$, and $\be$, such that 
\[
c_- \le c_\be \le c_+,
\]
where $c_-, c_+ > 0$ are constants depending only on $n$, $s$, $\lambda$, and $\Lambda$, but independent of $\be$. 
\end{lem}
\begin{proof}
We compute $\L u (x)$ by using polar coordinates 
\[K(dy) = \frac{dr}{r^{1+2s}}\,K(d\theta),\]
with $r>0$ and $\theta\in \mathbb S^{n-1}$, where we used the homogeneity of the kernel $K$.
Then, we get
\[
\begin{split}
\L  u & (x)  = \frac{1}{2} \int_{\S^{n-1}}\int_{0}^{\infty} \big(2u(x) - u(x+r\theta) - u(x-r\theta)\big) \frac{dr}{r^{1+2s}}{K}(d\theta)  \\
& = \frac{1}{4} \int_{\S^{n-1}}\int_{-\infty}^{\infty} \big(2u(x) - u(x+r\theta) - u(x-r\theta)\big) \frac{dr}{|r|^{1+2s}}{K}(d\theta)  \\
& = \frac{1}{4} \int_{\S^{n-1}}\int_{-\infty}^{\infty} \big(2\bar u(x\cdot \be) - \bar u((x+r\theta)\cdot\be ) - \bar u((x-r\theta)\cdot \be )\big) \frac{dr}{r^{1+2s}}{K}(d\theta)  \\
& = c_s\int_{\S^{n-1}} \big({(-\Delta)^s_\R}\bar u\big)(x\cdot \be+r \theta\cdot \be)\big|_{r = 0}{K}(d\theta)   \\
& = c_s\big({(-\Delta)^s_\R}\bar u\big)(x\cdot \be)\int_{\S^{n-1}} |\theta\cdot\be|^{2s} {K}(d\theta)   = c_s \A(\be)({(-\Delta)^s_\R}\bar u)(x\cdot \be),
\end{split}
\]
for some $c_\be = c_s\A(\be)> 0$, as we wanted to see. 
\end{proof}

As a consequence we obtain the following: 

\begin{lem}
\label{lem:1dseps}
Let $s\in (0, 1)$, and let $\L\in \GLh$, and let $\eps\in (0, s)$. 
Then,  for any $\be\in \S^{n-1}$ the function 
\[u(x) = (x\cdot \be)^{s+\eps}_+\]
 satisfies 
\[
\left\{
\begin{array}{rcll}
\L u & = & - c_\be (x\cdot \be)^{\eps - s}_+& \quad\text{in}\quad \{x\cdot \be > 0 \},\\
u & = & 0  & \quad\text{in}\quad \{x\cdot \be \le 0 \},
\end{array}
\right.
\]
for some $c_\be > 0$ depending only on $n$, $s$, $\eps$, $\be$, $\lambda$, and $\Lambda$, such that 
\[
c_- \le c_\be \le c_+,
\]
for $c_-, c_+ > 0$ depending only on $n$, $s$, $\eps$, $\lambda$, and $\Lambda$.
\end{lem}

\begin{figure}
\centering
\makebox[\textwidth][c]{\includegraphics[scale = 1]{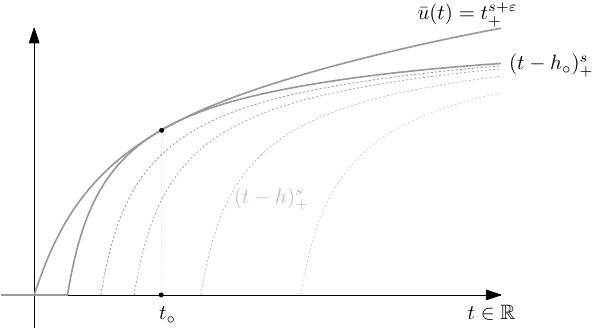}}
\caption{\label{fig:10} We slide the function $t^s_+$ from the right until it touches $\bar u$.}
\end{figure}

\begin{proof}
By Lemma~\ref{lem:L1d_2}, if we denote $\bar u (t) = t^{s+\eps}_+$,
\[
\L u (x) = c_\be'{(-\Delta)^s_\R} \bar u (x\cdot \be)\quad\text{for}\quad x\in \R^n. 
\]

The fractional Laplacian is a $2s$-homogeneous operator (see Lemma~\ref{lem:invariances_s}) so we get 
\[
{(-\Delta)^s_\R} \bar u (t) = c_* t^{\eps - s}_+\quad\text{for}\quad t \in \R,
\]
where $c_* = {(-\Delta)^s_\R} \bar u(1) \in \R$ by Lemma~\ref{lem:laplu_s}. We just need to check the sign of $c_*$. In order to do that, we ``slide'' the function $t^s_+$ from the right until we touch $\bar u$ (see Figure~\ref{fig:10}). Namely, consider $(t-h)^s_+$ for $h > 0$, which for $h > 0$ large enough is below $\bar u$. We now make $h$ small until they touch at some point $t_\circ > 0$, at $h_\circ > 0$. Then, 
\[
c_* (t_\circ^{\eps - s})_+ = {(-\Delta)^s_\R} \bar u (t_\circ) < {(-\Delta)^s_\R} [(t-h_\circ)^s_+](t_\circ) = 0,
\]
where we have used that ${(-\Delta)^s_\R} \left(\bar u (t) - (t-h_\circ)^s_+\right)  (t_\circ) < 0$ because $t_\circ$ is a strict minimum, and that ${(-\Delta)^s_\R} t_+^s=0$ for $t>0$; see Proposition \ref{prop:onedbarrier}. 
Thus, $c_* < 0$ and the result follows. 
\end{proof}

 As well as:

\begin{lem}
\label{lem:1dseps2}
Let $s\in (0, 1)$, and let $\L\in \GLh$, and let $\eps\in (0, s)$. 
Then,  for any $\be\in \S^{n-1}$ the function 
\[u(x) = (x\cdot \be)^{s-\eps}_+\]
 satisfies 
\[
\left\{
\begin{array}{rcll}
\L u & = &  c_\be (x\cdot \be)^{-\eps - s}_+& \quad\text{in}\quad \{x\cdot \be > 0 \},\\
u & = & 0  & \quad\text{in}\quad \{x\cdot \be \le 0 \},
\end{array}
\right.
\]
for some $c_\be > 0$ depending only on $n$, $s$, $\eps$, $\be$, $\lambda$, and $\Lambda$, such that 
\[
c_- \le c_\be \le c_+,
\]
for $c_-, c_+ > 0$ depending only on $n$, $s$, $\eps$, $\lambda$, and $\Lambda$.
\end{lem}
\begin{proof}
As in the proof of Lemma~\ref{lem:1dseps} it is enough to check that $c_* = {(-\Delta)^s_\R} \bar u(1) >0$, where $\bar u (t) = t^{s-\eps}_+$. To do that, we now  ``slide'' the function $t^s_+$ from the left until we touch $\bar u$: we consider $(t+h)^s_+$ for $h > 0$, which for $h > 0$ large enough is above $\bar u$, and make it small until they touch at some point $t_\circ > 0$, for $h = h_\circ > 0$. Then, 
\[
c_* (t_\circ^{-\eps - s})_+ = {(-\Delta)^s_\R} \bar u (t_\circ) > {(-\Delta)^s_\R} [(t+h_\circ)^s_+](t_\circ) = 0,
\]
using that ${(-\Delta)^s_\R} \left(\bar u (t) - (t+h_\circ)^s_+\right)  (t_\circ)> 0$ since $t_\circ$ is a strict maximum, and that ${(-\Delta)^s_\R} t_+^s=0$ for $t>0$; see Proposition \ref{prop:onedbarrier}. 
Thus, $c_* > 0$. 
\end{proof}

\section{Barriers for stable operators}

We next construct appropriate barriers in $C^1$ and $C^{1,\alpha}$ domains.
Recall that, given any $C^{1,\alpha}$ domain $\Omega\subset \R^n$,  we say that $\varrho$ is a $C^{1,\alpha}$-radius for $\Omega$  if, at scale $\varrho$,  the domain looks like the graph of a $C^{1,\alpha}$ function with $C^{1,\alpha}$ norm bounded by 1; see Definition~\ref{defi:varrho}. 
Notice also that, by \cite[Theorem 1.3]{Lieberman}, we can take the regularized distance $\dr$ so that $\|\dr\|_{C^{1,\alpha}(\overline\Omega)}$ depends only on $n$ and $\varrho$.

\subsection*{An approximate solution}

The key result of this section is the following.

\begin{prop}\label{lem11}
Let $s\in (0, 1)$, and let $\L\in \GLh$.
 Let $\alpha\in(0,s)$, let $\Omega$ be any $C^{1,\alpha}$ domain with $C^{1,\alpha}$-radius $\varrho_\circ > 0$, and  let $ \dr = \dr_\Omega$ be given by Definition \ref{defi:distance}.
Then,  $\dr^s\in H^s(\R^n)$ and
 
\begin{enumerate}[leftmargin=*, label=(\roman*)]
\item \label{it:lem11_i} We have
\[
|\L( \dr ^s)| \leq C\dr^{s\alpha-s}\qquad \textrm{in}\quad \Omega.
\]

\item \label{it:lem11_ii} If, in addition, the kernel satisfies \eqref{eq:boundabovekernel}, then
\begin{equation}\label{eq-psi-10}
|\L( \dr^s)| \leq C\dr^{\alpha-s}\qquad \textrm{in}\quad \Omega.
\end{equation}
\end{enumerate}

\noindent The constant $C$ depends only on $n$, $s$, $\alpha$, $\varrho_\circ$, $\lambda$, and $\Lambda$.
\end{prop}

For this, we need a few technical lemmas.
The first one reads as follows.

\begin{lem}\label{lem13}
Let $\Omega$ be any $C^{1,\alpha}$ domain with $C^{1,\alpha}$-radius $\varrho_\circ>0$.
Then, for each $x_\circ\in \Omega$ we have
\[\left| \dr(x_\circ+y)-\bigl( \dr(x_\circ)+\nabla  \dr(x_\circ)\cdot y\bigr)_+\right|\leq C|y|^{1+\alpha}\qquad \textrm{for}\quad  y\in \R^n.\]
The constant $C$ depends only on $n$ and $\varrho_\circ$.
\end{lem}

\begin{proof}
Let us consider $\tilde  \dr$, a $C^{1,\alpha}(\R^n)$ extension of $ \dr|_\Omega$ satisfying $\tilde \dr\leq 0$ in $\R^n\setminus\Omega$ (see, e.g., Lemma~\ref{lem:extensionthm}).
Then, since $\tilde  \dr\in C^{1,\alpha}(\R^n)$ we   have
\[\left|\tilde \dr(x)- \dr(x_\circ)-\nabla \dr(x_\circ)\cdot (x-x_\circ)\right|\leq C|x-x_\circ|^{1+\alpha}\]
in all of $\R^n$. Here we used $\tilde  \dr(x_\circ)= \dr(x_\circ)$ and $\nabla \tilde  \dr(x_\circ)=\nabla \dr(x_\circ)$.

Now, using that $|a_+-b_+|\leq |a-b|$, combined with $(\tilde \dr)_+= \dr$, we find
\[\left| \dr(x)-\bigl( \dr(x_\circ)+\nabla  \dr(x_\circ)\cdot (x-x_\circ)\bigr)_+\right|\leq C|x-x_\circ|^{1+\alpha}\]
for all $x\in \R^n$, and  the lemma follows.
\end{proof}

\begin{lem}\label{lem13_2}
Let $\Omega$ be any bounded $C^{1}$ domain. Then, there exists a modulus of continuity $\omega$ depending only on $\Omega$ such that, for each $x_\circ\in \Omega$ we have
\[\left| \dr(x_\circ+y)-\bigl( \dr(x_\circ)+\nabla  \dr(x_\circ)\cdot y\bigr)_+\right|\leq \omega(|y|)|y|\qquad \textrm{for}\quad  y\in \R^n.\]
\end{lem}

\begin{proof}
The proof is analogous to the one of Lemma~\ref{lem13}. We consider $\tilde  \dr$, a $C^{1}(\R^n)$ extension of $ \dr|_\Omega$ satisfying $\tilde \dr\leq 0$ in $\R^n\setminus\Omega$.
By Taylor's theorem (since $\tilde\dr\in C^1(\R^n)$) we know that for any $x\in \R^n$,
\[
\tilde\dr(x) - \tilde\dr(x_\circ) = \nabla \tilde \dr\left(\lambda x + (1-\lambda) x_\circ)\right)\cdot (x-x_\circ)
\]
for some $\lambda\in [0, 1]$ (that depends on $x$). Hence, once again using that $\tilde\dr\in C^1(\R^n)$,
\[
\begin{split}
& \hspace{-1cm}\left|\tilde \dr(x)- \tilde \dr(x_\circ)-\nabla \tilde \dr(x_\circ)\cdot (x-x_\circ)\right| = \\
& \qquad =  \left|\left( \nabla \tilde \dr\left(\lambda x + (1-\lambda) x_\circ)\right)-\nabla \tilde \dr(x_\circ)\right)\cdot (x-x_\circ)\right| \\
&\qquad  \leq \omega(\lambda|x-x_\circ|)|x-x_\circ|\le \omega(|x-x_\circ|)|x-x_\circ|
\end{split}
\]
in all of $\R^n$, for some modulus $\omega$. Using now that $\tilde  \dr(x_\circ)= \dr(x_\circ)$, $\nabla \tilde  \dr(x_\circ)=\nabla \dr(x_\circ)$, $|a_+-b_+|\leq |a-b|$, and $(\tilde \dr)_+= \dr$, gives the desired result.
\end{proof}

The second one reads as follows.

\begin{lem}\label{lem15}
Let $\Omega$ be any Lipschitz  domain with Lipschitz constant $L$ and $C^{0,1}$-radius $\varrho_\circ > 0$ (see Definition~\ref{defi:varrho_Lip}). Let $x_\circ\in \Omega$, $\rho=d_\Omega(x_\circ)$,  $\gamma>-1$, and $\beta\neq\gamma$.
Then,
\[\int_{\Omega \cap (B_1\setminus B_{\rho/2})}d_\Omega^\gamma(x_\circ+y)\frac{dy}{|y|^{n+\beta}}\leq C\bigl(1+\rho^{\gamma-\beta}\bigr).\]
The constant $C$ depends only on $\gamma$, $\beta$, $\varrho_\circ$, and $L$.
\end{lem}

\begin{proof}
Let us assume $x_\circ = 0$. 
Notice that, since $\Omega$ is Lipschitz,  there is $\kappa_*>0$ such that for any $t\in(0,\kappa_*]$ the level set $\{d_\Omega=t\}$ is Lipschitz as well.
Since
\begin{equation}\label{kappa11}
\int_{(B_1\setminus B_{\rho/2})\cap \{d_\Omega\geq\kappa_*\}}d_\Omega^\gamma(y)\frac{dy}{|y|^{n+\beta}}\leq C,
\end{equation}
 we just have to bound the same integral in the set $\{0<d_\Omega<\kappa_*\}$.
Here we used that $B_r\cap \{d_\Omega\geq\kappa_*\}=\varnothing$ if $r\leq \kappa_*-2\rho$, which follows from the fact that $d_\Omega(0)=2\rho$.

We will use the following estimate for $t\in (0,\kappa_*)$:
\[\mathcal{H}^{n-1}\bigl(\{d_\Omega=t\}\cap (B_{2^{-k+1}}\setminus B_{2^{-k}})\bigr)\leq C(2^{-k})^{n-1},\]
which follows, for example, from the fact that $\{d=t\}$ is Lipschitz.
Note also that $\{d_\Omega=t\}\cap B_r=\varnothing$ if $t>r+2\rho$.

Let $M\geq0$ be such that $2^{-M}\leq \rho/2\leq 2^{-M+1}$.
Then, using the coarea formula,
\begin{equation}\label{kappa22}\begin{split}
\int_{(B_1\setminus B_{\rho/2})\cap \{0<d_\Omega<\kappa_*\}}& d_\Omega^\gamma(y)\frac{dy}{|y|^{n+\beta}}\leq \\
&\hspace{-1cm}\leq\ \sum_{k=1}^M \frac{1}{2^{-k(n+\beta)}}\int_{(B_{2^{-k+1}}\setminus B_{2^{-k}})\cap \{0 < d_\Omega<C2^{-k}\}}
    d^\gamma(y)|\nabla d_\Omega(y)|\,dy \\
&\hspace{-1cm}\leq\  \sum_{k=1}^M \frac{1}{2^{-k(n+\beta)}}\int_0^{C2^{-k}}t^\gamma dt\int_{(B_{2^{-k+1}}\setminus B_{2^{-k}})\cap \{d_\Omega=t\}}
    d\mathcal{H}^{n-1}(y)  \\
&\hspace{-1cm}\leq\  C\sum_{k=1}^M \frac{(2^{-k})^{\gamma+1}2^{-k(n-1)}}{2^{-k(n+\beta)}}= C\sum_{k=1}^M 2^{k(\beta-\gamma)}= C(1+\rho^{\gamma-\beta}).
\end{split}\end{equation}
Here we used that $\gamma\neq \beta$; in case $\gamma=\beta$ we would get $C(1+|\log \rho|)$.

Combining \eqref{kappa11} and \eqref{kappa22}, the lemma follows.
\end{proof}

We will also need the following:

\begin{lem}
\label{lem:distHs}
Let $\Omega\subset \R^n$ be any bounded Lipschitz domain, and let $\varphi\in C(\R^n)$ be such that
\[
\left\{
\begin{array}{rcll}
\varphi  & = & 0 & \quad\text{in}\quad \R^n\setminus \Omega,\\
|\varphi | & \le & Cd^\alpha & \quad\text{in}\quad  \Omega,\\
|\nabla \varphi | & \le & Cd^{\alpha-1} & \quad\text{in}\quad \Omega,\\
\end{array}
\right.
\]
with 
\[\alpha> s-{\textstyle \frac12}.\]
Then, $\varphi \in H^s(\R^n)$. 
\end{lem}

\begin{proof}
Let us denote $d(x) = d_\Omega(x)$. 
We will bound $[\varphi]_{H^s(\R^n)}$  (observe that $\varphi  \in L^2(\R^n)$ since $|\varphi |^2 \le d_\Omega^{2\alpha}$, $\Omega$ is Lipschitz, and $2\alpha > -1$). 

We separate the computation into two terms:
\[
[\varphi]_{H^s(\R^n)} = \iint_{\R^n\times\R^n} \frac{|\varphi(x) - \varphi(y)|^2}{|x-y|^{n+2s}}\, dx\, dy = I + II,
\]
where 
\[
I := 2 \int_\Omega\int_{\R^n\setminus \Omega}\frac{|\varphi(x)|^2}{|x-y|^{n+2s}}\, dy\, dx
\]
and 
\[
II := \int_\Omega\int_\Omega \frac{|\varphi(x) - \varphi(y)|^2}{|x-y|^{n+2s}}\, dx\, dy.
\]

For the first term, observe that $|x-y|\ge d_\Omega(x) +d_\Omega(y) \ge d_\Omega(x)$ for $x\in \Omega$ and $y\in \R^n\setminus \Omega$, 
and therefore
\[
I\le 2 \int_\Omega d^{2\alpha}(x) \int_{|x-y|\ge d(x)}\frac{dy}{|x-y|^{n+2s}}\,  dx \le c\int_\Omega d^{2\alpha-2s}(x)\, dx< \infty,
\]
where we have used that $\int_{B_r^c} |z|^{-n-2s}\, dz = cr^{-2s}$. The boundedness of the last term follows because $\Omega$ is Lipschitz and $2\alpha -2s > -1$, by assumption. 

For the second term, we can split the integral in $y$ into two parts: one integrating in $B_{d(x)/2}(x)$ and the other integrating in the rest. In the first integral, we can use that 
\[
|\varphi(x)-\varphi(y)|^2\le \|\nabla \varphi\|^2_{B_{d(x)/2}(x)}|x-y|^2 \le  Cd^{2\alpha-2}(x)|x-y|^2,
\]
to deduce
\[
\begin{split}
II& \le \int_\Omega\int_{B_{d(x)/2}(x)} d^{2\alpha-2}(x)|x-y|^{2-n-2s}\, dy\, dx\\
& \qquad +2\int_\Omega\int_{\Omega\setminus B_{d(x)/2}(x)} (\varphi^2(x) + \varphi^2(y))|x-y|^{-n-2s}\, dy\, dx.
\end{split}
\]
Now, since $ \int_{B_{d(x)/2}} |z|^{2-n-2s}\, dz\le C d^{2-2s}(x)$, the first term above is again of the form $\int_\Omega d^{2\alpha-2s}(x)\, dx < \infty$. On the other hand, we also have 
\[
\int_\Omega\int_{\Omega\setminus B_{d(x)/2}(x)}\varphi^2(x) |x-y|^{-n-2s}\, dy\, dx <\infty
\]
arguing as above, so that 
\[
II \le C + 2\int_\Omega\int_{\Omega\setminus B_{d(x)/2}(x)}d^{2\alpha}(y) |x-y|^{-n-2s}\, dy\, dx.
\]
Observe, now, that we can assume $\Omega\subset B_1$ (after scaling) and we can bound the last term using Lemma~\ref{lem15}:
\[
II \le C\left(1+\int_\Omega d^{2\alpha-2s}(x)\, dx\right)< \infty. 
\]
This completes the proof.
\end{proof}

We now give the:

\begin{proof}[Proof of Proposition \ref{lem11}]
The fact that $\dr^s\in H^s(\R^n)$ follows from Lemma \ref{lem:distHs}.
Let $x_\circ\in \Omega$ and $\rho=d_\Omega(x_\circ)$, which is comparable to $\dr(x_\circ)$.

Notice that when $\rho\geq \rho_\circ>0$ then $ \dr^s$ is smooth in a neighborhood of $x_\circ$, and hence $(\L \dr^s)(x_\circ)$ is bounded by a constant depending only on $\rho_\circ$.
Thus, we may assume that $\rho\in(0,\rho_\circ)$.

We start by proving \ref{it:lem11_ii}.
Let us denote
\[\ell(x)=\bigl( \dr(x_\circ)+\nabla  \dr(x_\circ)\cdot (x-x_\circ)\bigr)_+,\]
which satisfies
\[\L(\ell^s)=0\qquad\textrm{in}\quad \{\ell>0\};\]
see Proposition~\ref{prop:onedbarrier} and Lemma~\ref{lem:L1d}.

Now, notice that
\begin{equation}
\label{eq:noticedxcirc} \dr(x_\circ)=\ell(x_\circ)\qquad{\rm and}\qquad \nabla \dr(x_\circ)=\nabla \ell(x_\circ).
\end{equation}
By Lemma \ref{lem13} we have
\begin{equation}\label{poaerf}
\big| \dr(x_\circ+y)-\ell(x_\circ+y)\big|\leq C|y|^{1+\alpha}.
\end{equation}
Using $|a^s-b^s|\leq |a-b|({ a}^{s-1}+{ b}^{s-1})$ for $a,b\geq0$ (with the convention that $a^{s-1} = 0$ if $a = 0$) we find
\begin{equation}\label{bound-1}
\left| \dr^s(x_\circ+y)-\ell^s(x_\circ+y)\right|\leq C|y|^{1+\alpha}\left(\dr^{s-1}(x_\circ+y)+\ell^{s-1}(x_\circ+y)\right).
\end{equation}
(Again, assuming that $\dr^{s-1} = \dr^{s-1}\chi_{\dr > 0}$ and $\ell^{s-1} = \ell^{s-1}\chi_{\ell > 0}$.) On the other hand, since $ \dr\in C^{1,\alpha}(\overline\Omega)$, using Taylor's theorem for H\"older functions (cf. \eqref{eq:tayl_form}) we have
\[
\big|\dr(x) - \dr(x_\circ) - \nabla \dr(x_\circ)\cdot(x-x_\circ)\big|\le C \rho^{1+\alpha}
\]
for all $x\in B_{\rho/2}(x_\circ)$, and so, since $\dr(x)$ is comparable to $\rho$ in $B_{\rho/2}(x_\circ)$ we deduce 
\[\ell(x) \ge \dr(x_\circ) + \nabla \dr(x_\circ)\cdot(x-x_\circ)\ge \dr(x) - C\rho^{1+\alpha}>c{\rho}>0\quad \textrm{in}\quad B_{\rho/2}(x_\circ),\]
provided that $\rho_\circ$ is small and for some $c > 0$ (both depending only on $\Omega$).

We also have 
\[
\bigl| \dr^s-\ell^s\bigr|(x)\leq |\dr - \ell|(x)\left\|\dr^{s-1}+\ell^{s-1}\right\|_{L^\infty(B_{\rho/2}(x_\circ))}\le C \rho^{s-1} |\dr - \ell|(x)
\]
for $x\in B_{\rho/2}(x_\circ)$. We have used here that both $\dr$ and $\ell$ are comparable to $\rho$ in $B_{\rho/2}(x_\circ)$.  Thanks to the control on the second derivatives given by Lemma~\ref{lem:distance} and due to \eqref{eq:noticedxcirc}, we have 
\begin{equation}\label{bound-2}
\bigl| \dr^s-\ell^s\bigr|(x_\circ+y)\leq \|D^2 \dr\|_{L^\infty(B_{\rho/2}(x_\circ))}|y|^2 \rho^{s-1} \leq C\rho^{s+\alpha-2}|y|^2
\end{equation}
for $y\in B_{\rho/2}$. 

It follows then from \eqref{bound-1} and \eqref{bound-2} that
\[\bigl| \dr^s-\ell^s\bigr|(x_\circ+y) \le
\left\{\begin{array}{ll}
C \rho^{s+\alpha-2} |y|^2 & \mbox{in}\quad  B_{\rho/2}\\
C |y|^{1+\alpha}\left(\dr^{s-1}(x_\circ+y)+\ell^{s-1}(x_\circ+y)\right)  & \mbox{in}\quad  B_1\setminus B_{\rho/2} \\
C  |y|^s              &  \mbox{in}\quad \R^n \setminus B_1.
\end{array}\right.\]

Recalling that $\L(\ell^s)(x_\circ)=0$ and thanks to the assumption on the kernel \eqref{eq:boundabovekernel}, we find
\[
\begin{split}
|\L( \dr^s)(x_\circ)| &=  |\L\bigl( \dr^s-\ell^s)(x_\circ)|
\\
&\le C\int_{\R^n} \bigl| \dr^s-\ell^s\bigr|(x_\circ+y) \frac{\Lambda}{|y|^{n+2s}}\,dy
\\
&\le  C\int_{B_{\rho/2}} \rho^{s+\alpha-2}|y|^2\frac{dy}{|y|^{n+2s}}+C\int_{\R^n\setminus B_1}|y|^s \frac{dy}{|y|^{n+2s}}\,+\\
&\quad +C\int_{B_1\setminus B_{\rho/2}} |y|^{1+\alpha}\left(\dr^{s-1}(x_\circ+y)+\ell^{s-1}(x_\circ+y)\right)\frac{dy}{|y|^{n+2s}}
\\
&\le   C(\rho^{\alpha-s}+1)+C\int_{B_1\setminus B_{\rho/2}}\frac{ \dr^{s-1}(x_\circ+y)+\ell^{s-1}(x_\circ+y) }{|y|^{n+2s-1-\alpha}}\,dy.
\end{split}
\]
Thus, using Lemma \ref{lem15} twice (on $\dr$ and on $\ell$, and recalling that $\dr$ is comparable to $d_\Omega$), we find
\[|\L( \dr^s)(x_\circ)|\leq C\rho^{\alpha-s},\]
and \eqref{eq-psi-10} follows (recall that $\rho$ is comparable to $\dr(x_\circ)$).
Thus, case \ref{it:lem11_ii} is proved.

Let us now prove \ref{it:lem11_i}.
In this case,  using $|a^s-b^s|\leq  |a-b|^s$ for $a,b\geq0$, it follows from \eqref{poaerf} that
\[
\left| \dr ^s(x_\circ+y)-\ell^s(x_\circ+y)\right|\leq C|y|^{(1+\alpha)s}\quad\text{for}\quad y \in \R^n,
\]
instead of \eqref{bound-1}.
The bound \eqref{bound-2} remains valid, and hence we have
\[\bigl| \dr ^s-\ell^s\bigr|(x_\circ+y) \le
\begin{cases}
C \rho^{s+\alpha-2} |y|^2 & \mbox{for } y\in B_{\rho/2}\\
C |y|^{(1+\alpha)s} &  \mbox{for } y\in B_1\setminus B_{\rho/2} \\
C  |y|^s             &   \mbox{for } y\in \R^n \setminus B_1.
\end{cases}\]

Recalling that $\L(\ell^s)(x_\circ)=0$ (by Proposition \ref{prop:onedbarrier}), we now find
\[
\begin{split}
|\L(\dr ^s)(x_\circ)| &=  |\L\bigl( \dr^s-\ell^s)(x_\circ)|
\\
&\le C\int_{\R^n} \bigl| \dr^s-\ell^s\bigr|(x_\circ+y) K(dy) 
\\
&\le  C\int_{B_{\rho/2}} \rho^{s+\alpha-2}|y|^2K(dy) +C\int_{\R^n\setminus B_1}|y|^s K(dy) \\
&\quad +C\int_{B_1\setminus B_{\rho/2}} |y|^{(1+\alpha)s}K(dy) 
\\
&\le   C\rho^{\alpha-s}+C+ C\rho^{\alpha s-s} \leq C\rho^{\alpha s-s},
\end{split}
\]
where we have used \eqref{eq:prop_kernel2}-\eqref{eq:prop_kernel3}, and the lemma follows.
\end{proof}

When $\alpha>s$ the proof of Proposition~\ref{lem11}~\ref{it:lem11_ii} gives that 
\[\L(\dr^s)\in L^\infty(\Omega).\]
This is false, however, for operators $\L$ without the assumption \eqref{eq:boundabovekernel}; see~\cite{RS-stable}.
Concerning the higher regularity of $\L(\dr^s)$, it was proved in \cite{AR20} that
\[\left.\begin{array}{c}
\partial\Omega\in C^\beta \\
K|_{\mathbb S^{n-1}} \in C^{2\beta+1}(\mathbb S^{n-1})\\
\beta > 1+s,\quad \beta,\beta\pm s\notin\N
\end{array}\right\}
\qquad\Longrightarrow \qquad
\L(\dr^s)\in C^{\beta-1-s}(\overline\Omega).
\]
This is a crucial ingredient in the proof of higher regularity of the quotient $u/\dr^s$, and of Theorem \ref{thm:u1overu2}.

\subsection*{Sub- and supersolutions}

We next use Proposition \ref{lem11} to construct sub- and supersolutions that are comparable to $\dr^s$.

Moreover, we also want to bound $\L (\dr^{s+\eps})$ in a $C^{1,\alpha}$ domain. 
This will be used to construct sub- and supersolutions in $C^{1,\alpha}$ domains.

\begin{lem}\label{lem12}
Let $s\in (0, 1)$, and let $\L\in \GLh$. 
Let $\alpha\in(0,1)$, let $\Omega$ be any $C^{1,\alpha}$ domain with $C^{1,\alpha}$-radius $\varrho_\circ > 0$.
Then, for any $\eps\in(0,s)$, we have
\[
\L( \dr^{s+\eps})\leq -c\dr^{\eps-s}+C\qquad \textrm{in}\quad \Omega\cap B_{1},
\]
with $c > 0$ and $C > 0$ depending only on $n$, $s$, $\eps$, $\alpha$, $\varrho_\circ$, $\lambda$, and $\Lambda$. 
\end{lem}

\begin{proof}
Let $x_\circ\in \Omega$ and $\rho=\dr(x_\circ)$. 
Exactly as in Proposition~\ref{lem11}, one finds that
\[
\left| \dr^{s+\eps}(x_\circ+y)-\ell^{s+\eps}(x_\circ+y)\right|\leq C|y|^{(1+\alpha)(s+\eps)},
\]
and
\[
\bigl| \dr^{s+\eps}-\ell^{s+\eps}\bigr|(x_\circ+y)\leq C\rho^{s+\eps+\alpha-2}|y|^2
\]
for $y\in B_{\rho/2}$.
  Therefore, as in Proposition~\ref{lem11},
\[|\L\bigl( \dr^{s+\eps}-\ell^{s+\eps})(x_\circ)| \le  C(1+\rho^{\alpha (s+\eps)+\eps-s}).\]
Combined with Lemma~\ref{lem:1dseps}, we find
\[\L( \dr^{s+\eps})(x_\circ)\leq -c\rho^{\eps-s}+C(1+\rho^{\alpha(s+\eps)+\eps-s})\leq -\frac{c}{2}\rho^{s-\eps}+C,\quad \textrm{in}\quad \Omega\cap B_{1},
\]
as wanted.
\end{proof}

As a corollary we obtain the following supersolution in $C^{1,\alpha}$ domains near the boundary:

\begin{cor}[Supersolution]
\label{cor:supersol_domains}
Let $s\in (0, 1)$, and let $\L\in \GLh$. Let $\alpha\in(0,1)$, let $\Omega$ be any bounded $C^{1,\alpha}$ domain with $C^{1,\alpha}$-radius $\varrho_\circ > 0$.

Then, for any $\eps > 0$ there exists $\varphi\in C(\R^n)\cap C^\infty(\Omega)\cap H^s(\R^n)$ such that $\varphi = 0$ in $\R^n\setminus\Omega$ and
\[
\left\{
\begin{array}{rcccll}
&&\L \varphi & \ge & d_\Omega^{\eps-s}& \quad\text{in}\quad \{0 < d_\Omega(x) < \delta\}\\
d_\Omega^s & \le&  \varphi& \le & Cd_\Omega^s & \quad\text{in}\quad \Omega,
\end{array}
\right.
\]
for some $\delta,C > 0$ depending only on $n$, $s$, $\eps$, $\alpha$, $\varrho_\circ$, ${\rm diam}(\Omega)$, $\lambda$, and $\Lambda$. 
\end{cor}

\begin{proof}
Let us assume $\eps < \alpha s$, and consider 
\[
\varphi:= C\dr_\Omega^s - \dr_\Omega^{s+\eps'}
\]
with $\eps' =  \frac12 \eps$ and $c$ small enough.
Then, by Lemma \ref{lem:distHs}, we have $\varphi\in H^s(\R^n)$. 
Moreover, by definition $\varphi$ directly satisfies the second equation if $\delta$ and $c$ are small enough. On the other hand, thanks to  Proposition~\ref{lem11}-\ref{it:lem11_i} and Lemma~\ref{lem12} we know that 
\[
\L \varphi \ge -C\dr_\Omega^{s \alpha-s} + c\dr_\Omega^{\eps'-s}- C.
\]
Since $\eps' < \eps < \alpha s$, and we are done.
\end{proof}

We also get subsolutions in $C^{1,\alpha}$ domains:

\begin{cor}[Subsolution]
\label{cor:subsol_domains}
Let $s\in (0, 1)$, and let $\L\in \GLh$. 
Let $\alpha\in(0,1)$, let $\Omega$ be any bounded $C^{1,\alpha}$ domain with $C^{1,\alpha}$-radius $\varrho_\circ > 0$.

Then, for any $\eps > 0$ there exists $\varphi\in C(\R^n)\cap C^\infty(\Omega)\cap H^s(\R^n)$ such that $\varphi = 0$ in $\R^n\setminus\Omega$ and
\[
\left\{
\begin{array}{rcrcll}
&&\L \varphi & \le & -d_\Omega^{\eps-s}& \quad\text{in}\quad \{0 < d_\Omega(x) < \delta\}\\
d_\Omega^s& \le&  \varphi& \le & C d_\Omega^s  & \quad\text{in}\quad \Omega,
\end{array}
\right.
\]
for some $\delta,C > 0$ depending only on $n$, $s$, $\eps$, $\alpha$, $\varrho_\circ$, ${\rm diam}(\Omega)$,  $\lambda$, and $\Lambda$. 
\end{cor}

\begin{proof}
We proceed similarly to Corollary~\ref{cor:supersol_domains} by considering \[
\varphi:= \dr_\Omega^s + \dr_\Omega^{s+\eps'}
\]
with $\eps' < \eps < \alpha s$, which satisfies 
\[
\L \varphi \le C\dr_\Omega^{s(\alpha - 1)} - c \dr_\Omega^{\eps - s} +C  \le \dr_\Omega^{\eps - s}
\]
if $d_\Omega$ is small enough.  
\end{proof}

In $C^1$ domains we get, instead, the following supersolution. 

\begin{prop}\label{lem11_2}
Let $s\in (0, 1)$, $\alpha\in (0, s)$, and let $\L\in \GLh$.
 Let $\Omega$ be any bounded $C^{1}$ domain, and  let $ \dr = \dr_\Omega$ be given by Definition~\ref{defi:distance}.
Then, 
\[
\L( \dr^\alpha) \geq c\dr^{\alpha-2s}\qquad \textrm{in}\quad \left\{0 < d_\Omega(x) < \delta\right\},
\]
for some constants $c>0$ and $\delta>0$ depending only on $n$, $s$, $\alpha$, $\Omega$, $\lambda$, and $\Lambda$.
\end{prop}
\begin{proof}
The proof follows along the lines of  the proof of Proposition~\ref{lem11}. 

Let $x_\circ\in \Omega$ and $\rho=d_\Omega(x_\circ)$. When $\rho\geq \rho_\circ>0$ then $ \dr^\alpha$ is smooth in a neighborhood of $x_\circ$, and hence $(\L \dr^\alpha)(x_\circ)$ is bounded by a constant depending only on $\rho_\circ$. Thus, we may assume that $\rho\in(0,\rho_\circ)$.

Let us denote
\[\ell(x)=\bigl( \dr(x_\circ)+\nabla  \dr(x_\circ)\cdot (x-x_\circ)\bigr)_+ = a\big(\be\cdot (x-x_\circ + b)\big)_+,\]
with $ a =  |\nabla \dr(x_\circ) |$, $\be  =\frac{\nabla \dr(x_\circ)}{a}$, and $b = \dr(x_\circ) \nabla\dr(x_\circ)$, which by Lemma~\ref{lem:1dseps2}  (and the translation invariance of the operator) satisfies
\begin{equation}
\label{eq:boundellalpha}
\L(\ell^\alpha)(x_\circ) = c a (\be\cdot b)^{\alpha-2s}  = c' \rho^{\alpha -2s}> 0,
\end{equation}
for some $c'>0$ that depends only on $n$, $s$, $\Lambda$, $\lambda$, and $\Omega$.  
We have also used here that 
\[
\frac{1}{C}\le |\nabla \dr(x_\circ)| \le C\quad\text{in}\quad \{0 < d_\Omega(x) < \delta\},
\]
for some $C$ depending only on $\Omega$, which follows form the fact that it is true on $\partial\Omega$ by \eqref{eq:compddr} and $\dr\in C^1(\overline{\Omega})$.

By Lemma \ref{lem13_2} we have
\begin{equation}\label{poaerf2}
\big| \dr(x_\circ+y)-\ell(x_\circ+y)\big|\leq \omega(|y|) |y|,
\end{equation}
and using $|a^\alpha-b^\alpha|\leq |a-b|^\alpha$ for $a,b\geq0$ and $\alpha\in [0, 1]$, we find
\begin{equation}\label{bound-12}
\left| \dr^\alpha(x_\circ+y)-\ell^\alpha(x_\circ+y)\right|\leq C\omega(|y|)^\alpha |y|^\alpha.
\end{equation}

On the other hand, as in the proof of Proposition~\ref{lem11} we can estimate
\begin{equation}\label{bound-22}
\bigl| \dr^\alpha-\ell^\alpha\bigr|(x_\circ+y)\leq  
C\omega(\rho)\rho^{\alpha-2}|y|^2
\end{equation}
for $y\in B_{\rho/2}$. From \eqref{bound-12} and \eqref{bound-22} we have
\[\bigl| \dr^\alpha-\ell^\alpha\bigr|(x_\circ+y) \le
\left\{\begin{array}{ll}
C \omega(\rho) \rho^{\alpha-2} |y|^2 & \mbox{in}\quad  B_{\rho/2}\\
C \omega(|y|)^\alpha |y|^\alpha & \mbox{in}\quad  B_1\setminus B_{\rho/2} \\
C  |y|^\alpha              &  \mbox{in}\quad \R^n \setminus B_1.
\end{array}\right.\]

We now have (thanks to \eqref{eq:boundellalpha}), 
\[
\L (\dr^\alpha)(x_\circ) \ge c'\rho^{\alpha -2s} - I_1-I_2 -I_3, 
\]
where, 
\[
I_1 := \int_{\R^n\setminus B_1} (\dr^\alpha-\ell^\alpha)(x_\circ+y)K(dy) \le \int_{\R^n\setminus B_1} C|y|^\alpha K(dy)  \le C,
\]
using \eqref{eq:prop_kernel3} we get
\[
\begin{split}
I_2 & := \int_{B_1\setminus B_\rho} (\dr^\alpha-\ell^\alpha)(x_\circ+y)K(dy) \\
& \le C \int_{B_1\setminus B_\rho } \omega(|y|)^\alpha |y|^\alpha K(dy) \\
& \le C \int_{B_1\setminus B_{\sqrt{\rho}} } |y|^\alpha K(dy) +  C \omega(\sqrt{\rho})^\alpha \int_{B_{\sqrt{\rho}}\setminus B_{{\rho}} } |y|^\alpha K(dy) \\
& \le C \rho^{\frac{\alpha-2s}{2}} + C \omega(\sqrt{\rho})^\alpha \rho^{\alpha-2s},
\end{split}
\]
and thanks to \eqref{eq:prop_kernel2},
\[
I_3 = \int_{B_\rho} (\dr^\alpha-\ell^\alpha)(x_\circ+y)K(dy) \le C\omega(\rho) \rho^{\alpha-2} \int_{B_\rho} |y|^2 K(dy)  \le C \omega(\rho)\rho^{\alpha-2s}.
\]

Putting all together, we have 
\[
\L (\dr^\alpha)(x_\circ) \ge c'\rho^{\alpha -2s} - C\rho^{\alpha-2s}\left(\rho^{2s-\alpha}+\rho^{\frac{2s-\alpha}{2}} + \omega(\sqrt{\rho})^\alpha +\omega(\rho)\right).
\]
By choosing $\rho$ sufficiently small (depending on $n$, $s$, $\Lambda$, $\lambda$, and $\Omega$), we can ensure 
\[
\L(\dr^\alpha)(x_\circ) \ge \frac{c'}{2}\rho^{\alpha-2s},
\]
as we wanted to see. 
\end{proof}

\section{Barriers in Lipschitz  (and more general) domains}

We now by construct the following subsolution in quite general domains $\Omega\subset\R^n$.
Notice that the assumption on $\Omega$ is satisfied for any bounded Lipschitz domain.

\begin{lem} \label{subsolution-2s-eps}
Let $s\in (0, 1)$, and let $\L\in \LLL$ --- that is, of the form \eqref{eq:Lu1}-\eqref{eq:Kint1}-\eqref{eq:compdef}. 
Let $\Omega\subset \R^n$ be any open set such that, 
\[\textrm{for all}\quad r\in(0,\rho_\circ), \quad z\in\partial\Omega, \quad \textrm{there is a ball}\quad B_{\kappa r}(x_{r,z})\subset \Omega \cap B_r(z),\]
for some $\kappa,\rho_\circ>0$.

Then, there exist $C_1,\gamma_\circ > 0$ depending only on $n$, $s$, $\lambda$, $\Lambda$, $\kappa$, and $\rho_\circ$ such that 
\[
\L (\dr^{2s-\gamma}) \le C_1 - d^{-\gamma}  \quad\text{in}\quad \Omega
\]
for all $\gamma <\gamma_\circ$, where $\dr$ is given by Lemma \ref{lem:distance}.
\end{lem}

\begin{proof}
First notice that, since the function $\dr^{2s-\gamma}$ is smooth inside $\Omega$,   the inequality clearly holds in $\{d>\rho\}$ for some $C_1$ depending on $\rho>0$.
Thus, it suffices to show that there is a $\rho>0$ for which the inequality holds in $\{0<d<\rho\}$.

Let $\rho>0$ to be fixed later, and let $x_\circ\in \Omega$ with $d(x_\circ)=r<\rho$.
Consider the rescaled function
\[w(x):=\frac{\dr^{2s-\gamma}(x_\circ+rx)}{r^{2s-\gamma}}.\]
Since $|D^2\dr|\leq Cd^{-1}$ in $\Omega$, the function $w$ satisfies
\[\|w\|_{C^2(B_{1/2})} \leq C.\]

On the other hand, by the assumption on $\Omega$ (rescaled), we will have that for all $R\in(1,\rho_\circ/r)$ there is a ball $B_{\kappa R}(z_R)\subset B_{R+1} \cap\big(\frac1r(\Omega-x_\circ)\big)$.
In particular, we have
\[w(x) \geq c|x|^{2s-\gamma}\quad \textrm{in} \quad B_{\kappa R/2}(z_R),\]
with $c>0$, where we used that $\dr\geq d$.

Since $\|w\|_{C^2(B_{1/2})}\leq C$ and $w\geq0$ everywhere, we have
\[\tilde \L w(0) \leq C\Lambda - c\lambda\int_{\bigcup_{R\in(1,\rho_\circ/r)} B_{\frac{\kappa}{2} R}(z_R)} |y|^{2s-\gamma} \frac{dy}{|y|^{n+2s}}\]
for every $\tilde \L\in \LLL$.

Now  notice that, since $d(x_\circ)=r$, we have $|z_R|\geq \kappa R-r$, and therefore $B_{\kappa \bar R/2}(z_{\bar R}) \cap B_{R+1} =\varnothing$ if $\frac{\kappa}{2}\bar R \geq R+1+r$.
In particular, $B_{\kappa \bar R/2}(z_{\bar R}) \cap B_{\kappa R/2}(z_R)= \varnothing$ if $\bar R\geq \frac{5}{\kappa}R$.
Since in each of these balls we have
\[\int_{B_{\kappa R/2}(z_R)}|y|^{2s-\gamma} \frac{dy}{|y|^{n+2s}} \geq cR^{-\gamma}, \]
then by taking $\gamma$ and $r$ small enough we have
\[\tilde \L w(0) \leq C\Lambda - c\lambda \int_{\bigcup_{R\in(1,\rho_\circ/r)} B_{\frac{\kappa}{2} R}(z_R)}  \frac{dy}{|y|^{n+\gamma}} \leq - 1.\]
 
Rescaling back to $\dr$, we deduce that 
\[\L (\dr^{2s-\gamma})(x_\circ) < - d^{-\gamma}.\]
Since this holds for any $x_\circ\in \{0<d<\rho\}$ and for every $\L\in \LLL$, the lemma is proved.
\end{proof}

We next want to construct supersolutions in Lipschitz (and more general) domains.
For this, we need the following.

\begin{lem} \label{barrier-chi-Omega}
Let $s\in (0, 1)$ and   let $\L\in \LLL$ --- that is, of the form \eqref{eq:Lu1}-\eqref{eq:Kint1}-\eqref{eq:compdef}. 
Let $\Omega\subset \R^n$ be any open set satisfying
\[
\big|B_r(z)\cap \Omega^c\big| \geq \mu|B_r| \quad
\textrm{for all}\quad z\in\partial \Omega,\quad r>0,
\]
for some $\mu>0$.
Then,
\[
\L \chi_\Omega \ge c_\circ d^{-2s}\quad\text{in}\quad \Omega,
\]
for some $c_\circ> 0$ depending only on $n$, $s$, $\lambda$, and $\mu$.
\end{lem}

\begin{proof}
For any $x_\circ\in \Omega$ we have
\[
\L \chi_\Omega(x_\circ)  =  \int_{\R^n} \chi_{\Omega^c}(x_\circ+y)K(y) \, dy.
\]
By assumption on $\Omega$, it follows that
\[
\big|B_r(x_\circ)\cap \Omega^c\big| \geq \mu|B_r| \quad
\textrm{for all}\quad x_\circ\in\Omega,\quad r\geq  2 d(x_\circ).\]
Thus, denoting $r_\circ= 2d(x_\circ)$, we have
\[
\L\chi_\Omega(x_\circ) \geq \lambda\int_{\R^n\setminus B_{r_\circ}} \chi_{\Omega^c}(x_\circ+y)\frac{dy}{|y|^{n+2s}} \geq \frac{c\mu\lambda}{r_\circ^{2s}},
\]
and the result follows.
\end{proof}

Using the previous lemma, we can now establish the following.

\begin{lem} \label{supersolution-d-eps}
Let $s\in (0, 1)$ and let $\L\in \LLL$ --- that is, of the form \eqref{eq:Lu1}-\eqref{eq:Kint1}-\eqref{eq:compdef}. 
Let $\Omega\subset \R^n$ be any open set satisfying \begin{equation}\label{assumption-Omega-appendix}
\textrm{for all}\quad r>0, \quad z\in\partial\Omega, \quad \textrm{there is a ball}\quad B_{\kappa r}(x_{r,z})\subset \Omega^c \cap B_r(z),
\end{equation}
for some $\kappa>0$.
Then, there exist $c_1,\eps_\circ > 0$ depending only on $n$, $s$, $\lambda$, $\Lambda$, and $\kappa$ such that 
\[
\L (\dr^{\eps}) \ge c_1d^{\eps-2s}  \quad\text{in}\quad \Omega
\]
for all $\eps <\eps_\circ$.
\end{lem}

\begin{proof}
Let $x_\circ\in \Omega$ and $2r_\circ=d(x_\circ)$. 
Since $|D^2\dr|\leq Cd^{-1}$ in $\Omega$ and $\dr\asymp d$,   a direct computation shows that
\[\|\dr^\eps - \dr^\eps(x_\circ)\chi_{\Omega}\|_{C^2(B_{r_\circ})} \leq C\varepsilon r_\circ^{\eps-2}.\]

Thus, we find
\[\big|\L(\dr^\eps-\dr^\eps(x_\circ)\chi_{\Omega})(x_\circ)\big| \leq C\eps r_\circ^{\eps-2s} 
+ \Lambda\int_{\R^n\setminus B_{r_\circ}} \big|\dr^\eps-\dr^\eps(x_\circ)\chi_{\Omega}\big|\frac{dx}{|x-x_\circ|^{n+2s}}. \]

In order to bound the last integral, let us fix $\eta,\gamma >0$ small, and $M>1$ large, to be chosen later, and define the following subsets of $\Omega\setminus B_{r_\circ}(x_\circ)$,
\[A_\circ := \big\{ (1-\eta)\dr^\eps(x_\circ) < \dr^\eps < (1+\eta)\dr^\eps(x_\circ)\big\}\setminus B_{r_\circ}(x_\circ),\]
\[D_\circ := \big\{ 0 < \dr^\eps < (1-\eta)\dr^\eps(x_\circ)\big\} 
\subset \big\{0<\dr < \gamma r_\circ\big\},\]
\[E_\circ := \big\{ \dr^\eps > (1+\eta)\dr^\eps(x_\circ)\big\}
\subset \Omega\setminus B_{Mr_\circ}(x_\circ),\]
where the inclusions hold provided that $\eps>0$ is small enough (depending only on $\eta$, $\gamma$, and $M$). In particular, in the last inclusion we are using that if $x\in E_\circ$, then 
\[
c(1+\eta)^\eps r_\circ < \dr(x) \le C d(x) \le C\left(|x-x_\circ|+ r_\circ\right),
\]
by the triangle inequality, so that taking $\eps > 0$ small enough (depending on $\eta$ and $M$) we can make sure that $E_\circ \subset \Omega\setminus B_{Mr_\circ}(x_\circ)$. 

Thus, we have
\[\int_{A_\circ} \frac{\big|\dr^\eps-\dr^\eps(x_\circ)\chi_{\Omega}\big|}{|x-x_\circ|^{n+2s}}\, dx
\leq \int_{\Omega\setminus B_{r_\circ}(x_\circ)} C\eta r_\circ^\eps\frac{dx}{|x-x_\circ|^{n+2s}} \leq C\eta r_\circ^{\eps-2s},\]
and
\[\int_{E_\circ} \frac{\big|\dr^\eps-\dr^\eps(x_\circ)\chi_{\Omega}\big|}{|x-x_\circ|^{n+2s}}\, dx 
\leq \int_{\Omega\setminus B_{M r_\circ}(x_\circ)} \dr^\eps(x)\frac{dx}{|x-x_\circ|^{n+2s}} \leq C(Mr_\circ)^{\eps-2s}.\]

It only remains to bound the integral over $D_\circ$.
For this, we use the assumption \eqref{assumption-Omega-appendix} on $\Omega$, which by Lemma \ref{GMT-d-eps} implies that
\[\big|\{0<\dr<\gamma\rho\}\cap B_\rho(x_\circ)\big| \leq C\gamma^\theta\rho^n\]
for all $\rho\geq r_\circ$, for some $\theta,C>0$ depending only on $\kappa$.

Using this, and the layer cake representation
\[\begin{split}
\int_U \frac{dx}{|x-x_\circ|^{n+2s}} & = \int_0^\infty \big|U\cap \{|x-x_\circ|^{-n-2s}>t\} \big| dt \\
& = (n+2s)\int_0^\infty \big|U\cap B_\rho(x_\circ)\big|  \frac{d\rho}{\rho^{n+2s+1}},
\end{split}\]
we find
\[\begin{split}
\int_{D_\circ} \frac{\big|\dr^\eps-\dr^\eps(x_\circ)\chi_{\Omega}\big|}{|x-x_\circ|^{n+2s}}\, dx 
& \leq \int_{\{0<\dr<\gamma r_\circ\}} Cr_\circ^\eps \frac{dx}{|x-x_\circ|^{n+2s}} \\
& =  Cr_\circ^\eps \int_{r_\circ}^\infty \big|\{0<\dr<\gamma r_\circ\} \cap B_\rho(x_\circ) \big| \frac{d\rho}{\rho^{n+2s+1}} \\
&\leq Cr_\circ^\eps \int_{r_\circ}^\infty \gamma^\theta \rho^n \frac{d\rho}{\rho^{n+2s+1}} \\
& =  C\gamma^\theta r_\circ^{\eps-2s}.
\end{split}\]

Combining all the previous inequalities, we obtain
\[\big|\L(\dr^\eps-\dr^\eps(x_\circ)\chi_{\Omega})(x_\circ)\big| \leq C\big(\eps+\eta+M^{\eps-2s}+\gamma^\theta\big) r_\circ^{\eps-2s},\]
where can still choose $\eps$, $\eta$, and $\gamma$ as small as we want, and $M$ as large as needed.

On the other hand, since $2r_\circ=d(x_\circ)$, by Lemma \ref{barrier-chi-Omega} we have 
\[\L(\dr^\eps(x_\circ)\chi_{\Omega})(x_\circ) > c_\circ r_\circ^{\eps-2s},\]
for some $c_\circ>0$.
Therefore, taking 
\[\eps+\eta+M^{\eps-2s}+\gamma^\theta \leq {\textstyle \frac12}c_\circ,\]
the result follows. 
\end{proof}

The following singular supersolution will also be useful.
Notice that this result is purely nonlocal, and that it is false for $s=1$.

\begin{lem} \label{barrier-singular-eps}
Let $s\in (0, 1)$ and let $\L\in \LLL$ --- that is, of the form \eqref{eq:Lu1}-\eqref{eq:Kint1}-\eqref{eq:compdef}. 
Let $\Omega\subset \R^n$ be any open set satisfying \eqref{assumption-Omega-appendix} for some $\kappa>0$.
Then, there exist $c_1,\eps_\circ > 0$ depending only on $n$, $s$, $\lambda$, $\Lambda$, and $\kappa$ such that 
\[
\L (\dr^{-\eps}\chi_\Omega) \ge c_1d^{-\eps-2s}  \quad\text{in}\quad \Omega
\]
for all $\eps <\eps_\circ$.
\end{lem}

\begin{proof}
The proof is analogous to that of Lemma \ref{supersolution-d-eps}.
Indeed, let $x_\circ\in \Omega$ and $2r_\circ=d(x_\circ)$.
Exactly as before, we have
\[\|\dr^{-\eps} - \dr^{-\eps}(x_\circ)\|_{C^2(B_{r_\circ})} \leq C\varepsilon r_\circ^{-\eps-2}\]
and
\[\begin{split}
\big|\L\big(\dr^{-\eps}\chi_{\Omega}-\dr^{-\eps}(x_\circ)& \chi_{\Omega}\big)(x_\circ)\big| 
\leq C\eps r_\circ^{-\eps-2s}  \\
& + \Lambda\int_{\R^n\setminus B_{r_\circ}} \big|\dr^{-\eps}\chi_{\Omega}-\dr^{-\eps}(x_\circ)\chi_{\Omega}\big|\frac{dx}{|x-x_\circ|^{n+2s}}. 
\end{split}\]

We fix $\eta,\gamma >0$ small, and $M>1$ large, to be chosen later, and define 
\[\tilde A_\circ := \big\{ (1-\eta)\dr^{-\eps}(x_\circ) < \dr^{-\eps} < (1+\eta)\dr^{-\eps}(x_\circ)\big\}\setminus B_{r_\circ}(x_\circ),\]
\[\tilde D_\circ := \big\{ 0 < \dr^{-\eps} < (1-\eta)\dr^{-\eps}(x_\circ)\big\} 
\subset \Omega\setminus B_{Mr_\circ}(x_\circ),\]
\[\tilde E_\circ := \big\{ \dr^{-\eps} > (1+\eta)\dr^{-\eps}(x_\circ)\big\}
\subset \big\{\dr < \gamma r_\circ\big\},\]
 where the inclusions hold provided that $\eps>0$ is small enough (depending only on $\eta$, $\gamma$, and $M$, and using the triangle inequality for $\tilde D_\circ$).

Thus, we have
\[\begin{split}
\int_{\tilde A_\circ} \frac{\big|\dr^{-\eps}-\dr^{-\eps}(x_\circ)\chi_{\Omega}\big|}{|x-x_\circ|^{n+2s}} \, dx
& \leq \int_{\Omega\setminus B_{r_\circ}(x_\circ)} \frac{C\eta r_\circ^{-\eps}\, dx}{|x-x_\circ|^{n+2s}}  \leq C\eta r_\circ^{-\eps-2s},
\end{split}\]
and
\[\begin{split}
\int_{\tilde D_\circ} \frac{\big|\dr^{-\eps}-\dr^{-\eps}(x_\circ)\chi_{\Omega}\big|}{|x-x_\circ|^{n+2s}} \, dx
& \leq \int_{\Omega\setminus B_{M r_\circ}(x_\circ)}\frac{ \dr^{-\eps}(x_\circ)\,dx}{|x-x_\circ|^{n+2s}} \leq C(Mr_\circ)^{-\eps-2s}.
\end{split}\]

To bound the integral over $\tilde E_\circ$ we use that
\[\frac{1}{|x-x_\circ|^{n+2s}} = (n+2s)\int_0^\infty \chi_{B_\rho(x_\circ)}(x) \frac{d\rho}{\rho^{n+2s+1}}\]
and 
\[\int_{\tilde U} \dr^{-\eps}(x)dx 
= \int_0^\infty \big|\tilde U\cap \{\dr^{-\eps}>t\}\big|dt = \eps \int_0^\infty \big|\tilde U\cap \{\dr<r\}\big| \frac{dr}{r^{\eps+1}}\]
to get 
\[\begin{split}
\int_U \dr^{-\eps}(x)& \frac{dx}{|x-x_\circ|^{n+2s}}  
= (n+2s)\int_0^\infty \int_{U \cap B_\rho(x_\circ)} \dr^{-\eps}(x)dx \frac{d\rho}{\rho^{n+2s+1}} \\
& = \eps(n+2s)\int_0^\infty \int_0^\infty \big|U \cap B_\rho(x_\circ) \cap \{\dr<r\}\big| \frac{dr}{r^{\eps+1}} \frac{d\rho}{\rho^{n+2s+1}}.
\end{split}\]
Hence, we find
\[\begin{split}
\int_{\tilde E_\circ} \big|\dr^{-\eps}-& \dr^{-\eps}(x_\circ) \chi_{\Omega}\big|  \frac{dx}{|x-x_\circ|^{n+2s}} 
 \leq \int_{\{0<\dr<\gamma r_\circ\}} \dr^{-\eps}(x) \frac{dx}{|x-x_\circ|^{n+2s}} \\
& =  C\eps \int_{r_\circ}^\infty  \int_0^\infty \big|\{0<\dr<\min\{r,\gamma r_\circ\}\} \cap B_\rho(x_\circ) \big| \frac{dr}{r^{\eps+1}}
 \frac{d\rho}{\rho^{n+2s+1}} .
\end{split}\]
Now, by Lemma \ref{GMT-d-eps}, we have that (since $\rho \ge r_\circ$)
\[\big|\{0<\dr<r\} \cap B_\rho(x_\circ) \big| 
\leq C\left(r/\rho\right)^\theta \rho^n = Cr^\theta \rho^{n-\theta} \qquad \textrm{for all}\quad r\in(0,\rho),\]
and thus
\[\begin{split}
\int_{\tilde E_\circ} \big|\dr^{-\eps}-\dr^{-\eps}(x_\circ)& \chi_{\Omega}\big|  \frac{dx}{|x-x_\circ|^{n+2s}} \leq 
  \\  
 & \leq  C\eps \int_{r_\circ}^\infty  \int_0^\infty \min\{r^\theta,(\gamma r_\circ)^\theta\}\rho^{n-\theta} \frac{dr}{r^{\eps+1}}
 \frac{d\rho}{\rho^{n+2s+1}} \\
 & = C(\gamma r_\circ)^{\theta-\eps} r_\circ^{-\theta-2s} \\
 & = C\gamma^{\theta-\eps} r_\circ^{-\eps-2s},
\end{split}\]
provided that $\eps<\theta$.

Combining all the previous inequalities, we find
\[\big|\L\big(\dr^{-\eps}-\dr^{-\eps}(x_\circ)\chi_{\Omega}\big)(x_\circ)\big| \leq C\big(\eps+\eta+M^{-\eps-2s}+\gamma^{\theta-\eps}\big) r_\circ^{-\eps-2s},\]
where can still choose $\eps$, $\eta$, $\gamma$ as small as we want, and $M$ as large as needed.

On the other hand, since $2r_\circ=d(x_\circ)$, by Lemma \ref{barrier-chi-Omega} we have 
\[\L\big(\dr^{-\eps}(x_\circ)\chi_{\Omega}\big)(x_\circ) > c_\circ r_\circ^{-\eps-2s},\]
for some $c_\circ>0$.
Therefore, taking $\eps$, $\eta$, $\gamma$ small enough, and $M$ large enough, the result follows.
\end{proof}

Finally, we show the GMT lemma that we used before.
 
\begin{lem} \label{GMT-d-eps}
Let $\Omega\subset \R^n$ be any open set satisfying \eqref{assumption-Omega-appendix} for some $\kappa>0$, and let $d(x):={\rm dist}(x,\Omega^c)$.
Then, there exist $\theta,C > 0$ depending only on $n$ and $\kappa$ such that 
\[\big|\{0<d<r\rho\}\cap B_\rho(z)\big| \leq Cr^\theta|B_\rho|\]
for all $z\in\partial\Omega$ and all $r \in(0,1)$, $\rho>0$.
\end{lem}

This result follows from classical results on Geometric Measure Theory; see Remark \ref{GMT-rem} below.
Still, for convenience of the reader, we provide a simple proof here. This proof is due to Riccardo Tione.

\begin{proof}
The result is invariant under rescalings, so it suffices to prove it for $\rho=1$ and $z=0$.
Indeed, up to considering the domain $\frac{1}{\rho}(\Omega-z)$, it suffices to show that
\[\big|\{0<d<r\}\cap B_1\big| \leq Cr^\theta,\]
for some $\theta,C>0$ depending only on $n$ and $\kappa$.

For this, we define for any $R>0$,
\[D_R := \left\{ x\in\R^n: {\rm dist}(x,\partial\Omega)\leq R \right\}.\]
Let $r\in(0,\frac12)$, and consider 
\[\mathcal E := \left\{ B_r(x) : x\in \partial\Omega\cap B_{1+r} \right\}.\]
The balls in $\mathcal E$ form a covering of the compact set $D_{\kappa r/2} \cap \overline{B_1}$, and hence  we can extract a finite subcovering, say $B_r(x_1)$, ..., $B_r(x_N)$.
Thanks to Vitali's covering lemma, \cite{EG}, we can take a subcollection of these balls, which we denote $B_r(z_1)$, ..., $B_r(z_m)$, that are pairwise disjoint and satisfy 
\[D_{\kappa r/2} \cap \overline{B_1} \subset \bigcup_{i=1}^m B_{3r}(z_i).\]

Now, for each $i=1,...,m$, we have the following dichotomy: either 
\begin{equation}\label{dichotomy-GMT}
B_r(z_i)\subset B_1,
\end{equation}
or 
\begin{equation}\label{dichotomy-GMT2}
B_{r}(z_i) \subset B_{1+2r} \setminus B_{1-2r}.
\end{equation}
In case \eqref{dichotomy-GMT} holds, by assumption on $\Omega$ there is a point $y_i$ such that
\[B_{\kappa r}(y_i) \subset \Omega^c \cap B_r(z_i).\]
Observe that $B_{\kappa r}(y_i)$ are pairwise disjoint, since $B_r(z_i)$ are.
Moreover,  they have the following property
\begin{equation}\label{claim-Riccardo}
B_{\kappa r/2}(y_i) \subset \big(D_r \setminus D_{\kappa r/2}\big) \cap B_1.
\end{equation}
Indeed, since $z_i\in \partial\Omega$, we immediately have that $B_{\kappa r/2}(y_i)\subset B_r(z_i)\subset D_r\cap B_1$. Moreover, 
\[B_{\kappa r/2}(y_i) \cap D_{\kappa r/2} = \varnothing,\]
since for any $p\in B_{\kappa r/2}(y_i)$, we have $B_{\kappa r/2}(p)\subset B_{\kappa r}(y_i)\subset \Omega^c$. That is, \eqref{claim-Riccardo} holds whenever \eqref{dichotomy-GMT} holds.

Combining the previous information, we have
\[
\big| D_{\kappa r/2} \cap B_1 \big|  \leq 
3^n \sum_{i=1}^m|B_r(z_i)| \leq 3^n\hspace{-5mm} \sum_{\eqref{dichotomy-GMT}\ {\rm holds}}|B_r(z_i)| + 3^n\hspace{-5mm} \sum_{\eqref{dichotomy-GMT2} \ {\rm holds}}|B_r(z_i)|,
\]
with
\[\sum_{\eqref{dichotomy-GMT}\ {\rm holds}}|B_r(z_i)| 
= \frac{2^n}{\kappa^n} \hspace{-2mm}\sum_{\eqref{claim-Riccardo} \ {\rm holds}}|B_{\kappa r/2}(y_i)| 
\leq \frac{2^n}{\kappa^n} \big|\big(D_r \setminus D_{\kappa r/2}\big) \cap B_1 \big|
\]
and
\[ \sum_{\eqref{dichotomy-GMT2} \ {\rm holds}}|B_r(z_i)| 
\leq \big| B_{1+2r} \setminus B_{1-2r}\big| \leq Cr.
\]
Notice that we used (twice) that the balls $B_r(z_i)$ are pairwise disjoint.

From  the previous inequalities, we find
\[
\big| D_{\kappa r/2} \cap B_1 \big| \leq C\big|\big(D_r \setminus D_{\kappa r/2}\big) \cap B_1 \big| + Cr,
\]
and thus 
\[
\big| D_{\kappa r/2} \cap B_1 \big| \leq \vartheta \big|D_r \cap B_1 \big| + r,
\]
with $\vartheta=\frac{C}{C+1}\in(0,1)$.
Since this holds for every $r\in(0,\frac12)$, a standard argument then shows that 
\[\big|D_r \cap B_1 \big| \leq Cr^\theta,\]
for some $\theta,C>0$, and we are done.
\end{proof}

\begin{rem} \label{GMT-rem}
The statement of Lemma \ref{GMT-d-eps} is invariant under rescalings, and therefore it suffices to prove it in case $\rho=1$.
If we denote $E:=\partial\Omega \cap B_1$, the assumption \eqref{assumption-Omega-appendix} on $\Omega$ implies that the set $E$ is $\kappa$-porous, in the sense that for any $z\in E$ and $r>0$ there exists a ball $B_{\kappa r}(x_{r,z}) \subset E^c \cap B_r(z)$.
It is then well-known that this implies that the Minkowski dimension of $E$ is strictly less than $n$, with a uniform bound on the $(n-\theta)$-dimensional upper Minkowski content, i.e., 
\[\big|\{0<d<r\}\cap E\big| \leq Cr^{\theta}\]
for some $\theta,C>0$ depending only on $n$ and $\kappa$; see \cite[Proof of Theorem~2.1]{KR}, and more precisely, the equation before (2.5) therein.
\end{rem}

\section{Barriers  for general elliptic operators}

For general kernels in the class $\GL$, we need a stronger assumption on the domain $\Omega$ in order to build a supersolution.
Indeed, the following result holds in $C^1$ domains, but it is false in general in Lipschitz domains.

\begin{lem} \label{barrier-chi-Omega-2}
Let $s\in (0, 1)$,   let $\L\in \GL$, and let $\Omega\subset \R^n$ be any bounded $C^1$ domain.
Then,
\[
\L \chi_\Omega \ge c_\circ d^{-2s}\quad\text{in}\quad \Omega,
\]
for some $c_\circ> 0$ depending only on $n$, $s$, $\lambda$, $\Lambda$, and $\Omega$.
\end{lem}

\begin{proof}
Let us define, for any $\nu\in \mathbb S^{n-1}$ and any $\delta>0$, the cone 
\[\mathcal C_\nu^\delta:=\big\{x\cdot \nu \geq \delta|x|\big\}.\]
Since $\Omega$ is a $C^1$ domain, it follows that for any $\delta>0$ there exists $\varrho_\circ>0$ small enough (depending only on $\delta$ and $\Omega$) so that for all $z\in\partial\Omega$
\begin{equation}
\label{eq:deltarcirc}
z+\mathcal C_\nu^\delta \cap B_{\varrho_\circ} \subset \Omega^c,
\end{equation}
where $\nu\in \mathbb S^{n-1}$ is the unit outward normal vector to $\partial\Omega$ at $z$. Let us fix
\begin{equation}
\label{eq:deltachosen}
\delta = \frac18 \left(\frac{\lambda}{\Lambda}\right)^{\frac12},
\end{equation}
and the corresponding $\varrho_\circ$ such that \eqref{eq:deltarcirc} holds at all $z\in \partial\Omega$ for this $\delta$ (so that $\varrho_\circ$ depends only on $\lambda$, $\Lambda$, and $\Omega$). 

For any $x_\circ\in \Omega$, since $\Omega$ is bounded, from the ellipticity condition on $K$ we have that
\[
\L \chi_\Omega(x_\circ)  =  \int_{\R^n} \chi_{\Omega^c}(x_\circ+y)K(dy)\geq c_\circ\lambda >0,
\]
and so the result follows at points far from $\partial\Omega$. Let us now also see that the statement holds at points close to $\partial\Omega$. In particular, let us assume that $d(x_\circ) < \frac{\varrho_\circ}{2M^2}$, for some $M > 1$ to be chosen later. 
 
For any $x_\circ$ close enough to $\partial\Omega$, we have 
\[\L \chi_\Omega(x_\circ)  \geq  \int_{(z-x_\circ)+\mathcal C_\nu^\delta \cap B_{\varrho_\circ}} K(dy),\]
where $z\in \partial\Omega$ is the closest point to $x_\circ$ on $\partial\Omega$. By denoting $\rho=Md(x_\circ)$ and rescaling the previous integral, we find
\[\L \chi_\Omega(x_\circ)  \geq  \rho^{-2s} \int_{\frac{z-x_\circ}{\rho}+\mathcal C_\nu^\delta \cap B_{\varrho_\circ/\rho}} K_\rho(dy),\]
where $K_\rho(dy)=\rho^{2s} K(\rho\, dy)$   satisfies the same ellipticity conditions as $K$.
Since $\varrho_\circ/\rho>2M$ and $|z-x_\circ|=d(x_\circ)=\rho/M$, this yields
\[\L \chi_\Omega(x_\circ)  \geq  \rho^{-2s} \int_{x_\rho+\mathcal C_\nu^\delta \cap B_{2M}} K_\rho(dy),\]
for some $x_\rho\in \partial B_{1/M}$.

Now notice that, for $M>1$ large enough (depending only on $\delta$ and $n$, which in turn depend only on $\lambda$, $\Lambda$, and $n$), we will have 
\[x_\rho+\mathcal C_\nu^\delta \cap B_{2M}(x_\rho)\supset \mathcal C_\nu^{2\delta} \cap B_{M}\setminus B_1,\]
and hence
\[\L \chi_\Omega(x_\circ)  \geq  \rho^{-2s} \int_{\mathcal C_\nu^{2\delta} \cap B_{M}\setminus B_1} K_\rho(dy).\]

On the other hand, by the ellipticity conditions  \eqref{eq:Kellipt_gen_L}-\eqref{eq:Kellipt_gen_l} on $K_\rho$, for any $\theta>0$ we have 
\[0<\lambda \leq \int_{B_2\setminus B_1} |y\cdot\nu|^2 K_\rho(dy) \leq \theta^2\Lambda + 4\int_{\{|y\cdot \nu|>\theta\}\cap B_2\setminus B_1} K_\rho(dy), \]
where we used that $|y\cdot\nu|^2\leq 4$.
This yields (also using the symmetry of the kernel)
\[\int_{\{y\cdot \nu>\theta\}\cap B_2\setminus B_1} K_\rho(dy) = \frac12\int_{\{|y\cdot \nu|>\theta\}\cap B_2\setminus B_1} K_\rho(dy) \geq \frac18(\lambda - \theta^2\Lambda)>0,\]
if $\theta>0$ is small enough.
Thanks to this, since $\mathcal C_\nu^{2\delta} \cap B_{M}\setminus B_1 \supset \{y\cdot \nu>4\delta\}\cap B_2\setminus B_1$, we deduce that 
\[\L \chi_\Omega(x_\circ)  \geq  \rho^{-2s} \int_{\mathcal C_\nu^{2\delta} \cap B_{M}\setminus B_1} K_\rho(dy) \geq \frac{\rho^{-2s}}{8}(\lambda - 16\delta^2\Lambda)>\frac{\lambda}{16}\rho^{-2s},\]
thanks to the choice of $\delta$, \eqref{eq:deltachosen}.
Recalling that $\rho=M d(x_\circ)$, the result follows.
\end{proof}

Using the previous lemma, we can now establish the following.

\begin{lem} \label{supersolution-d-eps-2}
Let $s\in (0, 1)$,   let $\L\in \GL$, and let $\Omega\subset \R^n$ be any bounded $C^1$ domain.
Then, there exist $c_1,\eps_\circ > 0$ depending only on $n$, $s$, $\lambda$, $\Lambda$, and $\Omega$ such that 
\[
\L (\dr^{\eps}) \ge c_1d^{\eps-2s}  \quad\text{in}\quad \Omega
\]
for all $\eps <\eps_\circ$.
\end{lem}

\begin{proof}
The proof is a modification of that of Lemma \ref{supersolution-d-eps}.
Let $x_\circ\in \Omega$ and $2r_\circ=d(x_\circ)$. 
Then, we have
\[\|\dr^\eps - \dr^\eps(x_\circ)\chi_{\Omega}\|_{C^2(B_{r_\circ})} \leq C\varepsilon r_\circ^{\eps-2}\]
and
\[\L\big(\dr^\eps-\dr^\eps(x_\circ)\chi_{\Omega}\big)(x_\circ) \geq -C\eps r_\circ^{\eps-2s} 
- \int_{\R^n\setminus B_{r_\circ}} \big(\dr^\eps-\dr^\eps(x_\circ)\chi_{\Omega}\big) K(-x_\circ+dy). \]

For any  $\eta >0$ small, and $M>1$ large, we define as in the proof of Lemma~\ref{supersolution-d-eps}
\[A_\circ := \big\{ (1-\eta)\dr^\eps(x_\circ) < \dr^\eps < (1+\eta)\dr^\eps(x_\circ)\big\}\setminus B_{r_\circ}(x_\circ),\]
\[D_\circ := \big\{ 0 < \dr^\eps < (1-\eta)\dr^\eps(x_\circ)\big\} 
\subset \big\{0<\dr <   r_\circ\big\},\]
\[E_\circ := \big\{ \dr^\eps > (1+\eta)\dr^\eps(x_\circ)\big\}
\subset \Omega\setminus B_{Mr_\circ}(x_\circ),\]
where the inclusions hold provided that $\eps>0$ is small enough (depending only on $\eta$,  and $M$).

Thus, we have
\[\int_{A_\circ} \big|\dr^\eps-\dr^\eps(x_\circ)\chi_{\Omega}\big| K(-x_\circ+dy)
\leq \int_{\Omega\setminus B_{r_\circ}(x_\circ)} C\eta r_\circ^\eps K(-x_\circ+dy) \leq C\eta r_\circ^{\eps-2s},\]
and
\[\int_{E_\circ} \hspace{-1mm} \big|\dr^\eps-\dr^\eps(x_\circ)\chi_{\Omega}\big| K(-x_\circ+dy)
\leq \int_{\Omega\setminus B_{M r_\circ}(x_\circ)}\hspace{-4mm} \dr^\eps(x) K(-x_\circ+dy) \leq C(Mr_\circ)^{\eps-2s}.\]
To bound the integral over $D_\circ$, simply notice that $\chi_\Omega = 1$  in $D_\circ$ and $\dr \le r_\circ \le d(x_\circ) \le \dr(x_\circ)$ in $D_\circ$  as well, so
\[
\int_{D_\circ} \big(\dr^\eps-\dr^\eps(x_\circ)\chi_{\Omega}\big) K(-x_\circ+dy)  \leq 0.
\]

Combining all the previous inequalities, we obtain
\[\L\big(\dr^\eps-\dr^\eps(x_\circ)\chi_{\Omega}\big)(x_\circ) \geq -C\big(\eps+\eta+M^{\eps-2s}\big) r_\circ^{\eps-2s},\]
where can still choose $\eps$ and $\eta$  as small as we want, and $M$ as large as needed.

On the other hand, since $2r_\circ=d(x_\circ)$, by Lemma \ref{barrier-chi-Omega} we have 
\[\L(\dr^\eps(x_\circ)\chi_{\Omega})(x_\circ) > c_\circ r_\circ^{\eps-2s},\]
for some $c_\circ>0$.
Therefore, taking 
\[\eps+\eta+M^{\eps-2s} \leq {\textstyle \frac12}c_\circ,\]
the result follows.
\end{proof}

\begin{rem}
\label{rem:lemma_ch2_toadd}
In Lemmas~\ref{barrier-chi-Omega-2} and \ref{supersolution-d-eps-2}, we require the domain $\Omega$ to be $C^1$. Observe that, in fact, we only used that it is Lipschitz with a small Lipschitz constant (on a sufficiently small scale), which we fix in \eqref{eq:deltachosen}. Hence, the results in  Lemmas~\ref{barrier-chi-Omega-2} and \ref{supersolution-d-eps-2} are also true in a Lipschitz domain with a universally small Lipschitz constant (depending only on $\lambda$, and $\Lambda$). The same happens with Proposition~\ref{lem11_2}. 
\end{rem}

\addtocontents{toc}{\protect\setcounter{tocdepth}{2}}

%% file: notation.tex
\chapter*{Notation}
\label{notation}  

Let us introduce some of the notation that is used throughout the book. 
\\[0.4cm]
\noindent {\bf Matrix notation.}\\[0.3cm]
\begin{tabular}{ l m{10cm} }
 $A = (a_{ij})_{ij}$ & Matrix with $(i, j)^{\rm th}$ entry denoted by $a_{ij}$.   \\[0.2cm] 
 ${\rm Id}$& Identity matrix.  \\[0.2cm]
 ${\rm tr}\, A$& Trace of the matrix $A$, i.e., ${\rm tr}\, A = a_{11}+\dots+a_{nn}$.  \\[0.15cm]
 ${\rm det}\, A$& Determinant of the matrix $A$. \\[0.2cm]
 $A^\top$& Transpose of the matrix $A$. \\[0.4cm]
\end{tabular}

\noindent {\bf Functional notation.}\\[0.3cm]
\begin{tabular}{ l m{11cm} }
 $u$ & Unless stated otherwise, $u$ denotes a function $u:\R^n\to \R$.   \\[0.2cm] 
  $u^+,u^-$ & Positive and negative part of a function, $u^+ = \max\{u, 0\}$, $u^- = \max\{-u, 0\}$.   \\[0.2cm]
    $\chi_E$ & Characteristic function of the set $E$, i.e., $\chi_E(x) = 1$ for $x\in E$, and $\chi_E(x) = 0$ for $x\notin E$.  \\[0.2cm]
    ${\rm supp}\, u$ & Support of $u$, ${\rm supp}\,u = \overline{\{x : u(x) \neq 0\}}$.          \\[0.2cm]
    $\ave_A$ & Average integral over the set $A$, i.e., $\ave_A f := \frac{1}{|A|}\int_A f$.    \\[0.2cm]
    $f \asymp g$ & Comparable functions: there exists $C > 0$ independent of $u$ and $v$ such that $C^{-1} g \le f \le C g$ in a given domain.  \\[0.2cm]
    $c$ or $C$ & Denotes a small ($c$) or big ($C$) constant, whose dependence is given by the context of the corresponding statement. Its value can change from line to line. 
\end{tabular}
 
\newpage
\noindent {\bf Differential notation.} Let $u:U \to \R$ be a function.\\[0.3cm]
\begin{tabular}{ l m{9cm} }
 $\de_i u, \de_{x_i} u , u_{x_i}$ & Partial derivative in the $i$-th coordinate direction, $\frac{\de u}{\de x_i}$. \\[0.15cm] 
  $\de_e u$ & Derivative in the $e\in \mathbb{S}^{n-1}$ direction. \\[0.15cm] 
  $\nabla u, Du$ & Gradient, $\nabla u = (\de_1 u, \dots, \de_n u)$. \\[0.15cm] 
   $\de_{ij} u, \de_{x_i x_j} u , u_{x_i x_j}$ & Second partial derivatives in the directions $e_i$ and $e_j$, $\frac{\de^2 u}{\de x_i\de x_j }$. \\[0.15cm] 
    $D^2 u$ & Hessian, $D^2 u  = (\de_{ij} u)_{ij} \in \mathcal{M}_n$. \\[0.15cm] 
    $D^k u$ & Higher derivatives forms, $D^k u := (\de_{i_1}\dots\de_{i_k} u)_{i_1,\dots ,i_k}$. \\[0.15cm] 
    $|D^k u(x)|$ & Norm of $D^k u(x)$ (any equivalent norm). \\[0.15cm] 
    $\|D^k u(x)\|_{\mathcal{F}}$ & Norm of $D^k u$, $\| |D^k u|\|_{\mathcal{F}}$.\\[0.15cm] 
    $\Delta u$ & Laplacian of $u$, $\Delta u = \de_{11} u + \dots + \de_{nn} u$.
\end{tabular}
\\[0.4cm]

\noindent {\bf Geometric and sets notation.}\\[0.25cm]
\begin{tabular}{ l  m{11cm} }
 $\N_0$ & $\N\cup\{0\}$.  \\[0.15cm] 
 $\R^n$, $\mathbb{S}^n$ & $n$-dimensional Euclidean space, $n$-sphere.  \\[0.0cm] 
 $e_i\in \mathbb{S}^{n-1}$ & $i^{\rm th}$ element of the base, $e_i = (0,\dots,0, \stackrel{(i)}{1},0,\dots 0)$. \\[0.15cm]  
 $x \in \R^n$& Typical point $x = (x_1, \dots, x_n)$.  \\[0.15cm]
  $|x| $& Modulus of the point $x$, i.e., $|x| = \sqrt{x_1^2 +\dots + x_n^2}$.  \\[0.15cm] 
    $|U| $& $n$-dimensional Lebesgue measure of a set $U\subset\R^n$.  \\[0.15cm] 
 $\R^n_+$& $\{x = (x_1,\dots, x_n)\in \R^n : x_n > 0\}$.  \\[0.15cm]
  $\mathring{U}$& Interior of the set $U\subset \R^n$. \\[0.15cm]
 $\partial U$& Boundary of the set $U\subset \R^n$. \\[0.15cm]
  $U^c$& Complement of the set $U\subset \R^n$, $U^c := \R^n\setminus U$. \\[0.15cm]
 $V\ssubset  U$& The set $V$ is compactly contained in $U$, that is $\overline{V}\subset U$. \\[0.15cm]
  $B_r(x)$& Ball of radius $r$ centered at $x$, i.e., $B_r(x) := \{y \in \R^n : |x-y|< r\}$. \\[0.15cm]
  $x\cdot y$ & For $x, y\in \R^n$, their scalar product, $x\cdot y = x_1y_1+\dots+x_n y_n$. 
\end{tabular}
\\[1cm]

\label{domainnotation}
\noindent {\bf Domains.} We say that $\Omega\subset\R^n$ is a domain if it is an open connected set.\\[0.2cm]
A domain $\Omega$ is said to be $C^{k, \alpha}$ (resp. $C^k$) if $\partial\Omega$ can be written locally as the graph of a $C^{k,\alpha}$ (resp. $C^k$) function.
\\[0.4cm]

\newpage

\noindent {\bf Function spaces.} Let $U\subset\R^n$ be an open set. \\[0.3cm]
\begin{tabular}{ l m{9.6cm} }
 $C(U), C^0(U)$ & Space of continuous functions $u: U\to \R$.   \\[0.2cm] 
  $C(\overline{U}), C^0(\overline{U})$ & Functions $u\in C(U)$  that are continuous up to the boundary.   \\[0.2cm] 
  $C^k(U), C^k(\overline{U})$ & Space of functions $k$ times continuously differentiable (resp. up to the boundary).   \\[0.2cm]
    $C^{k,\alpha}(U)$ & H\"older spaces, see Appendix \ref{app.A}.   \\[0.2cm]
   $C^\infty(U), C^\infty(\overline{U})$  & Set of functions in $C^k(U)$ or $C^k(\overline{U})$ for all $k\ge 1$.\\[0.2cm]
   $C_c(U), C^k_c(U)$  & Set of functions with compact support in $U$.\\[0.2cm]
   $C_0(U), C^k_0(U)$  & Set of functions with $u = 0$ on $\partial U$.\\[0.2cm]
      $  C^{k,\alpha}_{\rm loc}(U)$  & Set of functions in $C^{k,\alpha}(K)$ for any  $K\ssubset U$.\\[0.2cm]
    $L^p(U)$ & $L^p$ space.\\[0.2cm]
    $L^\infty(U)$ & $L^\infty$ space.\\[0.2cm]
        $L^p_{\rm loc}(U), L^\infty_{\rm loc}(U)$ & Set of functions in $L^p(K)$ (resp. $L^\infty(K)$) for any $K\ssubset U$.\\[0.2cm]
    ${\rm esssup}_\Omega u$ & {Essential supremum of $u$ in $\Omega$:  infimum of the essential upper bounds, ${\rm esssup}_\Omega u := \inf\{b > 0 : |\{u > b\}| = 0\}$.}\\[0.2cm]
    $W^{1,p}, W^{1,p}_0$ & Sobolev spaces.\\[0.2cm]
    $H^{1}, H^{1}_0$ & Sobolev spaces with $p = 2$.\\[0.2cm]
        $W^{s,p}, W^{s,p}_0$ & Fractional Sobolev spaces.\\[0.2cm]
    $H^{s}, H^{s}_0$ & Fractional Sobolev spaces with $p = 2$.\\[0.2cm]
    $\|\cdot\|_{\mathcal{F}}$ & Norm in the functional space $\mathcal{F} \in \{C^0, C^k, L^p, \dots\}$, defined when used for the first times.
\end{tabular}
\\[1cm]
\noindent {\bf Measures.}\label{not:measures} Given a measure $\mu(dx)$ in $\R^n$, and for a given $x_\circ\in \R^n$, and $r\in \R\setminus\{0\}$, we denote by $\mu(x_\circ + r\, dx)$ the measure $\tilde \mu(dx)$ such that  $\tilde\mu(B) = \mu(x_\circ + r B)$ for all Borel sets $B\subset \R^n$.

\newpage
\noindent {\bf Specific notation.}\\[0.3cm]
\begin{tabular}{ l m{9cm} }
 $\sqrtl$ & The square root of the Laplacian, see \eqref{eq:square_root_Laplacian}. \\[0.15cm] 
 $\fls$ & The fractional Laplacian, see \eqref{eq:fractional_Laplacian}. \\[0.15cm] 
  $L^1_\omega(\cdot)$ and $\|\cdot\|_{L^1_\omega(\cdot)}$ & See Definition~\ref{defi:L1w}. \\[0.15cm] 
   $L^1_{\omega_s}(\cdot)$ and $\|\cdot\|_{L^1_{\omega_s}(\cdot)}$ & See Definition~\ref{defi:L1omegas}. \\[0.15cm] 
            $L^\infty_{\tau}(\cdot)$ and $\|\cdot\|_{L^\infty_{\tau}(\cdot)}$ & See Definition~\ref{defi:Linfs}. \\[0.15cm] 
   $\GL$ & General elliptic operators, see Definition~\ref{defi:G}. \\[0.15cm] 
      $\GLh$ & General stable operators (i.e., operators in $\GL$ with homogeneous kernels), see Definition~\ref{defi:Gh}. \\[0.15cm] 
      $\G_s(\lambda, \Lambda; \mu)$ & Regular general elliptic operators, see Definition~\ref{defi:Gmu}. \\[0.15cm] 
            $[K]_\mu$ and $[\L]_\mu$ & See \eqref{Cmu-assumption}-\eqref{eq:LKmu}. \\[0.15cm] 
               $\LLL$ & Operators with kernels comparable to the fractional Laplacian, \eqref{eq:Lu1}-\eqref{eq:Kint1}-\eqref{eq:compdef}, see Definition~\ref{defi:LL}. \\[0.15cm]
                              $\LLh$ & Operators in $\LLL$ with homogeneous kernels, see Definition~\ref{defi:LLh}. \\[0.15cm]  
      $\LL_s(\lambda, \Lambda; \mu)$ & Regular operators in $\LLL$, see Definition~\ref{defi:LL}. \\[0.15cm] 
                  $[K]_{C^\mu}$ and $[\L]_{C^\mu}$ & See Definition~\ref{defi:LL}. \\[0.15cm] 
                                 $\III$ & Fully nonlinear integro-differential operators,  see Definition~\ref{defi:II}. \\[0.15cm] 
      $\II_s(\lambda, \Lambda; \mu)$ & Regular fully nonlinear operators in $\III$, see Definition~\ref{defi:II}. \\[0.15cm] 
                  $[\I]_{C^\mu}$ & See Definition~\ref{defi:II}. \\[0.15cm] 
                  $\langle\cdot, \cdot\rangle_K$ and  $\langle\cdot, \cdot\rangle_{K;\Omega}$ & See \eqref{eq:bilin} and \eqref{eq:prod_KO}. \\[0.15cm] 
                                    $d_\Omega(x) = d(x)$ & Distance to the exterior of $\Omega$,  see \eqref{eq:ddef}. \\[0.15cm] 
                                                       $\dr_\Omega(x) =  \dr(x)$ & Regularized distance, see Definition~\ref{defi:distance}. \\[0.15cm]                 
                  $C^{k,\alpha}$-radius   & See Definition~\ref{defi:varrho}. \\[0.15cm] 
                    $C^{0,1}$-radius  & See Definition~\ref{defi:varrho_Lip}. \\[0.15cm] 
                                   $\M^\pm_\LL$ & Extremal operators in the class $\LL$, see \eqref{eq:extremal_def}. \\[0.15cm]   
                                   $\MMpm$ & Extremal operators in the class $\LLL$, see \eqref{eq:Mpexplicit}-\eqref{eq:Mpexplicit2}. \\[0.15cm]         

\end{tabular}


%% file: biblio.tex

\bibliographystyle{amsalpha}

